\pdfoutput=1  

\documentclass[11pt]{article}
\usepackage[margin=.8in,left=.8in]{geometry}
\usepackage{amsmath}
\usepackage{amsfonts}
\usepackage{amssymb}
\usepackage{amsthm}
\usepackage{mathtools}
\usepackage{subcaption}
\usepackage{graphicx}
\usepackage{color}
\usepackage{xcolor}
\usepackage{xfrac}
\usepackage{stmaryrd}
\usepackage{rotating}
\usepackage{hyperref}
\usepackage[all]{xy}

\usepackage{lmodern} 
\usepackage{ulem}

\usepackage[titletoc]{appendix}

\usepackage{tikz}
\usetikzlibrary{decorations.pathmorphing}
\usepackage{tikz-cd}
\usetikzlibrary{calc}
\usetikzlibrary{matrix}
\usetikzlibrary{positioning,arrows,shapes}
\usetikzlibrary{intersections,shapes.arrows, calc}
\usetikzlibrary{patterns}
\usetikzlibrary{mindmap} 
\usetikzlibrary[mindmap] 

\usepackage{mathrsfs} 

\usepackage{wasysym}

\usepackage[safe]{tipa} 

\usepackage{lipsum}
\usepackage{adjustbox}

\definecolor{darkblue}{rgb}{0.05,0.25,0.65}
\definecolor{greenii}{RGB}{20,140,10}
\definecolor{lightgray}{rgb}{0.9,0.9,0.9}

\usepackage{multirow}

\usepackage{stmaryrd}    

\usepackage{enumerate} 

\usepackage{helvet}   

\usepackage{mathptmx}
\usepackage{amsmath}
\usepackage{graphicx}
\DeclareRobustCommand{\coprod}{\mathop{\text{\fakecoprod}}}
\newcommand{\fakecoprod}{%
  \sbox0{$\prod$}%
  \smash{\raisebox{\dimexpr.9625\depth-\dp0}{\scalebox{1}[-1]{$\prod$}}}%
  \vphantom{$\prod$}%
}

\usepackage{floatflt}  
\usepackage{array}     
\newcolumntype{L}[1]{>{\raggedright\let\newline\\\arraybackslash\hspace{0pt}}m{#1}}
\newcolumntype{C}[1]{>{\centering\let\newline\\\arraybackslash\hspace{0pt}}m{#1}}
\newcolumntype{R}[1]{>{\raggedleft\let\newline\\\arraybackslash\hspace{0pt}}m{#1}}

\usepackage[new]{old-arrows}   

\newdir{> }{{}*!/10pt/@{>}}

\usepackage{amssymb}

\def\acts{\raisebox{-1pt}{\;\;\;\begin{rotate}{90} $\curvearrowright$\end{rotate}}\;}

\makeatletter
\newif\if@sup
\newtoks\@sups
\def\append@sup#1{\edef\act{\noexpand\@sups={\the\@sups #1}}\act}%
\def\reset@sup{\@supfalse\@sups={}}%
\def\mk@scripts#1#2{\if #2/ \if@sup ^{\the\@sups}\fi \else%
  \ifx #1_ \if@sup ^{\the\@sups}\reset@sup \fi {}_{#2}%
  \else \append@sup#2 \@suptrue \fi%
  \expandafter\mk@scripts\fi}
\def\tensor#1#2{\reset@sup#1\mk@scripts#2_/}
\def\multiscripts#1#2#3{\reset@sup{}\mk@scripts#1_/#2%
  \reset@sup\mk@scripts#3_/}
\makeatother

\makeatletter
\newbox\slashbox \setbox\slashbox=\hbox{$/$}
\def\itex@pslash#1{\setbox\@tempboxa=\hbox{$#1$}
  \@tempdima=0.5\wd\slashbox \advance\@tempdima 0.5\wd\@tempboxa
  \copy\slashbox \kern-\@tempdima \box\@tempboxa}
\def\slash{\protect\itex@pslash}
\makeatother

\def\clap#1{\hbox to 0pt{\hss#1\hss}}
\def\mathllap{\mathpalette\mathllapinternal}
\def\mathrlap{\mathpalette\mathrlapinternal}
\def\mathclap{\mathpalette\mathclapinternal}
\def\mathllapinternal#1#2{\llap{$\mathsurround=0pt#1{#2}$}}
\def\mathrlapinternal#1#2{\rlap{$\mathsurround=0pt#1{#2}$}}
\def\mathclapinternal#1#2{\clap{$\mathsurround=0pt#1{#2}$}}

\let\oldroot\root
\def\root#1#2{\oldroot #1 \of{#2}}
\renewcommand{\sqrt}[2][]{\oldroot #1 \of{#2}}

\DeclareSymbolFont{symbolsC}{U}{txsyc}{m}{n}
\SetSymbolFont{symbolsC}{bold}{U}{txsyc}{bx}{n}
\DeclareFontSubstitution{U}{txsyc}{m}{n}

\DeclareSymbolFont{stmry}{U}{stmry}{m}{n}
\SetSymbolFont{stmry}{bold}{U}{stmry}{b}{n}

\DeclareFontFamily{OMX}{MnSymbolE}{}
\DeclareSymbolFont{mnomx}{OMX}{MnSymbolE}{m}{n}
\SetSymbolFont{mnomx}{bold}{OMX}{MnSymbolE}{b}{n}
\DeclareFontShape{OMX}{MnSymbolE}{m}{n}{
    <-6>  MnSymbolE5
   <6-7>  MnSymbolE6
   <7-8>  MnSymbolE7
   <8-9>  MnSymbolE8
   <9-10> MnSymbolE9
  <10-12> MnSymbolE10
  <12->   MnSymbolE12}{}

\makeatletter
\def\Decl@Mn@Delim#1#2#3#4{%
  \if\relax\noexpand#1%
    \let#1\undefined
  \fi
  \DeclareMathDelimiter{#1}{#2}{#3}{#4}{#3}{#4}}
\def\Decl@Mn@Open#1#2#3{\Decl@Mn@Delim{#1}{\mathopen}{#2}{#3}}
\def\Decl@Mn@Close#1#2#3{\Decl@Mn@Delim{#1}{\mathclose}{#2}{#3}}
\Decl@Mn@Open{\llangle}{mnomx}{'164}
\Decl@Mn@Close{\rrangle}{mnomx}{'171}
\Decl@Mn@Open{\lmoustache}{mnomx}{'245}
\Decl@Mn@Close{\rmoustache}{mnomx}{'244}
\makeatother

\makeatletter
\DeclareRobustCommand\widecheck[1]{{\mathpalette\@widecheck{#1}}}
\def\@widecheck#1#2{%
    \setbox\z@\hbox{\m@th$#1#2$}%
    \setbox\tw@\hbox{\m@th$#1%
       \widehat{%
          \vrule\@width\z@\@height\ht\z@
          \vrule\@height\z@\@width\wd\z@}$}%
    \dp\tw@-\ht\z@
    \@tempdima\ht\z@ \advance\@tempdima2\ht\tw@ \divide\@tempdima\thr@@
    \setbox\tw@\hbox{%
       \raise\@tempdima\hbox{\scalebox{1}[-1]{\lower\@tempdima\box
\tw@}}}%
    {\ooalign{\box\tw@ \cr \box\z@}}}
\makeatother

\makeatletter
\def\udots{\mathinner{\mkern2mu\raise\p@\hbox{.}
\mkern2mu\raise4\p@\hbox{.}\mkern1mu
\raise7\p@\vbox{\kern7\p@\hbox{.}}\mkern1mu}}
\makeatother

\newcommand{\gt}{>}

\def\1{{\bf 1}}

\def\<{\langle}
\def\>{\rangle}

\renewcommand{\(}{\begin{equation}}
\renewcommand{\)}{\end{equation}}
\newcommand{\bea}{\begin{eqnarray*}}
\newcommand{\eea}{\end{eqnarray*}}

\usepackage{cleveref}

\crefformat{section}{\S#2#1#3} 
\crefformat{subsection}{\S#2#1#3}
\crefformat{subsubsection}{\S#2#1#3}

\theoremstyle{italics}
\newtheorem{theorem}{Theorem}[section]
\newtheorem{lemma}[theorem]{Lemma}
\newtheorem{prop}[theorem]{Proposition}

\theoremstyle{definition}
\newtheorem{defn}[theorem]{Definition}
\newtheorem{notation}[theorem]{Notation}
\newtheorem{example}[theorem]{Example}
\newtheorem{examples}[theorem]{Examples}

\newtheorem{remark}[theorem]{Remark}
\newtheorem{note[theorem]}{Note}

\usepackage{amsfonts}
\usepackage{colortbl}

\renewcommand{\emph}{\textit}

\def\singular{\rotatebox[origin=c]{70}{$<$}}
\def\smooth{\rotatebox[origin=c]{70}{$\subset$}}
\def\orbisingular{\rotatebox[origin=c]{70}{$\prec$}}

\def\singularG{\raisebox{-3pt}{$\singular^{\hspace{-5.7pt}\raisebox{1pt}{\scalebox{.83}{$G$}}}$}}
\def\smoothG{\raisebox{-3pt}{$\smooth^{\hspace{-5.7pt}\raisebox{1pt}{\scalebox{.83}{$G$}}}$}}
\def\orbisingularG{\raisebox{-3pt}{$\orbisingular^{\hspace{-5.7pt}\raisebox{2pt}{\scalebox{.83}{$G$}}}$}}
\def\orbisingularE{\raisebox{-3pt}{$\orbisingular^{\hspace{-5.5pt}\raisebox{1pt}{\scalebox{.83}{$1$}}}$}}
\def\orbisingularGa{\raisebox{-3pt}{$\orbisingular^{\hspace{-5.7pt}\raisebox{2pt}{\scalebox{.83}{$G_{\mathrlap{1}}$}}}$}\,}
\def\orbisingularGb{\raisebox{-3pt}{$\orbisingular^{\hspace{-5.7pt}\raisebox{2pt}{\scalebox{.83}{$G_{\mathrlap{2}}$}}}$}\,}
\def\orbisingularGab{\;\raisebox{-2.4pt}{$\orbisingular^{\hspace{-15.7pt}\raisebox{2pt}{\scalebox{.83}{$G_1 \!\times\!G_2$}}}$}\,}
\def\orbisingularK{\raisebox{-3pt}{$\orbisingular^{\hspace{-5.7pt}\raisebox{2pt}{\scalebox{.83}{\hspace{1pt}$K$}}}$}}
\def\orbisingularHc{\raisebox{-3pt}{$\orbisingular^{\hspace{-5.7pt}\raisebox{2pt}{\scalebox{.83}{$H_{\mathrlap{c}}$}}}$}\,}

\begin{document}

\title{Proper Orbifold Cohomology}

\author{Hisham Sati, \quad Urs Schreiber}

\maketitle

\begin{abstract}
  The concept of \emph{orbifolds} should unify differential geometry
  with equivariant homotopy theory, so that
  orbifold cohomology should unify differential cohomology
  with proper equivariant cohomology theory.
  Despite the prominent role
  that orbifolds have come to play in mathematics and
  mathematical physics, especially in string theory,
  the formulation of a general theory of orbifolds
  reflecting this unification has remained an open problem.
  Here we present a natural theory argued to achieve this.
  We give both a general abstract axiomatization
  in higher topos theory, as well as concrete models
  for ordinary as well as for super-geometric
  and for higher-geometric orbifolds.
  Our first main result is a fully faithful embedding of the
  $2$-category of orbifolds into a singular-cohesive $\infty$-topos
  whose intrinsic cohomology theory is
  proper globally equivariant differential generalized cohomology,
  subsuming traditional orbifold cohomology,
  Chen-Ruan cohomology, and orbifold K-theory.
  Our second main result is a general construction of
  orbifold {\'e}tale cohomology which we show to naturally unify
  {\bf (i)} tangentially twisted cohomology of smooth but curved spaces
  with {\bf (ii)} RO-graded proper equivariant cohomology of flat but singular spaces.
  As a fundamental example we present
  J-twisted orbifold Cohomotopy theories with coefficients in
  shapes of generalized Tate spheres.
  According to {``Hypothesis H''} this includes the
  proper orbifold cohomology theory that controls non-perturbative
  string theory.
\end{abstract}

\tableofcontents

\newpage

\section{Introduction}

\subsection{Motivation}

The concept of \emph{orbifolds}
\cite{Satake56}\cite{Satake57}\cite{Thurston80}\cite{Haefliger84}
-- manifolds with singularities modeled on
fixed points of finite group actions
(review in \cite{MoerdijkMrcun03}\cite[\S 6]{Kapovich08}\cite[\S 4]{BoyerGalicki08}\cite{IKZ10})
--
has become commonplace in mathematics
(e.g. \cite{BLP05}\cite[\S 13]{Ratcliffe06}\cite{JiYau11}),
and plays a central role in theoretical physics
(see \cite{AMR02}), notably so in string theory
(\cite{DHVW85}\cite{DHVW86}\cite{BailinLove99}\cite{SS19a}).
However, the definition of the homotopy theory,
i.e. of the $\infty$-category (\cite{Lurie09}, see \cref{HigherToposTheory})
of orbifolds, hence in particular of \emph{orbifold cohomology},
is subtle, as witnessed by the convoluted history of the concept;
see \cite[Intro.]{Lerman08}\cite[\S 1]{IKZ10}.
In fact, the issue has remained open:

\medskip

\hypertarget{OrbifoldsAsGeometricGroupoids}{}
\noindent {\bf Orbifolds as {\'e}tale stacks?}
A proposal popular among Lie theorists \cite{MoerdijkPronk97}
(see \cite{Moerdijk02}\cite{Lerman08}\cite{Amenta12}) is to
regard an orbifold with local charts
$G_i \!\!\acts U_i$ \eqref{GAction} as
\begin{itemize}
\vspace{-2mm}
\item the {\'e}tale groupoid; specifically: Lie groupoid (see \cite{MoerdijkMrcun03}\cite{TX06})
or topological groupoid (see \cite{CPRST14});
\vspace{-.7cm}
\item equivalently, the {\'etale} geometric stack;
  specifically: differentiable or topological stack
  (\cite{Carchedi11}\cite{Carchedi12}\cite{Ginot13})
\end{itemize}
\vspace{-2mm}
obtained by gluing the corresponding
\emph{homotopy quotient stacks} $U_i \!\sslash\!\! G_i$ \eqref{HomotopyQuotientGeneralizedSpace}.

\vspace{.1cm}

\noindent This proposal is directly modeled (explictly so in \cite[\S 8]{Joyce12}) on the
concept of Deligne-Mumford stacks in
algebraic geometry (\cite{DeligneMumford69}, review in \cite{Kresch09})
and extends to a concept of general {\'e}tale $\infty$-stacks
\cite{Carchedi13}\cite{Carchedi15}. It
relies on the fact that {\'e}tale stacks,
in their role as homotopy-theoretic generalizations of sheaves,
fully capture geometric aspects
(via generalized sheaf cohomology \cite{Brown73}, see \cite{NSS12}),
while in their role as geometric refinements of classifying spaces
they support \emph{Borel equivariant cohomology} (see \cite{Tu11}).
However, Borel cohomology is
coarser than the
\emph{proper equivariant cohomology}
that is generally relevant in theory and in applications:

\medskip

\noindent {\bf Proper equivariant cohomology}\footnote{
  We follow \cite{DHLPS19} with this terminology, see
  Remark \ref{Proper} below.
},
formulated in equivariant homotopy theory
(review in \cite{Blu17}\cite{May96}),
is obtained by
refining the purely homotopy-theoretic nature of Borel cohomology
by the geometric (`cohesive'', see \cref{Results})
nature of fixed loci (see Example \ref{GeometricFixedPointSpacesDiffer}) of topological group actions
-- hence by the characteristic nature of orbifold geometry --
as encoded in the category of orbits of the equivariance group
(recalled in \cref{GEquivariantHomotopyTheory}).
The proper equivariant version of ordinary cohomology is
known as
\emph{Bredon cohomology}
\cite{Bredon67a}\cite{Bredon67b} (review in \cite[\S 1.4]{Blu17}\cite[\S 7]{tomDieck79});
beyond that, there is a wealth of
proper equivariant generalized cohomology theories
(Def. \ref{GenuineGEquivariantCohomology} below)
such as \emph{equivariant K-theory} \cite{Segal68}\cite{AtiyahSegal69}
(which is proper equivariant by \cite[A3.2]{AtiyahSegal04}\cite[A.5]{FHT07}\cite{DavisLuck98})
and \emph{equivariant Cohomotopy theory}
\cite{Segal71}\cite[\S 8]{tomDieck79}\cite{SS19a}\cite{SS19b}.

\medskip
However, if orbifolds are modeled just by {\'e}tale stacks,
then their proper equivariant cohomology remains, by and large, invisible.
This is true even for Chen-Ruan orbifold cohomology:

\medskip

\noindent {\bf Traditional orbifold cohomology and its shortcomings.}
Given an orbifold $\mathcal{X}\!$,
we write (see \cref{Results}) $\smooth \mathcal{X}$ for the
{\'e}tale stack underlying it, and\footnote{
The ``esh''-symbol ``\raisebox{1pt}{\textesh}'' stands for \emph{shape}
\cite[3.4.5]{dcct}\cite[9.7]{Shulman15},
following \cite{Borsuk75},
which for well-behaved topological spaces is another term for
their \emph{homotopy type} \cite[7.1.6]{Lurie09}\cite[4.6]{Wang17};
see Example \ref{SmoothInfinityGroupoids}.}
$\raisebox{1pt}{\textesh}\, \smooth \mathcal{X}$
for its geometric realization or classifying space
(often denoted $B \mathcal{X}$).
In the case that $\mathcal{X}$ is the global quotient orbifold
of a $G$-space $X$, this is the homotopy type of the
\emph{Borel construction}; so that we may generally call
$\raisebox{1pt}{\textesh}\, \smooth \mathcal{X}$
the \emph{Borel space} of the orbifold.
Now, traditional orbifold cohomology is
\cite[p. 38]{ALR07} just the ordinary cohomology
(e.g. singular cohomology)
of this Borel space, hence is \emph{Borel cohomology}:

\vspace{-.2cm}

\begin{equation}
  \label{TraditionalOrbifoldCohomology}
  \overset{
    \mathclap{
    \mbox{
      \tiny
      \color{darkblue}
      \bf
      \begin{tabular}{c}
        traditional
        \\
        orbifold cohomology
      \end{tabular}
    }
    }
  }{
    H^\bullet_{\mathrm{trad}}
    \big(
      \mathcal{X}
      \!
      ,\,
      A
    \big)
  }
  \;:=\;
  \overset{
    \mbox{
      \tiny
      \color{darkblue}
      \bf
      Borel cohomology $\phantom{AAA}$
    }
  }{
  \underset{
    \mathclap{
    \mbox{
      \tiny
      \color{darkblue}
      \bf
      \begin{tabular}{c}
        singular
        \\
        cohomology
      \end{tabular}
    }
    \;\;\;\;
    }
  }{
    H^\bullet_{\mathrm{sing}}
  }
  \big(
    \underset{
      \mathclap{
      \mbox{
        \tiny
        \color{darkblue}
        \bf
        \begin{tabular}{c}
          Borel
          \\
          space
        \end{tabular}
      }
      }
    }{
      \,\raisebox{1pt}{\textesh} \smooth \mathcal{X}
    }
    \!,\,
    A
  \big).
  }
\end{equation}
This can be considered with any kind of coefficients $A$,
notably in the generality of local coefficient systems \cite{MoerdijkPronk99},
but it always remains an invariant of just the Borel space.
Moreover, for a coefficient ring that inverts the order of
the isotropy groups of $\mathcal{X}$, hence in particular
for rational, real and complex number coefficients
$A \in \{\mathbb{Q}, \mathbb{R}, \mathbb{C}\}$,
the purely torsion cohomology of the
orbifold's finite isotropy groups
becomes invisible, and traditional orbifold cohomology reduces
further (e.g. \cite[Prop. 2.12]{ALR07}) to an invariant of just the shape
$\raisebox{1pt}{\textesh}\,\singular \mathcal{X}$
of the singular quotient space $\singular \mathcal{X}$
(the ``coarse moduli space'') underlying the orbifold
(often denoted $\left\vert \mathcal{X}\right\vert$):

\newpage
$\phantom{A}$
\vspace{-.8cm}

\begin{equation}
  \label{SatakeCohomology}
  \overset{
    \mbox{
      \tiny
      \color{darkblue}
      \bf
      \begin{tabular}{c}
       rational
       \\
       orbifold cohomology
      \end{tabular}
    }
  }{
  H^\bullet_{\mathrm{trad}}
  \big(
    \mathcal{X}
    \!,\,
    \mathbb{Q}
  \big)
  }
  \;\simeq\;
  \overset{
    \mbox{
      \tiny
      \color{darkblue}
      \bf
      ordinary cohomology
    }
  }{
  \underset{
    \mathclap{
    \mbox{
      \tiny
      \color{darkblue}
      \bf
      \begin{tabular}{c}
        singular
        \\
        cohomology
      \end{tabular}
    }
    \;\;\;\;
    }
  }{
    H^\bullet_{\mathrm{sing}}
  }
  \big(
    \underset{
      \mathclap{
      \mbox{
        \tiny
        \color{darkblue}
        \bf
        \begin{tabular}{c}
          naive/coarse
          \\
          quotient
          \\
          space
        \end{tabular}
      }
      }
    }{
      \raisebox{1pt}{\textesh}\,\singular\,\mathcal{X}
    }
    \!,\,
    \mathbb{Q}
  \big).
  }
\end{equation}

\vspace{-.1cm}
\noindent It is in this form that orbifold cohomology was
originally introduced in \cite[Thm. 1]{Satake56}
(following \cite{Bai54}, reviewed in \cite[2.1]{ALR07}).

\medskip

Of course it did not go unnoticed that this coarse
notion of orbifold cohomology
is insensitive to the actual nature of orbifolds.
In reaction to this
(and motivated by algebraic
constructions \cite{DHVW85}\cite{DHVW86}
on 2d conformal field theories interpreted
as describing strings propagating on orbifold spacetimes),
Chen and Ruan famously introduced a new orbifold cohomology theory in
\cite{CR04}. But in fact (see \cite[p. 4,7]{Clader14} for review)
Chen-Ruan cohomology of an orbifold
is just Satake's coarse cohomology \eqref{SatakeCohomology}
(typically considered with complex coefficients),
but applied to the corresponding ``inertia orbifold''
$\mathbf{Maps}\big(\raisebox{1pt}{\textesh}S^1, \mathcal{X}\big)$:

\vspace{-.3cm}

\begin{equation}
  \overset{
    \mathclap{
    \mbox{
      \tiny
      \color{darkblue}
      \bf
      \begin{tabular}{c}
        Chen-Ruan
        \\
        orbifold cohomology
      \end{tabular}
    }
    }
  }{
    H^\bullet_{\mathrm{CR}}
    \big(
      \mathcal{X}
    \big)
  }
  \;\;
  \simeq
  \;\;
  \overset{
    \mathclap{
    \mbox{
      \tiny
      \bf
      \color{darkblue}
      \begin{tabular}{c}
        traditional
        orbifold cohomology
      \end{tabular}
    }
    }
  }{
  H^\bullet_{\mathrm{trad}}
  \big(
    \underset{
      \mathclap{
       \mbox{
         \tiny
         \color{darkblue}
         \bf
         \begin{tabular}{c}
           inertia
           orbifold
         \end{tabular}
       }
       }
     }{
       \mathbf{Maps}
       (
         \raisebox{1pt}{\textesh}S^1,
         \mathcal{X}
       )
     }
    \!,\,
    \mathbb{C}
  \big).
  }
\end{equation}

\vspace{-.1cm}

\noindent Still, it turns out that, for global $G$-quotient orbifolds
$\mathcal{X} = \orbisingular( X \!\sslash\! G )$,
Chen-Ruan cohomology is equivalent to
a proper equivariant cohomology theory, namely
to Bredon cohomology with coefficient system given specifically by:

\vspace{-.6cm}

\begin{equation}
  \label{ChenRuanBredonCoefficientSystem}
  A_{\mathrm{CR}}
  \;:\;
  G/H
  \;\longmapsto\;
  \mathrm{ClassFunctions}(H, \mathbb{C})
  \,.
\end{equation}
This was observed in \cite[p. 18]{Moerdijk02}, using
\cite[Thm. 5.5]{Honk90} with \cite[Prop. 6.5 b)]{Honk88}:

\vspace{-.2cm}

\begin{equation}
  \label{ChenRuanCohomologyIsABredonCohomologyTheory}
  \overset{
    \mathclap{
    \mbox{
      \tiny
      \color{darkblue}
      \bf
      \begin{tabular}{c}
        Chen-Ruan
        \\
        orbifold cohomology
      \end{tabular}
    }
    }
  }{
  \mathbf{H}^\bullet_{\mathrm{CR}}
  \big(
    \underset{
      \mathclap{
      \mbox{
        \tiny
        \color{darkblue}
        \bf
        \begin{tabular}{c}
          global quotient
          \\
          orbifold
        \end{tabular}
      }
      }
    }{
      \orbisingular(X \!\sslash\! G)
    }
    \,,\,
    \mathbb{C}
  \big)
  }
  \;\simeq\;
  \overset{
    \mathclap{
    \mbox{
      \tiny
      \color{darkblue}
      \bf
      \begin{tabular}{c}
        Bredon cohomology
      \end{tabular}
    }
    }
  }{
  H_G^\bullet
    \big(
      X
      \,,\,
      \underset{
        \mathclap{
        \mbox{
          \tiny
          \color{darkblue}
          \bf
          \begin{tabular}{c}
            specific system
            \\
            of coefficients
            \eqref{ChenRuanBredonCoefficientSystem}
          \end{tabular}
        }
        }
      }{
        A_{\mathrm{CR}}
      }
    \big).
  }
\end{equation}

\vspace{-1mm}
\noindent Thus the success of Chen-Ruan cohomology
(surveyed in \cite[\S 4,5]{ALR07})
highlights
the relevance of proper equivariance in orbifold cohomology.
At the same time, this means that
to detect the full proper equivariant homotopy type
of orbifolds,
one needs an orbifold cohomology theory that
induces Bredon coefficient systems more general
than \eqref{ChenRuanBredonCoefficientSystem};
and, in fact, one that subsumes also generalized equivariant
cohomology theories such as equivariant K-theory.
In \cite{AdemRuan01} the authors \emph{define}
orbifold K-theory to be the equivariant K-theory of any global quotient
presentation (see also \cite{ARZ06}\cite{BecerraUribe09}\cite{HW11}):

\vspace{-.3cm}

\begin{equation}
  \label{TraditionalOrbifoldKTheory}
  \overset{
    \mathclap{
    \mbox{
      \tiny
      \color{darkblue}
      \bf
      \begin{tabular}{c}
        traditional
        \\
        orbifold K-theory
      \end{tabular}
    }
    }
  }{
    K^\bullet
    \big(
      \underset{
        \mathclap{
        \mbox{
          \tiny
          \color{darkblue}
          \bf
          \begin{tabular}{c}
            global quotient
            \\
            orbifold
          \end{tabular}
        }
        }
      }{
        \orbisingular
        (X \!\sslash\! G )
      }
    \big)
  }
  \;:=\;
  \overset{
    \mathclap{
    \mbox{
      \tiny
      \color{darkblue}
      \bf
      \begin{tabular}{c}
        equivariant
        \\
        K-theory
      \end{tabular}
    }
    }
  }{
    K^\bullet_G
    (
      X
    )\;.
  }
\end{equation}

\vspace{-1mm}
\noindent This works for the case of K-theory, because it has
been proven explicitly
\cite[Prop. 4.1]{PronkScull10}
that the right hand side of \eqref{TraditionalOrbifoldKTheory}
is independent of the choice of global quotient presentation.
However, in general, this approach of circumventing
an intrinsic definition of orbifold cohomology
by just defining it to be equivariant cohomology
of global quotient presentations is,
besides being somewhat unsatisfactory, in need of justification:

\medskip

\noindent {\bf Orbifolds in global equivariant homotopy theory?}
That orbifold cohomology should also capture proper equivariant cohomology
was suggested explicitly in
\cite{PronkScull10}.
However, the fundamental issue remained that
a quotient presentation $\mathcal{X} \simeq \orbisingular(X \!\sslash\! G)$
of an orbifold is not intrinsic to the orbifold, similarly to
a choice of coordinate atlas, while
in equivariant cohomology theory the equivariance group $G$ is
traditionally taken to be fixed.
But this suggests
\cite[Intro.]{Schwede17}\cite[p. ix-x]{Schwede18}
(details in \cite{Juran20})
that the right context for orbifold cohomology is
``global'' equivariant homotopy theory
\cite{Schwede18}
(following \cite{HenriquesGepner07} and originally motivated from patterns
seen in genuine equivariant stable homotopy theory
\cite{Segal71}\cite{LMS86})
where the equivariance group $G$ is allowed to vary in a prescribed
class of groups. On the other hand, plain global homotopy theory retains no
geometric information!

\medskip

\noindent {\bf The open problem} is thus to set up a mathematical theory
of \emph{proper orbifold cohomology} which unifies:
\begin{itemize}
\vspace{-2mm}
\item[{\bf (i)}] the higher \emph{geometric} (differential, {\'e}tale)
aspects of orbifolds
captured by geometric $\infty$-stack theory;
and
\vspace{-3mm}
\item[{\bf (ii)}] the \emph{singular} (equivariant) aspects of orbifolds
captured by proper and global equivariant homotopy theory.
\end{itemize}

\newpage

To achieve this, we look to higher topos theory
\cite{ToenVezzosi05}\cite{Lurie09}\cite{Rezk10}
(more pointers in \cref{HigherToposTheory} below)
as an ambient foundational homotopy theory of higher geometric
spaces \cite{dcct}\cite{Schreiber19}:

\medskip

\noindent {\bf $\infty$-Toposes as collections of generalized higher geometric spaces.}
Viewed from the outside (i.e., ``externally''), an $\infty$-topos is
a collection of geometric spaces of a given flavor, which may be:
\\
{\bf (a)} generalized geometric, $\;\,$ as well as $\;\;$  {\bf (b)} higher geometric.

\vspace{3mm}

\hspace{-.9cm}
\begin{tabular}{ll}

\begin{minipage}[left]{9.3cm}

\noindent {\bf (a)} Here ``{\bf generalized geometry}'' refers to what Grothedieck
called \emph{functorial geometry} \cite{Grothendieck65}
(review in \cite{DemazureGabriel80}),
which he urged in \cite{Grothendieck73} should supercede any point-set
(e.g. locally ringed)
definition of geometric spaces (further amplified by Lawvere, e.g. \cite{Lawvere05}\cite{Lawvere91}). In hindsight, the basic idea here is
just that of how physicists describe emergent spacetimes $\mathcal{X}$
in terms of what (classical) $p$-brane sigma-models on worldvolumes $\Sigma$
detect when probing $\mathcal{X}$ (see \cite{FSS13a}\cite{FSS13b}\cite{JSSW19}):
\end{minipage}

&

\begin{minipage}[left]{8cm}
\vspace{-.4cm}
\begin{equation}
  \label{ExamplesOfGrosSites}
  \begin{aligned}
  &
  \mathrm{Charts}
  \;=\;
  \\
  &
  \;\;\;\;\;
  \left\{
    \begin{array}{ll}
      \emph{CartesianSpaces} & \mbox{(Def. \ref{CartesianSpaces})}
      \\
      \emph{JetsOfCartesianSpaces} & \mbox{(Def. \ref{FormalCartesianSpaces})}
      \\
      \emph{SuperCartesianSpaces} & \mbox{(Def. \ref{SuperFormalCartesianSpaces})}
      \\
      \emph{Singularities} & \mbox{(Def. \ref{CategoryOfSingularities})}
      \\
      \emph{SingularCharts} & \mbox{(Lem. \ref{YonedaOnProductSite})}
      \\
      \ldots
    \end{array}
  \right.
  \end{aligned}
\end{equation}
\end{minipage}
\end{tabular}

\vspace{.2cm}

\noindent Given any category of local model spaces
(often: ``affine spaces'', here: ``charts''; see Def. \ref{ChartsForCohesion} below),
such as those shown in \eqref{ExamplesOfGrosSites},
one may encode a would-be generalized (``target''-)space $\mathcal{X}$ by
assigning to each $\Sigma \in \mathrm{Charts}$ the collection
\begin{equation}
  \label{AbstractProbes}
  \overset{
    \mathclap{
    \raisebox{3pt}{
      \tiny
      \color{darkblue}
      \bf
      probe space
    }
    }
  }{
    \Sigma
  }
   \;\; \xymatrix{
        \;   \ar@{|->}[r]
      &
}
  \overset{
    \mathclap{
    \raisebox{2pt}{
      \tiny
      \color{darkblue}
      \bf
      \begin{tabular}{c}
        collection of probes of
        generalized space $\mathcal{X}$ by $\Sigma$
      \end{tabular}
    }
    }
  }{
  \mathcal{X}(\Sigma)
  \;:=\;
  \big\{
      ``
      \;
      \Sigma \to\mathcal{X}
      \!\text{ ''}
  \big\}
  }
\end{equation}
of geometric (e.g. smooth, super-geometric, etc.) maps into
$\mathcal{X}$; where the quotation marks indicate that, at this point
of bootstrapping $\mathcal{X}$ into existence,
the category in which these maps are actual morphisms is yet to be specified.
To that end, one observes that a minimal set of consistency conditions on
such an abstract assignment \eqref{AbstractProbes} to be
anything like collections of maps into a space $\mathcal{X}$ are:

\vspace{.2cm}

\hspace{-.8cm}
\begin{tabular}{|l||l|}
\hline
$\phantom{A}$
\!\!\!
\begin{minipage}[left]{7.8cm}
\vspace{.1cm}
\hspace{-2mm}
\noindent {\bf (1) Functoriality of probes.}

For every morphism $\phi$ of $\mathrm{Charts}$
there is an operation of ``pre-composition of probe maps by $\phi$'':

\vspace{-.7cm}

\begin{equation}
  \label{IdeaOfFunctorialityOfProbes}
  \raisebox{20pt}{
 \xymatrix@C=10pt{
   \overset{
     \mathclap{
     \raisebox{4pt}{
       \tiny
       \color{darkblue}
       \bf
       \begin{tabular}{c}
         map of
         \\
         probe spaces
       \end{tabular}
     }
     }
   }{
     \Sigma_1
   }
   \ar[d]_-{ \phi }
    & &
   \overset{
     \mathclap{
     \raisebox{4pt}{
       \tiny
       \color{darkblue}
       \bf
       \begin{tabular}{c}
         pre-composition operation
         \\
         on collections of probes
       \end{tabular}
     }
     }
   }{
     \mathcal{X}(\Sigma_1)
   }
   \\
   \Sigma_2
   & &
   \mathcal{X}(\Sigma_2)
   \ar[u]|-{
     \scalebox{.6}{\!\!\! $
       \mathcal{X}(\phi) = \mbox{`` \ $(-)\circ \phi$\ ''}
     $}
   }
 }}
 \;\;
 \begin{array}{c}
   \mbox{such that}
   \\
   \phantom{a}
   \\
   \mathcal{X}(\phi_2) \circ \mathcal{X}(\phi_1)
   \\
   \;\simeq\;
   \\
   \mathcal{X}( \phi_2 \circ \phi_1 )
 \end{array}
\end{equation}

\end{minipage}

&

\begin{minipage}[left]{8.7cm}
\vspace{-.2cm}
\noindent
{\bf  (2) Gluing of probes.}

If $\{ \! \xymatrix@C=10pt{ U_i \ar[r] & \Sigma } \!\}_{i \in I}$
is a cover of $\Sigma \in \mathrm{Charts}$, then probes of
$\mathcal{X}$ by $\Sigma$ should be equivalent to those tuples of probes by the
$U_i$ which are coherently identified on intersections:
\begin{equation}
  \label{IdeaOfGluing}
  \mathcal{X}(\Sigma)
  \;\simeq\;
  \left\{
  \!\!\!\!\!\!
  \mbox{\small
    \begin{tabular}{c}
      tuples of probes $U_i \longrightarrow \mathcal{X}$
      \\
      identified on intersections $U_i \cap U_j$
      \\
      compatibly on $U_i \cap U_j \cap U_k$
      \\
      etc.
    \end{tabular}
  }
  \!\!\!\!\!\!\!
  \right\}
\end{equation}

\end{minipage}
\\
\hline
\end{tabular}

\vspace{.2cm}

\noindent In the jargon of topos theory (see \cite{MacLaneMoerdijk92}\cite{Johnstone02}),
condition \eqref{IdeaOfFunctorialityOfProbes}
says that the collection $\mathcal{X}(-)$ of probes of $\mathcal{X}$
is a \emph{pre-sheaf} on $\mathrm{Charts}$, while
condition \eqref{IdeaOfGluing} says that this is in fact a \emph{sheaf}.
Hence the category of generalized geometric spaces
probeable by $\mathrm{Charts}$ is the category of sheaves
(the \emph{Grothendieck topos}) on $\mathrm{Charts}$:
\begin{equation}
  \label{1CategoryOfGeneralizedSpaces}
  \mathllap{
    \mbox{
      \tiny
      \color{darkblue}
      \bf
      \begin{tabular}{c}
        topos of generalized
        \\
        geometric spaces
      \end{tabular}
    }
    \;\;\;
  }
  \mathrm{GeneralizedSpaces}
  \;:=\;
  \mathrm{Sheaves}
  \big(
    \mathrm{Charts}
  \big)
  \mathrlap{
    \;\;\;
    \mbox{
      \tiny
      \color{darkblue}
      \bf
      \begin{tabular}{c}
        category of sheaves
        \\
        on site of charts
      \end{tabular}
    }
  }
\end{equation}
\noindent Now, every $\Sigma \in \mathrm{Charts}$
is itself canonically regarded as a generalized space
$y(\Sigma) \in \mathrm{GeneralizedSpaces}$, by taking
its probes to be those given by morphisms of $\mathrm{Charts}$
(this is the \emph{Yoneda embedding}\footnote{Shown here for sub-canonical Grothendieck
topologies on $\mathrm{Charts}$, which is the case in all examples of interest here.},
recalled as Prop. \ref{InfinityYonedaEmbedding} below):
\begin{equation}
  \label{ChartsAsGeneralizedSpaces}
  \mathllap{
    \mbox{
      \tiny
      \color{darkblue}
      \bf
      \begin{tabular}{c}
        chart regarded as
        \\
        generalized space
      \end{tabular}
    }
    \;\;\;
  }
  y(\Sigma) \;:\;
  \Sigma'
  \xymatrix{
    \;
    \ar@{|->}[r]
    &
   }
   \big\{
     \Sigma' \to\Sigma
   \big\}
  \;=:\;
  \mathrm{Charts}(\Sigma', \Sigma)
  \mathrlap{
    \mbox{
      \tiny
      \color{darkblue}
      \bf
      \begin{tabular}{c}
        collection of its
        \\
        $\Sigma'$-shaped probes
      \end{tabular}
    }
  }
\end{equation}
Hence we have completed the bootstrap construction of generalized spaces
$\mathcal{X}$ in \eqref{AbstractProbes} if we may remove the quotation marks there,
hence if for
$\mathcal{X} \in \mathrm{GeneralizedSpaces}$ there is a natural equivalence
\begin{equation}
  \label{ChartsYonedaLemma}
  \mathllap{
    \mbox{
      \tiny
      \color{darkblue}
      \bf
      \begin{tabular}{c}
        collection of $\Sigma$-shaped
        \\
        probes of $\mathcal{X}$
      \end{tabular}
    }
    \;\;\;
  }
  \mathcal{X}(\Sigma)
  \;\simeq\;
  \big\{
       y(\Sigma)
    \to
      \mathcal{X}
  \big\}
  \;:=\;
  \mathrm{GeneralizedSpaces}
  \big(
    y(\Sigma), \mathcal{X}
  \big)
  \mathrlap{
    \;\;\;
    \mbox{
      \tiny
      \color{darkblue}
      \bf
      \begin{tabular}{c}
        collection of maps
        \\
        from $y(\Sigma)$ to $\mathcal{X}$
      \end{tabular}
    }
  }
  \,.
\end{equation}
That this is indeed the case is the statement of the
\emph{Yoneda lemma} (recalled as Prop. \ref{YonedaLemma} below),
which thus implies consistency and existence of generalized geometry.

\newpage

\noindent {\bf (b)} On the other hand, ``{\bf higher geometry}''
(see \cite{FSS13a}\cite{FSS19}\cite{JSSW19} for exposition and applications)
refers to
the refinement of the above theory of generalized geometric spaces,
where the collection of probes \eqref{AbstractProbes} of a generalized
space is not necessarily just a set, but may be a set equipped with
equivalences (gauge transformations) between its elements,
and with higher order equivalences (higher gauge transformations) between these,
etc. -- called an \emph{$\infty$-groupoid}
(e.g., modeled as a Kan simplicial set, see \cite[I.3]{GoerssJardine99}).
For example, for $X \in \mathrm{Sets}$ and $G$ a discrete group acting on $S$,
the corresponding \emph{action groupoid} (Example \ref{ActionGroupoids} below)
consists of the elements $x \in X$,
but equipped with an equivalence between $x_1$ and $x_2$
for every group element whose action takes $x_1$ to $x_2$:

\vspace{-.2cm}

\begin{equation}
\scalebox{.9}{
\raisebox{-85pt}{
\begin{tikzpicture}

\draw (0,0) node
{
$
  \begin{aligned}
  \overset{
  }{
    \mathllap{
      \mbox{
        \tiny
        \color{darkblue}
        \bf
        \begin{tabular}{c}
          homotopy
          \\
          quotient
        \end{tabular}
      }
      \;\;
    }
    X \!\sslash\! G
  }
  \;\;\simeq\;\;
  &
  \left\{
  \raisebox{-6pt}{
  \xymatrix{
    y \ar@(ul,ur)^{g_i}
  }
  }
  \;\; \;\;
  \raisebox{46pt}{
  \xymatrix@C=4pt@R=16pt{
    &&
    \scalebox{.9}{$
      g_1 \!\cdot\! x
    $}
    \ar[ddd]|<<<<<<<<{ \scalebox{.7}{$ g_3 \!\cdot\! g_2 $} }
    \ar[ddrr]^-{g_2}
    \\
    \\
    \scalebox{.8}{$x$}
    \ar[uurr]^-{ g_1}
    \ar[drr]_-{\,
      \scalebox{.7}{$
        g_3 \!\cdot\! g_2 \!\cdot\!  g_1
      $}
    }
    \ar[rrrr]|-{\phantom{AA}}^>>>>>>>>{\,
      \scalebox{.65}{$g_2 \!\cdot\! g_1$}
    }
    &&&&
    \scalebox{.8}{$
      g_2 \!\cdot\! g_1 \!\cdot\! x
    $}
    \ar[dll]^{  g_3}
    \\
    &&
    \scalebox{.85}{$ g_3 \!\cdot\! g_2 \!\cdot\! g_1 \!\cdot\! x$}
  }
  }
  \;\;
  \xymatrix@C=12pt{
    z
    \ar@/^1pc/[r]^-{g_1 }
    \ar@{<-}@/_1pc/[r]_-{ g_1^{-1} }
    &
    g_1 \!\cdot\! z
    \ar@/^1pc/[r]^-{ g_2 }
    \ar@{<-}@/_1pc/[r]_-{ g_2^{-1} }
    &
    g_2 \!\cdot\! g_1 \!\cdot\! z
    \ar@/^1pc/[r]^-{ g_3 }
    \ar@{<-}@/_1pc/[r]_-{ g_3^{-1} }
    &
    g_3 \!\cdot\! g_2 \!\cdot\! g_1 \!\cdot\! z
  }
  \;\;\;
  \cdots
  \right\}
  \\
  \\
   \overset{
  }{
    \mathllap{
      \mbox{
        \tiny
        \color{darkblue}
        \bf
        \begin{tabular}{c}
          plain
          \\
          quotient
        \end{tabular}
      }
      \;\;
    }
    X / G
  }
  \;\;\simeq\;\;
  &
  \Big\{
  \;\;\;
  [y]
  \;\;\;\;\;\;\;\;\;\;\;\;\;\;\;\;\;\;\;\;\;\;\;
  [x]
  \;\;\;\;\;\;\;\;\;\;\;\;\;\;\;\;\;\;\;\;\;\;\;\;\;\;\;\;\;\;\;\;\;\;\;\;\;
  \;\;\;\;\;\;\;\;\;\;\;\;\;\;\;\;\;\;\;\;
  [z]
  \;\;\;\;\;\;\;\;\;\;\;\;\;\;\;\;\;
  \;\;\;\;\;\;\;\;\;\;\;\;\;\;\;\;\;\;\;\;
  \cdots
  \Big\}
  \end{aligned}
$
};

\draw[->] (-6.5,.2) to node[left] {\small$\tau_0$} (-6.5,-1.4);

\end{tikzpicture}
}
}
\end{equation}

\vspace{-.1cm}

\noindent This is a model for the \emph{homotopy quotient} of $X$ by $G$,
which resolves the plain quotient $X/G$ (the set of equivalence classes)
by remembering not only \emph{that} but \emph{how} two elements are equivalent.
More precisely, the action groupoid remembers the
\emph{graph} and \emph{syzygies} of the $G$-action,
encoded in its Kan simplicial \emph{nerve} (Example \ref{Nerve} below):

\vspace{-.4cm}

\begin{equation}
  \label{ActionGroupoidAsSimplicialSet}
    \hspace{-1cm}
  X \sslash G
  \;
  \simeq
  \;
  \left(\!\!\!
  \xymatrix{
    \ar@<+0pt>@{..}[rr]
    &&
    \;X \times G \times G\;\;\;
    \ar@<+18pt>[rrrr]|-{\; (x,g_1, g_2) \mapsto (g_1 \!\cdot\! x,\, g_2) \;}
    \ar@<+8pt>@{<-}[rrrr]
    \ar[rrrr]|-{\; (x, g_1, g_2) \mapsto (x,\, g_2 \!\cdot\! g_1) \;}
    \ar@<-8pt>@{<-}[rrrr]
    \ar@<-16pt>[rrrr]|-{\; (x, g_1, g_2) \mapsto (x, g_1) \;}
    &&&&
    \;\;
    \overset{
      \mathclap{
      \raisebox{6pt}{
        \tiny
        \color{darkblue}
        \bf
        \begin{tabular}{c}
          set of
          \\
          morphisms
        \end{tabular}
      }
      }
    }{
      \;X \times G\;
    }
    \;\;
    \ar@<+8pt>[rrr]^-{ (x,g) \mapsto g \!\cdot\! x }
    \ar@{<-}[rrr]|-{\;  (x, e) \mapsfrom x \;}
    \ar@<-8pt>[rrr]_-{ (x,g) \mapsto x }
    &&&
    \;\;\;
    \overset{
      \mathclap{
      \raisebox{6pt}{
        \tiny
        \color{darkblue}
        \bf
        \begin{tabular}{c}
          set of
          \\
          objects
        \end{tabular}
      }
      }
    }{
      \;X
    }
  }
 \; \right)
\end{equation}

\vspace{1mm}
\hspace{-.9cm}
\begin{tabular}{ll}
\begin{minipage}[left]{10.5cm}
In particular, if an element $y \in X$ is fixed by the group action, then
in the homotopy quotient it appears as the one-object groupoid also
known as $K(G,1)$ or (since $G$ is assumed to be discrete here) as
$B G$.
\end{minipage}
&
\hspace{.0cm}
\begin{minipage}[left]{6.5cm}
\begin{equation}
  \label{bgquot}
 \hspace{-.3cm}
  K(G,1) \;\simeq\; B G \;\simeq\;
\Big\{
       \raisebox{-8pt}{
      \xymatrix{
        {}^\ast \ar@(ul,ur)|-{\;g\,}
      }
    }
    \;\;\;
    \vert
    \;\;
    g \in G
  \Big\}
\end{equation}
\end{minipage}
\end{tabular}

\vspace{.2cm}

More generally, if $X \in \mathrm{Charts}$ in the list \eqref{ExamplesOfGrosSites}
is equipped with the action of a discrete group $G$, then
we obtain a higher generalized space
$\mathcal{X} := X \!\sslash\! G$
whose $\infty$-groupoid of $\Sigma$-shaped probes \eqref{AbstractProbes}
is the action groupoid
of the induced action on the set of $\Sigma$-shaped probes of $X$
(the following formula is for contractible charts, Lemma \ref{HommingChartsIntoHomotopyFiberSequences}):

\vspace{-.2cm}

\begin{equation}
  \label{HomotopyQuotientGeneralizedSpace}
  \mathllap{
    \mbox{
      \tiny
      \color{darkblue}
      \bf
      \begin{tabular}{c}
        global quotient
        \\
        orbifold
      \end{tabular}
    }
    \;\;\;
  }
  X \!\sslash\! G
  \;\;:\;\;
  \Sigma
  \;\longmapsto\;
  \big(
    X \!\sslash\! G
  \big)(\Sigma)
  \;:=\;
   X(\Sigma) \!\sslash\! G
   \;=\;
   \mathrm{Charts}(\Sigma,X) \!\sslash\! G \;.
   \mathrlap{
     \;\;\;
     \mbox{
       \tiny
       \color{darkblue}
       \begin{tabular}{c}
         \tiny
         \color{darkblue}
         \bf
         \begin{tabular}{c}
           groupoid of its
           \\
           $\Sigma$-shaped probes
         \end{tabular}
       \end{tabular}
     }
   }
\end{equation}

\vspace{1mm}
\noindent Such a higher generalized space with collections of probes \eqref{AbstractProbes}
being groupoids, and satisfying the appropriate gluing condition \eqref{IdeaOfGluing},
may be called a \emph{$2$-sheaf}
or \emph{sheaf of groupoids} \cite{Brylinski93} on $\mathrm{Charts}$, in
generalization of \eqref{1CategoryOfGeneralizedSpaces}, but is
commonly known as a \emph{stack} \cite{DeligneMumford69}\cite{Giraud72}\cite{Jardine01}\cite{Hollander08},
following \emph{champ} \cite{Giraud71}.
Generally,
a higher generalized space with $\infty$-groupoids of probes
is thus an \emph{$\infty$-sheaf} or  \emph{$\infty$-stack}
on $\mathrm{Charts}$,
in generalization of \eqref{1CategoryOfGeneralizedSpaces}:

\vspace{-.2cm}

\begin{equation}
  \label{InfinityCategoryOfInfinityStacks}
  \overset{
    \mathclap{
    \raisebox{3pt}{
      \tiny
      \color{darkblue}
      \bf
      $\infty$-topos
    }
    }
  }{
    \mathbf{H}
  }
  \;\;:=\;\;
  \overset{
    \raisebox{2pt}{
      \tiny
      \color{darkblue}
      \bf
      \begin{tabular}{c}
        context for
        higher generalized geometry
      \end{tabular}
    }
  }{
    \mathrm{HigherGeneralizedSpaces}
  }
  \;\;:=\;\;
  \overset{
    \mathclap{
    \raisebox{2pt}{
      \tiny
      \color{darkblue}
      \bf
      \hspace{-14pt}
      \begin{tabular}{c}
        $\infty$-category
        of $\infty$-stacks
      \end{tabular}
    }
    }
  }{
    \mathrm{Sheaves}_\infty\big(
      \underset{
        \mathclap{
        \raisebox{-2pt}{
          \tiny
          \color{darkblue}
          \bf
          $\infty$-site of
          probe spaces
        }
        }
      }{
        \mathrm{Charts}
      }
    \big)
  }
  \,.
\end{equation}

\vspace{1mm}
The theory of $\infty$-stacks originates with \cite{Brown73},
developed in \cite{Jardine87}\cite{Jardine96} (survey in \cite{Jardine15}) and brought into
the more abstract form of $\infty$-topos theory in \cite{ToenVezzosi05}\cite{Lurie09}\cite{Rezk10}.
In fact, finitary constructions internal to
$\infty$-toposes behave so well that they may naturally be formulated
\cite{Shulman19} in a kind of programming language now known as
\emph{homotopy type theory} \cite{UFP13}.
While we will not dwell on this here,
we do focus on elegant internal constructions.
For some of these, a
homotopy type-theoretic formulation has already been explored in the
literature:
{\small
\begin{center}
\label{HoTTFormalizations}
\begin{tabular}{|c||ll|l|}
  \hline
  \begin{tabular}{c}
\bf    Theory internal
    \\
\bf    to an $\infty$-topos
  \end{tabular}
  &
  \multicolumn{2}{c|}{
    \begin{tabular}{c}
 \bf     Internal formulation in
      \\
  \bf    traditional mathematics
    \end{tabular}
  }
  &
  \begin{tabular}{c}
  \bf  Partial formulation in
    \\
  \bf  homotopy type theory
  \end{tabular}
  \\
  \hline
  \hline
  Galois theory
  &
  \cref{GaloisTheory}
  &
  \cite{NSS12}
  &
  \cite{BvDR18}
  \\
  \hline
  \begin{tabular}{c}
    modalities \&
    cohesion
  \end{tabular}
  &
  \cref{DifferentialCohomology}
  &
  \cite{SSS09}\cite{dcct}
  &
  \cite{RSS17}\cite{Shulman15}
  \\
  \hline
  {\'e}tale $\infty$-stacks
  &
  \cref{VFolds} & \cite{KS17}
  &
  \cite{Wellen18}\cite{CherubiniRijke20}
  \\
  \hline
  cohomology & \cref{OrbifoldCohomology}
  &
  \cite{SSS09}\cite{NSS12}\cite{dcct}
  &
  \cite{Cavallo15}\cite{BuchholtzHou18}
  \\
  \hline
\end{tabular}
\end{center}
}

\newpage

\noindent {\bf Differential topology in an $\infty$-topos.}
As a consequence of the above, every $\infty$-topos $\mathbf{H}$
behaves like a homotopy theory of generalized geometric spaces.
In order to narrow back in, among these generalized spaces,
on those with some minimum properties,
we may, following \cite{Lawvere91}\cite{Lawvere94}\cite{Lawvere07},
axiomatize qualities of geometric objects
(such as being \emph{discrete}, \emph{smooth},
\emph{\'etale}, \emph{reduced}, \emph{bosonic}, \emph{singular}, etc.)
via the systems of (co-)reflective sub-$\infty$-categories
$\mathbf{H}_{\Circle}, \mathbf{H}_{\Box}, \cdots \subset \mathbf{H}$,
that the objects with these properties (should) form
inside $\mathbf{H}$ \cite{SSS09}\cite{dcct}:
\vspace{-3mm}
\begin{equation}
  \label{AdjointTriples}
  \xymatrix{
    \mathllap{
      \mbox{
        \tiny
        \color{darkblue}
        \bf
        \begin{tabular}{c}
          ambient $\infty$-topos of
          \\
          generalized geometric spaces
        \end{tabular}
        \hspace{.3cm}
      }
    }
    \mathbf{H} \;\;
    \ar@{->}@<+14pt>[rr]^-{ i_! }_-{\bot}
    \ar@{<-^{)}}@<-0pt>[rr]|-{\; i^\ast \;}
    \ar@{->}@<-14pt>[rr]_-{ i_\ast }^-{\bot}
    &&
   \;\; \mathbf{H}_{\Box}
  }
  \phantom{AAAAA}
  \mbox{or}
  \phantom{AAAAA}
  \xymatrix{
    \mathbf{H} \;\;
    \ar@{<-^{)}}@<+14pt>[rr]^-{ i^\ast }_-{\bot}
    \ar@{->}@<-0pt>[rr]|-{\; i_\ast \;}
    \ar@{<-^{)}}@<-14pt>[rr]_-{ i^! }^-{\bot}
    &&
   \;\;
   \mathbf{H}_{\Circle}
   \mathrlap{
     \mbox{
       \hspace{.3cm}
       \color{darkblue}
       \tiny
       \bf
       \begin{tabular}{c}
         sub-$\infty$-category of
         \\
         objects of pure $\Circle$-nature
       \end{tabular}
     }
   }
  }
\end{equation}

\vspace{-3mm}
\noindent This induces systems of adjoint (co-)projection operators
$\Circle \dashv \Box : \mathbf{H} \to \mathbf{H}$,
the associated \emph{idempotent (co-)monads}:
\begin{equation}
  \label{AdjointModalities}
  \Circle
  \;:=\;
  i^\ast \circ i_!
  \;,\;\;
  \Box
  \;:=\;
  i^\ast \circ i_\ast
  \phantom{AAAA}
  \mbox{or}
  \phantom{AAAA}
  \Circle
  \;:=\;
  i^\ast \circ i_\ast
  \;,\;\;
  \Box
  \;:=\;
  i^! \circ i_\ast
  \,,
\end{equation}

\vspace{-1mm}
\noindent to which we refer as
\emph{modal operators} or just \emph{modalities}
\cite{dcct}\cite{RSS17}\cite{Corfield20}.
These are idempotent (Prop. \ref{IdempotentMonads})

\vspace{-.3cm}

\begin{equation}
  \Circle \Circle X
  \;\simeq\;
  \Circle X\;,
  \phantom{AAAA}
  \Box \Box X
  \;\simeq\;
  \Box X
  \,,
\end{equation}
which means that they act like \emph{projecting out}
certain qualitative aspects of generalized spaces,
while them being adjoint means that they project out
an \emph{opposite pair} of such qualities.
Therefore, their (co-)unit transformations
$\eta^{\Box}$ \eqref{AdjunctionUnit}
and
$\epsilon^{{\mbox{\tiny $\Circle$}}}$ \eqref{CounitOfAdjunction}
exhibit every $X \in \mathbf{H}$ as carrying a quality
intermediate to these two opposite extreme
aspects \cite[p. 245]{LawvereRosebrugh03}:
\vspace{-.35cm}
$$
  \xymatrix{
    \underset{
      \mathclap{
      \mbox{
        \tiny
        \color{darkblue}
        \bf
        \begin{tabular}{c}
          pure
          \\
          $\Circle$-aspect
        \end{tabular}
      }
      }
    }{
      \Circle X
    }
    \;\;
    \ar[rr]^-{\epsilon_X^{{\mbox{\tiny $\Circle$}}}}
    &&
    \;\;\;\;\;
    \underset{
      \mathclap{
      \mbox{
        \tiny
        \color{darkblue}
        \bf
        \begin{tabular}{c}
          generalized
          \\
          geometric space
        \end{tabular}
      }
      }
    }{
      X
    }
   \;\;\;\;\;\;\;
    \ar[rr]^-{ \eta^{\Box}_X }
       &&
    \;\;
    \underset{
      \mathclap{
      \mbox{
        \tiny
        \color{darkblue}
        \bf
        \begin{tabular}{c}
          pure
          \\
          $\Box$-aspect
        \end{tabular}
      }
      }
    }{
      \Box X
    }
  }
  \,.
  $$

\vspace{-.25cm}
\noindent For example, any adjoint modality $\flat \dashv \sharp$
(see Def. \ref{CohesiveTopos} below)
that contains the initial modality $\varnothing \dashv \ast$
(which globally projects to the initial and the terminal object, respectively)
acts like projecting out \emph{discrete} and
\emph{purely continuous} (co-discrete, chaotic) aspects of a space.
Consequently,  the existence of such a modality on $\mathbf{H}$
exhibits each space $X \in \mathbf{H}$ as carrying
quality intermediate to these extremes,
hence, in this example, as equipped with a kind of
\emph{topology} (see \cite[\S 3]{Shulman15}, following \cite{Lawvere94}).

\vspace{1mm}
We observe here that extending this basic example to a larger system of adjoint modalities
allows to abstractly encode the presence of differential geometric structure
(Def. \ref{ElasticInfinityTopos} below) and of super-geometric structure
(Def. \ref{SolidTopos} below) in $\infty$-toposes,
and hence on higher generalized spaces.

\medskip

\noindent {\bf Generalized cohomology in an $\infty$-topos.}
Following \cite{SSS09}\cite{NSS12}\cite{dcct},
we may regard the concept of
$\infty$-toposes $\mathbf{H}$ as the ultimate notion of
\emph{generalized cohomology theory},
subsuming and combining all of:

\vspace{-1mm}
\begin{center}
\begin{tabular}{lll}
\rowcolor{lightgray}
\emph{Sheaf hypercohomology}
&
in non-discrete $\infty$-toposes
&
\cite{Brown73}
\\
\emph{Non-abelian cohomology}
&
in general $\infty$-toposes
&
\cite{SSS09}\cite[3]{NSS12}
\\
\rowcolor{lightgray}
\emph{Twisted non-abelian cohomology}
&
in slice $\infty$-toposes
&
Prp. \ref{SliceInfinityTopos}, Rem. \ref{TwistedCohomology}
\\
\emph{Twisted abelian cohomology}
&
in tangent $\infty$-toposes
&
Exl. \ref{TangentInfinityTopos}, Rem. \ref{AbelianTwistedCohomology}
\\
\rowcolor{lightgray}
\emph{Differential cohomology}
&
in cohesive $\infty$-toposes
&
Def. \ref{CohesiveTopos},\hspace{4pt} Rem. \ref{DifferentialCohomologyTheory}
\\
\emph{{\'E}tale cohomology}
&
in elastic $\infty$-toposes
&
Def. \ref{ElasticInfinityTopos}, Def. \ref{EtaleCohomologyOfVFolds}
\\
\rowcolor{lightgray}
\emph{Superspace cohomology}
&
in solid $\infty$-toposes
&
Def. \ref{SolidTopos}, Rem. \ref{SuperCohomologyTheory}
\\
\emph{Proper equivariant cohomology}
&
in singular $\infty$-toposes
&
Def. \ref{SingularCohesiveInfinityTopos},
Rem. \ref{ProperEquivariantCohomologyTheory},
Thm. \ref{OrbifoldCohomologyEquivariant}
\end{tabular}
\end{center}

\vspace{-1mm}

\noindent
In all these cases, for $X, A \in \mathbf{H}$ any two objects, with
$X$ regarded as a domain ``space'' and $A$ as the ``coefficients'', the
\emph{$A$-cohomology of $X$} is embodied by the homomorphisms from $X$ to $A$:
\begin{enumerate}[{\bf (i)}]
 \vspace{-2mm}
  \item a morphism $\xymatrix{X \ar[r]^-{c} & A}$ is a
  \emph{cocycle};
  \vspace{-4mm}
  \item a homotopy
   $
     \xymatrix{
       X
       \ar@/^.8pc/[r]^-{c_1}_-{\ }="s"
       \ar@/_.8pc/[r]_-{c_2}^-{\ }="t"
       &
       A
       \ar@{=>} "s"; "t"
     }
   $
   is a \emph{coboundary};
   \vspace{-3mm}
   \item
     the homotopy groups of the cocycle space

     \vspace{-.5cm}
     \begin{equation}
       \label{IntrinsicCohomologyOfAnInfinityTopos}
       \phantom{AAA}
       H^{-n}(X,A)
       \;:=\;
       \pi_n \, \mathbf{H}(X,A)
       \;\simeq\;
       \pi_0 \, \mathbf{H}(X, \Omega^n A)
     \end{equation}

     \vspace{-2mm}
  \noindent
  are the \emph{cohomology sets} of $X$ with coefficients in $A$.
  (Here $\Omega^n(-)$ is the $n$-fold based loop space.)
\end{enumerate}

    \vspace{-2mm}
\noindent
This is the \emph{intrinsic cohomology theory} of the $\infty$-topos
$\mathbf{H}$ -- we discuss various examples below in \cref{OrbifoldCohomology}.

\newpage

\subsection{Results}
\label{Results}

$\,$
\vspace{-.8cm}

\hspace{-.9cm}
\begin{tabular}{ll}

\begin{minipage}[left]{7cm}

We survey the results presented below:

\medskip

\noindent {\bf Axiomatic orbifold geometry in modal \newline homotopy theory.}

Building on the above, therefore, to formulate proper orbifold cohomology
we ask for $\infty$-toposes \eqref{InfinityCategoryOfInfinityStacks}
equipped with a system of adjoint modalities \eqref{AdjointModalities}
that capture both aspects of proper orbifold cohomology:

\noindent {\bf (i)} the {\it geometric} (differential, {\'e}tale) aspect
and

\noindent {\bf (ii)} the {\it singular} (proper equivariant) aspect.

\end{minipage}

&
\hspace{2cm}
\raisebox{0pt}{
\scalebox{.8}{
\begin{tabular}{|c|c|c|}
  \multicolumn{3}{c}{
    \bf
    Modalities for Singular Super-Geometry (\cref{SingularCohesiveGeometry})
  }
  \\
  \hline
  \multicolumn{3}{|c|}{
    \begin{tabular}{c}
      $\tau_n$
      \\
      \bf
      $n$-groupoidal
    \end{tabular}
  }
  \\
  \hline
  \hline
  \begin{tabular}{c}
    $\mathclap{\phantom{\vert^{\vert^{\vert^{\vert^{\vert}}}}}}$
    $\raisebox{1pt}{\textesh}$
    \\
    \color{darkblue}
    \bf
    shaped
    $\mathclap{\phantom{\vert_{\vert_{\vert_{\vert_{\vert}}}}}}$
  \end{tabular}
  &
  \begin{tabular}{c}
    $\flat$
    \\
    \color{darkblue}
    \bf
    discrete
  \end{tabular}
  &
  \begin{tabular}{c}
    $\sharp$
    \\
    \color{darkblue}
    \bf
    continuous
  \end{tabular}
  \\
  \hline
  \begin{tabular}{c}
    $\mathclap{\phantom{\vert^{\vert^{\vert^{\vert^{\vert}}}}}}$
    $\Re$
    \\
    \color{darkblue}
    \bf
    reduced
    $\mathclap{\phantom{\vert_{\vert_{\vert_{\vert_{\vert}}}}}}$
  \end{tabular}
  &
  \begin{tabular}{c}
    $\Im$
    \\
    \color{darkblue}
    \bf
    {\'e}tale
  \end{tabular}
  &
  \begin{tabular}{c}
    $\mathcal{L}$
    \\
    \color{darkblue}
    \bf
    locally constant
  \end{tabular}
  \\
  \hline
  \begin{tabular}{c}
    $\mathclap{\phantom{\vert^{\vert^{\vert^{\vert^{\vert}}}}}}$
    $\rightrightarrows$
    \\
    \color{darkblue}
    \bf
    even
    $\mathclap{\phantom{\vert_{\vert_{\vert_{\vert_{\vert}}}}}}$
  \end{tabular}
  &
  \begin{tabular}{c}
    $\rightsquigarrow$
    \\
    \color{darkblue}
    \bf
    bosonic
  \end{tabular}
  &
  \begin{tabular}{c}
    $\mathrm{R}\!\mathrm{h}$
    \\
    \color{darkblue}
    \bf
    rheonomic
  \end{tabular}
  \\
  \hline
  \hline
  \begin{tabular}{c}
    $\mathclap{\phantom{\vert^{\vert^{\vert^{\vert^{\vert}}}}}}$
    $\singular$
    \\
    \color{darkblue}
    \bf
    singular
    $\mathclap{\phantom{\vert_{\vert_{\vert_{\vert_{\vert}}}}}}$
  \end{tabular}
  &
  \begin{tabular}{c}
    $\smooth$
    \\
    \color{darkblue}
    \bf
    smooth
  \end{tabular}
  &
  \begin{tabular}{c}
    $\mathclap{\phantom{\vert^{\vert^{\vert^{\vert^{\vert}}}}}}$
    $\orbisingular$
    \\
    \color{darkblue}
    \bf
    orbi-singular
    $\mathclap{\phantom{\vert_{\vert_{\vert_{\vert_{\vert}}}}}}$
  \end{tabular}
  \\
  \hline
\end{tabular}
}
}
\end{tabular}

\hspace{-.9cm}
\begin{tabular}{ll}
\begin{minipage}[left]{7cm}
\noindent {\bf 1. The geometric aspect of orbifold \newline theory.}
In order to formulate, internal to suitable $\infty$-toposes,
the {\bf (a)} differential topology, {\bf (b)} differential geometry,
and {\bf (c)} super-geometry of orbifolds
(hence of manifolds, super-manifolds, super-orbifolds, etc.)
in their smooth guise as {\'e}tale $\infty$-stacks \eqref{InfinityCategoryOfInfinityStacks},
we consider a
corresponding progression of adjoint modalities \eqref{AdjointModalities},
which starts out in the form of the ``axiomatic cohesion''
of \cite{Lawvere07}, on to a second layer that contains a
``de Rham shape'' operation $\Im$ as considered in \cite{Simpson96} \cite{SimpsonTeleman97},
and  then to a third layer which captures super-geometry in a new axiomatic way.
\end{minipage}
&
\hspace{1cm}
\begin{minipage}[left]{8cm}
  \xymatrix@C=1.5em@R=1.5em{
    \scalebox{1.2}{$\mathrm{id}$}
    \ar@{}[r]|-{\dashv}
    \ar@{}[d]|-{\vee}
    &
    \scalebox{1.2}{$\mathrm{id}$}
    \ar@{}[d]|-{\vee}
    \\
    \scalebox{1.2}{$\rightrightarrows$}
    \ar@{}[r]|-{\dashv}
    &
    \scalebox{1.2}{$\rightsquigarrow$}
    \ar@{}[r]|-{\dashv}
    \ar@{}[d]|-{\vee}
    &
    \scalebox{1.2}{$\mathrm{R}\!\mathrm{h}$}
    \ar@{}[d]|-{\vee}
    &&&
    \mbox{
      \tiny
      \color{darkblue}
      \bf
      \hspace{-.4cm}
      \begin{tabular}{c}
        for super-geometry in
        \\
        solid $\infty$-toposes
        (Def. \ref{SolidTopos})
      \end{tabular}
    }
    \\
    &
    \scalebox{1.2}{$\Re$}
    \ar@{}[r]|-{ \dashv }
    &
    \scalebox{1.2}{$\Im$}
    \ar@{}[r]|-{\dashv}
    \ar@{}[d]|-{\vee}
    &
    \scalebox{1.2}{$\mathcal{L}$}
    \ar@{}[d]|-{\vee}
    &&
    \mbox{
      \tiny
      \color{darkblue}
      \bf
      \hspace{-.4cm}
      \begin{tabular}{c}
        for differential geometry in
        \\
        elastic $\infty$-toposes
        (Def. \ref{ElasticInfinityTopos})
      \end{tabular}
    }
    \\
    &
    &
    \scalebox{1.2}{
    \raisebox{1pt}{\textesh}
    }
    \ar@{}[r]|-{\dashv }
    &
    \scalebox{1.2}{$\flat$}
    \ar@{}[r]|-{\dashv }
    \ar@{}[d]|-{\vee}
    &
    \scalebox{1.2}{$\sharp$}
    \ar@{}[d]|-{\vee}
    &
    \mbox{
      \tiny
      \color{darkblue}
      \bf
      \hspace{-.4cm}
      \begin{tabular}{c}
        for differential topology in
        \\
        cohesive $\infty$-toposes
        (Def. \ref{CohesiveTopos})
      \end{tabular}
    }
    \\
    &
    &
    &
    \scalebox{1.2}{$\varnothing$}
    \ar@{}[r]|-{\dashv}
    &
    \scalebox{1.2}{$\ast$}
  }
\end{minipage}
\end{tabular}

\vspace{3mm}
\noindent {\bf 2. The singular aspect of orbifold theory.}

\vspace{1mm}

\hspace{-.95cm}
\hypertarget{FigureD}{}
\begin{tabular}{ll}
  \begin{minipage}[l]{7cm}

Envision the picture of an orbifold singularity $\orbisingular$
and a mathematical magnifying glass held over the singular point.
Under this magnification, one sees resolved the singular point as
a {\it fuzzy fattened point}, to be denoted
$\orbisingularG$.
Removing the magnifying glass,
what one sees with the bare eye depends on how one squints:

\begin{enumerate}[{\bf (i)}]
\vspace{-4mm}
\item The physicists (see, e.g., \cite[\S 1.3]{BailinLove99}) and the
  classical geometers (see, e.g., \cite{IKZ10}\cite{Watts15}) say
  that they see an actual singular point,
  such as the tip of a cone $\singular$.
  This is the {\it plain quotient} \newline
  $
    \singularG := \ast / G = \ast
  $, a point.

\vspace{-2.5mm}
 \item
The higher geometers (see, e.g., \cite{MoerdijkPronk97} \cite{CPRST14})
say that they see the smooth $G$-action around that point,
   hence a smooth stacky geometry $\smooth$.
   This is the
   {\it homotopy quotient}
   \newline
   $
     \smoothG
       \;:=\;
     \ast \!\sslash\! G
       \;=\;
     B G
       \;=\;
     K(G,1)
   $
   \eqref{bgquot}.
   \end{enumerate}

  \end{minipage}
&
  \hspace{0cm}
  \scalebox{0.9}{
\begin{tabular}{c}
\begin{tabular}{c|c}
\bf Singular quotient $\;\;\;\;$
&
\hspace{-2cm}
\bf Smooth homotopy quotient
\\
\\
\includegraphics[width=.23\textwidth]{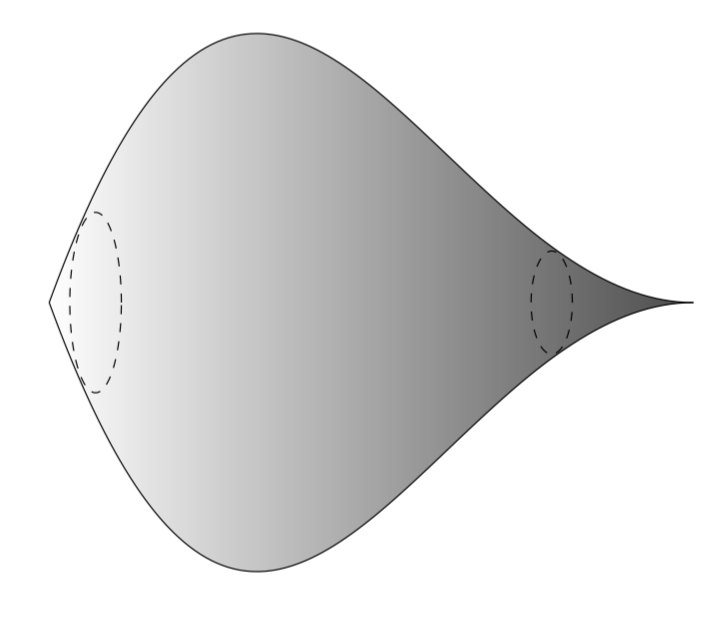}
\hspace{.6cm}
&
\hspace{0.2cm}
\includegraphics[width=.30\textwidth]{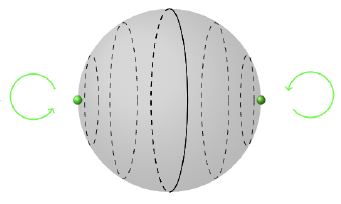}
\hspace{2cm}
\end{tabular}
\\
$$
  \hspace{-3cm}
  \raisebox{36pt}{
  \xymatrix@R=10pt{
    &
    \overset{
      \mbox{
        \color{darkblue}
        \bf
       \fbox{ orbi-singularity}
      }
    }{
      \scalebox{1.4}{$\orbisingularG$}
    }
    \ar@{|->}[dr]^-{
      \overset{
        \begin{rotate}{-32}
          \color{darkblue}
          \tiny
          \bf
          $
          \mathclap{
            \mbox{
              \begin{tabular}{c}
                project onto
                \\
                purely smooth aspect
                \\
                \phantom{$\vert_{\vert_{\vert}}$}
              \end{tabular}
            }
          }
          $
        \end{rotate}
      }{
        \smooth
      }
    }
    \ar@{|->}[dl]_-{
      \overset{
        \begin{rotate}{32}
          \color{darkblue}
          \tiny
          \bf
          $
          \mathclap{
            \mbox{
              \begin{tabular}{c}
                project onto
                \\
                purely singular aspect
                \\
                \phantom{$\vert_{\vert_{\vert}}$}
              \end{tabular}
            }
          }
          $
        \end{rotate}
      }{
        \singular
      }
    }
    \ar@{}[d]|{
      \mbox{\tiny \color{gray} \bf
        \begin{tabular}{c}
          \\
   opposite extreme
          \\
     aspects of orbifold singularities
        \end{tabular}
      }
    }
    \\
    \underset{
      \mathclap{
      \fbox{
        \color{darkblue}
        \bf
        \begin{tabular}{c}
          singular
          \\
          quotient
        \end{tabular}
      }
      }
    }{
      \scalebox{1.2}{$\singularG$} =  \ast / G
      \mathrlap{= \ast}
    }
    &&
    \underset{
      \mathclap{
      \fbox{
        \color{darkblue}
        \bf
        \begin{tabular}{c}
          smooth
          \\
          homotopy quotient
        \end{tabular}
      }
      }
    }{
      \scalebox{1.2}{$\smoothG$} = \ast \!\sslash\! G
      \mathrlap{= B G}
    }
  }
  }
$$
\end{tabular}
}
\end{tabular}

\medskip
\noindent We observe in \cref{SingularGeometry} that just this
is captured
by the cohesive structure on global equivariant homotopy theory that had been
observed in \cite{Rezk14}, but whose conceptual interpretation had remained
open \cite[Footnote 8]{Rezk14}.

\newpage

\noindent {\bf Differential geometry of {\'e}tale $\infty$-stacks.}
We present, in \cref{VFolds}, a general theory of higher differential geometry
formulated internally to these elastic $\infty$-toposes (\cref{DifferentialCohomology}).
This deals with {\'e}tale $\infty$-stacks
locally modeled on any group $\infty$-stack $V$
(``$V$-folds'', Def. \ref{VManifold}).
For the special case $V = (\mathbb{R}^n,+)$,
this  subsumes ordinary manifolds (Example \ref{OrdinaryManifolds})
and ordinary {\'e}tale Lie groupoids (Example \ref{EtaleLieGroupoidAsRnFold}).
For $V$ a super-symmetry group \eqref{EmbeddingOfSuperLieGroupsInGroupsOfJetsOfSuperGroupoids},
this produces a theory of super-orbifolds (Example \ref{GeneralSuperEtaleStacks}),
capturing, for instance,
those that appear as target spaces in
superstring theory (e.g. \cite{ParkReym04}\cite{GaberdielIsraelRabinovic08})
and M-theory \cite{HSS18},
or those that appear as moduli spaces of super-Riemann surfaces
\cite{Rabin87}\cite{LeBrunRothstein88}\cite{Witten12}\cite{CodogniViviani17}.

\vspace{-8mm}
$$
\hspace{-2mm}
  \xymatrix@R=4pt@C=6em{
    \ar@{<..>}@<-5pt>[rrr]|<<<<<<<<<<<<<{\mbox{\color{red} \tiny \bf \;\;coarser\;\;}}
   |>>>>>>>>>>>>>{\mbox{\color{greenii} \tiny \bf \;\;finer\;\;}}
    &&&
    \\
   \fbox{ $(X/G)_{\mathrm{top}}$ }
    \ar@{}[r]|-{ \longmapsfrom }
    &
   \fbox{$ X/G$}
    \ar@{}[r]|-{ \longmapsfrom }
    &
   \fbox{$X \!\sslash\! G$}
    \ar@{}[r]|-{ \longmapsfrom }
    &
  \fbox{$\orbisingular
    \big(
      X \!\sslash\! G
    \big)
    $}
    \\
    {\phantom{a}}
    \\
    \big\{
      \!\!\!\!\!\!\!\!
    \mbox{
      \tiny
    \raisebox{1pt}{  \bf \color{darkblue}
    \begin{tabular}{c}
        underlying
        \\
        topological spaces
      \end{tabular}
      }
    }
      \!\!\!\!\!\!\!
    \big\}
    \ar@{<-}[r]^-{ \scalebox{0.6}{$\mathrm{Dtplg}$} }_-{
      \mbox{
        \tiny
        Prop. \ref{AdjunctionTopDiff}
      }
    }
    &
    \big\{
      \!\!\!\!\!\!\!
    \mbox{
      \tiny
      \raisebox{1pt}{\bf \color{darkblue}
      \begin{tabular}{c}
        orbifolds as
        \\
        diffeological spaces
      \end{tabular}
    }
    }
      \!\!\!\!\!\!\!
    \big\}
    \ar@{<-}@/^1.8pc/[rr]^-{ \scalebox{0.6}{$\mathrm{Snglr}$} }_-{
      \mbox{
        \tiny
        Prop. \ref{OverSingularitiesCohesion},
        Prop. \ref{SingularQuotientOfGOrbiSingularSpaces}
      }
    }
    \ar@{<-}[r]_-{ \tau_0 }
    \ar@{^{(}->}[dd]_-{
      \tiny
    }
    &
    \big\{
      \!\!\!\!\!\!\!
    \mbox{
      \tiny
      \raisebox{1pt}{\bf \color{darkblue}
      \begin{tabular}{c}
        orbifolds as
        \\
        {\'e}tale stacks
      \end{tabular}
    }
    }
      \!\!\!\!\!\!\!
    \big\}
    \ar@{<-}@<+0pt>[r]^-{ \scalebox{0.6}{$\mathrm{Smth}$} }_-{
      \mbox{
        \tiny
        Prop. \ref{OverSingularitiesCohesion},
        Prop. \ref{ShapeOfOrbiSingularSpacesAsPresheafOnSingularities}
      }
    }
    \ar@<-10pt>[r]_-{  }
    \ar@{^{(}->}[dd]_-{
      \mbox{
        \tiny
        Def. \ref{VManifold}
      }
    }
    &
    \big\{
      \!\!\!\!\!\!\!
    \mbox{
      \tiny
      \raisebox{1pt}{\bf \color{darkblue}
      \begin{tabular}{c}
        orbifolds as
        \\
        cohesive orbi-spaces
      \end{tabular}
      }
    }
      \!\!\!\!\!\!\!
    \big\}
    \ar@{^{(}->}[dd]^-{
      \mbox{
        \tiny
        Def. \ref{OrbiVFolds}
      }
    }
    \\
    \\
    &
    \big\{
      \!\!\!\!\!\!\!
    \mbox{
      \tiny
      \raisebox{1pt}{\bf \color{darkblue}
      \begin{tabular}{c}
        sheaves on
        \\
        $\mathrm{CartesianSpaces}$
      \end{tabular}
      }
    }
      \!\!\!\!\!\!\!
    \big\}
    \ar@{^{(}->}[r]
    &
    \big\{
      \!\!\!\!\!\!\!\!
    \mbox{
      \tiny
      \raisebox{1pt}{\bf \color{darkblue}
      \begin{tabular}{c}
        $\infty$-stacks on
        \\
        $\mathrm{CartesianSpaces}$
      \end{tabular}
      }
    }
      \!\!\!\!\!\!\!
    \big\}
    \ar@{^{(}->}[r]^-{ \scalebox{0.6}{$\mathrm{OrbSnglr}$} }_-{
      \mbox{
        \tiny
        Prop. \ref{OverSingularitiesCohesion},
        Def. \ref{OrbiVFolds}
      }
    }
    \ar@{=}[d]
    &
    \big\{
      \!\!\!\!\!\!\!\!
    \mbox{
      \tiny
      \raisebox{1pt}{\bf \color{darkblue}
      \begin{tabular}{c}
        $\infty$-stacks on
        \\
        $\mathrm{SingularCartesianSpaces}$
      \end{tabular}
      }
    }
      \!\!\!\!\!\!\!
    \big\}
    \ar@{=}[d]
    \\
    &
    &
    \underset{
      \mbox{
        \tiny
        Example \ref{SmoothInfinityGroupoids}
      }
    }{
    {\scalebox{.75}{ \bf \color{darkblue}   $
      \mathrm{SmoothGroupoids}_\infty
      $}
    }
    }
    \ar@{^{(}->}[r]
    &
    \underset{
      \mbox{
        \tiny
        Example \ref{OrbiSingularSmoothInfinityGroupoids}
      }
    }{
    {\scalebox{.75}{\bf \color{darkblue}  $
      \mathrm{SingularSmoothGroupoids}_\infty
      $}
    }
    }
  }
$$

\vspace{-.5mm}

\noindent {\bf Gerbes on {\'e}tale $\infty$-stacks.}
With orbifolds, in their incarnation as {\'e}tale stacks, thus embedded into
a fully-fledged $\infty$-topos, the general theory of
$\infty$-bundles \cite{NSS12}\cite{NSS12b} (see \cref{GaloisTheory} below)
applies to provide the theory of
fiber bundles on orbifolds (e.g. \cite{LGTX04}\cite{Seaton06}\cite{BoyerGalicki08})
and of gerbes on orbifolds
(e.g. \cite{LU04}\cite{Carchedi10}\cite{BehrendXu11}\cite{TT14})
naturally generalized to higher, to non-abelian and to twisted gerbes
on orbifolds.

\vspace{1mm}
\noindent {\bf Differential cohomology of {\'e}tale $\infty$-stacks.}
Moreover, since the intrinsic cohomology theory of
cohesive $\infty$-toposes
is differential cohomology (Remark \ref{DifferentialCohomologyTheory}),
this realization of {\'e}tale $\infty$-stacks
within cohesive $\infty$-toposes immediately provides their differential
cohomology theory (see \cite{SSS09}\cite{FSS13a}\cite{FSS15}).
This includes, in particular (as made explicit in \cite{ParkRedden19}),
the Borel-equivariant/orbifold Deligne cohomology considered in
\cite{KubelThom18} (which, for finite groups, coincides with \cite{LupercioUribe03}\cite{Gomi05}),
given, in low degrees, by:
\begin{enumerate}[{\bf (i)}]
\vspace{-1.5mm}
\item  gerbes with connection,
hence including what in string theory is known as the \emph{discrete torsion}
classification
of the B-field on orbifolds
\cite{Vafa86}\cite{Vafa-Witten}\cite{Sharpe00a}\cite{Sharpe00c}\cite{Sharpe02}\cite{Ruan};
or
\vspace{-.25cm}
\item
 2-gerbes with connection, hence including
what in M-theory is known as the \emph{discrete torsion}
classification of the C-field on
orbifolds \cite{Sharpe00b}\cite{Seki01}\cite[\S 4.6.2]{dBDHKMMS02}.
\end{enumerate}

\vspace{0mm}
\noindent {\bf Geometric enhancement of global equivariant homotopy theory.}
We enhance all of the above to a theory of properly orbi-singular spaces, formulated
internally to ``singular-elastic'' $\infty$-toposes
(\cref{SingularGeometry}),
where a natural notion of orbi-singularization $\orbisingular$
(Prop. \ref{OverSingularitiesCohesion})
promotes (Def. \ref{OrbiVFolds})
any such $\infty$-category of {\'e}tale $\infty$-groupoids
faithfully to its \emph{proper} orbifold version
(Remark \ref{Proper}).
This detects geometric fixed point spaces
(Def. \ref{GeometricAndHomotopyFixedPoints}) in the sense of
proper equivariant homotopy theory.
We show
(Prop. \ref{ShapeOfOrbisingularizedTopologicalGroupoidIsOrbispace},
Lemma \ref{ShapeOfOrbiSingularSpacesAsPresheafOnSingularities})
that the cohesive shape (Def. \ref{CohesiveTopos})
of the orbi-singularization of an {\'e}tale groupoid is its incarnation
as an orbispace in global equivariant homotopy theory,
in the sense of \cite{HenriquesGepner07}\cite{Rezk14}\cite{Koerschgen16}\cite{Schwede17}
(Remark \ref{TraditionalWayOfRegardingTopologicalGroupoidsAsOrbifolds}).

\vspace{1mm}
\hypertarget{FigureG}{}
\hspace{-.9cm}
\begin{tabular}{ll}

\begin{minipage}[l]{7.8cm}
\noindent {\bf The proper 2-category of orbifolds.}
One model for the axioms of singular-cohesive homotopy theory
is the $\infty$-topos of
\emph{singular-smooth $\infty$-groupoids}
(Examples \ref{SmoothInfinityGroupoids}, \ref{OrbiSingularSmoothInfinityGroupoids} below).
In this model, the proper 2-category (Rem. \ref{Proper})
of orbifolds $\mathcal{X}\!$,
either Lie theoretically (Example \ref{FrechetSmoothGOrbifolds})
or topologically (Example \ref{TopologicalGOrbiSingularSpaces}),
is equivalent, via passage to
\begin{enumerate}[{\bf (i)}]
\vspace{-3.3mm}
\item  their purely
smooth aspect $\smooth \mathcal{X}$,
to the traditional 2-category of {\'e}tale stacks
(Prop. \ref{SingularQuotientOfGOrbiSingularSpaces}),
\vspace{-.8cm}
\item while their  purely singular aspect
$\singular \mathcal{X}$ gives the underlying singular
coarse quotient space (Prop. \ref{PropertiesOfUniversalCoveringSpaces}).
\end{enumerate}
  \end{minipage}

  &
  \hspace{.6cm}
  \raisebox{28pt}{
  \xymatrix@C=40pt@R=-12pt{
    &
    \overset{
    \scalebox{.9}{
      \fbox{
        \begin{tabular}{c}
          \color{darkblue}
          \bf
          proper
          \\
          \color{darkblue}
          \bf
          orbifold
          \\
          \cref{OrbifoldGeometry}
        \end{tabular}
      }
      }
    }{
      \phantom{\mathclap{\vert^{\vert^{\vert}}}}
      \orbisingular
      (X \!\sslash\! G)
    }
    \ar@{|->}[dr]^-{
              \begin{rotate}{-31.5}
          \color{darkblue}
          \tiny
          \bf
          $
          \mathclap{
            \mbox{
              \begin{tabular}{c}
                project onto
                \\
                purely smooth aspect
                \\
                \phantom{$\vert_{\vert_{\vert}}$}
                \\
                \phantom{$\vert_{\vert_{\vert}}$}
              \end{tabular}
            }
          }
          $
        \end{rotate}
      }_-{\;\;\;\;\;\;\;
        \smooth
      }
    \ar@{|->}[dl]_-{
          \begin{rotate}{32}
          \color{darkblue}
          \tiny
          \bf
          $
          \mathclap{
            \mbox{
              \begin{tabular}{c}
                project onto
                \\
                purely singular aspect
                \\
                \phantom{$\vert_{\vert_{\vert}}$}
                \\
                \phantom{$\vert_{\vert_{\vert}}$}
              \end{tabular}
            }
          }
          $
        \end{rotate}
      }^-{\!\!\!\!
        \singular
      }
    \ar@{}[d]|{
      \mbox{\tiny \color{gray} \bf
        \begin{tabular}{c}
          \\
          opposite extreme
          \\
          aspects of orbifolds
        \end{tabular}
      }
    }
    \\
    \underset{
      \mathclap{
      \scalebox{.9}{
      \fbox{
        \begin{tabular}{c}
          \color{darkblue} \bf
          singular
          \\
          \color{darkblue} \bf
          space
          \\
          \cite{IKZ10}
        \end{tabular}
      }
      }
      }
    }{
      X/G
    }
    &&
    \underset{
      \mathclap{
      \scalebox{.9}{
      \fbox{
        \begin{tabular}{c}
          \color{darkblue} \bf
          smooth
          \\
          \color{darkblue} \bf
          {\'e}tale stack
          \\
          \cite{MoerdijkPronk97}\cite{CPRST14}
        \end{tabular}
      }
      }
      }
    }{
      X \!\sslash\! G
    }
  }
  }
\end{tabular}

\noindent {\bf Cartan geometry of {\'e}tale $\infty$-stacks.}
In this internal formulation we find
all fundamental phenomena of differential geometry
naturally generalized to {\'e}tale $\infty$-stacks,
hence in particular to orbifolds:

\medskip
{\small
\renewcommand\arraystretch{1.17}
\hspace{-.2cm}
\begin{tabular}{llll}
  &
  \begin{tabular}{l}
     \cref{VFolds}
  \end{tabular}
  &
    \hspace{-.2cm}
    \bf
    \begin{tabular}{l}
      Cartan geometry
      \\
      for {\'e}tale $\infty$-stacks
    \end{tabular}
  &
  Discussion for ordinary orbifolds, e.g. in
  \\
  \hline
  \hline
  {\bf (i)}
  &
  \begin{tabular}{l}
    Def. \ref{TangentBundleOfVFoldIsFiverBundle}
  \end{tabular}
  &
  \hspace{-5pt}
  \raisebox{-2pt}{
  \begin{tikzpicture}
    \clip (-1.1,4pt) rectangle (1.1,-5pt);
    \draw (0,0) node
    {
      \colorbox{lightgray}{\textbf{\textit{Frame bundles}}}
    };
  \end{tikzpicture}
  }
  &
  \cite[p. 42]{MoerdijkMrcun03}
  \\
  \hline
  {\bf (ii)}
  &
  \begin{tabular}{l}
    Def. \ref{GStructures}
  \end{tabular}
  &
  \hspace{-5pt}
  \raisebox{-2pt}{
  \begin{tikzpicture}
    \clip (-.95,4pt) rectangle (.95,-5pt);
    \draw (0,0) node
    {
      \colorbox{lightgray}{\textbf{\textit{G-structures}}}
    };
  \end{tikzpicture}
  }
  &
  \multirow{2}{*}{
    \hspace{-.2cm}
        \cite{Wolak16}\cite{Zu06}\cite{BagaevZhukova03}
  }
  \\
  &
  \begin{tabular}{l}
    Def. \ref{IntegrableGStructure},
    \\
    Def. \ref{LocallyIntegrableGStructures}
  \end{tabular}
  &
  \begin{tabular}{l}
    -{\it locally}
    \\
    -{\it globally}
  \end{tabular}
  {\it integrable}
  &
  \\
    {\bf (ii.a)}
  &
  \begin{tabular}{l}
    Ex. \ref{TorsionFreeGStructureOnSmoothManifolds}
  \end{tabular}
  &
  {\textbf{\textit{Geometric structures}}}
  &
   \cite[\S 1.8]{Apanasov00}\cite{Wolak16}
  \\
  & &
  - {\it Riemannian structure}
  &
  \hspace{-.3cm}
  \begin{tabular}{l}
   \cite{Borzellino92}\cite{HM04}\cite{Ratcliffe06}\cite{BagaevZhukova07}\cite{Hepworth09a}\cite{Hepworth09b}
  \\
   \cite{Akutagawal12}\cite{Kankaanrinta13}\cite{BDP17}\cite{Lange18}
  \end{tabular}
  \\
  & &
  - {\it Flat structure}
  &
    \cite{BDP17}\cite{Reffert06}\cite[\S 8]{IbanezUranga12}\cite{SS19a}
  \\
  & &
  - {\it Complex structure}
  &
   \cite{SteerWren99}\cite{FrenkelSzczesny07}
  \\
  & &
  - {\it Symplectic structure}
  &
   \cite{Verbitsky00}\cite{Godinho01}\cite{DE05}\cite{HM12}\cite{CP14}\cite{Chen17}\cite{RojoCarulli19}
  \\
  & &
  - {\it Lorentzian structure}
  &
  \cite{HorowitzSteif91}\cite{Nekrasov02}\cite{LiuMooreSeiberg02a}\cite{LiuMooreSeiberg02b}\cite{BerkoozReichmann07}\cite{ZhukovaRogozhina12}
  \\
& &
  - {\it Pseudo-Riemannian structure}
  &
\cite{Melnick}\cite{Zhukova18}\cite{BagaevZhukova19}
  \\
  & &
  - {\it Conformal structure}
  &
 \cite{Apanasov98}\cite{Apanasov00}
  \\
  & &
  - {\it CR-structure}
  &
 \cite{DragomirMasamune02}
  \\
& &
  - {\it Hypercomplex structure}
  &
 \cite{BoyerGalickiMann98}
  \\
    & & - \ldots
  \\
  {\bf (ii.b)}
  &
  \begin{tabular}{l}
    Ex. \ref{TorsionFreeGStructureOnSmoothManifolds}
  \end{tabular}
  &
  {\textbf{\textit{Special holonomy}}}
  &
 \cite{Jo00}\cite{CheegerTian05}
  \\
  & &
  - {\it K{\"a}hler structure}
  &
 \cite{Fujiki83}\cite{Jeffres97}\cite{Abreu01}\cite{BBFMT16}
  \\
  & &
  - {\it Calabi-Yau structure}
  &
         \hspace{-.3cm}
  \begin{tabular}{l}
 \cite{Roan91}\cite{Joyce98}\cite{Joyce99}\cite{Joyce99b}\cite[\S 6.5.1]{Jo00}\cite{Stapleton10}\cite{RuanZhang11}
     \\
    \cite{CDR16}
  \end{tabular}
  \\
  & &
  - {\it Quaternionic K{\"a}hler}
  &
 \cite{GalickiLawson88}\cite[\S 7.5.2]{Jo00}
  \\
  & &
  - {\it Hyper-K{\"a}hler struc.}
  &
  \cite{BD00}
  \\
  & &
  - {\it $G_2$-structure}
  &
 \cite[\S 11]{Jo00}\cite{Reidegeld15}
  \\
  & &
  - {\it $\mathrm{Spin}(7)$-structure}
  &
 \cite[\S 13]{Jo00}\cite{Bazaikin07}
  \\
  \hline
  {\bf (iii)}
  &
  Def. \ref{Isometries}
  &
  \hspace{-5pt}
  \raisebox{-2pt}{
  \begin{tikzpicture}
    \clip (-1.22,4pt) rectangle (1.22,-5pt);
    \draw (0,0) node
    {
      \colorbox{lightgray}{\textbf{\textit{Local isometries}}}
    };
  \end{tikzpicture}
  }
  &
 \cite{BagaevZhukova07}
  \\
  \hline
  {\bf (iv)}
  &
  Def. \ref{HaefligerGroupoid}
  &
  \hspace{-5pt}
  \raisebox{-2.5pt}{
  \begin{tikzpicture}
    \clip (-1.19,5pt) rectangle (1.19,-5pt);
    \draw (0,0) node
    {
      \colorbox{lightgray}{\textbf{\textit{Haefliger stacks}}}
    };
  \end{tikzpicture}
  }
  &
  \cite{Haefliger71}\cite{Haefliger84}\cite[\S 2.2, \S 3]{Carchedi12}\cite{Carchedi15}).
  \\
  \hline
  {\bf (v)}
  &
  Def. \ref{TangentialStructures}
  &
  \hspace{-5pt}
  \raisebox{-3pt}{
  \begin{tikzpicture}
    \clip (-1.55,5pt) rectangle (1.55,-5pt);
    \draw (0,0) node
    {
      \colorbox{lightgray}{\textbf{\textit{Tangential structures}}}
    };
  \end{tikzpicture}
  }
  &
\cite{Weelinck18}\cite{Pardon20}
  \\
  {\bf (v.a)}
  &
  Ex. \ref{HigherSpinStructure}
  &
  {\textbf{\textit{Higher Spin-structures}}}
  \\
  & &
  - {\it Orientation}
  &
\cite{Druschel94}
  \\
  & &
  - {\it Spin structure}
  &
   \cite{Ve}\cite{Acosta01}\cite{BGR}\cite{DLM}
  \\
  & &
  - {\it Spin${}^c$ structure}
  &
 \cite[\S 14]{Du}
  \\
  & &
  - {\it String structure}
  &
  \cite{PW88}\cite{LU04}\cite{LU06}
  \\
  & &
  - {\it Fivebrane structure}
  &
 \cite{BotvinnikLabbi12} (cf. \cite{SSS08}\cite{SSS09})
  \\
  & &
  \ldots
\end{tabular}
}

\vspace{3mm}
\noindent {\bf Orbifold {\'e}tale cohomology.}
Based on this, we give a natural general definition of
{\it {\'e}tale cohomology of $V$-{\'e}tale $\infty$-stacks }
(Def. \ref{EtaleCohomologyOfVFolds})
hence in particular of {\it orbifold {\'e}tale cohomology},
which is sensitive to the above (integrable) $G$-structures,
and hence to geometry/special holonomy on orbifolds.
For example, in the case of complex structure,
this orbifold {\'e}tale cohomology
subsumes traditional notions of complex-geometric orbifold
cohomology such as orbifold Dolbeault cohomology
\cite{Bai54}\cite{Bai56}\cite{CR04}\cite{Fernandez03}
and orbifold Bott-Chern cohomology \cite{Angella12}\cite{Ma05}.
Abstractly, orbifold {\'e}tale cohomology is the
intrinsic cohomology \eqref{IntrinsicCohomologyOfAnInfinityTopos}
of integrably $G$-structured
{\'e}tale $\infty$-stacks when regarded in the slice
$\infty$-topos (Prop. \ref{SliceInfinityTopos})
over the $G$-Haefliger stack (Def. \ref{HaefligerGroupoid})
via the classifying map of their $G$-structure (Prop. \ref{GStructuredHaefligerGroupoidClassifiesIntegrability}).
As such, orbifold {\'e}tale cohomology is the progenitor of
\emph{tangentially twisted} (proper) orbifold cohomology
(Def. \ref{TangentiallyTwistedCohomology}, Def. \ref{OrbifoldCohomologyTangentiallyTwisted}),
to which we turn next.

\newpage

\noindent{\bf Proper equivariant cohomology.}
While the proper 2-category of orbifolds is equivalent to
the traditional one of orbifolds as {\'e}tale stacks,
its full embedding into an ambient singular-cohesive $\infty$-topos
(\cref{SingularGeometry})
provides for more general coefficient objects, each of which is
guaranteed to produce a proper orbifold Morita-class invariant
(Remark \ref{Proper}).
Our {\bf first main Theorem} \ref{OrbifoldCohomologyEquivariant} shows
that the intrinsic cohomology \eqref{IntrinsicCohomologyOfAnInfinityTopos}
of orbifolds, regarded in singular-cohesive homotopy theory (Def. \ref{SingularCohesiveInfinityTopos}),
subsumes all proper $G$-equivariant cohomology theories:
Bredon cohomology with any coefficient system, as well as
proper equivariant generalized cohomology theories.

\medskip

\noindent {\bf Traditional orbifold cohomology.}
In particular, Prop. \ref{ShapeOfOrbisingularizedTopologicalGroupoidIsOrbispace} and
Theorem \ref{OrbifoldCohomologyEquivariant}
imply, via \cite{Juran20} (Remark \ref{TraditionalWayOfRegardingTopologicalGroupoidsAsOrbifolds}),
that proper orbifold cohomology in singular-cohesive homotopy theory
subsumes Chen-Ruan orbifold cohomology, via \eqref{ChenRuanCohomologyIsABredonCohomologyTheory},
and Adem-Ruan orbifold K-theory, via \eqref{TraditionalOrbifoldKTheory}. Hence it
also subsumes Freed-Hopkins-Teleman orbifold K-theory \cite{FHT07}
(reviewed in \cite[\S 3.2.2]{Nuiten13})
and Jarvis-Kaufmann-Kimura's ``full orbifold K-theory'' \cite{JKK05}\cite{GHHK08}
for orbifolds with global quotient presentations
(by \cite[Prop. 3.5]{FHT07} and \cite[3]{JKK05},
respectively).
Moreover, singular-cohesion provides a natural transformation
$\!\!
  \xymatrix@C=3em{
    \smooth
    \mathcal{X}
    \ar[r]|-{\;
      \epsilon^{\scalebox{.5}{\smooth}}_{\mathcal{X}}
    \;}
    &
    \orbisingular
    \mathcal{X}
  }\!\!
$
which restricts this proper orbifold cohomology to the
underlying {\'e}tale stack, where it
reduces to traditional Borel orbifold cohomology
\eqref{TraditionalOrbifoldCohomology} and, in particular,
to Satake cohomology \eqref{SatakeCohomology}
(see also, e.g., \cite{ADG11}\cite{BNSS18}).

\medskip
\noindent {\bf Twisted orbifold cohomology.} All these
cohomology theories generalize to their twisted cohomology versions
(e.g., local coefficients for ordinary cohomology, as in \cite{MoerdijkPronk99},
or twisted K-theory, as in \cite{AdemRuan01}),
by passage to slice $\infty$-toposes
of the ambient singular-cohesive $\infty$-topos (Remark \ref{TwistedCohomology}).
In particular, slicing of orbifolds over $\orbisingular^{\hspace{-5pt}\scalebox{.6}{$\mathbb{Z}_2$}}$
via their orientation bundle promotes them
(Example \ref{OrientifoldCohomology}) to \emph{orientifolds}
\cite{DFM11}\cite[4.4]{FSS15}\cite{SS19a}.

\hspace{-.9cm}
\begin{tabular}{ll}
\begin{minipage}[l]{7.4cm}
\noindent {\bf Proper orbifold {\'e}tale cohomology.}
Finally, we promote (Def. \ref{OrbifoldCohomologyTangentiallyTwisted})
orbifold {\'e}tale cohomology,
in its guise as tangentially twisted cohomology, to a
\emph{proper} orbifold cohomology theory
in the above sense (Remark \ref{Proper}).
Our {\bf second main theorem} \ref{OrbifoldCohomologyTangentiallyTwistedReduces}
shows  that this \emph{proper orbifold {\'etale} cohomology} unifies:

\vspace{-.3cm}

\begin{enumerate}[{\bf (i)}]
 \vspace{-2mm}
\item {($\smooth$)} \hspace{.1cm}
{\'e}tale cohomology (Def. \ref{EtaleCohomologyOfVFolds})
of smooth $V$-folds (Def. \ref{VManifold}).
 \vspace{-3mm}
\item {($\,\flat\,$)} \hspace{.1cm} proper equivariant cohomology
(Def. \ref{ProperEquivariantCohomologyInSingularCohsion})
of flat orbifolds (Def. \ref{FlatVFolds}),  i.e., of their flat frame bundles
(Prop. \ref{FlatFrameBundlesAreVFolds}).
\end{enumerate}
\end{minipage}
&
\hspace{1.6cm}
\raisebox{35pt}{
  \xymatrix@R=-4pt@C=-5pt{
  &
  \overset{
    \mathclap{
    \raisebox{20pt}{
    \scalebox{.9}{
      \fbox{
      \begin{tabular}{c}
        \color{darkblue}
        \bf
        proper
        \\
        \color{darkblue}
        \bf
        orbifold {\'e}tale cohomology
        \\
        Def. \ref{OrbifoldCohomologyTangentiallyTwisted}
      \end{tabular}
      }
    }
    }
    }
  }{
  \mathclap{\phantom{\vert^{\vert^{\vert}}}}
  H^{
    \scalebox{.7}{$ \raisebox{1pt}{\rm\textesh}\orbisingular\tau $}
  }
  \big(
    \mathcal{X}
    \!,\,
    \mathcal{A}
  \big)
  }
  \ar@{|->}[ddl]^-{
    (i_{\scalebox{.6}{$\smooth$}})^\ast
  }_-{ \!\!\!\!                  \begin{rotate}{50}
           \hspace{-.9cm}
    \color{darkblue}
        \bf
       \scalebox{.6}{ smooth orbifolds}
        \end{rotate}
  }
  \ar@{|->}[ddr]_-{
    (i_{\scalebox{.6}{$\flat$}})^\ast
  }^-{
  \!
      \begin{rotate}{-38}
           \hspace{-.7cm}
    \color{darkblue}
        \bf
       \scalebox{.6}{flat orbifolds}
        \end{rotate}
  }
  \\
  \\
  \underset{
    \mathclap{
    \scalebox{.9}{
      \fbox{
      \hspace{-5pt}
      \begin{tabular}{c}
        \color{darkblue}
        \bf
        smooth
        \\
        \color{darkblue}
        \bf
        {\'e}tale cohomology
        \\
        Def. \ref{TangentiallyTwistedCohomology}
      \end{tabular}
      \hspace{-5pt}
      }
    }
    }
  }{
  \mathclap{\phantom{\vert_{\vert_{\vert}}}}
  H^{
    \scalebox{.7}{$ \raisebox{0pt}{\rm\textesh}\tau $}
  }
  \big(
    X
    \!,\,
    A
  \big)
  }
  &&
  \underset{
    \mathclap{
    \scalebox{.9}{
      \fbox{
      \hspace{-12pt}
      \begin{tabular}{c}
        \color{darkblue}
        \bf
        proper
        \\
        \color{darkblue}
        \bf
        equivariant cohomology
        \\
        Def. \ref{ProperEquivariantCohomologyInSingularCohsion}
      \end{tabular}
      \hspace{-12pt}
      }
    }
    }
  }{
    \mathclap{\phantom{\vert_{\vert_{\vert}}}}
    H_{\flat G}
    \big(
      (\flat G)\mathrm{Frames}(X)
      ,
      A
    \big)
  }
  }
  }
\end{tabular}

\vspace{.1cm}

\hspace{-.9cm}
\begin{tabular}{ll}
\label{JTwistedOrbifoldCohomotopyDiagram}
\begin{minipage}[l]{7.4cm}
\noindent {\bf J-twisted orbifold Cohomotopy.}
We construct a fundamental class of examples of such
proper orbifold {\'e}tale cohomology theories,
which we call \emph{J-twisted orbifold Cohomotopy theories}
(Def. \ref{JTwistedOrbifoldCohomotopyTheory}).
Their coefficients are \emph{Tate spheres} (Def. \ref{TheVSphere}),
in the sense of (unstable) motivic homotopy theory
(Example \ref{TateSphereInMotivicHomotopyTheory}), with twisting via
an intrinsic \emph{Tate J-homomorphism} (Def. \ref{TateJHomomorphism}).
Specified to ordinary orbifolds
(Example \ref{JTwistedOrbifoldCohomotopyOnOrdinaryOrbifolds}),
Theorem \ref{OrbifoldCohomologyTangentiallyTwistedReduces}
shows that J-twisted orbifold Cohomotopy subsumes:
\begin{enumerate}[{\bf (i)}]
 \vspace{-3mm}
\item {($\smooth$)} \hspace{.05cm} J-twisted Cohomotopy theory
of smooth but curved spaces, as introduced in \cite{FSS19b}\cite{FSS19c}.

 \vspace{-3mm}
\item
{($\flat$)} \hspace{.05cm} RO-graded equivariant Cohomotopy theory
of flat orbifolds, as discussed in \cite{SS19a}\cite{SS19b}.
\end{enumerate}
\end{minipage}
&
\hspace{1.6cm}
\raisebox{43pt}{
  \xymatrix@R=4pt@C=10pt{
  &
  \overset{
    \mathclap{
    \raisebox{20pt}{
    \scalebox{.9}{
      \fbox{
      \begin{tabular}{c}
        \color{darkblue}
        \bf
        J-twisted
        \\
        \color{darkblue}
        \bf
        orbifold Cohomotopy
        \\
        Def. \ref{JTwistedOrbifoldCohomotopyTheory}
      \end{tabular}
      }
    }
    }
    }
  }{
  \mathclap{\phantom{\vert^{\vert^{\vert}}}}
  \pi^{
    \scalebox{.7}{$ \raisebox{1pt}{\rm\textesh}\orbisingular\tau $}
  }
  \big(
    \mathcal{X}
  \big)
  }
    \ar@{|->}[ddl]^-{
    (i_{\scalebox{.6}{$\smooth$}})^\ast
  }_-{ \!\!\!\!                  \begin{rotate}{55}
           \hspace{-.9cm}
    \color{darkblue}
        \bf
       \scalebox{.6}{ smooth orbifolds}
        \end{rotate}
  }
  \ar@{|->}[ddr]_-{
    (i_{\scalebox{.6}{$\flat$}})^\ast
  }^-{
  \!
      \begin{rotate}{-43}
           \hspace{-.7cm}
    \color{darkblue}
        \bf
       \scalebox{.6}{flat orbifolds}
        \end{rotate}
  }
  \\
  \\
  \underset{
    \mathclap{
    \scalebox{.9}{
      \fbox{
      \hspace{-5pt}
      \begin{tabular}{c}
        \color{darkblue}
        \bf
        J-twisted
        \\
        \color{darkblue}
        \bf
        Cohomotopy theory
        \\
        \cite{FSS19b}\cite{FSS19c}
      \end{tabular}
      \hspace{-5pt}
      }
    }
    }
  }{
  \mathclap{\phantom{\vert_{\vert_{\vert}}}}
  \pi^{
    \scalebox{.7}{$ \raisebox{0pt}{\rm\textesh}\tau $}
  }
  \big(
    X
  \big)
  }
  &&
  \underset{
    \mathclap{
    \scalebox{.9}{
      \fbox{
      \hspace{-12pt}
      \begin{tabular}{c}
        \color{darkblue}
        \bf
        RO-graded
        \\
        \color{darkblue}
        \bf
        equivariant Cohomotopy
        \\
        \cite{SS19a}\cite{SS19b}
      \end{tabular}
      \hspace{-12pt}
      }
    }
    }
  }{
    \mathclap{\phantom{\vert_{\vert_{\vert}}}}
    \pi^{}_{\flat G}
    \big(
      (\flat G)\mathrm{Frames}(X)
    \big)
  }
  }
  }
\end{tabular}

\vspace{.15cm}

\noindent We conclude with a Remark \ref{HypothesisH}
on the impact of this unification.

\newpage

\noindent {\bf Other approaches and outlook.}
We briefly comment on relations of our constructions to other approaches
in the literature, further discussion of which is beyond the scope of this
article.

\vspace{1mm}
\noindent {\bf Proper $\infty$-categories of general {\'e}tale $\infty$-stacks.}
Another general theory of {\'e}tale $\infty$-stacks
has been presented in \cite{Carchedi13},
generalizing an elegant characterization of {\'e}tale 1-stacks
due to \cite{Carchedi12}
by following the discussion of derived Deligne-Mumford stacks
conceived as structured $\infty$-toposes in \cite{Lurie09b}.
This approach proceeds externally via characterizing the \emph{sites}
(recalled below as Prop. \ref{ToposLexReflection})
which present $\infty$-toposes of {\'e}tale $\infty$-stacks; and is thus
complementary to the internal perspective proceeding
from inside an ambient $\infty$-topos which we are presenting here.
We briefly indicate the relation between the two:

\begin{itemize}
\vspace{-2mm}
\item[$\circ$] The approach in \cite{Carchedi13} is to pick
an $\infty$-site of $\mathrm{ModelSpaces}$ (denoted ``$\mathcal{L}$'' there)
which is equipped with a suitable notion of which of its 1-morphisms
qualify as being {\'e}tale maps (the external version of our notion
Def. \ref{FormallyEtaleMorphism}).
The inclusion $i$ of the wide subcategory on these {\'e}tale
morphisms
induces, by left Kan extension,
a pair of adjoint $\infty$-functors $(i_! \dashv i^\ast)$
between the corresponding $\infty$-stack $\infty$-toposes,
and the {\'e}tale $\infty$-stacks are then
characterized as those in the essential
image of the left adjoint $i_!$.
This is shown on the right of the following diagram:

\vspace{-.4cm}
\begin{equation}
  \label{SubsumingMoreGeneralEtaleStacks}
  \hspace{-1.4cm}
  \raisebox{4pt}{\footnotesize
  \xymatrix@R=12pt@C=2.5em{
{\small     \mathrm{Sheaves}_\infty
    \big(
      \mathrm{ModelSpaces}
      \times
      \mathrm{Singularities}
    \big)
    }
    \ar@{->}@<+8pt>[rr]|-{ \,\mathrm{Smth} \,}
    \ar@{<-^{)}}@<-9pt>[rr]|-{
      \;  \mathrm{OrbSinglr} \,
    }
      ^-{ \raisebox{3pt}{\scalebox{.7}{$\bot$}} }
    \ar@{}@<-27pt>[rr]|-{
        \mbox{
          \tiny
          \color{darkblue}
          \bf
          Prop. \ref{OverSingularitiesCohesion}
        }
    }
    &&
 {\small      \mathrm{Sheaves}_\infty
    \big(
      \mathrm{ModelSpaces}
    \big)
    }
    \ar@{<-}@<-8pt>[r]|-{\; i_!  \,}
    \ar@<+8pt>[r]|-{\; i^\ast \,}_-{ \raisebox{-8pt}{\;\scalebox{.7}{$\top$}} }
    &
 {\small      \mathrm{Sheaves}_\infty
    \big(
      \mathrm{LocalModelSpaces}^{\mbox{\tiny{\'e}t}}
    \big)
    }
    \ar@{->>}[dl]
    \\
    \underset{
      \mbox{
        \tiny
        \color{darkblue}
        \bf
        \begin{tabular}{c}
          proper $\infty$-category
          \\
          of higher orbifolds
          \\
          (Remark \ref{Proper})
        \end{tabular}
      }
    }{
   {\small      \mathrm{OrbSinglr}
      \big(
        \mbox{{\'E}taleStacks}_\infty
      \big)
      }
    }
    \ar@{_{(}->}[u]
    &&
    \underset{
      \mbox{
        \tiny
        \color{darkblue}
        \bf
        \begin{tabular}{c}
          $\infty$-category
          \\
          of {\'e}tale $\infty$-stacks
          \\
          {\cite{Carchedi13}}
        \end{tabular}
      }
    }{
 \small      \mbox{{\'E}taleStacks}_\infty
    }
    \ar[ll]^-{ \simeq }
    \ar@{_{(}->}[u]
  }
  }
\end{equation}

\vspace{-6mm}
\noindent
\vspace{-1mm}
\item[$\circ$]Following Remark \ref{Proper}, we may and should enhance this
construction to the
\emph{proper $\infty$-category of higher orbifolds}
Def. \ref{SingularCohesiveInfinityTopos},  Def. \ref{OrbiVFolds},
as shown on the left in \eqref{SubsumingMoreGeneralEtaleStacks}.

\vspace{-3mm}
\item[$\circ$] In fact, the archetypical example of $\mathrm{ModelSpaces}$ considered in
\cite{Carchedi13} is $\mathrm{SmoothManifolds}$ (Def. \ref{SmoothManifolds}), in which case
the left hand side of \eqref{SubsumingMoreGeneralEtaleStacks} is
the singular-cohesive $\infty$-topos of
our Examples \ref{SmoothInfinityGroupoids}, \ref{OrbiSingularSmoothInfinityGroupoids},
containing the proper (Def. \ref{OrbiVFolds})
$\infty$-category of orbi-$\mathbb{R}^n$-folds in our Example \ref{EtaleLieGroupoidAsRnFold}.

\vspace{-3mm}
\item[$\circ$] On the other hand, a general $\infty$-topos
$\mathrm{Sheaves}_\infty(\mathrm{ModelSpaces})$ is
not going to be cohesive (Def. \ref{CohesiveTopos})
or even elastic (Def. \ref{ElasticInfinityTopos}).
This means that various nice geometric properties,
which we derive here, of
objects in the proper $\infty$-category of higher orbifolds,
are not guaranteed to exist in the general setup of \cite{Carchedi13}.
Notably the theory of frame bundles
on orbifolds, according to Prop. \ref{TangentBundleOfVFoldIsFiverBundle},
and the main theorem on the induced {\'e}tale cohomology
of orbifolds (Theorem \ref{OrbifoldCohomologyTangentiallyTwistedReduces})
crucially uses the internal modal logic of
singular-cohesive and singular-elastic $\infty$-toposes
as in \cref{SingularCohesiveGeometry},
which may not exist, or not exist completely, for any given
site of $\mathrm{ModelSpaces}$ as in \cite{Carchedi13}.
\end{itemize}

\noindent {\bf Proper orbifold differential cohomology.}
While
\begin{itemize}
\vspace{-2mm}
\item[{\bf (i)}] generalized differential cohomology on smooth manifolds
\cite{HopkinsSinger05} is well-understood (see \cite{Bunke12})
and
\vspace{-2mm}
\item[{\bf (ii)}] plain global equivariant cohomology has been established
\cite{Schwede18} and understood to provide proper orbifold
cohomology (\cite{Juran20}, see Remark \ref{TraditionalWayOfRegardingTopologicalGroupoidsAsOrbifolds} below),
\end{itemize}

\vspace{-2mm}
\noindent their combination to (generalized, global)
\emph{proper equivariant differential cohomology} has remained
elusive.
Explicit constructions have been explored for the case
of equivariant/orbifold differential K-theory
\cite{SzaboValentino07}\cite{BunkeSchick09} \cite{Ortiz09}, but
even these theories do not seem to be well-understood yet \cite[p. 47]{BunkeSchick12}.
What has been missing is a coherent theoretical framework
for proper equivariant differential cohomology: Since
\begin{itemize}
\vspace{-2mm}
\item[{\bf (a)}] differential cohomology is the intrinsic cohomology
\eqref{IntrinsicCohomologyOfAnInfinityTopos} of
cohesive $\infty$-toposes (by Remark \ref{DifferentialCohomologyTheory})
and
\vspace{-3mm}
\item[{\bf (b)}] proper equivariant cohomology is the intrinsic cohomology
of $\infty$-toposes over a (global) orbit category (by Remark \ref{ProperEquivariantCohomologyTheory}),
\end{itemize}

\vspace{-2mm}
\noindent proper equivariant differential cohomology should be the intrinsic
cohomology of $\infty$-toposes that combine these two properties.
This is exactly what our notion of singular-cohesive $\infty$-toposes
expresses (Def. \ref{SingularCohesiveInfinityTopos}),
as confirmed by Theorem \ref{OrbifoldCohomologyEquivariant}.
For example, in singular-cohesive $\infty$-toposes there exists
the (global) proper equivariant
version of twisted differential non-abelian cohomology \cite{FSS20c},
now given by homotopy fiber
products parametrized over $\mathrm{Singularities}$ (Def. \ref{CategoryOfSingularities}).
Hence singular-cohesive $\infty$-toposes constitute a coherent framework in which
to discuss
proper equivariant/orbifold differential cohomology in general.
We will develop this elsewhere.

\newpage

\section{Preliminaries}

We recall basics of higher topos theory in \cref{HigherToposTheory}
and lay out in \cref{GaloisTheory}
the \emph{internal formulation}, in $\infty$-toposes,
of group actions and the classification of fiber bundles.

\subsection{Topos theory}
\label{HigherToposTheory}

We briefly record basics of $\infty$-topos theory \cite{ToenVezzosi05}\cite{Lurie09}\cite{Rezk10}
(review is in \cite{Rezk19}, exposition with an eye towards
differential geometric applications is in \cite{FSS13a}).
This is to set up our notation and
to highlight some less widely used aspects that we need further below.

\medskip

\noindent {\bf Categories.}
We make free use of the language and the basic facts
of category theory and homotopy theory
(see \cite{GoerssJardine99}\cite{Riehl14}\cite{Richter20})
as well as of $\infty$-category theory
(see \cite{Joyal08}\cite{Joyal08b}\cite{Lurie09}\cite{Riehl14}\cite{Cisinski19}).

\noindent {\bf (i)} We write $\mathrm{Categories}_\infty$ for the (``very large'')
$\infty$-category of (large) $\infty$-categories
\cite{Rezk98}\cite{Bergner05}\cite[Ch. 3]{Lurie09}, though
we only use this for declaring $\infty$-categories.
Inside $\mathrm{Categories}_\infty$,
there is the sequence of full sub-$\infty$-categories (Def. \ref{FullyFaithfulFunctor})
of $n$-categories (i.e.: $(n,1)$-categories)
as well as of $n$-groupoids (see Def. \ref{InfinityGroupoids})
for all $n \in \mathbb{N}$, denoted thus:
\begin{equation}
  \label{InfinityCategoryOfInfinityCategories}
  \xymatrix@R=7pt{
    \mathrm{Sets}
   \; \ar@{^{(}->}[r]
    &
    \mathrm{Categories}_1
   \; \ar@{^{(}->}[r]
    &
    \mathrm{Categories}_2
   \; \ar@{^{(}->}[r]
    &
    \cdots
   \; \ar@{^{(}->}[r]
    &
    \mathrm{Categories}_\infty
    \\
    \mathrm{Sets}
  \;  \ar@{^{(}->}[r]
    \ar@{=}[u]
    &
    \mathrm{Groupoids}_1
  \;  \ar@{^{(}->}[r]
    \ar@{^{(}->}[u]
    &
    \mathrm{Groupoids}_2
  \;  \ar@{^{(}->}[r]
    \ar@{^{(}->}[u]
    &
    \cdots
  \;  \ar@{^{(}->}[r]
    &
    \mathrm{Groupoids}_\infty
    \ar@<+5pt>@{^{(}->}[u]_-{ \dashv }
    \ar@<-6pt>@{<-}[u]_-{\mathrm{Core}}
  }
\end{equation}

\noindent {\bf (ii)} Here $\mathrm{Core}(\mathcal{C})$ denotes
the maximal $\infty$-groupoid inside an $\infty$-category $\mathcal{C}$.

\noindent {\bf (iii)} For $\mathcal{C} \;\in\; \mathrm{Categories}_\infty$ and
for $X, Y \in \mathcal{C}$ a pair of objects, we write
\vspace{-2mm}
\begin{equation}
  \label{HomInfinityGroupoids}
  \mathcal{C}(X,Y) \;:=\; \mathrm{Hom}_{\mathcal{C}}(X,Y)
  \;\;\in\;
  \mathrm{Groupoids}_\infty
  \end{equation}

\vspace{-2mm}
\noindent
for the \emph{hom-$\infty$-groupoid}, i.e. the $\infty$-groupoid of morphisms between them,
and higher homotopies between these
(see \cite[1.2.2]{Lurie09}\cite{DuggerSpivak09}).
This is well-defined, up to equivalence of $\infty$-groupoids,
independently of
which model for $\infty$-categories is used,
since these are all equivalent to each other \cite{Bergner06}\cite{Bergner14}.
We have no need to specify any particular model for $\infty$-categories
(except for the construction of examples, in \cref{ModelCategoryPresentations}).

\begin{defn}[Fully faithful functor {\cite[1.2.10]{Lurie09}}]
\label{FullyFaithfulFunctor}
For $\mathcal{C}, \mathcal{D} \in \mathrm{Categories}_\infty$
\eqref{InfinityCategoryOfInfinityCategories},
a functor $F \;:\; \xymatrix@C=12pt{ \mathcal{C} \ar[r] & \mathcal{D} }$
is called \emph{fully faithful}, to be denoted

\vspace{-.9cm}

\begin{equation}
  \label{FullSubcategoryInclusion}
  \xymatrix{
    \mathcal{C}
    \; \ar@{^{(}->}[r]^-{ F }
    &
    \mathcal{D}
    \,,
  }
\end{equation}
if it is an equivalence on all
hom-$\infty$-groupoids \eqref{HomInfinityGroupoids}:

\vspace{-.3cm}

$$
  \underset{
    X,Y \in \mathcal{C}
  }{\forall}
  \;\;
  \xymatrix{
    \mathcal{C}(X\,,\,Y)
    \ar[rr]^-{ F_{X,Y} }_-{ \simeq }
    &&
    \mathcal{D}
    \big(
      F(X)
      \,,\,
      F(Y)
    \big)
    \,.
  }
$$
In this case we also say that \eqref{FullSubcategoryInclusion}
exhibits a \emph{full sub-$\infty$-category inclusion}.
\end{defn}

\medskip

\noindent {\bf Topology.} The category of
\emph{$\Delta$-generated} or \emph{D-topological} spaces
(Remark \ref{EuclideanGeneratedIsDeltaGenerated})
is both: a convenient foundation for homotopy theory (Prop. \ref{NiceDTopologicalSpaces})
as well as pivotal for our key example context (Example \ref{SmoothInfinityGroupoids}):

\begin{defn}[Topological spaces]
  \label{TopologicalSpaces}
  \label{DTopologicalSpaces}
  We write
  \vspace{-1mm}
 \begin{equation}
    \label{CWComplexes}
    \xymatrix{
      \mathrm{CWComplexes}
      \;
      \ar@{^{(}->}[r]
      &
      \;\mathrm{DTopologcalSpaces}\;
      \ar@{^{(}->}[r]
      &
      \;\mathrm{TopologicalSpaces}
   }
   \;\;\;
   \in \mathrm{Categories}_1
  \end{equation}
   \vspace{-.6cm}

  \noindent
  for (from right to left):

  \noindent {\bf (i) } the category of all topological spaces
  with continuous functions between them;

  \noindent {\bf (ii)} the full subcategory
  on those spaces whose topology coincides with the final topology
  on the set of continuous functions out of a Euclidean space
  $\mathbb{R}^n$, hence whose open subsets
  coincide with those subsets whose pre-images under
  all continuous functions $\mathbb{R}^n \to X$
  are open in $\mathbb{R}^n$, for all $n \in \mathbb{N}$;

  \noindent {\bf (iii)}
  the further full subcategory on those
  that admit the structure of a CW-cell complex, hence
  that are homeomorphic to topological spaces which are obtained,
  starting with the empty space, by gluing on standard $n$-disks along
  their $(n-1)$-sphere boundaries, iteratively for $n \in \mathbb{N}$.
\end{defn}

\newpage

\begin{remark}[D-topological is $\Delta$-generated]
  \label{EuclideanGeneratedIsDeltaGenerated}
  $\phantom{A}$
  \vspace{-2mm}
\item {\bf (i)}   Since the topological $n$-simplex
  $\Delta^n_{\mathrm{top}}$ is a retract of
  the Euclidean space $\mathbb{R}^n$, the condition
  on $X \in \mathrm{TopologicalSpaces}$
  of being D-topological (Def. \ref{DTopologicalSpaces})
  is equivalent to being \emph{$\Delta$-generated}, in that the
  open subsets of $X$ are precisely those whose pre-images
  under all continuous functions of the
  form $\Delta^n_{\mathrm{top}} \to X$ are open.

  \vspace{-1mm}
  \item {\bf (ii)}  The concept of $\Delta$-generated spaces is due to
  \cite{Smith}\cite{Dugger03}; and independently due to
  \cite{SYH10}, where they are called \emph{numerically generated}.

  \vspace{-1mm}
  \item {\bf (iii)}   We say \emph{D-topological} to better bring out their
  conceptual role, in view of Prop. \ref{AdjunctionTopDiff} below.
\end{remark}

\begin{prop}[D-topological spaces are convenient]
  \label{NiceDTopologicalSpaces}
  The category of $\mathrm{DTopologicalSpaces}$ (Def. \ref{DTopologicalSpaces})
  is a \emph{convenient category of topological spaces} in the sense of \cite{Steenrod67} in that it:
  \begin{enumerate}[{\bf (i)}]
  \vspace{-2mm}
    \item contains all CW-complexes \eqref{CWComplexes} \cite[Cor. 4.4]{SYH10};
    \vspace{-2mm}
    \item  has all small limits and colimits \cite[Prop. 3.4]{SYH10};
  \vspace{-2mm}
    \item  is locally presentable \cite[Cor. 3.7]{FajstrupRosicky07};
   \vspace{-2mm}
    \item  is Cartesian closed \cite[Cor. 4.6]{SYH10}:
      the \emph{mapping space}
      between $X, Y \in \mathrm{DTopologicalSpaces}$
      is the reflection \eqref{AdjunctionBetweenTopologicalAndDiffeologicalSpaces}
      of
      the internal mapping space
      $\mathbf{Maps}$ \eqref{InternalHomAdjunction}
      of $\mathrm{DiffeologicalSpaces}$ \cite[Prop. 4.7]{SYH10}:
       \vspace{-2mm}
     \begin{equation}
        \label{InternalHomInDTopologicalSpaces}
        \mathrm{Maps}(X,Y)
        \;=\;
        \mathrm{Dtplg}
        \Big(
          \mathbf{Maps}
          \big(
            \mathrm{Cdfflg}(X),
            \,
            \mathrm{Cdfflg}(Y)
          \big)
        \Big).
      \end{equation}
  \end{enumerate}
\end{prop}

\medskip

\noindent {\bf Differential topology.}
D-topological spaces lend themselves to differential topology
via their joint (co-)reflection (Prop. \ref{AdjunctionTopDiff})
both into all topological spaces and into diffeological spaces
(Def. \ref{DiffeologicalSpaces}):
\begin{defn}[Cartesian spaces]
  \label{CartesianSpaces}
  We write
   \vspace{-1.5mm}
  $$
    \xymatrix{
      \mathrm{CartesianSpaces}
      \;
      \ar@{^{(}->}[r]
      &
      \mathrm{SmoothManifolds}
    }
    \;\;\in\;
    \mathrm{Categories}_1
  $$

   \vspace{-1.5mm}
\noindent
  for the category whose objects are the natural numbers
  $n \in \mathbb{N}$, thought of as representing the
  Cartesian spaces $\mathbb{R}^n$, and whose morphisms
  are the \emph{smooth} functions between these.
  We regard this category as equipped with the coverage
  (Grothendieck pre-topology)
  whose covers are the differentially good open covers
  (i.e., such that all non-empty finite intersections of patches
  are
  \emph{diffeomorphic} to a Cartesian space \cite[6.3.9]{FSS12}).
\end{defn}

\begin{defn}[Diffeological spaces]
  \label{DiffeologicalSpaces}
{\bf (i)} The category of \emph{diffeological spaces}
  (\cite{Souriau80}\cite{Souriau84}\cite{IZ85}, see \cite{BaezHoffnung08} \cite{IZ13})
  is the full subcategory of sheaves on
  $\mathrm{CartesianSpaces}$ (Def. \ref{CartesianSpaces})
  \vspace{-2mm}
  \begin{equation}
    \label{FullSubcategoryOfDiffeologicalSpaces}
    \xymatrix{
    \mathrm{DiffeologicalSpaces}
      \;\ar@{^{(}->}[r] &
    \mathrm{Sheaves}(\mathrm{CartesianSpaces})
    }
  \end{equation}

  \vspace{-2mm}
\noindent
  on those $X \in \mathrm{Sheaves}(\mathrm{CartesianSpaces})$
  which are \emph{concrete sheaves} \cite{Dubuc79}
  supported on their \emph{underlying set}
  \vspace{-2mm}
  $$
    X_s \; := \; \mathrm{Sheaves}(\mathrm{SmthMfd})(\ast, X)
  $$

\vspace{-2mm}
\noindent
  in that the canonical function
  \vspace{-4mm}
  $$
  \xymatrix{
    X(U)   \;\ar@{^{(}->}[r] & \mathrm{Set}(U_s, X_s)
    }
  $$

  \vspace{-1mm}
\noindent
  is an injection,
  for all $U \in \mathrm{CartesianSpaces}$,
  with $U_s$ denoting their underlying set of $U$.

  \vspace{-2mm}
 \item {\bf (ii)}  We call
  \begin{equation}
    \label{PlotsOfADiffeologicalSpace}
    X(U)
    \;\;
    \underset{
      \mathclap{
      \mbox{
        \tiny
        Prop. \ref{YonedaLemma}
      }
      }
    }{\simeq}
    \;\;
    \mathrm{DiffeologicalSpaces}(U,X)
    \;\;\in\;\;
    \mathrm{Set}
  \end{equation}
  the set of \emph{$U$-plots} of
  the diffeological space $X$.
\end{defn}

\begin{prop}[Topological/diffeological adjunction]
  \label{AdjunctionTopDiff}
  {\bf (i)}
  There is an adjunction \cite[Prop. 3.1]{SYH10}
  \vspace{-2mm}
  \begin{equation}
    \label{AdjunctionBetweenTopologicalAndDiffeologicalSpaces}
    \xymatrix{
      \mathrm{TopologicalSpaces}
      \;\;
      \ar@{<-}@<+6pt>[rr]^-{ \mathrm{Dtplg} }
      \ar@<-6pt>[rr]_{ \mathrm{Cdfflg} }^-{\bot}
      &&
     \;\; \mathrm{DiffeologicalSpaces}
    }
  \end{equation}

  \vspace{-3mm}
\noindent   between the categories of topological spaces (Def. \ref{TopologicalSpaces})
  and of diffeological spaces (Def. \ref{DiffeologicalSpaces}),
  where
  \begin{itemize}
    \vspace{-2mm}
    \item the right adjoint $\mathrm{Cdfflg}$ sends a topological space to the
    same underlying set equipped with the {\bf topological diffeology}
    whose plots \eqref{PlotsOfADiffeologicalSpace} are precisely the continuous functions;
     \vspace{-2mm}
    \item the left adjoint $\mathrm{Dtplg}$ sends a diffeological space
    to the same underlying set equipped with the
    {\bf diffeological topology} (``D-topology'' \cite[2.38]{IZ13}\cite{CSW13}),
    which is the
    final topology with respect to all plots \eqref{PlotsOfADiffeologicalSpace},
    hence such that a subset is
    open precisely if its pre-image under all plots is open.
  \end{itemize}

\vspace{-2mm}
  \noindent {\bf (ii)} The fixed points
  $X \in \mathrm{TopologicalSpaces}$
  of this  adjunction
  are the D-topological spaces (Remark \ref{EuclideanGeneratedIsDeltaGenerated})
  \vspace{-2mm}
   $$
    X
    \;\;
    \mbox{is D-topological}
    \phantom{AA}
    \Leftrightarrow
    \phantom{AA}
    \xymatrix{
      \mathrm{Dtplg}
      \big(
        \mathrm{Cdfflg}(X)
      \big)
      \ar[r]^-{ \epsilon_X }_-{\simeq}
      &
      X\;.
    }
  $$

\vspace{-2mm}
  \noindent {\bf (iii)} The adjunction is idempotent \cite[Lemma 3.3]{SYH10},
  hence factors through the category
  of D-topological spaces,
  exhibiting them as a co-reflective subcategory of
  $\mathrm{TopologicalSpaces}$ and a reflective subcategory
  of
 \\  $\mathrm{DiffeologicalSpaces}$:
 \vspace{-4mm}
  \begin{equation}
    \label{AdjunctionBetweenTopologicalAndDiffeologicalSpaces}
    \xymatrix{
      \mathrm{TopologicalSpaces}
    \;  \ar@{<-^{)}}@<+6pt>[rr]^-{ }
     \; \ar@<-6pt>[rr]_-{ \mathrm{Cdfflg} }^-{\bot}
      &&
      \;\; \mathrm{DTopologicalSpaces}
     \; \ar@{<-}@<+6pt>[rr]^-{ \mathrm{Dtplg} }
     \; \ar@{^{(}->}@<-6pt>[rr]_{ }^-{\bot}
      &&
     \; \mathrm{DiffeologicalSpaces}\;.
    }
  \end{equation}
\end{prop}

The following Prop. \ref{ModelStructureOnDTopologicalSpaces}
is due to \cite[Thm. 3.3]{Haraguchi13}.

\begin{prop}[Model structure on D-topological spaces]
  \label{ModelStructureOnDTopologicalSpaces}
 $\phantom{A}$

 \noindent {\bf (i)}
  The standard cell inclusions define a cofibrantly generated model category
  structure on $\mathrm{DTopologicalSpaces}$ (Def. \ref{DTopologicalSpaces}).

  \noindent {\bf (ii)} With respect to this model structure and the
  standard model structure on $\mathrm{TopologicalSpaces}$,
  the co-reflection \eqref{AdjunctionBetweenTopologicalAndDiffeologicalSpaces}
  becomes a Quillen equivalence:
  \vspace{-1mm}
  $$
    \xymatrix{
      \mathrm{TopologicalSpaces}
      \;\;
      \ar@<+8pt>@{<-^{)}}[rr]
      \ar@<-8pt>[rr]_-{ \mathrm{Cdfflg} }^-{
        \simeq_{{}_{\mathrlap{\mathrm{Quillen}}}}
      }
      &&
      \;\;\mathrm{DTopologicalSpaces}
    }.
  $$
\end{prop}

\noindent {\bf Differential geometry.}
\begin{defn}[Smooth Manifolds]
  \label{SmoothManifolds}
  We write
   \vspace{-2mm}
  \begin{equation}
    \label{CategoryOfSmoothManiolds}
    \mathrm{SmoothManifolds}
    \;\in\;
    \mathrm{Categories}
  \end{equation}

  \vspace{-2mm}
\noindent
  for the category of finite-dimensional paracompact smooth manifolds
  with smooth functions between them.
  We regard this as a site
  with the Grothendieck topology of open covers.
\end{defn}

\begin{prop}[Cartesian spaces are dense in the site of manifolds]
  \label{CartesianSpacesIsDenseSubsite}
  With respect to the coverages in Def. \ref{SmoothManifolds}
  and Def. \ref{CartesianSpaces}, the inclusion
  $\mathrm{CartesianSpaces} \overset{i}{\hookrightarrow} \mathrm{SmoothManifolds}$
  is a \emph{dense sub-site}, in that it induces an equivalence
  of categories of sheaves
  \vspace{-3mm}
  $$
    \xymatrix{
      \mathrm{Sheaves}(\mathrm{CartesianSpaces})
      \;
      \ar@{<-}@<+5pt>[rr]^-{ i^\ast }
      \ar@<-5pt>[rr]_-{ i_\ast }^-{ \simeq }
      &&
     \; \mathrm{Sheaves}(\mathrm{SmoothManifolds})\;.
    }
  $$
\end{prop}

\begin{prop}[Smooth manifolds inside diffeological spaces]
  \label{SmoothManifoldsInsideDiffeologicalSpaces}
  Every $X \in \mathrm{SmoothManifolds}$ \eqref{CategoryOfSmoothManiolds}
  becomes a diffeological space
  (Def. \ref{DiffeologicalSpaces}) on its underlying set
  by taking its plots \eqref{PlotsOfADiffeologicalSpace}
  of shape $U \in \mathrm{CartesianSpaces}$ to be the
  ordinary smooth functions:

  \vspace{-.9cm}

  $$
    X(U) \;:=\; \mathrm{SmoothManifolds}(U,X)
    \,.
  $$
  More generally, every
  possibly infinite-dimensional Fr{\'e}chet manifold
  (e.g. \cite[2.2]{KS17})
  becomes a diffeological space this way.
  Moreover, this constitutes fully faithfull embeddings
  (Def. \ref{FullyFaithfulFunctor})
  into the category of Diffeological spaces \cite[Thm. 3.1.1]{Losik94}:
  \begin{equation}
  \label{EmbeddingOfSmoothManifoldsIntoDiffeologicalSpaces}
  \xymatrix{
    \underset{
      \mathclap{
      \raisebox{-3pt}{
        \tiny
        \color{darkblue}
        \bf
        finite-dimensional
      }
      }
    }{
      \mathrm{SmoothManifolds}
    }
    \;
    \ar@{^{(}->}[r]
    &
    \underset{
      \mathclap{
      \raisebox{-3pt}{
        \tiny
        \color{darkblue}
        \bf
        possibly infinite-dimensional
      }
      }
    }{
      \mbox{\rm Fr{\'e}chetSmoothManifolds}
    }
    \;
    \ar@{^{(}->}[r]
    &
    \mathrm{DiffeologicalSpaces}
    \,.
  }
  \end{equation}
\end{prop}

\medskip

\noindent {\bf Homotopy theory.}
\begin{defn}[$\infty$-Groupoids]
  \label{InfinityGroupoids}
 {\bf (i)}
   We write
   \vspace{-2mm}
  \begin{equation}
    \label{InfinityCategoryOfInfinityGroupoids}
    \mathrm{Groupoids}_\infty
      \;\in\;
    \mathrm{Categories}_\infty
  \end{equation}

\vspace{-2mm}
\noindent
  for the $\infty$-category which is presented by the
  topologically enriched category whose objects are
  the CW-complexes \eqref{CWComplexes} and whose
  hom-spaces are the mapping spaces \eqref{InternalHomInDTopologicalSpaces}.

\vspace{-1mm}
 \item {\bf (ii)}  The full sub-$\infty$-category (Def. \ref{FullyFaithfulFunctor})
 on the homotopy $n$-types is
 that of \emph{$n$-groupoids}
 \vspace{-1mm}
  $$
    \mathrm{Groupoids}_n \xymatrix{\ar@{^{(}->}[r]&}\mathrm{Groupoids}_\infty
    \,.
  $$
\end{defn}

\newpage
\begin{defn}[Topological shape]
  \label{ShapeOfTopologicalSpaces}
 {\bf (i)}  We write
  $$
    \mathrm{Shp}_{\mathrm{Top}}
    \;:\;
    \xymatrix{
      \mathrm{CWComplexes}
      \;\;
      \ar[r]
      &
      \;\;
      \mathrm{Groupoids}_\infty
    }
  $$
  for the $\infty$-functor from the 1-category
  of CW-complexes \eqref{CWComplexes}
  to the $\infty$-category of $\infty$-groupoids
  (Def. \ref{InfinityGroupoids})
  which, as a topologically enriched functor, is
  the identity on objects, and is on hom-spaces
  the continuous map given by the identity
  function from the discrete set of continuous maps to the
  mapping space \eqref{InternalHomInDTopologicalSpaces}.

\vspace{-1mm}
 \item {\bf (ii)}    For any choice of CW-approximation functor
 \vspace{-3mm}
  $$
    \xymatrix{
      \mathrm{TopologicalSpaces}
      \ar[rr]^-{ (-)_{\mathrm{cof}} }
      &&
      \mathrm{CWComplexes}
    }
  $$
  we get the corresponding functor on all topological
  spaces (Def. \ref{TopologicalSpaces}),
  hence on D-toplogical spaces (Def. \ref{DTopologicalSpaces})
  which we denote by the same symbol:
  \vspace{-3mm}
  \begin{equation}
    \label{ShapeOfTopologicalSpaces}
    \mathrm{Shp}_{\mathrm{Top}}
    \;:\;
    \xymatrix{
      \mathrm{TopologicalSpaces}
      \ar[rr]^-{ (-)_{\mathrm{cof}} }
      &&
      \mathrm{CWComplexes}
      \ar[rr]^-{ \mathrm{Shp}_{\mathrm{Top}} }
      &&
      \mathrm{Groupoids}_\infty \;.
    }
  \end{equation}
\end{defn}

\begin{example}[Delooping groupoids]
  \label{DeloopingGroupoids}
  For $G \in \mathrm{Groups}^{\mathrm{fin}}$,
  consider the groupoid
  with a single object $\ast$,
  and with $G$ as its set of morphisms,
  whose composition is given by the product in the group:
  \vspace{-2mm}
  \begin{equation}
    \label{DeloopingGroupoidConvention}
    \xymatrix@R=.5em@C=3em{
      &
      \ast
      \ar[dr]^-{g_2}
      \\
      \ast
      \ar[ur]^-{ g_1 }
      \ar[rr]_-{ g_2 \cdot g_1 }
      &&
      \ast
    }
  \end{equation}

  \vspace{-2mm}
\noindent
This groupoid is the
topological shape \eqref{ShapeOfTopologicalSpaces} of
the Eilenberg-MacLane space $K(G,1)$ as well as
(since $G$ is assumed to be finite) the classifying space
$B G$.
More intrinsically, this groupoid is, equivalently,
the homotopy quotient of the point by the trivial $G$-action:

\vspace{-.9cm}

$$
  \ast \!\sslash\! G
  \;\in\;
  \mathrm{Groupoids}_1
  \xymatrix{
        \;   \ar@{^{(}->}[r]
      &
}
  \mathrm{Groupoids}_\infty
  \,.
$$
\end{example}
More generally:
\begin{example}[Action groupoids]
  \label{ActionGroupoids}
  For $G \in \mathrm{Groups}^{\mathrm{fin}}$ a finite group and for $X \in \mathrm{Set}$
  a set equipped with a $G$-action
    \vspace{-3mm}
    \begin{equation}
    \label{ActionOnSet}
    \xymatrix@R=-2pt{
      G \times X
      \ar[rr]^-{ \rho }
      &&
      X
      \\
      (g,x)
      \ar@{|->}[rr]
      &&
      g \cdot x
    }
  \end{equation}

   \vspace{-1mm}
\noindent
  the corresponding \emph{action groupoid} has
  as objects the elements of $X$ and
  its morphisms and their composition are given as follows:
   \vspace{-2mm}
  \begin{equation}
    \label{ActionGroupoidConvention}
    \xymatrix@R=.5em@C=3em{
      &
      g_1 \cdot x
      \ar[dr]^-{g_2}
      \\
      x
      \ar[ur]^-{ g_1 }
      \ar[rr]_-{ g_2 \cdot g_1 }
      &&
      g_2 \cdot g_1 \cdot x
    }
  \end{equation}
  This action groupoid is a model for the
  homotopy quotient of $X$ by its $G$-action
   \vspace{-1mm}
  $$
    X \!\sslash\! G
    \;\in\;
    \mathrm{Groupoids}_1
       \xymatrix{
        \;   \ar@{^{(}->}[r]
      &
}
    \mathrm{Groupoids}_\infty\;.
  $$
\end{example}

The following elementary example plays a pivotal role in
later constructions (Lemma \ref{ShapeOfOrbiSingularSpacesAsPresheafOnSingularities}):
\begin{example}[Hom-groupoid into action groupoid]
 \label{HomGroupoidFromBGIntoActionGroupoids}
Let $G \in \mathrm{Groups}^{\mathrm{fin}}$,
$X \in \mathrm{Set}$ equipped with a $G$-action \eqref{ActionOnSet},
hence with action groupoid/homotopy quotient
$X \!\sslash\! G \;\in\; \mathrm{Groupoids}_1$ (Example \ref{ActionGroupoids}).
Let $K \in \mathrm{Groups}^{\mathrm{fin}}$ be any finite group,
with $\ast \!\sslash\! K \;\in\; \mathrm{Groupoids}_1$
its delooping groupoid (Example \ref{DeloopingGroupoids}).
Then the hom-groupoid (functor groupoid) of morphisms (functors)
$\ast \!\sslash\! K \longrightarrow X \!\sslash\! G$
is, equivalently, the action groupoid of $G$ acting on
the set of pairs consisting of a group homomorphism
$\phi : K \to G$ and a point in $X$ fixed by the image of $\phi$:
 \vspace{-1mm}
\begin{equation}
  \label{GroupoidHomsBetweenQuotientGroupoidsOfSets}
  \mathrm{Groupoids}_1
  \big(
     \ast \!\sslash\! K
     \,
     ,
     \;
     X \!\sslash\! G
  \big)
  \;\simeq\;
  \left( \;\;\;\;\;\;\;
 \raisebox{5pt}{$\underset{{}_{
        \mathclap{
          \phi \in \mathrm{Groups}(K,G)
        }}
      }{
        \bigsqcup
      }
      \;\;\;\;\;
      X^{\phi(K)}
     $ }
  \right)
  \sslash G
  \,.
\end{equation}

\vspace{-2mm}
\noindent Here
\begin{itemize}
\vspace{-2mm}
  \item
    $\phi(K) \subset G$ denotes the subgroup of $G$
    which is image of the
    group homomorphism $\phi : K \to G$;
 \vspace{-2mm}
  \item
  $X^{\phi(K)}
    =
    \left\{
    x \in X
    \; \left|\;
    \underset{h \in \phi(K)}{\forall}
    h\cdot x = x
    \right.
    \right\}
  $
  denotes the $\phi(K)$-fixed-point set in $X$;
  \vspace{-2mm}
  \item
  the $G$-action by which the homotopy quotient is taken
  is the conjugation action on $\phi$,
  hence $g \cdot \phi \;:=\;  \mathrm{Ad}_g \circ \phi$,
  and the given $G$-action
  on $x \in X$.
\end{itemize}

\newpage
\vspace{-2mm}
\noindent This follows by direct unwinding of the definition of
functors and of natural transformations between the groupoids
\eqref{DeloopingGroupoidConvention} and \eqref{ActionGroupoidConvention}.
\end{example}

\begin{defn}[Simplicial-topological shape]
  \label{SimplicialTopologicalShape}
  Let
  \vspace{-2mm}
  \begin{equation}
    \label{SimplicialTopologicalSpace}
    X_\bullet
    \;:\;
    \xymatrix{
      \Delta^{\mathrm{op}}
      \ar[r]
      &
      \mathrm{TopologicalSpaces}
    }
  \end{equation}

  \vspace{-2mm}
\noindent
  be a simplicial topological space, for instance the nerve
  of a topological groupoid. Then we say that its
  \emph{simplicial-topological shape}
  is the homotopy colimit (Prop. \ref{LimitsAndColimitsAsAdjoints})
  of its degreewise topological shape (Def. \ref{ShapeOfTopologicalSpaces}):
  \vspace{-2mm}
  \begin{equation}
    \label{SimplicialTopologicalShape}
    \mathrm{Shp}_{\mathrm{sTop}}
    \big(
      X_\bullet
    \big)
    \;:=\;
    \underset{\longrightarrow}{\mathrm{lim}}
    \big(
      \mathrm{Shp}_{\mathrm{Top}}
      (X)
    \big)_\bullet
    \;\;
    \in
    \mathrm{Groupoids}_\infty
    \,.
  \end{equation}
\end{defn}

The following Prop. \ref{SimplicialTopologicalShapeOfDegreewiseCofibrantSimplicialSpaces}
appears as \cite[4.3,. 4.4]{Wang18}:

\begin{prop}[Simplicial-topological shape of degreewise cofibrant spaces is fat geometric realization]
  \label{SimplicialTopologicalShapeOfDegreewiseCofibrantSimplicialSpaces}
  If $X_\bullet$ is a simplicial topological space
  \eqref{SimplicialTopologicalSpace} which degreewise
  admits the structure of a retract of a cell complex
  (for instance: degreewise a CW-complex \eqref{CWComplexes}),
  then its simplicial topological shape \eqref{SimplicialTopologicalShape}
  is equivalent to its \emph{fat geometric realization}
  $\left\Vert -\right\Vert$ \cite{Segal74} (see \cite[2.3]{HenriquesGepner07}):
  \vspace{-2mm}
  $$
    \underset{
      \mathclap{
      \mbox{
        \tiny
        \color{darkblue}
        \bf
        \begin{tabular}{c}
          degreewise cofibrant
          \\
          simplicial topological spaces
        \end{tabular}
      }
      }
    }{
      X_\bullet
    }
    \;\in\;
    \big(
      \mathrm{TopologicalSpaces}_{\mathrm{cof}}
    \big)^{\Delta_{\mathrm{op}}}
    \;\;\;\;\;\;\;\;
    \Rightarrow
    \;\;\;\;\;\;\;\;
    \underset{
      \mathclap{
      \mbox{
        \tiny
        \color{darkblue}
        \bf
        \begin{tabular}{c}
          simplicial
          \\
          topological shape
        \end{tabular}
      }
      }
    }{
      \mathrm{Shp}_{\mathrm{sTop}}
      \big(
        X_\bullet
      \big)
    }
    \;\;\simeq\;\;\;
    \underset{
      \mathclap{
      \mbox{
        \tiny
        \color{darkblue}
        \bf
        \begin{tabular}{c}
          fat
          geometric
          \\
          realization
        \end{tabular}
      }
      }
    }{
      \left\Vert X_\bullet \right\Vert
    }
    \,.
  $$
\end{prop}

\begin{defn}[Diffeological simplices]
  \label{SmoothSimplices}
  \label{DiffeologicalSingularSimplices}
{\bf (i)}  We write
  \vspace{-3mm}
  $$
    \xymatrix@R=-4pt{
      \Delta
      \ar[rr]^-{ \Delta^\bullet_{\mathrm{smth}} }
      &&
      \mathrm{DiffeologicalSpaces}
      \\
      [n]
      \ar@{|->}[rr]
      &&
      \Big\{
        \vec x \in \mathbb{R}^{n+1}
        \,\vert\,
        \underset{i}{\sum} x^i  = 1
      \Big\}    }
  $$

  \vspace{-2mm}
\noindent  for the \emph{diffeological extended simplicies}, hence for the
  simplicial object in diffeological spaces (Def. \ref{DiffeologicalSpaces})
  (in fact in smooth manifolds, under Prop. \ref{SmoothManifoldsInsideDiffeologicalSpaces})
  which in degree $n$ is the extended
  $n$-simplex in $\mathbb{R}^{n+1}$, regarded with its sub-diffeology,
  and whose face and degeneracy maps are the
  standard ones (see \cite[Def. 4.3]{ChristensenWu14}\cite[p. 1]{Pavlov19}).

\vspace{-1mm}
\item {\bf (ii)}
The induced nerve/realization construction is a pair of adjoint functors
(Def. \ref{AdjointInfinityFunctors})
\vspace{-3mm}
  \begin{equation}
    \label{DiffeologicalSingularSimplicialSetAdjunction}
    \xymatrix{
      \mathrm{DiffeologicalSpaces}
     \; \ar@{<-}@<+6pt>[rr]^-{
        \left\vert - \right\vert_{\mathrlap{\mathrm{diff}}}
      }
      \ar@<-6pt>[rr]_{
        \mathrm{Sing}_{\mathrlap{\mathrm{diff}}}
      }^-{\bot}
      &&
   \;  \mathrm{SimplicialSets}
    }
  \end{equation}

  \vspace{-2mm}
\noindent   between the categories of simplicial sets
  and of diffeological spaces (Def. \ref{DiffeologicalSpaces}),
  where
 the right adjoint $\mathrm{Sing}_{\mathrm{diff}}$
     sends $X \in \mathrm{DiffeologicalSpaces}$ to its
     \emph{smooth singular simplicial set}
     $$
       \mathrm{Sing}_{\mathrm{diff}}(X)_\bullet
       \;:=\;
       \mathrm{DiffeologicalSpaces}
       \big(
         \Delta^\bullet_{\mathrm{diff}},
         X
       \big).
     $$
\end{defn}

The following Prop. \ref{DiffeologicalSingularSimpliciaSetOfContinuousDiffeology}
is due to \cite[Prop. 4.14]{ChristensenWu14}:

\begin{prop}[Diffeological singular simplicial set of continuous Diffeology]
  \label{DiffeologicalSingularSimpliciaSetOfContinuousDiffeology}
  For all $X_{\mathrm{top}} \in \mathrm{TopologicalSpaces}$
  there is a weak homotopy equivalence between
  the diffeological singular simplicial set
  (Def. \ref{DiffeologicalSingularSimplices}) of its
  continuous diffeology (Def. \ref{AdjunctionTopDiff})
  and its ordinary singular simplicial set:
  \vspace{-2mm}
  $$
    \mathrm{Sing}
    \big(
      X_{\mathrm{top}}
    \big)
    \;\simeq_{\mathrm{whe}}\;
    \mathrm{Sing}_{\mathrm{diff}}
    \big(
      \mathrm{Cdfflg}(X_{\mathrm{top}})
    \big)
    \,.
  $$

  \vspace{-2mm}
\noindent
  Equivalently this means, in the terminology to be introduced
  in a moment,
  that the topological shape \eqref{ShapeOfTopologicalSpaces}
  of topological spaces is equivalent to the cohesive shape
  (Def. \ref{CohesiveTopos})
  of their incarnation as continuous-diffeological spaces
  (see Example \ref{SmoothInfinityGroupoids} below):
  \vspace{-2mm}
  $$
    \mathrm{Shp}_{\mathrm{Top}}
    \big(
      X_{\mathrm{top}}
    \big)
    \;\simeq\;
    \mathrm{Shp}
    \big(
      \mathrm{Cdfflg}(X_{\mathrm{top}})
    \big)
    \;\;\;\;
    \in \mathrm{Groupoids}_\infty\;.
  $$
\end{prop}

\medskip

\noindent {\bf Universal constructions.}
All diagrams we consider now are homotopy-coherent, even if
we do not notationally indicate the higher cells, unless some are to be highlighted.
Similarly, all universal constructions we consider now are $\infty$-categorical,
even if this is not further pronounced by the terminology. In particular,
we say ``colimit'' $\underset{\longrightarrow}{\mathrm{lim}}$
for ``homotopy colimit'',
``limit'' $\underset{\longleftarrow}{\mathrm{lim}}$ for ``homotopy limit''
(see Prop. \ref{LimitsAndColimitsAsAdjoints}),
``Cartesian square'' for ``homotopy Cartesian square'', etc.:
\begin{notation}[Cartesian squares]
  \label{CartesianSquares}
  We say a square in an $\infty$-category is \emph{Cartesian},
  to be denoted
   \vspace{-2mm}
  \begin{equation}
    \label{PullbackSquare}
    \xymatrix@R=12pt{
      X \times_B Y
      \ar[d]_-{ f^\ast g }
      \ar[rr]
      \ar@{}[drr]|-{ \mbox{\tiny(pb)} }
      &&
      Y
      \ar[d]^-{g}
      \\
      X
      \ar[rr]|-{\, f \,}
      &&
      B
    }
  \end{equation}

   \vspace{-2mm}
\noindent
  if it is an limit cone over the diagram consisting of $f$ and $g$.
  We also say this is the \emph{pullback square} of $g$ along $f$.
\end{notation}
\begin{example}[Pullback of equivalence is equivalence]
  \label{PullbackOfEquivalenceIsEquivalence}
  Let $\mathcal{C} \in \mathrm{Categories}_\infty$. Then a square in $\mathcal{C}$
  whose right vertical morphism is an equivalence
  is Cartesian (Notation \ref{CartesianSquares}) precisely if the
  left vertical morphism is also an equivalence:

  \vspace{-.7cm}

  \begin{equation}
    \raisebox{16pt}{
    \xymatrix@R=15pt{
      A
        \ar@{}[drr]|-{\mbox{\tiny(pb)}}
        \ar[rr]
        \ar[d]
        &&
      B
        \ar[d]^-{\simeq}
      \\
      C
        \ar[rr]
        &&
      D
    }
    }
    \phantom{AAAA}
    \Leftrightarrow
    \phantom{AAAA}
    \raisebox{16pt}{
    \xymatrix@R=15pt{
      A
      \ar[d]^-{\simeq}
      \\
      C
    }
    }
  \end{equation}
  hence precisely if $\xymatrix@C=11pt{C \ar[r] & D}$ is equivalent to
  $\xymatrix@C=11pt{A \ar[r] & B}$ in $\mathcal{C}^{\Delta^1}$.
\end{example}
\begin{prop}[Pasting law {\cite[Lemma 4.4.2.1]{Lurie09}}]
  \label{PastingLaw}
  In any $\infty$-category,
  consider a diagram of the form
   \vspace{-2mm}
  $$
    \xymatrix@R=1.2em@C=3em{
      A
      \ar[d]\ar[r]
      \ar@{}[dr]|{
        \begin{rotate}{-45}
          $\mathclap{
            \raisebox{-3pt}{$\Downarrow$}
          }$
        \end{rotate}
      }
      &
      B
      \ar[d] \ar[r]
      \ar@{}[dr]|{
        \begin{rotate}{-45}
          $\mathclap{
            \raisebox{-3pt}{$\Downarrow$}
          }$
        \end{rotate}
      }
      &
      C \ar[d]
      \\
      D \ar[r] & E \ar[r] & F
    }
  $$

   \vspace{-2mm}
\noindent
  such that the right square is Cartesian (Notation \ref{CartesianSquares}).
  Then the left square is Cartesian
  if and only if the total rectangle is Cartesian.
\end{prop}

\begin{defn}[Adjoint $\infty$-functors {\cite[5.2.2.7, 5.2.2.8]{Lurie09}\cite[4.4.4]{RiehlVerity13}}]
  \label{AdjointInfinityFunctors}
  Let $\mathcal{C}, \mathcal{D} \;\in\; \mathrm{Categories}_\infty$
  \eqref{InfinityCategoryOfInfinityCategories} and
  $
    L
    :
    \xymatrix{
      \mathcal{C}
      \ar@{<->}[r]
      &
      \mathcal{D}
    }
    : R
  $
  two functors between them, back and forth.
  This is an \emph{adjoint pair} with $L$ \emph{left adjoint}
  and $R$ \emph{right adjoint}, to be denoted $(L \dashv R)$:
  \vspace{-2mm}
\noindent
  \begin{equation}
    \label{AdjointFunctors}
    \xymatrix{
      \mathcal{D}
        \ar@{<-}@<+6pt>[rr]^-{ L }
        \ar@<-6pt>[rr]_-{ R }^-{
          \raisebox{2pt}{
            \scalebox{.7}{$
              \bot
            $}
          }
        }
        &&
      \mathcal{C}
    }
  \end{equation}

  \vspace{-1mm}
\noindent  if there is a natural equivalence of hom-$\infty$-groupoids
  \eqref{HomInfinityGroupoids} of the form
  \begin{equation}
    \label{AdjunctionHomEquivalence}
    \mathcal{D}
    \big(
      L(-)
      \,,\,
      -
    \big)
    \;\;\;\;\simeq\;\;\;\;
    \mathcal{C}
    \big(
      -
      \,,\,
      R(-)
    \big)
  \end{equation}
  (This is unique when it exists \cite[Prop. 5.2.1.3, 5.2.6.2]{Lurie09}).
  In this case, one says:
  \begin{itemize}
  \vspace{-3mm}
    \item[{\bf (i)}]
      The \emph{adjunction unit} is the
      natural transformation
      \vspace{-2mm}
      \begin{equation}
        \label{AdjunctionUnit}
        \xymatrix{
          X
          \ar[r]^-{ \eta_X }
          &
          R \circ L (X)
        }
      \end{equation}

      \vspace{-2mm}
    \noindent
      which is  the (pre-)image under
      \eqref{AdjunctionHomEquivalence} of the identity on $R(X)$.

      \vspace{-3mm}
    \item[{\bf (ii)}]
      The \emph{adjunction co-unit} is the
      natural transformation
      \vspace{-2mm}
      \begin{equation}
        \label{CounitOfAdjunction}
        \xymatrix{
          L \circ R (X)
          \ar[r]^-{ \epsilon_X }
          &
          X
        }
      \end{equation}

      \vspace{-2mm}
   \noindent
      which is the image under
      \eqref{AdjunctionHomEquivalence} of the identity on $L(X)$.
  \end{itemize}
\end{defn}
As in the classical situation of 1-category theory, it follows that:

\begin{prop}[Triangle identities]
  \label{TriangleIdentities}
  Let
  $
    \xymatrix{
      \mathcal{D}
        \ar@{<-}@<+7pt>[r]|-{\, L \,}
        \ar@<-5pt>[r]|-{\, R \,}^-{
          \raisebox{1pt}{
            \scalebox{.6}{$
              \bot
            $}
          }
        }
        &
      \mathcal{C}
    }
  $
  be a pair of adjoint $\infty$-functors (Def. \ref{AdjointInfinityFunctors}).
  Then their adjunction unit $\eta$ \eqref{AdjunctionUnit}
  and counit $\epsilon$ \eqref{CounitOfAdjunction}
  satisfy the following natural equivalences:

  \noindent
  {\bf (i)} for all $c \in \mathcal{C}$,
  $$
    \xymatrix@R=2pt@C=3em{
      &
      L \circ R \circ L (c)
      \ar[dr]^-{ \epsilon_{L(c)} }
      \\
      L(c)
      \ar[ur]^-{ L(\eta_c) }
      \ar@{=}[rr]
      &&
      L(c)
      \,;
    }
  $$

  \noindent
  {\bf (ii)} for all $d \in \mathcal{D}$,
  $$
    \xymatrix@R=2pt@C=3em{
      &
      R \circ L \circ R (d)
      \ar[dr]^-{ R(\epsilon_d) }
      \\
      R(d)
      \ar[ur]^-{ \eta_{R(d)} }
      \ar@{=}[rr]
      &&
      R(d)
      \,.
    }
  $$
\end{prop}

\begin{prop}[Right/left adjoints preserve limits/colimits {\cite[5.2.3.5]{Lurie09}}]
  \label{AdjointsPreserveCoLimits}
  Let
 $
    \xymatrix{
      \mathcal{D}
        \ar@{<-}@<+7pt>[r]|-{\, L \,}
        \ar@<-5pt>[r]|-{\, R \,}^-{
          \raisebox{1pt}{
            \scalebox{.6}{$
              \bot
            $}
          }
        }
        &
      \mathcal{C}
    }
  $
   be a pair of adjoint $\infty$-functors
  (Def. \ref{AdjointInfinityFunctors}) and let
  $\mathcal{I} \in \mathrm{Categories}_\infty$\,.
  \begin{itemize}
  \vspace{-2mm}
    \item[{\bf (i)}]
      If $X_\bullet : \!\! \xymatrix@C=10pt{\mathcal{I} \ar[r] & \mathcal{D}}$
      is a diagram whose limit exists, then this limit is preserved by
      the right adjoint $R$:
      \vspace{-1mm}
      \begin{equation}
        \label{RightAdjointPreservesLimit}
        R\big(
          \underset{\longleftarrow}{\mathrm{lim}}
          \,
          X_\bullet
        \big)
        \;\simeq\;
          \underset{\longleftarrow}{\mathrm{lim}}
          \,
          R X_\bullet
      \end{equation}

      \vspace{-6mm}
    \item[{\bf (ii)}]
      If $X_\bullet : \!\! \xymatrix@C=10pt{\mathcal{I} \ar[r] & \mathcal{C}}$
      is a diagram whose colimit exists, then this colimit is preserved by
      the left adjoint $L$:
      \vspace{-1mm}
\noindent
      \begin{equation}
        \label{PreservedColimitLeftAdjoint}
        L\big(
          \underset{\longrightarrow}{\mathrm{lim}}
          \,
          X_\bullet
        \big)
        \;\simeq\;
          \underset{\longrightarrow}{\mathrm{lim}}
          \,
          L X_\bullet
      \end{equation}
  \end{itemize}
\end{prop}
Conversely:

\begin{prop}[Adjoint $\infty$-functor theorem {\cite[5.5.2.9]{Lurie09}}]
  \label{AdjointFunctorTheorem}
  Let $\mathcal{C}_{1,2} \in \mathrm{Categories}_\infty$
  be presentable (e.g. $\infty$-toposes, Def. \ref{InfinityTopos}),
  then an $\infty$-functor $\xymatrix{\mathcal{C}_1 \ar[r] & \mathcal{C}_2}$ is a:

  \noindent {\bf (i)} right adjoint (i.e., has a left adjoint, Def. \ref{AdjointInfinityFunctors})
  precisely if it preserves limits \eqref{RightAdjointPreservesLimit};

  \noindent {\bf (ii)}  left adjoint (i.e., has a right adjoint, Def. \ref{AdjointInfinityFunctors})
  precisely if it preserves colimits \eqref{PreservedColimitLeftAdjoint}.
\end{prop}

\begin{prop}[Fully faithful adjoints {\cite[5.2.7.4]{Lurie09}}]
  \label{CharacterizationOfFullyFaithfulAdjoints}
 For adjoint $\infty$-functors (Def. \ref{AdjointInfinityFunctors})
 $
    \xymatrix{
      \mathcal{D}
        \ar@{<-}@<+7pt>[r]|-{\, L \,}
        \ar@<-5pt>[r]|-{\, R \,}^-{
          \raisebox{1pt}{
            \scalebox{.6}{$
              \bot
            $}
          }
        }
        &
      \mathcal{C}
    }
  $,
  \begin{itemize}
   \vspace{-2mm}
    \item[{\bf (i)}]  $L$ is fully faithful
    $\xymatrix{\mathcal{D} \ar@{<-^{)}}[r]^-{L} & \,\mathcal{C} }$
     (Def. \ref{FullyFaithfulFunctor})
    iff
      the adjunction unit $\eta$ \eqref{AdjunctionUnit} is an equivalence:
      $\xymatrix@C=12pt{
        \mathrm{id}
        \ar[r]^-{\eta}_-{ \simeq }
        &
        R \circ L
      }$;
       \vspace{-5mm}
    \item[{\bf (ii)}]  $R$ is fully faithful
    $\xymatrix{\mathcal{D} \; \ar@{^{(}->}[r]^-{R} & \mathcal{C} }$
    (Def. \ref{FullyFaithfulFunctor})
    iff
      the adjunction counit $\epsilon$ \eqref{CounitOfAdjunction}
      is an equivalence
      $\xymatrix@C=12pt{
        L \circ R \ar[r]^-{\epsilon}_-{ \simeq } & \mathrm{id}
      }$.
  \end{itemize}
\end{prop}

\begin{prop}[Idempotent Monads and Comonads]
  \label{IdempotentMonads}
  For
  $
    \xymatrix{
      \mathcal{D}
        \ar@{<-}@<+7pt>[r]|-{\, L \,}
        \ar@<-5pt>[r]|-{\, R \,}^-{
          \raisebox{1pt}{
            \scalebox{.6}{$
              \bot
            $}
          }
        }
        &
      \mathcal{C}
    }
  $
  a pair of adjoint $\infty$-functors (Def. \ref{AdjointInfinityFunctors}):

\vspace{1mm}
  \noindent {\bf (i)} If $R$ is fully faithful (Def. \ref{FullyFaithfulFunctor})
  then $\Circle := R \circ L$ is idempotent, exhibited by the
  $\Circle$-image of the adjunction unit $\eta$ \eqref{AdjunctionUnit}:

  \vspace{-.8cm}

  \begin{equation}
    \label{IdempotencyOfComonad}
    \xymatrix{
      \Circle (c)
      \ar[rr]^-{ \Circle( \eta_{L(c)} ) }_-{ \simeq }
      &&
      \Circle \circ \Circle (c)\;.
    }
  \end{equation}

  \noindent {\bf (ii)} If $L$ is fully faithful (Def. \ref{FullyFaithfulFunctor})
  then
  $\Box := L \circ R$ is idempotent, exhibited by the
  $\Box$-image of the adjunction counit $\eta$ \eqref{CounitOfAdjunction}:

  \vspace{-.8cm}

  \begin{equation}
    \label{IdempotencyOfMonad}
    \xymatrix{
      \Box \circ \Box (d)
      \ar[rr]^{ \Box( \epsilon_{R(d)} ) }_-{ \simeq }
      &&
      \Box(d)\;.
    }
  \end{equation}
\end{prop}
\begin{proof}
  Consider case {\bf (i)}, the other case is formally dual.
  Since $R$ is fully faithful, by assumption,
  the condition that
  $\Circle (\eta_{L(c)}) := R \circ L(\eta_{L(c)})$
  is an equivalence is equivalent to $L(\eta_{L(c)})$
  being an equivalence.
  But, by the triangle identity (Prop. \ref{TriangleIdentities}),
  we have that the composite $\epsilon_{L(L(c))} \circ L(\eta_{L(c)})$
  is an equivalence, while by Prop. \ref{CharacterizationOfFullyFaithfulAdjoints}
  the counit $\epsilon$ is a natural equivalence. By cancellation,
  this implies
  that $L(\eta_{L(c)})$ is an equivalence.
\hfill \end{proof}

\medskip

\noindent {\bf $\infty$-Toposes.}
\noindent For our purposes, we take the following characterization to be the definition
of $\infty$-toposes. This is due to Rezk and Lurie \cite[6.1.6.8]{Lurie09};
we follow the presentation in \cite[Prop. 2.2]{NSS12}:

\begin{defn}[$\infty$-topos]
  \label{InfinityTopos}
  An $\infty$-topos $\mathbf{H}$ is a presentable $\infty$-category
  with the following properties:
  \begin{enumerate}[{\bf (i)}]
   \vspace{-2mm}
    \item {\bf Universal colimits.}
    For all morphisms $f : X \longrightarrow B$
    and all small diagrams $A : I \longrightarrow \mathbf{H}_{/B}$,
    there is an equivalence:
     \vspace{-2mm}
    \begin{equation}
      \label{PreserveColimitsByPullback}
      \underset{
        \underset{
          i
        }{\longrightarrow}
      }{\lim}
      f^\ast A_i
      \;\simeq\;
      f^\ast
      \big(
        \underset{
          \underset{
            i
          }{\longrightarrow}
        }{\lim}
        A_i
      \big)
    \end{equation}

     \vspace{-2mm}
\noindent
    between the pullback \eqref{PullbackSquare} of the colimit and the colimit
    over the pullbacks of its components.
     \vspace{-3mm}
    \item {\bf Univalent universes.} For every sufficiently large
    regular cardinal $\kappa$, there exists a morphism
    $\widehat{\mathrm{Objects}}_\kappa \longrightarrow \mathrm{Objects}_\kappa$
    in $\mathbf{H}$, such that for every object $X \in \mathbf{H}$,
    pullback \eqref{PullbackSquare}
    along morphisms $X \longrightarrow \mathrm{Objects}_\kappa$
    constitutes an equivalence
    \vspace{-3mm}
    \begin{equation}
      \label{ObjectClassifier}
      \xymatrix@R=-2pt{
        \mathrm{Core}
        \big(
          \mathbf{H}_{/_{\!\!\kappa} X}
        \big)
        \ar@{}[r]|-{\simeq}
        &
        \mathbf{H}\big(X, \mathrm{Objects}_\kappa \big)
        \\
        E \ar@{}[r]|-{ \longmapsto } & \vdash E
      }
      \phantom{AAAAA}
      \raisebox{20pt}{
      \xymatrix@C=5em@R=1em{
        E
          \ar@{}[dr]|-{\mbox{\tiny\rm(pb)}}
          \ar[r]
          \ar[d]
          &
        \widehat{\mathrm{Objects}}_\kappa
        \ar[d]
        \\
        X
          \ar[r]_-{ \vdash E }
          &
        \mathrm{Objects}_\kappa
      }
      }
    \end{equation}

     \vspace{-2mm}
\noindent
    between the $\infty$-groupoid core \eqref{InfinityCategoryOfInfinityCategories}
    of bundles (Notation \ref{SliceCategory})
    which are $\kappa$-small over $X$,
    and the hom-$\infty$-groupoid \eqref{HomInfinityGroupoids}
    of morphisms from $X$ to
    the \emph{object classifier} $\mathrm{Objects}_\kappa$.
  \end{enumerate}
\end{defn}

\begin{example}[Internal mapping space in an $\infty$-topos]
  Let $\mathbf{H}$ be an $\infty$-topos (Def. \ref{InfinityTopos})
  and $X \in \mathbf{H}$ an object. As a special case of
  universality of colimits \eqref{PreserveColimitsByPullback}, we
  have that the functor $X \times (-)$ of Cartesian product with $X$
  preserves all colimits. Hence, by the adjoint $\infty$-functor theorem
  (Prop. \ref{AdjointFunctorTheorem}),
  this functor has a right adjoint, to be denoted $\mathbf{Maps}(X,-)$,
  the \emph{internal hom}- or \emph{internal mapping space}-
  or \emph{mapping stack}-functor:
  \vspace{-2mm}
  \begin{equation}
    \label{InternalHomAdjunction}
    \xymatrix{
      \mathbf{H}\;\;
      \ar@{<-}@<+5pt>[rrr]^-{ X \times (-) }
      \ar@<-5pt>[rrr]_-{
        \underset{
          \mathclap{
          \raisebox{-3pt}{
            \tiny
            \color{darkblue}
            \bf
            \begin{tabular}{c}
              internal
              mapping space
            \end{tabular}
          }
          }
        }{
          \mathbf{Maps}(X,-)
        }
      }^-{ \bot }
      &&&
    \;\;  \mathbf{H}\;.
    }
  \end{equation}

  \vspace{-2mm}
\noindent
  By adjointness, the probes of the internal mapping space
  over any $U \in \mathbf{H}$ are given by
  \begin{equation}
    \label{ProbesOfInternalHom}
    \mathbf{H}
    \big(
      U, \mathbf{Maps}(X,Y)
    \big)
    \;\simeq\;
    \mathbf{H}
    \big(
      U \times X
      \,,\,
      Y
    \big).
  \end{equation}
\end{example}


\begin{prop}[Colimits and equifibered transformations {\cite[6.1.3.9(4)]{Lurie09}\cite[6.5]{Rezk10}}]
  \label{ColimitsOfEquifiberedTransformations}
  Let $\mathbf{H}$ be an $\infty$-topos (Def. \ref{InfinityTopos}),
  $\mathcal{I}$ a small $\infty$-category,
  $X_\bullet, Y_\bullet : \xymatrix@C=10pt{\mathcal{I} \ar[r] & \mathbf{H}}$
  two $\mathcal{I}$-shaped diagrams.

  \noindent {\bf (i)}
  If
  $\xymatrix@C=25pt{X_\bullet \ar@{=>}[r]|-{\;f_\bullet\;} & Y_\bullet}$
  is a natural transformation which
  is \emph{equifibered} \cite[p. 9]{Rezk10}, in that
  its value on
  all morphisms $\xymatrix@C=25pt{i_1 \ar[r]|-{\;\phi\;} & i_2}$ in $\mathcal{Y}$
  is a Cartesian square (Notation \ref{CartesianSquares}),
  then the value of $\underset{\longrightarrow}{\mathrm{lim}} \, f_\bullet$
  on all colimit component morphisms is also Cartesian:
  \vspace{-4mm}
  \begin{equation}
    \label{ColimitOfEquifiberedTransformationIsEquifibered}
    \underset{
      i_1 \overset{\phi}{\to} i_2
    }{\forall}
    \;\;\;\;\;
    \raisebox{20pt}{
    \xymatrix@R=1em@C=4em{
      X_{i_1}
      \ar[r]^-{ f_{i_1} }
      \ar[d]_-{X_\phi}
      \ar@{}[dr]|-{\mbox{\tiny\rm(pb)}}
      &
      Y_{i_1}
      \ar[d]^-{Y_\phi}
      \\
      X_{i_2}
      \ar[r]_-{ f_{i_2} }
      &
      Y_{i_2}
    }
    }
    \;\;\;\;\;\;\;\;\;\;\;\;\;
    \Rightarrow
    \;\;\;\;\;\;\;\;\;\;\;\;\;
    \underset{
      i
    }{\forall}
    \;\;\;\;\;\;
    \raisebox{22pt}{
    \xymatrix@R=1em@C=4em{
      X_{i}
      \ar[r]^-{ f_{i} }
      \ar[d]_-{ q_{X_i} }
      \ar@{}[dr]|-{\mbox{\tiny\rm(pb)}}
      &
      Y_{i_1}
      \ar[d]^-{ q_{Y_i} }
      \\
      \underset{\longrightarrow}{\mathrm{lim}}
      \, X_\bullet
      \ar[r]_-{
        \underset{\longrightarrow}{\mathrm{lim}}
        \,
        f_\bullet
      }
      &
      \underset{\longrightarrow}{\mathrm{lim}}
      \, Y_\bullet
    }
    }
  \end{equation}

  \vspace{-4mm}
  \noindent {\bf (ii)}
  Let $X^\rhd_\bullet : \xymatrix@C=10pt{\mathcal{I}^{\rhd} \ar[r] & \mathbf{H}}$
  be a cocone under $X_\bullet$, with tip $\mathcal{X} \in \mathbf{H}$,
  and let $Y^\rhd_\bullet : \xymatrix@C=10pt{\mathcal{I}^{\rhd} \ar[r] & \mathbf{H}}$
  denote the colimiting cocone under $Y_\bullet$ with tip
  $\underset{{\longrightarrow}}{\mathrm{lim}}Y_\bullet$.
  If $\xymatrix@C=12pt{X^{\rhd}_\bullet \ar@{=>}[r]^-{f^{\rhd}_\bullet} & Y^{\rhd}_\bullet}$
  is a natural transformation of cocone diagrams which is equifibered, then
  $X^\rhd_\bullet$ is a colimiting cocone:
  \vspace{-2mm}
  \begin{equation}
    \label{EquifiberedCoconeTransformationIntoColimitExhibitsColimit}
    \underset{
      i_1 \overset{\phi}{\to} i_2
    }{\forall}
    \;\;\;
    \raisebox{20pt}{
    \xymatrix@R=1.3em@C=4em{
      X_{i_1}
      \ar[r]^-{ f_{i_1} }
      \ar[d]_-{X_\phi}
      \ar@{}[dr]|-{\mbox{\tiny\rm(pb)}}
      &
      Y_{i_1}
      \ar[d]^-{Y_\phi}
      \\
      X_{i_2}
      \ar[r]_-{ f_{i_2} }
      &
      Y_{i_2}
    }
    }
    \;\;\;\;\;
    \mathrm{and}
    \;\;\;\;\;
    \underset{
      i
    }{\forall}
    \;\;\;
    \raisebox{22pt}{
    \xymatrix@R=1.3em@C=4em{
      X_{i}
      \ar[r]^-{ f_{i} }
      \ar[d]_-{ q_{X_i} }
      \ar@{}[dr]|-{\mbox{\tiny\rm(pb)}}
      &
      Y_{i_1}
      \ar[d]^-{ q_{Y_i} }
      \\
      \mathcal{X}
      \ar[r]_-{
        \underset{\longrightarrow}{\mathrm{lim}}
        \,
        f_\bullet
      }
      &
      \underset{\longrightarrow}{\mathrm{lim}}
      \, Y_\bullet
    }
    }
    \;\;\;\;\;\;\;\;\;\;\;\;
    \Rightarrow
    \;\;\;\;\;\;\;\;\;\;\;\;
    \mathcal{X}
    \;\simeq\;
    \underset{\longrightarrow}{\mathrm{lim}}\, X_\bullet
    \,.
  \end{equation}
\end{prop}
\begin{example}[Initial object in $\infty$-topos is empty object {\cite[p. 16]{Rezk19}}]
  \label{InitialObjectInInfinityToposIsEmpty}
  Let $\mathbf{H}$ be an $\infty$-topos (Def. \ref{InfinityTopos}).
Applying the implication
  \eqref{EquifiberedCoconeTransformationIntoColimitExhibitsColimit}
  in Prop. \ref{ColimitsOfEquifiberedTransformations} to the colimit
  over the empty diagram, which is the initial object, shows that
  any object with a morphism to the initial object is
  itself equivalent to the initial object.
  Hence if we write

  \vspace{-.6cm}

  \begin{equation}
    \label{InInfinityToposInitialObject}
    \varnothing \;\in\; \mathbf{H}
    \phantom{AAA}
    \mbox{s.t.}
    \phantom{AAA}
    \underset{X \in \mathbf{H}}{\forall}
    \big(
      \mathbf{H}(
        \varnothing
        \,,\,
        X
      )
      \;\simeq\;
      \ast
    \big)
  \end{equation}

  \vspace{-2mm}
\noindent
  for the initial object, this means that

  \vspace{-5mm}

  \begin{equation}
    \label{MorphismIntoInitialObjectImpliesDomainIsInitialObject}
    \xymatrix{
      X \ar[r]^-{\exists} & \varnothing
    }
    \phantom{AAAA}
    \Rightarrow
    \phantom{AAAA}
    X \;\simeq\; \varnothing
    \,.
  \end{equation}
\end{example}

\begin{prop}[Tensoring of $\infty$-toposes over $\infty$-groupoids]
  \label{TensoringOfInfinityToposesOverInfinityGroupoids}
  Let $\mathbf{H}$ be an $\infty$-topos (Def. \ref{InfinityTopos})
  with inverse base geometric morphism (Prop. \ref{BaseGeometricMorphism})
  denoted
  $\Delta : \mathrm{Groupoids}_\infty \longrightarrow \mathbf{H}$\,.
  Then, for $S \in \mathrm{Groupoids}_\infty$ and $X,Y \in \mathbf{H}$,
  there is a natural equivalence of $\infty$-groupoids

  \vspace{-.4cm}

  \begin{equation}
    \label{InfinityTensoring}
    \mathbf{H}
    \big(
      \Delta(S) \times X
      \,,\,
      Y
    \big)
    \;\simeq\;
    \mathrm{Groupoids}_{\infty}
    \big(
      S
      \,,\,
      \mathbf{H}(X,Y)
    \big).
  \end{equation}
\end{prop}
\begin{proof}
  By \cite[Cor. 4.4.4.9]{Lurie09} we have, for
  $S \in \mathrm{Groupoids}_\infty \hookrightarrow \mathrm{Categories}_\infty$
  and $X, Y \in \mathbf{H}$, natural equivalences
  \vspace{-1.5mm}
  \begin{equation}
    \label{InfinityGroupoidIsHomotopyColimitOverItself}
    \underset{
      \underset{S}{\longrightarrow}
    }{\mathrm{lim}}
    \;
    \mathrm{const}_{\ast}
    \;\simeq\;
    S
    \phantom{AAA}
    \mbox{and}
    \phantom{AAA}
    \mathbf{H}
    \Big(
      \underset{
        \underset{
          S
        }{\longrightarrow}
      }{\mathrm{lim}}
      \,
      \mathrm{const}_X
      \,,\,
      Y
    \Big)
    \;\simeq\;
    \mathrm{Groupoids}_\infty
    \big(
      S
      \,,\,
      \mathbf{H}(X,Y)
    \big).
  \end{equation}

   \vspace{-3mm}
\noindent
  This implies the statement in the form \eqref{InfinityTensoring}
  by using {\bf (a)} that $\Delta$ preserves all colimits as
  well as finite limits (Prop. \ref{BaseGeometricMorphism})
  and {\bf (b)} that Cartesian products
  may be taken inside colimits, as a special case of \eqref{PreserveColimitsByPullback}:
\vspace{-.5mm}
  $$
    \begin{aligned}
      \mathbf{H}
      \big(
        \Delta(S)
        \times X
        \,,\,
        Y
      \big)
      & \simeq
      \mathbf{H}
      \big(
        \Delta
        \big(
          \underset{
            \underset{S}{\longrightarrow}
          }{\mathrm{lim}}
          \,
          \ast
        \big)
        \times X
        \,,\,
        Y
      \big)
      \; \simeq
      \mathbf{H}
      \big(
        \big(
          \underset{
            \underset{S}{\longrightarrow}
          }{\mathrm{lim}}
          \,
          \underset{
            \mathclap{
              \simeq \, \ast
            }
          }{
          \underbrace{
            \Delta(\ast)
          }
          }
        \,\big)
        \times X
        \,,\,
        Y
      \big)
      \\
      & \simeq
      \mathbf{H}
      \Big(
        \big(
          \underset{
            \underset{S}{\longrightarrow}
          }{\mathrm{lim}}
          \,
          \underset{
            \simeq \, X
          }{
          \underbrace{
            (\ast \times X)
          }
          }
        \,\big)
        \,,\,
        Y
      \Big)
      \; \simeq
      \mathrm{Groupoids}_\infty
      \big(
        S
        \,,\,
        \mathbf{H}(X,Y)
      \big).
    \end{aligned}
  $$
  \vspace{-2.5mm}
  \noindent The composite equivalence is \eqref{InfinityTensoring}.
\hfill \end{proof}

\vspace{1mm}
\noindent {\bf Sheaves.}
\begin{notation}[$\infty$-Presheaves]
  \label{Presheaf}
  For $\mathcal{C}$ a small $\infty$-category, we write
  \vspace{-1.5mm}
  \begin{equation}
    \label{Presheaves}
    \mathrm{PreSheaves}_\infty(\mathcal{C})
    \;:=\;
    \mathrm{Functors}_\infty
    \big(
      \mathcal{C}^{\mathrm{op}}
      \,,\,
      \mathrm{Groupoids}_\infty
    \big)
  \end{equation}

  \vspace{-2mm}
\noindent
  for the $\infty$-category of $\infty$-presheaves on $\mathcal{C}$.
  More generally, if $\mathbf{H}$ is any $\infty$-topos
  (Def. \ref{InfinityTopos}) we also write
\vspace{-1.5mm}
  \begin{equation}
    \label{HValuedPresheaves}
    \mathrm{PreSheaves}_\infty
    \big(
      \mathcal{C}, \mathbf{H}
    \big)
    \;:=\;
    \mathrm{Functors}_\infty
    \big(
      \mathcal{C}^{\mathrm{op}}
      \,,\,
      \mathbf{H}
    \big).
  \end{equation}
\end{notation}

\begin{prop}[Limits and colimits in an $\infty$-topos {\cite[Lem. 4.2.4.3]{Lurie09}}]
  \label{LimitsAndColimitsAsAdjoints}
  Let $\mathbf{H}$ be an $\infty$-topos (Def. \ref{InfinityTopos})
  and $\mathcal{C}$ a small $\infty$-category.
  Then the $\infty$-functor which sends an object in $\mathbf{H}$
  to the $\mathbf{H}$-valued presheaf \eqref{HValuedPresheaves}
  constant on this object has a right- and a left-adjoint
  (Def. \ref{AdjointInfinityFunctors}),
  given by the limit and colimit construction,
  respectively:
  \vspace{-2mm}
  \begin{equation}
    \label{LimitAndColimitAsAdjoints}
    \xymatrix{
      \mathrm{Functors}_{\infty}
      \big(
        \mathcal{C}
        \,,\,
        \mathbf{H}
      \big)
      \;\;\;
      \ar@<+11pt>[rr]^-{ \underset{\longrightarrow}{\mathrm{lim}} }
      \ar@{<-}[rr]|-{\; \mathrm{const}\, }^-{ \raisebox{1pt}{\scalebox{.7}{$\bot$}} }_-{ \raisebox{-1pt}{\scalebox{.7}{$\bot$}} }
      \ar@<-11pt>[rr]_-{ \underset{\longleftarrow}{\mathrm{lim}} }
      &&
    \;\;\;  \mathbf{H}
    }
  \end{equation}
\end{prop}

\begin{prop}[$\infty$-Yoneda embedding {\cite[Lemma 5.5.2.1]{Lurie09}}]
  \label{InfinityYonedaEmbedding}
  Let $\mathcal{C}$ be an $\infty$-category. Then the
  $\infty$-functor from $\mathcal{C}$ to its $\infty$-presheaves
  \eqref{Presheaves} which assigns
  \emph{representable presheaves}
  \vspace{-2mm}
  \begin{equation}
    \label{TheInfinityYonedaEmbedding}
    \xymatrix@R=-4pt{
      \mathcal{C}
      \;
      \ar@{^{(}->}[r]^-{ y }
      & \;
      \mathrm{PreSheaves}_\infty(\mathcal{C})
      \\
      c \ar@{|->}[r]
      &
      \mathcal{C}(-,c)
    }
  \end{equation}

  \vspace{-2mm}
  \noindent  is fully faithful (Def. \ref{FullyFaithfulFunctor}).
\end{prop}
\begin{prop}[$\infty$-Yoneda lemma {\cite[Lemma 5.5.2.1]{Lurie09}}]
  \label{YonedaLemma}
  Let $\mathcal{C}$ be an $\infty$-category.
  Then for $X \in \mathrm{PreSheaves}_\infty(\mathcal{C})$
  \eqref{Presheaves} and
  $c \in \mathcal{C}$, there is a natural equivalence
  \vspace{-2mm}
  $$
    \mathrm{PreSheaves}_\infty
    \big(
      y(c), X
    \big)
    \;\simeq\;
    X(c)
    \,,
  $$

  \vspace{-2mm}
\noindent
  where $y$ is the Yoneda embedding \eqref{TheInfinityYonedaEmbedding}
  from Prop. \ref{InfinityYonedaEmbedding}.
\end{prop}

\begin{prop}[(Co-)Limits of presheaves are computed objectwise {\cite[Cor. 5.1.2.3]{Lurie09}}]
  \label{LimitsAndColimitsOfPresheavesComputedObjectwise}
  Let $\mathbf{H}$ be an $\infty$-topos,
  let $\mathcal{C}$ and $\mathcal{D}$  be small $\infty$-categories, and let
  \vspace{-2mm}
  $$
    I
      \;:\;
    \mathcal{D}   \xymatrix{
           \ar[r]
      &
}\mathrm{PreSheaves}_\infty(\mathcal{C}, \mathbf{H})
  $$

  \vspace{-2mm}
\noindent
  be a diagram of $\mathbf{H}$-valued $\infty$-presheaves over $\mathcal{C}$.
  Then the limit and colimit over $I$ exist and are given objectwise
  over $c \in \mathcal{C}$
  by the limit and colimit of the components in $\mathrm{Groupoids}_\infty$:
  \vspace{-1mm}
  \begin{align*}
    \big(\underset{\longrightarrow}{\mathrm{lim}}\,I\big)
   & \;:\;
    c
    \;\longmapsto\;
    \big(\underset{\longrightarrow}{\mathrm{lim}}\,I_c\big)
    \,,
 \\
    \big(\underset{\longleftarrow}{\mathrm{lim}}\,I\big)
  &  \;:\;
    c
    \;\longmapsto\;
    \big(\underset{\longleftarrow}{\mathrm{lim}}\,I_c\big)
    \,.
  \end{align*}
\end{prop}

\begin{lemma}[Colimit of representable functor is contractible]
  \label{InfinityColimitOverRepresentableInfinityFunctorIsContractible}
  Let $\mathcal{C}$ be a small $\infty$-category,
  and consider an
  $\infty$-functor
  $y C : \mathcal{C}^{\mathrm{op}} \longrightarrow \mathrm{Groupoids}_\infty$
  to the $\infty$-category of $\infty$-groupoids
  \eqref{InfinityCategoryOfInfinityGroupoids}, which is
  representable, hence which is in the essential image of the
  $\infty$-Yoneda embedding \eqref{TheInfinityYonedaEmbedding}.
  Then the colimit of this functor is contractible:
  \vspace{-1mm}
  \begin{equation}
    \underset{
      \underset{\mathcal{C}}{\longrightarrow}
    }{
      \mathrm{lim}
    }
    (y C)
    \;\simeq\;
    \ast \;.
  \end{equation}
\end{lemma}
\begin{proof}
  The terminal $\ast \in \mathrm{Groupoids}_\infty$
  is characterized by the fact that for $S \in \mathrm{Groupoids}_\infty$
  there is a natural equivalence
  \vspace{-1mm}
  $$
    S \;\simeq\;
    \mathrm{Groupoids}_\infty
    \big(
      \ast
      \,,\,
      S
    \big)
    \,.
  $$
  Hence it is sufficient to see that
  $\underset{\longrightarrow}{\mathrm{\lim}}\,(y C)$ satisfies the
  same property. But we have the following sequence of natural
  equivalences:
  $$
    \begin{aligned}
      \mathrm{Groupoids}_\infty
      \left(
        \underset{\longrightarrow}{\mathrm{lim}}\,
        (y C)
        \,,\,
        S
      \right)
      & \;\simeq\;
      \mathrm{Functors}_\infty
      \big(
        \mathcal{C}^{op}
      \big)
      \big(
        y C
        \,,\,
        \mathrm{const}\,
      \big)
      \\
      & \;\simeq\;
      (
        \mathrm{const}\, S
      )
      (C)
      \; \simeq\;
      S
      \,.
    \end{aligned}
  $$

  \vspace{-1mm}
\noindent
  Here the first step is the adjunction \eqref{LimitAndColimitAsAdjoints}, while the
  second step is the $\infty$-Yoneda lemma (Prop. \ref{YonedaLemma}).
\hfill \end{proof}

\begin{prop}[Topos is accessibly lex reflective in presheaves over site {\cite[6.1.0.6]{Lurie09}}]
  \label{ToposLexReflection}
  Let $\mathbf{H}$ be an $\infty$-topos (Def. \ref{InfinityTopos}).

 \noindent {\bf (i)} Then there exists an \emph{$\infty$-site}
  for $\mathbf{H}$, namely a small $\mathcal{C} \in \mathrm{Categories}_\infty$
  equipped with a pair of adjoint $\infty$-functors (Def. \ref{AdjointInfinityFunctors})
  between $\mathbf{H}$ and $\mathrm{PreSheaves}_\infty(\mathcal{C})$
  (Notation \ref{Presheaf}):
  \vspace{-2mm}
  \begin{equation}
    \label{ToposReflectionInPresheaves}
    \xymatrix{
      \mathbf{H} \;\;
      \ar@{<-}@<+8pt>[rr]^-{L}
      \ar@<-8pt>@{^{(}->}[rr]^-{
        \raisebox{2pt}{
          \scalebox{.9}
          {$
            \bot
          $}
        }
      }
      &&
     \; \mathrm{PreSheaves}_\infty
      \big(
        \mathcal{C}
      \big)
    }
  \end{equation}
  such that
  {\bf (a)} the right adjoint is
  accessible and fully faithful (Def. \ref{FullyFaithfulFunctor})
  and {\bf (b)} the left adjoint preserves finite limits
  (in addition to preserving all colimits, by Prop. \ref{AdjointsPreserveCoLimits}).

 \noindent {\bf (ii)}  Conversely, any such accessibly embedded lex reflective sub-$\infty$-category
  of an $\infty$-category of $\infty$-presheaves is an $\infty$-topos.
\end{prop}
\begin{defn}[Sheaf $\infty$-topos {\cite[6.2]{Lurie09}}]
  \label{InfinityToposOfInfinitySheaves}
  An $\infty$-topos $\mathbf{H}$ (Def. \ref{InfinityTopos})
  is called an \emph{$\infty$-category of $\infty$-sheaves}
  or \emph{of $\infty$-stacks}, or just a \emph{sheaf topos} for short,
  to be denoted
  \begin{equation}
    \label{SheafTopos}
    \mathbf{H}
    \;\simeq\;
    \mathrm{Sheaves}_\infty
    \big(
      \mathcal{C}
    \big)
  \end{equation}
  if there exists a \emph{site} $\mathcal{C}$,
  namely a small
  $\mathcal{C} \in \mathrm{Categories}_\infty$
  with a reflection $(L\mathrm{const} \dashv \Gamma)$ \eqref{ToposReflectionInPresheaves}
  as in Prop. \ref{ToposLexReflection}, such that
  $L\mathrm{const}$ exhibits localization at a set
  $$
    \Big\{
    \!\!
    \underset{
      \mathclap{
      \raisebox{+1pt}{
        \tiny
        \color{darkblue}
        \bf
        covering sieves
      }
      }
    }{
    \xymatrix{
        U
        \;
        \ar@{^{(}->}[r]
        &
        \;y(c)
    }
    }
    \!\!
    \Big\}
    \;\subset\;
    \underset{c \in \mathcal{C}}{\sqcup}
    \mathrm{SubObjects}
    \big(
      y(c)
    \big)
  $$
  of monomorphisms
  (Def. \ref{Monomorphism}) into representable presheaves \eqref{TheInfinityYonedaEmbedding}.
\end{defn}

\begin{prop}[Base geometric morphism {\cite[6.3.4.1]{Lurie09}}]
  \label{BaseGeometricMorphism}
  Let $\mathbf{H}$ be an $\infty$-topos (Def. \ref{InfinityTopos}).
  There is an essentially unique pair
  of adjoint $\infty$-functors (Def. \ref{AdjointInfinityFunctors})
  between $\mathbf{H}$ and $\mathrm{Groupoids}_\infty$
  (Def. \ref{InfinityGroupoids})
    \vspace{-2mm}
  \begin{equation}
    \label{BaseGeometricMorphismAdjunction}
    \xymatrix{
      \mathbf{H} \;\;
      \ar@{<-}@<+8pt>[rr]^-{L\mathrm{const}}
      \ar@<-8pt>[rr]^-{
        \raisebox{2pt}{
          \scalebox{.9}
          {$
            \bot
          $}
        }
      }_-{ \Gamma }
      &&
     \; \mathrm{Groupoids}_\infty
    }
  \end{equation}

    \vspace{-2mm}
\noindent  such that the left adjoint $L\mathrm{const}$ preserves
  finite limits (in addition to preserving all colimits, by Prop. \ref{AdjointsPreserveCoLimits}).
\end{prop}

\begin{example}[Base geometric morphism via site]
  Let $\mathbf{H}$ be an $\infty$-topos (Def. \ref{InfinityTopos})
  and $\mathcal{C}$ a site (Prop. \ref{ToposLexReflection}).
  Then the composite of pairs of adjoint $\infty$-functors
  (Def. \ref{AdjointInfinityFunctors})
    \vspace{-2mm}
  \begin{equation}
    \label{BaseGeometricMorphismViaFactorizationThroughPresheavesOnSite}
    \xymatrix{
      \mathbf{H} \;\;
      \ar@{<-}@<+8pt>[rr]^-{L}
      \ar@<-8pt>@{^{(}->}[rr]^-{
        \raisebox{2pt}{
          \scalebox{.9}
          {$
            \bot
          $}
        }
      }
      &&
      \; \mathrm{PreSheaves}_\infty
      \big(
        \mathcal{C}
      \big)
      \ar@{<-}@<+8pt>[rr]^-{\mathrm{const}}
      \ar@<-8pt>[rr]_-{
        \underset{\longleftarrow}{\mathrm{\lim}}
      }^-{
        \raisebox{2pt}{
          \scalebox{.9}
          {$
            \bot
          $}
        }
      }
      &&
     \; \mathrm{Groupoids}_\infty
    }
  \end{equation}

    \vspace{-2mm}
\noindent  of {\bf (a)} the reflection into presheaves over the site
  (Prop. \ref{ToposLexReflection})
  with {\bf (b)} the limit-construction on presheaves (Prop. \ref{LimitsAndColimitsAsAdjoints})
  is such that the composite left adjoint
  $L\mathrm{const}$
  preserves finite limits
  (since $L$ does by Prop. \ref{ToposLexReflection} and
  $\mathrm{const}$ does by Prop. \ref{AdjointsPreserveCoLimits} with
  Prop. \ref{LimitsAndColimitsAsAdjoints}).
  Hence, by the essential uniqueness of Prop. \ref{BaseGeometricMorphism},
  the composite
  \eqref{BaseGeometricMorphismViaFactorizationThroughPresheavesOnSite}
  is a factorization of the base geometric morphism of $\mathbf{H}$.
\end{example}

\medskip

\noindent {\bf Bundles.}
\begin{notation}[Bundles and slicing.]
  \label{SliceCategory}
  Let $\mathbf{H}$ an $\infty$-topos (Def. \ref{InfinityTopos})
  and $X \in \mathbf{H}$ an object. We write:

  \noindent {\bf (i)} $(X,p) \;\in\; \mathbf{H}_{/X}$
  for objects in the slice $\infty$-category of $\mathbf{H}$
  over $X$, corresponding to morphisms $p$ to $X$ in $\mathbf{H}$
  (\emph{bundles} over $X$):
  \vspace{-4mm}
  $$
    \xymatrix@R=11pt{
      E
      \ar[d]^-{p}
      \\
      X
    }
  $$

  \vspace{-2mm}
  \noindent {\bf (ii)} $(f, \alpha) \in \mathbf{H}_{/_{X}}
  \big(
    (E_1, p_1)
    \,,\,
    (E_2, p_2)
  \big)$
  for morphisms in the slice $\infty$-category, corresponding
  to diagrams in $\mathbf{H}$ of the form
  \vspace{-2mm}
  \begin{equation}
    \label{MorphismInSliceCategory}
    \raisebox{10pt}{
    \xymatrix@R=1em@C=3em{
      E_1
      \ar[dr]_-{p_1}^{\ }="s"
      \ar[rr]^-{ f }_-{\ }="t"
      &&
      E_2
      \ar[dl]^-{ p_2 }
      \\
      &
      X
      \ar@{=>}_\alpha "s"; "t"
    }
    }
  \end{equation}
\end{notation}

\begin{prop}[Slice $\infty$-topos {\cite[Prop. 6.3.5.1 (1)]{Lurie09}}]
  \label{SliceInfinityTopos}
  Let $\mathbf{H}$ be an $\infty$-topos (Def. \ref{InfinityTopos})
  and $X \in \mathbf{H}$ an object. Then the slice $\infty$-category
  $\mathbf{H}_{/X}$ (Notation \ref{SliceCategory}) is also an $\infty$-topos.
\end{prop}

\begin{example}[Iterated slice $\infty$-topos]
  \label{IteratedSliceTopos}
  Let $\mathbf{H}$ be an $\infty$-topos (Def. \ref{InfinityTopos}),
  $X \in \mathbf{H}$ and $(Y,p) \in \mathbf{H}_{/X}$ an object in the
  slice, hence (Notation \ref{SliceCategory}) a morphism
  $\xymatrix@C=16pt{Y\!\! \ar[r]|<<<{p} & X}$.
    Then $\mathbf{H}_{/X}$ is itself an $\infty$-topos, by
  Prop. \ref{SliceInfinityTopos}, and we may slice again to obtain the
  iterated slice $\infty$-topos

  \vspace{-.8cm}

  \begin{equation}
    \label{IteratedSliceCategory}
    \left(
      \mathbf{H}_{/X}
    \right)_{/(Y,p)}
    \;\;
    \in
    \mathrm{Categories}_\infty\;.
  \end{equation}
  \vspace{-3mm}
  \item   {\bf (i)} an object in \eqref{IteratedSliceCategory}
  is a diagram in $\mathbf{H}$ of this form:
  $ \;
    \raisebox{10pt}{
    \xymatrix@R=-2pt{
      Z
      \ar[dddr]^-{\ }="s"
      \ar[drr]_>>>>{\ }="t"
      \\
      &&
      Y
      \ar[ddl]^{p}
      \\
      \\
      & X
      \ar@{=>} "s"; "t"
    }
    }
  $
\item
  \hspace{-.4cm}
  \begin{tabular}{l}
    {\bf (ii)} a morphism in \eqref{IteratedSliceCategory}
    is a diagram in $\mathbf{H}$ of this form:
    \\
    \phantom{\bf (ii)} (This is furthermore filled by a 3-morphism,
    \\
    \phantom{\bf (ii)} which we notationally suppress, for readability.)
  \end{tabular}
 $
    \raisebox{20pt}{
    \xymatrix@R=20pt@C=5em{
      Z_1
      \ar[ddr]^>>>>>{\ }="s1"
      \ar@{}[ddr]|-{
      }^>>>>>>>{\ }="s4"
      \ar@/_.46pc/[drr]|>>>>>>>{}_->>>>>{\ }="t1"^-<<<<<<{\ }="t2"
      \ar[r]^-{}|-{}
      &
      Z_2
      \ar@/^.55pc/[dr]|>>>>>>{}_>>>>>>{\ }="s2"
      \ar[dd]|<<<{\phantom{A\vert}}|<<<<<<{\phantom{AA}}|>>>>>{\phantom{AA} \atop \phantom{a}}^<<<{\ }="s3"_-{\ }="t4"
      \\
      &&
      Y
      \\
      &
      X
      \ar@{<-}[ur]_-{ p }^-{\ }="t3"
      \ar@/_.23pc/@{=>} "s1"; "t1"|<<<<<<<<<<<<<<<<<{}
      \ar@{=>}_<<<<<<<<<<<<<{} "s2"; "t2"
      \ar@<+2pt>@{=>} "t3"+(-10,4); "t3"+(-1.5,-1)|>>>>>>{\phantom{AA}}
      \ar@{=>} "s4"; "t4"
    }
    }
$

\end{example}

\begin{prop}[Hom-$\infty$-groupoids in slices {\cite[Prop. 5.5.5.12]{Lurie09}}]
  \label{HomsInSliceInfinityCategory}
  Let $\mathbf{H}$ be an $\infty$-topos (Def. \ref{InfinityTopos})
  and $B \in \mathbf{H}$ an object.
  Then for $(X_1, p_1), (X_2, p_2) \in \mathbf{H}_{/B}$
  two objects in the slice over $B$ (Prop. \ref{SliceInfinityTopos})
  the hom-$\infty$-groupoid between them is given by the following
  homotopy fiber-product of hom-$\infty$-groupoids of $\mathbf{H}$:
  \vspace{-2mm}
  \begin{equation}
    \label{HomGroupoidInSlice}
    \mathbf{H}_{/B}
    \big(
      (X_1, p_1)
      \,,\,
      (X_2, p_2)
    \big)
    \;\;\simeq\;\;
    \{p_1\}
    \underset{
      \mathbf{H}(X_1,B)
    }{\times}
    \mathbf{H}(X_1, X_2)
  \end{equation}

  \vspace{-2mm}
\noindent
  hence by the $\infty$-groupoid given by the following
  Cartesian square (Notation \ref{CartesianSquares}):
  \vspace{-2mm}
  $$
    \xymatrix@R=1.5em{
      \mathbf{H}_{/B}
      \big(
        (X_1, p_1)
        \,,\,
        (X_2, p_2)
      \big)
      \ar[d]
      \ar[rr]
      \ar@{}[drr]|-{\mbox{\tiny\rm(pb)}}
      &&
      \mathbf{H}(X_1, X_2)
      \ar[d]^-{ p_2 \circ (-) }
      \\
      \ast
      \ar[rr]_-{ \vdash p_1 }
      &&
      \mathbf{H}(X_1, B)
    }
  $$
\end{prop}

\vspace{-4mm}
\begin{prop}[Base change {\cite[HTT 6.3.5]{Lurie09}}]
  \label{BaseChange}
  Let $\mathbf{H}$ be an $\infty$-topos (Def. \ref{InfinityTopos}).
  Then for every morphism $X \overset{f}{\to} Y$ in $\mathbf{H}$
  there is an induced \emph{base change} adjoint triple
  (Def. \ref{AdjointInfinityFunctors}) between
  the corresponding slice $\infty$-toposes (Prop. \ref{SliceInfinityTopos}):
  \vspace{-4mm}
  \begin{equation}
    \label{AdjointTripleBaseChange}
    \xymatrix{
      \mathbf{H}_{/X}
     \;\; \ar@<+14pt>[rr]^-{ f_! }
      \ar@{<-}[rr]|-{ f^\ast }^-{\raisebox{4pt}{\footnotesize $\bot$}}
      \ar@<-14pt>[rr]_-{ f_\ast }^-{\footnotesize \bot}
      &&
   \;\;   \mathbf{H}_{/Y}
    }
  \end{equation}

    \vspace{-1mm}
\noindent
  where,   in $\mathbf{H}$,
 $f_!$ is given by postcomposition with $f$ while
$f^\ast$ is given by pullback along $f$.
\end{prop}

\begin{example}[Bundle morphisms covering base morphisms]
  \label{BundleMorphismsCoveringBaseMorphisms}
  For $\mathbf{H}$ an $\infty$-topos (Def. \ref{InfinityTopos}),
  the system of all its slice $\infty$-toposes (Prop. \ref{SliceInfinityTopos})

  \vspace{-.6cm}

  \begin{equation}
    \label{SystemOfSliceToposes}
    \xymatrix@R=-2pt{
      \mathbf{H}^{\mathrm{op}}
      \ar[rr]^-{ \mathbf{H}_{\!/(-)} }
      &&
      \mathrm{Categories}_\infty
      \\
      X
      \ar@{}[rr]|-{\longmapsto}
      &&
      \mathbf{H}_{\!/X}
    }
  \end{equation}

  \vspace{-1mm}
\noindent  related via contravariant base change \eqref{AdjointTripleBaseChange}
  arranges into the ``arrow $\infty$-topos'' \cite[2.4.7.12]{Lurie09}

  \vspace{-.3cm}

  \begin{equation}
    \label{InfinityToposOfBundles}
    \mathrm{Bundles}(\mathbf{H})
    \;:=\;
    \int_X \mathbf{H}_{/X}
    \;\simeq\;
    \mathbf{H}^{\Delta[1]}
    \,,
  \end{equation}
  which, in view of Notation \ref{SliceCategory}, may be thought of
  as the $\infty$-category of bundles in $\mathbf{H}$,
  but now with
  bundle morphisms allowed to cover non-trivial base morphisms.
\end{example}

\begin{example}[Spectral bundles and tangent $\infty$-topos]
  \label{TangentInfinityTopos}
  Let $\mathbf{H}$ be an $\infty$-topos (Def. \ref{InfinityTopos}).
  Instead of the system
  \eqref{SystemOfSliceToposes}
  of its plain slices, consider the corresponding system
  of \emph{stabilized} slices
  (stabilized under the suspension/looping adjunction
  on pointed objects, e.g. \cite[1.4]{Lurie07}):
  \vspace{-1mm}
  \begin{equation}
    \label{SystemOfSliceToposes}
    \xymatrix@R=-2pt{
      \mathbf{H}^{\mathrm{op}}
      \ar[rr]^-{ \mathrm{Stab}\left(\mathbf{H}_{\!/(-)}\right) }
      &&
      \mathrm{Categories}_\infty
      \\
      X
      \ar@{}[rr]|-{\longmapsto}
      &&
      \mathrm{Stab}\big(\mathbf{H}_{\!/X}\big)
    }
  \end{equation}

  \vspace{-.7cm}

 \noindent  The resulting total $\infty$-category
  \begin{equation}
    \label{InfinityToposOfSpectralBundles}
    \mathrm{SpectralBundles}(\mathbf{H})
    \;:=\;
    \int_X \mathrm{Stab}\big(\mathbf{H}_{/X}\big)
    \,,
  \end{equation}
  is that of \emph{bundles of spectra} in $\mathbf{H}$
  (parametrized spectrum objects).
  Remarkably, this is itself an $\infty$-topos
  \cite[35.5]{Joyal08b}\cite[6.1.1.11]{Lurie17},
  also called the \emph{tangent $\infty$-topos} $T \mathbf{H}$ of $\mathbf{H}$
  (e.g. \cite{Lurie07}\cite{BM19}).
\end{example}

\begin{example}[Base change along terminal morphism]
  \label{BaseChangeAlongTerminalMorphism}
  Let $\mathbf{H}$ be an $\infty$-topos (Def. \ref{InfinityTopos})
  and $X \in \mathbf{H}$ any object.
  With $\mathbf{H} \simeq \mathbf{H}_{/\ast}$ regarded
  as its own slice (Prop. \ref{SliceInfinityTopos}) over the terminal object,
  base change (Prop. \ref{BaseChange}) along the terminal
  morphism $X \to \ast$ is of the form
  \vspace{-4mm}
  \begin{equation}
    \label{BaseChangeToGlobalContext}
    \xymatrix{
      \mathbf{H}_{/X}
     \;\; \ar@<+15pt>[rr]^-{ \mathrm{dom} }
      \ar@{<-}[rr]|-{\, X \times (-) \,}^-{\raisebox{4pt}{\footnotesize $\bot$}}
      \ar@<-15pt>[rr]_-{  }^-{\bot}
      &&
   \;\;   \mathbf{H}
    }
  \end{equation}
  where {\bf (a)} the top functor
  sends a morphism $Y \to X$ to its domain object $Y$,
  and {\bf (b)} the middle functor is Cartesian product with $X$.
  In particular, it follows that:

  \noindent
  {\bf (i)} The base geometric morphism (Prop. \ref{AdjointInfinityFunctors})
  of the slice $\infty$-topos
  $\mathbf{H}_{/X}$ (Prop. \ref{SliceInfinityTopos})
  is given by
  \begin{equation}
    \label{BaseGeometricMorphismForSliceTopos}
    \big(
      \Delta \,\dashv\, \Gamma
    \big)
    \;\simeq\;
    \big(
      (X \to \ast)^\ast \,\dashv\, (X \to \ast)_\ast
    \big)
  \end{equation}
  (since $(X \to \ast)^\ast$ is a left adjoint that also preserves finite
  limits, as it is also a right adjoint, Prop. \ref{AdjointsPreserveCoLimits}).

  \noindent {\bf (ii)}
  The forgetful functor $\mathrm{dom} : \mathbf{H}_{/X} \to \mathbf{H}$
  is a left adjoint $(X \to \ast)_!$ and hence preserves all colimits
  (Prop. \ref{AdjointsPreserveCoLimits}).
\end{example}
While $\mathrm{dom}$ \eqref{BaseChangeToGlobalContext}
does not preserve all limits,
it does preserve fiber products:
\begin{prop}[Fiber products in slice $\infty$-toposes]
  \label{HomotopyFiberProductsInOverToposes}
  Let $\mathbf{H}$ be an $\infty$-topos (Def. \ref{InfinityTopos}),
  $B \in \mathbf{H}$, $\mathbf{H}_{/B}$ the slice $\infty$-topos
  (Prop. \ref{SliceInfinityTopos}) and
  $\xymatrix@C=15pt{\mathbf{H}_{/B} \ar[r]^{\scalebox{0.6}{{\rm dom}}} & \mathbf{X} }$
  its forgetful functor \eqref{BaseChangeToGlobalContext}
  from Example \ref{BaseChangeAlongTerminalMorphism}.

\vspace{-2mm}
  \noindent
  {\bf (i)} Given a cospan
  $\!\!\!\xymatrix@C=11pt{ (Y, \phi_Y) \ar[r] & (X, \phi_X) \ar@{<-}[r] & (Z,\phi_Z) }\!\!$
  in $\mathbf{H}_{/B}$, the underlying object of its fiber product is
  the fiber product of its underlying objects:
  \vspace{-5mm}
  \begin{equation}
    \label{DomainOfFiverProductInSlice}
    \mathrm{dom}
    \left(
      (Y,\phi_Y)
      \underset{
        (X,\phi_X)
      }{\times}
      (Z,\phi_Z)
    \right)
    \;\;\simeq\;\;
    Y \underset{X}{\times} Z
    \,.
  \end{equation}

  \vspace{-1mm}
  \noindent {\bf (ii)} In particular, since $(X,\mathrm{id}_X)$ is the terminal
  object in $\mathbf{H}_{/X}$, so that the plain product in the slice is
  \vspace{-1mm}
  $$
    (Y,\phi_Y) \times (Z,\phi_Z)
    \;\;=\;\;
    (Y,\phi_Y)
      \underset{
        (X,\mathrm{id}_X)
      }{\times}
    (Z,\phi_Z)
    \,,
  $$

  \vspace{-1mm}
\noindent
  we have the that product in $\mathbf{H}_{/X}$ is given by the
  fiber product over $X$ in $\mathbf{H}$:
  \vspace{-2mm}
  $$
    \mathrm{dom}
    \Big(
      (Y,\phi_Y) \times (Z,\phi_Z)
    \Big)
    \;\;\simeq\;\;
    Y \underset{X}{\times} Z
    \,.
  $$
\end{prop}
\begin{proof}
  Generally, limits in $\mathbf{H}_{/X}$
  are given by limits in $\mathbf{H}$ over
  the underlying co-cone diagram.
  Specifically: for $Y : \xymatrix@C=9pt{ \mathcal{I} \ar[r] & \mathbf{H} }$
  we have
  $
    \mathrm{dom}\big(
      \underset{\longleftarrow}{\mathrm{lim}}
      \,
      Y_\bullet
    \big)
    \simeq
    \underset{\longleftarrow}{\mathrm{lim}}
    \,
    \big(Y/X \big)_\bullet
    \,.
  $
  With this, the claim follows from the fact that the canonical inclusion of
  diagram categories
  \vspace{-4mm}
  $$
    \big\{
      \xymatrix@C=15pt{
        y \ar[r] & b \ar@{<-}[r] & z
      }
    \big\}
     \;    \xymatrix{
        \;   \ar@{^{(}->}[r]
      &
}
          \Bigg\{\!\!\!
      \raisebox{15pt}{
      \xymatrix@C=30pt@R=6pt{
        & t
        \ar[d]
        \ar[dl]
        \ar[dr]
        \\
        y \ar[r] & b \ar@{<-}[r] & z
      }
      }
   \!\!\Bigg\}
  $$

  \vspace{-1mm}
\noindent
  is an initial functor (i.e., under $(-)^{\mathrm{op}}$ it is
  a final functor).
\hfill \end{proof}

\begin{prop}[Terminal right base change of bare $\infty$-groupoids]
  \label{BaseChangeFromBareInfinityGroupoids}
  In the base $\infty$-topos $\mathbf{H} = \mathrm{Groupoids}_\infty$
  \eqref{InfinityCategoryOfInfinityGroupoids},
  the right base change along the terminal morphism (Example \ref{BaseChangeAlongTerminalMorphism})
  of an object $X \in \mathrm{Groupoids}_\infty$
  is given by the hom-$\infty$-groupoid out of $X$,
  regarded as the terminal object in the slice:
  \vspace{-2mm}
  $$
    (X \to \ast)_\ast
    \;\simeq\;
    \mathbf{H}_{/X}\big(X, - \big)
    \;\;:\;\;
    \xymatrix{
      \big(
        \mathrm{Groupoids}_\infty
      \big)_{/X}
      \ar[r]
      &
      \mathrm{Groupoids}_\infty\;.
    }
  $$
\end{prop}

\vspace{-9mm}
$\phantom{A}$
\begin{proof}
  We have the following chain of natural equivalences:
  \vspace{-1mm}
  \begin{equation}
    \label{TowardsUnderstandingRightBaseChangeFromInfinityGroupoid}
    \begin{aligned}
      &
      \mathrm{Groupoids}_\infty
      \big(
        A,
        (\mathrm{Groupoids}_\infty)_{/X}(X,B)
      \big)
      \;
      \simeq
      \;
      (\mathrm{Groupoids}_\infty)_{/X}
      \big(
        \Delta(A) \times_X X , B
      \big)
      \\
      & \simeq\;
      (\mathrm{Groupoids}_\infty)_{/X}
      \big(
        \Delta(A) , B
      \big)
      \;\simeq\;
      (\mathrm{Groupoids}_\infty)_{/X}
      \big(
        (X \to \ast)^\ast(A), B
      \big).
    \end{aligned}
  \end{equation}

  \vspace{-1mm}
  \noindent
  Here the first step observes that
  the slice $(\mathrm{Groupoids}_\infty)_{/X}$
  is itself an $\infty$-topos by Prop. \ref{SliceInfinityTopos},
  so that the tensoring equivalence of Prop. \ref{TensoringOfInfinityToposesOverInfinityGroupoids}
  applies. The second step uses the fact that $X$ is regarded as the terminal object
  in its own slice, so that forming Cartesian product with it
  is equivalently the identity operation. The last step observes
  that for the slice
  $\infty$-topos $\Delta \simeq (X \to \ast)^\ast$ \eqref{BaseGeometricMorphismForSliceTopos}
  by Example \ref{BaseChangeAlongTerminalMorphism}.
  In summary,
  the total equivalence of \eqref{TowardsUnderstandingRightBaseChangeFromInfinityGroupoid} is
  the hom-equivalence that characterizes $\mathbf{H}_{/X}(X,-)$ as
  a right adjoint to $(X \to \ast)^\ast$.
\hfill \end{proof}

\begin{prop}[Base change along effective epi is conservative {\cite[3.15]{NSS12}} ]
  \label{BaseChangeAlongEffectiveEpimorphismsIsConservative}
  Let $\mathbf{H}$ be an $\infty$-topos (Def. \ref{InfinityTopos}).
  For $\!\xymatrix@C=15pt{Y \ar@{->>}[r] & X}\!$ an effective epimorphism (Def. \ref{EffectiveEpimorphisms})
  in $\mathbf{H}$, the induced base change (Prop. \ref{BaseChange})
   \vspace{-2mm}
  $$
    \xymatrix{
      \mathbf{H}_{/X}
      \ar[r]^-{ i^\ast }
      &
      \mathbf{H}_{/Y}
    }
  $$

   \vspace{-2mm}
\noindent
  is a conservative $\infty$-functor, meaning that a morphism
  $f \in \mathbf{H}_{/X}$ is an equivalence if
  its base change $i^\ast(f)$ in $\mathbf{H}_{/Y}$ is an equivalence.
\end{prop}

\begin{prop}[Colimits of classifying maps are classifying maps of colimits]
  \label{ColimitsOfClassifyingMapsAreClassifyingMapsOfColimits}
  Let $\mathbf{H}$ be an $\infty$-topos (Def. \ref{InfinityTopos}),
  $\mathcal{I}$ a small $\infty$-category,
  $X_\bullet : \mathcal{I} \to \mathbf{H}$
  a diagram and
  $
    (\vdash E)_\bullet :
      X_\bullet
 \to
      \mathrm{const}_{\mathrm{Objects}_\kappa}
  $
  a transformation to the diagram constant on the
  object classifier \eqref{ObjectClassifier}, thus classifying a
  diagram $E_\bullet : \mathcal{I} \to \mathbf{H}$
  of bundles over $X_\bullet$.
  Then the colimit of $(\vdash E)_\bullet$ formed in the slice
  $\mathbf{H}_{/_{\mathrm{Objects}_\kappa}}$ (Prop. \ref{SliceInfinityTopos})
  is the colimit
  of $X_\bullet$ equipped with the classifying map for the colimit
  of $E_\bullet$:
  \vspace{-1mm}
  $$
    \underset{\longrightarrow}{\mathrm{lim}}\, (\vdash E)_\bullet
    \;\;\simeq\;\;
    \vdash\!\big( \underset{\longrightarrow}{\mathrm{lim}}\, E_\bullet \big).
  $$
\end{prop}

\vspace{-4mm}

\hspace{-.9cm}
\begin{tabular}{ll}
\begin{minipage}[left]{10cm}
\begin{proof}
Since underlying the colimit
$\underset{\longrightarrow}{\mathrm{lim}}(\vdash E)_\bullet$
in the slice $\infty$-topos $\mathbf{H}_{/\mathrm{Objects}_\kappa}$
is the colimit $\underset{\longrightarrow}{\mathrm{lim}}\, X_\bullet$
in $\mathbf{H}$ (by Example \ref{BaseChangeAlongTerminalMorphism})
we are dealing with a situation as shown in the diagram on the right
(where a simplicial diagram shape is shown just for definiteness of
illustration).
We need to demonstrate that the front square in this diagram is Cartesian.
Observe that
\begin{itemize}
\vspace{-3.5mm}
\item[{\bf (a) }] the vertical squares over
each $\vdash E_i$ are Cartesian by assumption,
whence
\vspace{-3.5mm}
\item[{\bf (b)}] also the solid vertical squares over each
$\xymatrix@C=10pt{X_i \ar[r] & X_j}$ are Cartesian, by
the pasting law (Prop. \ref{PastingLaw}).
\end{itemize}

\vspace{-3.5mm}
\noindent
This means that the assumption of Prop. \ref{ColimitsOfEquifiberedTransformations}
is satisfied for the left part of the diagram
(regarded as a transformation of diagrams from top to bottom)
implying that the dashed square is Cartesian.

This implies, together with {\bf (a)}, that the front square is Cartesian,
again the pasting law (Prop. \ref{PastingLaw}).
\hfill \end{proof}
\end{minipage}
&
\hspace{1cm}
\begin{minipage}[left]{8cm}
 \raisebox{50pt}{
  \xymatrix@R=8pt@C=11pt{
    && &
    \ar@<-10pt>@{..>}[ddl]
    \ar@<-5pt>@{<..}[ddl]
    \ar@{..>}[ddl]
    \ar@<+5pt>@{<..}[ddl]
    \ar@<+10pt>@{..>}[ddl]
    \\
    \\
    & &
    E_1
    \ar[dddd]|<<<<<<<<<<<<<{\phantom{AA\vert}}|>>>{\phantom{A\vert}}
    \ar[dddrrr]
    \ar@<-3pt>[ddl]
    \ar@{<-}[ddl]
    \ar@<+3pt>[ddl]
    \\
    \\
    &
    E_0
    \ar[dddd]|<<<<<<<<<<<<<<<<{\phantom{AA\vert}}
    \ar[drrrr]
    \ar@{-->}[ddl]
    \\
    &&
    &&&
    \widehat{\mathrm{Objects}}_\kappa
    \ar[dddd]
    \\
    \underset{\longrightarrow}{\mathrm{\lim}}\, E_\bullet
    \ar[dddd]
    \ar[urrrrr]
    &
    &
    X_1
    \ar[dddrrr]|-{\;
      \vdash E_1 \;\;\;\;
    }
    \ar@<-3pt>[ddl]
    \ar@{<-}[ddl]
    \ar@<+3pt>[ddl]
    \\
    \\
    &
    X_0
    \ar[drrrr]|-{\;
      \vdash E_0 \,
    }
    \ar@{-->}[ddl]
    \\
    &
    &&&& \mathrm{Objects}_\kappa
    \\
    \underset{\longrightarrow}{\mathrm{lim}}\;
    X_\bullet
    \ar[urrrrr]_{ \underset{\longrightarrow}{\mathrm{lim}}\,(\vdash E_\bullet) }
  }
  }
\end{minipage}
\end{tabular}

\medskip

\noindent {\bf $n$-Truncation.}
\begin{defn}[$n$-truncated objects {\cite[Def. 5.5.6.1]{Lurie09}}]
  \label{nTruncatedObjects}
  Let $n \in \{-2,-1,0,1,2,\cdots\}$.

  \vspace{-2mm}
  \item {\bf (i)}  An $\infty$-groupoid is called \emph{$n$-truncated}
  for $n \geq 0$
  if all its homotopy groups of degree $> n$ are trivial.
  It is called \emph{$(-1)$-truncated} if it is either
  empty or contractible, and \emph{(-2)-truncated} if it is
  (non-empty and) contractible.

    \vspace{-1mm}
\item {\bf (ii)}   Let $\mathcal{C}$ be an $\infty$-category.
  Then an object $X \in \mathcal{C}$ is \emph{$n$-truncated}
  if for all objects $U \in \mathcal{C}$ the
  hom-$\infty$-groupoid $\mathcal{C}(U,X)$ is $n$-truncated,
  in the above sense.
\end{defn}

\begin{defn}[$n$-truncated morphisms  {\cite[Def. 5.5.6.8]{Lurie09}}]
  \label{nTruncatedMorphism}
  Let $n \in  \{-2,-1,0,1,2,\cdots\}$.

    \vspace{-2mm}
  \item {\bf (i)}
  A morphism of $\infty$-groupoids is called \emph{$n$-truncated}
  if all its homotopy fibers are $n$-truncated $\infty$-groupoids
  according to Def. \ref{nTruncatedObjects}.

  \vspace{-2mm}
  \item {\bf (ii)}   Let $\mathcal{C}$ be an $\infty$-category. A morphism
  $X \overset{f}{\longrightarrow} Y$ in $\mathcal{C}$
  is called \emph{$n$-truncated} if for all objects
  $U \in \mathcal{C}$ the induced morphism of hom-$\infty$-groupoids
  $\xymatrix{
    \mathcal{C}(U,X)
    \ar[rr]|-{\;
      \mathcal{C}(U,f)
   \; }
    &&
    \mathcal{C}(U,Y)
    }
  $
  is $n$-truncated in the above sense.
\end{defn}

\begin{defn}[Monomorphisms]
  \label{Monomorphism}
  A (-1)-truncated morphism
  (Def. \ref{nTruncatedMorphism})
  is also called a
  \emph{monomorphism}, to be denoted
  \vspace{-2mm}
  \begin{equation}
    \label{NotationMonomorphism}
    \xymatrix{X \; \ar@{^{(}->}[r] & Y\;.}
  \end{equation}
\end{defn}

\begin{prop}[Monomorphisms are preserved by pushout {\cite[p. 21]{Rezk19}}]
  Let $\mathbf{H}$ be an $\infty$-topos (Def. \ref{InfinityTopos}). Then the class
  of monomorphisms in $\mathbf{H}$ (Def. \ref{Monomorphism})
  is closed under
      {\bf (i)} pullback
      and
      {\bf (ii)} composition.
\end{prop}

\begin{defn}[Poset of subobjects]
  \label{PosetOfSubobjects}
  Let $\mathbf{H}$ be an $\infty$-topos and
  $X \in \mathbf{H}$ any object. Then the \emph{poset of subobjects}
  of $X$ is the sub-$\infty$-category (Def. \ref{nTruncationModality})
  of $(-1)$-truncated objects of the slice over $X$:
  \vspace{-2mm}
  \begin{equation}
    \label{ThePosetOfSubobjects}
    \xymatrix{
      \mathrm{SubObjects}(X)
     \; \ar@{^{(}->}[rr]
      &&
      \mathbf{H}_{/_X}
    }
  \end{equation}

  \vspace{-2mm}
  \noindent
  whose objects are equivalently the monomorphisms (Def. \ref{Monomorphism})
  $U \hookrightarrow X$.
\end{defn}

\begin{prop}[$n$-Trucation modality {\cite[5.5.6.18]{Lurie09}}]
  \label{nTruncationModality}
  If $\mathbf{H}$ is an $\infty$-topos (Def. \ref{InfinityTopos}),
  for all $n \in \{-1,0,1,2, \cdots\}$,
  its full sub-$\infty$-category (Def. \ref{FullyFaithfulFunctor})
  of $n$-truncated objects (Def. \ref{nTruncatedObjects})
  is reflective, in that the inclusion functor has a left adjoint
  (Def. \ref{AdjointInfinityFunctors}):
  \vspace{-2mm}
  \begin{equation}
    \label{nTruncationAdjunction}
    \xymatrix@R=2pt{
      \mathbf{H}
   \;\;   \ar@{->}@<+6pt>[rr]^-{ \tau_n }_-{ \bot }
      \ar@{<-^{)}}@<-6pt>[rr]_-{ i_n }
      &&
    \;\;  \mathbf{H}_{n}
      \\
      \mathclap{
        \mbox{
          \tiny
          \rm
          \color{darkblue} \bf
          $\infty$-topos
        }
      }
      &&
      \mathclap{
        \mbox{
          \tiny
          \rm
          \color{darkblue} \bf
          \begin{tabular}{c}
          sub-$\infty$-category
          \\
          of $n$-truncated objects
          \end{tabular}
        }
      }
    }
  \end{equation}

  \vspace{-2mm}
  \noindent  We write
  for the induced \emph{$n$-truncation modality} \eqref{AdjointModalities}:
  \begin{equation}
    \label{nTruncatedModality}
    \big(
      \underset{
        \mbox{
          \tiny
          \color{darkblue} \bf
          \begin{tabular}{c}
            ``$n$-truncated''
          \end{tabular}
        }
      }{
      {\pmb\tau}_n
      \;:=\;
      i_n \circ \tau_n
      }
    \big)
    \;:\;
    \mathbf{H}
      \longrightarrow
    \mathbf{H}\;.
  \end{equation}

\end{prop}
\begin{defn}[Effective epimorphisms {\cite[Cor. 6.2.3.5]{Lurie09}}]
  \label{EffectiveEpimorphisms}
  Let $\mathbf{H}$ be an $\infty$-topos. A morphism in
  $\mathbf{H}$ is called an \emph{effective epimorphism},
  to be denoted
  \vspace{-2mm}
  \begin{equation}
    \label{NotationEffectiveEpi}
    \xymatrix{ Y \ar@{->>}[r]^-{f} & Z }
  \end{equation}

  \vspace{-1mm}
\noindent
  if,
  when regarded as an object of the slice over $X$ (Prop. \ref{SliceInfinityTopos}),
  its $(-1)$-truncation (Prop. \ref{nTruncationModality}) is the terminal object
  \vspace{-2mm}
    $$
    \tau_{(-1)} (f) \;\simeq\; \ast
    \;\;\;
    \in \mathbf{H}_{/X}
    \,.
  $$
  We write

  \vspace{-.5cm}

  \begin{equation}
    \label{CategoryOfEffectiveEpimorphisms}
    \mathrm{EffectiveEpimorphisms}(\mathbf{H})
    \;\subset\;
    \mathbf{H}^{(0 \to 1)}
    \;\in\;
    \mathrm{Categories}_\infty
  \end{equation}
  for the full sub-$\infty$-category (Def. \ref{FullyFaithfulFunctor})
  of the arrow-category
  of $\mathbf{H}$ on those that are effective epimorphisms.
\end{defn}

\begin{defn}[$n$-Connected morphisms {\cite[Prop. 6.5.1.12]{Lurie09}}]
  \label{nConnectedMorphisms}
  Let $\mathbf{H}$ be an $\infty$-topos (Def. \ref{InfinityTopos})
  and $n \in \{-1, 0, 1, 2 \cdots\}$.
  Then a morphism $\xymatrix{X \ar[r]^f & Y}$ in $\mathbf{H}$
  is called \emph{$n$-connected} if, regarded as an object
  in the slice over $X$ (Prop. \ref{SliceInfinityTopos}),
  its $n$-truncation (Def. \ref{nTruncationModality}) is the terminal object:
  \vspace{-2mm}
  $$
    \xymatrix{
      Y \ar[r]^-f & X
      \;\;\;
      \mbox{
        is $n$-truncated
      }
    }
    \phantom{AA}
    \Leftrightarrow
    \phantom{AA}
    \tau_n(f) \;\simeq\; \ast
    \;
    \in \mathbf{H}_{/X}
    \,.
  $$

  \vspace{-2mm}
\noindent
  Hence the $(-1)$-connected morphisms are equivalently
  the effective epimorphisms (Def. \ref{EffectiveEpimorphisms}).
\end{defn}

\begin{lemma}[Effective epimorphisms are preserved by pullback {\cite[6.2.3.15]{Lurie09}}]
  \label{EffectiveEpimorphismsArePreservedByPullback}
  Let $\mathbf{H}$ be an $\infty$-topos (Def. \ref{InfinityTopos}). Then the class
  of effective epimorphisms in $\mathbf{H}$ (Def. \ref{EffectiveEpimorphisms})
  is closed under
  {\bf (i)} pullback
  and
  {\bf (ii)} composition.
\end{lemma}

\newpage

\noindent {\bf $n$-Image factorization.}
\begin{prop}[Connected/truncated factorization system {\cite[Ex. 5.2.8.16]{Lurie09}\cite[Prop. 5.8]{Rezk10}}]
  \label{nConnectednTruncatedFactorizationSystem}
  Let $\mathbf{H}$ be an $\infty$-topos.
  Then, for all $n \in \{-1,0,1,2, \cdots\}$,
  the pair of classes
  of $n$-connected/$n$-truncated morphisms
  (Def. \ref{nConnectedMorphisms},
  Def. \ref{nTruncatedMorphism})
  forms an orthogonal factorization system:

  \noindent {\bf (i)}  every morphism $f$ in $\mathbf{H}$ factors essentially uniquely as
  \vspace{-3mm}
  \begin{equation}
    \label{nImageFactorization}
   \hspace{1cm}
    \xymatrix@R=.5em@C=3em{
      X
      \ar[dr]_-{
        \mbox{
          \tiny
          \rm
          $n$-connected \phantom{A}
        }
      }
      \ar[rr]^-{ f }
      &&
      Y
      \\
      &
      \mathrm{im}_n(f)
      \ar[ur]_-{
        \mbox{
          \tiny
          \rm
          $n$-truncated
        }
      }
    }
  \end{equation}

  \vspace{-2mm}
  \noindent {\bf (ii)}
  every commuting square as follows has an essentially unique dashed lift:
 \vspace{-2mm}
  \begin{equation}
    \label{ConnectedTruncatedLiftingProperty}
  \hspace{1cm}
   \raisebox{20pt}{
    \xymatrix@R=.4em@C=2em{
      X
      \ar[rr]
      \ar[dd]_{
        \mbox{
          \tiny
          \rm
          $n$-connected
        }
      }
      & &
      A
      \ar[dd]^{
        \mbox{
          \tiny
          \rm
          $n$-truncated
        }
      }
      \\
      \\
      Y
      \ar@{-->}[uurr]
      \ar[rr]
      &&
      B
    }
    }
  \end{equation}
\end{prop}

\begin{example}[Epi/mono factorization]
  \label{EffectiveEpimorphismMonomorphismFactorization}
  For $n = -1$, the connected/truncated factorization system
  (Prop. \ref{nConnectednTruncatedFactorizationSystem})
  has as left class the effective epimorphisms (Def. \ref{EffectiveEpimorphisms})
  and as right class the monomorphisms (Def. \ref{Monomorphism}).
  Hence, with the notation from \eqref{NotationEffectiveEpi}
  and \eqref{NotationMonomorphism}:

 \noindent  {\bf (i)} the (-1)-image factorization \eqref{nImageFactorization}
  reads:
  \vspace{-1cm}
  \begin{equation}
    \label{EpiMonoFactorization}
   \hspace{1cm}
    \xymatrix@R=.3em@C=3em{
      X
      \ar@{->>}[dr]_-{
      }
      \ar[rr]^-{ f }
      &&
      Y
      \\
      &
      \mathrm{im}_{-1}(f)
      \ar@{^{(}->}[ur]_-{
      }
    }
  \end{equation}

  \vspace{.3cm}

  \noindent {\bf (ii)} the lifting property \eqref{ConnectedTruncatedLiftingProperty}
  for $n = -1$ reads:

  \vspace{-1.3cm}

  \begin{equation}
   \label{EpiMonoLifting}
   \raisebox{40pt}{
    \xymatrix@R=.4em@C=2em{
      X
      \ar[rr]
      \ar@{->>}[dd]_{
        \mbox{
          \tiny
          \rm
        }
      }
      & &
      A\mathclap{\phantom{\vert}}
      \ar@{^{(}->}[dd]^{
        \mbox{
          \tiny
          \rm
        }
      }
      \\
      \\
      Y
      \ar@{-->}[uurr]
      \ar[rr]
      &&
      B
    }
    }
  \end{equation}
\end{example}

\medskip

\noindent {\bf Groupoids and Stacks.}
\begin{defn}[Groupoids internal to an $\infty$-topos {\cite[6.1.2.7]{Lurie09}}]
  \label{InHigherToposGroupoids}
  Let $\mathbf{H}$ be an $\infty$-topos (Def. \ref{InfinityTopos}).

  \noindent {\bf (i)} A \emph{groupoid in}  $\mathbf{H}$ is a simplicial diagram
  \begin{equation}
    \label{SimplicialDiagramUnderlyingGroupoidObject}
    X_\bullet
    :
    \xymatrix@C=14pt{
      \Delta^{\mathrm{op}}
      \ar[r]
      &
      \mathbf{H}
    }
  \end{equation}
  which satisfies the \emph{groupoidal Segal condition}:
  For all $n \in \mathbb{N}$ and for all
  partitions of the set of $n+1$ elements by two subsets
  that share a unique element, the corresponding image
  under $X_\bullet$ is a Cartesian square (Notation \ref{CartesianSquares}):
  \vspace{-2mm}
  \begin{equation}
    \label{GroupoidSegalCondtion}
    \raisebox{24pt}{
    \xymatrix@C=20pt@R=4pt{
      &
      \{0,1,\cdots, n\}
      \\
      S_1
      \ar@{^{(}->}[ur]
      &&
      S_2
      \ar@{_{(}->}[ul]
      \\
      &
      \ast
      \ar[ul]
      \ar[ur]
      \ar@{}[uu]|-{\mbox{\tiny(po)}}
    }
    }
    \;\;\;\;\;\;\;\;
    \overset{
      X_\bullet
    }{\longmapsto}
    \;\;\;\;\;\;\;\;
    \raisebox{24pt}{
    \xymatrix@C=25pt@R=4pt{
      &
      \;\;
        X_n
      \;\;
      \ar[dl]
      \ar[dr]
      \ar@{}[dd]|-{\mbox{\tiny(pb)}}
      \\
      X_{\left\vert S_1\right\vert-1}
      \ar[dr]
        &&
      X_{\left\vert S_2 \right\vert-1}
      \ar[dl]
      \\
      &
      X_0
    }
    }
  \end{equation}

\vspace{-1mm}
  \noindent {\bf (ii)} We write

  \vspace{-.9cm}

  \begin{equation}
    \label{InInfinityToposGroupoids}
    \xymatrix{
      \mathrm{Groupoids}(\mathbf{H})
      \;
      \ar@{^{(}->}[r]
      &
      \mathbf{H}^{(\Delta^{\mathrm{op}})}
    }
    \;\in\;
    \mathrm{Categories}_\infty
  \end{equation}
  for the full sub-$\infty$-category of that of simplicial
  diagrams in $\mathbf{H}$ on those that are groupoids.
\end{defn}
\begin{example}[Nerves]
  \label{Nerve}
  Let $\mathbf{H}$ be an $\infty$-topos (Def. \ref{InfinityTopos})
  and $\xymatrix{X \ar[r]^-{f} & \mathcal{X}}$ a morphism in $\mathbf{H}$.
  Its \emph{nerve} is the simplicial diagram
  \vspace{-2mm}
  of its iterated homotopy fiber products:
  \begin{equation}
    \label{NerveOfAnAtlas}
    \xymatrix@R=-10pt{
      \mathllap{
        \mathrm{Nerve}_\bullet(f)
        \;:\;
        \;
      }
      \Delta^{\mathrm{op}}
      \ar[rr]
      &&
      \mathbf{H}
      \\
      [n]
      \ar@{}[rr]|-{\longmapsto}
      &&
      \underset{
        \mbox{
          \tiny
          $n$ factors
        }
      }{
      \underbrace{
        X \times_{{}_{\mathcal{X}}} X \times_{{}_{\mathcal{X}}}
        \cdots
        \times_{{}_{\mathcal{X}}} X
      }
      }
    }
  \end{equation}

   \vspace{-2mm}
\noindent
  with face maps the projections and degeneracy maps the diagonals.
  This is evidently a groupoid object according to Def. \ref{InHigherToposGroupoids}:
  \vspace{-3mm}
  $$
    \mathrm{Nerve}_\bullet(f)
    \;\in\;
    \mathrm{Groupoids}(\mathbf{H})
    \,.
  $$
\end{example}

\vspace{1mm}
\begin{prop}[Groupoids equivalent to stacks with atlases {\cite[6.2.3.5]{Lurie09}}]
  \label{GroupoidsEquivalentToEffectiveEpimorphisms}
$\phantom{A}$

  \vspace{0cm}

  \hspace{-.9cm}
  \begin{tabular}{ll}
  \begin{minipage}[l]{9.8cm}
  Let $\mathbf{H}$ be an $\infty$-topos (Def. \ref{InfinityTopos}).
  Then the $\infty$-functor sending
  $X_\bullet \in \mathrm{Groupoids}(\mathbf{H})$
  (Def. \ref{InHigherToposGroupoids}) to the $X_0$-component
  of its colimiting cocone
  \begin{itemize}
  \vspace{-3mm}
\item[{\bf (i)}] lands in effective epimorphisms \eqref{CategoryOfEffectiveEpimorphisms}
  and
  \vspace{-3mm}
  \item[{\bf (ii)}]
  constitutes an equivalence of $\infty$-categories
  whose inverse is given by the construction of nerves
  (Example \ref{Nerve}):
  \end{itemize}
  \vspace{-3mm}
  \begin{equation}
    \label{EquivalenceBetweenGroupoidsAndEffectiveEpimorphisms}
    \xymatrix@R=-2pt{
      \mathrm{Groupoids}(\mathbf{H})
      \ar[r]^-{ \simeq }
      &
      \mathrm{EffectiveEpimorphisms}(\mathbf{H})
      \\
      X_\bullet
      \ar@{}[r]|-{\longmapsto}
      &
      \big(
        X_0
          \twoheadrightarrow
        \underset{\longrightarrow}{\mathrm{lim}}\, X_\bullet
      \big)
      \\
      \mathrm{Nerve}_\bullet(a)
      \ar@{}[r]|<<<<<<<<{\longmapsfrom}
      &
      \big(
        X
          \overset{a}{\twoheadrightarrow}
        \;\mathcal{X}\;
      \big)
    }
  \end{equation}
  \end{minipage}
  &
  \begin{minipage}[l]{7cm}
  \begin{equation}
  \label{StacksAtlasesAndGroupoids}
  \raisebox{40pt}{
  \xymatrix@C=10pt{
    \ar@<-12pt>@{..>}[d]
    \ar@<-6pt>@{<..}[d]
    \ar@{..>}[d]
    \ar@<+6pt>@{<..}[d]
    \ar@<+12pt>@{..>}[d]
    &
    \ar@<-12pt>@{..>}[d]
    \ar@<-6pt>@{<..}[d]
    \ar@{..>}[d]
    \ar@<+6pt>@{<..}[d]
    \ar@<+12pt>@{..>}[d]
    \\
    X \times_{\mathcal{X}} X
    \ar@{}[r]|-{\simeq}
    \ar@<-6pt>[d]_-{\mathrm{pr}_1}
    \ar@{<-}[d]|-{\Delta}
    \ar@<+6pt>[d]^{\mathrm{pr}_2}
    &
    X_1
    \ar@<-6pt>[d]_-{s}
    \ar@{<-}[d]|-{e}
    \ar@<+6pt>[d]^-{t}
    &
    \mbox{
      \footnotesize
      \color{darkblue}
      \bf
      ``groupoid''
    }
    \\
    X_0
    \ar@{->>}[d]_-{a}
    \ar@{=}[r]
    &
    X_0
    \ar@{->>}[d]
    &
    \ar@{}[d]|-{
      \mbox{
        \footnotesize
        \color{darkblue}
        \bf
        ``atlas''
      }
    }
    \\
    \mathcal{X}
    \ar@{}[r]|-{\simeq}
    &
    \underset{\longrightarrow}{\mathrm{lim}}\, X_\bullet
    &
    \mbox{
      \footnotesize
      \color{darkblue}
      \bf
      ``stack''
    }
  }
  }
  \end{equation}
  \end{minipage}
 \end{tabular}
\end{prop}
\begin{remark}[Internal groupoids with prescribed properties]
  Often one considers
  $X_\bullet \in \mathrm{Groupoids}(\mathbf{H})$ (Def. \ref{InHigherToposGroupoids})
  whose simplicial component diagram \eqref{SimplicialDiagramUnderlyingGroupoidObject} is inside
  a chosen  sub-$\infty$-category of $\mathbf{H}$.
  Key examples are {\'e}tale groupoids (Def. \ref{EtaleGroupoids} below)
  and $V$-{\'e}tale groupoids (Remark \ref{VFoldsAndVEtaleGroupoids} below).
\end{remark}

\begin{remark}[Morita morphisms of groupoids]
  \label{MoritaMorphismsOfGroupoids}
  A morphism between
  stacks $\mathcal{X} := \underset{\longrightarrow}{\mathrm{lim}}\,X_\bullet$
  underlying groupoids $X_\bullet$
  (according to Prop. \ref{GroupoidsEquivalentToEffectiveEpimorphisms})
  \emph{without} (i.e., disregarding) the corresponding atlas
  is also known as a \emph{Morita morphism}
  (in particular, a \emph{Morita equivalence} if it is an
  equivalence),
  or a \emph{Hilsum-Skandalis morphism}
  \cite{HilsumSkandalis87}\cite{Pradines89},
  or a
  \emph{groupoid bibundle} \cite{Blohmann07}\cite[Prop. 2.2.34]{Nuiten13}
  between the corresponding groupoids:
  \vspace{-1mm}
  $$
    \xymatrix@R=2pt@C=2.5em{
      \mathrm{Groupoids}(\mathbf{H})
      \ar[r]^-{ \simeq }
      &
      \mathrm{EffectiveEpimorphisms}(\mathbf{H})
      \ar[r]^-{ \mathrm{codom} }
      &
      \mathbf{H}
      \\
      \mathllap{
        \mbox{
          \tiny
          \color{darkblue}
          \bf
          groupoid
        }
        \;\;\;
      }
      X_\bullet \ar@{}[r]|-{ \longmapsto }
      &
      (X_0 \twoheadrightarrow \mathcal{X})
      \ar@{}[r]|-{\longmapsto}
      &
      \mathcal{X}
      \mathrlap{
        \;\;
        \mbox{
          \tiny
          \color{darkblue}
          \bf
          ``stack''
        }
      }
      \ar[dd]_-{
        f}
        ^>>>{
          \!\!\!\!
          {
            \tiny
            \color{darkblue}
            \bf
            \begin{tabular}{l}
              morphism of underlying stacks =
              \\
              ``Morita morphism'' of groupoids
            \end{tabular}
          }
        }
      \\
      \\
      \mathllap{
        \mbox{
          \tiny
          \color{darkblue}
          \bf
          groupoid
        }
        \;\;\;
      }
      Y_\bullet \ar@{}[r]|-{ \longmapsto }
      &
      (Y_0 \twoheadrightarrow \mathcal{Y})
      \ar@{}[r]|-{\longmapsto}
      &
      \mathcal{Y}
      \mathrlap{
        \;\;
        \mbox{
          \tiny
          \color{darkblue}
          \bf
          ``stack''
        }
      }
    }
  $$

  \vspace{0mm}
\noindent
  Hence whether or not there is a conceptual distinction between
  ``geometric groupoids'' and ``stacks'' depends on whether
  morphisms of groupoids are taken to be their plain morphisms or their
  Morita morphisms. In practice, one is typically interested
  in the latter case. Indeed, the groupoid atlas of a stack,
  whose preservation restricts Morita morphisms to plain morphisms
  of groupoids,
  by Prop. \ref{GroupoidsEquivalentToEffectiveEpimorphisms},
  is, in practice, typically required to exist with a certain property,
  but not required to be preserved by morphisms
  (this is so notably for $V$-{\'e}tale groupoids, Remark \ref{VFoldsAndVEtaleGroupoids} below).
  In particular, the
  $\mathrm{SmoothGroupoids}_\infty$ of Example \ref{SmoothInfinityGroupoids}
  and the $\mathrm{JetsOfSmoothGroupoids}_\infty$ of Example \ref{FormalSmoothInfinityGroupoids} below
  are $\infty$-groupoids with Morita morphisms understood,
  hence could also be called (\emph{jets of}) \emph{smooth $\infty$-stacks}.
\end{remark}

\begin{prop}[Equifibered morphisms of groupoids]
  \label{GroupoidMorphismsCartesian}
  Let $\mathbf{H}$ be an $\infty$-topos (Def. \ref{InfinityTopos})
  and  $X_\bullet, Y_\bullet \;\in\; \mathrm{Groupoids}(\mathbf{H})$
  (Def. \ref{InHigherToposGroupoids}).
  Then, under the equivalence \eqref{EquivalenceBetweenGroupoidsAndEffectiveEpimorphisms}
  between groupoids and their stacks with atlases (Prop. \ref{GroupoidsEquivalentToEffectiveEpimorphisms}),
  we have that equifibered morphisms of groupoids
  correspond to Cartesian squares between their atlases:
   \vspace{-3mm}
  $$
    \xymatrix@R=10pt{
      X_\bullet
      \ar@{=>}[r]^-{f_\bullet}
      &
      Y_\bullet
    }
    \;\;
    \mbox{such that}
    \;\;\;\;
    \underset{
      [n_1] \overset{\phi}{\to} [n_2]
    }{\forall}
    \;\;\;\;
    \raisebox{15pt}{
    \xymatrix@R=15pt@C=50pt{
      X_{n_1}
      \ar[r]^-{ f_{n_1} }
      \ar[d]_-{X_\phi}
      \ar@{}[dr]|-{\mbox{\tiny\rm(pb)}}
      &
      Y_{n_1}
      \ar[d]^-{Y_\phi}
      \\
      X_{n_2}
      \ar[r]_-{ f_{n_2} }
      &
      Y_{n_2}
    }
    }
    \phantom{AAAA}
    \Leftrightarrow
    \phantom{AAAA}
    \raisebox{15pt}{
    \xymatrix@R=15pt@C=25pt{
      X_0
      \ar@{->>}[d]_-{ a_X }
      \ar[rr]^-{ f_0 }
      \ar@{}[drr]|-{ \mbox{\tiny\rm(pb)} }
      &&
      Y_0
      \ar@{->>}[d]^-{ a_Y }
      \\
      \mathcal{X}
      \ar[rr]_{
        \underset{
          \longrightarrow
        }{\mathrm{lim}}\, f_\bullet
      }
      &&
      \mathcal{Y}
    }
    }
  $$
\end{prop}

\vspace{-4mm}
\begin{proof}
  From right to left this follows by the pasting law
  (Prop. \ref{PastingLaw}), while from left to right this
  is Prop. \ref{ColimitsOfEquifiberedTransformations}.
\hfill \end{proof}

\medskip

\subsection{Galois theory}
\label{GaloisTheory}

We discuss here the \emph{internal} formulation in $\infty$-toposes
of the theory of \emph{groups}, \emph{group actions}, and
\emph{fiber bundles}, following
\cite{NSS12}\cite{SSS09} (see \cite{FSS13a} for exposition).
Externally, these concepts are known as
\emph{grouplike $A_\infty$-algebras} or equivalently:
\emph{grouplike $E_1$-algebras}
(here: in $\infty$-stacks) and as their \emph{$A_\infty$-modules} etc.,
and are traditionally presented by simplicial techniques
\cite{May72}\cite{Lurie17}.
But internally the theory becomes finitary and elementary,
with all concepts emerging naturally from pastings
of a few Cartesian squares. Accordingly,
much of the following
constructions may readily be expressed fully formally
in homotopy type theory \cite{BvDR18} (see p. \pageref{HoTTFormalizations}).
Thus, the following elegant characterizations of

\begin{itemize}
 \vspace{-2.5mm}
\item[$\circ$] groups (Prop. \ref{LoopingAndDelooping}),
 \vspace{-2.5mm}
\item[$\circ$] group actions (Prop. \ref{InfinityAction}),
 \vspace{-2.5mm}
\item[$\circ$] principal bundles (Prop. \ref{ClassificationPrincipalInfinityBundles}),
 \vspace{-2.5mm}
\item[$\circ$] fiber bundles (Prop. \ref{FiberBundlesClassified}),
\end{itemize}

 \vspace{-2.5mm}
\noindent in an $\infty$-topos $\mathbf{H}$ may be taken to be
the \emph{definition} of these notions for all purposes of
internal constructions.

\medskip

\noindent {\bf Groups.} The following characterization of group $\infty$-stacks
(Prop. \ref{LoopingAndDelooping}) is the time-honored
\emph{May recognition theorem} \cite{May72} generalized from
$\mathrm{Groupoids}_\infty$ to general $\infty$-toposes
\cite[7.2.2.11]{Lurie09}\cite[6.2.6.15]{Lurie17}:
\begin{prop}[Groups {\cite[Thm. 2.19]{NSS12}}]
  \label{LoopingAndDelooping}
  Let $\mathbf{H}$ be an $\infty$-topos (Def. \ref{InfinityTopos}).
    Then the operation of sending an
  $\infty$-group $G$ to the homotopy quotient
  of its action on a point constitutes an
  equivalence of $\infty$-categories:
  \vspace{-2mm}
  \begin{equation}
    \label{LoopingDeloopingEquivalence}
    \xymatrix@R=2pt{
      \mathrm{Groups}
      \big(
        \mathbf{H}
      \big)
      \ar@{<-}@<+6pt>[rr]^-{ \Omega }
      \ar@{->}@<-6pt>[rr]_-{ \mathbf{B} }^-{
        \raisebox{1pt}{$\simeq$}
      }
      &&
      \mathbf{H}^{\ast/}_{\geq 1}
      \\
      G
        \ar@{|->}[rr]
      &&
      \ast \!\sslash\! G
    }
  \end{equation}

   \vspace{-2mm}
\noindent  between the $\infty$-category of $\infty$-group objects
  and the $\infty$-category of pointed and
  connected objects in $\mathbf{H}$.
 The inverse equivalence is given by forming the loop space object
  \vspace{-1mm}
 \begin{equation}
   \label{GroupIsLooping}
   \raisebox{10pt}{
   \xymatrix@C=3em@R=1em{
     \mathllap{
       G \simeq \;
     }
     \Omega \mathbf{B}G
       \ar[d]
       \ar[r]
       \ar@{}[dr]|-{\mbox{\tiny\rm(pb)}}
     &
     \ast
     \ar[d]
     \\
     \ast
     \ar[r]
     &
     \mathbf{B}G
   }
   }
 \end{equation}
\end{prop}
\begin{example}[Point in delooping is an effective epi]
  \label{PointInDeloopingIsEffectiveEpi}
  For $G \in \mathrm{Groups}(\mathbf{H})$,
  the essentially unique morphism that exhibits
  its delooping as a pointed object (Prop. \ref{LoopingAndDelooping})
   \vspace{-2mm}
  \begin{equation}
    \label{PointInclusionIntoDelooping}
    \xymatrix{
      \ast
      \ar@{->>}[r]
      &
      \mathbf{B}G
    },
  \end{equation}

 \vspace{-2mm}
\noindent
  is an effective epimorphism (Def. \ref{nTruncationModality}).
  Thus, Prop. \ref{GroupoidsEquivalentToEffectiveEpimorphisms}
  says here that

\noindent   {\bf (i)} groups in $\mathbf{H}$ are, equivalently, the
  groupoids in $\mathbf{H}$ (Def. \ref{InHigherToposGroupoids}) that
  admit an atlas by the point
  and,

  \noindent {\bf (ii)} with \eqref{GroupIsLooping}, we have
  \begin{equation}
    \label{DeloopingAsHomotopyColimit}
    \mathbf{B}G
    \;\simeq\;
    \underset{\longrightarrow}{\mathrm{lim}}
    G^{\times_\bullet}
    \;\;
    \in
    \;
    \mathbf{H}\;.
  \end{equation}
\end{example}

\begin{example}[Neutral element]
  \label{NeutralElement}
  Let $\mathbf{H}$ be an $\infty$-topos.
  Given a group $G \in$

 \vspace{-3mm}
  \hspace{-.9cm}
  \begin{tabular}{ll}
  \begin{minipage}[l]{12cm}
    $\mathrm{Groups}(\mathbf{H})$
  in the form of a pointed connected object
  $\ast \to \mathbf{B}G$, according to
  Prop. \ref{LoopingAndDelooping}, its
  \emph{neutral element}
  $
\ast \overset{e}{\longrightarrow}
      G
      $
  is the diagonal morphism into the
  defining homotopy fiber product \eqref{GroupIsLooping},
  hence the canonical morphism induced by the universal property
  of the homotopy fiber product from the equivalence with itself
  of the point inclusion into $\mathbf{B}G$ \eqref{PointInclusionIntoDelooping}.
  \end{minipage}
  &
   $$
   \raisebox{40pt}{
    \xymatrix@C=4em@R=4pt{
      &
      \ast
      \ar@{=}@/^1pc/[dddr]
      \ar@{=}@/_1pc/[dddl]
      \ar@{..>}[dd]^-e
      \\
      \\
      &
      G
      \ar[dl]
      \ar[dr]
      \ar@{}[dd]|-{\mbox{\tiny\rm(pb)}}
      \\
      \ast
      \ar[dr]
      &&
      \ast
      \ar[dl]
      \\
      &
      \mathbf{B}G
    }
    }
  $$
  \end{tabular}
\end{example}
\begin{example}[Group division/shear map]
  \label{GroupOperation}
  Let $\mathbf{H}$ be an $\infty$-topos.
  Given a group $G \in \mathrm{Groups}(\mathbf{H})$
  in the form of a pointed connected object
  $\ast \longrightarrow \mathbf{B}G$, according to
  Prop. \ref{LoopingAndDelooping},
  the group division operation
  \vspace{-2mm}
  $$
    \xymatrix@C=2.3em{
      G \times G
       \ar[rr]^-{\, (-) \cdot (-)^{-1} } && G
    }
  $$

   \vspace{-2mm}
\noindent
  is exhibited by the universal morphism shown dashed in the following
  diagram:
   \vspace{-2mm}
  \begin{equation}
    \label{GroupOperationFromLooping}
    \raisebox{30pt}{
    \xymatrix@R=16pt{
      \ar@<-6pt>@{..}[d]
      \ar@{..}[d]
      \ar@<+6pt>@{..}[d]
      &&
      \ar@<-6pt>@{..}[d]
      \ar@{..}[d]
      \ar@<+6pt>@{..}[d]
      \\
      G \times G
      \ar@<-3pt>[d]
      \ar@<+3pt>[d]
      \ar@{-->}[rr]^-{
        (-) \cdot (-)^{-1}
      }
      &&
      G
      \ar@<-3pt>[d]
      \ar@<+3pt>[d]
      \\
      G
      \ar[d]
      \ar[rr]
      &&
      \ast
      \ar[d]
      \\
      \ast
      \ar[rr]
      &&
      \mathbf{B}G
    }
    }
    {\phantom{AAAAAAAAAA}}
    \raisebox{20pt}{
    \xymatrix@C=12pt@R=5pt{
      &
      G \times G
      \ar[ddl]
      \ar[ddr]
      \ar@{-->}[drrr]^{ (-) \cdot (-)^{-1} }
      \\
      && & &
      G
      \ar[ddl]
      \ar[ddr]_<<<<{\ }="s3"
      \\
      G
      \ar[ddr]
      \ar[drrr]_>>>>>>>>{\ }="s1"
      &&
      G
      \ar[ddl]|<<<{\phantom{AA}}
      \ar[drrr]|>>>>>>>>>>>>>>>>{\phantom{AA}}|>>>>>>>>>{\phantom{AA}}_>>>>>{\ }="s2"
      \\
      &&&
      \,\ast\,
      \ar[ddr]^>>>>{\ }="t3"
      &&
      \,\ast\,
      \ar[ddl]
      \\
      &
      \ast
      \ar[drrr]^<<<<<<<<<<<<<{\ }="t1"^<<<<<<<<<<<<<{\ }="t2"
      \\
      &
      && &
      \mathbf{B}G
      \ar@{=>} "s1"; "t1"
      \ar@{=>}|<<<<{\phantom{AAA}}|>>>>>>>>>>{\phantom{AA}}  "s2"; "t2"
      \ar@{=>} "s3"; "t3"
    }
    }
  \end{equation}

   \vspace{-2mm}
\noindent
  On the left, we are showing this as part of a morphism of
  {\v C}ech nerve augmented simplicial diagrams.
    On the right, the situation is shown in more detail:
  Here the right and the two bottom squares
  are all the looping relation \eqref{GroupIsLooping},
  while the left square exhibits the plain product of $G$ with itself.
  With this, the universal property of the right square implies the
  essentially unique dashed morphism making the total diagram homotopy-commute.
  Notice:

 \vspace{-1mm}
\item {\bf (i)} The two top squares are also Cartesian:
  This follows from the pasting law (Prop. \ref{PastingLaw})
  using, for the top front square, that the left and right and the
  bottom rear squares are Cartesian; and similarly for the
  top rear square.

 \vspace{-1mm}
\item {\bf (ii)} The total homotopy
  filling the top and the right faces in \eqref{GroupOperationFromLooping}
  is, by commutativity, equivalent to the total homotopy filling the
  left and the bottom faces. But, in performing the composition this way,
  the direction of one of the two bottom homotopies gets reversed.
  This is why this construction gives the division map $(-)\cdot (-)^{-1}$
  (shear map) instead of the plain group product.
\end{example}

\begin{prop}[Mayer-Vietoris sequence {\cite[Prop. 3.6.142]{dcct}}]
  \label{MayerVietorisSequence}
  Let $\mathbf{H}$ be an $\infty$-topos (Def. \ref{InfinityTopos}),
  $G \in \mathrm{Groups}(\mathbf{H})$ (Prop. \ref{LoopingAndDelooping})
  and $(X,f), (Y,g) \in \mathbf{H}_{/_G}$ two
  objects in the slice (Prop. \ref{SliceInfinityTopos}) over the underlying object of $G$.
    Then their homotopy fiber product
    \vspace{-3mm}
  $$
    \xymatrix@C=5em@R=1.5em{
      X \underset{G}{\times} Y
      \ar[d]_-{ \mathrm{pr}_X }
      \ar[r]^-{ \mathrm{pr}_Y }
      \ar@{}[dr]|<<<<<<<{\mbox{\tiny\rm(pb)}}
      &
      Y
      \ar[d]^-g
      \\
      X \ar[r]_-f
      &
      G
    }
  $$

     \vspace{-2mm}
\noindent
  is equivalently exhibited by the
  following \emph{Mayer-Vietoris homotopy fiber sequence}
  \vspace{-2mm}
  \begin{equation}
    \label{MayerVietorisSquare}
    \xymatrix@C=5em@R=1.5em{
      X \times_G Y
      \ar[d]_-{ (\mathrm{pr}_X, \mathrm{pr}_Y) }
      \ar@{}[drr]|-{ \mbox{\tiny\rm(pb)} }
      \ar[rr]
      &&
      \ast
      \ar[d]
      \\
      X \times Y
      \ar[r]^-{ (f,g) }
      \ar@/_1.2pc/[rr]_-{ f \cdot g^{-1} }
      &
      G \times G
      \ar[r]^-{ (-)\cdot (-)^{-1} }
      &
      G
      \,,
    }
  \end{equation}

     \vspace{-2mm}
\noindent
  where the morphism on the bottom right is the
  group division map \eqref{GroupOperationFromLooping}.
\end{prop}

\medskip

\noindent {\bf Group actions.}
\begin{prop}[Group actions {\cite[4.1]{NSS12}}]
  \label{InfinityAction}
  Let $\mathbf{H}$ an $\infty$-topos and
  $G \in \mathrm{Groups}(\mathbf{H})$ (Prop. \ref{LoopingAndDelooping}).
\vspace{-2mm}
  \item {\bf (i)} An \emph{action} $(X,\rho)$ of $G$
  is an object $X \in \mathbf{H}$
  and homotopy fiber sequence in $\mathbf{H}$ of the form
  \vspace{-2mm}
  \begin{equation}
    \label{InfinityActionHomotopyFiberSequence}
    \xymatrix@R=1em@C=3em{
      X \ar[r]^-{\mathrm{fib}(\rho)} &
      X \!\sslash\! G
      \ar[d]^-{\rho}
      \\
      &
      \mathbf{B}G
      \,,
    }
  \end{equation}

   \vspace{-2mm}
\noindent
  where $\mathbf{B}G$ is the delooping of $G$ \eqref{LoopingAndDelooping}.

 \vspace{-1mm}
    \item {\bf (ii)} The object $X \!\sslash\! G$ appearing in
  \eqref{InfinityActionHomotopyFiberSequence} is,
  equivalently, the
  \emph{homotopy quotient} of the action of $G$ on $V$:
  \vspace{-2mm}
  \begin{equation}
    \label{HomotopyQuotientAsColimit}
    X \!\sslash\! G
    \;\simeq\;
    \underset{\longrightarrow}{\mathrm{lim}}
    \left(
      \xymatrix{
        \ar@{}[r]|-{\cdots}
        &
        X \times G \times G \;
        \ar@<+7pt>[r]
        \ar@{<-}@<+3.5pt>[r]
        \ar[r]
        \ar@{<-}@<-3.5pt>[r]
        \ar@<-7pt>[r]
        &
     \;   X \times G \;
        \ar@<+3.5pt>[r]
        \ar@{<-}[r]
        \ar@<-3.5pt>[r]
        &
      \;  X
      }
    \right).
  \end{equation}

\vspace{-3mm}
 \item {\bf (iii)}  Hence the $\infty$-category of $G$-actions
 is, equivalently, the slice $\infty$-topos (Prop. \ref{SliceInfinityTopos}) of $\mathbf{H}$ over
 $\mathbf{B}G$:
\vspace{-2mm}
  \begin{equation}
    \label{CategoryOfGActions}
    G\mathrm{Actions}(\mathbf{H})
    \;\simeq\;
    \mathbf{H}_{/\mathbf{B}G}
    \mathrlap{
      \;\;\;\;\in\;
      \mathrm{Categories}_\infty\;.
    }
  \end{equation}
  \end{prop}

We record the following immediate but important aspect of this characterization:

\begin{lemma}[Homotopy quotient maps are effective epimorphisms]
  \label{HomotopyQuotientMapsAreEffectiveEpimorphisms}
  Let $\mathbf{H}$ be an $\infty$-topos,
  $G \in \mathrm{Groups}(\mathbf{H})$ (Prop. \ref{LoopingAndDelooping}),
  and $(X,\rho) \in G\mathrm{Actions}(\mathbf{H})$ (Prop. \ref{InfinityAction}).
 Then the quotient morphism from $X$ to its homotopy quotient \eqref{HomotopyQuotientAsColimit}
  is an effective epimorphism (Def. \ref{EffectiveEpimorphisms}):
  \vspace{-2mm}
  $$
    \xymatrix{
      X \ar@{->>}[rr]^{
        \mathrm{fib}(\rho)
      } && X \!\sslash\! G
    }.
  $$
\end{lemma}
\begin{proof}
  By \eqref{InfinityActionHomotopyFiberSequence} in Prop. \ref{InfinityAction},
  the quotient map sits in a homotopy pullback square of the form
   \vspace{-2mm}
  $$
    \xymatrix@R=1.2em{
      X
      \ar[d]
      \ar[rr]^{\mathrm{fib}(\rho)}
      \ar@{}[drr]|-{\mbox{\tiny\rm(pb)}}
      &&
      X \!\sslash\! G
      \ar[d]^-{\rho}
      \\
      \ast \ar[rr] && \mathbf{B}G
    }
  $$

   \vspace{-2mm}
\noindent
  The bottom morphism is an effective epimorphism (Example \ref{PointInDeloopingIsEffectiveEpi}).
  Since these are preserved by pullback
  (Lemma \ref{EffectiveEpimorphismsArePreservedByPullback}), the claim follows.
\hfill \end{proof}

\begin{example}[Left multiplication action]
  \label{LeftMultiplicationAction}
  Let $\mathbf{H}$ be an $\infty$-topos (Def. \ref{InfinityTopos})
  and $G \in \mathrm{Groups}(\mathbf{H})$ (Prop. \ref{LoopingAndDelooping}).
  The defining looping relation \eqref{GroupIsLooping}
  exhibits, by comparison with \eqref{InfinityActionHomotopyFiberSequence},
  an action of $G$ on itself:
   \vspace{-2mm}
  $$
    \xymatrix@R=12pt{
      G \ar[rr]^{ \mathrm{fib}(\rho_\ell) }
      &&
      \ast \ar[d]^-{ \rho_\ell }
      \\
      && \mathbf{B}G
    }
  $$

   \vspace{-2mm}
\noindent
  This is the \emph{left multiplication action} with
  $
    G \sslash G \;\simeq\; \ast
    \,.
  $
\end{example}
\begin{example}[Adjoint action]
  \label{AdjointAction}
  Let $\mathbf{H}$ be an $\infty$-topos (Def. \ref{InfinityTopos})
  and $G \in \mathrm{Groups}(\mathbf{H})$ (Prop. \ref{LoopingAndDelooping}).
    Then the free loop space object
  $\mathcal{L}\mathbf{B}G$ of the delooping $\mathbf{B}G$ \eqref{LoopingDeloopingEquivalence},
  defined by the Cartesian square
 \vspace{-1mm}
  $$
    \xymatrix@R=17pt{
      \mathcal{L}\mathbf{B}G
      \ar[rr]
      \ar[d]_-{ \rho_{\mathrm{ad}} }
      \ar@{}[drr]|-{ \mbox{\tiny\rm(pb)} }
      &&
      \mathbf{B}G
      \ar[d]^-{ \Delta }
      \\
      \mathbf{B}G
      \ar[rr]_-{\Delta}
      &&
      \mathbf{B}G
      \times
      \mathbf{B}G
    }
  $$

   \vspace{-3mm}
\noindent
  sits in a homotopy fiber sequence of the form
   \vspace{-3mm}
  $$
    \xymatrix@R=14pt{
      G
      \ar[rr]^-{ \mathrm{fib}(\rho_{\mathrm{ad}}) }
      &&
      \mathcal{L}\mathbf{B}G
      \ar[d]^-{ \rho_{\mathrm{ad}} }
      \\
      &&
      \mathbf{B}G
    }
  $$

   \vspace{-2mm}
\noindent
  By comparison with \eqref{InfinityActionHomotopyFiberSequence}, this
  exhibits an action of $G$ on itself. This is the \emph{adjoint action}
  with
  $ G \!\sslash_{\!\mathrm{ad}}\! G \;\simeq\; \mathcal{L}\mathbf{B}G$.
\end{example}
\begin{defn}[Equivariant maps]
  \label{Equivariance}
  By the functoriality/universality of the homotopy fiber construction
  in \eqref{InfinityActionHomotopyFiberSequence}
  and using the equivalence \eqref{CategoryOfGActions},
  we have the $\infty$-functor that assigns the underlying objects of
  the $G$-actions in Def. \ref{InfinityAction}:
  \begin{equation}
    \xymatrix{
      G\mathrm{Actions}(\mathbf{H})
      \simeq
      \mathbf{H}_{/_{\mathbf{B}G}}
      \ar[r]^-{ \mathrm{fib} }
      &
      \mathbf{H}
    }.
  \end{equation}
  With two $G$-actions $(X_i, \rho_i)$
  given, we say that a morphism $X_1 \to X_2  \; \in \mathbf{H}$
  between their underlying objects \emph{is equivariant} if it
  lifts through this functor, hence if it is the image of a morphism
  $(X_1, \rho_1) \to (X_2, \rho_2) \;\in \mathrm{G}\mathrm{Actions}(\mathbf{H})$.
\end{defn}

\begin{example}[Group division is equivariant under diagonal left and adjoint action]
  Let $\mathbf{H}$ be an $\infty$-topos (Def. \ref{InfinityTopos}) and
  $G \in \mathrm{Groups}(\mathbf{H})$ (Prop. \ref{LoopingAndDelooping}).
    Then the group division operation (Example \ref{GroupOperation})
  is equivariant (Def. \ref{Equivariance}) with respect to
  the diagonal
  left multiplication action $\rho_\ell$ (Example \ref{LeftMultiplicationAction})
  on its domain and the adjoint action $\rho_{\mathrm{ad}}$ (Example \ref{AdjointAction})
  on its codomain:
  \vspace{-3mm}
  \begin{equation}
    \label{EquivarianceOfGroupDivisionUnderDiagonalLeftAndAdjointAction}
    \xymatrix{
      (G, \rho_{\ell})
      \times
      (G, \rho_\ell)
      \ar[rr]^-{ (-) \cdot (-)^{-1} }
      &&
      (G, \rho_{\mathrm{ad}})
    }
    \mathrlap{
      \;\;\;\;\;\;\;\;
      \in
      G\mathrm{Actions}(\mathbf{H})\;.
    }
  \end{equation}
\end{example}
\begin{proof}
  Observe the following pasting of Cartesian squares:
   \vspace{-2mm}
$$
  \xymatrix{
    G \times G
    \ar[rr]^-{
      (-) \cdot (-)^{-1}
    }
    \ar[d]_-{
      (-)^{-1}\cdot (-) \circ \sigma
    }
    &&
    G
    \ar[rr]
    \ar[d]
    &&
    \ast
    \ar[d]
    \\
    G
    \ar[d]
    \ar[rr]
    &&
    \mathcal{L}\mathbf{B}G
    \ar[rr]
    \ar[d]
    &&
    \mathbf{B}G
    \ar[d]^-{\Delta}
    \\
    \ast
    \ar[rr]
    &&
    \mathbf{B}G
    \ar[rr]^-{\Delta}
    &&
    \mathbf{B}G \times \mathbf{B}G
  }
$$
The middle horizontal composite, regarded as a morphism in the
slice over $\mathbf{B}G$ and hence
as a morphism of $G$-actions \eqref{InfinityActionHomotopyFiberSequence},
gives \eqref{EquivarianceOfGroupDivisionUnderDiagonalLeftAndAdjointAction}.
\hfill \end{proof}

\begin{prop}[Restricted and induced group actions]
  \label{PullbackAction}
  Let $\mathbf{H}$ be an $\infty$-topos.
  Then, for
  $\phi : H \to G$ a morphism
  in $\mathrm{Groups}(\mathbf{H})$ (Prop. \ref{LoopingAndDelooping}),
  there is a
  triple of adjoint $\infty$-functors (Def. \ref{AdjointInfinityFunctors}) between
  the corresponding $\infty$-categories of group actions
  (Prop. \ref{InfinityAction})
  \vspace{-2mm}
  \begin{equation}
    \label{InducedGroupActionAdjunction}
    \xymatrix{
      H\mathrm{Actions}(\mathbf{H})
      \;\;
      \ar@<+14pt>[rr]^{
        \overset{
          \mbox{
            \tiny
            \color{darkblue}
            \bf
            ``left-induced''
          }
        }{
          \mathbf{B}\phi_!
        }
      }
      \ar@{<-}@<+0pc>[rr]|{ \; \mathbf{B}\phi^\ast \; }^-{
        \raisebox{2pt}{\scalebox{.8}{\small $\bot$}}
      }      \ar@<-14pt>[rr]^-{
        \raisebox{-1pt}{\scalebox{.8}{\small $\bot$}}
      }
_{
        \underset{
          \mbox{
            \tiny
            \color{darkblue}
            \bf
            ``right-induced''
          }
        }{
          \mathbf{B}\phi_\ast
        }
      }
      &&
     \;\; G\mathrm{Actions}(\mathbf{H})
    }
  \end{equation}
  such that $\mathbf{B}\phi^\ast$ preserves the object being acted on
  (``restricted action'').
\end{prop}
\begin{proof}
  By \eqref{CategoryOfGActions} in
  Prop. \ref{InfinityAction}, an adjoint triple (Def. \ref{AdjointInfinityFunctors})
  of the form \eqref{InducedGroupActionAdjunction}
  is given by base change (Prop. \ref{BaseChange})
  of homotopy quotients \eqref{HomotopyQuotientAsColimit}
  along the delooped morphism
  $\mathbf{B}\phi$ (Prop. \ref{LoopingAndDelooping}).
  This means that $\mathbf{B}\phi^\ast$ is
  given by sending the homotopy fiber sequence \eqref{InfinityActionHomotopyFiberSequence}
  corresponding to a $G$-action to
  the following homotopy pullback
  (Prop. \ref{LoopingAndDelooping}):
  \vspace{-2mm}
  \begin{equation}
    \label{PullbackOfInfinityActions}
    \xymatrix@C=7em{
      X
      \ar[r]_-{\mathrm{fib}(\phi^\ast \rho)}
      \ar@/^1pc/[rr]^-{ \mathrm{fib}(\phi) }
      &
      X \!\sslash\! H
      \ar@{}[dr]|-{ \mbox{\tiny (pb) } }
      \ar[d]_-{ \phi^\ast \rho }
      \ar[r]
      &
      X \!\sslash\! G
      \ar[d]^-{ \rho }
      \\
      &
      \mathbf{B}H
      \ar[r]_-{ \mathbf{B}\phi }
      &
      \mathbf{B}G
    }
  \end{equation}

   \vspace{-2mm}
\noindent  That this preserves the object $X$ being acted on, as indicated,
  follows by the pasting law (Prop. \ref{PastingLaw}).
\hfill \end{proof}

\begin{defn}[Automorphism group]
  \label{AutomorphismGroup}
  Let $\mathbf{H}$ be an $\infty$-topos and
  $F \in \mathbf{H}$ an object. Then the
  \emph{automorphism group}
  $\mathrm{Aut}(F) \in \mathrm{Groups}(\mathbf{H})$
  of $F$ is the looping (Prop. \ref{LoopingAndDelooping})
  of the (-1)-image \eqref{nImageFactorization}
  of the classifying map \eqref{ObjectClassifier} of $F$:
  \vspace{-1mm}
  $$
    \xymatrix@C=3em{
      \ast
      \ar@{->>}[rr]^-{ \mbox{\tiny (-1)-conn.} }
      \ar@/_1.3pc/[rrrr]_-{ \vdash F }
      &&
      \mathbf{B} \mathrm{Aut}(F)
     \; \ar@{^{(}->}[rr]^-{ \mbox{\tiny (-1)-trunc.} }
      &&
      \mathrm{Objects}_\kappa\;.
    }
  $$
  The canonical action of this group (Prop. \ref{InfinityAction})
  on $V$ is exhibited, via \eqref{InfinityActionHomotopyFiberSequence},
  by the left square of the
  following pasting composite of Cartesian squares:
     \vspace{-2mm}
  \begin{equation}
    \label{AutomorphismAction}
    \raisebox{20pt}{
    \xymatrix@C=5em@R=1.5em{
      F
      \ar[d]
      \ar[r]^-{ \mathrm{fib}(\rho_{\mathrm{Aut}}) }
      \ar@{}[dr]|-{ \mbox{\tiny(pb)} }
      &
      F \!\sslash\! \mathrm{Aut}(F)
      \ar[r]
      \ar[d]^-{\rho_{\mathrm{Aut}} }
      \ar@{}[dr]|-{ \mbox{\tiny(pb)} }
      &
      \widehat{\mathrm{Objects}}_\kappa
      \ar[d]
      \\
      \ast
      \ar@/_1.5pc/[rr]_{ \vdash F }
      \ar@{->>}[r]
      &
      \mathbf{B} \mathrm{Aut}(F)
      \ar@^{{(}->}[r]
      &
      \mathrm{Objects}_\kappa
      \;,
    }
    }
  \end{equation}

     \vspace{-2mm}
\noindent
  where we use the pasting law (Prop. \ref{PastingLaw})
  to identify $F$ as the homotopy fiber of $\rho_{\mathrm{Aut}}$.
\end{defn}

\begin{prop}[Automorphism group is universal]
\label{AutomorphimsGroupIsUniversal}
Let $\mathbf{H}$ be an $\infty$-topos,
$G \in \mathrm{Groups}(\mathbf{H})$ (Prop. \ref{LoopingAndDelooping}),
and $(X, \rho) \in G\mathrm{Actions}(\mathbf{H})$ (Def. \ref{InfinityAction}).
Then there is a group homomorphism from $G$
to the automorphism group (Def. \ref{AutomorphismGroup})
 \vspace{-2mm}
$$
  \xymatrix{
    G
    \ar[r]^-{ i_\rho }
    &
    \mathrm{Aut}(X)
  }
$$

 \vspace{-1mm}
\noindent
such that the action $\rho$
is the restricted action (Prop. \ref{PullbackAction}) along $i_\rho$
of the canonical automorphism action \eqref{AutomorphismAction},
i.e., such that there is a Cartesian square of this form:
$$
  \xymatrix{
    X \!\sslash\! G
    \ar[rr]
    \ar@{}[drr]|-{ \mbox{\tiny\rm(pb)} }
    \ar[d]_-{ \rho }
    &&
    X \!\sslash\! \mathrm{Aut}(X)
    \ar[d]^-{ \rho_{\mathrm{Aut}} }
    \\
    \mathbf{B} G
    \ar[rr]_-{ \mathbf{B} i_\rho }
    &&
    \mathbf{B} \mathrm{Aut}(X)
  }
$$
\end{prop}
\begin{proof}
  Let $\kappa$ be a regular cardinal such that $X$ is
  $\kappa$-small, and consider the following solid diagram
  of classifying maps
  \eqref{ObjectClassifier}
  for $\rho$, $\rho_{\mathrm{Aut}}$
  and for $X$:
   \vspace{-2mm}
  $$
  \xymatrix@R=12pt{
    X
    \ar[rr]
    \ar[dr]
    \ar[dd]
    &&
    X \!\sslash\! \mathrm{Aut}(X)
    \ar[dr]
    \ar[dd]|-{\phantom{AA}}
    \\
    &
    X \!\sslash\! G
    \ar[dd]
    \ar[rr]
    \ar@{-->}[ur]
    &&
    \widehat{\mathrm{Objects}}_\kappa
    \ar[dd]
    \\
    \ast
    \ar[dr]|-{
      \mbox{\tiny(-1)-connected}
    }
    \ar[rr]|-{\phantom{AA}}
    &&
    \mathbf{B}\mathrm{Aut}(X)
    \ar@{^{(}->}[dr]|-{\mbox{\tiny (-1)-truncated}}
    \\
    &
    \mathbf{B} G
    \ar@{-->}[ur]|-{\mathbf{B} i_\rho}
    \ar[rr]_{ \vdash \rho }
    &&
    \mathrm{Objects}_\kappa
  }
$$

 \vspace{-2mm}
\noindent
Here the bottom square homotopy-commutes by the essential uniqueness
of the classifying map $\vdash  X$ \eqref{ObjectClassifier}.
Hence the dashed lift exists
essentially uniquely \eqref{ConnectedTruncatedLiftingProperty},
by the connected/truncated factorization system
(Prop. \ref{nConnectednTruncatedFactorizationSystem}).
\hfill \end{proof}

\medskip

\noindent {\bf Principal bundles.}
\begin{prop}[Principal bundles {\cite[Thm. 3.17]{NSS12}}]
  \label{ClassificationPrincipalInfinityBundles}
  Let $\mathbf{H}$ be an $\infty$-topos,
  $X \in \mathbf{H}$,  and
  $G \in \mathrm{Groups}(\mathbf{H})$ (Prop. \ref{LoopingAndDelooping}).
  Then $G$-principal $\infty$-bundles $P \to X$
  over $X$ are, equivalently, given by
  \emph{classifying maps} $\vdash P : X \to \mathbf{B}G$.
  Forming their homotopy fibers
  \vspace{-2mm}
  $$
    \xymatrix@R=1.5em{
      P
      \ar[d]_-{ \mathrm{fib}(\vdash P) }
      \\
      X
      \ar[r]_-{ \vdash P }
      &
      \mathbf{B}G
    }
  $$

   \vspace{-2mm}
\noindent
  constitutes
  an equivalence of $\infty$-groupoids:
  \vspace{-3mm}
  \begin{equation}
    \label{ClassificationOfPrincipalBundles}
    \xymatrix@R=-2pt{
      G\mathrm{Bundles}_X(\mathbf{H})
      \ar@{<-}[rr]^-{ \mathrm{fib} }_-{ \simeq }
      &&
      \mathbf{H}
      (
        X, \mathbf{B}G
      )\;.
      \\
      P
      &\longmapsto&
      \vdash P
    }
  \end{equation}
\end{prop}

\begin{remark}[Principal base spaces are homotopy quotients]
  \label{PrincipalBaseSpacesAreHomotopyQuotients}
  Comparison of the abstract
  characterization of {\bf (i)} group actions (Prop. \ref{InfinityAction})
  and {\bf (ii)} principal bundles
  (Prop. \ref{ClassificationPrincipalInfinityBundles}),
  reveals that these are about one and the same
  abstract concept, just viewed from two different perspectives:
  In an $\infty$-topos, every $G$-principal bundle
  is a $G$-action whose homotopy quotient is the given base space;
  and, conversely, every $G$-action is that of a principal bundle
  over its homotopy quotient:
\vspace{-1mm}
  $$
    \xymatrix@R=1em@C=1pt{
      \mbox{
        \tiny
        \color{darkblue}
        \bf
        \begin{tabular}{c}
                 principal
          \\
          G-bundle
        \end{tabular}
      }
          &
     \hspace{1cm}P  \raisebox{8pt}{\;\;\;\begin{rotate}{270} $\curvearrowright$\end{rotate}}\;\;\;\;\; G
      \ar[d]
      &&
      \mbox{
        \tiny
        \color{darkblue}
        \bf
        \begin{tabular}{c}
          G-action
        \end{tabular}
      }
      \\
      \mbox{
        \tiny
        \color{darkblue}
        \bf
        \begin{tabular}{c}
          base
          \\
          space
        \end{tabular}
      }
      &
      X
      \ar@{}[r]|-{ \simeq }
      &
      P \!\sslash\! G
      &
      \mbox{
        \tiny
        \color{darkblue}
        \bf
        \begin{tabular}{c}
          homotopy
          \\
          quotient
        \end{tabular}
      }
    }
  $$

 \vspace{-1mm}
\noindent
  Notice (see \cite[3.1]{NSS12} for exposition) that it is the higher
  geometry inside an $\infty$-topos that makes this work.
\end{remark}

\begin{defn}[Atiyah groupoid]
  \label{TheAtiyahGroupoid}
  Let $\mathbf{H}$ be an $\infty$-topos (Def. \ref{InfinityTopos}),
  $X \in \mathbf{H}$,
  $G \in \mathrm{Groups}(\mathbf{H})$ (Prop. \ref{LoopingAndDelooping}),
  and $P \in G\mathrm{Bundles}_X$ (Prop. \ref{ClassificationPrincipalInfinityBundles}).
  Then the \emph{Atiyah groupoid} of $P$
  is the groupoid
  $
    \mathrm{At}_\bullet(P)
    \;\in\;
    \mathrm{Groupoids}(\mathbf{H})
  $
  (Def. \ref{InHigherToposGroupoids})
  whose corresponding stack with atlas (via Prop. \ref{EtaleGroupoidsAndEtaleAtlases})
  is the (-1)-image projection (Example \ref{EffectiveEpimorphismMonomorphismFactorization})
  of the bundle's classifying map $\vdash P$ \eqref{ClassificationOfPrincipalBundles}:
  \begin{equation}
    \label{AtiyahGroupoidFactorization}
    \xymatrix{
      X
      \;
      \ar@{->>}[rr]
      \ar@/_1.3pc/[rrrr]_-{
        \vdash P
      }
      &&
      \;\mathcal{A}\!t(P)\;
      \ar@{^{(}->}[rr]
      &&
      \mathbf{B}G
    }.
  \end{equation}

\end{defn}

\medskip
\noindent {\bf Fiber bundles.}
\begin{defn}[Fiber bundle]
  \label{FiberBundle}
  Let $\mathbf{H}$ be an $\infty$-topos (Def. \ref{InfinityTopos}).
  \vspace{-1mm}
\item {\bf (i)}    morphism $Y \overset{p}{\longrightarrow} X$
  in $\mathbf{H}$
  is a \emph{fiber bundle} with \emph{typical fiber}
  $F \in \mathbf{H}$ if there exists an effective epimorphism
  $\xymatrix{ U \ar@{->>}[r]^i & X }$ (Def. \ref{EffectiveEpimorphisms})
  and a Cartesian square (Notation \ref{CartesianSquares}) of the form
  $$
    \xymatrix@C=5em@R=1.5em{
      U \times F
      \ar[d]
      \ar[r]
      \ar@{}[dr]|-{\mbox{\tiny(pb)}}
      &
      Y
      \ar[d]^-{ p }
      \\
      U
      \ar@{->>}[r]_-{ i }
      &
            X
    }
  $$
 \vspace{-6mm}
\item {\bf (ii)} We write
  $$
    F\mathrm{FiberBundles}_{X}(\mathbf{H})
    \;\subset\;
    \mathrm{Core}\big(\mathbf{H}_{/X}\big)
    \;\in\;
    \mathrm{Groupoids}_\infty
  $$
  for the full $\infty$-groupoid of
  the core \eqref{InfinityCategoryOfInfinityCategories}
  of the slice $\mathbf{H}_{/X}$ over $X$ (Prop. \ref{SliceInfinityTopos})
  on the $F$-fiber bundles.
\end{defn}

\begin{prop}[Classification of fiber bundles {\cite[Prop. 4.10]{NSS12}}]
  \label{FiberBundlesClassified}
  Let $\mathbf{H}$ be an $\infty$-topos (Def. \ref{InfinityTopos})
  and $X, F \in \mathbf{H}$.
  Then fiber bundles over $X$ (Def. \ref{FiberBundle})
  with typical fiber $F$ are equivalent to
  morphisms $X \longrightarrow \mathbf{B} \mathrm{Aut}(F)$
  from $X$ to the delooping (Prop. \ref{LoopingAndDelooping})
  of the automorphism group (Def. \ref{AutomorphismGroup})
  of $F$:
   \vspace{-4mm}
  \begin{equation}
    \label{FiberBundleClassification}
    \xymatrix@R=-2pt{
      \\
      F \mathrm{FiberBundles}_{X}(\mathbf{H})
      \ar[rr]^-{ \simeq }
      &&
      \mathbf{H}
      \big(
        X
        \,,\,
        \mathbf{B}\mathrm{Aut}(F)
      \big)
      \\
      E
      \ar@{}[rr]|-{ \longmapsto }
      &&
      \vdash\, E
    }
  \end{equation}
\end{prop}
\begin{proof}
  Let $\kappa$ be a regular cardinal such that $F$ is
  $\kappa$-small. Then, by assumption, we have the following solid diagram
  of classifying maps \eqref{ObjectClassifier}:
   \vspace{-2mm}
  $$
    \xymatrix@C=4em@R=.8em{
      U \times F
      \ar[dd]
      \ar[rr]^{ \mathrm{pr}_2 }
      \ar[dr]
      &&
      F \!\sslash\! \mathrm{Aut}(F)
      \ar[dd]|-{
        \phantom{AA}
      }
      \ar[dr]
      \\
      &
      E
      \ar[dd]
      \ar[rr]
      &&
      \widehat{\mathrm{Objects}}_\kappa
      \ar[dd]
      \\
      U
      \ar[rr]|-{
        \phantom{AA}
      }
      \ar@{->>}[dr]|-{
        \mbox{
          \tiny
          \rm
          (-1)-connected
        }
        \;\;\;\;\;
      }
      \ar@{}[rr]|{
        {\phantom{A}}
      }
      &&
      \mathbf{B} \mathrm{Aut}(F)
      \;\;
      \ar@{^{(}->}[dr]|-{
        \mbox{
          \tiny
          \rm
          (-1)-truncated
        }
      }
      \\
      &
      X
      \ar[rr]_-{
      }
      \ar@{-->}[ur]^-{ \vdash E }
      &&
      \mathrm{Objects}_{\kappa}
    }
  $$
  Now the (-1)-connected/(-1)-truncated factorization system
  (Prop. \ref{nConnectednTruncatedFactorizationSystem}) implies
  that the dashed morphism exists essentially uniquely \eqref{ConnectedTruncatedLiftingProperty}.

  It just remains to see that this assignment is independent of the
  choice of $U$: For $\xymatrix@C=1.5em{U' \! \ar@{->>}[r] & \! X}$
  any other effective epimorphism with $(\vdash E)'$ the associated
  classifying map as above, observe that the fiber product
  $\xymatrix@C=1.5em{ U \times_X U' \! \ar@{->>}[r] & \! X}$ is again an effective
  epimorphism, since the class of effective epimorphisms
  is closed under pullbacks as well as under composition
  (Lemma \ref{EffectiveEpimorphismsArePreservedByPullback}).
  Therefore $\vdash E$ and $(\vdash E)'$ are jointly lifts in
  a diagram as above but with $U \times_X U'$ in the top left.
  Hence, by the essential uniqueness of lifts in the
  connected/truncated orthogonal factorization system, they are equivalent,
  $(\vdash E) \simeq (\vdash E)'$,
  in an essentially unique way.
\hfill \end{proof}

\begin{notation}[Associated bundles]
We say that

\noindent {\bf (i)}
the morphism $\vdash E$ in \eqref{FiberBundleClassification}
is the \emph{classifying map} of $E$
and

\noindent {\bf  (ii)} that $E$ is \emph{associated} to the
$\mathrm{Aut}(F)$-principal bundle which is classified by $\vdash E$
according to Prop. \ref{ClassificationPrincipalInfinityBundles}.
\end{notation}

\begin{remark}[Twisted cohomology in slice $\infty$-toposes]
  \label{TwistedCohomology}
  Prop. \ref{FiberBundlesClassified} implies
  (together with the universal property of the
  pullback) that sections $\sigma$ of $A$-fiber bundles
  $E$ over some $X$ are, equivalently, lifts $c$ of the classifying
  map $c \;:=\; \vdash E$ \eqref{FiberBundleClassification} through
  $\rho_{\mathrm{Aut}}$ \eqref{AutomorphismAction}:

  \vspace{-.5cm}

\begin{equation}
  \label{SectionOfFiberBundleAsLiftOfClassifyingMap}
  \raisebox{20pt}{
  \xymatrix{
    &&
    A \!\sslash\! \mathbf{Aut}(A)
    \ar[d]^-{ \rho_{\mathrm{Aut}} }
    \\
    X
    \ar[rr]_-{
      \underset{
        \raisebox{-3pt}{
          \tiny
          \color{darkblue}
          \bf
          classifying map
        }
      }{
        \tau \;:=\; \vdash E
      }
    }
    \ar@{-->}[urr]^{
      \overset{
        \mathllap{
        \raisebox{-1pt}{
          \tiny
          \color{darkblue}
          \bf
          \begin{tabular}{c}
            lift of
            \\
            classifying map
          \end{tabular}
        }
        \!\!\!\!
        }
      }{
        c
      }
    }
    &&
    \mathbf{B}
    \mathbf{Aut}(A)
  }
  }
  \;\;\;\;\;\;\simeq\;\;\;\;\;\;\;
  \raisebox{20pt}{
  \xymatrix@C=3em{
    &
    \overset{
      \mathclap{
      \raisebox{3pt}{
        \tiny
        \color{darkblue}
        \bf
        \begin{tabular}{c}
          associated
          bundle
        \end{tabular}
      }
      }
    }{
      E
    }
    \ar[d]_-p
    \ar[r]
    \ar@{}[dr]|-{ \mbox{\tiny(pb)} }
    &
    A \!\sslash\! \mathbf{Aut}(A)
    \ar[d]^-{ \rho_{\mathrm{Aut}} }
    \\
    X
    \ar@{=}[r]
    \ar@{-->}[ur]^{
      \overset{
        \mathllap{
        \raisebox{4pt}{
          \tiny
          \color{darkblue}
          \bf
          \begin{tabular}{c}
            section
          \end{tabular}
        }
        }
      }{
        \sigma
      }
    }
    &
    X
    \ar[r]_-{ \tau }
    &
    \mathbf{B}
    \mathbf{Aut}(A)
  }
  }
\end{equation}

\vspace{-.2cm}
\item {\bf (i)} If $A$ is regarded here as a coefficient object for $A$-cohomology
\eqref{IntrinsicCohomologyOfAnInfinityTopos}, then such a section
$\sigma$ is a locally $A$-valued cocycle, which is ``twisted''
over $X$ according to the classifying map $\tau$.
Hence such a $\sigma$ is a cocycle in
(non-abelian)
\emph{$\tau$-twisted cohomology} \cite[4.2]{NSS12}.
But the left hand side of \eqref{SectionOfFiberBundleAsLiftOfClassifyingMap}
is, equivalently, a morphism \eqref{MorphismInSliceCategory}
in the slice $\infty$-topos (Prop. \ref{SliceInfinityTopos})
$\mathbf{H}_{/\mathbf{B}\mathbf{Aut}(A)}$.
It follows that
{\it twisted cohomology is
the intrinsic cohomology \eqref{IntrinsicCohomologyOfAnInfinityTopos}
of slice $\infty$-toposes}:

\vspace{-.8cm}

\begin{equation}
  \label{TwistedCohomologySliceHom}
  \overset{
    \mathclap{
    \raisebox{5pt}{
      \tiny
      \color{darkblue}
      \bf
      \begin{tabular}{c}
        $\mathbf{\tau}$-twisted
        \\
        cohomology
      \end{tabular}
    }
    }
  }{
    H^\tau
    \big(
      X
      \,,\,
      A
    \big)
  }
  \;:=\;
  \pi_0\,
  \mathbf{H}_{/\mathbf{B} \mathbf{Aut}(A)}
  \Big(
    (X,\tau)
    \,,\,
    (A \!\sslash\! \mathbf{Aut}(A), \rho_{\mathrm{Aut}})
  \Big)
  \;\;
  \simeq
  \;\;
  \left\{
  \raisebox{23pt}{
  \xymatrix@C=14pt{
    X
    \ar[dr]_-{\tau}^-{\ }="t"
    \ar@{-->}[rr]^-{
      \overset{
        \raisebox{3pt}{
          \tiny
          \color{darkblue}
          \bf
          cocycle
        }
      }{
        c
      }
    }_>>>>>>{\ }="s"
    &&
    A \!\sslash\! \mathbf{Aut}(A)
    \ar[dl]^-{ \rho_{\mathrm{Aut}} }
    \\
    &
    \mathbf{B}
    \mathbf{Aut}(A)
    \ar@{=>} "s"; "t"
  }
  }
  \right\}_{\!\!\!\!\big/\sim}
\end{equation}

\vspace{-.3cm}

\item {\bf (ii)} By the universality of $\mathbf{Aut}(A)$
(Prop. \ref{AutomorphimsGroupIsUniversal}), this holds
for slicing over \emph{any} pointed connected object $\mathbf{B}G$
\eqref{LoopingDeloopingEquivalence}.

  \vspace{-1.5mm}
\item {\bf (iii)} If the base object is not connected, the intrinsic cohomology
of its slice may be thought of as a mixture of twisted and parametrized
cohomology. We encounter an example of this in Def. \ref{EtaleCohomologyOfVFolds}
below.
\end{remark}

\begin{remark}[Twisted cohomology as global sections]
  \label{TwistedCohomologyAsGlobalSections}
  The $\infty$-groupoid of sections of the associated bundle
  $\xymatrix@C=12pt{E := \tau^\ast(A \!\sslash\! G) \ar[r]^-p & X}$
  in \eqref{SectionOfFiberBundleAsLiftOfClassifyingMap},
  is equivalently its image $\Gamma_X(E)$ under the base geometric morphism
  (Prop. \ref{BaseGeometricMorphism})

  \vspace{-.4cm}
  $$
    \xymatrix{
      \mathbf{H}_{/X}
        \ar@{<-}@<+6pt>[rr]^-{ \Delta_X }
        \ar@<-6pt>[rr]_-{ \Gamma_X }^-{ \bot }
        &&
      \mathrm{Groupoids}_\infty
    }
  $$
  \vspace{-.4cm}

  \noindent
  of the slice $\infty$-topos $\mathbf{H}_X$ (Prop. \ref{SliceInfinityTopos}),
  in that
  (by Prop. \ref{TensoringOfInfinityToposesOverInfinityGroupoids})
  $
    \Gamma_X(E)
    \;\simeq\;
    \mathbf{H}_X
    \big(
      \mathrm{id}_X
      \,,\,
      p
    \big)
    \,.
  $
  Hence the $\tau$-twisted cohomology \eqref{TwistedCohomologySliceHom}
  of $X$ is equivalently the set of connected components of the
  $\infty$-groupoid of global sections:

  \vspace{-.4cm}
  \begin{equation}
    \label{TwistedCohomologyIsGlobalSectionOfAssociatedBundle}
    \mathbf{H}^\tau
    \big(
      X;
      \,
      A
    \big)
    \;\simeq\;
    \pi_0
    \,
    \Gamma_X\big( \tau^\ast(A \!\sslash\! G) \big)
    \,.
  \end{equation}
\end{remark}

\begin{remark}[Twisted abelian cohomology in tangent $\infty$-toposes]
  \label{AbelianTwistedCohomology}
  Let $\mathbf{H}$ be an $\infty$-topos (Def. \ref{InfinityTopos}).

  \vspace{-1.5mm}
  \item {\bf (i)} Notice that
  the intrinsic cohomology \eqref{IntrinsicCohomologyOfAnInfinityTopos}
  of $\mathrm{Bundles}(\mathbf{H})$ (Example \ref{BundleMorphismsCoveringBaseMorphisms})
  is still twisted cohomology
  as in Remark \ref{TwistedCohomology}, just up to a change in perspective:
  now the twisting $\tau$ is encoded not in the domain object,
  but in the cocycles
  on these (a morphism
  of the form
  $\xymatrix@C=11pt{\!\mathrm{id}_X \ar[r] & \rho_{\mathrm{Aut}} \!}$
  in $\mathrm{Bundles}(\mathbf{H})$ is still manifestly given by the
  diagrams in \eqref{SectionOfFiberBundleAsLiftOfClassifyingMap}).

    \vspace{-2mm}
 \item {\bf (ii)}   Therefore, similarly, the intrinsic cohomology
  \eqref{IntrinsicCohomologyOfAnInfinityTopos}
   in the tangent $\infty$-topos
  $\mathrm{SpectralBundles}(\mathbf{H})$ (Example \ref{TangentInfinityTopos})
  is twisted cohomology
  with local coefficients being spectra
  \cite[4.1]{dcct}\cite{ABGHR14}\cite{GS19a}\cite{GS19b},
  hence is \emph{twisted abelian cohomology}.

  \vspace{-1.5mm}
 \item {\bf (iii)}  In the case that $\mathbf{H} = \mathrm{Groupoids}_\infty$,
  the base tangent $\infty$-topos
    \vspace{-1mm}
  \begin{equation}
    \label{BundlesOfPlainSpectra}
    T \mathrm{Groupoids}_\infty
    \;=\;
    \mathrm{SpectralBundles}
    \big(
      \mathrm{Groupoids}_\infty
    \big)
  \end{equation}

    \vspace{-2mm}
\noindent   is the topic of traditional parametrized stable homotopy theory
  \cite{James95}\cite{MaySigurdsson06}\cite[2]{ABGHR14}\cite{BM19}
  and its intrinsic cohomology theory \eqref{IntrinsicCohomologyOfAnInfinityTopos} is traditional twisted generalized cohomology
  \cite{Douglas05}\cite{AndoBlumbergGepner10}.
\end{remark}

\medskip

\noindent {\bf Fixed points and fixed loci.}
\begin{defn}[Fixed points and fixed loci]
  \label{FixedPoints}
  Let $\mathbf{H}$ be an $\infty$-topos,
  $G \in \mathrm{Groups}(\mathbf{H})$
  (Prop. \ref{LoopingAndDelooping})
  and $(X,\rho) \in G\mathrm{Actions}(\mathbf{H})$ (Prop. \ref{InfinityAction}).

  \vspace{-.3cm}

  \begin{equation}
    \label{HomotopyFixedPoints}
    \hspace{-3.2cm}
 \begin{minipage}[left]{7.5cm}
     \noindent {\bf (i)}
    A \emph{fixed point} of $(X,\rho)$
    is an element
    $\xymatrix{\ast \ar[r]^-x & X}$
    induced from a section $x \!\sslash\! G$
    of $\rho$ in \eqref{InfinityActionHomotopyFiberSequence},
    as shown on the right
    (where we are using the pasting law, Prop. \ref{PastingLaw},
     and Example \ref{PullbackOfEquivalenceIsEquivalence}
     to identify the top square as Cartesian).
 \end{minipage}
   \hspace{2.5cm}
    \raisebox{32pt}{
    \xymatrix@C=3em@R=1.5em{
      \ast
      \ar[rr]
      \ar[d]_-{ x }
      \ar@{}[drr]|-{ \mbox{\tiny(pb)} }
      &&
      \mathbf{B}G
      \ar[d]_-{ x \!\sslash\! G }
      \ar@{=}@/^1.4pc/[dd]
      \\
      X
      \ar[rr]|-{\, \mathrm{fib}(\rho) \, }
      \ar[d]
      \ar@{}[drr]|-{\mbox{\tiny(pb)}}
      &&
      X \!\sslash\! G
      \ar[d]_-{ \rho }
      \\
      \ast
      \ar[rr]
      &&
      \mathbf{B}G\,,
    }
    }
  \end{equation}

  \vspace{-.3cm}

  \noindent   {\bf (ii)} The \emph{$G$-fixed locus} of $(X,\rho)$ is the object

 \vspace{-2mm}
  \begin{equation}
    \label{HomotopyFixedLocus}
    X^G
    \;:=\;
    \mathbf{B}(G \to \ast)_\ast
    \big(
      (X,\rho)
    \big)
    \;\;\;
    \in
    1\mathrm{Action}(\mathbf{H})
    \;\simeq\;
    \mathbf{H}
    \,,
  \end{equation}
   that is right induced (Prop. \ref{PullbackAction})
   along the unique morphism to the trivial group.
\end{defn}
\begin{example}[Global points of fixed loci are homotopy fixed points]
  The global points of a homotopy-fixed locus $X^G$
  \eqref{HomotopyFixedLocus} are indeed, equivalently,
  the fixed points \eqref{HomotopyFixedPoints}.
  By the adjunction \eqref{InducedGroupActionAdjunction},
  we have the hom-equivalence \eqref{AdjunctionHomEquivalence}
   \vspace{-2mm}
  $$
   \big(
    \xymatrix{
      \ast \ar[r]
      &
      X^G
      \; = \mathbf{B}(G \to 1)_\ast (X,\rho)
    }
    \big)
    \;\;\;\;
    \leftrightarrow
    \;\;\;\;
    \big(
      \xymatrix{
        \mathbf{B}(G \to 1)^\ast (\ast)
        \ar[r]
        &
        (X,\rho)
      }
    \big)
  $$

   \vspace{-2mm}
\noindent
  and, by Prop. \ref{InfinityAction}, the latter morphisms are equivalent
  to homotopy-commuting diagrams of the form
   \vspace{-2mm}
  $$
    \xymatrix@C=3em@R=1.5em{
      \mathbf{B}G
      \ar[dr]_-{
        \mathbf{B}(G \to 1)^\ast(\ast)
      \phantom{A}}^-{\simeq}
      \ar[rr]^{ x \!\sslash\! G }
      &&
      X \!\sslash\! G
      \ar[dl]^-{\rho}
      \\
      & \mathbf{B}G
    }
  $$

   \vspace{-2mm}
\noindent
  This is just the type of diagram
  characterizing homotopy fixed points.
  as seen vertically on the right in \eqref{HomotopyFixedPoints}.
\end{example}
\begin{example}[Fixed loci in $\infty$-groupoids]
  \label{FixedLociInInfinityGroupoids}
  Consider $\mathbf{H} := \mathrm{Groupoids}_\infty$,
  $G \in \mathrm{Groups}(\mathrm{Groupoids}_\infty)$
  and  $(X,\rho) \in
  G\mathrm{Actions}
  \big(
    \mathrm{Groupoids}_\infty
  \big)$. Then the $G$-fixed locus (Def. \ref{FixedPoints})
  is given (due to Prop. \ref{BaseChangeFromBareInfinityGroupoids}) by
  $$
    X^G
    \;\simeq\;
    \mathbf{H}_{/ \ast \!\sslash\! G}
    \big(
      \ast \!\sslash\! G
      \,,\,
      X \!\sslash\! G
    \big)
    \;\;\;\;
    \in
    \mathrm{Groupoids}_\infty
    \,.
  $$
\end{example}

\begin{defn}[Pointed-automorphism group]
  \label{PointedAutomorphismGroup}
  Let $\mathbf{H}$ be an $\infty$-topos and
  $
    {\xymatrix{\ast \ar[r]^x & X }} \in \mathbf{H}^{\Delta^1}
  $
  a pointed object in
  $\mathbf{H}$. Then its \emph{pointed-automorphism} group
  $\mathrm{Aut}_\ast(X) \in \mathrm{Groups}(\mathbf{H})$
  is its automorphism group, according to Def. \ref{AutomorphismGroup},
  formed in the arrow $\infty$-topos $\mathbf{H}^{\Delta^1}$.
  This is characterized by a diagram in $\mathbf{H}$ of the form
 \vspace{-1mm}
  \begin{equation}
    \label{PointedAutomorphismAction}
    \raisebox{30pt}{
    \xymatrix@R=8pt@C=12pt{
      \;\ast\;
      \ar[dr]
      \ar[rr]
      \ar[dd]
      &&
      \ast \!\sslash\! \mathrm{Aut}_\ast(X)
      \ar[dr]
      \ar@{=}[dd]|-{ \phantom{AA} }
      \\
      &
      X
      \ar[rr]
      \ar[dd]
      &&
      X \sslash \mathrm{Aut}_\ast(X)
      \ar[dd]^-{ \rho_{\mathrm{Aut}_\ast} }
      \\
      \;\ast\;
      \ar[dr]
      \ar[rr]|-{\phantom{AA}}
      &&
      \mathbf{B} \mathrm{Aut}_\ast(X)
      \ar@{=}[dr]
      \\
      &
      \;\ast\;
      \ar[rr]
      &&
      \mathbf{B}\mathrm{Aut}_\ast(X)
    }
    }
  \end{equation}

   \vspace{-1mm}
\noindent
  where the front, rear, top and bottom squares are
  Cartesian: the bottom face trivially, the
  front face exhibiting the action on $X$,
  the top face exhibiting the given base point as a homotopy
  fixed point (Def. \ref{FixedPoints}) and the rear square
  exhibiting the trivial action on that point.
\end{defn}

\begin{defn}[Group-automorphism group]
  \label{GroupOFGroupAutomorphisms}
  Let $\mathbf{H}$ be an $\infty$-topos and
  $G \in \mathrm{Groups}(\mathbf{H})$ (Prop. \ref{LoopingAndDelooping}).
  Then the group of group-automorphisms of $G$
  is the group of pointed-automorphisms (Def. \ref{PointedAutomorphismGroup})
  of its delooping $\mathbf{B}G$ \eqref{LoopingDeloopingEquivalence}:
  \vspace{-2mm}
  $$
    \mathrm{Aut}_{\mathrm{Grp}}(G)
    \;:=\;
    \mathrm{Aut}_\ast
    (
      \mathbf{B}G
    )
    \;\;\;
    \in
    \mathrm{Groups}(\mathbf{H})
    \,.
  $$
\end{defn}

\begin{prop}[Canonical action of group-automorphism group]
  \label{CanonicalActionOfGroupAutomomorphisms}
  Let $\mathbf{H}$ be an $\infty$-topos and
  $G \in \mathrm{Groups}(\mathbf{H})$ (Prop. \ref{LoopingAndDelooping}).
  The
  group-automorphism group of $G$ (Def. \ref{GroupOFGroupAutomorphisms})
  has a canonical action (Prop. \ref{InfinityAction})
  \vspace{-1mm}
  $$
    (G, \rho_{\mathrm{Aut}_{\mathrm{Grp}}})
    \;\in\;
    \mathrm{Aut}_{\mathrm{Grp}}(G)\mathrm{Actions}(\mathbf{H})
  $$

  \vspace{-2mm}
\noindent
  on the underlying object $G \in \mathbf{H}$, which is
  such that

  \begin{itemize}
  \vspace{-.3cm}
  \item[{\bf (i)}] The neutral element $\xymatrix{\ast \ar[r]^e & G}$
  (Example \ref{NeutralElement})
  is a fixed point of the action (Def. \ref{FixedPoints}).

  \vspace{-.3cm}
  \item[{\bf (ii)}]
  Together with the defining action
  on the delooping $\mathbf{B}G$ of $G$ \eqref{LoopingDeloopingEquivalence},
  the looping equivalence \eqref{GroupIsLooping}
    $
    \xymatrix{
      G
      \ar[r]^-{ \simeq }
      &
      \Omega \mathbf{B}G
    }
  $
  is $\mathrm{Aut}_{\mathrm{Grp}}(G)$-equivariant (Def. \ref{Equivariance}).
  \end{itemize}
\end{prop}
\begin{proof}
  First consider item {\bf (ii)}:
  Write  $G \sslash \mathrm{Aut}_{\mathrm{Grp}}(G)$
  for the homotopy fiber product in the following
  pullback square
  \vspace{-2mm}
  \begin{equation}
    \label{HomotopyQuotientOfGByGroupAutomorphisms}
    \raisebox{20pt}{
    \xymatrix@R=1.5em{
      G \!\sslash\! \mathrm{Aut}_{\mathrm{Grp}}(G)
      \ar[rr]
      \ar[d]
      \ar@{}[drr]|-{ \mbox{\tiny\rm(pb)} }
      &&
      \ast \!\sslash\! \mathrm{Aut}_{\mathrm{Grp}}(G)
      \ar[d]
      \\
      \ast \!\sslash\! \mathrm{Aut}_{\mathrm{Grp}}(G)
      \ar[rr]
      &&
      (\mathbf{B}G) \!\sslash\! \mathrm{Aut}_{\mathrm{Grp}}(G)
    }
    }
  \end{equation}
  We need to show that this really is the homotopy quotient
  of the canonical group-automorphism action with the
  claimed property, in that it makes the
  total solid rear rectangle of the following
  diagram be Cartesian:
  \begin{equation}
    \label{TowardsUnderstandingGroupAutomorphismActionOnG}
    \raisebox{80pt}{
    \xymatrix@R=10pt@C=12pt{
      &
      \ast
      \ar[rr]
      \ar@{-->}[dd]_-{ e }
        &&
      \ast \sslash \mathrm{Aut}_{\mathrm{Grp}}(G)
      \ar@{-->}[dd]^-{ (\mathrm{id},\, \mathrm{id}) }
      \\
      \\
      &
      G
      \ar[dr]
      \ar[rr]
      \ar[dd]
      &&
      G \!\sslash\! \mathrm{Aut}_{\mathrm{Grp}}(G)
      \ar[dr]
      \ar[dd]|-{\phantom{AA}}
      \\
      &
      &
      \ast
      \ar[rr]
      \ar[dd]
      &&
      \ast \sslash \mathrm{Aut}_{\mathrm{Grp}}(G)
      \ar[dd]
      \\
      &
      \;\ast\;
      \ar[dr]
      \ar[rr]|-{\phantom{AA}}
      \ar[dd]
      &&
      \ast \!\sslash\! \mathrm{Aut}_{\mathrm{Grp}}(G)
      \ar[dr]
      \ar@{=}[dd]|-{ \phantom{AA} }
      \\
      &
      &
      \mathbf{B}G
      \ar[rr]
      \ar[dd]
      &&
      (\mathbf{B}G) \sslash \mathrm{Aut}_{\mathrm{Grp}}(\mathbf{B}G)
      \ar[dd]^-{ \rho_\ast }
      \\
      &
      \;\ast\;
      \ar[dr]
      \ar[rr]|-{\phantom{AA}}
      &&
      \mathbf{B} \mathrm{Aut}_{\mathrm{Grp}}(G)
      \ar@{=}[dr]
      \\
      &
      &
      \;\ast\;
      \ar[rr]
      &&
      \mathbf{B} \mathrm{Aut}_{\mathrm{Grp}}(G)
    }
    }
  \end{equation}

  \vspace{-.3cm}

  \noindent Here:
   \vspace{-.25cm}
  \begin{itemize}
   \item
   the bottom part is the diagram \eqref{PointedAutomorphismAction}
  (for $X = \mathbf{B}G$) which exhibits the pointed-automorphism action on
   $\mathbf{B}G$;
   \vspace{-.25cm}
   \item the top front square is Cartesian and exhibits the
    base point being a homotopy-fixed point;
    as in \eqref{PointedAutomorphismAction},
   \vspace{-.25cm}
   \item
   the top left square is Cartesian and exhibits the looping/delooping relation
  \eqref{GroupIsLooping};
  \vspace{-.25cm}
  \item the top right square is \eqref{HomotopyQuotientOfGByGroupAutomorphisms}
  and this Cartesian by definition.
\end{itemize}
\vspace{-.2cm}
Hence the solid top rear square and thus the total
  solid rear square are Cartesian, by the
  pasting law (Prop. \ref{PastingLaw}).

\vspace{1mm}
Finally to see item {\bf (i)}: Observe that there is the dashed morphism
shown in the top right of \eqref{TowardsUnderstandingGroupAutomorphismActionOnG},
this being the diagonal morphism induced from the Cartesian property of the
top right square, by the above. This means, by construction, that the total
vertical morphism on the right is an equivalence.
Now define the dashed top square to be a pullback square.
Then, by the pasting law (Prop. \ref{PastingLaw}), the pullback
object in the top left of the dashed square is equivalently
the pullback of the total rear diagram, hence the pullback of
an equivalence to a point, hence is itself equivalent to the point, as shown.
Since the point is terminal, the top left dashed morphism is
thus a cone over the Cartesian square on the top left.
By the universal property of the homotopy fiber product, this
means that the top left dashed morphism must be the neutral element
(Example \ref{NeutralElement}). The top dashed square hence exhibits
this as a homotopy fixed point.
\hfill \end{proof}

\begin{prop}[Group disivion is equivariant under group-automorphisms]
  \label{GroupDivisionIsEquivariantUnderGroupAutomorphisms}
  Let $\mathbf{H}$ be an $\infty$-topos
  and $G \in \mathrm{Groups}(\mathbf{H})$ (Prop. \ref{LoopingAndDelooping}).
  Then the group division morphism
  $\xymatrix{ G \times G \ar[rr]|-{\, (-) \cdot (-)^{-1}\, }  && G}$
  (Example \ref{GroupOperation})
  is equivariant (Def. \ref{Equivariance}) with respect
  to the canonical group-automorphism action
  (Prop. \ref{CanonicalActionOfGroupAutomomorphisms})
  of the group-automorphism group $\mathrm{Aut}_{\mathrm{Grp}}(G)$
  (Def. \ref{GroupOFGroupAutomorphisms}) acting on all three copies of
  $G$:
   \vspace{-3mm}
  $$
  \hspace{-4cm}
    \xymatrix{
      (G, \rho_{\mathrm{Aut}_{\mathrm{Grp}}})
      \times
      (G, \rho_{\mathrm{Aut}_{\mathrm{Grp}}})
      \ar[rr]^-{
        (-) \cdot (-)^{-1}
      }
      &&
      (G, \rho_{\mathrm{Aut}_{\mathrm{Grp}}})
    }
    \mathrlap{
      \;\;\;\;\;\;
      \in
      \mathrm{Aut}_{\mathrm{Grp}}(G)
      \mathrm{Actions}(\mathbf{H})\;.
    }
  $$
\end{prop}
\begin{proof}
  By \eqref{GroupOperationFromLooping} the
  group division morphism is a universal morphism induced
  from pasting of copies of the looping square \eqref{GroupIsLooping}.
  Thus the claim follows by Prop. \ref{CanonicalActionOfGroupAutomomorphisms}.
\hfill \end{proof}

\newpage
\section{Singular geometry}
 \label{SingularCohesiveGeometry}

Here we establish foundations of
a geometric homotopy theory of orbifolds
which unifies:
\begin{enumerate}[{\bf (i)}]
\vspace{-3mm}
\item
 \cref{DifferentialCohomology} --
 the cohesive geometric homotopy theory due to \cite{SSS09}\cite{Schreiber19},
 which reflects the {\it geometric aspects} of orbifolds;

\vspace{-3mm}
\item
  \cref{SingularGeometry} --
  the cohesive global-equivariant homotopy theory due to \cite{Rezk14},
  understood as reflecting the {\it singular aspects} of orbifolds,
  as in \hyperlink{FigureD}{Figure D}.
\end{enumerate}

\vspace{-5mm}
\subsection{Geometry}
\label{DifferentialCohomology}
\label{EtaleCohomology}

We present axioms internal to $\infty$-toposes for

\begin{itemize}
\vspace{-.3cm}
\item[-] \cref{DifferentialTopology} -- Differential topology
\vspace{-.2cm}
\item[-] \cref{DifferentialGeometry} -- Differential geometry
\vspace{-.2cm}
\item[-] \cref{SuperGeometry} -- Super-geometry
\end{itemize}

\vspace{-3mm}
\noindent
This is to provide, in \cref{OrbifoldGeometry} below,
a general abstract theory of geometric aspects of orbi-singular spaces
and of {\'e}tale $\infty$-stacks.

\subsubsection{Differential topology}
 \label{DifferentialTopology}

We present a formulation of differential topology
internal to $\infty$-toposes which we call \emph{cohesive} \cite{dcct}.
In $\infty$-categorical generalization of \cite{Lawvere94}\cite{Lawvere07},
this involves an abstract \emph{shape} operation $\mbox{\textesh}$ that
relates higher geometric spaces
to their bare underlying homotopy type.

\begin{defn}[Cohesive $\infty$-topos]
  \label{CohesiveTopos}
 {\bf (i)} An $\infty$-topos $\mathbf{H}$ (Def. \ref{InfinityTopos})
  is called \emph{cohesive} if its base geometric morphism
  (Prop. \ref{BaseGeometricMorphism}),
  to be denoted
  $\mathrm{Pnts}
    :
   \xymatrix@C=11pt{\mathbf{H} \ar[r] & \mathrm{Groupoids}_\infty}
  $,
  is part of an adjoint quadruple of $\infty$-functors (Def. \ref{AdjointInfinityFunctors})
  \vspace{-2mm}
  \begin{equation}
    \label{AdjunctionCohesion}
 \hspace{4cm}
    \xymatrix@R=1em{
      \mathbf{H} \;\;
        \ar@{->}@<+28pt>[rr]|<\times|-{
          \mathllap{
            \scalebox{.8}{
              \color{darkblue} \bf
              ``shape''
              \hspace{2.8cm}
            }
          }
          \;\mathrm{Shp}\;
        }_-{\raisebox{-6pt}{\tiny $\bot$}}
        \ar@{<-^{)}}@<+14pt>[rr]|-{
          \mathllap{
            \scalebox{.8}{
              \color{darkblue} \bf
              ``discrete''
              \hspace{2.8cm}
            }
          }
         \; \mathrm{Disc} \;
        }_-{\raisebox{-6pt}{\tiny $\bot$}}
        \ar@<+0pt>@<-0pt>[rr]|-{
          \mathllap{
            \scalebox{.8}{
              \color{darkblue} \bf
              ``points''
              \hspace{2.83cm}
            }
          }
         \; \mathrm{Pnts} \;
        }_-{\raisebox{-6pt}{\tiny $\bot$}}
        \ar@{<-^{)}}@<-14pt>[rr]|-{
          \mathllap{
            \scalebox{.8}{
              \color{darkblue} \bf
              ``chaotic''
              \hspace{2.8cm}
            }
          }
         \; \mathrm{Chtc} \;
        }
      &&
     \;\;  \mathbf{B}
      \\
      \mathclap{
        \mbox{
          \tiny
          \color{darkblue} \bf
          \begin{tabular}{c}
            cohesive
            \\
            $\infty$-topos
          \end{tabular}
        }
     \;\;\;\; }
      &&
      \mathclap{\;\;\;\;\;\;\;
        \mbox{
          \tiny
          \color{darkblue} \bf
          \begin{tabular}{c}
            discrete
            \\
            sub-topos
          \end{tabular}
        }
      }
    }
  \end{equation}

  \vspace{-3mm}
  \noindent such that
  {\bf (a)} $\mathrm{Disc}$ and $\mathrm{Chtc}$
  are fully faithful (Def. \ref{FullyFaithfulFunctor}),
  and {\bf (b)} such that $\mathrm{Shp}$ preserves finite products.

\vspace{-1.5mm}
  \item {\bf (ii)} We write
  \begin{equation}
    \label{CohesiveModalitiesFromAdjointQuadruple}
    \underset{
      \mbox{
        \small
        \color{darkblue} \bf
        ``shape''
      }
    }{
    \big(
      \;
      \raisebox{1pt}{\textesh} := \mathrm{Disc} \circ \mathrm{Shp}
      \;
    \big)
    }
    \dashv
    \underset{
      \mbox{
        \small
        \color{darkblue} \bf
        ``discrete''
      }
    }{
      \big(
        \;
        \flat := \mathrm{Disc} \circ \mathrm{Pnts}
        \;
      \big)
    }
    \dashv
    \underset{
      \mbox{
        \small
        \color{darkblue} \bf
        ``continuous''
      }
    }{
    \big(
      \;
      \sharp := \mathrm{Chtc} \circ \mathrm{Pnts}
      \;
    \big)
    }
    \;:\;
    \mathbf{H} \longrightarrow \mathbf{H}
  \end{equation}
  for the induced adjoint triple (Def. \ref{AdjointInfinityFunctors}) of modalities \eqref{AdjointModalities}
  (\emph{cohesive modalities}).
\end{defn}
\noindent The following direct consequence may serve to illustrate how
these axioms are put to work:

\vspace{-1mm}
\begin{prop}[Composite cohesive modalities]
  \label{CompositeCohesiveModalities}
  The cohesive modalities (Def. \ref{CohesiveTopos}) satisfy:
\vspace{-2mm}
  $$
    \mbox{\rm\textesh} \circ \flat
    \;\simeq\;
    \flat
    \phantom{AA}
    \mbox{and}
    \phantom{AA}
    \flat \circ \sharp
    \;\simeq\;
    \flat\;.
  $$
\end{prop}
\begin{proof}
  That $\mathrm{Disc}$ and $\mathrm{Chtc}$
  in \eqref{AdjunctionCohesion} are fully faithful
  means, equivalently (Prop. \ref{CharacterizationOfFullyFaithfulAdjoints}),
  that the co-unit morphisms \eqref{CounitOfAdjunction}
   \vspace{-2mm}
  $$
    \xymatrix{
      \mathrm{Shp} \circ \mathrm{Disc}
      \ar[r]^-{ \simeq }
      &
      \mathrm{id}
    }
    \;,
    \phantom{AAA}
    \xymatrix{
      \mathrm{Pnts} \circ \mathrm{Chtc}
      \ar[r]^-{ \simeq }
      &
      \mathrm{id}
    }
  $$

 \vspace{-1mm}
\noindent
  are natural equivalences.
 Hence  the image under
  $\mathrm{Disc} \circ (-) \circ \mathrm{Pnts}$
  of the first of these
  is a natural equivalence
  of the form
  $$
    \mbox{\textesh}
    \circ
    \flat
    \;=
    \xymatrix{
      \mathrm{Disc} \circ \mathrm{Shp}
      \circ
      \mathrm{Disc} \circ \mathrm{Pnts}
      \ar[r]^-{\simeq}
      &
      \mathrm{Disc} \circ \mathrm{Pnts}
    }
    =
    \;
    \flat
    \,.
  $$
  while the image of the second
  is of the form
   \vspace{-4mm}
   $$
    \flat \circ \sharp
    \;=
    \xymatrix{
      \mathrm{Disc} \circ \mathrm{Pnts}
      \circ
      \mathrm{Chtc} \circ \mathrm{Pnts}
      \ar[r]^-{ \simeq }
      &
      \mathrm{Disc} \circ \mathrm{Pnts}
    }
    =
    \;
    \flat
    \,.
  $$

  \vspace{-.9cm}

\hfill\end{proof}

\begin{lemma}[Only the empty object has empty shape]
  \label{OnlyTheEmptyObjectHasEmptyShape}
  Let $\mathbf{H}$ be a cohesive $\infty$-topos
  (Def. \ref{CohesiveTopos}).
  Then $X \in \mathbf{H}$
  is empty, i.e., equivalent to the initial object $\varnothing$ \eqref{InInfinityToposInitialObject},
  precisely if its shape
  \eqref{CohesiveModalitiesFromAdjointQuadruple} is empty:
  $$
    X \;\simeq\; \varnothing
    \phantom{AAA}
    \Leftrightarrow
    \phantom{AAA}
    \raisebox{1pt}{\rm\textesh}\,X
    \;\simeq\;
    \varnothing
    \;.
  $$
\end{lemma}
\begin{proof}
  In one direction, assume that $X \simeq \varnothing$.
  Noticing that $\varnothing$ is the initial colimit and that
  colimits are preserved by $\mbox{\textesh}$, this being a left adjoint
  (Prop. \ref{AdjointsPreserveCoLimits}),
  it follows that $\raisebox{1pt}{\textesh}(\varnothing)
  \;\simeq\; \varnothing$.

  In the other direction, assume that the shape of $X$ is empty.
  Then the shape unit \eqref{AdjunctionUnit} is a morphism of the form
  $$
    \xymatrix{
      X
      \ar[r]^-{\eta^{\scalebox{.6}{\textesh}}_X}
      &
      \raisebox{1pt}{\textesh}\, X
      \;\simeq\;
      \varnothing
    }
  $$
  and thus $X \simeq \varnothing$ follows
  as in \eqref{MorphismIntoInitialObjectImpliesDomainIsInitialObject},
  by universality of colimits (Example \ref{InitialObjectInInfinityToposIsEmpty}).
  \hfill
\end{proof}

\medskip

\noindent {\bf Cohesive $\infty$-group actions.}
The condition that $\mathrm{Shp}$ preserves finite products
implies the following properties.

\vspace{-1mm}
\begin{prop}[Shape preserves groups, actions and their homotopy quotients]
  \label{ShapePreservesHomotopyColimitsByDiscreteGroups}
  Let $\mathbf{H}$ be a cohesive $\infty$-topos (Def. \ref{CohesiveTopos}),
  $G \in \mathrm{Groups}\big( \mathbf{H}\big)$ \eqref{LoopingDeloopingEquivalence} and
  and $(X,\rho) \in G\mathrm{Actions}(\mathbf{H})$
  (Prop. \ref{InfinityAction}).

  \vspace{-1mm}
 \item {\bf (i)}  Then the shape
  $\raisebox{1pt}{\rm\textesh} X$
  \eqref{CohesiveModalitiesFromAdjointQuadruple}
  of $X$
  is equipped with an induced $\mbox{\rm\textesh} G$-action,
  such that
  the shape of the homotopy quotient \eqref{HomotopyQuotientAsColimit}
  is the homotopy quotient of the shapes.
  The analogous statement holds for $\flat$ \eqref{CohesiveModalitiesFromAdjointQuadruple}:
  \vspace{-1mm}
  $$
    \raisebox{1pt}{\rm\textesh}\,
    \big(
      X \!\sslash\! G
    \big)
    \;\simeq\;
    \big(
      \raisebox{1pt}{\rm\textesh}\, X
    \big)
    \!\sslash\!
    \big(
      \raisebox{1pt}{\rm\textesh}\, G
    \big)
    \phantom{AA}
    \mbox{and}
    \phantom{AA}
    \flat\,
    \big(
      X \!\sslash\! G
    \big)
    \;\simeq\;
    \big(
      \flat\, X
    \big)
    \!\sslash\!
    \big(
      \flat\, G
    \big)
    \,.
  $$

  \vspace{-2mm}
\noindent
  \item {\bf (ii)}  In particular, both $\mbox{\textesh}$ and $\flat$ preserve
  group objects and their deloopings (Prop. \ref{LoopingAndDelooping}):
  \vspace{-1mm}
  $$
    \raisebox{1pt}{\rm\textesh}\, \mathbf{B}G
    \;\simeq\;
    \mathbf{B} \,\raisebox{1pt}{\rm\textesh}\, G
    \phantom{AA}
    \mbox{and}
    \phantom{AA}
    \flat\, \mathbf{B}G
    \;\simeq\;
    \mathbf{B} \,\flat\, G
    \,.
  $$
\end{prop}
\begin{proof}
  The homotopy quotient of $X$ by $G$ is, equivalently, a colimit
  over a simplicial diagram of finite Cartesian products of copies of
  $X$ and $G$ \eqref{HomotopyQuotientAsColimit}.
  Hence the statement follows for every $\infty$-functor that
  commutes with simplicial colimits and with finite products.
   But, since $\mbox{\textesh}$ is a left adjoint, it commutes with
  all colimits (Prop. \ref{AdjointsPreserveCoLimits})
  and also with finite products,
  by assumption on $\mathrm{Shp}$ and since $\mathrm{Disc}$
  is a right adjoint.
  Similarly, $\flat$ is both left and right adjoint, and hence
  preserves all colimits and all limits (again Prop. \ref{AdjointsPreserveCoLimits}).
  That preservation of homotopy quotients implies
  preservation of $\infty$-groups follows by the delooping theorem
  (Prop. \ref{LoopingAndDelooping}).
\hfill \end{proof}

\begin{lemma}[Cohesive shape presverves homotopy fiber products]
  \label{CohesiveShapePresvervesHomotopFiberProductsOverDiscreteObjects}
  In a cohesive $\infty$-topos $\mathbf{H}$ (Def. \ref{CohesiveTopos}),
  the shape functor
  $\mathrm{Shp}$ \eqref{AdjunctionCohesion}
  preserves homotopy fiber products over
  cohesively discrete objects. That is,
  for $B \in \! \xymatrix{\mathbf{B} \; \ar@{^{(}->}[r]^-{ \scalebox{.6}{$\mathrm{Disc}$}}
  & \mathbf{H} }$ and $X, Y \in \mathbf{H}_{/B}$, we have a natural
  equivalence
   \vspace{-2mm}
  $$
    \mathrm{Shp}
    \big(
      X
      \underset{B}{\times}
      Y
    \big)
    \;\simeq\;
    \mathrm{Shp}(X)
    \underset{B}{\times}
    \mathrm{Shp}(Y)
    \,.
  $$
\end{lemma}
\begin{proof}
  This is proven in \cite[Thm. 3.8.19]{dcct} under
  the assumption that $\mathbf{H}$ admits an
  $\infty$-cohesive site of definition. This assumption
  was shown to be unnecessary in \cite[Lemma 3.10]{BP19}.
\hfill \end{proof}
\begin{lemma}[Shape of $\eta^{\scalebox{.6}{\textesh}}$-induced action]
  \label{ShapeOfAnInducedAction}
  Let $\mathbf{H}$ be a cohesive $\infty$-topos (Def. \ref{CohesiveTopos}),
  $G \in \mathrm{Groups}(\mathbf{H})$ (Prop. \ref{LoopingAndDelooping})
  and $(X,\rho) \in G\mathrm{Actions}(\mathbf{H})$ (Prop. \ref{InfinityAction}).

 \noindent  {\bf (i)} The left-induced action (Prop. \ref{PullbackAction})
  \vspace{-2mm}
  $$
 \hspace{1cm}
    (\widetilde X, \widetilde \rho)
    \;:=\;
    \mathbf{B}
    \big(
      \eta^{\scalebox{.7}{$\mbox{\rm\textesh}$}}_{G}
    \big)_!
    \,
    (X,\rho)
    \;\in\;
    (\raisebox{1pt}{\rm\textesh} G)\mathrm{Actions}(\mathbf{H})
  $$

   \vspace{-2mm}
\noindent
  along the shape unit morphism \eqref{AdjunctionUnit}
  $
    \xymatrix@C=3em{
      G
      \ar[r]|-{\;
        \eta_G^{\scalebox{.7}{\rm\textesh}}
     \; }
      &
      \raisebox{0pt}{\rm\textesh}\,G
    }
  $
  acts on an object
  whose shape \eqref{CohesiveModalitiesFromAdjointQuadruple}
  is that of $X$:
  \vspace{-2mm}
  $$
    \raisebox{1pt}{\rm\textesh}\, \widetilde X
    \;\simeq\;
    \raisebox{1pt}{\rm\textesh}\, X \;,
  $$

  \vspace{-2mm}
\noindent
  whence
  \begin{equation}
    \label{ShapedAction}
    \big(
      \raisebox{1pt}{\rm\textesh}X \;, \;
      \raisebox{1pt}{\rm\textesh}\rho
    \big)
    \;\in\;
    \big(
      \raisebox{1pt}{\rm\textesh}G
    \big)
    \mathrm{Actions}(\mathbf{H})\;.
  \end{equation}

\noindent  {\bf (ii)} Similarly, the restricted-induced action (Prop. \ref{PullbackAction})
  $$
    (\widetilde X, \widetilde \rho)
    \;:=\;
    \mathbf{B}
    \big(
      \raisebox{0pt}{\rm\textesh}
      \epsilon^\flat
    \big)^\ast
    \circ
    \mathbf{B}
    \big(
      \eta^{\scalebox{.7}{$\mbox{\rm\textesh}$}}_{G}
    \big)_!
    \,
    (X,\rho)
    \;\in\;
    (\flat G)\mathrm{Actions}(\mathbf{H})
  $$

   \vspace{-2mm}
\noindent
  along the pair of group homomorphisms (using Prop. \ref{ShapePreservesHomotopyColimitsByDiscreteGroups})
  $
    \xymatrix@C=3em{
      G
      \ar[r]|-{\;
        \eta_G^{\scalebox{.7}{\rm\textesh}}
      \;}
      &
      \raisebox{0pt}{\rm\textesh}\,G
      \ar@{<-}[r]|-{\;
        \scalebox{.7}{$\mbox{\rm\textesh}
        \,
        \epsilon_G^\flat$}
      \;}
      &
      \flat\,G
    }
  $
  acts on an object
  whose shape \eqref{CohesiveModalitiesFromAdjointQuadruple}
  is that of $X$:
  \vspace{-1mm}
  $$
    \raisebox{1pt}{\rm\textesh}\, \widetilde X
    \;\simeq\;
    \raisebox{1pt}{\rm\textesh}\, X \;.
  $$
\end{lemma}
\begin{proof}
By Prop. \ref{InfinityAction} and Prop. \ref{PullbackAction},
the object $\widetilde X$ sits in a diagram
of Cartesian squares (Notation \ref{CartesianSquares})
as shown on the left in the following
(the full square in case {\bf (i)}, the pasting decomposition for case {\bf (ii)}):
\begin{equation}
  \label{ShapeOfInducedAction}
  \hspace{-.8cm}
  \raisebox{40pt}{
  \xymatrix@R=1.5em@C=2.8em{
    \widetilde X
    \ar[dd]
    \ar[r]
    \ar@{}[rdd]|-{
      \mbox{\tiny(pb)}
    }
    &
    \widetilde X
      \!\sslash\!
    \flat G
    \ar[dd]
    \ar[r]
     \ar@{}[ddr]|-{
      \mbox{\tiny(pb)}
     }
    &
    X \!\sslash\! G
    \ar[d]^-{\rho }
    &
    &
    \raisebox{1pt}{\textesh}\, \widetilde X
    \ar[dd]
    \ar@{}[ddr]|-{
      \mbox{\tiny(pb)}
    }
    \ar[r]
    &
    \big(
      \raisebox{1pt}{\textesh}\, \widetilde X
    \big)
    \!\sslash\!
    \big(
      \flat G
    \big)
    \ar[r]
    \ar[dd]
    \ar@{}[ddr]|-{
      \mbox{\tiny(pb)}
    }
    &
    \big(
      \raisebox{1pt}{\textesh}\,
      X
    \big)
    \!\sslash\!
    \big(
      \raisebox{1pt}{\textesh}\,
      G
    \big)
    \ar[dd]^-{
      \mbox{\textesh}\,\rho
    }
    \\
    &
    &
    \mathbf{B}G
    \ar[d]^-{
      \mathbf{B}
      \,
      \eta_G^{\scalebox{.5}{$\mbox{\rm\textesh}$}}
    }
    &
    {\phantom{AA}}
    \overset{
      \mbox{\textesh}
    }{\longmapsto}
    {\phantom{AA}}
    &
    \\
    \ast
    \ar[r]
    &
    \mathbf{B}\,\flat\,G
    \ar[r]_-{
      \scalebox{.7}{$
        \mathbf{B}
        \,\mbox{\rm\textesh}\,
        \epsilon_G^{\flat}
      $}
    }
    &
    \mathbf{B}
    \,\raisebox{1pt}{\rm\textesh}\,
    G
    &&
    \ast
    \ar[r]
    &
    \mathbf{B}\,\flat\,G
    \ar[r]_-{
      \scalebox{.7}{$
        \mathbf{B}
        \,\mbox{\rm\textesh}\,
        \epsilon_G^{\flat}
      $}
    }
    &
    \mathbf{B}
    \,\raisebox{1pt}{\rm\textesh}\,
    G
  }
  }
\end{equation}
But, since the objects in the bottom row
$
  \mathbf{B}
  \,\raisebox{1pt}{\textesh}\,
  G
  \;\simeq\;
  \raisebox{1pt}{\textesh}\,
  \mathbf{B}G
$
and
$
  \mathbf{B}
  \,\flat\,
  G
  \;\simeq\;
  \flat\,
  \mathbf{B}G
$
(equivalences by Prop. \ref{ShapePreservesHomotopyColimitsByDiscreteGroups})
are both cohesively discrete, Lemma
\ref{CohesiveShapePresvervesHomotopFiberProductsOverDiscreteObjects}
says that the image of these squares under shape are still Cartesian.
This is shown on the right in \eqref{ShapeOfInducedAction}, where we have identified the
shape of the various objects by using Prop. \ref{ShapePreservesHomotopyColimitsByDiscreteGroups}
and idempotency of the modality (Prop. \ref{IdempotentMonads}).
With this, the pasting law (Prop. \ref{PastingLaw}) implies
that the outer right square in \eqref{ShapeOfInducedAction} is itself Cartesian, hence that
$\raisebox{1pt}{\textesh}\widetilde X$ is the homotopy fiber of
$\mbox{\textesh}\rho$. This implies the claim, by Prop. \ref{InfinityAction}.
\hfill \end{proof}

\begin{prop}[Automorphisms along shape-unit]
  \label{AutomorphismsAlongShapeUnit}
  Let $\mathbf{H}$ be a cohesive $\infty$-topos (Def. \ref{CohesiveTopos}),
  $G \in \mathrm{Groups}(\mathbf{H})$ (Prop. \ref{LoopingAndDelooping})
  and $(X,\rho) \in G\mathrm{Actions}(\mathbf{H})$ (Prop. \ref{InfinityAction}).
    There is a canonical homomorphism
     \vspace{-2mm}
  \begin{equation}
    \label{AutomorphismsAlongShapeUnit}
    \xymatrix{
      \mathrm{Aut}(X)
      \ar[rr]^{
        \mathrm{Aut}
        \big(
          \eta^{\scalebox{.7}{\rm\textesh}}_X
        \big) }
      &&
      \mathrm{Aut}(\raisebox{1pt}{\rm\textesh}X)
    }
  \end{equation}

   \vspace{-2mm}
\noindent
  from the automorphism group (Def. \ref{AutomorphismGroup})
  of $X$ to that of the shape \eqref{CohesiveModalitiesFromAdjointQuadruple} of $X$, which is such that
  the shape unit $\eta^{\scalebox{.7}{\rm\textesh}}_X$ \eqref{AdjunctionUnit}
  is equivariant (Def. \ref{Equivariance})
  with respect to the canonical automorphism action \eqref{AutomorphismAction}
  on $X$ and the
  restriction (Prop. \ref{PullbackAction})
  along this morphism \eqref{AutomorphismsAlongShapeUnit}
  of the canonical automorphism action on $\raisebox{1pt}{\rm\textesh}X$:
   \vspace{-2mm}
  $$
    \xymatrix{
      (X,\rho_{\mathrm{Aut}(X)})
      \ar[rr]^-{
        \eta^{\scalebox{.7}{\rm\textesh}}_X
      }
      &&
      \mathrm{Aut}(\eta^{\scalebox{.7}{\rm\textesh}}_X)^\ast
      \big(
        \raisebox{1pt}{\rm\textesh}X,
        \rho_{
          \scalebox{.7}{$
           \mathrm{Aut}(\raisebox{1pt}{\rm\textesh}X)
          $}
        }
      \big)
      \mathrlap{
        \;\;\;\;\;\;
        \in
        \mathrm{Aut}(X)\mathrm{Actions}(\mathbf{H})
        \,.
      }
    }
  $$
\end{prop}
\begin{proof}
  Take the morphism \eqref{AutomorphismsAlongShapeUnit} to be the
  composite
   \vspace{-2mm}
\noindent
  $$
    \xymatrix@R=7pt@C=40pt{
      \mathrm{Aut}(X)
      \ar[rr]^-{
        \mathrm{Aut}\big(\eta^{\scalebox{.7}{\textesh}}_X\big)
      }
      \ar[dr]_-{
        \eta^{\scalebox{.5}{\textesh}}_{\mathrm{Aut}(X)}
      }
      &&
      \mathrm{Aut}
      \big(
        \raisebox{1pt}{\textesh}X
      \big)
      \\
      &
      \raisebox{1pt}{\textesh}
      \big(\mathrm{Aut}(X)\big)
      \ar[ur]_{ \;\;
        \scalebox{.6}{$
          \Omega \vdash \raisebox{0pt}{\textesh} \rho_{\mathrm{Aut}}
        $}
      }
    }
  $$

  \vspace{-2mm}
\noindent   where {\bf (a)} the left morphism is the shape unit \eqref{AdjunctionUnit},
  using Prop. \ref{ShapePreservesHomotopyColimitsByDiscreteGroups},
  while {\bf (b)} the right morphism is that which exhibits, via
  Prop. \ref{AutomorphimsGroupIsUniversal},
  the $\raisebox{1pt}{\textesh}\mathrm{Aut}(X)$-action
  $\raisebox{1pt}{\textesh}\rho_{\mathrm{Aut}}$
  \eqref{ShapedAction}
  on $\raisebox{1pt}{\textesh}X$ from Lemma \ref{ShapeOfAnInducedAction}.
    Then consider the following diagram of homotopy fiber sequences:
     \vspace{-2mm}
  $$
  \hspace{-.5cm}
    \xymatrix@R=1em@C=5em{
      &&
      \raisebox{1pt}{\textesh}X
      \ar[rr]
      &&
      (\raisebox{1pt}{\textesh}X)
        \!\sslash\!
      \mathrm{Aut}
        \big(
          \raisebox{1pt}{\textesh}X
        \big)
      \ar[dd]^-{
        \scalebox{.7}{$
          \rho_{
            \scalebox{.7}{$
              \mathrm{Aut}(\raisebox{1pt}{\textesh}(X))
            $}
          }
        $}
      }
      \\
      &
      \;\raisebox{1pt}{\textesh}X\;
      \ar[rr]
      \ar@{=}[ur]
      &&
      (\raisebox{1pt}{\textesh}X)
      \!\sslash\!
      \big(
        \raisebox{1pt}{\textesh}\mathrm{Aut}(X)
      \big)
      \ar[dd]|-{
        \scalebox{.7}{\textesh}
        \rho_{\mathrm{Aut}(X)}
      }
      \ar[ur]
      \ar@{}[dr]|-{\mbox{\tiny(pb)}}
      \\
      \;\;\;X\;\;\;
      \ar[rr]
      \ar[ur]|-{\;
        \scalebox{.7}{$
          \eta^{\scalebox{.7}{\textesh}}_X
        $}
      }
      &&
      X \!\sslash\! \mathrm{Aut}(X)
      \ar[dd]^-{
        \rho_{\mathrm{Aut}(X)}
      }
      \ar[ur]|-{\;
        \scalebox{.7}{$
          \eta^{\scalebox{.7}{\textesh}}_{X \sslash \mathrm{Aut}(X)}
        $}
      }
      &&
      \mathbf{B}
      \mathrm{Aut}
      \big(
        \raisebox{1pt}{\textesh}X
      \big)
      \\
      &
      &&
      \mathbf{B} \raisebox{1pt}{\textesh}\mathrm{Aut}(X)
      \ar[ur]|-{ \;\;
        \scalebox{.7}{$
          \vdash
          \raisebox{0pt}{\textesh}
          \rho_{\mathrm{Aut}}
        $}
      }
      \\
      &&
      \mathbf{B} \mathrm{Aut}(X)
      \ar[ur]|-{\;
        \scalebox{.7}{$
          \eta^{\scalebox{.7}{\textesh}}_{\mathbf{B}\mathrm{Aut}(X)}
        $}
      }
      \ar@/_2pc/[uurr]|-{\;\;
        \scalebox{.7}{$
          \mathrm{Aut}
          \big(
            \eta^{\scalebox{.7}{\textesh}}_X
          \big)
        $}
      }
    }
  $$
  Here {\bf (i)} the fiber squence in the middle is that
  from the right of \eqref{ShapeOfInducedAction},
  {\bf (ii)}
  the right part is the defining pullback from Prop. \ref{AutomorphimsGroupIsUniversal},
  while {\bf (iii)} the left part exists by the naturality of
  $\eta^{\scalebox{.7}{\textesh}}$. By the commutativity of the
  total front square it factors through the coresponding pullback square,
  thus implying the claim.
\hfill \end{proof}

\medskip

\noindent {\bf Concrete cohesive objects.}
\begin{defn}[Concrete objects]
  \label{ConcreteObjects}
  Let $\mathbf{H}$ be a cohesive $\infty$-topos (Def. \ref{CohesiveTopos}).

  \vspace{-1mm}
\item {\bf (i)}   For $X \in \mathbf{H}_0 \hookrightarrow \mathbf{H}$
  0-truncated (Def. \ref{nTruncatedObjects}),
  we say that $X$ is a \emph{concrete object}
  or \emph{concrete cohesive space} if the unit
  $\eta^\sharp_X$ \eqref{AdjunctionUnit}  of the
  $\sharp$-modality \eqref{CohesiveModalitiesFromAdjointQuadruple}
  is (-1)-truncated (Def. \ref{nTruncatedMorphism}),
  hence a monomorphism. By the $0$-image factorization \eqref{nImageFactorization},
  $$
    \xymatrix{
      X
      \ar@{->>}[rr]^-{\mbox{\tiny (-1)-conn.}}
      \ar@/_1pc/[rrrr]_-{
        \underset{
          \mathclap{
          \mbox{
            \tiny
            \color{darkblue}
            \bf
            unit morphism of $\sharp$-modality
            }
          }
        }{
          \eta^\sharp_X
        }
      }
      &&
      \overset{
        \mathclap{
        \mbox{
          \tiny
          \color{darkblue}
          \bf
         \!\!\!  image factorization
        }
        }
      }{
       \;\;\;\;\; \sharp_1 X \;\;\;\;\;\;\;\;\;\;
      }
      \ar@{^{(}->}[rr]^-{\mbox{\tiny (-1)-trunc.}}
      &&
     \; \sharp X
    }
  $$
  this means equivalently that $X$ is equivalent to its
  0-image under the $\sharp$-unit \eqref{AdjunctionUnit}:
  \vspace{-3mm}
  \begin{equation}
    \label{ConcreteObjectCondition}
    X \in \mathbf{H}_{0}
    \;:\;
    \;\;\;\;\;\;\;\;
    \mbox{
      $X$
      is concrete
    }
    \phantom{AAA}
    \Leftrightarrow
    \phantom{AAA}
    \xymatrix{
      X
      \; \ar@{^{(}->}[r]^-{ \eta^\sharp_X }
      &
      \sharp X
    }
    \phantom{AAA}
    \Leftrightarrow
    \phantom{AAA}
    \sharp_1 X
    \;\simeq\;
    X\;.
  \end{equation}

   \vspace{-4mm}
\item {\bf (ii)}   We write
  \begin{equation}
    \label{CategoryOfConcrete0TruncatedObjects}
    \mathbf{H}_{0,\sharp_1}
    \longhookrightarrow
    \mathbf{H}_0
    \longhookrightarrow
    \mathbf{H}
  \end{equation}
  for the full subcategory of the 0-truncated objects
  on those which are concrete.

  \vspace{-1.5mm}
 \item {\bf (iii)}  Moreover, for $n \in \mathbb{N}$ we define, recursively,
  full sub-$\infty$-categories of concrete $(n+1)$-truncated objects
  (Def. \ref{nTruncatedObjects})
  \begin{equation}
    \label{CategoryOfConcreteNTruncatedObjects}
    \mathbf{H}_{n+1,\sharp_1}
      \longhookrightarrow
    \mathbf{H}_{n+1}
      \longhookrightarrow
    \mathbf{H}
  \end{equation}
  by declaring that $X \in \mathbf{H}_{n+1}$ is \emph{concrete}
  if:

  \begin{itemize}
  \vspace{-4mm}
\item it admits a \emph{concrete atlas},
  namely an effective epimorphism out of a
  concrete 0-truncated object \eqref{ConcreteObjectCondition},

    \vspace{-3mm}
\item  such that
 the homotopy fiber product of the atlas
  with itself (which is an $n$-truncated object) is a concrete:
   \vspace{-2mm}
  \begin{equation}
    \label{ConcreteNTruncatedObject}
    X \in \mathbf{H}_{n+1}
    \;:\;
    \;\;\;\;
    \mbox{$X$ is concrete}
    \;\;\;\;\;\;\;
    \Leftrightarrow
    \;\;\;\;\;\;\;
    \underset{
      \mathclap{
       X_0 \in \mathbf{H}_{0, \sharp_1}
       }
    }{\exists}
    \;\;\;\;\;
    :
    \;\;
    \xymatrix{
      X_0
      \ar@{->>}[rr]^-{
        \mbox{
          \tiny
          $(-1)$-trunc.
        }
      }
      &&
      X
    }
    \;\;\;
    \mbox{and}
    \;\;\;
    X_0 \underset{X}{\times} X_0
    \;\in\;
    \mathbf{H}_{n,\sharp_1} \;.
  \end{equation}
\end{itemize}
\end{defn}

\noindent {\bf Cohesive charts.}
\begin{defn}[Charts]
  \label{ChartsForCohesion}
  Let $\mathbf{H}$ be a cohesive $\infty$-topos
  (Def. \ref{CohesiveTopos}). We say that
  an \emph{$\infty$-category of cohesive charts} for $\mathbf{H}$ is
  an $\infty$-site $\mathrm{Charts}$ for $\mathbf{H}$
  (Prop. \ref{ToposLexReflection})
  \vspace{-3mm}
  $$
    \xymatrix{
      \mathbf{H} \;
      \ar@{<-}@<+8pt>[rr]^-{L}
      \ar@<-8pt>@{^{(}->}[rr]^-{
        \raisebox{2pt}{
          \scalebox{.9}
          {$
            \bot
          $}
        }
      }
      &&
     \; \mathrm{PreSheaves}_\infty
      (
        \mathrm{Charts}
      )
    }
  $$
  \vspace{-2mm}
  all of whose objects
  (under the $\infty$-Yoneda embedding $y$, Prop. \ref{InfinityYonedaEmbedding})
  have contractible shape \eqref{CohesiveModalitiesFromAdjointQuadruple}:
  \vspace{-1mm}
  \begin{equation}
    \label{GeometricallyContractibleGenerators}
    \hspace{-1cm}
    \raisebox{6pt}{
    \xymatrix@R=-1pt{
      \mathrm{Charts}
     \; \ar@{^{(}->}[r]^-{ y }
      &
     \; \mathbf{H} \;
      \ar[r]^-{ \mathrm{Shp} }
      &
      \mathrm{Groupoids}_\infty
      \\
      U
      \ar@{|->}[r]
      &
     \; U \;
      \ar@{|->}[r]
      &
      \mathrm{Shp}(U)
      \;
      \mathrlap{
        \simeq \ast
      }
    }
    }
    \phantom{AAAAA}
    \Leftrightarrow
    \phantom{AAAAA}
    \raisebox{6pt}{
    \xymatrix@R=-1pt{
      \mathrm{Charts}
    \;  \ar@{^{(}->}[r]^-{ y }
      &
      \; \mathbf{H} \;
      \ar[r]^-{ \raisebox{1pt}{\textesh} }
      &
      \mathbf{H}
      \\
      U
      \ar@{|->}[r]
      &
     \;  U
     \; \ar@{|->}[r]
      &
      \raisebox{1pt}{\textesh}(U)
      \;
      \mathrlap{
        \simeq \ast
      }
    }
    }
  \end{equation}
\end{defn}

\begin{lemma}[Charts are cohesively connected]
  \label{ChartsAreCohesivelyConnected}
  Let $\mathbf{H}$ be a cohesive $\infty$-topos (Def. \ref{CohesiveTopos})
  with a site of $\mathrm{Charts}$ (Def. \ref{ChartsForCohesion}).
  Then, for $U \in \mathrm{Charts}$ and
  $\big\{ X_i \in \mathbf{H}\big\}_{i \in I}$
  an indexed set of objects of $\mathbf{H}$,
  we have that every morphism from $U$ into the coproduct of the
  $X_i$ factors through one of the $X_i$:
\vspace{-3mm}
  $$
    \xymatrix{
      U \ar[r]^-{f}
      &
      \underset{i \in I}{\sqcup}
      X_i
    }
    \phantom{AAAA}
    \Leftrightarrow
    \phantom{AAAA}
    \underset{i_0 \in I}{\exists}
    \;\;
    \xymatrix@C=4em@R=1.5em{
      U
      \ar@/_1pc/[rr]_-{f}
      \ar@{-->}[r]
      &
      X_{i_0}
      \ar[r]^-{ q_{X_{i_0}} }
      &
      \underset{i \in I}{\sqcup}
      X_i
      \,.
    }
  $$
\end{lemma}
\begin{proof}
  Consider the pullbacks $\xymatrix@C=13pt{U_i \ar[r]^{q_{U_i}} & U}$
  along $f$ of the canonical inclusions of  the $X_i$ into their coproduct,
  given by these Cartesian squares (Notation \ref{CartesianSquares}):
  \vspace{-2mm}
  \begin{equation}
    \label{ConnectedComponentsOfPreimageOfMapOutOfChart}
    \raisebox{20pt}{
    \xymatrix@R=1.2em{
      U_i
      \ar[d]_-{ q_{U_i} }
      \ar[rr]
      \ar@{}[drr]|-{\mbox{\tiny(pb)}}
      &&
      X_i
      \ar[d]^-{ q_{X_i} }
      \\
      U
      \ar[rr]_-{f}
      &&
      \underset{i \in I}{\bigsqcup}\;
      X_i
    }
    }
  \end{equation}

  \newpage

  \noindent By Prop. \ref{ColimitsOfEquifiberedTransformations},
  this is such that
  \begin{equation}
    \label{ChartAsDisjointUnionOfConnectedPreimages}
    U \;\simeq\; \underset{i \in I}{\bigsqcup}\; U_i\;.
  \end{equation}
  The image of \eqref{ChartAsDisjointUnionOfConnectedPreimages} under
  shape \eqref{CohesiveModalitiesFromAdjointQuadruple}
  is
  \vspace{-1mm}
  $$
    \ast
    \;\simeq\;
    \raisebox{1pt}{\textesh}\,U
    \;\simeq\;
    \underset{i \in I}{\bigsqcup}\; \raisebox{1pt}{\textesh}\, U_i
    \;\;\;\;\;\;\;\;
    \in
    \xymatrix{
      \mathrm{Groupoids}_\infty
      \; \ar@{^{(}->}[r]^-{ \mathrm{Disc} }
      &
      \mathbf{H}
    }
    \,,
  $$

  \vspace{-2mm}
\noindent
  where on the left we used the defining property \eqref{GeometricallyContractibleGenerators}
  of charts and on the right we used that the shape operation,
  being a left adjoint, preserves coproducts (Prop. \ref{AdjointsPreserveCoLimits}).
  But, since $\ast \in \mathrm{Groupoids}_\infty$ is connected,
  this implies that there is $i_0 \in I$ with
  $$
    \raisebox{0.5pt}{\textesh}\, U_i
    \;\simeq\;
    \left\{
    \!
    \begin{array}{lcl}
      \varnothing &\vert& i \neq i_0
      \\
      \ast &\vert& i = i_0
    \end{array}
    \right.
  $$
  From this, Lemma \ref{OnlyTheEmptyObjectHasEmptyShape} implies that
  $
    U_i
    \;\simeq\;
    \varnothing
       \mbox{
      for
      $i \neq i_0$
    }
  $
  and, with \eqref{ChartAsDisjointUnionOfConnectedPreimages}, this implies
  \vspace{-1mm}
  $$
    \xymatrix{
      U_{i_0}
      \ar[r]_-{\simeq}^-{ q_{U_{i_0}} }
      &
      U
    }
    \,.
  $$

    \vspace{-1mm}
\noindent  Using this in \eqref{ConnectedComponentsOfPreimageOfMapOutOfChart}
  gives the desired factorization.
\hfill \end{proof}

\begin{lemma}[Quotient by cohesively discrete $\infty$-group]
  Let $\mathbf{H}$ be a cohesive $\infty$-topos (Def. \ref{CohesiveTopos})
  which admits a site of $\mathrm{Charts}$ (Def. \ref{ChartsForCohesion}).
  Then, for
  \vspace{-2mm}
  \begin{equation}
    \label{AGeometricallyDiscreteGroup}
    G
    \;\in
    \xymatrix{
      \mathrm{Groups}(\mathrm{Groupoids}_\infty)
     \; \ar@{^{(}->}[r]^-{ \mathrm{Disc} }
      &
      \mathrm{Groups}(\mathbf{H})
    }
  \end{equation}
  a cohesively discrete $\infty$-group \eqref{LoopingDeloopingEquivalence}
  and $U \in \mathrm{Charts}$, we have an equivalence
   \vspace{-1mm}
  \begin{equation}
    \label{HomOfChartIntoDiscreteDeloopingGroupoid}
    \mathbf{H}(U, \ast \sslash G)
    \;\simeq\;
    \ast \sslash G
    \;\;\;
    \in
    \;
    \mathrm{Groupoids}_\infty
    \,.
  \end{equation}
\end{lemma}
\begin{proof}
  Since $\mathrm{Disc}$ is both a left and a right
  adjoint, it preserves (Prop. \ref{AdjointsPreserveCoLimits})
  the homotopy quotient  that
  corresponds to the effective epimorphism
  $\xymatrix{ \ast \ar@{->>}[r] & \ast \!\sslash\! G }$
  (Prop. \ref{GroupoidsEquivalentToEffectiveEpimorphisms})
  so that
   \vspace{-2mm}
  $$
    \ast \!\sslash\! G
    \;\in\;
    \xymatrix{
      \mathrm{Groupoids}_\infty
      \;\ar@{^{(}->}[rr]^-{ \mathrm{Disc} }
      &&
      \mathbf{H}
    }
  $$
  is a cohesively discrete object.
  With this, we have the following sequence of natural
  equivalences:
  $$
    \mathbf{H}
    \big(
      U,
      \,
      \ast \!\sslash\! G
    \big)
    \;\simeq\;
    \mathbf{H}
    \big(
      U,
      \,
      \mathrm{Disc}(\ast \!\sslash\! G)
    \big)
    \;\simeq\;
    \mathrm{Groupoids}_\infty
    \big(
      \mathrm{Shp}(U),
      \,
      \ast \!\sslash\! G
    \big)
    \;\simeq\;
    \mathrm{Groupoids}_\infty
    \big(
      \ast,
      \,
      \ast \!\sslash\! G
    \big)
    \;\simeq\;
    \ast \!\sslash\! G
    \,
  $$
  where the second step is the
  hom-equivalence \eqref{AdjunctionHomEquivalence} of the
  $\mathrm{Shp} \dashv \mathrm{Disc}$-adjunction and the
  third step is the condition that the chart $U$ has
  contractible shape.
\hfill \end{proof}

\begin{lemma}[Homming Charts into quotients by discrete groups]
  \label{HommingChartsIntoHomotopyFiberSequences}
  Let $\mathbf{H}$ be a cohesive $\infty$-topos
  (Def. \ref{CohesiveTopos}) which admits $\mathrm{Charts}$
  (Def. \ref{ChartsForCohesion}).
  Then, for $X \in \mathbf{H}$ an object equipped with an
  $\infty$-action (Prop. \ref{InfinityActionHomotopyFiberSequence})
  by a geometrically
  discrete $\infty$-group $G$
  \eqref{AGeometricallyDiscreteGroup},
  the homotopy quotient $X \!\sslash\! G$
  \eqref{HomotopyQuotientAsColimit}
  is given as an $\infty$-sheaf on $\mathrm{Charts}$, by
  assigning to $U \in \mathrm{Charts}$ the homotopy quotient
  of the $\infty$-groupoid of $U$-shapes plots of $X$:
   \vspace{-1mm}
  $$
    X \!\sslash\! G
    \;:\;
    U \;\longmapsto\;
    \mathbf{H}(U,X) \!\sslash\! G\;.
  $$
\end{lemma}
\begin{proof}
  Consider the image of the homotopy fiber sequence
  that characterizes the given $\infty$-action (Prop. \ref{InfinityAction})
  under homming the chart $U$ into it:
   \vspace{-2mm}
  \begin{equation}
    \label{HommingAChartIntoAHomotopyFiberSequence}
    \hspace{-3cm}
    \raisebox{20pt}{
    \xymatrix{
      X \ar[rr]^-{\mathrm{fib}(p)}
      &&
      X \!\sslash\! G
      \ar[d]^-{ p }
      \\
      &&
      \ast \!\sslash\! G
    }
    }
    \phantom{AAAA}
    \overset{
      \mathbf{H}(U,-)
    }{\longmapsto}
    \phantom{AAAA}
    \raisebox{20pt}{
    \xymatrix{
      \mathbf{H}(U,X)
      \ar[rr]^-{
        \mathrm{fib}
        \left(
          \mathbf{H}(U,p)
        \right)
      }
      &&
      \mathbf{H}(U,X) \!\sslash\! G
      \mathrlap{
         \;\;\simeq\;
         \mathbf{H}\big(
           U,
           X \!\sslash\! G
         \big)
      }
      \ar[d]^-{ \mathbf{H}(U,p) }
      \\
      &&
      \ast \!\sslash\! G
      \mathrlap{
        \;\;\simeq\;
        \mathbf{H}(U, \ast \!\sslash\! G)
      }
    }
    }
  \end{equation}

   \vspace{-2mm}
\noindent
  Since the hom-functor $\mathbf{H}(U,-)$ preserves
  limits, the result is again a homotopy fiber sequence,
  as shown on the right of \eqref{HommingAChartIntoAHomotopyFiberSequence}.
  Moreover, by the assumption that $G$ is geometrically discrete
  and that $U$ is geometrically contractible,
  we have the equivalence \eqref{HommingAChartIntoAHomotopyFiberSequence}
  shown on the bottom right.
  This means that the fiber sequence on the right of
  \eqref{HommingAChartIntoAHomotopyFiberSequence}
  exhibits $\mathbf{H}(U, X \!\sslash\! G )$
  as the homotopy quotient
  $\mathbf{H}(U,X) \!\sslash\! G$
  of an $\infty$-action by $G$ on
  $\mathbf{H}(U,X)$.
\hfill \end{proof}

\begin{lemma}[Fixed locus in 0-truncated objects for discrete groups]
  \label{FixedLocusIn0TruncatedObjectsForDiscreteGroups}
  Let $\mathbf{H}$ be a cohesive $\infty$-topos (Def. \ref{CohesiveTopos})
  with a site of $\mathrm{Charts}$ (Def. \ref{ChartsForCohesion}).
  Let $G \in \mathrm{Groups}(\mathbf{H})$ (Prop. \ref{LoopingAndDelooping})
  be discrete
  $G \simeq \flat G$ and 0-truncated, $G \simeq \tau_0 G$,
  and let $(X,\rho) \in G \mathrm{Actions}(\mathbf{H})$
  (Prop. \ref{InfinityAction}) with $X \simeq \tau_0 X$ also
  0-truncated.
  Then the $G$-fixed locus $X^G \in \mathbf{H}$ (Def. \ref{FixedPoints})
  is itself 0-truncated and
  such that, for $U \in \mathrm{Charts}$, we have a natural
  equivalence
  \begin{equation}
    \label{AsPresheafFixedLocusOf0TruncatedByDiscreteGroup}
    \mathbf{H}
    \big(
      U,
      X^G
    \big)
    \;\simeq\;
    \mathbf{H}(U,X)^G
    \;:=\;
    \Big\{
      \phi \in \mathbf{H}(U,X)
      \;\vert\;
      \underset{g \in G}{\forall}
      \;
      g\cdot \phi = \phi
     \Big\}
  \end{equation}

   \vspace{-2mm}
\noindent
  between {\bf (a)} the hom-set from $U$ to $X^G$
  and {\bf (b)} the naive set of fixed points in the hom-set from $U$ to
  $X$, with respect to the
  restriction (Prop. \ref{PullbackAction})
  along $K \hookrightarrow G$ of the
  induced $G$-action \eqref{HommingAChartIntoAHomotopyFiberSequence}
  on the latter.
\end{lemma}
\begin{proof}
  We claim that we have the following sequence of natural equivalences:
  \begin{equation}
    \label{TowardsUnderstandingTheFixedLocusInA0TruncatedObject}
    \begin{aligned}
      \mathbf{H}(U, X^G)
      & =
      \mathbf{H}\left(
        U,
        \mathbf{B}(G \to \ast)_\ast
        \big(
          (X,\rho)
        \big)
      \right)
      \\
      & \simeq
      \mathbf{H}_{/\mathbf{B}G}
      \left(
        \mathbf{B}(G \to \ast)_\ast
        \big(
          U
        \big),
        X \!\sslash\! G
      \right)
      \\
      & \simeq
      \mathbf{H}_{/\mathbf{B}G}
      \big(
        (\ast \!\sslash\! G)
        \times
        U
        \,,\,
        X \!\sslash\! G
      \big)
      \\
      & \simeq
      \mathbf{H}
      \big(
        (\ast \!\sslash\! G) \times U
        \,,\,
        X \!\sslash\! G
      \big)
      \underset{\scalebox{.6}{$
        \mathbf{H}
        \big(
          (\ast \!\sslash\! G) \times U
          \,,\,
          \ast \!\sslash\! G
        \big)$}
      }{\times}
      \big\{
        \mathrm{pr}_1
      \big\}
      \\
      & \simeq
      \mathrm{Groupoids}
      \Big(
        \ast \!\sslash\! G
        \,,\,
      \mathbf{H}
      \big(
        U
        \,,\,
        X \!\sslash\! G
      \big)
      \Big)
      \underset{\scalebox{.6}{$
        \mathrm{Groupoids}
        \Big(
          \ast \!\sslash\! G
          \,,\,
        \mathbf{H}
        \big(
          U
          \,,\,
          \ast \!\sslash\! G
        \big)
        \Big)
        $}
      }{\times}
      \big\{
        \widetilde{\mathrm{pr}_1}
      \big\}
      \\
      & \simeq
      \mathrm{Groupoids}
      \Big(
        \ast \!\sslash\! G
        \,,\,
      \mathbf{H}
      \big(
        U
        \,,\,
        X
      \big) \!\sslash\! G
      \Big)
      \underset{\scalebox{.6}{$
        \mathrm{Groupoids}
        \big(
          \ast \!\sslash\! G
          \,,\,
          \ast \!\sslash\! G
        \big)
        $}
      }{\times}
      \big\{
        \mathrm{id}
      \big\}
      \\
      & \simeq
      \mathbf{H}(U,X)^G
      \,.
    \end{aligned}
  \end{equation}
  Here the first three lines are the definition of
  fixed loci \eqref{HomotopyFixedLocus} and the
  hom-equivalences \eqref{AdjunctionHomEquivalence} of the resulting
  adjunction \eqref{BaseChangeToGlobalContext}.
  The fourth line is the characterization \eqref{HomGroupoidInSlice} of
  hom-$\infty$-groupoids in
  slices (Prop. \ref{HomsInSliceInfinityCategory}),
  the fifth line uses the tensoring \eqref{InfinityTensoring} of $\mathbf{H}$
  over $\mathrm{Groupoids}_\infty$ (Prop. \ref{TensoringOfInfinityToposesOverInfinityGroupoids}),
  and the sixth line follows by Prop. \ref{HommingChartsIntoHomotopyFiberSequences}.

  To see the last step in \eqref{TowardsUnderstandingTheFixedLocusInA0TruncatedObject},
  use the explicit presentation of the
  groupoid $\mathbf{H}(U,X) \!\sslash\! G$ as an action groupoid,
  by Example \ref{ActionGroupoids}.
  This way the projection
  map in the fiber product in the sixth line
  in \eqref{TowardsUnderstandingTheFixedLocusInA0TruncatedObject}
  is presented by a Kan fibration,
  whence this homotopy fiber product may be computed
  equivalently as
  a 1-categorical fiber product
  of sets of objects and of sets of morphisms,
  separately. Moreover, since $\{\mathrm{id}\}$ has no non-trivial morphisms
  and since the projection functor itself is faithful,
  there are in fact no non-trivial morphisms in this fiber product,
  which is hence just the set
  whose elements are precisely those functors of action groupoids
  which are equal to the identity on labels in $G$:
  \vspace{-2mm}
  $$
      \mathrm{Groupoids}
      \Big(
        \ast \!\sslash\! G
        \,,\,
      \mathbf{H}
      \big(
        U
        \,,\,
        X
      \big) \!\sslash\! G
      \Big)
      \!
      \underset{\scalebox{.6}{$
        \mathrm{Groupoids}
        \big(
          \ast \!\sslash\! G
          \,,\,
          \ast \!\sslash\! G
        \big)
        $}
      }{\times}
      \!
      \big\{
        \mathrm{id}
      \big\}
      \;\;\simeq\;
      \left\{\!\!
      \raisebox{28pt}{
    \xymatrix@R=4pt{
      \ast \!\sslash\! G
      \ar[rr]
      &&
      \mathbf{H}(U,X) \!\sslash\! G
      \\
      \ast
      \ar@{}[rr]|-{ \longmapsto }
      \ar[dd]_-{g \in G}
      &&
      \phi
      \ar[dd]^-{ g }
      \\
      \\
      \ast
      \ar@{}[rr]|-{ \longmapsto }
      &&
      g \cdot \phi
    }
    }
   \!\! \right\}
    \;\;\simeq\;\;
    \mathbf{H}(U,X)^G.
  $$

\vspace{-6mm}
\hfill\end{proof}

\begin{lemma}[$n$-Truncated morphisms via $n$-truncated homotopy fibers]
 \label{nTruncatedMorphismViannTruncatedHomotopyFiber}
Let $\mathbf{H}$ be an $\infty$-topos
which is cohesive (Def. \ref{CohesiveTopos}).
Let $G$ be a finite group in $\mathbf{H}$ \eqref{InclusionOfFiniteGroups}.
Then, for every $n \in \{-2, -1, 0,1,\cdots\}$ and for any morphism in $\mathbf{H}$
to its delooping groupoid (Example \ref{DeloopingGroupoids})
$\mathcal{X} \overset{p}{\longrightarrow} \ast \!\sslash\! G$,
the following are equivalent
 \vspace{-1mm}
  \item {\bf (i)} $p$ is an $n$-truncated morphism (Def. \ref{nTruncatedMorphism});
  \vspace{-1mm}
  \item {\bf (ii)} the homotopy fiber of $p$ (over the essentially unique
  point of $\ast \!\sslash\! G$) is an $n$-truncated object (Def. \ref{nTruncatedObjects}).
\end{lemma}
\begin{proof}
Let $U \in \mathrm{Charts}$ and consider homming it into the
homotopy fiber sequence in question:
 \vspace{-2mm}
$$
  \raisebox{20pt}{
  \xymatrix@C=4em@R=1.5em{
    X
    \ar@{}[dr]|-{\mbox{\tiny (pb)}}
    \ar[d]
    \ar[r]
    &
    \mathcal{X}
    \ar[d]^-{p}
    \\
    \ast
    \ar[r]
    &
    \ast \!\sslash\! G
  }
  }
  \;\;\;\;\;
  \Rightarrow
  \;\;\;\;\;
  \raisebox{20pt}{
  \xymatrix@C=4em@R=1.5em{
    \mathbf{H}(U, X)
    \ar@{}[dr]|-{\mbox{\tiny (pb)}}
    \ar[d]
    \ar[r]
    &
    \mathbf{H}(U, \mathcal{X})
    \ar[d]^-{ \mathbf{H}(U,\,p) }
    \\
    \ast
    \ar[r]
    &
    \mathbf{H}(U, \ast \!\sslash\! G)
    \mathrlap{
      \;\;
      \simeq
      \ast \!\sslash\! G
    }
  }
  }
$$

 \vspace{-2mm}
\noindent
Since the hom-functor $\mathbf{H}(U,-)$ preserves limits,
the square on the right is again a homotopy pullback. Since
$U$ is a chart and $G$ is discrete, we have
the equivalence  \eqref{HomOfChartIntoDiscreteDeloopingGroupoid}
shown on the bottom right.
Since $\ast \!\sslash\! G$ has an essentially unique point,
the square on the right exhibits the essentially unique homotopy fiber
of the morphism $\mathbf{H}(U,p)$.
Since the charts $U$ are generators of $\mathbf{H}$ (objects of an $\infty$-site
of definition), the morphism $p$ is $n$-truncated  (Def. \ref{nTruncatedMorphism})
precisely if for each chart $U$ the homotopy fiber of $\mathbf{H}(U,p)$ is
$n$-truncated. But the square on the right shows that this homotopy fiber is $\mathbf{H}(U,X)$, and hence
this means, equivalently, that $X$ is an $n$-truncated object
(according to Def. \ref{nTruncatedObjects}).
\hfill \end{proof}

\medskip

\noindent {\bf Examples of cohesive $\infty$-toposes.}
We indicate some examples of cohesive $\infty$-toposes (Def. \ref{CohesiveTopos}),
following \cite{dcct}.
For full details of the constructions see \cite{SS20}.
\vspace{-1mm}
\begin{example}[Discrete cohesion]
  The base $\infty$-topos $\mathrm{Groupoids}_\infty$
  is trivially a cohesive $\infty$-topos (Def. \ref{CohesiveTopos})
  with all operations being identities:
  \begin{equation}
    \xymatrix{
      \mathrm{Groupoids}_\infty
      \;\;
      \ar@{->}@<+30pt>[rr]|<\times|-{\;
          \mathrm{id}
        \;}_-{\raisebox{-6pt}{\tiny $\bot$}}
        \ar@{<-^{)}}@<+15pt>[rr]|-{
         \; \mathrm{id} \;
        }_-{\raisebox{-6pt}{\tiny $\bot$}}
        \ar@<+0pt>@<-0pt>[rr]|-{\;
          \mathrm{id}
       \; }_-{\raisebox{-6pt}{\tiny $\bot$}}
        \ar@{<-^{)}}@<-15pt>[rr]|-{
         \; \mathrm{id} \;
        }
      &&
    \;\;  \mathrm{Groupoids}_\infty
    }
  \end{equation}

  For emphasis we also call this the $\infty$-topos
  of \emph{geometrically discrete $\infty$-groupoids}.
\end{example}

\begin{defn}[Site for homotopical cohesion]
  \label{SiteForHomotopicalCohesion}
  A small $\infty$-site \eqref{SheafTopos} is
  an  \emph{$\infty$-site for homotopical cohesion}
  if
   \vspace{-1.5mm}
\item {\bf (i)}  its Grothendieck topology is trivial
  and
   \vspace{-1.5mm}
  \item {\bf (ii)} the underlying $\infty$-category has finite products,
  i.e., has a terminal object and binary Cartesian products.
\end{defn}

\begin{example}[Homotopical cohesion]
  \label{InfinityCohesivePresheafSite}
  The
  $\infty$-topos of $\infty$-sheaves (Def. \ref{InfinityToposOfInfinitySheaves})
  over an $\infty$-site
  $\mathcal{C}$ for homotopical cohesion
  (Def. \ref{SiteForHomotopicalCohesion})
  is cohesive (Def. \ref{CohesiveTopos}):
  \vspace{-2mm}
  \begin{equation}
    \xymatrix{
      \mathbf{H}
        :=
        \mathrm{Sheaves}_\infty(\mathcal{C})
        \;\;
          \ar@{->}@<+30pt>[rr]|<\times|-{\;
          \underset{\longrightarrow}{\mathrm{lim}}
        \;}_-{\raisebox{-6pt}{\tiny $\bot$}}
        \ar@{<-^{)}}@<+15pt>[rr]|-{
         \; \mathrm{const} \;
        }_-{\raisebox{-6pt}{\tiny $\bot$}}
        \ar@<+0pt>@<-0pt>[rr]|-{\;
          \underset{\longleftarrow}{\mathrm{lim}}
       \; }_-{\raisebox{-6pt}{\tiny $\bot$}}
        \ar@{<-^{)}}@<-15pt>[rr]|-{
         \; \mathrm{Chtc} \;
        }
      &&
    \;\;  \mathrm{Groupoids}_\infty
    }
  \end{equation}

 \item  {\bf (i)} The operation
  $\mathrm{Pnts}
    \;\simeq\;
    \underset{\longleftarrow}{\mathrm{lim}}
    $
    forms the limit of $\infty$-presheaves
    regarded as $\infty$-functors on $\mathcal{C}^{\mathrm{op}}$
    (by Prop. \ref{LimitsAndColimitsAsAdjoints});
    but since $\mathcal{C}$ is assumed to have a terminal object,
    this is equivalently just the evaluation on that object:
     \vspace{-1mm}
    $$
      \mathrm{Pnts}(X)
      \;\simeq\;
      X(\ast)
      \;\simeq\;
      \mathbf{H}(\ast, X)
      \,,
    $$

     \vspace{-1mm}
\noindent
    where on the right we used the $\infty$-Yoneda lemma (Prop. \ref{YonedaLemma}).
    This makes manifest how $\mathrm{Pnts}(X)$ is the
    ``underlying $\infty$-groupoid of points of $X$''.

  \item {\bf (ii)} The operation
    $\mathrm{Shp} \simeq \underset{\longrightarrow}{\lim}$
    is the colimit of $\infty$-presheaves
    regarded as $\infty$-functors (by Prop. \ref{LimitsAndColimitsAsAdjoints}).
    Since the colimit of any representable functor
    is the point (Lemma \ref{InfinityColimitOverRepresentableInfinityFunctorIsContractible})
   \vspace{-2mm}
    $$
      \xymatrix{
        \mathcal{C} \;
        \ar@/_1pc/[rrrr]_-{ \mathrm{const}_\ast }
        \ar@{^{(}->}[rr]^-{ y }
        &&
        \mathrm{Sheaves}_\infty(\mathcal{C})
        \ar[rr]^-{ \mathrm{Shp} }
        &&
        \mathrm{Groupoids}_\infty\;,
      }
    $$

     \vspace{-2mm}
\noindent
      this means that $\mathcal{C}$ serves itself
    as a category of $\mathrm{Charts}$
    in this case (Def. \ref{ChartsForCohesion}).
\end{example}

\begin{example}[Smooth cohesion]
  \label{SmoothInfinityGroupoids}
  The $\infty$-sheaf $\infty$-topos (Def. \ref{InfinityToposOfInfinitySheaves}) over the
  site of $\mathrm{SmoothManifolds}$
  (Def. \ref{SmoothManifolds}, see \cite[App.]{FSS12}),
  which we call the $\infty$-topos of
  \emph{smooth $\infty$-groupoids}
   \vspace{-2mm}
  $$
    \mathrm{SmoothGroupoids}_\infty
    \;:=\;
    \mathrm{Sheaves}_\infty(\mathrm{SmoothManifolds})
    \,,
  $$

   \vspace{-2mm}
\noindent
  is cohesive (Def. \ref{CohesiveTopos}):
  The adjoint quadruple \eqref{AdjunctionCohesion}
  arises as in Example \ref{InfinityCohesivePresheafSite},
  which here happens to descend from $\infty$-presheaves to $\infty$-sheaves.

  \noindent In this case we have:
\vspace{-1mm}
 \item {\bf (i)} A category of $\mathrm{Charts}$ (Def. \ref{ChartsForCohesion})
  is given (Prop. \ref{CartesianSpacesIsDenseSubsite})
  by $\mathrm{CartesianSpaces}$
  (Def. \ref{CartesianSpaces})
  \vspace{-2mm}
  \begin{equation}
    \label{ChartsForSmoothGroupoids}
    \xymatrix@R=-2pt{
      \mathrm{CartesianSpaces}
     \; \ar@{^{(}->}[r]^-{ y }
      &
      \mathrm{Sheaves}_\infty(\mathrm{CartesianSpaces})
      \ar[rr]^-{\simeq}
      &&
      \mathrm{SmoothGroupoids}_\infty
      \\
      & &&
      \raisebox{1pt}{\textesh}(y(\mathbb{R}^n))
      \;\simeq\;
      \ast
    }
  \end{equation}

\vspace{-.3cm}
\item {\bf (ii)} The concrete 0-truncated objects
  (Def. \ref{ConcreteObjects})
  are equivalently the
  \emph{diffeological spaces} (Def. \ref{DiffeologicalSpaces}),
  including the \emph{D-topological spaces}\footnote{
    These are the \emph{$\Delta$-generated spaces} of
    \cite{Smith}\cite{Dugger03}; see Remark \ref{EuclideanGeneratedIsDeltaGenerated}.
  }
  (Def. \ref{DTopologicalSpaces})
  as well as smooth and possibly infinite-dimensional
  Fr{\'e}chet manifolds
  (Prop. \ref{SmoothManifoldsInsideDiffeologicalSpaces})
  as further full subcategories
  \eqref{AdjunctionBetweenTopologicalAndDiffeologicalSpaces}:

  \vspace{-3mm}
  \begin{equation}
    \label{ConcreteSmoothInfinityGroupoidsAreDiffeologicalSpaces}
    \hspace{-.5cm}
    \raisebox{+14pt}{
    \xymatrix@C=20pt@R=-2pt{
      \mathrm{TopologicalSpaces}
      \ar[r]^-{   \scalebox{.6}{$ \mathrm{Cdfflg}$} }
      &
      \mathrm{DTopologicalSpaces}
       \; \ar@{^{(}->}[rd]
      \\
      &
      &
      \mathrm{DiffeologicalSpaces}
      \;  \ar@{^{(}->}[rr]^-{
          \mbox{
            \tiny
            \begin{tabular}{c}
              concrete
              \\
              0-truncated
              \\
              objects
            \end{tabular}
          }
        }
        &&
      \mathrm{SmoothGroupoids}_\infty
      \\
      &
      \mbox{Fr{\'e}chetManifolds}
      \ar@{^{(}->}[ru]
    }
    }
  \end{equation}

  \noindent {\bf (iii)}
  The concrete 1-truncated objects (Def. \ref{ConcreteObjects})
  form the $(2,1)$-category of diffeological groupoids with
  Morita/Hilsum-Skandalis morphisms (Remark \ref{MoritaMorphismsOfGroupoids})
  between them,
  which includes, by \eqref{ConcreteSmoothInfinityGroupoidsAreDiffeologicalSpaces},
  the $(2,1)$-categories
  of D-topological groupoids and of
  (possibly infinite-dimensional Fr{\'e}chet-)Lie groupoids:
  \vspace{-2mm}
  \begin{equation}
    \label{1TrunConcreteSmoothInfinityGroupoidsAreDiffeologicalSpaces}
    \scalebox{.95}{
    \hspace{-.5cm}
    \xymatrix@R=-2pt@C=2em{
      \mathrm{TopologicalGroupoids}
      \ar[r]^-{    \scalebox{.6}{$ \mathrm{Cdfflg}$} }
      &
      \mathrm{DTopologicalGroupoids}
       \; \ar@{^{(}->}[rd]
      \\
      &
      &
      \mathrm{DiffeologicalGroupoids}
      \;  \ar@{^{(}->}[r]^-{
          \mbox{
            \tiny
            \begin{tabular}{c}
              concrete
              \\
              1-truncated
              \\
              objcts
            \end{tabular}
          }
        }
        &
      \mathrm{SmoothGroupoids}_\infty
      \\
      &
      \mbox{Fr{\'e}chetLieGroupoids}
      \ar@{^{(}->}[ru]
    }
    }
  \end{equation}

\vspace{-.1cm}
  \noindent {\bf (iv)}
  The cohesive shape \eqref{CohesiveModalitiesFromAdjointQuadruple}
  is given equivalently \cite{Pavlov}\cite{Pavlov19} by the
  \emph{smooth $\infty$-path $\infty$-groupoid}:
  \begin{equation}
    \label{SmoothShapeViaPathInfinityGroupoid}
    \raisebox{1pt}{\textesh}
    \,
    X
    \;\simeq\;
    \underset{\longrightarrow}{\mathrm{lim}}
    \,
    \mathbf{Maps}
    \big(
      \Delta^\bullet_{\mathrm{smth}},
      X
    \big)
    \;\in\;
    \mathrm{SmoothGroupoids}_\infty
    \,,
    \phantom{AA}
    \mbox{hence}
    \;\;\;
    \mathrm{Shp}(X)
    \;\simeq\;
    \underset{\longrightarrow}{\mathrm{lim}}
    \,
    X(\Delta^\bullet_{\mathrm{smth}})
    \;\in\;
    \mathrm{Groupoids}_\infty
  \end{equation}

  \vspace{-2mm}
\noindent   where $\Delta^\bullet_{\mathrm{smth}}$ is the simplicial smooth manifold of
  extended simplices (Def. \ref{SmoothSimplices}) and
  $\mathbf{Maps}(-,-)$ denotes the internal hom \eqref{InternalHomAdjunction}
  in $\mathrm{SmoothGroupoids}_\infty$.

  \noindent {\bf (v)}
  The cohesive shape \eqref{CohesiveModalitiesFromAdjointQuadruple}
  of
  {\bf (a)} any topological space
  and
  {\bf (b)} any finite-dimensional smooth manifold
  regarded, respectively, as
  smooth $\infty$-groupoids via \eqref{ConcreteSmoothInfinityGroupoidsAreDiffeologicalSpaces}
  is equivalently
  (by \eqref{SmoothShapeViaPathInfinityGroupoid}
  with Prop. \ref{DiffeologicalSingularSimpliciaSetOfContinuousDiffeology},
  and by \cite[4.3.29]{dcct}, respectively)
  its standard
  topological homotopy type $\mathrm{Shp}_{\mathrm{Top}}$
  \eqref{ShapeOfTopologicalSpaces}:
 \vspace{-2mm}
  \begin{equation}
    \label{CohesiveShapeOfTopologicalSpaces}
   \;\;\;\;\;\; \;\;\;\;  \mathllap{
     \mbox{\bf (a)}
      \;\;\;\;
    }
    \xymatrix@C=12pt{
      \mathrm{TopologicalSpaces}
      \ar@/_2pc/[rrrrrr]_{\mathrm{Shp}_{\mathrm{Top}}}^{
        \raisebox{6pt}{
          $\Downarrow \mathrlap{\scalebox{.7}{$\simeq$}}$
        }
      }
      \ar[rr]^-{    \scalebox{.6}{$\mathrm{Cdfflg}$} }
      &&
      \mathrm{DiffeologicalSpaces}
     \; \ar@{^{(}->}[rr]
      &&
      \mathrm{SmoothGroupoids}_\infty
      \ar[rr]^-{ \mathrm{Shp} }
      &&
      \mathrm{Groupoids}_\infty
    }
  \end{equation}
 \vspace{-3mm}
  \begin{equation}
    \label{CohesiveShapeOfSmoothManifolds}
     \;\;\;\;\;\; \;\;\;\;  \mathllap{
      \mbox{\bf (b)}
      \;\;\;\;
    }
    \xymatrix@C=12pt{
      \mathrm{SmoothManifolds}
      \ar@/_2pc/[rrrrrr]_{
        \mathrm{Shp}_{\mathrm{Top}}
        \circ
        \mathrm{Dtplg}
      }^{
        \raisebox{6pt}{
          $\Downarrow \mathrlap{\scalebox{.7}{$\simeq$}}$
        }
      }
      \ar@{^{(}->}[rr]^-{   }
      &&
      \mathrm{DiffeologicalSpaces}
     \; \ar@{^{(}->}[rr]
      &&
      \mathrm{SmoothGroupoids}_\infty
      \ar[rr]^-{ \mathrm{Shp} }
      &&
      \mathrm{Groupoids}_\infty
    }
  \end{equation}

  \noindent {\bf (vi)} The cohesive shape \eqref{CohesiveModalitiesFromAdjointQuadruple}
  of a topological groupoid,
  when regarded, via its coreflection \eqref{AdjunctionBetweenTopologicalAndDiffeologicalSpaces},
  as a D-topological groupoid and hence as a smooth $\infty$-groupoid
  \eqref{1TrunConcreteSmoothInfinityGroupoidsAreDiffeologicalSpaces}
  is equivalently
  (by \eqref{CohesiveShapeOfTopologicalSpaces}, and since
  $\raisebox{1pt}{\textesh}$ is left adjoint and hence preserves
  homotopy colimits, Prop. \ref{AdjointsPreserveCoLimits}) its simplicial-topological shape (Def. \ref{SimplicialTopologicalShape}):
  \vspace{-2mm}
  \begin{equation}
    \label{CohesiveShapeOfTopologicalGroupoids}
    \xymatrix@C=12pt{
      \mathrm{TopologicalGroupoids}
      \ar@/_2pc/[rrrrrr]_{\mathrm{Shp}_{\mathrm{sTop}}}^{
        \raisebox{6pt}{
          $\Downarrow \mathrlap{\scalebox{.7}{$\simeq$}}$
        }
      }
      \ar[rr]^-{    \scalebox{.6}{$ \mathrm{Cdfflg}$} }
      &&
      \mathrm{DiffeologicalGroupoids}
     \; \ar@{^{(}->}[rr]
      &&
      \mathrm{SmoothGroupoids}_\infty
      \ar[rr]^-{ \mathrm{Shp} }
      &&
      \mathrm{Groupoids}_\infty
    }
  \end{equation}
\end{example}

\begin{example}[Spectral cohesion]
  \label{CohesiveTangentInfinityTopos}
  Let $\mathbf{H}$ be a cohesive $\infty$-topos
  (Def. \ref{CohesiveTopos}).
  Then its tangent $\infty$-topos
  $T \mathbf{H} \;=\; \mathrm{SpectralBundles}(\mathbf{H})$
  (Example \ref{TangentInfinityTopos}) is cohesive
  \cite[4.1.9]{dcct}
  over the base tangent $\infty$-topos \eqref{BundlesOfPlainSpectra}:

  \vspace{-.1cm}
  \begin{equation}
    \label{TangentCohesion}
 \hspace{1cm}
    \xymatrix@R=1em{
      T\mathbf{H} \;\;
        \ar@{->}@<+28pt>[rr]|<\times|-{
          \;T\mathrm{Shp}\;
        }_-{\raisebox{-6pt}{\tiny $\bot$}}
        \ar@{<-^{)}}@<+14pt>[rr]|-{
         \; T\mathrm{Disc} \;
        }_-{\raisebox{-6pt}{\tiny $\bot$}}
        \ar@<+0pt>@<-0pt>[rr]|-{
         \; T\mathrm{Pnts} \;
        }_-{\raisebox{-6pt}{\tiny $\bot$}}
        \ar@{<-^{)}}@<-14pt>[rr]|-{
         \; T\mathrm{Chtc} \;
        }
      &&
     \;\;
     T\mathrm{Groupoids}_\infty
    }
  \end{equation}
 \end{example}

 \newpage

\begin{remark}[Differential cohomology in cohesive $\infty$-toposes]
  \label{DifferentialCohomologyTheory}
  The intrinsic cohomology theory \eqref{IntrinsicCohomologyOfAnInfinityTopos}
  of a cohesive $\infty$-topos (Def. \ref{CohesiveTopos}) is \emph{differential cohomology}
  \cite{dcct}.

\vspace{-1mm}
\item {\bf (i)}  In the case when $\mathbf{H} := \mathrm{SmoothGroupoids}_\infty$
  (Example \ref{SmoothInfinityGroupoids}),
  this is a non-abelian differential cohomology theory
  generalizing the theory of Cartan-Ehresmann connections
  on smooth fiber bundles to $\infty$-connections on
  smooth $\infty$-bundles \cite{SSS09}\cite{FSS12}\cite{NSS12}.

\vspace{-1mm}
  \item {\bf (ii)} In the case when
  $\mathbf{H} := T \mathrm{SmoothGroupoids}_{\infty}$
  is the cohesive tangent
  $\infty$-topos (Example \ref{CohesiveTangentInfinityTopos})
  to that of smooth $\infty$-groupoids (Example \ref{SmoothInfinityGroupoids}),
  the intrinsic cohomology
  furthermore
  subsumes abelian Hopkins-Singer differential cohomology
  theories and variants \cite{BNV13},
  as well as the twisted versions of these
  (Remark \ref{AbelianTwistedCohomology}),
  such as twisted differential KU-theory \cite{GS19a}
  and twisted differential KO-theory \cite{GS19b}.
\end{remark}

\subsubsection{Differential geometry}
\label{DifferentialGeometry}

We present a formulation of differential geometry
internal to $\infty$-toposes which we call \emph{elastic} \cite{dcct},
adjoining to the plain \emph{shape} operation $\mbox{\textesh}$ of \cref{DifferentialTopology}
a \emph{de Rham shape} operation $\Im$,
in generalization of \cite{Simpson96}\cite{SimpsonTeleman97}.

\begin{defn}[Elastic $\infty$-topos]
  \label{ElasticInfinityTopos}
  $\phantom{A}$
\vspace{-1mm}
 \item {\bf (i)} An \emph{elastic $\infty$-topos}
  over $\mathbf{B} = \mathrm{Groupoids}_\infty$
  is an $\infty$-topos $\mathbf{H}$ (Def. \ref{InfinityTopos})
  whose base geometric
  morphism (Prop. \ref{BaseGeometricMorphism}),
  to be denoted
  $\mathrm{Pnts} : \xymatrix@C=11pt{\mathbf{H} \ar[r] & \mathrm{Groupoids}_\infty}\,,$
  is equipped with
  a factorization as follows,
  having adjoints (Def. \ref{AdjointInfinityFunctors}) as shown:
\vspace{0mm}
  \begin{equation}
    \label{ElasticityAdjunctions}
    \hspace{2cm}
       \xymatrix@C=4em{
     \ar@{}@<+14pt>[r]|-{ \mathrm{Shp} \; : }
    \ar@{}@<-14pt>[r]|-{ \mathrm{Pnts} \; : }
      &
      \mathbf{H} \;
      \ar@{<-^{)}}@<+28pt>[rr]|-{
          \mathllap{
            \scalebox{.8}{
              \color{darkblue} \bf
              ``reduced''
              \hspace{4.1cm}
            }
          }
          \; \mathrm{Rdcd} \;
        }_-{\raisebox{-6pt}{\tiny $\bot$}}
        \ar@{->}@<+14pt>[rr]|-{
          \mathllap{
            \scalebox{.8}{
              \color{darkblue} \bf
              ``infinitesimal shape''
              \hspace{4cm}
            }
          }
        \;  \mathrm{Shp}_{\mathrm{inf}}  \;
        }_-{\raisebox{-6pt}{\tiny $\bot$}}
        \ar@{<-^{)}}@<+0pt>[rr]|-{
          \mathllap{
            \scalebox{.8}{
              \color{darkblue} \bf
              ``infinitesimally discrete''
              \hspace{4cm}
            }
          }
      \;     \mathrm{Disc}_{\mathrm{inf}} \;
        }_-{\raisebox{-6pt}{\tiny $\bot$}}
        \ar@<-14pt>[rr]|-{
          \mathllap{
            \scalebox{.8}{
              \color{darkblue} \bf
              ``infinitesimal points''
              \hspace{4cm}
            }
          }
          \; \mathrm{Pnts}_{\mathrm{inf}} \;
        }
        \ar@{<-^{)}}@<-28pt>[rrrr]|-{ \;\mathrm{Chtc} \;}
      &
      \ar@<-14pt>@{}[rr]|-{\raisebox{-.6pt}{$---$}}_-{\raisebox{-6pt}{\tiny $\bot$}}
      &
    \;  \mathbf{H}_{\Re}
        \ar@{->}@<+14pt>[rr]|<\times|-{ \; \mathrm{Shp}_\Re \;}_-{\raisebox{-6pt}{\tiny $\bot$}}
        \ar@{<-^{)}}@<+0pt>[rr]|-{\; \mathrm{Disc}_{\Re} \;}_-{\raisebox{-6pt}{\tiny $\bot$}}
        \ar@<+0pt>@<-14pt>[rr]|-{\; \mathrm{Pnts}_\Re  \;}
      &&
    \;\;  \mathbf{B}
     \!\! \ar@{}@<+0pt>[r]|-{ : \; \mathrm{Disc} }
      &
      \\
      &
      \mathclap{
        \mbox{
          \tiny
          \color{darkblue} \bf
          \begin{tabular}{c}
            elastic
            \\
            $\infty$-topos
          \end{tabular}
        }
      }
      &&
      \mathclap{
        \mbox{
          \tiny
          \color{darkblue} \bf
          \begin{tabular}{c}
            reduced
            \\
            sub-topos
          \end{tabular}
        }
      }
      &&
      \mathclap{
        \mbox{
          \tiny
          \color{darkblue} \bf
          \begin{tabular}{c}
            discrete
            \\
            sub-topos
          \end{tabular}
        }
      }
    }
  \end{equation}

  \vspace{-3mm}
  \item {\bf (ii)} Hence the elastic $\infty$-topos $\mathbf{H}$ is, in particular,
  a cohesive $\infty$-topos over $\mathbf{B}$, according to Def. \ref{CohesiveTopos},
  and so is its sub-$\infty$-topos $\mathbf{H}_\Re$ of reduced objects.

  \vspace{-2mm}
   \item {\bf (iii)} We write
   \vspace{-2mm}
  \begin{equation}
    \label{ElasticityModalities}
    \underset{
      \mbox{
        \footnotesize
        \color{darkblue}
        \bf
        ``reduced''
      }
    }{
      \big(
        \;
        \Re
          \;:=\;
        \mathrm{Rdcd}
          \circ
        \mathrm{Shp}_{\mathrm{inf}}
        \;
      \big)
    }
    \;\dashv\;
    \underset{
      \mbox{
        \footnotesize
        \color{darkblue} \bf
        ``{\'e}tale''
      }
    }{
    \big(
      \;
      \Im \;:=\; \mathrm{Disc}_{\mathrm{inf}} \circ \mathrm{Shp}_{\mathrm{inf}}
      \;
    \big)
    }
    \;\dashv\;
    \underset{
      \mbox{
        \footnotesize
        \color{darkblue} \bf
        ``locally constant''
      }
    }{
    \big(
      \;
      \mathcal{L} \;:=\; \mathrm{Disc}_{\mathrm{inf}} \circ \mathrm{Pnts}_{\mathrm{inf}}
      \;
    \big)
    }
    \;:\;
    \mathbf{H} \longrightarrow \mathbf{H}
  \end{equation}

  \vspace{-.2cm}

  for the further induced modalities \eqref{AdjointModalities}
  (\emph{elastic modalities}),
  accompanying the cohesive modalities of \eqref{CohesiveModalitiesFromAdjointQuadruple}.
\end{defn}

\medskip

\noindent {\bf Examples of elastic $\infty$-toposes.}
We indicate some examples of elastic $\infty$-toposes (Def. \ref{ElasticInfinityTopos}),
following \cite{dcct}.
For full details on the constructions, see \cite{SS20}.
\begin{defn}[Jets of Cartesian spaces]
  \label{FormalCartesianSpaces}
 Let $k \in \mathbb{N}$.

 \vspace{-1mm}
\item {\bf (i)}   We write
  \vspace{-2mm}
  \begin{equation}
    \label{FiniteJetsOfCartesianSpaces}
    \xymatrix@R=-2pt@C=3.5em{
      k\mathrm{JetsOfCartesianSpaces}
      \; \ar@{^{(}->}[r]^-{ C^\infty(-) }
      &
      \mathrm{CommutativeAlgebras}^{\mathrm{op}}_{\mathbb{R}}
      \\
      \mathbb{R}^n \times \mathbb{D}_W
      \ar@{|->}[r]
      &
      C^\infty(\mathbb{R}^n)
      \otimes_{\mathbb{R}}
      (\mathbb{R} \oplus W)
    }
  \end{equation}

  \vspace{-1mm}
  \noindent
  for the full subcategory of that of commutative
  $\mathbb{R}$-algebras on those which are tensor products
  of {\bf (a)} the algebra of real-valued smooth functions on a
  Cartesian space $\mathbb{R}^n$, with
  {\bf (b)} a finite-dimensional real algebra
  with a maximal ideal $W$ that is
  nilpotent of order $k+1$, in that $W^{k+1} = 0$.

 \vspace{-1mm}
 \item {\bf (ii)}
  We write
  \vspace{-4mm}
  \begin{equation}
    \label{InfiniteJetsOfCartesianSpaces}
    \xymatrix@R=-4pt{
      \infty\mathrm{JetsOfCartesianSpaces}
      \ar@{}[r]|-{ := }
      &
      \underset{k \in \mathbb{N}}{\bigcup}
      k\mathrm{JetsOfCartesianSpaces}
      \;
      \ar@{^{(}->}[rr]^-{ C^\infty(-) }
      &&
      \mathrm{CommutativeAlgebras}^{\mathrm{op}}_{\mathbb{R}}
      \\
      \mathbb{R} \times \mathbb{D}_W
      \ar@{|->}[rrr]
      &&&
      C^\infty\big( \mathbb{R}^n \big)
      \otimes_{\mathbb{R}}
      W
    }
  \end{equation}

  \vspace{-2mm}
\noindent
  for the analogous full subcategory where each $W$ is
  (finite dimensional and) nilpotent
  of some finite order.

   \vspace{-1mm}
   \item {\bf (iii)}
   We regard these categories as equipped with the coverage
   (Grothendieck pre-topology)
   whose covers are the families of morphisms of the form
  \vspace{-1mm}
  $$
    \big\{ \!\!
      \xymatrix{
        \mathbb{R}^n \times \mathbb{D}
        \ar[r]^{f_i \times \mathrm{id}}
        &
        \mathbb{R}^n \times \mathbb{D}
      }
    \!\! \big\}_{i \in I}
    \phantom{AA}
    \mbox{such that
     $
     \phantom{A}
     \big\{ \!\!
       \xymatrix{
         \mathbb{R}^n
         \ar[r]^{f_i}
         &
         \mathbb{R}^n
       }
    \!\!  \big\}_{i \in I}
     $
     \begin{tabular}{l}
       is a cover in
       $\mathrm{CartesianSpaces}$ (Def. \ref{CartesianSpaces}).
     \end{tabular}
     }
  $$
\end{defn}

\newpage

\begin{lemma}[Coreflections of jets of Cartesian spaces]
  \label{CoreflectionsOfJetsOfCartesianSpaces}
  Consider the $k\mathrm{JetsOfCartesianSpaces}$ from Def. \ref{FormalCartesianSpaces}.

  \noindent {\bf (i)} For $k = 0$, this is equivalently the category
  of plain Cartesian spaces of Def. \ref{CartesianSpaces}:

  \vspace{-2mm}
  $$
    0\mathrm{JetsOfCartesianSpaces}
    \;\simeq\;
    \mathrm{CartesianSpaces}\;.
  $$

\vspace{0mm}
  \noindent {\bf (ii)} For any $k \in \mathbb{N}$, the evident full inclusion
  of $k\mathrm{JetsOfCartesianSpaces}$ into
  $(k+1)\mathrm{JetsOfCartesianSpaces}$ is co-reflective
    \vspace{-3mm}
  $$
  \hspace{-2mm}
    \scalebox{.95}{
    \xymatrix{
      \infty\mathrm{JetsOfCartesianSpaces}
      \ar@<+6pt>@{<-^{)}}[r]^-{\scalebox{0.6}{$ \mathrm{Rdcd}_{\infty}$} }
      \ar@<-6pt>@{->}[r]_-{\scalebox{0.6}{$
        \mathrm{Shp}_{\mathrm{inf},\infty}
        $}
      }^-{
        \raisebox{1pt}{
          \scalebox{.8}{$
            \bot
          $}
        }
      }
      &
      \cdots
      \ar@<+6pt>@{<-^{)}}[r]^-{\scalebox{0.6}{$  \mathrm{Rdcd}_2 $}}
      \ar@<-6pt>@{->}[r]_-{\scalebox{0.6}{$
        \mathrm{Shp}_{\mathrm{inf},2}
        $}
      }^-{
        \raisebox{1pt}{
          \scalebox{.8}{$
            \bot
          $}
        }
      }
      &
      2\mathrm{JetsOfCartesianSpaces}
      \ar@<+6pt>@{<-^{)}}[r]^-{ \scalebox{0.6}{$ \mathrm{Rdcd}_1$} }
      \ar@<-6pt>@{->}[r]_-{\scalebox{0.6}{$
        \mathrm{Shp}_{\mathrm{inf},1}
        $}
      }^-{
        \raisebox{1pt}{
          \scalebox{.8}{$
            \bot
          $}
        }
      }
      &
      1\mathrm{JetsOfCartesianSpaces}
      \ar@<+6pt>@{<-^{)}}[r]^-{\scalebox{0.6}{$  \mathrm{Rdcd}$} }
      \ar@<-6pt>@{->}[r]_-{\scalebox{0.6}{$
        \mathrm{Shp}_{\mathrm{inf}}
        $}
      }^-{
        \raisebox{1pt}{
          \scalebox{.8}{$
            \bot
          $}
        }
      }
      &
      \mathrm{Cartesian Spaces}
    }
    }
  $$

 \vspace{-2mm}
\noindent   with
  \begin{equation}
    \label{InfinitesimalShapeInTermsOfFunctionAlgebras}
    C^\infty
    \left(
      \mathrm{Shp}_{\mathrm{inf},k}
      \big(
        \mathbb{R}^n \times \mathbb{D}_W
      \big)
    \right)
    \;\simeq\;
    C^\infty(\mathbb{R}^n)
    \otimes_{\mathbb{R}}
    (\mathbb{R} \oplus W)/W^{k+1}.
  \end{equation}
\end{lemma}
\begin{proof}
  Statement {\bf (i)} follows as a special case of the general fact, sometimes known
  as \emph{Milnor's exercise}
  (since the key idea is hinted at in \cite[Prob. 1-C]{MilnorStasheff74}),
  that passage to their real algebras of smooth
  functions embeds smooth manifolds fully faithfully into the opposite
  or real algebras (a general proof is in \cite[35.10]{KMS93},
  see also \cite{Gra05};
  for general perspective see \cite[6]{Nestruev03}) :
  \vspace{-2mm}
  $$
    \xymatrix{
      \mathrm{SmoothManifolds}\;
      \ar@{^{(}->}[rr]^-{ C^\infty(-) }
      &&
      \;\mathrm{CommutativeAlgebras}_{\mathbb{R}}^{\mathrm{op}}\;.
    }
  $$

  \vspace{-1mm}
\noindent
  Statement {\bf (ii)} follows readily from the definition,
  using the fact that algebra homomorphisms preserve order of nilpotency.
\hfill \end{proof}

\begin{example}[Jets of smooth $\infty$-groupoids]
  \label{FormalSmoothInfinityGroupoids}
  For $k \in \mathbb{N} \sqcup \{\infty\}$,
  the $\infty$-sheaf $\infty$-topos (Def. \ref{InfinityToposOfInfinitySheaves})
  over the site of $k$-jets of Cartesian spaces
  (Def. \ref{FormalCartesianSpaces})
   \vspace{-1mm}
  $$
    k\mathrm{JetsOfSmoothGroupoids}_{\infty}
    \;:=\;
    \mathrm{Sheaves}_\infty(k\mathrm{JetsOfCartesianSpaces})
  $$

   \vspace{-1mm}
\noindent
  is elastic (Def. \ref{ElasticInfinityTopos}),
  with
  $( \mathrm{Rdcd} \dashv \mathrm{Shp}_{\mathrm{inf}} )$ in
  \eqref{ElasticityAdjunctions}
  given by Kan extension of the co-reflections of sites
  from Lemma \ref{CoreflectionsOfJetsOfCartesianSpaces}:
  $$
    \xymatrix@C=4em{
      k\mathrm{JetsOfSmoothGroupoids}_\infty
      \ar@{<-^{)}}@<+28pt>[rr]|-{
          \; \mathrm{Rdcd} \;
        }_-{\raisebox{-6pt}{\tiny $\bot$}}
        \ar@{->}@<+14pt>[rr]|-{
          \; \mathrm{Shp}_{\mathrm{inf}}  \;
        }_-{\raisebox{-6pt}{\tiny $\bot$}}
        \ar@{<-^{)}}@<+0pt>[rr]|-{
          \;
          \mathrm{Disc}_{\mathrm{inf}}
          \;
        }_-{\raisebox{-6pt}{\tiny $\bot$}}
        \ar@<-14pt>[rr]|-{
          \; \mathrm{Pnts}_{\mathrm{inf}} \;
        }
        \ar@{<-^{)}}@<-28pt>[rrrr]|-{ \;\mathrm{Chtc} \;}
      &
      &
      \;
      \mathrm{SmoothGroupoids}_\infty
        \ar@{->}@<+14pt>[rr]|<\times|-{ \; \mathrm{Shp}_\Re \;}_-{\raisebox{-6pt}{\tiny $\bot$}}
        \ar@{<-^{)}}@<+0pt>[rr]|-{\; \mathrm{Disc}_{\Re} \;}_-{\raisebox{-6pt}{\tiny $\bot$}}
        \ar@<+0pt>@<-14pt>[rr]|-{\; \mathrm{Pnts}_\Re  \;}
      &&
      \;\;
      \mathrm{Groupoids}_\infty
    }
  $$

 \noindent {\bf (i)} Here for $k = 1$ we will, for short, abbreviate
   \vspace{-1mm}
  \begin{equation}
    \label{1JetsOfSmoothGroupoids}
    \mathrm{JetsOfSmoothGroupoids}_{\infty}
    \;:=\;
    1\mathrm{JetsOfSmoothGroupoids}_{\infty}
    \,.
  \end{equation}

   \vspace{-1mm}
\noindent {\bf (ii)}  For the case $k = \infty$,
  the underlying 1-topos is the
  ``Cahiers topos'' \cite{Dubuc79a}\cite{Kock86}\cite{KS17}.

  \noindent {\bf (iii)} For any $k$, we have:

\begin{itemize}
  \vspace{-2mm}
  \item[{\bf (a)}] The full sub-$\infty$-topos of reduced objects
  \eqref{ElasticityAdjunctions} is
  (by Lemma \ref{CoreflectionsOfJetsOfCartesianSpaces})
  that of smooth $\infty$-groupoids
  from Example \ref{SmoothInfinityGroupoids}
   \vspace{-1mm}
 \begin{equation}
    \label{SmoothGroupoidsAsReduced}
    \xymatrix{
      k\mathrm{JetsOfSmoothGroupoids}_\infty
      \ar@{<-^{)}}[rr]^-{ \mathrm{Disc}_{\mathrm{inf}} }
      &&
   \;   \mathrm{SmoothGroupoids}_\infty
    }
  \end{equation}

 \vspace{-1mm}
  \item[{\bf (b)}] the  0-truncated concrete objects
  (Def. \ref{ConcreteObjects}) are still equivalently the
  \emph{diffeological spaces} (Def. \ref{DiffeologicalSpaces})
  as was the case
  in \eqref{ConcreteSmoothInfinityGroupoidsAreDiffeologicalSpaces}
   \vspace{-3mm}
  \begin{equation}
    \label{ConcreteSmoothFormalInfinityGroupoidsAreDiffeologicalSpaces}
    \hspace{-1cm}
    \raisebox{+14pt}{
    \xymatrix@R=-4pt{
      \mathrm{DTopologicalSpaces}
       \; \ar@{^{(}->}[rrd]
      \\
      &&
      \mathrm{DiffeologicalSpaces}
      \;  \ar@{^{(}->}[rr]^-{
          \mbox{
            \tiny
            \begin{tabular}{c}
              0-truncated
              \\
              concrete
              \\
              objects
            \end{tabular}
          }
        }
        &&
      k\mathrm{JetsOfSmoothGroupoids}_\infty
      \\
      \mbox{Fr{\'e}chetManifolds}
      \ar@{^{(}->}[rru]
    }
    }
  \end{equation}

   \vspace{-1mm}
\noindent
  and, more generally, the 1-truncated concrete
  objects are still the \emph{diffeological groupoids},
  as was the case in \eqref{1TrunConcreteSmoothInfinityGroupoidsAreDiffeologicalSpaces}:
   \vspace{-3mm}
  \begin{equation}
    \label{1TrunConcreteFormalSmoothInfinityGroupoidsAreDiffeologicalSpaces}
    \hspace{-1cm}
    \raisebox{+14pt}{
    \xymatrix@R=-4pt{
      \mathrm{DTopologicalGroupoids}
       \; \ar@{^{(}->}[rrd]
      \\
      &&
      \mathrm{DiffeologicalGroupoids}
      \;  \ar@{^{(}->}[rr]^-{
          \mbox{
            \tiny
            \begin{tabular}{c}
              1-truncated
              \\
              concrete
              \\
              objects
            \end{tabular}
          }
        }
        &&
      k\mathrm{JetsOfSmoothGroupoids}_\infty
      \\
      \mbox{Fr{\'e}chetLieGroupoids}
      \ar@{^{(}->}[rru]
    }
    }
  \end{equation}

 \vspace{-2mm}
  \item[{\bf (c)}] A category of charts (Def. \ref{ChartsForCohesion})
  for $\mathrm{JetsOfSmoothGroupoids}_\infty$
  is given by $k\mathrm{JetsOfCartesianSpaces}$
  (Def. \ref{FormalCartesianSpaces}) itself.
  \end{itemize}
\end{example}

\newpage

\noindent {\bf {\'E}tale geometry.}
\begin{defn}[{\'E}tale-over-$X$ modality]
  \label{EtaleOverXModality}
  Let $\mathbf{H}$ be an elastic $\infty$-topos (Def. \ref{ElasticInfinityTopos})
  and $X \in \mathbf{H}$ an object.
  We say that the
  \emph{{\'e}tale-over-$X$} modality on the slice $\infty$-topos over $X$
  (Def. \ref{SliceInfinityTopos}) is the $\infty$-functor
  \vspace{-2mm}
  $$
    \raisebox{20pt}{
    \xymatrix@R=6pt{
      \mathbf{H}_{/X}
      \ar[rr]^-{ \Im_X }
      &&
      \mathbf{H}_{/X}
      \\
      Y
      \ar[dd]_-{f}
      &&
      Y \times_{\Im X} \Im{Y}
      \ar[dd]^-{    \scalebox{.6}{$(\eta^\Im_X)^\ast (\Im f)$} }
      \\
      \ar@{}[rr]|-{\longmapsto}
      &&
      \\
      X
      &&
      X
    }
    }
    {\phantom{AAAAAA}}
    \raisebox{22pt}{
    \xymatrix@R=7pt@C=10pt{
      X
      \ar@/^1pc/[drrrrrr]|-{\;\eta^\Im_X}
      \ar@/_1pc/[dddrr]_-f
      \ar@{-->}[drr]
      \\
      &&
      X \times_{\Im X} \Im Y
      \ar[dd]
      \ar[rrrr]
      \ar@{}[ddrrrr]|-{\mbox{\tiny(pb)}}
      &&&&
      \Im X
      \ar[dd]^-{ \Im f }
      \\
      \\
      &&
      Y
      \ar[rrrr]|-{\; \eta^\Im_Y }
      &&&&
      \Im Y
    }
    }
  $$
  which sends any morphism $f$ into $X$ to the pullback of
  its image under the plain {\'e}tale modality $\Im$ \eqref{ElasticityModalities}
  along its unit morphism \eqref{AdjunctionUnit}, hence to the left vertical morphism
  in the Cartesian square shown on the right.
\end{defn}

\begin{defn}[Local diffeomorphism]
  \label{FormallyEtaleMorphism}
  Let $\mathbf{H}$ be an elastic $\infty$-topos (Def. \ref{ElasticInfinityTopos}). We say that a morphism
  $Y \overset{f}{\to} X$ in $\mathbf{H}$ is
  a \emph{local diffeomorphism}
  if it is {\'e}tale-over-$X$ (Def. \ref{EtaleOverXModality})
   \vspace{-1mm}
  $$
    \Im_X(f)
    \;\simeq\;
    X\;,
  $$
 \vspace{-1mm}
\noindent
  hence (see Prop. \ref{EtaleTopos} for this implication)
  if the naturality square of the unit \eqref{AdjunctionUnit} of the $\Im$-modality
  \eqref{ElasticityModalities} is a Cartesian square:
  \vspace{-2mm}
  \begin{equation}
    \label{PullbackSquareForLocalDiffeomorphisms}
    \xymatrix@C=5em@R=1.5em{
      Y
      \ar[d]_-{f}^-{ \mbox{\tiny{\'e}t} }
      &
      \ar@{}[d]|-{ \mbox{$\Leftrightarrow$} }
      &
      Y
        \ar@{}[dr]|-{\mbox{\tiny (pb)}}
        \ar[d]_-{f}
        \ar[r]^-{\eta^{\Im}_Y}
      &
      \Im Y
        \ar[d]^-{\Im f}
      \\
      X
      &&
      X \ar[r]_-{\eta^{\Im}_X}
      &
      \Im X
    }
  \end{equation}
\end{defn}
\begin{lemma}[Closure of class of local diffeomorphisms]
  \label{ClosureOfLocalDiffeomorphisms}
  Let $\mathbf{H}$ be an elastic $\infty$-topos (Def. \ref{ElasticInfinityTopos}).
  The class of local diffeomorphisms in
  $\mathbf{H}$ (Def. \ref{FormallyEtaleMorphism})

  \noindent {\bf (i)} satisfies left-cancellation: given
  a pair of composable morphisms $f, g$
  where $g$ is a local diffeomorphism, then $f$ is so precisely if
  the composite $g \circ f$ is:
  \vspace{-2mm}
  \begin{equation}
    \label{LocalDiffeoComposition}
    \raisebox{20pt}{
    \xymatrix@R=1.5em{
      Z
      \ar[dr]_-{
        g \circ f
      }
      \ar[rr]^-{ f }
      &&
      Y
      \ar[dl]|-{\mbox{\tiny\rm{\'e}t}}^{ g }
      \\
      & X
    }
    }
    \;\;\;\;\;\;\;
    \Rightarrow
    \;\;\;\;\;
    \Big(
      \mbox{$f$ is a local diffeo}
      \;\;\;\;
      \Leftrightarrow
      \;\;\;\;
      \mbox{$g \circ f$ is a local diffeo}
    \Big).
  \end{equation}

  \noindent {\bf (ii)} is closed under pullbacks:
  if in a Cartesian square the right vertical morphism is
  a local diffeomorphism, then so is the left morphism
  $$
    \raisebox{20pt}{
    \xymatrix{
      Y' \times_X Y
      \ar@{}[drr]|-{ \mbox{\tiny\rm(pb)} }
      \ar[rr]
      \ar[d]_-{ g^\ast f }
      &&
      Y
      \ar[d]_-{\mbox{\tiny\rm{\'e}t}}^-{ f }
      \\
      Y'
      \ar[rr]_-{g}
      &&
      X
    }
    }
    \;\;\;\;\;\;\;\;
    \Rightarrow
    \;\;\;\;\;\;\;
    \mbox{$g^\ast f$ is a local diffeo}.
  $$
\end{lemma}
\begin{proof}
  This is a routine argument:
  {\bf (i)}
  For two composable morphisms, consider the pasting of
  their $\eta^\Im$-naturality squares
   \vspace{-2mm}
 $$
    \raisebox{20pt}{
    \xymatrix@R=14pt{
      Z
      \ar[d]_-{f}
      \ar[rr]^-{ \eta^\Im_Z }
      \ar@{}[drr]|-{ \mbox{\tiny(pb)} }
      &&
      \Im Z
      \ar[d]^-{ \Im f }
      \\
      Y
      \ar[d]_-{g}
      \ar[rr]|-{ \;\eta^\Im_Y }
      \ar@{}[drr]|-{ \mbox{\tiny(pb)} }
      &&
      \Im Y
      \ar[d]^-{ \Im g }
      \\
      X
      \ar[rr]_-{ \eta^\Im_X }
      &&
      \Im X
    }
    }
  $$

   \vspace{-3mm}
\noindent
  By the functoriality of $\Im$, the total rectangle is the
  $\eta^\Im$-naturality square of $g \circ f$. But, by the
  pasting law (Prop. \ref{PastingLaw})
  and the assumption that the bottom square is Cartesian,
  the total rectangle is Cartesian precisely if so is the top square.

\noindent  {\bf (ii)} For two morphisms with the same codomain, consider
  the pasting of their pullback square with the
  $\eta^\Im$-naturality square of one of them, as shown on the left
  here:
   \vspace{-3mm}
  $$
    \raisebox{20pt}{
    \xymatrix{
      Y' \times_X Y
      \ar[d]_-{ g^\ast f }
      \ar[rr]^-{ f^\ast g }
      \ar@{}[drr]|-{\mbox{\tiny\rm(pb)}}
      &&
      Y
      \ar[d]^-{f}
      \ar[rr]^{ \eta^\Im }
      \ar@{}[drr]|-{\mbox{\tiny\rm(pb)}}
      &&
      \Im Y
      \ar[d]^-{\Im f}
      \\
      Y'
      \ar[rr]_-{ g }
      &&
      X
      \ar[rr]_-{ \eta^{\Im}_X }
      &&
      \Im X
    }
    }
    \phantom{AAA}
    \simeq
    \phantom{AAA}
    \raisebox{20pt}{
    \xymatrix{
      Y \times_X Y'
      \ar[d]_-{ g^\ast f }
      \ar[rr]^-{ \eta^\Im_{(Y \times_X Y')} }
      &&
      \Im( Y \times_X Y' )
      \ar[rr]^-{ \Im( f^\ast f ) }
      \ar[d]_-{ \Im ( g^\ast f ) }
      \ar@{}[drr]|-{ \mbox{\tiny\rm(pb)} }
      &&
      \Im Y'
      \ar[d]^-{ \Im g }
      \\
      Y'
      \ar[rr]_-{ \eta^\Im_{Y'} }
      &&
      \Im Y'
      \ar[rr]_-{ \Im(f') }
      &&
      \Im X
    }
    }
  $$

   \vspace{-2mm}
\noindent
  By the naturality of $\eta^\Im$, this pasting diagram
  on the left is equivalent
  to that shown on the right. Moreover, if $f$ is a local
  diffeomorphisms, it follows that three of the squares are pullbacks
  (the rightmost one by using that $\Im$ is right adjoint and
  thus preserves pullbacks, Prop. \ref{AdjointsPreserveCoLimits}), as shown. With that, the pasting law
  (Prop. \ref{PastingLaw}) implies, first, that the total rectangle
  on the left is a pullback, hence also that on the left,
  and then that the remaining square on the right is a pullback.
  This means that $g^\ast f$ is a local diffeomorphism.
\hfill \end{proof}

\begin{defn}[Local neighborhood]
  \label{InfinitesimalNeighborhood}
  Let $\mathbf{H}$ be an elastic $\infty$-topos
  (Def. \ref{ElasticInfinityTopos}).
For $Y \overset{f}{\longrightarrow} X$
  a morphism in $\mathbf{H}$, we say that the
  corresponding
  \emph{local neighborhood}
  of $Y$ in $X$ is the purely {\'e}tale aspect of $f$,
  hence is the
  object $N_{{}_{\!f}} X \in \mathbf{H}_{\!/X}$
  given by $\Im_{\!/X}(f) \simeq  (\eta_X^\Im)^\ast(\Im f) $,
  hence given by the following homotopy pullback square:
   \vspace{-2mm}
  $$
    \xymatrix@C=5em@R=1.5em{
      N_{{}_{\! f}} X
      \ar@{}[dr]|-{
        \mbox{\tiny(pb)}
      }
      \ar[r]
      \ar[d]_{ \Im_{\!/X}(f) }
      &
      \Im X
      \ar[d]^-{ \Im f }
      \\
      Y
      \ar[r]_-{ \eta_X^\Im }
      &
      \Im Y
    }
  $$
\end{defn}

\begin{defn}[Tangent bundle]
  \label{InfinitesimalTangentBundle}
   Let $\mathbf{H}$ be an elastic $\infty$-topos
  (Def. \ref{ElasticInfinityTopos}).
Then for $X \in \mathbf{H}$ any object, we say that
  its \emph{infinitesimal tangent bundle} is
  $$
    T X
      \;:=\;
    X
      \underset{\Im X}{\times}
    X
    \;\in\;
    \mathbf{H}_{\!/X}
    \,,
  $$
  hence the left morphism in
  this Cartesian square:
   \vspace{-2mm}
  \begin{equation}
    \label{TangentBundlePullbackDefinition}
    \xymatrix@C=5em@R=1.5em{
      T X
      \ar[r]
      \ar[d]_-{
        (\eta_X^\Im)^\ast (\eta^\Im_X)_!
        (\mathrm{id}_X)
      }
      \ar@{}[dr]|-{
        \mbox{\tiny(pb)}
      }
      &
      X
      \ar[d]^-{ \eta_X^\Im }
      \\
      X
      \ar[r]_-{ \eta_X^\Im }
      &
      \Im X
    }
  \end{equation}
\end{defn}

\begin{example}[Local neighborhood of a point]
  \label{LocalNeighbourhoodOfAPoint}
  Let $\mathbf{H}$ be an elastic $\infty$-topos
  (Def. \ref{ElasticInfinityTopos}).
  For $X \in \mathbf{H}$ any object and
  $\ast \overset{x}{\longrightarrow} X$ any point,
  the homotopy fiber of the tangent bundle
  (Def. \ref{InfinitesimalTangentBundle}) over $x$
  is equivalent to
  the local neighborhood of $x$ (Def. \ref{InfinitesimalNeighborhood}):
  \begin{equation}
    \label{LocalNeighbourhoodOfAPointIsFiverOfTangentBundle}
    T_x X \;\simeq\; N_x X
    \,.
  \end{equation}
  This follows immediately from the definitions, by the
  pasting law (Prop. \ref{PastingLaw}):
   \vspace{-1mm}
  $$
    \xymatrix@C=5em@R=1.5em{
      \mathllap{
        N_x X
        \simeq
        \;
      }
      T_x X
      \ar@{}[dr]|-{ \mbox{\tiny(pb)} }
      \ar[d]
      \ar[r]
      &
      T X
      \ar@{}[dr]|-{ \mbox{\tiny(pb)} }
      \ar[r]
      \ar[d]
      &
      X
      \ar[d]^-{ \eta_X^\Im }
      \\
      \ast
      \ar[r]_-{ x }
      &
      X
      \ar[r]_-{
        \eta_X^\Im
      }
      &
      \Im X
    }
  $$
\end{example}

\begin{prop}[Pullback along local diffeomorphisms preserves tangent bundles]
  \label{PullbackAlongLocalDiffeomorphismsPreservesTangentBundles}
  In an elastic $\infty$-topos
  (Def. \ref{ElasticInfinityTopos}),
  pullback along a local diffeomorphism
  $\xymatrix{ Y \ar[r]^-{f}_-{\mbox{\tiny{\'e}t}} & X }$
  (Def. \ref{FormallyEtaleMorphism})
  preserves tangent bundles (Def. \ref{InfinitesimalTangentBundle})
  in that
 \vspace{-3mm}
  $$
    f^\ast( T X )
    \;\simeq\;
    T Y
    \phantom{AAA}
    \mbox{via:}
    \phantom{AAA}
    \raisebox{20pt}{
    \xymatrix@R=1.5em@C=5em{
      T Y
      \ar@{}[dr]|-{ \mbox{\tiny\rm(pb)} }
      \ar[d]
      \ar[r]^-{ T f }
      &
      T X
      \ar[d]
      \\
      Y \ar[r]_-{f}^-{ \mbox{\tiny{\'e}t} }
      &
      X
    }
    }
  $$
\end{prop}
\begin{proof}
  Consider the pasting of the defining Cartesian
  squares, shown on the left here:
   \vspace{-2mm}
  $$
    \raisebox{20pt}{
    \xymatrix@R=1.5em@C=5em{
      f^\ast T X
      \ar@{}[dr]|-{ \mbox{\tiny (pb)} }
      \ar[d]
      \ar[r]
      &
      T X
      \ar[d]
      \ar[r]
      \ar@{}[dr]|-{ \mbox{\tiny (pb)} }
      &
      X
      \ar[d]^-{ \eta_X^\Im }
      \\
      Y
      \ar[r]_-{f}^-{\mbox{\tiny{\'e}t}}
      &
      X
      \ar[r]_-{ \eta_X^\Im }
      &
      \Im X
    }
    }
    \phantom{AA}
    \simeq
    \phantom{AA}
    \raisebox{20pt}{
    \xymatrix@R=1.5em@C=5em{
      T Y
      \ar[d]
      \ar[r]
      \ar@{}[dr]|-{\mbox{\tiny(pb)}}
      &
      Y
      \ar[d]^-{ \eta_X^\Im }
      \ar[r]^-{ f }
      \ar@{}[dr]|-{ \mbox{\tiny(pb)} }
      &
      X
      \ar[d]^-{ \eta_X^\Im }
      \\
      Y
      \ar[r]_-{ \eta^\Im_X }
      &
      \Im Y
      \ar[r]_-{ \Im f }
      &
      \Im X
    }
    }
  $$

  \vspace{-2mm}
\noindent
By the pasting law (Prop. \ref{PastingLaw}),
  the total rectangle on the left
  is itself Cartesian. Moreover, the
  bottom composite morphism on the left is equivalent
  to the bottom composite morphism on the right,
  by the naturality of $\eta_X^\Im$.
  Therefore, using again the pasting law (Prop. \ref{PastingLaw})
  the total rectangle on the left is equivalent to the
  pasting of the two consecutive Cartesian squares
  shown on the right. These identify, in the top row,
  the middle object $Y$ by
   \eqref{PullbackSquareForLocalDiffeomorphisms}
  and thus the left object $T Y$ by \eqref{TangentBundlePullbackDefinition}.
\hfill \end{proof}

\medskip

\noindent {\bf {\'E}tale toposes.}
\begin{defn}[{\'E}tale topos]
  \label{EtaleTopos}
  Let $\mathbf{H}$ be an elastic $\infty$-topos (Def. \ref{ElasticInfinityTopos})
  and $X \in \mathbf{H}$.
  Then we say that the \emph{{\'e}tale $\infty$-topos} of $X$,
  to be denoted $\mbox{\bf{\'E}t}_{X}$, is the full sub-$\infty$-category
  (Def. \ref{FullyFaithfulFunctor})
  of the slice $\infty$-topos over $X$ (Prop. \ref{SliceInfinityTopos})
  on those morphisms that are local diffeomorphisms (Def. \ref{FormallyEtaleMorphism}):
  \vspace{-2mm}
  \begin{equation}
    \label{EtaleToposInsideSliceTopos}
    \mbox{\bf{\'E}t}_{X}
    \;:=\;
    \xymatrix{
      \big(\mathbf{H}_{/X}\big)_{\Im_X}
      \ar@{^{(}->}[r]
      &
      \mathbf{H}_{/X}
    }.
  \end{equation}
\end{defn}

\begin{prop}[Reflections of {\'e}tale toposes]
  \label{PropertiesOfEtaleToposes}
  Let $\mathbf{H}$ be an elastic $\infty$-topos (Def. \ref{ElasticInfinityTopos})
  and $X \in \mathbf{H}$ an object. Then the {\'e}tale topos
  $\mbox{\bf{{\'E}t}}_{X}$ from Def. \ref{EtaleTopos}:

  \noindent {\bf (i)} is indeed an $\infty$-topos (Def. \ref{InfinityTopos});

  \noindent {\bf (ii)} its defining full inclusion \eqref{EtaleToposInsideSliceTopos}
  has both a left- and a right-adjoint (Def. \ref{AdjointInfinityFunctors}):
  \vspace{-2mm}
  \begin{equation}
    \label{AdjointsToInclusionOfEtaleSlice}
    \xymatrix{
      \mbox{\bf{\'E}t}_{X}
      \;
      \ar@<+14pt>@{<-}[rr]^-{
        \mathrm{Etl}_{X}
      }_-{ \bot }
      \ar@{^{(}->}[rr]|-{\;
        i_X
      \,}
      \ar@<-14pt>@{<-}[rr]_-{
        \mathrm{LcllCnstnt}_X
      }^-{ \bot }
      &&
      \;
      \mathbf{H}_{/X}
    }
  \end{equation}

\vspace{-2mm}
  \noindent {\bf (iii)} whose induced adjoint modality \eqref{AdjointModalities}
  \vspace{-2mm}
  \begin{equation}
    \label{ElateOverXModality}
    \big(
    \;
    \underset{
      \raisebox{-3pt}{
        \footnotesize
        \color{darkblue}
        \bf
        ``{\'e}tale over $X$''
      }
    }{
      \Im_X
      \;:=\;
      i_X \circ \mbox{\rm{\'E}tl}_X
    }
    \;
    \big)
    \;\;\dashv\;\;
    \big(
    \;
    \underset{
      \mathclap{
      \raisebox{-3pt}{
        \footnotesize
        \color{darkblue}
        \bf
        ``locally constant over $X$''
      }
      }
    }{
      \mathcal{L}_X
      \;:=\;
      i_X \circ \mbox{\rm{LcllCnstnt}}_X
    }
    \;
    \big)
  \end{equation}

  \vspace{-2mm}
\noindent
  is on the left that of Def. \ref{EtaleOverXModality}:
  \vspace{-3mm}
\begin{equation}
  \label{EtlX}
  \mbox{\rm\'Etl}_X
  \;:\;
  \raisebox{20pt}{
  \xymatrix{
    Y
    \ar[d]^-p
    \\
    X
  }
  }
  \;\;\;\longmapsto\;\;\;
  \raisebox{20pt}{
  \xymatrix{
    (\eta_X^\Im)^\ast (\Im Y)
    \ar[d]_-{
      (\eta_X^\Im)^\ast (\Im p)
    }
    \\
    X
  }
  }
  \phantom{AA}
  \mbox{i.e.:}
  \phantom{AA}
  \raisebox{20pt}{
  \xymatrix{
    (\eta_X^\Im)^\ast (\Im Y)
    \ar[d]_-{
      (\eta_X^\Im)^\ast (\Im p)
    }
    \ar[rr]^-{
      (\Im p)^\ast ( \eta_X^\Im )
    }
    \ar@{}[drr]|-{
      \mbox{\tiny (pb)}
    }
    &&
    \Im Y
    \ar[d]^-{ \Im p }
    \\
    X \ar[rr]_-{ \eta_X^\Im }
    &&
    \Im X
  }
  }
\end{equation}
\end{prop}
\begin{proof}
First to see that \eqref{EtlX} is well-defined as a functor to $\mbox{\bf{\'E}t}_X$
(this proceeds as in
\cite[3.3]{CassidyHebertKelly85}\cite[3]{CarboniJanelidzeKellyPare97}\cite[7.3]{CherubiniRijke20}):
We need to check that $(\eta_X^\Im)^\ast (\Im p)$ is a local diffeomorphism
(Def. \ref{FormallyEtaleMorphism}). For this, it is sufficient to have
equivalences
\begin{equation}
  \label{EquivalentConditionForEtlXToTakeValuesInLocalDiffeos}
  \Im\big(
    (\eta_X^\Im)^\ast (\Im p)
  \big)
  \;\simeq\;
  \Im p
  \,,
\end{equation}
and
\begin{equation}
  \label{SecondEquivalentConditionForEtlXToTakeValuesInLocalDiffeos}
  (\Im p)^\ast ( \eta_X^\Im )
  \;\simeq\;
  \eta_X^\Im
\end{equation}
because then the Cartesian square on the right of \eqref{EtlX}
exhibits this property.

But \eqref{EquivalentConditionForEtlXToTakeValuesInLocalDiffeos} follows
by applying $\Im$ to the square on the right of \eqref{EtlX},
by idempotency (Prop. \ref{IdempotentMonads}) and since equivalences are preserved by
pullback (Example \ref{PullbackOfEquivalenceIsEquivalence}).
With this, \eqref{SecondEquivalentConditionForEtlXToTakeValuesInLocalDiffeos}
follows from the naturality of the $\Im$-unit,
by the universal factorization shown dashed in the following diagram:

\newpage

\begin{equation}
  \label{EtlxXTowards}
  \raisebox{20pt}{
  \xymatrix@C=2.5em@R=10pt{
    (\eta_X^\Im)^\ast (\Im Y)
    \ar@{-->}[dr]^-{\simeq}
    \ar@/^1pc/[drrrrr]^-{ \eta^\Im_{(\eta_X^\Im)^\ast (\Im Y)} }
    \ar@/_1pc/[dddr]_-{ (\eta_X^\Im)^\ast (\Im p) }
    \\
    &
    (\eta_X^\Im)^\ast (\Im Y)
    \ar[dd]|-{
      (\eta_X^\Im)^\ast (\Im p)
    }
    \ar[rrrr]|-{\;
      (\Im p)^\ast ( \eta_X^\Im )
   \; }
    \ar@{}[ddrrrr]|-{
      \mbox{\tiny (pb)}
    }
    &&&&
    \Im Y
    \ar[dd]^-{ \Im p }
    \\
    \\
    &
    X \ar[rrrr]_-{ \eta_X^\Im }
    &&&&
    \Im X
  }
  }
\end{equation}
Notice that, similarly, there is a natural transformation
  \vspace{-2mm}
\begin{equation}
  \label{UnitForEtlX}
  \xymatrix@R=6pt@C=3em{
    Y
    \ar[dr]_-{p}
    \ar[rr]^-{ \eta^{\mathrm{Etl}_X}_Y }
    &&
    \mathrm{Etl}_X(Y)
    \ar[dl]^-{ \mbox{\tiny{\'et}} }
    \\
    &
    X
  }
\end{equation}
induced as the universal factorization shown dashed in the following diagram:
  \vspace{-2mm}
\begin{equation}
  \label{EtlxAdjunctionUnit}
  \raisebox{20pt}{
  \xymatrix@C=2.8em@R=10pt{
    Y
    \ar@{-->}[dr]^-{\simeq}
    \ar@/^1pc/[drrrrr]^-{ \eta^\Im_{Y} }
    \ar@/_1pc/[dddr]_-{ p }
    \\
    &
    (\eta_X^\Im)^\ast (\Im Y)
    \ar[dd]|-{
      (\eta_X^\Im)^\ast (\Im p)
    }
    \ar[rrrr]|-{\;
      (\Im p)^\ast ( \eta_X^\Im )
    \;}
    \ar@{}[ddrrrr]|-{
      \mbox{\tiny (pb)}
    }
    &&&&
    \Im Y
    \ar[dd]^-{ \Im p }
    \\
    \\
    &
    X \ar[rrrr]_-{ \eta_X^\Im }
    &&&&
    \Im X
  }
  }
\end{equation}

  \vspace{-2mm}
\noindent and notice that this in an $\Im$-equivalence:
\begin{equation}
  \label{EtlXUnitIsImEquivalence}
  \Im
  \big(
    \eta^{\mathrm{Etl}_X}_{Y_1}
  \big)
  \phantom{AA}
  \mbox{is an equivalence}
  \,.
\end{equation}
Condition \eqref{EtlXUnitIsImEquivalence} follows by applying $\Im$ to the whole left part of
the diagram on the right of \eqref{TowardsAdjointEtlX}, using
idempotency (Prop. \ref{IdempotentMonads}) and that equivalences are preserved by
pullback (Example \ref{PullbackOfEquivalenceIsEquivalence}).

\medskip

\noindent Second, to see that \eqref{EtlX} defines a left adjoint to the inclusion:
We need to check the corresponding hom-equivalence \eqref{AdjunctionHomEquivalence},
shown on the left here:

\vspace{-.6cm}

\begin{equation}
  \label{TowardsAdjointEtlX}
  \hspace{-4mm}
  \raisebox{20pt}{
  \xymatrix@C=.5em{
    \mbox{\'Etl}_X(Y_1)
    \ar[dr]_-{ \scalebox{.6}{$\mbox{\'Etl}_X(p)$} }
    \ar[rr]^-{
      \widetilde f
    }
    &&
    Y_2
    \ar[dl]^-{\mbox{\tiny{\'e}t}}
    \\
    & B
  }
  }
  \Leftrightarrow
  \raisebox{20pt}{
  \xymatrix@C=.5em{
    Y_1
    \ar[dr]_-{p}
    \ar[rr]^-{f}
    &&
    Y_2
    \ar[dl]^-{\mbox{\tiny{\'e}t}}
    \\
    & B
  }
  }
  \simeq
    \raisebox{62pt}{
  \xymatrix@C=1.5em@R=10pt{
    &&&
    \Im Y_1
    \ar[dddrr]|>>>>>>>>>>{\phantom{AA} \atop \phantom{AA}}
    \ar[rrrr]^-{ \Im f }
    && &&
    \;\;\;\;\Im Y_2\;\;\;\;
    \ar[dddll]
    \\
    \\
    Y_1
    \ar[dddrr]_-{p}
    \ar@/^1pc/[uurrr]^-{ \eta^\Im_{Y_1} }
    \ar@{-->}[rr]|-{ \eta^{\mathrm{Etl}_X}_{Y_1} }
    &
    &
    (\eta^\Im_X)^\ast (\Im X)
    \ar[uur]|-{\scalebox{.7}{$
      \eta^\Im_{(\eta^\Im_X)^\ast (\Im X)}
      $}
    }
    \ar[ddd]|-{
      \mathrm{Etl}_X(p)
    }
    \ar@{-->}[rr]|-{\; \widetilde f \;}
    &&
    Y_2
    \ar[dddll]|-{\;\mbox{\tiny{\'e}t}\;}
    \ar[uurrr]|-{ \eta^\Im_{Y_2} }
    \\
    & && &&
    \Im X
    \\
    \\
    &&
    X
    \ar[uurrr]_-{\eta^\Im_X}
  }
  }
\end{equation}
On the right of \eqref{TowardsAdjointEtlX} we show an
induced factorization:
The square sub-diagram on the right of \eqref{TowardsAdjointEtlX} is Cartesian by the assumption that we are homming into a local diffeomorphism,
while the square in the middle is Cartesian by \eqref{EtlxXTowards}.
Thus, given $f$, the morphism $\widetilde f$ is induced by
the universal property of the right Cartesian square.
Conversely, given $\widetilde f$, precomposition with the
$\eta^{\mathrm{Etl}_X}_{Y_1}$ \eqref{EtlxAdjunctionUnit} gives a morphism $f$.
To see that this correspondence is an equivalence, we just
need to observe that $\Im (\widetilde f) \simeq \Im f$.
This follows by \eqref{EtlXUnitIsImEquivalence}.

\medskip

Thus we have established the existence of the left adjoint
$\mbox{\'Etl}_X$. With this, to see the
right adjoint $\mathrm{LcllCnst}_X$ as well as the fact that
$\mbox{\bf{\'E}t}$ is an $\infty$-topos, it is now sufficient
to show that $\xymatrix@C=12pt{\mbox{\bf{\'E}t}_X \ar@{^{(}->}[r]^{i_X} & \mathbf{H}_{/X} }$
preserves colimits:
Because, by the reflection $\mbox{\'Etl}_X$ this
implies, first, that $\mbox{\bf{\'E}t}_X$ is a presentable $\infty$-category, in fact an $\infty$-topos
(by Prop. \ref{ToposLexReflection},
since it is then an accessibly embedded reflective subcategory of the slice
$\mathbf{H}_{/X}$, which is an $\infty$-topos by Prop. \ref{SliceInfinityTopos});
and thus, second, the existence of the right adjoint by the
adjoint $\infty$-functor theorem (Prop. \ref{AdjointFunctorTheorem}).

\newpage

So to see that $i_X$ preserves colimits, consider any small $\mathcal{I} \in \mathrm{Categories}_\infty$
and a diagram
\vspace{-2mm}
\begin{equation}
  \label{DiagramInTheEtaleSlice}
  Y_\bullet
  :
  \xymatrix@C=24pt{
    \mathcal{I}
    \ar[r]
    &
    \mbox{\bf{\'E}t}_X
    \; \ar@{^{(}->}[r]^-{ i_X }
    &
    \mathbf{H}_{/X} \;.
  }
\end{equation}

\vspace{-2mm}
\noindent
Since $i_X$ is fully faithful by construction,
it is sufficient to show that the colimit of this diagram formed in
$\mathbf{H}_{/X}$ is itself in the image of $i_X$.
This colimit, in turn, is computed in $\mathbf{H}$
(by Example \ref{BaseChangeAlongTerminalMorphism}) with its
morphism $q$ to $X$ universally induced, and this we need to
show to be a local diffeomorphism (Def. \ref{FormallyEtaleMorphism}).
Hence we need to show that
the following square on the left is Cartesian:
\vspace{-2mm}
$$
  \raisebox{20pt}{
  \xymatrix{
    \underset{\longrightarrow}{\mathrm{lim}}\, Y_\bullet
    \ar[d]_-{ q }
    \ar[rr]^-{
      \eta^\Im_{\underset{\longrightarrow}{\mathrm{lim}} Y_\bullet}
    }
    \ar@{}[drr]|-{ \mbox{\tiny(pb)} }
    &&
    \Im
    \big(
      \underset{\mathclap{\longrightarrow}}{\mathrm{lim}} Y_\bullet
    \big)
    \ar[d]^-{ \Im q }
    \\
    X
    \ar[rr]_-{ \eta^\Im_X }
    &&
    \Im X
  }
  }
  \phantom{AAA}
  \Leftrightarrow
  \phantom{AAA}
  \raisebox{20pt}{
  \xymatrix{
    \underset{\longrightarrow}{\mathrm{lim}} Y_\bullet
    \ar[d]_-{ q }
    \ar[rr]^-{
      (\eta^\Im_{Y_\bullet})
    }
    \ar@{}[drr]|-{\mbox{\tiny(pb)}}
    &&
    \underset{\longrightarrow}{\mathrm{lim}}
    \big(
      \Im Y_\bullet
    \big)
    \ar[d]^-{ \Im q }
    \\
    X
    \ar[rr]_-{ \eta^\Im_X }
    &&
    \Im X
  }
  }
  \phantom{AAA}
  \Leftrightarrow
  \phantom{AAA}
  \underset{i \in \mathcal{I}}{\forall}
  \raisebox{20pt}{
  \xymatrix{
    Y_i
    \ar[d]_-{ q_i }
    \ar[rr]^-{
      \eta^\Im_{Y_i}
    }
    \ar@{}[drr]|-{\mbox{\tiny(pb)}}
    &&
    \Im Y_i
    \ar[d]^-{ \Im q_i }
    \\
    X
    \ar[rr]_-{ \eta^\Im_X }
    &&
    \Im X
  }
  }
$$

\vspace{-2mm}
\noindent
But, since $\Im$ is a left adjoint and hence preserves colimits
(Prop. \ref{AdjointsPreserveCoLimits}), this is equivalent to the
square on in middle being Cartesian. Finally, by universality
of colimits \eqref{PreserveColimitsByPullback} in the $\infty$-topos
$\mathbf{H}$, this is equivalent to all the squares on the right being
Cartesian. This is the case, by the assumption \eqref{DiagramInTheEtaleSlice}.
\hfill \end{proof}

\begin{remark}[Local and global $\infty$-section functors.]
  \label{GlobalSections}
  Let $\mathbf{H}$ be an elastic $\infty$-topos (Def. \ref{ElasticInfinityTopos})
  and $X \in \mathbf{H}$.
  Then we may think of the {\'e}tale $\infty$-topos
  $\mbox{\bf{\'E}t}_{X}$ (Def. \ref{EtaleTopos}, Prop. \ref{PropertiesOfEtaleToposes})
  as the internal construction of the
  $\infty$-topos of \emph{$\infty$-sheaves over $X$}. Under this
  interpretation:

  \noindent
  {\bf i)} the $\infty$-functor $\mathrm{LcllCnst}$
  \eqref{AdjointsToInclusionOfEtaleSlice} has the interpretation
  of sending any $\infty$-bundle
  $\xymatrix@C=12pt{ E \ar[r] & X }$ (Notation \ref{SliceCategory})
  to its \emph{$\infty$-sheaf of local sections}
  $\underline{E} := \mathrm{LcllCnst}_X(E)$;

  \noindent
  {\bf ii)} the direct image of the base geometric morphism
  \eqref{BaseGeometricMorphismAdjunction} has the interpretation of
  sending any $\infty$-sheaf to its $\infty$-groupoid of
  global sections:

  \vspace{-.4cm}
  \begin{equation}
    \xymatrix@C=44pt{
      \overset{
        \mathclap{
        \raisebox{3pt}{
          \tiny
          \color{darkblue}
          \bf
          \begin{tabular}{c}
            $\infty$-bundles
            \\
            over $X$
          \end{tabular}
        }
        }
      }{
        \mathbf{H}_{\!/X}
      }
      \ar@{<-^{)}}@<+6pt>[rr]^-{ i_X }
      \ar@<-6pt>[rr]_-{
        \underset{
          \mathclap{
            \raisebox{3pt}{
              \tiny
              \color{darkblue}
              \bf
              \begin{tabular}{c}
                form $\infty$-sheaf
                of local sections
              \end{tabular}
            }
          }
        }{
        \underline{(-)}
        \,:=\,
        \mathrm{LcllCnstnt}_X
        }
      }^-{ \bot }
      \ar@<-24pt>@/_1.6pc/[rrrr]_-{ \Gamma_X }
      &&
      \overset{
        \mathclap{
        \raisebox{3pt}{
          \tiny
          \color{darkblue}
          \bf
          \begin{tabular}{c}
            $\infty$-sheaves
            \\
            on $X$
          \end{tabular}
        }
        }
      }{
        \mbox{\bf{\'E}t}_{X}
      }
      \ar@<+6pt>@{<-}[rr]^-{ \Delta_X }
      \ar@<-6pt>[rr]_-{
        \underset{
          \mathclap{
          \raisebox{-3pt}{
            \tiny
            \color{darkblue}
            \bf
            \begin{tabular}{c}
              form $\infty$-groupoid
              of global sections
            \end{tabular}
          }
          }
        }{
          \Gamma_X
        }
      }^-{ \bot }
      &&
      \mathrm{Groupoids}_\infty
    }
  \end{equation}
  \vspace{-.3cm}

  \noindent
  Notice that the global sections of the $\infty$-sheaf
  of local sections of an $\infty$-bundle $E$ is the
  global sections of that $\infty$-bundle
  (as in Remark \ref{TwistedCohomology}):

  \vspace{-.7cm}
  $$
    \Gamma_X
    \big(
      \underline{E}
    \big)
    \;\simeq\;
    \Gamma_X(E)
  $$
  \vspace{-.4cm}

  \noindent
  (by the essential uniqueness of the
  base geometric morphism (Prop. \ref{BaseGeometricMorphism})
  and the fact that the base geometric morphism on
  $\infty$-bundles forms global sections,
  Remark \ref{TwistedCohomologyAsGlobalSections}).
\end{remark}

\medskip

\noindent {\bf {\'E}tale groupoids.}
\begin{defn}[{\'E}tale groupoid]
  \label{EtaleGroupoids}
  Let $\mathbf{H}$ be an elastic $\infty$-topos (Def. \ref{ElasticInfinityTopos}).

\noindent {\bf (i)}   We say that
$X_\bullet \in \mathrm{Groupoids}(\mathbf{H})$ (Def. \ref{InHigherToposGroupoids})
  is an {\it {\'e}tale groupoid} if all its face maps are local diffeomorphisms
  (Def. \ref{FormallyEtaleMorphism}):
  \vspace{-2mm}
  $$
    \mbox{
      $X_\bullet$
      is {\'etale groupoid}
    }
    \phantom{AAA}
    \Leftrightarrow
    \phantom{AAA}
    \underset{
      {n \in \mathbb{N}}
      \atop
      {0 \leq i \leq n}
    }{\forall}
    \;\;
    \xymatrix{
      X_{n+1}
      \ar[rr]^-{d_i}_-{ \mbox{\tiny{\'e}t} }
      &&
      X_{n}
    }.
  $$

  \vspace{-2mm}
\noindent {\bf (ii)}   We write
  \begin{equation}
    \label{CategoryOfEtaleGroupoids}
    \mbox{\'EtaleGroupoids}(\mathbf{H})
      \xymatrix{
        \;   \ar@{^{(}->}[r]
      &
}
    \mathrm{Groupoids}(\mathbf{H})
    \;\;\;
    \in
    \;
    \mathrm{Categories}_\infty
  \end{equation}
  for the full sub-$\infty$-category of that of all groupoids
  \eqref{InInfinityToposGroupoids} on those that are
  {\'e}tale groupoids.
\end{defn}
As a variant of Prop. \ref{GroupoidsEquivalentToEffectiveEpimorphisms} we have:
\begin{prop}[{\'E}tale groupoids are equivalent to stacks with {\'e}tale atlases]
  \label{EtaleGroupoidsAndEtaleAtlases}
  $\;$\\

  \vspace{-.8cm}

  \hspace{-.9cm}
  \begin{tabular}{llll}

  \vspace{-3cm}

  \begin{minipage}[left]{8cm}
  Let $\mathbf{H}$ be an elastic $\infty$-topos (Def. \ref{ElasticInfinityTopos})
  and $X_\bullet \in \mathrm{Groupoids}(\mathbf{H})$ (Def. \ref{InHigherToposGroupoids}).
  Then the following conditions are equivalent:

 \noindent  {\bf (i)} The groupoid $X_\bullet$ is an {\'e}tale groupoid (Def. \ref{EtaleGroupoids}).

  \noindent  {\bf (ii)} The associated atlas $\xymatrix@R=8pt{X_0 \ar@{->>}[r]^-a & \mathcal{X} }$
  (via Prop. \ref{GroupoidsEquivalentToEffectiveEpimorphisms})
    is a local diffeomorphism (Def. \ref{FormallyEtaleMorphism}).
  \end{minipage}
  &&&
  \begin{minipage}[left]{8cm}
    \begin{equation}
  \label{EtaleAtlasesAndGroupoids}
  \raisebox{40pt}{
  \xymatrix@C=15pt@R=1.5em{
    \ar@<-12pt>@{..>}[d]
    \ar@<-6pt>@{<..}[d]
    \ar@{..>}[d]
    \ar@<+6pt>@{<..}[d]
    \ar@<+12pt>@{..>}[d]
    &
    \ar@<-12pt>@{..>}[d]
    \ar@<-6pt>@{<..}[d]
    \ar@{..>}[d]
    \ar@<+6pt>@{<..}[d]
    \ar@<+12pt>@{..>}[d]
    \\
    X \times_{\mathcal{X}} X
    \ar@{}[r]|-{\simeq}
    \ar@<-6pt>[d]_-{\mbox{\tiny{\'e}t}}
    \ar@{<-}[d]|-{\mbox{\tiny{\'e}t}}
    \ar@<+6pt>[d]^-{\mbox{\tiny{\'e}t}}
    &
    X_1
    \ar@<-6pt>[d]_-{s}
    \ar@{<-}[d]|-{e}
    \ar@<+6pt>[d]^-{t}
    &
    \mbox{
      \footnotesize
      \color{darkblue}
      \bf
      ``{\'e}tale groupoid''
    }
    \\
    X_0
    \ar@{->>}[d]_-{a}^-{ \mbox{\tiny{\'e}t} }
    \ar@{=}[r]
    &
    X_0
    \ar@{->>}[d]
    &
    \ar@{}[d]|-{
      \mbox{
        \footnotesize
        \color{darkblue}
        \bf
        ``{\'e}tale atlas''
      }
    }
    \\
    \mathcal{X}
    \ar@{}[r]|-{\simeq}
    &
    \underset{\longrightarrow}{\mathrm{lim}}\, X_\bullet
    &
    \mbox{
      \footnotesize
      \color{darkblue}
      \bf
      ``{\'e}tale stack''
    }
  }
  }
  \end{equation}
  \end{minipage}
  \end{tabular}
\end{prop}

\vspace{2.5cm}

\begin{proof}
  By definition of local diffeomorphisms, we need to demonstrate
  the logical equivalence shown on the left:
   \vspace{-3mm}
  \begin{equation}
    \label{TowardsEtaleGroupoids}
    \hspace{-3mm}
    \underset{
      n_1 \overset{\phi}{\to} n_2
    }{\forall}
    \;\;\;\;\;
    \raisebox{20pt}{
    \xymatrix@C=4em@R=1.5em{
      X_{n_1}
      \ar[r]^-{ \eta^\Im_{X_{n_1}} }
      \ar[d]_-{X_\phi}
      \ar@{}[dr]|-{\mbox{\tiny\rm(pb)}}
      &
      \Im X_{n_1}
      \ar[d]^-{\Im X_\phi}
      \\
      X_{n_2}
      \ar[r]_-{ \eta^\Im_{X_{n_2}} }
      &
      \Im X_{n_2}
    }
    }
    \;\;\;\;\;
    \Leftrightarrow
    \;\;\;\;\;
    \raisebox{22pt}{
    \xymatrix@C=4em@R=1.5em{
      X_{0}
      \ar[r]^-{ \eta^\Im_{X_0} }
      \ar[d]_-{ a }
      \ar@{}[dr]|-{\mbox{\tiny\rm(pb)}}
      &
      \Im X_0
      \ar[d]^-{ \Im a }
      \\
      \underset{\longrightarrow}{\mathrm{lim}}
      \, X_\bullet
      \ar[r]_-{
        \eta^\Im_{
          \underset{\longrightarrow}{\mathrm{lim}}
          \,
          X_\bullet
        }
      }
      &
      \Im
      \underset{\longrightarrow}{\mathrm{lim}}
      \, X_\bullet
    }
    }
    \;\;\;\;\;
    \Leftrightarrow
    \;\;\;\;\;
    \raisebox{22pt}{
    \xymatrix@C=4em@R=1.5em{
      X_{0}
      \ar[r]^-{ \eta^\Im_{X_0} }
      \ar[d]_-{ a }
      \ar@{}[dr]|-{\mbox{\tiny\rm(pb)}}
      &
      \Im X_0
      \ar[d]^-{ \Im a }
      \\
      \underset{\longrightarrow}{\mathrm{lim}}
      \, X_\bullet
      \ar[r]_-{
        \underset{\longrightarrow}{\mathrm{lim}}\,
        \eta^\Im_{
          X_\bullet
        }
      }
      &
      \underset{\longrightarrow}{\mathrm{lim}}
      \, \Im X_\bullet
    }
    }
  \end{equation}

   \vspace{-2mm}
\noindent
  But since $\Im$ preserves all limits and colimits
  (being a left and a right adjoint, Prop. \ref{AdjointsPreserveCoLimits}), we have
  {\bf (a)} also the logical equivalence shown on the right
  of \eqref{TowardsEtaleGroupoids};
  and {\bf (b)} that $\Im X_\bullet$ is itself a groupoid
  with atlas $\Im a$, and that
  $\xymatrix@C=3em{ X_\bullet \ar[r]|-{\, \eta^\Im_{X_\bullet}}  &  \Im X_\bullet }$
  is a morphism in
  $\mathrm{Groupoids}(\mathbf{H})$ \eqref{InInfinityToposGroupoids}.
  By {\bf (a)}, it is now sufficient to
  prove the composite logical equivalence in \eqref{TowardsEtaleGroupoids}.
  By {\bf (b)}, this follows with Prop. \ref{GroupoidMorphismsCartesian}.
\hfill \end{proof}

\begin{prop}[Tangent stacks]
  \label{TangentGroupoids}
  Let $\mathbf{H}$ be an elastic $\infty$-topos
  (Def. \ref{ElasticInfinityTopos})
  and $X_\bullet \in \mbox{\rm{\'E}taleGroupoids}(\mathbf{H})$
  (Def. \ref{EtaleGroupoids}) with {\'e}tale atlas
  $\xymatrix@C15pt{X \ar[r]^{\mbox{\tiny{\'e}t}} & \mathcal{X} }$
  (via Prop. \ref{EtaleGroupoidsAndEtaleAtlases}).
  Then:

  \noindent {\bf (i)} the system of tangent bundles $T X_\bullet$
  (Def. \ref{InfinitesimalTangentBundle})
  is itself an {\'e}tale groupoid (Def. \ref{EtaleGroupoids}),
  the \emph{tangent groupoid};

  \noindent {\bf (ii)}
  its atlas (under Prop. \ref{EtaleGroupoidsAndEtaleAtlases})
  is the differential $\xymatrix{ T X_0 \ar[r]^{T a} & T \mathcal{X} }$
  of the given atlas,
  hence the \emph{tangent stack} is:
  \vspace{-2mm}
  \begin{equation}
    \label{TangentStack}
    T \mathcal{X}
    \;\simeq\;
    \underset{\longrightarrow}{\mathrm{lim}} T X_\bullet
  \end{equation}
\end{prop}

\begin{proof}
{\bf (i)} That $TX_\bullet$ is itself a groupoid
(Def. \ref{InHigherToposGroupoids}) follows because both
the tangent bundle construction $T(-)$ \eqref{TangentBundlePullbackDefinition}
as well as the groupoid Segal conditions \eqref{GroupoidSegalCondtion}
are pullback constructions, hence limits, which commute over each other.
To see that $T X_\bullet$ is an {\'e}tale groupoid, consider
the following diagram:
\begin{equation}
  \label{TowardsTangentGroupoids}
  \hspace{-2cm}
  \raisebox{40pt}{
  \xymatrix@R=3pt@C=15pt{
    & & &
    \ar@{..>}@<-10pt>[dl]
    \ar@{<..}@<-5pt>[dl]
    \ar@{..>}[dl]
    \ar@{<..}@<+5pt>[dl]
    \ar@{..>}@<+10pt>[dl]
    &&&&
    \ar@{..>}@<-10pt>[dl]
    \ar@{<..}@<-5pt>[dl]
    \ar@{..>}[dl]
    \ar@{<..}@<+5pt>[dl]
    \ar@{..>}@<+10pt>[dl]
    \\
    & &
    T
    \big(
      X_0 \times_{\mathcal{X}} X_0
    \big)
    \ar[dddd]|<<<<{\phantom{AA\vert}}|<<<<<<<<<<<<{\phantom{AA\vert}}
    \ar[rrrr]|<{\phantom{AAA}}
    \ar@<-5pt>[dl]
    \ar@{<-}[dl]
    \ar@<+5pt>[dl]
    &&&&
    X_0 \times_{\mathcal{X}} X_0
    \ar[dddd]
    \ar@<-5pt>[dl]
    \ar@{<-}[dl]
    \ar@<+5pt>[dl]
    \\
    &
    \;T X_0\;
    \ar[rrrr]|<{\phantom{AAAA}}
    \ar[dddd]|<<<<<<{\phantom{AA}}
    \ar[dl]
    &&&&
    \;\;X_0\;\;
    \ar[dddd]
    \ar[dl]
    \\
    \underset{\longrightarrow}{\mathrm{lim}}
    T X_\bullet
    \ar[rrrr]
    \ar@{-->}[dddd]
    &&&&
    \;\;\mathcal{X}\;\;
    \ar@{-->}[dddd]
    \\
    &&&
    \ar@{..>}@<-10pt>[dl]
    \ar@{<..}@<-5pt>[dl]
    \ar@{..>}[dl]
    \ar@{<..}@<+5pt>[dl]
    \ar@{..>}@<+10pt>[dl]
    &&&&
    \ar@{..>}@<-10pt>[dl]
    \ar@{<..}@<-5pt>[dl]
    \ar@{..>}[dl]
    \ar@{<..}@<+5pt>[dl]
    \ar@{..>}@<+10pt>[dl]
    \\
    & &
    X_0 \times_{\mathcal{X}} X_0
    \ar[rrrr]|<{\phantom{AAAAAA}}|<<<<<<<<<<<<<<<<<<<<<{\phantom{AA}}|>>>>>>>>>>>{\phantom{AA}}
     \ar@<-5pt>[dl]
    \ar@{<-}[dl]
    \ar@<+5pt>[dl]
    &&&&
    \Im
    \big(
      X_0 \times_{\mathcal{X}} X_0
    \big)
    \ar@<-5pt>[dl]
    \ar@{<-}[dl]
    \ar@<+5pt>[dl]
    \\
    &
    \;\;X_0\;\;
    \ar[rrrr]|<{\phantom{AAAAA}}|>>>>>>>>>>>{\phantom{AA}}
    \ar@{-->}[dl]
    &&&&
    \;\Im X_0\;
    \ar@{-->}[dl]
    \\
    \;\;\mathcal{X}\;\;
    \ar@{-->}[rrrr]
    &&&&
    \;\Im \mathcal{X}\;
  }
  }
\end{equation}
Here the simplicial sub-diagram in the top right consists
of local diffeomorphism by the assumption that $X_\bullet$
is {\'e}tale. But this implies that all the horizontal squares
in the top of \eqref{TowardsTangentGroupoids} are Cartesian,
by Prop \ref{PullbackAlongLocalDiffeomorphismsPreservesTangentBundles},
hence that also all morphisms of the simplicial sub-diagram in the top
left are local diffeomorphisms, by Lemma \ref{ClosureOfLocalDiffeomorphisms}.

\noindent {\bf (ii)} To see \eqref{TangentStack} we need to show that the
front square in \eqref{TowardsTangentGroupoids} is Cartesian.
Observe:
\begin{itemize}
 \vspace{-2mm}
\item[{\bf (a)}] All horizontal squares in \eqref{TowardsTangentGroupoids}
are Cartesian: the top ones by the above argument for {\bf (i)},
the bottom ones by the assumption that $X_\bullet$ is {\'e}tale.

 \vspace{-3mm}
\item[{\bf (b)}] All solid vertical  squares in \eqref{TowardsTangentGroupoids}
are also Cartesian, by definition \eqref{TangentBundlePullbackDefinition}
of tangent bundles.

 \vspace{-3mm}
\item[{\bf (c)}] The object $\mathcal{X}$ in the bottom front left of
\eqref{TowardsTangentGroupoids} is not just the colimit of the
simplicial sub-diagram in the bottom left, but in fact of the
full left sub-diagram (because of the colimit of the top left sub-diagram
in the front top left). Similarly, the object $\Im \mathcal{X}$
is in fact the colimit over the full right sub-diagram in
\eqref{TowardsTangentGroupoids}
(using that $\Im$ preserves colimits, being a left adjoint, Prop. \ref{AdjointsPreserveCoLimits}).
\end{itemize}

 \vspace{-2mm}
\noindent
Now {\bf (a)} and {\bf (b)} verify the assumption of
Prop. \ref{ColimitsOfEquifiberedTransformations}
applied to the diagram \eqref{TowardsTangentGroupoids}, regarded as a
natural transformation from its left part to its right part;
and with {\bf (c)}, the conclusion of Prop. \ref{ColimitsOfEquifiberedTransformations}
says that the front square in \eqref{TowardsTangentGroupoids} is
Cartesian.
\hfill \end{proof}

\begin{lemma}[Degreewise local diffeomorphisms of {\'e}tale groupoids]
  \label{DegreewiseLocalDiffeomorphismsOfEtaleGroupoids}
  Let $\mathbf{H}$ be an elastic $\infty$-topos (Def. \ref{ElasticInfinityTopos})
  and $X_\bullet, Y_\bullet \in \mbox
  {\'EtaleGroupoids}(\mathbf{H})$ (Def. \ref{EtaleGroupoids}).
   If a morphism
  $\xymatrix@C=25pt{ X_\bullet \ar[r]|{\, f_\bullet } & Y_\bullet  }$
  is such that
  for all $n \in \mathbb{N}$, the component
    $\xymatrix@C=25pt{ X_n \ar[r]|{\, f_n } & Y_n  }$
  is a local diffeomorphism (Def. \ref{FormallyEtaleMorphism}),
  then induced morphism on stacks
    $\xymatrix@C=40pt{ \mathcal{X} \ar[r]|{\;\;\scalebox{.6}{$\underset{\longrightarrow}{\mathrm{lim}}\,f_\bullet$}\;} & \mathcal{Y} }$
    is also a local diffeomorphism (Def. \ref{EtaleGroupoidsAndEtaleAtlases}).
\end{lemma}
\begin{proof}
  Consider the following diagram:
 \vspace{-7mm}
  $$
  \hspace{2cm}
    \xymatrix@R=.8em{
      \ar@{..>}@<-10pt>[d]
      \ar@{<..}@<-5pt>[d]
      \ar@{..>}[d]
      \ar@{<..}@<+5pt>[d]
      \ar@{..>}@<+10pt>[d]
      &&
      \ar@{..>}@<-10pt>[d]
      \ar@{<..}@<-5pt>[d]
      \ar@{..>}[d]
      \ar@{<..}@<+5pt>[d]
      \ar@{..>}@<+10pt>[d]
      \\
      X_1
      \ar@<-5pt>[dd]
      \ar@{<-}[dd]
      \ar@<+5pt>[dd]
      \ar[rr]
      \ar[dr]|-{ \eta^\Im_{X_1} }
      &
      \ar@{..>}@<-10pt>[d]
      \ar@{<..}@<-5pt>[d]
      \ar@{..>}[d]
      \ar@{<..}@<+5pt>[d]
      \ar@{..>}@<+10pt>[d]
      &
      Y_1
      \ar@<-5pt>[dd]|-{\phantom{AA}}
      \ar@{<-}[dd]|-{\phantom{AA}}
      \ar@<+5pt>[dd]|-{\phantom{AA}}
      \ar[dr]|-{ \eta^\Im_{Y_1} }
      &
      \ar@{..>}@<-10pt>[d]
      \ar@{<..}@<-5pt>[d]
      \ar@{..>}[d]
      \ar@{<..}@<+5pt>[d]
      \ar@{..>}@<+10pt>[d]
      \\
      &
      \Im X_1
      \ar@<-5pt>[dd]
      \ar@{<-}[dd]
      \ar@<+5pt>[dd]
      \ar[rr]
      &&
      \Im Y_1
      \ar@<-5pt>[dd]
      \ar@{<-}[dd]
      \ar@<+5pt>[dd]
      \\
      X_0
      \ar[rr]|-{\phantom{AAAA}}
      \ar[dd]
      \ar[dr]|-{ \eta^\Im_{X_0} }
      &&
      Y_0
      \ar@{-->}[dd]|-{\phantom{AA}}
      \ar[dr]|-{ \eta^\Im_{Y_0} }
      \\
      &
      \Im X_0
      \ar[dd]
      \ar[rr]
      &&
      \Im Y_0
      \ar@{-->}[dd]
      \\
      \mathcal{X}
      \ar@{-->}[rr]|-{\phantom{AA}}
      \ar[dr]|-{ \eta^\Im_{\mathcal{X}} }
      &&
      \mathcal{Y}
      \ar@{-->}[dr]|-{ \eta^\Im_{\mathcal{Y}} }
      \\
      &
      \Im \mathcal{X}
      \ar@{-->}[rr]
      &&
      \Im \mathcal{Y}
    }
  $$

  \vspace{-1mm}
\noindent
  Observe that:
  \begin{itemize}
  \vspace{-2mm}
\item[{\bf (a)}] all solid $\eta^\Im$-naturality squares in this diagram
  are Cartesian, by the assumption that the rear part of the diagram
  is a degreewise local diffeomorphism of {\'e}tale groupoids.
    \vspace{-3mm}
 \item[{\bf (b)}] $\mathcal{Y}$ is not just the colimit
  of the partial diagram $Y_\bullet$ in the rear right, but in
  fact is also the colimit of the full non-dashed rear part of the
  diagram (using that $\mathcal{X}$ is the colimit of the rear left part).
  Similarly, $\Im \mathcal{Y}$ is the colimit of the non-dashed front
  part of the diagram
  (using that $\Im$ preserves limits and colimits, being a left and
  a right adjoint, Prop. \ref{AdjointsPreserveCoLimits}).
  \end{itemize}
    \vspace{-2mm}
 \noindent  Hence if we regard the diagram as a natural transformation
  from its rear to its front part, then Prop. \ref{ColimitsOfEquifiberedTransformations}
  applies and says that also the bottom dashed square is Cartesian,
  and hence that $\mathcal{X} \to \mathcal{Y}$ is a local diffeomorphism.
\hfill \end{proof}

\begin{defn}[{\'E}talification of groupoids]
\label{EtalificationOfGroupoids}
Let $\mathbf{H}$ be an elastic $\infty$-topos (Def. \ref{ElasticInfinityTopos})
and $X_\bullet \in \mathrm{Groupoids}(\mathbf{H})$ (Def. \ref{InHigherToposGroupoids}).
Notice that, by Prop. \ref{GroupoidsEquivalentToEffectiveEpimorphisms}
for all $n \in \mathbb{N}$ we have
for all $0 \leq i \leq n$ that
all face maps $\xymatrix{ X_{n+1} \ar[r]|-{\;d_i\,} & X_n }$
are in fact equivalent to each other,
being related by an automorphism of
$X_{n+1}$ given by permutation of fiber product factors
\eqref{NerveOfAnAtlas}
 \vspace{-1mm}
\begin{equation}
  \label{GroupoidAsDiagramOverDegree0Space}
  X_\bullet
  \;\;\;\simeq\;\;\;
  \raisebox{30pt}{
  \xymatrix@R=9pt@C=9pt{
    &&& &&&
    X_2
    \ar[dr]^-{\simeq}
    &&&&
    \\
    &&&
    X_1
    \ar[dr]^-{\simeq}
    &&&&
    X_2
    \ar[dr]^-{\simeq}
    &&&&
    \\
    X_0
    \ar@{<-}[rrrr]|-{\;d_0}
    \ar@{<-}[urrr]|-{\;d_1}
    &&
    &&
    X_1
    \ar@{<-}[rrrr]|-{\;d_0}
    \ar@{<-}[urrr]|-{\;d_1}
    \ar@{<-}[uurr]|-{\;d_2}
    && &&
    X_2
    \ar@{<..}[rrrr]
    \ar@{<..}[urrr]
    \ar@{<..}[uurr]
    &&&&
  }
  }
\end{equation}
(and similarly for the degeneracy maps).
Therefore, we may regard $X_\bullet$
as a diagram in the slice $\mathbf{H}_{X_0}$.
and apply $\mathcal{L}_{X_0}$ \eqref{ElateOverXModality}
to this diagram \eqref{GroupoidAsDiagramOverDegree0Space}
to obtain
 \vspace{-3mm}
\begin{equation}
  \label{TowardsConstructingTheEtalifiedAtlas}
  X_\bullet^{\mbox{\tiny{\'et}}}
  \;\;\;\simeq\;\;\;
  \raisebox{30pt}{
  \xymatrix@R=9pt@C=9pt{
    &&& &&&
    \mathcal{L}_{X_0} X_2
    \ar[dr]^-{\simeq}
    &&&&
    \\
    &&&
    \mathcal{L}_{X_0} X_1
    \ar[dr]^-{\simeq}
    &&&&
    \mathcal{L}_{X_0} X_2
    \ar[dr]^-{\simeq}
    &&&&
    \\
    X_0
    \ar@{<-}[rrrr]|-{\;\mbox{\tiny{\'e}t}\,}
    \ar@{<-}[urrr]|-{\;\mbox{\tiny{\'e}t}\,}
    &&
    &&
    \mathcal{L}_{X_0} X_1
    \ar@{<-}[rrrr]|-{\;\mbox{\tiny{\'e}t}\,}
    \ar@{<-}[urrr]|-{\;\mbox{\tiny{\'e}t}\,}
    \ar@{<-}[uurr]|-{\;\mbox{\tiny{\'e}t}\,}
    && &&
    \mathcal{L}_{X_0} X_2
    \ar@{<..}[rrrr]
    \ar@{<..}[urrr]
    \ar@{<..}[uurr]
    &&&&
  }
  }
\end{equation}

\vspace{-2mm}
\noindent
Observe that:
 \begin{itemize}
  \vspace{-3mm}
\item[{\bf (a)}]
 the simplicial diagram \eqref{TowardsConstructingTheEtalifiedAtlas}
is again a groupoid, since the right adjoint functor
$\mathcal{L}_{X_0}$ preserves the characterizing fiber products \eqref{GroupoidSegalCondtion}
(by Prop. \ref{AdjointsPreserveCoLimits});

  \vspace{-3mm}
\item[{\bf (b)}] this groupoid is {\'e}tale (Def. \ref{EtaleGroupoids}),
since the morphisms of the form
$\mathcal{L}_{X_0} X_n \longrightarrow X_0$
in \eqref{TowardsConstructingTheEtalifiedAtlas}
are local diffeomorphisms by construction, whence all other morphisms
$\mathcal{L}_{X_0} X_{n_1} \longrightarrow \mathcal{L}_{X_0} X_{n_2}$
are local diffeomorphisms by the left-cancellation property
\eqref{LocalDiffeoComposition}.
\end{itemize}

\vspace{-3mm}
\noindent
 Hence we say that:
 \begin{itemize}
  \vspace{-2mm}
\item[{\bf (i)}] The simplicial diagram \eqref{TowardsConstructingTheEtalifiedAtlas}
 is the \emph{\'etalification} of the groupoid $X_\bullet$.
 \vspace{-2mm}
\begin{equation}
  \label{EtalifiedGroupoid}
  X_\bullet^{\mbox{\tiny{\'et}}}
  \;\in\;
  \mbox{\'EtaleGroupoids}(\mathbf{H})\;.
\end{equation}

 \vspace{-4mm}
\item[{\bf (ii)}]  If the corresponding atlas of $X_\bullet$
(via Prop. \ref{GroupoidsEquivalentToEffectiveEpimorphisms}) is denoted
$\!
  \xymatrix{ X_0 \ar@{->>}[r] & \mathcal{X} }
 \!\!
$,
then we write
 \vspace{-3mm}
\begin{equation}
  \label{AtlasOfEtalifiedGroupoid}
  \xymatrix{
    X_0
    \ar@{->>}[r]^-{ \mbox{\tiny{\'et}} }
    &
    \mathcal{X}^{\scalebox{.6}{\'et}}
  }
\end{equation}

 \vspace{-3mm}
\noindent
for the corresponding {\'e}tale atlas (via Prop. \ref{EtaleGroupoidsAndEtaleAtlases})
of the {\'e}talified groupoid \eqref{EtalifiedGroupoid}.
\end{itemize}
\end{defn}

\subsubsection{Super-geometry}
\label{SuperGeometry}

We present a formulation of super-geometry
internal to $\infty$-toposes which we call \emph{solid} \cite{dcct}.

\medskip

\noindent {\bf Super-geometry.}
\begin{defn}[Solid $\infty$-topos]
  \label{SolidTopos}

 \vspace{-1mm}
  \item {\bf (i)} An $\infty$-topos $\mathbf{H}$ (Def. \ref{InfinityTopos})
  over $\mathbf{B} = \mathrm{Groupoids}_\infty$
  is a \emph{solid $\infty$-topos} if
  its base geometric morphism (Prop. \ref{AdjointInfinityFunctors}),
  to be called
  $\mathrm{Pnts} :
  \xymatrix@C=11pt{\mathbf{H} \ar[r] & \mathbf{B}}$, is equipped
  with a factorization as follows, with adjoints (Def. \ref{AdjointInfinityFunctors})
  as shown:
  \begin{equation}
    \label{SolidityAdjunctions}
    \hspace{2.2cm}
    \xymatrix@C=3.2em{
      \ar@{}@<+14pt>[r]|-{ \mathrm{Shp} \; : }
      \ar@{}@<-14pt>[r]|-{ \Gamma \; : }
      &
      \mathbf{H}
 \;\;        \ar@{->}@<+42pt>[rr]|<\times|-{
          \mathllap{
            \scalebox{.8}{
              \color{darkblue} \bf
              ``even''
              \hspace{4cm}
            }
          }
         \; \mathrm{Evn}  \;
        }_-{\raisebox{-6pt}{\;\tiny $\bot$}\;}
        \ar@{<-^{)}}@<+28pt>[rr]|-{
          \mathllap{
            \scalebox{.8}{
              \color{darkblue} \bf
              ``bosonic''
              \hspace{4cm}
            }
          }
          \; \mathrm{Bsnc} \;
        }_-{\raisebox{-6pt}{\tiny $\bot$}}
        \ar@{->}@<+14pt>[rr]|-{
          \mathllap{
            \scalebox{.8}{
              \color{darkblue} \bf
              ``super shape''
              \hspace{3.9cm}
            }
          }
          \mathrm{Shp}_{\mathrm{sup}}
        }_-{\raisebox{-6pt}{\tiny $\bot$}}
        \ar@{<-^{)}}@<+0pt>[rr]|-{
          \mathllap{
            \scalebox{.8}{
              \color{darkblue} \bf
              ``super discrete''
              \hspace{3.8cm}
            }
          }
          \; \mathrm{Disc}_{\mathrm{sup} \;}
        }
        \ar@<+0pt>@<-14pt>[rrrr]|>>>>>>>>>>>>>>>>{
          \; \mathrm{Pnts}_{\mathrm{inf}} \;
        }^-{\raisebox{+1pt}{\tiny $\bot$}}
        \ar@{<-^{)}}@<-28pt>[rrrrrr]|-{ \; \mathrm{Chtc} \;}^-{\raisebox{+2pt}{\tiny $\bot$}}
      &&
      \; \mathbf{H}_{\rightsquigarrow}
        \ar@{<-^{)}}@<+28pt>[rr]|-{
        \;  \mathrm{Rdcd} \;
        }_-{\raisebox{-6pt}{\tiny $\bot$}}
        \ar@{->}@<+14pt>[rr]|-{ \;\mathrm{Shp}_{\mathrm{inf}} \;}_-{\raisebox{-6pt}{\tiny $\bot$}}
        \ar@{<-^{)}}@<+0pt>[rr]|-{ \;\mathrm{Disc}_{\mathrm{inf}} \;}
      &
      \ar@<-14pt>@{}[rr]|-{\raisebox{-.6pt}{$---$}}
      &
      \; \mathbf{H}_{\Re}
        \ar@{->}@<+14pt>[rr]|<\times|-{ \;\mathrm{Shp}_\Re \; }_-{\raisebox{-6pt}{\tiny $\bot$}}
        \ar@{<-^{)}}@<+0pt>[rr]|-{ \;\mathrm{Disc}_{\Re} \; }_-{\raisebox{-6pt}{\tiny $\bot$}}
        \ar@<+0pt>@<-14pt>[rr]|-{ \; \mathrm{Pnts}_\Re  \;}
      &&
     \;\;  \mathbf{B}
      \ar@{}@<+0pt>[r]|-{ : \; \mathrm{Disc} }
      &
      \\
      &
      \mathclap{
        \mbox{
          \tiny
          \color{darkblue} \bf
          \begin{tabular}{c}
            solid
            \\
            $\infty$-topos
          \end{tabular}
        }
      }
      &&
      \mathclap{
        \mbox{
          \tiny
          \color{darkblue} \bf
          \begin{tabular}{c}
            bosonic
            \\
            sub-topos
          \end{tabular}
        }
      }
      &&
      \mathclap{
        \mbox{
          \tiny
          \color{darkblue} \bf
          \begin{tabular}{c}
            reduced
            \\
            sub-topos
          \end{tabular}
        }
      }
      &&
      \mathclap{\;\;\;\;\;\;
        \mbox{
          \tiny
          \color{darkblue} \bf
          \begin{tabular}{c}
            discrete
            \\
            sub-topos
          \end{tabular}
        }
      }
    }
  \end{equation}
  \item {\bf (ii)}  In particular, a solid $\infty$-topos is also
  an elastic $\infty$-topos (Def. \ref{ElasticInfinityTopos}),
  as its is sub-$\infty$-topos $\mathbf{H}_{\rightsquigarrow}$
  of bosonic objects.

    \item {\bf (iii)} We write
  \vspace{-2mm}
  \begin{equation}
    \label{ModeSupergeometric}
    \underset{
      \mbox{
        \footnotesize
        \color{darkblue} \bf
        ``even''
      }
    }{
    \big(
      \;
      \rightrightarrows
      \;\;:=\; \mathrm{Bsn} \circ \mathrm{Evn}
      \;
    \big)
    }
    \;\dashv\;
    \underset{
      \mbox{
        \footnotesize
        \color{darkblue} \bf
        ``bosonic''
      }
    }{
    \big(
      \;
      \rightsquigarrow
      \;\;:=\; \mathrm{Bsn} \circ \mathrm{Shp}_{\mathrm{sup}}
      \;
    \big)
    }
    \;\dashv\;
    \underset{
      \mbox{
        \footnotesize
        \color{darkblue} \bf
        ``rheonomic''
      }
    }{
    \big(
      \;
      \mathrm{R}\!\mathrm{h} \;:=\; \mathrm{Disc}_{\mathrm{sup}} \circ \mathrm{Shp}_{\mathrm{sup}}
      \;
    \big)
    }
    \;:\;
    \mathbf{H} \longrightarrow \mathbf{H}
  \end{equation}
  for the further induced modalities \eqref{AdjointModalities}
  (\emph{solid modalities})
  accompanying the elastic modalities \eqref{ElasticityModalities}
  and the cohesive modalities \eqref{CohesiveModalitiesFromAdjointQuadruple}.
\end{defn}

\medskip

\noindent{\bf Examples of solid $\infty$-toposes.}
We indicate an example of a solid $\infty$-topos (Def. \ref{SolidTopos}).
For full details on the construction see \cite{SS20}.
In generalization of Def. \ref{FormalCartesianSpaces}
we have the following:
\begin{defn}[$\infty$-Jets of super Cartesian spaces]
  \label{SuperFormalCartesianSpaces}
  $\phantom{A}$
   \vspace{-2mm}
 \item {\bf (i)}  Write
  \vspace{-2mm}
  \begin{equation}
    \label{SupcategorySuperCartesianSpacesInsideSuperAlgebras}
    \xymatrix@R=-2pt{
      \infty\mathrm{JetsOfSuperCartesianSpaces}
      \;\ar@{^{(}->}[rr]^-{C^\infty(-)}
      &&
      \mathrm{CommutativeSuperAlgebras}_{\mathbb{R}}^{\mathrm{op}}
      \\
      \mathbb{R}^{n\vert q} \times \mathbb{D}_W
      \ar@{|->}[rr]
      &&
      C^\infty\big(\mathbb{R}^n\big)
        \;\otimes_{\mathbb{R}}\;
      \wedge_{\mathbb{R}}^\bullet\big(\mathbb{R}^q\big)
        \;\otimes_{\mathbb{R}}\;
      (\mathbb{R} \oplus \mathrm{W})
    }
  \end{equation}
  for (as in \cite{KonechnySchwarz97}\cite{KonechnySchwarz00})
  the full subcategory of the opposite of super-commutative
  super-algebras over the real numbers on those which are tensor products of
  \begin{itemize}
  \vspace{-2mm}
      \item[{\bf (a)}]   algebras $C^\infty(\mathbb{R}^n)$ of smooth functions on
      a Cartesian space $\mathbb{R}^n$, for $d \in \mathbb{N}$;
    \vspace{-2mm}
    \item[{\bf (b)}]   Grassmann algebras $\wedge^\bullet_{\mathbb{R}} \mathbb{R}^q$
    on $q \in \mathbb{N}$ generators in odd degree;
   \vspace{-2mm}
    \item[{\bf (c)}] finite dimensional
    $\mathbb{R} \oplus W \;\in\; \mathrm{CommutativeAlgebras}$
    with a single nilpotent maximal ideal $W$.
\end{itemize}

 \vspace{-3mm}
\item {\bf (ii)}
  We regard this as a site via the
  the coverage (i.e., a Grothendieck pre-topology) whose covers are of the form

  \vspace{-.5cm}

  $$
    \Big\{ \!\!
      \xymatrix@C=15pt{
        \underset{
          \mathbb{R}^{n \vert q}
        }{
          \underbrace{
            \mathbb{R}^{n} \times \mathbb{R}^{0\vert q}
          }
        }
        \times \mathbb{D}
        \ar[rr]^-{f_i \times \mathrm{id} \times \mathrm{id}}
        &&
        \mathbb{R}^n \times \mathbb{R}^{0\vert q} \times \mathbb{D}
      }
    \!\! \Big\}_{i \in I}
     \;\;  \mbox{such that
     $
     \;
     \Big\{ \!\!
       \xymatrix{
         \mathbb{R}^n
         \ar[r]^{f_i}
         &
         \mathbb{R}^n
       }
    \!\!  \Big\}_{i \in I}
     \!\!\!\!
     $
     \begin{tabular}{l}
       is a cover in
       $\mathrm{CartesianSpaces}$ (Def. \ref{CartesianSpaces}).
     \end{tabular}
     }
  $$
\end{defn}

\begin{lemma}[Reflections of super-commutative algebras into commutative algebras]
  \label{ReflectionsOfSuperCommutativeAlgebras}
  The canonical inclusion of $\infty\mathrm{JetsOfCartesianSpaces}$
  (Def. \ref{FormalCartesianSpaces})
  into $\infty\mathrm{JetsOfSuperCartesianSpaces}$
  (Def. \ref{SuperFormalCartesianSpaces}) has a left and a right
  adjoint (Def. \ref{AdjointInfinityFunctors})
  \begin{equation}
    \label{AdjunctionBetweenSuperCartesianAndJetsOfCartesian}
    \xymatrix{
      \infty\mathrm{JetsOfSuperCartesianSpaces}
      \;\;
      \ar@<+16pt>[rr]|-{ \;\mathrm{Evn} \;}
      \ar@{<-^{)}}[rr]|-{ \;\mathrm{Bsnc} \;}^-{
        \raisebox{2pt}{
          \scalebox{.8}{$
            \bot
          $}
        }
      }_-{
        \raisebox{-2pt}{
          \scalebox{.8}{$
            \bot
          $}
        }
      }
      \ar@<-16pt>[rr]|-{\; \mathrm{Shp}_{\mathrm{sup}}\; }
      &&
    \;\;  \infty\mathrm{JetsOfCartesianSpaces}
    }
  \end{equation}
  where:

  \noindent {\bf (i)} The left adjoint $\mathrm{Evn}$
  in \eqref{AdjunctionBetweenSuperCartesianAndJetsOfCartesian} is given
  in terms of
  super-algebras of smooth functions \eqref{SupcategorySuperCartesianSpacesInsideSuperAlgebras}
  by passage to the sub-algebra of even-graded elements:
  \begin{equation}
    \begin{aligned}
    C^\infty
    \left(
      \mathrm{Evn}
      \big(
        \mathbb{R}^{n\vert q} \times \mathbb{D}
      \big)
    \right)
    & \simeq\;
    C^\infty\big(
      \mathbb{R}^{n\vert q} \times \mathbb{D}
    \big)_{\mathrm{even}}
    \\
    & \simeq\;
    C^\infty
    \big(
      \mathbb{R}^n \times \mathbb{D}
    \big)
      \otimes_{\mathbb{R}}
    C^\infty
    \big(
      \mathbb{R}^{0\vert q}
    \big)_{\mathrm{even}}\;.
    \end{aligned}
  \end{equation}

  \noindent {\bf (ii)} The right adjoint
  $\mathrm{Shp}_{\mathrm{sup}}$
  in \eqref{AdjunctionBetweenSuperCartesianAndJetsOfCartesian}
  is given in terms of
  super-algebras of smooth functions \eqref{SupcategorySuperCartesianSpacesInsideSuperAlgebras}
  by passage to the quotient algebra by the ideal of odd-graded elements:
  \begin{equation}
    \label{SuperShapeInTermsOfFunctionAlgebras}
    \begin{aligned}
    C^\infty
    \Big(
      \mathrm{Shp}_{\mathrm{sup}}
      \big(
        \mathbb{R}^{n\vert q} \times \mathbb{D}
      \big)
    \Big)
    & \simeq\;
    C^\infty\big(
      \mathbb{R}^{n\vert q} \times \mathbb{D}
    \big) / C^\infty\big(
      \mathbb{R}^{n\vert q} \times \mathbb{D}
    \big)_{\mathrm{odd}}
    \\
    & \simeq\;
    C^\infty
    \big(
      \mathbb{R}^n \times \mathbb{D}
    \big)
    \otimes_{\mathbb{R}}
    \underset{
      \simeq
      \,
      \mathbb{R}
    }{
      \underbrace{
        C^\infty
        \big(
          \mathbb{R}^{0 \vert q}
        \big) /
        C^\infty
        \big(
          \mathbb{R}^{0 \vert q}
        \big)_{\mathrm{odd}}
      }
    }
    \\
    & \simeq
    C^\infty
    \big(
      \mathbb{R}^n \times \mathbb{D}
    \big)
  \end{aligned}
\end{equation}
and hence directly by
\begin{equation}
  \mathrm{Shp}_{\mathrm{sup}}
  \big(
    \mathbb{R}^{n\vert q}
    \times
    \mathbb{D}
  \big)
  \;\simeq\;
  \mathbb{R}^{n} \times \mathbb{D}
  \,.
\end{equation}
\end{lemma}
\begin{proof}
  By regarding the situation under the defining
  embedding as being in $\mathrm{CommutativeSuperAlgebras}_{\mathbb{R}}$
  (Def. \ref{SuperFormalCartesianSpaces}), it is equivalent
  to the statement that the canonical inclusion of commutative algebras
  into super-commutative super-algebras
  has a right and a left adjoint given by passage to the even
  sub-algebra and to the quotient by the odd ideal, respectively:
  \vspace{-2MM}
  $$
    \xymatrix@C=4em{
      \mathrm{CommutativeSuperAlgebras}_{\mathbb{R}}
      \;
      \ar@<+16pt>[rr]|-{ \; A \;\mapsto\; A/A_{\mathrm{odd}}  \;}
      \ar@{<-^{)}}[rr]|-{  }^-{
        \raisebox{2pt}{
          \scalebox{.8}{$
            \bot
          $}
        }
      }_-{
        \raisebox{-2pt}{
          \scalebox{.8}{$
            \bot
          $}
        }
      }
      \ar@<-16pt>[rr]|-{\; A \;\mapsto\; \mathrm{A}_{\mathrm{even}} \;}
      &&
    \;\;  \mathrm{CommutativeAlgebras}_{\mathbb{R}}\;.
    }
  $$
  This follows readily by inspection
  from the fact that
  homomorphisms of super-algebras preserve super-degree, by definition.
  One place where this adjoint triple has been made explicit before is
  \cite[below Example 3.18]{CarchediRoytenberg12}.
\hfill \end{proof}

\begin{example}[Jets of super-geometric $\infty$-groupoids]
  \label{JetsOfSupergeometricGroupoids}
  The
  $\infty$-category of $\infty$-sheaves (Def. \ref{InfinityToposOfInfinitySheaves})
  $$
    \infty\mathrm{JetsOfSupergeometricGroupoids}_\infty
    \;:=\;
    \mathrm{Sheaves}_\infty
    \big(
      \infty\mathrm{JetsOfSuperCartesianSpaces}
    \big)
  $$
  over the site from Def. \ref{SuperFormalCartesianSpaces}
  is a solid $\infty$-topos (Def. \ref{SolidTopos}).

  \noindent {\bf (i)}
  Its bosonic \eqref{ModeSupergeometric} sub-$\infty$-topos is that of
  $\infty\mathrm{JetsOfSmoothGroupoids}$ (Example \ref{FormalSmoothInfinityGroupoids})
  and its reduced \eqref{ElasticityAdjunctions} sub-$\infty$-topos
  that of $\mathrm{SmoothGroupoids}_\infty$ (Example \ref{SmoothInfinityGroupoids}):
  $$
  \hspace{-2mm}
    \scalebox{.87}{
    \xymatrix@C=46pt@R=9pt{
      \infty\mathrm{JetsOfSupergeometricGroupoids}_\infty
      \ar@<+48pt>[r]|-{\; \mathrm{Evn}\; }
      \ar@<+30pt>@{<-^{)}}[r]|-{\; \mathrm{Bsnc} \;}^-{
        \raisebox{2pt}{
          \scalebox{.7}{$
            \bot
          $}
        }
      }
      \ar@<+15pt>[r]|-{\; \mathrm{Shp}_{\mathrm{sup}} \;}^-{
        \raisebox{2pt}{
          \scalebox{.7}{$
            \bot
          $}
        }
      }
      \ar@{<-^{)}}[r]|-{\;
        \mathrm{Disc}_{\mathrm{sup}} \;
      }^-{
        \raisebox{2pt}{
          \scalebox{.7}{$
            \bot
          $}
        }
      }
      &
      \;
      \infty\mathrm{JetsOfSmoothGroupoids}_\infty
      \ar@{<-^{)}}[r]^-{
        \mathrm{Disc}_{\mathrm{inf}}
      }
      &
      \;
      \mathrm{SmoothGroupoids}_\infty
      \ar@{<-^{)}}[r]^-{
        \mathrm{Disc}
      }
      &
      \;
      \mathrm{Gropoids}_\infty
      \\
      &
      \vdots
      \ar@{^{(}->}[u]
      \\
      &
      2\mathrm{JetsOfSmoothGroupoids}_\infty
      \ar@{^{(}->}[u]
      \\
      &
      \mathrm{JetsOfSmoothGroupoids}_\infty
      \ar@{^{(}->}[u]
      \ar@/_2pc/@{<-_{)}}[uuur]
    }
    }
  $$
  where the adjoint triple
  $\big(\mathrm{Evn} \dashv \mathrm{Bsnc} \dashv \mathrm{Shp}_{\mathrm{sup}}\big)$
  arises by left Kan extension from
  that of Lemma \ref{ReflectionsOfSuperCommutativeAlgebras}.

  \noindent {\bf (ii)} The full inclusion of $\mathrm{SmoothManifolds}$, inherited from
  \eqref{ConcreteSmoothInfinityGroupoidsAreDiffeologicalSpaces},
  extends to a full inclusion of super-manifolds
  (as in \cite[4.6]{CCF11}\cite[2]{HKST11}):

  \vspace{-.4cm}

  \begin{equation}
    \label{EmbeddingOfSupermanifoldsInSuperGroupoids}
    \scalebox{.95}{\xymatrix{
      \mathrm{SmoothManifolds}
      \;
      \ar@{^{(}->}[rr]^-{ \mathrm{Disc}_{\mathrm{sup}} }
      &&
      \;
      \mathrm{SuperManifolds}
      \;
      \ar@{^{(}->}[rr]
      &&
      \;
      \infty\mathrm{JetsOfSupergeometricGroupoids}_\infty
    }
    }
  \end{equation}

  \noindent {\bf (iii)} Accordingly, super-Lie groups
  (e.g. \cite{Yagi93}\cite[7]{CCF11}) embed
  faithfully into all group objects (Prop. \ref{LoopingAndDelooping}):
  \begin{equation}
    \label{EmbeddingOfSuperLieGroupsInGroupsOfJetsOfSuperGroupoids}
    \hspace{-4mm}
   \scalebox{.95}{ \xymatrix@C=.9em{
      \underset{
        \raisebox{-1pt}{
          \tiny
          \color{darkblue}
          \bf
          Lie groups
        }
      }{
        \mathrm{Groups}
        \big(
          \mathrm{SmoothManifolds}
        \big)
      }
      \;
      \ar@{^{(}->}[rr]^-{ \mathrm{Disc}_{\mathrm{sup}} }
      &&
      \;
      \underset{
        \raisebox{-1pt}{
          \tiny
          \color{darkblue}
          \bf
          super Lie groups
        }
      }{
        \mathrm{Groups}
        \big(
          \mathrm{SuperManifolds}
        \big)
      }
      \;
      \ar@{^{(}->}[r]
      &
      \;
      \mathrm{Groups}
      \big(
        \infty\mathrm{JetsOfSupergeometricGroupoids}_\infty
      \big)
    }
    }
  \end{equation}

  \noindent {\bf (iv)}
  In particular, for $d \in \mathbb{N}$
  and $\mathbf{N} \in \mathrm{Spin}(d,1)\mathrm{Representations}_{\mathbb{R}}$,
  the corresponding \emph{supersymmetry} groups, i.e., the
  \emph{super-Poincar{\'e} group}
  and its underlying translational \emph{super-Minkowski group}
  (e.g. \cite[\S 3]{Freed99}) are group objects
  \begin{equation}
    \label{SupersymmetryGroupsInSupergeometricGroupoids}
    \xymatrix{
      \underset{
        \mathclap{
        \raisebox{-6pt}{
          \tiny
          \color{darkblue}
          \bf
          \begin{tabular}{c}
            super-Minkowski
            \\
            super Lie group
          \end{tabular}
        }
        }
      }{
        \mathbb{R}^{d,1\vert \mathbf{N}}
      }
      \;\;\;\;\;
      \ar@{^{(}->}[r]
      &
      \underset{
        \mathclap{
        \raisebox{-2pt}{
          \tiny
          \color{darkblue}
          \bf
          \begin{tabular}{c}
            super-Poincar{\'e}
            \\
            super Lie group
          \end{tabular}
        }
        }
      }{\;
        \mathrm{Iso}
        \big(
          \mathbb{R}^{d,1\vert \mathbf{N}}
        \big)
      }
      \ar@{->>}[r]
      &
      \mathrm{Spin}(d,1)
    \;\;\;
    \in
    \mathrm{Groups}
    \big(
      \infty\mathrm{JetsOfSupergeometricGroupoids}
    \big).
    }
  \end{equation}

\end{example}
\begin{remark}[Superspace cohomology theory in solid $\infty$-toposes]
  \label{SuperCohomologyTheory}
  The intrinsic cohomology \eqref{IntrinsicCohomologyOfAnInfinityTopos}
  in the solid $\infty$-topos of
  $\infty\mathrm{JetsOfSupergeometricGroupoids}_\infty$
  (Example \ref{JetsOfSupergeometricGroupoids})

   \vspace{-1.5mm}
  \item {\bf (i)} includes
  the
  super-rational cohomology of super-Minkowski spacetimes
  \eqref{SupersymmetryGroupsInSupergeometricGroupoids}
  that governs the fundamental ($\kappa$-symmetric)
  super $p$-brane sigma-models of string/M-theory
  \cite{FSS13b}\cite{FSS16a}\cite{FSS16b}, review in \cite{FSS19}.

 \vspace{-1.5mm}
 \item {\bf (ii)}  Its enhancement to
 \emph{twisted} super-rational cohomology
  of super-Minkowski spacetimes \eqref{SupersymmetryGroupsInSupergeometricGroupoids},
  which happens (by Remark \ref{AbelianTwistedCohomology})
  in the intrinsic cohomology of the
  tangent $\infty$-topos
  $T\big(\infty\mathrm{JetsOfSupergeometricGroupoids}_\infty\big)$ (Example \ref{TangentInfinityTopos}),
  encodes the double dimensional reduction from fundamental M-branes to D-branes
  \cite{BSS18}.

  \vspace{-1.5mm}
\item {\bf (iii)}   Its enhancement to \emph{proper equivariant} super-rational cohomology
  of super-Minkowski spacetimes \eqref{SupersymmetryGroupsInSupergeometricGroupoids},
  which happens (by Remark \ref{ProperEquivariantCohomologyTheory} and Theorem \ref{OrbifoldCohomologyEquivariant} below)
  in the intrinsic cohomology
  of the singular-solid $\infty$-topos
  $\mathrm{Singular}\infty\mathrm{JetsOfSupergeometricGroupoids}_\infty$
  (Example \ref{ExamplesOfSingularCohesiveToposes} below),
  encodes also the black (solitonic) super $p$-branes
  \cite{HSS18}.
\end{remark}

\begin{lemma}[In super-geometric groupoids {\'e}tale implies bosonic]
  \label{InSuperGeometricGroupoidsImModalObjectsAreBosonic}
  $\,$\\

  \vspace{-.5cm}

  \noindent
  In the solid $\infty$-topos of
  $\infty\mathrm{JetsOfSupergeometricGroupoids}$ (Example \ref{SuperFormalCartesianSpaces})
  we have a natural equivalence
   \vspace{-1mm}
  \begin{equation}
    \label{ImModalIsBosonic}
    \rightsquigarrow
    \circ
    \,
    \Im
    \;\;\simeq\;\;
    \Im
  \end{equation}
  saying that $\Im$-modal objects \eqref{ElasticityModalities}
  are bosonic \eqref{ModeSupergeometric}.
\end{lemma}
\begin{proof}
  Observe that on
  $\infty\mathrm{JetsOfSuperCartesianSpaces}
  \overset{y}{\longhookrightarrow} \infty\mathrm{JetsOfSupergeometricGroupoids}_\infty$
  (Def. \ref{SuperFormalCartesianSpaces}),
  we have a natural equivalence

  \vspace{-.5cm}

  \begin{equation}
    \label{ImModalIsBosonicDual}
    \Re \, \circ \rightrightarrows
    \;\;\simeq\;\;
    \Re
  \end{equation}
  saying that the reduction \eqref{ElasticityModalities}
  of the even aspect \eqref{ModeSupergeometric}
  of the space is equivalently the reduced aspect.

  To see this, consider $\mathbb{R}^{n \vert q} \times \mathbb{D}_W
  \;\in\; \infty\mathrm{JetsOfSuperCartesianSpaces}$
  and use, by Example \ref{JetsOfSupergeometricGroupoids}
  with Lemma \ref{ReflectionsOfSuperCommutativeAlgebras},
  the operation $\Re \, \circ \rightrightarrows$ is
  given in terms of the defining super-algebras of functions
  \eqref{SuperFormalCartesianSpaces} by passage to the
  reduced algebra of the even subalgebra:
   \vspace{-2mm}
  $$
    \begin{aligned}
      C^\infty
      \Big(
        \Re \, \circ \rightrightarrows
        \big(
          \mathbb{R}^{n \vert q}
          \times
          \mathbb{D}_W
        \big)
      \Big)
      &
      \simeq
      \;
      \Big(
      \big(
        C^\infty
        \big( \mathbb{R}^n\big)
          \otimes_{\mathbb{R}}
        \big(
          \wedge^\bullet_{\mathbb{R}} \mathbb{R}^q
        \big)
          \otimes_{\mathbb{R}}
        (\mathbb{R} \oplus W)
      \big)_{\mathrm{even}}
      \Big)_{\mathrm{red}}
      \\
      & \simeq \;
      \Big(
        C^\infty
        \big( \mathbb{R}^n\big)
          \otimes_{\mathbb{R}}
        \underset{
          \mathclap{
            \simeq
            \,
            \mathbb{R}
              \oplus
            \wedge^2 \mathbb{R}^q
              \oplus
            \wedge^4 \mathbb{R}^q
              \oplus
            \cdots
          }
        }{
        \underbrace{
          \big(
            \wedge^\bullet_{\mathbb{R}} \mathbb{R}^q
          \big)_{\mathrm{even}}
        }
        }
          \otimes_{\mathbb{R}}
        (\mathbb{R} \oplus W)
      \Big)_{\mathrm{red}}
      \\
      & \simeq \;
      \Big(
        C^\infty(\mathbb{R}^n)
          \otimes_{\mathbb{R}}
        \big(
          \mathbb{R}
          \oplus
          (W \oplus \wedge^2\mathbb{R}^q \oplus \wedge^4 \mathbb{R}^q \oplus \cdots)
        \big)
      \Big)_{\mathrm{red}}
      \\
      & \simeq \;
        C^\infty(\mathbb{R}^n)
          \otimes_{\mathbb{R}}
        \underset{
          \simeq
          \,
          \mathbb{R}
        }{
        \underbrace{
          \big(
            \mathbb{R}
            \oplus
            (W \oplus \wedge^2\mathbb{R}^q \oplus \wedge^4 \mathbb{R}^q \oplus \cdots)
          \big)_{\mathrm{red}}
        }
        }
      \\
      & \simeq
      C^\infty(\mathbb{R}^n)
      \,.
    \end{aligned}
  $$

   \vspace{-2mm}
\noindent  Here in the last step we used that every non-unit element
  in the Grassmann algebra is nilpotent.
 But, by \eqref{SuperShapeInTermsOfFunctionAlgebras} and
  \eqref{InfinitesimalShapeInTermsOfFunctionAlgebras}, we also have
   \vspace{-2mm}
  $$
    \begin{aligned}
      C^\infty
      \big(
        \Re( \mathbb{R}^{n\vert q} \times \mathbb{D}_W )
      \big)
      & \simeq \;
      C^\infty
      \big(
        \mathrm{Shp}_{\mathrm{inf}}
        \circ
        \mathrm{Shp}_{\mathrm{sup}}
        ( \mathbb{R}^{n\vert q} \times \mathbb{D}_W )
      \big)
      \\
      & \simeq
      C^\infty
      \big(
        \mathrm{Shp}_{\mathrm{inf}}
        ( \mathbb{R}^{n} \times \mathbb{D}_W )
      \big)
      \\
      & \simeq
      C^\infty
      \big(
        \mathbb{R}^{n}
      \big)
      \,,
    \end{aligned}
  $$
  where in the first step we used the elastic structure
  \eqref{SolidityAdjunctions}
  $\Re \;:=\;
    \mathrm{Bsnc}
    \circ
    \mathrm{Rdcd}
    \circ
    \mathrm{Shp}_{\mathrm{inf}}
    \circ
    \mathrm{Shp}_{\mathrm{sup}}
  $
  leaving the two full embeddings on the left notationally implicit.
  Since all these equivalences are natural, this implies
  \eqref{ImModalIsBosonicDual}.
  With this, we have the following sequence of natural equivalences
  for general $X \in \mathbf{H} := \infty\mathrm{JetsOfSupergeometricGroupoids}_\infty$:
  $$
    \begin{aligned}
      \mathbf{H}
      \big(
        \mathbb{R}^{n\vert q} \times \mathbb{D}
        \,,\,
        \rightsquigarrow \circ \, \Im( X )
      \big)
      &
      \simeq \;
      \mathbf{H}
      \Big(
        \Re
        \circ
        \rightrightarrows
        \big(
          \mathbb{R}^{n\vert q} \times \mathbb{D}
        \big)
        \,,\,
        X
      \Big)
      \\
      & \simeq \;
      \mathbf{H}
      \Big(
        \Re
        \big(
          \mathbb{R}^{n\vert q} \times \mathbb{D}
        \big)
        \,,\,
        X
      \Big)
      \\
      & \simeq
      \mathbf{H}
      \big(
        \mathbb{R}^{n\vert q} \times \mathbb{D}
        \,,\,
        \Im X
      \big)
      \,,
    \end{aligned}
  $$
  where the first and the last steps are the defining
  hom-equivalences \eqref{AdjunctionHomEquivalence}
  while the middle step is \eqref{ImModalIsBosonicDual}.
  Thus the statement \eqref{ImModalIsBosonic} follows, by the
  $\infty$-Yoneda lemma (Prop. \ref{YonedaLemma}).
\hfill \end{proof}

\subsection{Singularities}
\label{SingularGeometry}

Given a cohesive $\infty$-topos $\mathbf{H}_{\tiny\smooth}$ as in
\cref{DifferentialTopology}, we construct here a
new $\infty$-topos $\mathbf{H}$ (Def. \ref{SingularCohesiveInfinityTopos} below),
to be called \emph{singular-cohesive}, with the following properties:

\begin{enumerate}
\vspace{-3mm}
  \item
    $\mathbf{H}$ contains (\eqref{SingularitiesUnderYoneda} below)
    for each finite group $G$,
    an object $\orbisingularG \in \mathbf{H}$,
    to be thought of as the generic
    $G$-orbi-singularity (\hyperlink{FigureD}{\it Figure D}).

\vspace{-3mm}
  \item $\mathbf{H}$ carries (Prop. \ref{OverSingularitiesCohesion} below)
   an adjoint triple
   of modalities \eqref{AdjointModalities} to be read as follows
   $$
     \underset{
       \mathclap{
       \mbox{
         \small
         \color{darkblue}
         \bf
         ``singular''
       }
       }
     }{
       \singular
     }
     \;\;\;\;\;\;\;\;\;\;\;\;
     \dashv
     \;\;\;\;\;\;\;\;\;\;\;\;
     \underset{
       \mathclap{
       \mbox{
         \small
         \color{darkblue}
         \bf
         ``smooth''
       }
       }
     }{
       \smooth
     }
     \;\;\;\;\;\;\;\;\;\;\;\;
     \dashv
     \;\;\;\;\;\;\;\;\;\;\;\;\;\;\;\;\;
     \underset{
       \mathclap{
       \mbox{
         \small
         \color{darkblue}
         \bf
         ``orbi-singular''
       }
       }
     }{
       \orbisingular
     }
   $$
   with $\mathbf{H}_{\tiny\smooth}$
   being the full sub-$\infty$-category of smooth objects in $\mathbf{H}$,

\vspace{-3mm}
   \item such that (Prop. \ref{OrbiSingularitiesAreOrbiSingular} below):

\hspace{-5mm}
   \begin{tabular}{lll}
     $
       \phantom{\mathclap{\vert_{\vert_{\vert_\vert}}}}
       \singular\big(\orbisingularG\big) \;\simeq\; \ast
     $
     &
     ``The purely singular aspect of an orbi-singularity
     is the quotient of a point, hence a point.''
     \\
     $
       \phantom{\mathclap{\vert_{\vert_{\vert_\vert}}}}
       \smooth\big(\orbisingularG\big) \;\simeq\; \ast \!\sslash\! G
     $
     &
     ``The purely smooth aspect of an orbi-singularity
     is the homotopy quotient of a point.''
     \\
     $\orbisingular\big(\orbisingularG\big) \;\simeq\; \orbisingularG$
     &
     ``An orbi-singularity is purely orbi-singular.''
   \end{tabular}

\end{enumerate}

Essentially this list of conditions might completely characterize
$\mathbf{H}$ to be as in Def. \ref{SingularCohesiveInfinityTopos} below.
Here we leave a fully axiomatic characterization of
singular cohesion as an open problem and are content with
making the following definitions:

\medskip

\noindent {\bf Singular cohesive geometry.}

\begin{defn}[The 2-site of singularities]
  \label{CategoryOfSingularities}
  $\phantom{A}$
  \vspace{-1.5mm}
 \item {\bf (i)}  We write
  \begin{equation}
    \label{Singularities}
    \mathrm{Singularities}
    \;:=\;
    \mathrm{Groupoids}_{\leq 1, \mathrm{cn},\mathrm{fin}}
     \xymatrix{\ar@{^{(}->}[r] &}
    \mathrm{Groupoids}_\infty
  \end{equation}
  for the full sub-$\infty$-category of $\infty$-groupoids
  on the connected 1-truncated objects whose $\pi_1$ is finite.
  \item {\bf (ii)}   A skeleton of this $(2,1)$-category has,
  of course, as objects
  the delooping groupoids (Example \ref{DeloopingGroupoids})
  $\ast \!\sslash\! G$
  that are presented by a single object
  and a \emph{finite} group $G$ of automorphisms of that object.
 \item {\bf (iii)}    When regarded as objects of $\mathrm{Singularities}$
 in \eqref{Singularities},
  we will denote these by ``$\orbisingular$'' attached
  to the symbol for the group:
  \begin{equation}
    \label{AnOrbifoldSingularity}
    \xymatrix@R=1em{
      \orbisingularG
      \ar@{|->}[d]
      \ar@{}[r]|-{ \in }
      &
      \mathrm{Singularities}
      \ar@{^{(}->}[d]
      \\
      \ast \!\sslash\! G
      \ar@{}[r]|-{ \in }
      &
      \mathrm{Groupoids}_\infty
    }
  \end{equation}
  \vspace{-5mm}
  \item {\bf (iv)}   The hom-$\infty$-groupoids between these singularities are,
  equivalently, the
  action groupoids (Example \ref{ActionGroupoids})
  whose objects are group homomorphisms and whose
  morphisms are conjugation actions on these:
  \begin{equation}
    \label{SingularitiesHomGroupoids}
    \begin{aligned}
    \mathrm{Singularities}
    \big(
      \orbisingularGa,
      \orbisingularGb
    \big)
    & :=
    \mathrm{Groupoids}_\infty
    \big(
      \ast \!\sslash\! G_1,
      \ast \!\sslash\! G_2
    \big)
    \\
    &
    \phantom{:}
    \simeq
    \mathrm{Groups}(G_1,G_2) \!\sslash\!_{\mathrm{conj}} G_2
    \end{aligned}
  \end{equation}

  \vspace{-3mm}
   \item {\bf (v)}  We regard $\mathrm{Singularities}$ as an $\infty$-site with
  trivial Grothendieck topology, so that $\infty$-sheaves on
  $\mathrm{Singularities}$ are $\infty$-presheaves \eqref{HValuedPresheaves}.
\end{defn}

\begin{remark}[The global orbit category]
 \label{GlobalOrbitCategory}
 The category $\mathrm{Singularities}$ in Def. \ref{CategoryOfSingularities}
 is sometimes known in the literature
 as the ``global orbit category''
 (though at other times this term is used for its wide but non-full subcategory
 on the faithful morphisms).
 It has elsewhere been denoted:
 ``$\mathrm{Orb}$ case \textcircled{\raisebox{-1pt}{$1$}}''
 (in \cite[4.1]{HenriquesGepner07}),
 ``$\mathrm{Glob}$'' (in \cite[2.2]{Rezk14}),
 ``$\mathrm{Orb}$'' (in \cite[2.1]{Koerschgen16}\cite[3.2]{Juran20})
 and (up to equivalence)
 ``$\mathbf{O}_{\mathrm{gl}}$'' (in \cite{Schwede17}\cite[2.2]{Koerschgen16}).
 The terminology in Def. \ref{CategoryOfSingularities}
 is meant to be more suggestive of the role this category
 plays in the theory, from the perspective of cohesive homotopy theory.
 In fact, the (global) orbit category is often taken to contain not
 just all finite groups, but all compact Lie groups, with the
 hom-spaces then being the geometric realization of the
 topological mapping groupoids.
 We restrict to discrete groups (hence finite if compact)
 for reasons explained in
 Remark \ref{NeedForDiscreteGroupsInSingularities} below.
 This restriction is also amplified in \cite{DHLPS19}.
\end{remark}

\begin{defn}[Singular-cohesive $\infty$-topos]
  \label{SingularCohesiveInfinityTopos}
  Consider a cohesive $\infty$-topos (Def. \ref{CohesiveTopos}),
   now to be denoted with ``$\smooth$''-subscripts
      \vspace{-1mm}
  \begin{equation}
    \label{NonSingularCohesiveTopos}
    \xymatrix{
      \mathbf{H}_{\tiny\smooth}
      \;\;  \ar@{->}@<+28pt>[rr]|<\times|-{\; \mathrm{Shp}\; }_-{\raisebox{-6pt}{\tiny $\bot$}}
        \ar@{<-^{)}}@<+14pt>[rr]|-{\; \mathrm{Dsc} \;}_-{\raisebox{-6pt}{\tiny $\bot$}}
        \ar@<+0pt>@<-0pt>[rr]|-{\; \mathrm{Pnts}\; }_-{\raisebox{-6pt}{\tiny $\bot$}}
        \ar@{<-^{)}}@<-14pt>[rr]|-{\; \mathrm{coDsc} \; }
      &&
     \;\;
     \mathbf{B}_{\tiny\smooth}
     \; := \; \mathrm{Groupoids}_\infty
    }
  \end{equation}

  \vspace{1mm}
\noindent   and assumed to have a site of $\mathrm{Charts}$
  (Def. \ref{ChartsForCohesion}).
  The corresponding \emph{singular-cohesive} $\infty$-topos is
  that of $\mathbf{H}_{\tiny\smooth}$-valued
  $\infty$-sheaves \eqref{HValuedPresheaves}
  over the site of $\mathrm{Singularities}$
  (Def. \ref{CategoryOfSingularities}):
  \vspace{-1mm}
  \begin{equation}
    \label{InfinitToposOverSingularities}
    \raisebox{20pt}{
    \xymatrix@C=5em@R=1.3em{
      \mathllap{
        \mathbf{H} :=
      }
      \mathrm{Sheaves}_\infty
      \big(
        \mathrm{Singularities},
        \,
        \mathbf{H}_{\tiny\smooth}
      \big)
      \ar@<+10pt>@{<-}[rr]^-{ \mathrm{NnOrbSnglr} }
      \ar@<-4pt>[rr]_-{ \mathrm{Smth}   }^-{ \bot }
      \ar@{<-}@<+34pt>[dd]_-{ \mathrm{Disc} }
      \ar@<+46pt>[dd]^-{ \mathrm{Pnts} }_-{ \dashv }
      &&
      \;\;\;
      \mathbf{H}_{\tiny\smooth}
      \ar@{<-}@<-2pt>[dd]_-{ \mathrm{Disc} }
      \ar@<+10pt>[dd]^-{ \mathrm{Pnts} }_-{ \dashv }
      \\
      \\
      \mathllap{
        \mathbf{B} := \;
      }
      \mathrm{Sheaves}_\infty
      \big(
        \mathrm{Singularities},
        \,
        \mathbf{B}_{\tiny\smooth}
      \big)
      \ar@<+10pt>@{<-^{)}}[rr]^-{ \mathrm{NnOrbSnglr} }
      \ar@<-4pt>[rr]_-{ \mathrm{Smth}  }^-{ \bot }
      &&
      \;\;\;\mathbf{B}_{\tiny\smooth}
    }
    }
  \end{equation}

  \vspace{-1mm}
\noindent
  where horizontally we are showing the base
  geometric morphisms (Prop. \ref{BaseGeometricMorphism})
  of sheaves over the site $\mathrm{Singularities}$,
  while vertically we are showing
  the base geometric morphism \eqref{AdjunctionCohesion}
  of
  $\mathbf{H}_{\tiny\smooth}$ over $\mathbf{B}_{\tiny\smooth}$
  extended objectwise over $\mathrm{Singularities}$, by functoriality.
\end{defn}

 \begin{lemma}[Singularities is 2-site for homotopical cohesion]
   \label{SingularitiesIsSiteForHomotopicalCohesion}
   The 2-site $\mathrm{Singularities}$ (Def. \ref{CategoryOfSingularities})
   is an $\infty$-site for homotopical cohesion, in the sense of
   Def. \ref{SiteForHomotopicalCohesion}.
 \end{lemma}
 \begin{proof}
  It is immediately checked that
  \begin{enumerate}
  \vspace{-2mm}
    \item the terminal object
    is given by the trivial group:

    \vspace{-1.1cm}

    \begin{equation}
      \label{SingTerm}
      \ast
        \;\simeq\;
      \orbisingularE
    \end{equation}

    \vspace{-.4cm}

    \item
    Cartesian product
    is direct product of groups:
    \hspace{13pt}
    $
      \orbisingularGa
      \,\times
      \orbisingularGb
      \;\simeq\;\;
      \orbisingularGab
      \,.
    $
\end{enumerate}

\vspace{-.8cm}

\hfill \end{proof}

  \begin{prop}[Singular cohesion]
    \label{OverSingularitiesCohesion}
    A singular-cohesive $\infty$-topos
    (Def. \ref{SingularCohesiveInfinityTopos})
    \vspace{-2mm}
    $$
      \xymatrix@C=5em@R=.1em{
        &
        \mathbf{H}
        \ar[dr]^-{ \mathrm{Pnts} }
        \ar[dl]_{ \mathrm{Smth} }
        \\
        \mathbf{H}_{\tiny\smooth}
        \ar[dr]_-{ \mathrm{Pnts} }
        &&
        \mathbf{B}
        \ar[dl]^-{ \mathrm{Smth} }
        \\
        &
        \mathbf{B}_{\tiny\smooth}
      }
    $$

     \vspace{-2mm}
\noindent
  is itself cohesive (Def. \ref{CohesiveTopos}) in two ways:
  \vspace{-1mm}
  \item {\bf (i)}
      over the singular-base $\infty$-topos $\mathbf{B}$
      by the cohesion of
      $\mathbf{H}_{\tiny\smooth} \to \mathbf{B}_{\tiny\smooth}$
      \eqref{AdjunctionCohesion}
      applied object-wise over all $\mathrm{Singularities}$
  \begin{equation}
    \label{CohesionOfSingularCohesive}
    \xymatrix@C=3em{
      \mathbf{H}
   \;\;     \ar@{->}@<+28pt>[rr]|<\times|-{
            \mathllap{
              \scalebox{.8}{
                \color{darkblue}
                \bf
                ``shape''
                \hspace{3.4cm}
              }
           }
           \;\mathrm{Shp} \;
        }_-{\raisebox{-6pt}{\tiny $\bot$}}
        \ar@{<-^{)}}@<+14pt>[rr]|-{
            \mathllap{
              \scalebox{.8}{
                \color{darkblue}
                \bf
                ``discrete''
                \hspace{3.4cm}
              }
           }
           \; \mathrm{Disc} \;
        }_-{\raisebox{-6pt}{\tiny $\bot$}}
        \ar@<+0pt>@<-0pt>[rr]|-{
            \mathllap{
              \scalebox{.8}{
                \color{darkblue}
                \bf
                ``points''
                \hspace{3.4cm}
              }
           }
          \;\mathrm{Pnts}\;
        }_-{\raisebox{-6pt}{\tiny $\bot$}}
        \ar@{<-^{)}}@<-14pt>[rr]|-{
            \mathllap{
              \scalebox{.8}{
                \color{darkblue}
                \bf
                ``chaotic''
                \hspace{3.4cm}
              }
           }
          \;\mathrm{Chtc} \;
        }
      &&
  \;\;    \mathbf{B}\;;
    }
  \end{equation}
    \vspace{-2mm}
      \item {\bf (ii)}
    over the non-singular cohesive base $\infty$-topos
    $\mathbf{H}_{\tiny\smooth}$
    (Def. \ref{CohesiveTopos})
    in that
    the global section geometric morphism
    $\mathbf{H}
      \overset{ \mathrm{Smth} }{\longrightarrow}
      \mathbf{H}_{\tiny\smooth}
    $
    of \eqref{InfinitToposOverSingularities}
    is part of a cohesive adjoint quadruple, to be denoted
    \begin{equation}
      \label{CohesionSingular}
      \hspace{1cm}
      \xymatrix@C=4em{
        \mathbf{H}
       \;\; \ar@{->}@<+31pt>[rr]|<\times|-{
          \mathllap{
            \scalebox{.8}{
              \color{darkblue}
              \bf
              ``singular''
              \hspace{3.4cm}
            }
          }
         \; \mathrm{Snglr} \;
        }_-{\raisebox{-7pt}{\tiny $\bot$}}
        \ar@{<-^{)}}@<+14pt>[rr]|-{
          \mathllap{
            \scalebox{.8}{
              \color{darkblue} \bf
              ``not orbi-singular''
              \hspace{2.85cm}
            }
          }
         \; \mathrm{NnOrbSnglr} \;
        }_-{\raisebox{-6pt}{\tiny $\bot$}}
        \ar@<+0pt>@<-0pt>[rr]|-{
          \mathllap{
            \scalebox{.8}{
              \color{darkblue} \bf
              ``smooth''
              \hspace{3.45cm}
            }
          }
        \;  \mathrm{Smth} \;
        }_-{\raisebox{-6pt}{\tiny $\bot$}}
        \ar@{<-^{)}}@<-14pt>[rr]|-{
          \mathllap{
            \scalebox{.8}{
              \color{darkblue} \bf
              ``orbi-singular''
              \hspace{3.2cm}
            }
          }
      \;    \mathrm{OrbSnglr} \;
        }
        &&
        \;\; \mathbf{H}_{\tiny\smooth}
      }
      \,.
    \end{equation}
\end{prop}
  \begin{proof}
    The first statement is immediate.
    The second statement follows via Lemma \ref{SingularitiesIsSiteForHomotopicalCohesion}
    by Example \ref{InfinityCohesivePresheafSite}.
  \hfill
\end{proof}

\begin{notation}[Singular-elastic/solid $\infty$-topos]
  \label{SingularSolidInfinityTopos}
  Let $\mathbf{H}$ be a singular-cohesive $\infty$-topos
  (Def. \ref{SingularCohesiveInfinityTopos})
  with underlying smooth cohesive $\infty$-topos
  $\mathbf{H}_{\tiny\smooth} \hookrightarrow \mathbf{H}$.
Then

\vspace{-1mm}
\item {\bf (i)} if $\mathbf{H}_{\tiny\smooth}$
   is in fact an elastic $\infty$-topos (Def. \ref{ElasticInfinityTopos}),
   we say that $\mathbf{H}$ is a \emph{singular-elastic $\infty$-topos};

\vspace{-1mm}
  \item {\bf (ii)}  if $\mathbf{H}_{\tiny\smooth}$
   is in fact a solid $\infty$-topos (Def. \ref{SolidTopos}),
   we say that $\mathbf{H}$ is a \emph{singular-solid $\infty$-topos}.
\end{notation}

\begin{defn}[Singular-cohesive modalities]
  \label{SincularCohesionModalies}
Given a singular cohesive $\infty$-topos (Def. \ref{SingularCohesiveInfinityTopos}),
with its singular cohesion from Prop. \ref{OverSingularitiesCohesion},
we write
  \begin{equation}
    \label{SingularityModalities}
    \underset{
      \mbox{\footnotesize \color{darkblue} \bf ``singular''}
    }{
    {
                \big(
      \singular
      \;:=\;
      \mathrm{NnOrbSnglr}
        \circ
      \mathrm{Snglr}
    \big)
    }}
    \;\;\dashv\;\;
    \underset{
      \mbox{\footnotesize\color{darkblue} \bf
        \begin{tabular}{c}
          ``smooth''
        \end{tabular}
      }
    }{
    {
    \big(
      \smooth
      \;:=\;
      \mathrm{NnOrbSnglr}
        \circ
      \mathrm{Smth}
    \big)
    }
    }
    \;\;\dashv\;\;
    \underset{
      \mbox{ \footnotesize\color{darkblue} \bf ``orbi-singular'' }
    }{
    {
      \big(
        \orbisingular
        \;:=\;
        \mathrm{OrbSnglr}
          \circ
        \mathrm{Smth}
      \big)
    }
    }
  \end{equation}

  \vspace{-1mm}
 \noindent  for the adjoint triple of modalities
 $\mathbf{H} \to \mathbf{H}$
  induced \eqref{AdjointModalities}
  via \eqref{CohesionSingular};
  accompanying the
  cohesive modalities \eqref{CohesiveModalitiesFromAdjointQuadruple}
  induced via \eqref{CohesionOfSingularCohesive}.
\end{defn}

\noindent
The above terminology reflects the difference
  (see \hyperlink{FigureD}{Figure D}) between
  a plain singularity $\singular$ (singular but not orbi-singular)
  as opposed to its enhancement to an actual orbifold singularity
  $\orbisingular$.
We record the following elementary but important consequence:
\begin{prop}[Smooth orbi-singular is smooth]
  \label{SmoothOrbiSingularIsSmooth}
  The singularity modalities (Def. \ref{SincularCohesionModalies})
  satisfy:
  $$
    \singular \circ \smooth
    \;\simeq\;
    \smooth
    \phantom{AA}
    \mbox{and}
    \phantom{AA}
    \smooth \circ \orbisingular
    \;\simeq\;
    \smooth\;.
  $$
\end{prop}
\begin{proof}
  As in Prop. \ref{CompositeCohesiveModalities}.
\hfill \end{proof}

\begin{lemma}[Objectwise application of singularity modalities]
  \label{ObjectwiseApplicationOfSingularityModalities}
  The singular-modalities in \eqref{CohesionSingular}
  are computed objectwise
  over $\mathrm{Charts}$, as in Example \ref{InfinityCohesivePresheafSite},
  followed by $\infty$-sheafification $L_{\mathrm{Charts}}$ \eqref{ToposReflectionInPresheaves}:
 \vspace{-2mm}
  $$
  \hspace{-1mm}
    \scalebox{.82}{\xymatrix@C=1.9em{
      \mathrm{Sheaves}_\infty
      \big(
        \mathrm{Singularities}
         \!\times\!
        \mathrm{Charts}
      \big)
      \ar@/^3pc/[rrrrr]^-{ \mathrm{Snglr} }
      \ar@{<-}@/_2.8pc/[rrrrr]_-{
        \mathrm{NnOrbSnglr}
      }
      \ar@{^{(}->}[r]
      &
      \mathrm{PreSheaves}_\infty
      \big(
        \mathrm{Singularities}
        \!\times\!
        \mathrm{Charts}
      \big)
      \ar@<+8pt>[rr]^-{
        \underset{
          \underset{
            \mathrm{Singularities}
          }{\longrightarrow}
        }{\mathrm{lim}}
      }
      \ar@{<-^{)}}@<-8pt>[rr]^-{ \raisebox{2pt}{\scalebox{0.8}{$\bot$}} }_-{
        \mathrm{const}_{\mathrm{Singularities}}
      }
      &&
      \mathrm{PreSheaves}_\infty
      \big(
        \mathrm{Charts}
      \big)
      \ar@<+8pt>[rr]|-{ L_{\mathrm{Charts}} }
      \ar@{<-^{)}}@<-8pt>[rr]^-{ \raisebox{2pt}{\scalebox{0.8}{$\bot$}} }
      &&
      \mathrm{Sheaves}_\infty
      \big(
        \mathrm{Charts}
      \big)
    }
    }
  $$
\end{lemma}
\begin{proof}
  By essential uniqueness of adjoints \eqref{AdjunctionHomEquivalence}.
\end{proof}

 \medskip

 \noindent {\bf Examples of singular-cohesive $\infty$-toposes.}
  \label{ExamplesOfSingularCohesiveToposes}
  \begin{example}[Singular $\infty$-groupoids]
    \label{SingularInfinityGroupoids}
    For $\mathbf{H}_{\tiny\smooth} := \mathrm{Groupoids}_{\infty}$ the base
    $\infty$-topos of plain $\infty$-groupoids \eqref{InfinityCategoryOfInfinityGroupoids},
    the singular-cohesive $\infty$-topos from
    Def. \ref{SingularCohesiveInfinityTopos}
    \vspace{-2mm}
    $$
      \mathrm{SingularGroupoids}_\infty
      \;:=\;
      \mathrm{Sheaves}_\infty
      \big(
        \mathrm{Singularities},
        \,
        \mathrm{Groupoids}_\infty
      \big)
    $$

    \vspace{-2mm}
\noindent
    is that of traditional unstable global homotopy theory
    \cite[\S 1s]{Schwede18},
    as discussed in this form in \cite[\S 4.1]{Rezk14}
    (here with evaluation on all finite groups instead of
    all compact Lie groups).
  \end{example}
\begin{example}[Singular-smooth $\infty$-groupoids]
  \label{OrbiSingularSmoothInfinityGroupoids}

  {\bf (i)} We call the singular-cohesive $\infty$-topos (Def. \ref{SingularCohesiveInfinityTopos})
  over those of smooth $\infty$-groupoids (Example
  \ref{SmoothInfinityGroupoids})
  the $\infty$-topos of
  \emph{singular-smooth $\infty$-groupoids}:
  \begin{equation}
    \label{SingularSmoothGroupoids}
    \begin{aligned}
    \mathrm{SingularSmoothGroupoids}_\infty
    & \;:=\;
    \mathrm{Sheaves}_\infty
    \big(
      \mathrm{Singularities},
      \,
      \mathrm{SmoothGroupoids}_\infty
    \big)
    \\
    &
    \;\;\simeq\;
    \mathrm{Sheaves}_\infty
    \big(
      \mathrm{CartesianSpaces}
      \times
      \mathrm{Singularities}
    \big)
    \,.
    \end{aligned}
  \end{equation}

\noindent  {\bf (ii)} We call the singular-elastic $\infty$-topos (Def. \ref{SingularSolidInfinityTopos})
  over $\mathrm{JetsOfSmoothGroupoids}_\infty$ (Example \ref{FormalSmoothInfinityGroupoids})
  \begin{equation}
    \label{SingularSmoothGroupoids}
    \begin{aligned}
    \mathrm{SingularJetsOfSmoothGroupoids}_\infty
    & \;:=\;
    \mathrm{Sheaves}_\infty
    \big(
      \mathrm{Singularities},
      \,
      \mathrm{JetsOfSmoothGroupoids}_\infty
    \big)
    \\
    &
    \;\;\simeq\;
    \mathrm{Sheaves}_\infty
    \big(
      \mathrm{JetsOfCartesianSpaces}
      \times
      \mathrm{Singularities}
    \big)
    \,.
    \end{aligned}
  \end{equation}

\noindent  {\bf (iii)} We call the singular-solid $\infty$-topos (Def. \ref{SingularSolidInfinityTopos})
  over $\infty\mathrm{JetsOfSupergeometricGroupoids}_\infty$
  (Example \ref{JetsOfSupergeometricGroupoids})
  \begin{equation}
    \label{SingularSupergeometricGroupoids}
    \hspace{-3mm}
    \begin{aligned}
    \mathrm{Singular}\infty\mathrm{JetsOfSupergeometricGroupoids}_\infty
    & :=
    \mathrm{Sheaves}_\infty
    \big(
      \mathrm{Singularities},
      \,
      \infty\mathrm{JetsOfSupergeometricGroupoids}_\infty
    \big)
    \\
    &
    \simeq
    \mathrm{Sheaves}_\infty
    \big(
      \infty\mathrm{JetsOfSuperCartesianSpaces}
      \times
      \mathrm{Singularities}
    \big)
    \,.
    \end{aligned}
  \end{equation}

  \noindent For the second lines of \eqref{SingularSmoothGroupoids},
  \eqref{SingularSmoothGroupoids}, and \eqref{SingularSupergeometricGroupoids},
  see Lemma \ref{YonedaOnProductSite}.
\end{example}

\medskip

\noindent {\bf Basic properties of singular cohesion.}
\begin{defn}[Orbi-singularities]
  \label{OrbiSingularities}
  Let $\mathbf{H}$ be singular-cohesive $\infty$-topos
  (Def. \ref{SingularCohesiveInfinityTopos}).

  \noindent {\bf (i)} We regard the objects
  $\orbisingularG \;\in\; \mathrm{Singularities}$
  \eqref{AnOrbifoldSingularity} as objects of $\mathbf{H}$ under the
  $\infty$-Yoneda-embedding (Prop. \ref{InfinityYonedaEmbedding})
  and the inclusion \eqref{InfinitToposOverSingularities}
  of discrete objects:
   \vspace{-2mm}
  \begin{equation}
    \label{SingularitiesUnderYoneda}
    \xymatrix{
    \orbisingularG
    \;\in\;
    \mathrm{Singularities}
    \; \ar@{^{(}->}[r]^-{y}  &
    \mathrm{Sheaves}_\infty
    \big(
      \mathrm{Singularities},
      \mathbf{B}_{\tiny\smooth}
    \big)
    \; \ar@{^{(}->}[r]^-
    {   \scalebox{.5}{$
         \mathrm{Disc}$}
        }
        &
    \mathrm{Sheaves}_\infty\big(
      \mathrm{Singularities},
      \mathbf{H}_{\tiny\smooth}
    \big)
    \;=\;
    \mathbf{H}
    }.
  \end{equation}

   \vspace{-2mm}
  \noindent {\bf (ii)} More generally, for
  \vspace{-2mm}
  $$
    G
    \;\in\;
    \xymatrix{
      \mathrm{Groups}(\mathrm{Groupoids}_\infty)
      \ar[rr]^-
      {   \scalebox{.5}{$
          \mathrm{Groups} ( \mathrm{Disc})$}
        }
      &&
      \mathrm{Groups}(\mathbf{H}_{\scalebox{.6}{$\smooth$}})
    }
  $$
  any discrete $\infty$-group \eqref{CohesionOfSingularCohesive},
  we also write
  \begin{equation}
    \label{GeneralOrbiSingularities}
    \orbisingularG \;:=\;
    \orbisingular( \mathbf{B}G  )
    \;\in\;
    \mathbf{H}
  \end{equation}
  for the orbi-singularization \eqref{CohesionSingular}
  of its delooping \eqref{LoopingDeloopingEquivalence}.
  \end{defn}
  Lemma \ref{ImagesAndPreimagesOfOrbisingularities} shows that the
  two notations in Def. \ref{OrbiSingularities} are consistent with each other.

  \begin{remark}[Finite groups in singular cohesion]
    \label{DiscreteObjectsInSingularCohesion}
    Given a singular-cohesive $\infty$-topos (Def. \ref{SingularCohesiveInfinityTopos}),
    the images of a finite group $G$ under the following sequence of inclusions are naturally
    all denoted by the same symbol:
    \vspace{-2mm}
    \begin{equation}
      \label{InclusionOfFiniteGroups}
      \hspace{-4mm}
      \xymatrix@R=-1pt@C=20pt{
        \mathrm{Groups}^{\mathrm{fin}}
       \; \ar@{^{(}->}[r]
        &
        \mathrm{Groups}
        (
          \mathrm{Set}
        )
      \;  \ar@{^{(}->}[r]
        &
        \mathrm{Groups}
        (
          \mathrm{Groupoids}_{\infty}
        )
        \; \ar@{^{(}->}[rr]^-{   \scalebox{.5}{$
          \mathrm{Groups} ( \mathrm{Disc})$}
        }
        &&
        \mathrm{Groups}
        \big(
          \mathbf{H}_{\tiny\smooth}
        \big)
     \;   \ar@{^{(}->}[rr]^-{   \scalebox{.5}{$
          \mathrm{Grp}( \mathrm{NnOrbSnglr})$}
        }
        &&
        \mathrm{Groups}
        \big(
          \mathbf{H}
        \big)
        \\
        G
        \ar@{|->}[r]
        &
        G
        \ar@{|->}[r]
        &
        G
        \ar@{|->}[rr]
        &&
        G
        \ar@{|->}[rr]
        &&
        G
      }
    \end{equation}

    \vspace{-1mm}
\noindent
    With this understood, we also have identifications as
    follows (where now the ambient $\infty$-categories are implicit
    from the context):
    \begin{equation}
      \label{DiscreteEquivalences}
      \ast \!\sslash\! G
      \;\simeq\;
      \mathrm{Disc}(\ast \!\sslash\! G)
      \phantom{AA}
      \mbox{and}
      \phantom{AA}
      \orbisingularG
      \;\simeq\;
      \mathrm{Disc}\big(
        \orbisingularG
      \big)
    \end{equation}
    where on the right we are recalling the definition \eqref{SingularitiesUnderYoneda}.
\end{remark}
Similarly:

\begin{remark}[Smooth charts in singular cohesion]
  \label{SmoothModelSpacesInSingularCohesion}
  Consider a singular-cohesive $\infty$-topos (Def. \ref{SingularCohesiveInfinityTopos})
  with an $\infty$-site $\mathrm{Charts}$ of charts (Def. \ref{ChartsForCohesion}).
  Then images of the charts $U \in \mathrm{Charts}$ under
    the $\infty$-Yoneda embedding (Prop. \ref{InfinityYonedaEmbedding}),
    and further under
    $\mathrm{NnOrbSnglr}$ \eqref{InfinitToposOverSingularities},
   are naturally denoted by the same symbol:
 \vspace{-2mm}
  \begin{equation}
    \label{SmoothChartsInSingularCohesion}
    \xymatrix@R=-1pt{
      \mathcal{S}
      \ar[rr]^-{ y }
      &&
      \mathbf{H}_{\tiny\smooth}
      \ar[rr]^-
      {   \scalebox{.5}{$
            \mathrm{NnOrbSnglr}$}
        }
      &&
      \mathbf{H}
      \\
      U
      \ar@{|->}[rr]
      &&
      U
      \ar@{|->}[rr]
      &&
      U
    }
  \end{equation}
\end{remark}

\begin{lemma}[$\infty$-Yoneda on product site]
  \label{YonedaOnProductSite}
  Consider a singular-cohesive $\infty$-topos $\mathbf{H}$
  (Def. \ref{SingularCohesiveInfinityTopos})
  with an $\infty$-site $\mathrm{Charts}$ of cohesive charts
  (Def. \ref{ChartsForCohesion}) for $\mathbf{H}_{\tiny\smooth}$.

  \vspace{-1mm}
 \item {\bf (i)} Then a site (Def. \ref{InfinityToposOfInfinitySheaves})
 for the full singular-cohesive $\mathbf{H}$ is the
  Cartesian product site
  \vspace{-2mm}
  \begin{equation}
      \label{SingularCharts}
    \mathrm{SingularCharts}
     \;:=\;
    \mathrm{Charts} \times \mathrm{Singularities}
  \end{equation}

  \vspace{-2mm}
\noindent
  in that
\vspace{-1mm}
 \begin{equation}
    \label{SingularCohesionAsSheavesOnProductSite}
    \mathbf{H}
    \;\simeq\;
    \mathrm{Sheaves}_\infty
    \big(
      \mathrm{Charts}
        \times
      \mathrm{Singularities}
    \big).
  \end{equation}

  \vspace{-2mm}
  \item {\bf (ii)} Under the $\infty$-Yoneda embedding (Prop. \ref{InfinityYonedaEmbedding})
  objects in the product site map to the
  Cartesian product of their prolonged Yoneda embeddings
  (in the sense of Remark \ref{DiscreteObjectsInSingularCohesion} and
  Remark \ref{SmoothModelSpacesInSingularCohesion}):
\vspace{-2mm}
  $$
    \xymatrix@R=-2pt{
      \mathrm{Charts}
      \times
      \mathrm{Singularities}
      \ar[rr]^-{ y }
      &&
      \mathbf{H}
      \\
      \big(
        U,
        \orbisingularG
      \big)
      \ar@{|->}[rr]
      &&
      U  \times \orbisingularG
      \,,
    }
  $$

  \vspace{-3mm}
\noindent  where on the right we are using the abbreviated notation from
  \eqref{SingularitiesUnderYoneda} and \eqref{SmoothChartsInSingularCohesion}.
\end{lemma}

\begin{proof}
  On the one hand, we have
  a natural equivalence
  \vspace{-2mm}
  \begin{equation}
    \label{YonedaEmbeddingOfObjectsInProductSite}
    \begin{aligned}
      \mathbf{H}
      \Big(
        y\big( U_1, \orbisingularGa \, \big) \; , \;
        y\big( U_2, \orbisingularGb \, \big)
      \Big)
      &
      \; \simeq \;
      \mathrm{Charts}( U_1, U_2)
      \times
      \mathrm{Singularities}
      \big(
        \ast \!\sslash\! G_1
        \; ,  \,
        \ast \!\sslash\! G_2
      \big)
    \end{aligned}
  \end{equation}

  \vspace{-2mm}
\noindent
  by fully-faithfulness of the $\infty$-Yoneda embedding
  (Prop. \ref{InfinityYonedaEmbedding})
  and by the definition of product sites.
  On the other hand, we have a sequence of natural equivalences
  \vspace{-1mm}
  \begin{equation}
    \label{YonedaImageOfCartesianProductsOfChartsWithSingularities}
    \begin{aligned}
            \mathbf{H}
      \Big(
        y
        \big(
          U_1,
          \orbisingularGa
        \big)
        \; , \;
        U_2
        \times
        \orbisingularGb
      \Big)
       &=
      \mathbf{H}
      \Big(
        y
        \big(
          U_1,
          \orbisingularGa
        \big),
        \,
        \mathrm{NnOrbSnglr}(U_2)
        \times
        \mathrm{Disc}\big(
          \orbisingularGb
        \big)
      \Big)
      \\
      & \simeq
      \mathbf{H}
      \Big(
        y
        \big(
          U_1,
          \orbisingularGa
        \big),
        \,
        \mathrm{NnOrbSnglr}(U_2)
      \Big)
      \times
      \mathbf{H}
      \Big(
        y
        \big(
          U_1,
          \orbisingularGa
        \big),
        \,
        \mathrm{Disc}\big(
          \orbisingularGb
        \big)
      \Big)
      \\
      & \simeq
      \mathbf{H}_{\tiny\smooth}
      \Big(
        \mathrm{Snglr}
        \big(
        y
        \big(
          U_1,
          \orbisingularGa
        \big)
        \big),
        \,
        U_2
      \Big)
      \times
      \mathbf{B}
      \Big(
        \raisebox{1pt}{\textesh}
        \big(
        y
        \big(
          U_1,
          \orbisingularGa
        \big)
        \big),
        \,
        \orbisingularGb
      \Big)
      \\
      & \simeq
      \mathbf{H}_{\tiny\smooth}
      \big(
        U_1,
        U_2
      \big)
      \times
      \mathbf{B}
      \big(
        \orbisingularGa
        \; , \;
        \orbisingularGb
    \,  \big)
      \\
      & \simeq
      \mathrm{Charts}
      \big(
        U_1,
        U_2
      \big)
      \times
      \mathrm{Singularities}
      \big(
        \orbisingularGa
        \; , \;
        \orbisingularGb
     \, \big).
    \end{aligned}
  \end{equation}

  \vspace{-1mm}
\noindent
  Here the first step is by definition, the second
  step is the universal property of the Cartesian product,
  and the third step is the hom-equivalence \eqref{AdjunctionHomEquivalence}
  of the adjunctions
  $\mathrm{Snglr} \dashv \mathrm{NnOrbSnglr}$
  and
  $\raisebox{1pt}{\textesh} \dashv \mathrm{Disc}$, respectively.
  In the fourth step, we use
  \eqref{GeometricallyContractibleGenerators} and \eqref{SinglrG},
  respectively. The last step is the fully-faithfulness
  of the $\infty$-Yoneda embedding (Prop. \ref{InfinityYonedaEmbedding}).
  Since both \eqref{YonedaEmbeddingOfObjectsInProductSite} and
  \eqref{YonedaImageOfCartesianProductsOfChartsWithSingularities}
  are natural in $\big( U', (\ast \!\sslash\! G)_{\tiny\orbisingular} \big)$,
  and since their right hand sides coincide, it follows by the
  $\infty$-Yoneda embedding (Prop. \ref{InfinityYonedaEmbedding})
  that also the representatives of the
  left hand sides coincide:
  $
    y\big(
      U_2,
      \orbisingularGb
    \big)
    \;\simeq\;
    U_2
    \times
    \orbisingularGb
  $.

  \vspace{-7mm}
\hfill
\end{proof}

  \begin{lemma}[Images and pre-images of orbi-singularities]
    \label{ImagesAndPreimagesOfOrbisingularities}
    Let $\mathbf{H}$ be a singular-cohesive $\infty$-topos
    (Def. \ref{SingularCohesiveInfinityTopos}).
    Then the images and pre-images of the
    generic singularities
    $\orbisingularG$
    \eqref{SingularitiesUnderYoneda} under
    the functors \eqref{CohesionSingular}
    exhibiting the singular cohesion
    are as follows (see \hyperlink{FigureD}{\it Figure D}):
    \vspace{-2mm}
    \begin{equation}
      \label{ImagedAndPreimagesofBasicOrbiSngularities}
      \raisebox{20pt}{
      \xymatrix@C=5em@R=-2pt{
        &
        \orbisingularG
        \ar@{|->}[dl]_{
          \mathrm{Snglr}
        }
        \ar@<-4pt>@{|->}[dr]_-{
          \mathrm{Smth}
        }
        \ar@<+4pt>@{<-|}[dr]^-{
          \mathrm{OrbSnglr}
        }
        & &
        \ast \!\sslash\! G
        &
        \mathrlap{\in \mathbf{H}}
        \\
        \mathllap{\ast = \;}
        \ast \!/ G
        &&
        \ast \!\sslash\! G
        \ar@{|->}[ur]_-{\;\;\;\;\;
          \mathrm{NnOrbSnglr}
        }
        &
        &
        \mathrlap{\in \mathbf{H}_{\tiny\smooth}}
      }
      }
    \end{equation}
  \end{lemma}
  \begin{proof}
  By the singular cohesion established in the
  proof of Prop. \ref{OverSingularitiesCohesion} we have that:
    \begin{enumerate}
    \vspace{-2mm}
    \item  the functor
    $
      \mathrm{Snglr}
      \;\simeq\;
      \underset{\longrightarrow}{\lim}
    $
    is the colimit functor (Prop. \ref{LimitsAndColimitsAsAdjoints}),
    \vspace{-4mm}
    \item the functor
    $
      \mathrm{Smth}
      \;\simeq\;
      \mathrm{Singularities}
      \big(
        \orbisingularE,
        -
      \big)
    $
    is the hom-functor \eqref{HomInfinityGroupoids}
    out of the terminal object \eqref{SingTerm}.
    \end{enumerate}
     \vspace{-2mm}
    Using this, we deduce the claim:
    \begin{enumerate}
    \vspace{-2mm}
    \item
      Since
      colimits of representable $\infty$-functors are
      equivalent to the point
      (Lemma \ref{InfinityColimitOverRepresentableInfinityFunctorIsContractible})
      we have
      \vspace{-2mm}
      \begin{equation}
        \label{SinglrG}
        \mathrm{Snglr}
        \big(
          \orbisingularG
        \big)
        \;\simeq\;
        \ast
        \;\simeq\;
        \ast / G
        \,.
      \end{equation}
    \vspace{-8mm}
    \item
      Observing that \eqref{SingularitiesHomGroupoids} reduces
    to
    $
      \mathrm{Singularities}
      \big(
        \orbisingularE,
        \orbisingularG
      \big)
      \;\simeq\;
      \ast \!\sslash\! G
    $
    we have
      \vspace{-2mm}
    $$
      \mathrm{Smth}
      \big(
        \orbisingularG
      \big)
      \;\simeq\;
      \ast \!\sslash\! G
      \,.
    $$
     \vspace{-9mm}
    \item
    With this and by the various adjunctions we have,
    for $U \in \mathbf{H}_{\tiny\smooth}$
    any geometically contractible generator \eqref{GeometricallyContractibleGenerators}
    and $K \in \mathrm{Groups}^{\mathrm{fin}}$ any
    finite group,
    the following sequence of
    natural equivalences:
     \vspace{-2mm}
    $$
    \hspace{-2cm}
      \begin{aligned}
        \mathbf{H}
        \Big(
          U
            \times
          \orbisingularK
          ,
          \mathrm{OrbSnglr}
          (\ast \!\sslash\! G)
        \Big)
        &
        \simeq
        \mathbf{H}_{\tiny\smooth}
        \Big(\,
          \underset{
            \simeq
            U
            \times
            \mathrm{Smth}
            \big(
              \orbisingularK
            \big)
          }{
          \underbrace{
            \mathrm{Smth}
            \big(
              U
                \times
              \orbisingularK
            \big)
          }
          },
          \ast \!\sslash\! G
        \Big)
        \;\simeq\;
        \mathbf{H}_{\tiny\smooth}
        \big(
          U
          \times
          (\ast \!\sslash\! K)
          ,
          \underset{
            \simeq \mathrm{Disc}
            (\ast \sslash G)
          }{
          \underbrace{
            \ast \!\sslash\! G
          }
          }
        \big)
        \\
        & \simeq
        \mathrm{Groupoids}_{\infty}
        \big(\,
          \underset{
            \simeq \ast
          }{
          \underbrace{
            \mathrm{Shp}(U)
          }
          }
          \times
          (\ast \!\sslash\! K),
          \ast \!\sslash\! G
        \big)
        \;\simeq\;
        \mathrm{Singularities}
        \big(
          \orbisingularK
          ,
          \orbisingularG
        \big)
        \\
        & \simeq
        \mathbf{B}
        \Big(
          \underset{
            \simeq \ast
          }{
            \underbrace{
              \mathrm{Shp}(U)
            }
          }
          \times
          \orbisingularK
          ,
          \orbisingularG
        \Big)
        \;\simeq\;
        \mathbf{B}
        \Big(
          \mathrm{Shp}
          \big(
            U
            \times
            \orbisingularK
          )
          ,
          \orbisingularG
        \Big)
        \\
        & \simeq
        \mathbf{H}
        \Big(
          U
          \times
          \orbisingularK
          ,
          \mathrm{Disc}
          \big(
            \orbisingularG
          \big)
        \Big)
        \;\simeq\;
        \mathbf{H}
        \Big(
          U
          \times
          \orbisingularK
          ,
          \orbisingularG
        \Big)
        \,,
      \end{aligned}
    $$

     \vspace{-3mm}

    \noindent where in several steps we recognized geometric discreteness,
    by \eqref{DiscreteEquivalences} in
    Remark \ref{DiscreteObjectsInSingularCohesion}.
   \end{enumerate}

 \vspace{-3mm}
\noindent    But, by Lemma \ref{YonedaOnProductSite}, this
    chain of natural equivalences in total
    a natural equivalence of the form
     \vspace{-2mm}
    $$
      \mathbf{H}
      \Big(
        y
        \big(
          U,
          \orbisingularK
        \big),
        \mathrm{OrbSnglr}
        \big(
          \ast \!\sslash\! G
        \big)
      \Big)
      \;\simeq\;
      \mathbf{H}
      \Big(
        y
        \big(
          U,
          \orbisingularK
        \big),
        \orbisingularG
      \Big)
      \,.
    $$

     \vspace{-2mm}
\noindent
    From this, the $\infty$-Yoneda embedding (Prop. \ref{InfinityYonedaEmbedding})
    implies that
    $
      \mathrm{OrbSnglr}
      \big(
        \ast \sslash G
      \big)
      \;\simeq\;
      \orbisingularG
      \,.
    $
  \hfill  \end{proof}

  It is useful to re-express this in terms of the modalities:
  \begin{prop}[Orbi-singularities are orbi-singular]
    \label{OrbiSingularitiesAreOrbiSingular}
    Let $\mathbf{H}$ be a singular-cohesive $\infty$-topos
    (Def. \ref{SingularCohesiveInfinityTopos})
    and consider a finite group $G \in \mathrm{Groups}^{\mathrm{fin}}$
    \eqref{InclusionOfFiniteGroups}.
    Then the images of the generic orbi-singularity
    $\orbisingularG \in \mathbf{H}$
    \eqref{SingularitiesUnderYoneda}
    under the modalities \eqref{SingularityModalities} are
    (see \hyperlink{FigureD}{Figure D}):
    \begin{equation}
      \label{GenericOrbisingularitiesUnderSingularModalities}
      \singular
      \big(
        \orbisingularG
      \big)
      \;\simeq\;
      \ast
      \,,
      \phantom{AAAA}
      \smooth
      \big(
        \orbisingularG
      \big)
      \;\simeq\;
      \ast \!\sslash\! G
      \,,
      \phantom{AAAA}
      \orbisingular
      \big(
        \orbisingularG
      \big)
      \;\simeq\;
      \orbisingularG
      \,.
    \end{equation}
  \end{prop}
\begin{proof}
This follows directly with Lemma \ref{ImagesAndPreimagesOfOrbisingularities} and the definition \eqref{SingularityModalities}.
For example:
 \vspace{-2mm}
    $$
        \orbisingular
        \big(
          \orbisingularG
        \big)
        \;\simeq\;
        \underset{
          \simeq
          \orbisingularG
        }{
          \underbrace{
          \mathrm{OrbSnglr}
            \circ
          \underset{
            \simeq
            \,
            \ast \!\sslash\! G
          }{
            \underbrace{
              \mathrm{Smth}
              \big(
                \orbisingularG
              \big)
            }
          }
          }
        }
    $$

\vspace{-.4cm}

\hfill \end{proof}

In the same vein, we also have the following immediate but important property:

\begin{prop}[Orbi-singularities are geometrically discrete]
    \label{OrbiSingularitiesAreGeometricallyDiscrete}
    Let $\mathbf{H}$ be a singular-cohesive $\infty$-topos
    (Def. \ref{SingularCohesiveInfinityTopos})
    and consider a finite group $G \in \mathrm{Groups}^{\mathrm{fin}}$
    \eqref{InclusionOfFiniteGroups}.

    \vspace{-1mm}
  \item {\bf (i)}   Then the basic orbi-singularity
    $\orbisingularG
     \in
    \mathbf{H}$
    \eqref{SingularitiesUnderYoneda}
    is geometrically discrete
    \eqref{CohesiveModalitiesFromAdjointQuadruple}
    and thus also pure shape:
    \begin{equation}
      \label{OrbiSingularityIsGeometricallyDiscrete}
        \flat
        \orbisingularG
        \;\simeq\;
        \orbisingularG
        \;,
   \qquad
        \mbox{\rm\textesh}
        \orbisingularG
        \;\simeq\;
        \orbisingularG
        \;.
    \end{equation}
    \item {\bf (ii)}   The same is true for
    $\mathrm{Smth}(\ast \!\sslash\! G)_{\tiny\orbisingular} \simeq \ast \!\sslash\! G$:
    \begin{equation}
      \label{OrbiSingularityIsGeometricallyDiscrete}
      \begin{aligned}
        \flat
        (\ast \!\sslash\! G)
        \;\simeq\;
        \ast \!\sslash\! G
        \;,
        \qquad
        \mbox{\rm\textesh}
        (\ast \!\sslash\! G)
        \;\simeq\;
        \ast \!\sslash\! G
        \;.
      \end{aligned}
    \end{equation}
\end{prop}
\begin{proof}
  Both statements follow immediately from the definitions and the
  fact that $G$ is finite and hence geometrically discrete
  \eqref{InclusionOfFiniteGroups}.
\hfill \end{proof}

\begin{remark}[Need for discrete/finite groups in $\mathrm{Singularities}$]
  \label{NeedForDiscreteGroupsInSingularities}
  It is to make
  Lemma \ref{ImagesAndPreimagesOfOrbisingularities}
  and hence Prop. \ref{OrbiSingularitiesAreOrbiSingular}
  true
  that Def. \ref{CategoryOfSingularities} requires the
  global orbit category $\mathrm{Singularities}$ to consist of
  finite groups, instead of more general compact Lie groups
  (Remark \ref{GlobalOrbitCategory}):
  If $\mathrm{Singularities}$ were to contain non-discrete
  compact Lie groups $G$, then the same argument as in Lemma
  \ref{ImagesAndPreimagesOfOrbisingularities} would give in
  \eqref{GenericOrbisingularitiesUnderSingularModalities}
  the following more general formula:
  $$
    \smooth \orbisingularG \;\;\simeq\;\; \ast \!\sslash\! \flat G
  $$
  (where on the right we think of the Lie group $G$
  as being cohesive via \eqref{ConcreteSmoothInfinityGroupoidsAreDiffeologicalSpaces}).
  Since the condition $G \simeq \flat G$
  characterizes discrete groups, this
  would break Prop. \ref{ShapeOfOrbisingularizedTopologicalGroupoidIsOrbispace}
  below, in that then the shape of the orbi-singularization of
  a topological groupoid would take non-traditional values on
  non-discrete groups in the global orbit category.
\end{remark}

The following lemma further illustrates the nature of orbi-singular cohesion:

\begin{lemma}[Smooth 0-truncated objects are orbi-singular]
  \label{Smooth0TruncatedObjectsAreOrbiSingular}
  Let $\mathbf{H}$ be a singular-cohesive $\infty$-topos
  (Def. \ref{SingularCohesiveInfinityTopos}).
  Then if $X \in \mathbf{H}_{\tiny\smooth,0}$
  is smooth \eqref{SingularityModalities}
  and 0-truncated
  (Def. \ref{nTruncatedObjects}),
  it is also orbi-singular \eqref{SingularityModalities}:
  \begin{equation}
    \tau_0(X)
    \simeq
    X
    \;\;\;
    \mbox{and}
    \;\;\;
    \smooth(X)
    \simeq
    X
    \phantom{AAAA}
    \Rightarrow
    \phantom{AAAA}
    \orbisingular(X)
    \;\simeq\;
    X
    \,.
  \end{equation}
\end{lemma}
\begin{proof}
  Since $X$ is smooth, there exists
  $X_{\tiny\smooth} \in \mathbf{H}_{\tiny\smooth}$ such that
  $
    X \simeq \mathrm{Smth}(X_{\tiny\smooth})
    \,.
  $
  Observe that $X$ being
  0-truncated implies that $X_{\tiny\smooth}$ is 0-truncated,
  (by using in Def.
  \ref{nTruncatedObjects}
  the hom-equivalence \eqref{AdjunctionHomEquivalence}
  of the right adjoint $\mathrm{Smth}$).

  Now let $\mathcal{S}$ be any site
  \eqref{ToposReflectionInPresheaves}
  for $\mathbf{H}_{\tiny\smooth}$.
  Then, for
  $U \in \mathcal{S} \overset{}{\hookrightarrow} \mathbf{H}_{\tiny\smooth}$
  and $G \in \mathrm{Groups}^{\mathrm{fin}}$, we have the
  following sequence of natural equivalences,
  using the various adjoint functors, their idempotency
  and respect for products:

$\phantom{A}$
\vspace{-4mm}
  $$
    \begin{aligned}
      (\orbisingular X)
      \big(
        \mathrm{Smth}(U) \times \orbisingularG
      \big)
      &
      \simeq
      (\orbisingular \mathrm{Smth}(X_{\tiny\smooth}))
      \big(
        \mathrm{Smth}(U) \times \orbisingularG
      \big)
      \\
      &
      \simeq
      \mathbf{H}
      \big(
        \mathrm{Smth}(U) \times \orbisingularG,
        \orbisingular \mathrm{Smth}(X_{\tiny\smooth})
      \big)
      \\
      & \simeq
      \mathbf{H}
      \big(
        \mathrm{Smth}(U) \times \smooth (\orbisingularG),
        \mathrm{Smth}(X_{\tiny\smooth})
      \big)
      \\
      & \simeq
      \mathbf{H}
      \big(
        \mathrm{Smth}(U \times (\ast \!\sslash\! G) ),
        \mathrm{Smth}(X_{\tiny\smooth})
      \big)
      \\
      & \simeq
      \mathbf{H}_{\tiny\smooth}
      \big(
        U \times (\ast \!\sslash\! G),
        X_{\tiny\smooth}
      \big)
      \\
      & \simeq
      \mathrm{Groupoids}_\infty
      \Big(
        \ast \!\sslash\! G,
        \,
        \mathbf{H}_{\tiny\smooth}
        \big(
          U,
          \,
          X_{\tiny\smooth}
        \big)
      \Big)
      \\
      & \simeq
      \mathrm{Groupoids}_\infty
      \Big(
        \ast,
        \,
        \mathbf{H}_{\tiny\smooth}
        \big(
          U,
          \,
          X_{\tiny\smooth}
        \big)
      \Big)
      \\
      & \simeq
      \mathbf{H}_{\tiny\smooth}
      \big(
        U,
        \,
        X_{\tiny\smooth}
      \big)
      \\
      & \simeq
      \mathbf{H}
      \big(
        \mathrm{Smth}(U),
        \,
        \mathrm{Smth}(X_{\tiny\smooth})
      \big)
      \\
      & \simeq
      \mathbf{H}
      \big(
        \mathrm{Smth}(U) \times \orbisingularG,
        \mathrm{Smth}(X_{\tiny\smooth})
      \big)
      \\
      & \simeq
      \mathbf{H}
      \big(
        \mathrm{Smth}(U) \times \orbisingularG,
        X
      \big)
      \\
      & \simeq
      X
      \big(
        \mathrm{Smth}(U) \times \orbisingularG
      \big).
    \end{aligned}
  $$
  Here the first and the last step
  use the $\infty$-Yoneda embedding (Prop. \ref{InfinityYonedaEmbedding}),
  while
  the middle step uses the fact that $X_{\tiny\smooth}$
  is 0-truncated,
  hence that $\mathbf{H}_{\tiny\smooth}(U,X_{\tiny\smooth})$
  is 0-truncated (i.e. a set), to find that there is in fact no dependency
  on $G$.
  Hence
  the claim follows by the $\infty$-Yoneda embedding (Prop. \ref{InfinityYonedaEmbedding}), in view of Lemma \ref{YonedaOnProductSite}.
\hfill \end{proof}

\begin{remark}[Degenerate case of orbi-singular]
  The natural language statement of
  Lemma \ref{Smooth0TruncatedObjectsAreOrbiSingular}
  shows that the modality $\orbisingular$ ``orbi-singular''
  \eqref{CohesionSingular} really means:
  ``All singularities \emph{that are present} are orbi-singularities.'',
  which becomes a trivially satisfied condition
  when there are no singularities,
  such as for smooth and 0-truncated objects.
\end{remark}

\medskip

\noindent {\bf Interplay between geometric and singular cohesion.}

\begin{lemma}[Smooth commutes with shape]
  \label{SmoothCommutesWithShape}
  In a singular-cohesive $\infty$-topos (Def. \ref{SingularCohesiveInfinityTopos})
  the smooth-modality \eqref{SingularityModalities}
  commutes with all three
  cohesive modalities
  \eqref{CohesiveModalitiesFromAdjointQuadruple}
  (as per Prop. \ref{OverSingularitiesCohesion}):
  \vspace{-2mm}
  $$
    \smooth \circ \mbox{\textesh}
    \;\simeq\;
    \mbox{\textesh} \circ \smooth
    \,,
    \phantom{AAA}
    \smooth \circ \flat
    \;\simeq\;
    \flat \circ \smooth
    \,,
    \phantom{AAA}
    \smooth \circ \sharp
    \;\simeq\;
    \sharp \circ \smooth\;.
  $$
\end{lemma}
\begin{proof}
Under the defining identification
$\mathbf{H} \simeq \mathrm{Sheaves}_\infty\big( \mathrm{Singularities}, \mathbf{H}_{\tiny\smooth}\big)$, let $\mathcal{X} \in \mathbf{H}$
be any object regarded as a $\mathbf{H}_{\tiny\smooth}$-valued
$\infty$-presheaf on $\mathrm{Singularities}$:
\vspace{-2mm}
$$
  \mathcal{X}
  \;:\;
  \orbisingularK
  \;\longmapsto\;
  \mathcal{X}(\orbisingularK)
  \;\in\;
  \mathbf{H}_{\tiny\smooth}
  \,.
$$

\vspace{-2mm}
\noindent
Observe then
(by Example \ref{InfinityCohesivePresheafSite}
via Lemma \ref{SingularitiesIsSiteForHomotopicalCohesion})
that $\smooth$
turns such a presheaf into the
constant presheaf on its value at
the terminal object $\orbisingularE$:
$$
  \big(
    \smooth \mathcal{X}
  \big)
  \;:\;
  \orbisingularK
  \;\longmapsto\;
  \mathcal{X}(\orbisingularE)
  \,.
$$
On the other hand, the geometric modalities operate
objectwise over $\mathrm{Singularities}$ (Remark \ref{ObjectwiseApplicationOfSingularityModalities}):
$$
  \big( \mbox{\textesh} \mathcal{X} \big)
  \;:\;
  \orbisingularK
  \;\longmapsto\;
  \mbox{\textesh}
  \big(
    \mathcal{X}(\orbisingularK)
  \big)
  \,.
$$

\vspace{-2mm}
\noindent
With this, we have the following sequence of
natural equivalences for
$\mathcal{X} \;\in\; \mathbf{H}$ and
$\orbisingularK \;\in\; \mathrm{Singularities}$:
\vspace{-2mm}
$$
  \begin{aligned}
    \big(
      \smooth \mbox{\textesh} \mathcal{X}
    \big)(\orbisingularK)
    &
    \simeq
    \big(
      \mbox{\textesh} \mathcal{X}
    \big)(\orbisingularE)
    \\
    & \simeq
      \mbox{\textesh}
      \big(
        \mathcal{X}(\orbisingularE)
      \big)
    \\
    & \simeq
      \mbox{\textesh}
      \big(
        (\smooth \mathcal{X})(\orbisingularK)
      \big)
    \\
    & \simeq
      \big(
        \mbox{\textesh} \smooth \mathcal{X}
      \big)(\orbisingularK)
      \,.
  \end{aligned}
$$

\vspace{-2mm}
\noindent
Hence the claim follows by the $\infty$-Yoneda embedding
(Prop. \ref{InfinityYonedaEmbedding}).
The argument for $\flat$ and $\sharp$ is analogous.
\hfill \end{proof}

\begin{remark}[Dichotomy between naive and proper orbifold cohomology via singular-cohesion]
  \label{RootOfProperAsComparedToBorelEquivariantCohomology}
  In contrast to Lemma \ref{SmoothCommutesWithShape},
  the orbi-singular modality $\orbisingular$ \eqref{SingularityModalities}
  does \emph{not} commute with the cohesive shape modality
  $\mbox{\textesh}$ \eqref{CohesiveModalitiesFromAdjointQuadruple},
  in general.
  This phenomenon is the very source of the \emph{proper equivariant}
  structure seen in singular-cohesive $\infty$-toposes,
  reflected in the following dichotomy between geometric- and
  homotopy fixed points of an orbi-space
  and in the distinction between proper- and Borel-equivariant
  cohomology:

\vspace{-2mm}
{\small
  \begin{center}
  \begin{tabular}{l|c||c|l}
    \hline
    &
    $\orbisingular \circ \mbox{\textesh}$
    &
    $\mbox{\textesh} \circ \orbisingular$
    &
    \\
    \hline
    \hline
    Def. \ref{GeometricAndHomotopyFixedPoints} {\bf (i)}
    &
    \begin{tabular}{c}
      Homotopy
      \\
      fixed-points
    \end{tabular}
    &
    \begin{tabular}{c}
      Geometric
      \\
      fixed-points
    \end{tabular}
    &
    Def. \ref{GeometricAndHomotopyFixedPoints} {\bf (ii)}
    \\
    \cline{2-3}
    Def. \ref{BorelEquivariantCohomologyInSingularCohsion}
    &
    \begin{tabular}{c}
      Borel-equivariant
      \\
      cohomology
    \end{tabular}
    &
    \begin{tabular}{c}
      Proper equivariant
      \\
      cohomology
    \end{tabular}
    &
    Def. \ref{ProperEquivariantCohomologyInSingularCohsion}
    \\
    \cline{2-3}
    Def. \ref{TangentiallyTwistedCohomology}
    &
    \begin{tabular}{c}
      Tangentially twisted
      \\
      cohomology
    \end{tabular}
    &
    \begin{tabular}{c}
      Tangentially twisted
      \\
      proper orbifold cohomology
    \end{tabular}
    &
    Def. \ref{OrbifoldCohomologyTangentiallyTwisted}
    \\
    \hline
  \end{tabular}
  \end{center}
  }

\end{remark}

\begin{defn}[Geometric- and homotopy-fixed points]
  \label{GeometricAndHomotopyFixedPoints}
  Let $\mathbf{H}$ be a singular-cohesive $\infty$-topos
  (Def. \ref{SingularCohesiveInfinityTopos}),
  $G \in \mathrm{Groups}(\mathbf{H})$
  (Prop. \ref{LoopingAndDelooping})
  being discrete $G \simeq \flat G$ and
  0-truncated $G \simeq \tau_0 G$, and
  $(X, \rho) \in G\mathrm{Actions}(\mathbf{H})$
  (Prop. \ref{InfinityAction})
  with smooth $X \simeq \smooth X$, hence
   \vspace{-2mm}
  $$
    X
    \;\in\;
    \xymatrix{
      \mathbf{H}_{\smooth}
    \;\;  \ar@{^{(}->}[rr]^-{
       \scalebox{.6}{$\mathrm{NnOrbSinglr}$ }
       }
      &&
      \mathbf{H}\;.
    }
  $$

 \vspace{-2mm}
  \noindent For any subgroup $K \subset G$,
  the $\infty$-groupoid of
  $\orbisingularK$-points
  in the slice (Prop. \ref{SliceInfinityTopos}) over $\orbisingularG$
  \eqref{GeneralOrbiSingularities}...

       \item {\bf (i)}
      ...of the orbi-singularization \eqref{CohesionSingular}
      of the shape \eqref{CohesionOfSingularCohesive}
      of $X \!\sslash\! G$
      is the \emph{homotopy fixed point space} of $X$
       \vspace{-2mm}
      \begin{equation}
        \label{HomotopyFixedPointSpaces}
        \mathrm{HmtpFxdPntSpc}^K(X)
        \;:=\;
        \mathbf{H}_{\!\big/\scalebox{.7}{$\orbisingularG$} }
        \left(
          \orbisingularK
          ,
          \;
          \orbisingular
          \,
          \raisebox{1pt}{\rm\textesh}
          \,
          (X \!\sslash\! G)
        \right).
      \end{equation}
      \vspace{-4mm}

    \item {\bf (ii)}
      ...of the shape \eqref{CohesionOfSingularCohesive}
      of the orbi-singularization \eqref{CohesionSingular}
      of $X \!\sslash\! G$
      is the \emph{geometric fixed point space} of $X$
       \vspace{-2mm}
      \begin{equation}
        \label{GeometricFixedPointSpaces}
        \mathrm{GmtrcFxdPntSpc}^K(X)
        \;:=\;
        \mathbf{H}_{\!\big/ \scalebox{.7}{$\orbisingularG$} }
        \left(
          \orbisingularK
          ,
          \;
          \raisebox{1pt}{\rm\textesh}
          \,
          \orbisingular
          \,
          (X \!\sslash\! G)
        \right).
      \end{equation}

        \vspace{-2mm}
   \noindent  On the right we are using
  Prop. \ref{OrbiSingularitiesAreOrbiSingular} and
  Prop. \ref{OrbiSingularitiesAreGeometricallyDiscrete}
  to see that both expressions indeed live in the slice over
  $\orbisingularG$.
\end{defn}

\medskip

\begin{prop}[Homotopy-fixed point spaces are fixed loci in shapes]
  \label{HomotopyFixedPointSpacesAreFixedLociInShapes}
  The homotopy-fixed point spaces
  \eqref{HomotopyFixedPointSpaces}
  of the $G$-space $X$ in Def. \ref{GeometricAndHomotopyFixedPoints}
  are, equivalently, the fixed-loci (Def. \ref{FixedPoints})
  of the shape $\mathrm{Shp}(X) \;\in\; \mathrm{Groupoids}_\infty$
  \eqref{AdjunctionCohesion} of X :
     \vspace{-2mm}
  \begin{equation}
    \label{HomotopyFixedPointSpaceAsFixedLocusInShape}
    \hspace{-2.5cm}
    \mathrm{HmtpFxdPntSpcs}^K(X)
    \;\simeq\;
    \big(
      \mathrm{Shp}(X)
    \big)^K
    \mathrlap{
      \;\;\;\;\;
      \in
      \mathrm{Groupoids}_\infty
    }
  \end{equation}
  with respect to the induced $G \simeq \raisebox{1pt}{\rm\textesh} G$-action
  (using Prop. \ref{ShapePreservesHomotopyColimitsByDiscreteGroups},
  discreteness of $G$ and cohesion in the form of Prop. \ref{CompositeCohesiveModalities}).
\end{prop}
\begin{proof}
  We claim a sequence of natural equivalences as follows:
  \vspace{-1mm}
  \begin{equation}
    \label{TowardsUnderstandingHomotopyFixedPointSpaces}
    \begin{aligned}
      \mathrm{HmtpFxdPntSpc}^K(X)
      & :=\;
      \mathbf{H}_{\!\big/\scalebox{.7}{$\orbisingularG$} }
      \left(
        \orbisingularK
        ,
        \;
        \orbisingular
        \,
        \raisebox{1pt}{\rm\textesh}
        \,
        (X \!\sslash\! G)
      \right)
      \\
      & \simeq
      \mathbf{H}_{\!\big/\scalebox{.7}{$\orbisingularG$} }
      \left(
        \orbisingularK
        ,
        \;
        \orbisingular
        \,
        (\raisebox{1pt}{\rm\textesh} X)
        \!\sslash\!
        G
      \right)
      \\
      & \simeq\;
      \mathbf{H}_{\!\big/\scalebox{.7}{$\mathrm{OrbSnglr}(\ast \!\sslash\! G)$} }
      \Big(
        \mathrm{OrbSinglr}
        \big(
          \ast \!\sslash\! K
        \big)
        ,
        \;
        \mathrm{OrbSnglr}
        \big(
          (\raisebox{1pt}{\rm\textesh} X)
          \!\sslash\!
          G
        \big)
      \Big)
      \\
      & \simeq
      \big(
        \mathbf{H}_{\tiny\smooth}
      \big)_{\!\big/\scalebox{.7}{$\ast \!\sslash\! G$} }
      \big(
        \ast \!\sslash\! K
        ,
        \;
        (\raisebox{1pt}{\rm\textesh} X)
        \!\sslash\!
        G
      \big)
      \\
      & \simeq
      \big(
        \mathrm{Groupoids}_\infty
      \big)_{\!\big/\scalebox{.7}{$\ast \!\sslash\! G$} }
      \big(
        \ast \!\sslash\! K
        ,
        \;
        \mathrm{Shp}(X)
        \!\sslash\!
        G
      \big)
      \\
      & \simeq
      \big(
        \mathrm{Groupoids}_\infty
      \big)_{\!\big/\scalebox{.7}{$\ast \!\sslash\! K$} }
      \big(
        \ast \!\sslash\! K
        ,
        \;
        \mathrm{Shp}(X)
        \!\sslash\!
        K
      \big)
      \\
      & \simeq
      \big(\mathrm{Shp}(X))^K
    \end{aligned}
  \end{equation}
  Here the first step is the definition \eqref{HomotopyFixedPointSpaces},
  and the second step uses Prop. \ref{ShapePreservesHomotopyColimitsByDiscreteGroups},
  discreteness of $G$ and cohesion in the form of Prop. \ref{CompositeCohesiveModalities}.
  In the third step we observe with
  $\orbisingularK \;\simeq\; \orbisingular( \ast \!\sslash\! K )$
  (Lemma \ref{ImagesAndPreimagesOfOrbisingularities})
  and $\orbisingular \;:=\; \mathrm{OrbSinglr} \circ \mathrm{Smth}$
  \eqref{SingularityModalities}
  that all objects and morphisms are in the image of
  $\mathrm{OrbSnglr}$, and in the fourth step we use that this
  functor is fully faithful, by Prop. \ref{OverSingularitiesCohesion}.
  In the fifth step, we similarly observe that
  all objects and morphisms are, in fact, furthermore in the image
  of $\mathrm{Disc}$ (by assumption on $G$ and by definition
  of $\raisebox{1pt}{\textesh} := \mathrm{Disc} \circ \mathrm{Shp}$
  \eqref{CohesiveModalitiesFromAdjointQuadruple}), which is fully faithful
  by the axioms of cohesion \eqref{AdjunctionCohesion}.
  The sixth step observes the universal factorization through the
  pullback
   \vspace{-1mm}
  $$
    \raisebox{20pt}{
    \xymatrix@R=1em{
      \ast \!\sslash\! K
      \ar[dr]
      \ar@{-->}[rr]
      &&
      \mathrm{Shp}(X) \!\sslash\! G
      \ar[dl]
      \\
      & \ast \!\sslash\! G
    }
    }
    \phantom{AAA}
    \simeq
    \phantom{AAA}
    \raisebox{30pt}{
    \xymatrix@C=15pt@R=.5em{
      \ast \!\sslash\! K
      \ar@{-->}[rr]
      \ar@{=}[dr]
      &&
      \mathrm{Shp}(X) \!\sslash\! K
      \ar[dr]
      \ar[dl]
      \ar@{}[dd]|-{\mbox{\tiny(pb)}}
      \\
      &
      \ast \!\sslash\! K
      \ar[dr]
      &&
      \mathrm{Shp}(X) \!\sslash\! G
      \ar[dl]
      \\
      &
      & \ast \!\sslash\! G
    }
    }
  $$

   \vspace{-2mm}
\noindent
  The pullback, in turn, is the homotopy quotient of the
  restricted action, as shown, by Prop. \ref{PullbackAction}.
  With this,
  the last step follows by Example \ref{FixedLociInInfinityGroupoids}.
    In summary, the composite of the sequence of equivalences
    \eqref{TowardsUnderstandingHomotopyFixedPointSpaces}
  gives the statement \eqref{HomotopyFixedPointSpaceAsFixedLocusInShape}.
\hfill \end{proof}

\begin{example}[Geometric fixed points generally differ from homotopy fixed points]
 \label{GeometricFixedPointSpacesDiffer}
 As in Example \ref{OrbiSingularSmoothInfinityGroupoids},
  let $\mathbf{H} := \mathrm{SingularSmoothGroupoids}_\infty$.
      For $n \in \mathbb{N}$, $n \geq 1$, consider
  the Cartesian space
  $
    \mathbb{R}^n
    \in
         \mathrm{SmoothManifolds}
     \longhookrightarrow
      \mathbf{H}
     $,
  via \eqref{ConcreteSmoothInfinityGroupoidsAreDiffeologicalSpaces},
  and regard it as equipped with the additive translation action of
  $\mathbb{Z}^n$ induced from the left action of
  the additive group $(\mathbb{R}^n, +)$ on itself,
  under the canonical inclusion
  $(\mathbb{Z}^n,+) \hookrightarrow (\mathbb{R}^n, +)$:
   \vspace{-1mm}
  \begin{equation}
    \label{ActionOnRnByZn}
    (\mathbb{R}^n, \rho_{\ell})
    \;\in\;
    \mathbb{Z}^n\mathrm{Actions}(\mathbf{Hs})\;.
  \end{equation}

   \vspace{-2mm}
\noindent   So the quotient of this action
  $
    \mathbb{R}^n \!\sslash\! \mathbb{Z}^n
    \;\simeq\;
    \mathbb{R}^n / \mathbb{Z}^n
    \;\simeq\;
    \mathbb{T}^n
     \in
          \mathrm{SmoothManifolds}
    \longhookrightarrow
      \mathbf{H}
     $
  is the standard $n$-torus.
  We then have for the two notions of fixed-point spaces
  from Def. \ref{GeometricAndHomotopyFixedPoints}:

  \noindent
  {\bf (i)} The {\it Homotopy-fixed point space}
  \eqref{HomotopyFixedPointSpaces} of the action \eqref{ActionOnRnByZn}
  is equivalently the point
  (by Prop. \ref{HomotopyFixedPointSpacesAreFixedLociInShapes} and
  \eqref{ChartsForSmoothGroupoids}):
  \vspace{-1mm}
  $$
         \mathrm{HmtpFxdPntSpc}^{\mathbb{Z}^n}(\mathbb{R}^n)
                \simeq
        \big(
          \underset{
            \simeq \, \ast
          }{
          \underbrace{
            \raisebox{2pt}{\textesh}
            \mathbb{R}^n
          }
          }
        \big)^{\mathbb{Z}^n}
     \simeq \, \ast
    $$

 \vspace{-2mm}
  \noindent {\bf (ii)} The {\it geometric fixed point space} \eqref{GeometricFixedPointSpaces}
  of the action \eqref{ActionOnRnByZn} is empty
   \vspace{-1mm}
  $$
          \mathrm{GmtrcFxdPntSpc}^{\mathbb{Z}^n}(\mathbb{R}^n)
            \simeq\;
      \big(
        \mathbb{R}^n
      \big)^{\mathbb{Z}^n}
\simeq
      \;
      \varnothing
  $$

   \vspace{-1mm}
\noindent
  This follows by Lemma \ref{ShapeOfOrbiSingularSpacesAsPresheafOnSingularities},
  using that no element of the set underlying $\mathbb{R}^n$ is
  fixed by the action of $\mathbb{Z}^n$.
\end{example}

\medskip

\newpage

\section{Orbifold geometry}
\label{OrbifoldGeometry}

Within an ambient context of singular-cohesive homotopy theory
(\cref{SingularCohesiveGeometry}),
we now formulate the two geometric aspects of orbifolds:

- \cref{Orbispaces} -- as cohesive spaces with orbi-singularities,

- \cref{VFolds} -- as cohesive spaces locally equivalent to a given model space.

\noindent In the end, we combine both aspects to form the
\emph{proper $\infty$-categories of orbifolds}:
this is Def. \ref{OrbiVFolds} below.

\subsection{Orbispaces}
\label{Orbispaces}

We observe (Prop. \ref{ShapeOfOrbisingularizedTopologicalGroupoidIsOrbispace}) that the
shape of the orbi-singularization of a topological groupoid,
regarded in singular-smooth homotopy theory (Example \ref{OrbiSingularSmoothInfinityGroupoids}),
is the corresponding \emph{orbispace} in
global equivariant homotopy theory.

\begin{remark}[Orbispaces in topology and in global equivariant homotopy theory]
  \label{TraditionalWayOfRegardingTopologicalGroupoidsAsOrbifolds}
$\,$
\\

\vspace{-.5cm}

\noindent {\bf (i) Orbispaces in topology.}
The term \emph{orbispace}
was originally introduced \cite{Haefliger90} to mean
the topological version of orbifolds,
i.e., Satake's original concept \cite{Satake56} but
disregarding any differentiable structure. From the perspective
of {\'e}tale groupoids/stacks, this means to consider
topological groupoids/stacks instead of Lie groupoids/differentiable stacks.
So this usage of the term ``orbispace'' serves to complete the
following table:
\begin{center}
\begin{tabular}{|c||c|l}
  \cline{1-2}
  Smooth manifold & Topological manifold
  \\
  \cline{1-2}
  \cline{1-2}
  orbifold & {\color{darkblue} orbispace}
  &
  \scalebox{.8}{
    ({\it geometric sense})
  }
  \\
  \cline{1-2}
  Lie groupoid & topological groupoid
  \\
  \cline{1-2}
  differentiable stack & topological stack
  \\
  \cline{1-2}
\end{tabular}
\end{center}
In this sense, orbispaces have been discussed, e.g., in
\cite{Haefliger84}\cite[\S 5]{Haefliger91}\cite{Chen01}\cite{Henriques01}.

\vspace{1mm}
\noindent {\bf (ii) Orbispaces in global equivariant homotopy theory.}
In \cite{HenriquesGepner07} it was suggested to change perspective
and to instead regard these topological groupoids
$\mathcal{X}_{\mathrm{top}}$ via the
systems of homotopy types of all their geometric fixed point spaces,
by the following formula \cite[4.2]{HenriquesGepner07}
(beware the differing conventions, as per Remark \ref{GlobalOrbitCategory}):
\vspace{-2mm}
\begin{equation}
  \label{TopologicalGroupoidAsPresheafOnOrbitCategory}
  \xymatrix{
    G
    \ar@{}[rr]|-{
      \longmapsto
    }
    &&
    \overset{
      \mathclap{
      \raisebox{2pt}{
        \tiny
        \color{darkblue}
        \bf
        homotopy type of
        (fat) geometric realization of
      }
      }
    }{
    \left\Vert
      \underset{
        \mathclap{
        \mbox{
          \tiny
          \color{darkblue}
          \bf
          \begin{tabular}{c}
            topological
            mapping
            groupoid
          \end{tabular}
        }
        }
      }{
      \mathbf{Maps}
      \big(
        \mathbf{B}G
        \,,\,
        \mathcal{X}_{\mathrm{top}}
      \big)
      }
    \right\Vert
    }
  }
  \phantom{AAAAA}
  \mbox{
    \begin{tabular}{c}
      \color{darkblue}
      orbispace
      \\
      \scalebox{.8}{({\it equivariant homotopical sense})}
    \end{tabular}
  }
\end{equation}
This is a global-equivariant version of how
topological $G$-spaces are incarnated in $G$-equivariant homotopy
theory via Elmendorf's theorem (recalled as Prop. \ref{ElmendorfTheorem}),
and has served to motivate the development of
global equivariant homotopy theory \cite{Schwede18}.

In the course of this development,
homotopy theorists adopted the term ``orbispace''
to refer not to the topological groupoid $\mathcal{X}_{\mathrm{top}}$
(as \cite{Haefliger90} originally did)
but rather to the global equivariant homotopy type that is represented via
\eqref{TopologicalGroupoidAsPresheafOnOrbitCategory}.
Usage of the term \emph{orbispace} in this sense of global homotopy theory is,
after \cite{HenriquesGepner07},
in \cite{Rezk14}\cite{Koerschgen16}\cite{Schwede17}\cite[3]{Lurie19}\cite{Juran20}.
In \cite[3.15]{Juran20} formula \ref{TopologicalGroupoidAsPresheafOnOrbitCategory}
is used (following suggestions in \cite[Introd.]{Schwede17}\cite[p. ix-x]{Schwede18})
to define (abelian, non-geometric) cohomology of orbifolds
with coefficients in global equivariant spectra.

Our Prop. \ref{ShapeOfOrbisingularizedTopologicalGroupoidIsOrbispace}
below shows that these two different meanings of the term ``orbispace''
in the literature are disentangled as well as unified by the notion
of singular cohesion (Def. \ref{SingularCohesiveInfinityTopos}), in
that orbispaces in the sense {\bf (ii)} are
the shape $\raisebox{1pt}{\textesh}$ \eqref{AdjunctionCohesion}
of the orbi-singularization $\orbisingular$ \eqref{SingularityModalities}
of the topological groupoids in {\bf (i)}:
\begin{equation}
  \label{FromTopologicalGroupoidsToOrbispaces}
  \xymatrix@R=0pt{
    \mathrm{TopologicalGroupoids}
    \ar[rr]^-{
      \scalebox{.85}{$
        \raisebox{1pt}{\textesh}
        \,
        \circ
        \orbisingular
      $}
    }
    &&
    \mathrm{Orbispaces}
    \\
    \mathcal{X}_{\mathrm{top}}
    \ar[rr]
    &&
    \Big(
      \orbisingularG
      \;\mapsto\;
      \left\Vert
        \mathbf{Maps}
        \big(
          \mathbf{B}G,\mathcal{X}_{\mathrm{top}}
        \big)
      \right\Vert
    \Big)
  }
\end{equation}

\vspace{-2mm}
\noindent Hence Prop. \ref{ShapeOfOrbisingularizedTopologicalGroupoidIsOrbispace} below
means that, before passing to their pure shape, we may
think of the orbi-singularizations of objects in singular-cohesive
$\infty$-toposes as \emph{cohesive orbispaces}, lifting the
concept of plain orbispaces in the sense {\bf (ii)}
from plain homotopy theory to
geometric (differential, {\'e}tale) homotopy theory, hence back to
sense {\bf (i)} and beyond.
\end{remark}

\newpage

The crucial fact underlying the phenomenon \eqref{FromTopologicalGroupoidsToOrbispaces},
both in Prop. \ref{ShapeOfOrbisingularizedTopologicalGroupoidIsOrbispace}
and in Lemma \ref{ShapeOfOrbiSingularSpacesAsPresheafOnSingularities} below,
is that the probe of an orbi-singular object $\orbisingular \mathcal{X}_{\tiny\smooth}$
by a generic orbi-singularity $\orbisingularK$ \eqref{AnOrbifoldSingularity}
is, by adjunction \eqref{SingularityModalities},
equivalently the probe of the underlying smooth object
by the smooth aspect of $\orbisingularK$,
hence is, by \eqref{GenericOrbisingularitiesUnderSingularModalities} in
Prop. \ref{OrbiSingularitiesAreOrbiSingular},
the geometric $G$-fixed locus in $\mathcal{X}_{\tiny\smooth}$:
\vspace{-1mm}
\begin{equation}
  \label{MapsFromOrbisingularityIntoOrbisingularSpace}
  \xymatrix{
    \orbisingularG
    \ar[r]
    &
    \orbisingular \mathcal{X}_{\tiny\smooth}
  }
  \phantom{AAA}
  \overset{
    \raisebox{3pt}{
      \tiny
      \eqref{SingularityModalities}
    }
  }{
    \Leftrightarrow
  }
  \phantom{AAA}
  \xymatrix{
    \smooth \orbisingularG
    \ar[r]
    &
    \mathcal{X}_{\tiny\smooth}
  }
  \phantom{AAA}
  \overset{
    \raisebox{3pt}{
      \tiny
      \eqref{GenericOrbisingularitiesUnderSingularModalities}
    }
  }{
    \Leftrightarrow
  }
  \phantom{AAA}
  \xymatrix{
    \ast \!\sslash\! G
    \ar[r]
    &
    \mathcal{X}_{\tiny\smooth}
  }.
\end{equation}

\vspace{-2mm}
\noindent
Equivalently, since $\orbisingularG \;\simeq\; \orbisingular (\ast \!\sslash\! G)$
(Lemma \ref{ImagesAndPreimagesOfOrbisingularities}) the composite
corespondence \eqref{MapsFromOrbisingularityIntoOrbisingularSpace} is
fully-faithfulness of $\orbisingular$.

\medskip

\noindent {\bf Example: Topological groupoids as cohesive orbispaces.}
\begin{prop}[Shape of orbi-singularized topological groupoid is orbispace]
  \label{ShapeOfOrbisingularizedTopologicalGroupoidIsOrbispace}
  Let $\mathbf{H} := \mathrm{SingularSmoothGroupoids}_\infty$
  (Example \ref{OrbiSingularSmoothInfinityGroupoids}), and
  let
  \vspace{-3mm}
  $$
    \xymatrix@R=-2pt{
      \mathrm{TopologicalGroupoids}
      \;
      \ar[r]^-{\scalebox{.6}{$
        \mathrm{Cdfflg}$}
      }
      &
      \;\mathrm{SmoothGroupoids}_\infty\;
      \ar[rr]^-{\scalebox{.6}{$ \mathrm{NnOrbSnglr}$} }
      &&
      \mathbf{H}
      \\
      \mathcal{X}_{\mathrm{top}}
      \ar@{|->}[rrr]
      &&&
      \mathcal{X}_{\tiny\smooth}
    }
  $$

  \vspace{-2mm}
\noindent
  be a topological groupoid, regarded via the embeddings
  \eqref{1TrunConcreteSmoothInfinityGroupoidsAreDiffeologicalSpaces}
  and \eqref{CohesionSingular}.
    If $\mathcal{X}_{\tiny\smooth}$
  is such that both its space of objects
  and of morphisms are retracts of cell complexes
  (for instance: both are CW-complexes \eqref{CWComplexes})
  then the shape \eqref{CohesionOfSingularCohesive} of its
  orbi-singularization \eqref{SingularityModalities}
  is,
  as an $\infty$-presheaf \eqref{InfinitToposOverSingularities}
  of $\infty$-groupoids on $\mathrm{Singularities}$ \eqref{CategoryOfSingularities}
  (i.e., on the global orbit category, Remark \ref{GlobalOrbitCategory})
  \vspace{-3mm}
  $$
    \raisebox{1pt}{\rm\textesh} \orbisingular \mathcal{X}_{\tiny\smooth}
    \;\in\;
    \xymatrix{
      \mathrm{Sh}_\infty\big(
        \mathrm{Singularities}
      \big)
      \; \ar@{^{(}->}[rr]^-{ \mathrm{Disc} }
      &&
      \mathbf{H}
    }
  $$

  \vspace{-2mm}
\noindent
  given by the assignment \eqref{MapsFromOrbisingularityIntoOrbisingularSpace}
  \vspace{-0mm}
  \begin{equation}
    \label{FormulaForShapeOfOrbisingularizationOfTopologicalGroupoids}
    \raisebox{1pt}{\rm\textesh} \orbisingular \mathcal{X}_{\tiny\smooth}
    \;\;:\;\;
    \orbisingularG
    \;\longmapsto\;
    \left\Vert
      \mathbf{Maps}
      \big(
        \mathbf{B} G, \mathcal{X}_{\mathrm{top}}
      \big)
    \right\Vert
    \,,
  \end{equation}
  \noindent
  where on the right we have the fat geometric realization of the
  topological functor groupoid
  \cite{Segal74} (see \cite[2.3]{HenriquesGepner07}),
  with $\mathbf{B}G \;\simeq\; \ast \!\sslash\! G$ (Example \ref{DeloopingGroupoids})
  regarded as a finite topological groupoid.
\end{prop}
\begin{proof}
Recall from \eqref{ChartsForSmoothGroupoids} in
Example \ref{SmoothInfinityGroupoids}
that $\mathrm{Charts} := \mathrm{CartesianSpaces}$
(Def. \ref{CartesianSpaces}) is a site of cohesive charts
(Def. \ref{ChartsForCohesion})
for $\mathrm{SmoothGroupoids}_\infty$.
We claim that
for $\mathbb{R}^n \in \mathrm{CartesianSpaces}$
and $\orbisingularG \in \mathrm{Singularities}$ (Def. \ref{CategoryOfSingularities}),
hence $\mathbb{R}^n \times \orbisingularG \in
\mathrm{CartesianSpaces} \times \mathrm{Singularities}$ (Lemma \ref{YonedaOnProductSite}),
we have the following sequence of natural equivalences:
\begin{equation}
  \label{ValueOfGOrbiSingularSpacesOnRepresentables}
  \hspace{-1cm}
  \begin{aligned}
    \mathbf{H}
    \big(
      \mathbb{R}^n \times \orbisingularG,
      \,
      \orbisingular \mathcal{X}_{\tiny\smooth}
    \big)
    & =
    \mathbf{H}
    \big(
      \mathbb{R}^n \times \orbisingularG,
      \,
      \mathrm{ OrbSnglr }
      \big(
        \mathcal{X}_{\tiny\smooth}
      \big)
    \big)
    \\
    & \simeq
    \mathbf{H}_{\tiny\smooth}
    \big(
      \underset{
        \simeq \, \mathbb{R}^n \times \mathbf{B}G
      }{\;
      \underbrace{
        \mathrm{Smth}
        \big(
          \mathbb{R}^n \times \orbisingularG
        \big)
      }
      }
      \;,
      \,
      \mathcal{X}_{\tiny\smooth}
    \big)
    \\
    & \simeq
    \mathbf{H}_{\tiny\smooth}
    \Big(
      \mathbb{R}^n
      \,,\,
      \mathbf{Maps}
      \big(
        \mathbf{B}G
        \,,\,
        \mathcal{X}_{\tiny\smooth}
      \big)
    \Big)
    \\
    &
    \simeq
    \mathbf{H}_{\tiny\smooth}
    \Big(
      \mathbb{R}^n
      \,,\,
      \mathrm{Cdfflg}
      \mathbf{Maps}
      \big(
        \mathbf{B}G
        \,,\,
        \mathcal{X}_{\mathrm{top}}
      \big)
    \Big).
  \end{aligned}
\end{equation}
Here the first step is \eqref{GOrbiSingularSpaceAsOrbSnglrOfHomotopyQuotient},
the
second is the hom-equivalence \eqref{AdjunctionHomEquivalence}
of the adjunction $\mathrm{Smth} \dashv \mathrm{OrbSnglr}$
\eqref{CohesionSingular} and using under the brace that
$\mathrm{Smth}$ preserves products (by Prop. \ref{AdjointsPreserveCoLimits}),
that $\mathbb{R}^n$ is already smooth,
and that
$\mathrm{Smth}\big( \orbisingularG \big)
\simeq (\ast \!\sslash\! G )$ by \eqref{ImagedAndPreimagesofBasicOrbiSngularities}.
The third step is
Lemma \ref{MappingStackFromDiscreteGroupToTopologicalStack}.

\noindent
Since also the composite of all these natural equivalences
is thus natural, the $\infty$-Yoneda lemma (Prop. \ref{YonedaLemma})
implies that
$$
  \orbisingular \mathcal{X}_{\tiny\smooth}
  \;:\;
  \orbisingularK
  \;\longmapsto\;
  \mathrm{Cdfflg}
  \,
  \mathrm{Maps}
  \big(
    \mathbf{B}G
    \,,\,
    \mathcal{X}_{\tiny\smooth}
  \big)
  \,.
$$
Now, since
$\raisebox{1pt}{\textesh}$ acts objectwise over $\orbisingularK$
\eqref{CohesionOfSingularCohesive}, we find from this that
  \begin{eqnarray*}
    \raisebox{1pt}{\textesh} \orbisingular \mathcal{X}_{\tiny\smooth}
    \;:\;
    \orbisingularK
    &\longmapsto &
    \raisebox{1pt}{\textesh}
    \mathrm{Cdfflg}
    \,
    \mathbf{Maps}
    \big(
      \mathbf{B}G
      \,,\,
      \mathcal{X}_{\mathrm{top}}
    \big)
    \\
    & \simeq &
    \mathrm{Shp}_{\mathrm{sTop}}
    \Big(
    \mathbf{Maps}
    \big(
      \mathbf{B}G
      \,,\,
      \mathcal{X}_{\mathrm{top}}
    \big)
    \Big)
    \\
    & \simeq\; &
    \left\Vert
    \mathrm{Maps}
    \big(
      \mathbf{B}G
      \,,\,
      \mathcal{X}_{\mathrm{top}}
    \big)
    \right\Vert\;.
  \end{eqnarray*}
Here the first step is
\eqref{CohesiveShapeOfTopologicalGroupoids} and
the last step follow by
Prop. \ref{SimplicialTopologicalShapeOfDegreewiseCofibrantSimplicialSpaces}.
\hfill \end{proof}

\medskip

\noindent {\bf Cohesive $G$-orbispaces.}
We now discuss in more detail
the analogue of Prop. \ref{ShapeOfOrbisingularizedTopologicalGroupoidIsOrbispace}
in {\bf (a)} the special case of global quotient
stacks $\mathcal{X}_{\tiny\smooth} \simeq X \!\sslash\! G$
by a discrete group $G$, but {\bf (b)} in the full generality of
$X$ being any 0-truncated cohesive space
(not necessarily a topological space, but for instance
a smooth manifold or diffeological space \eqref{ConcreteSmoothInfinityGroupoidsAreDiffeologicalSpaces}
or even a non-concrete object).

\begin{remark}[Good orbifolds and good cohesive orbispaces]
  \label{GoodOrbifolds}
  The traditional orbifolds that arise as global quotients
  $\mathcal{X}_{\tiny\smooth} \simeq X \!\sslash\! G$
  of a smooth manifold $X$ by the action of a \emph{discrete} group $G$
  are called \emph{good orbifolds} (e.g. \cite[6]{Kapovich08}). Therefore, the cohesive
  $G$-orbispaces discussed now (Def. \ref{GOrbiSpace}) could be called
  (after forgetting their slicing over $\orbisingularG$)
  the \emph{good cohesive orbispaces}.
\end{remark}

\begin{defn}[Cohesive $G$-orbispace]
  \label{GOrbiSpace}
\  Let $\mathbf{H}$ be a singular-cohesive $\infty$-topos
  (Def. \ref{SingularCohesiveInfinityTopos})
  and $G \in \mathrm{Groups}(\mathbf{H})$ (Prop. \ref{LoopingAndDelooping})
  discrete $G \simeq \flat G$.
  We say that a \emph{cohesive $G$-orbispace} is an object
  \vspace{-2mm}
  $$
    \raisebox{15pt}{
    \xymatrix@R=.8em{
      \mathcal{X}
      \ar[d]^-{p}
      \\
      \orbisingularG
    }
    }
    \;\;\in\;\;
    \mathbf{H}_{/_{\!\!\raisebox{2pt}{\scalebox{.7}{$\orbisingularG$}}}}
  $$

  \vspace{-3mm}
  \noindent
  in the slice over the $G$-orbi-singularity \eqref{GeneralOrbiSingularities}
  that is:

\vspace{-2mm}
 \begin{equation}
   \label{ConditionsDefiningGOrbiSingularSpaces}
   \hspace{0cm}
   \mbox{
   \begin{tabular}{llll}
  {\bf (a)} & orbi-singular: &
   $
     \phantom{\mathclap{\vert_{\vert_{\vert_{\vert_{\vert_{\vert}}}}}}}
     \orbisingular(p) \simeq p
   $
   &
   (Def. \ref{SincularCohesionModalies})
   ,
   \\
     {\bf (b)} & 0-truncated: &
   $
     \phantom{\mathclap{\vert_{\vert_{\vert_{\vert_{\vert_{\vert}}}}}}}
     (\tau_0)_{/\scalebox{.7}{$\orbisingularG$}}(p) \simeq p
   $
   &
   (Def. \ref{nTruncatedObjects})
   .
 \end{tabular}
  }
  \end{equation}
\end{defn}

\begin{defn}[Universal covering space of a $G$-orbi-singular space]
  \label{CoveringSpaceOfGOrbiSingularSpace}
  Given a Cohesive $G$-orbispace
  $\mathcal{X} \in \mathbf{H}_{/\scalebox{.7}{$\orbisingularG$}}$
  (Def. \ref{GOrbiSpace}),
  we say that its \emph{universal covering space}
  $X \in \mathbf{H}$
  the homotopy fiber of the defining morphism to $\orbisingularG$
  over its essentially unique point:
  \vspace{-2mm}
  \begin{equation}
    \label{CohesiveCoveringSpace}
    \raisebox{13pt}{
    \xymatrix@R=1em{
      X \ar[r]^-{\mathrm{fib}(p)}
      &
      \mathcal{X}
      \ar[d]^-p
      \\
      & \orbisingularG
    }
    }
  \end{equation}
\end{defn}

\medskip

\begin{prop}[Properties of universal covering spaces]
  \label{PropertiesOfUniversalCoveringSpaces}
  Let $\mathbf{H}$ be a singular-cohesive $\infty$-topos (Def. \ref{SingularCohesiveInfinityTopos}).
  Given a $G$-orbi-singular space
  $\mathcal{X} \in \mathbf{H}_{/\scalebox{.7}{$\orbisingularG$}}$
  (Def. \ref{GOrbiSpace}),
  its universal covering space $X$ (Def. \ref{CoveringSpaceOfGOrbiSingularSpace})
    \begin{enumerate}[{\bf (i)}]
  \vspace{-2mm}
 \item is:

\vspace{-2mm}
   \hspace{4.2cm}
   \begin{tabular}{llll}
    {\bf (a)}
    &
    0-truncated:
    &
    $
      \phantom{\mathclap{\vert_{\vert_{\vert_{\vert_{\vert}}}}}}
      \tau_0(X) \simeq X
    $
    &
    (Def. \ref{nTruncatedObjects})\;,
    \\
    {\bf (b)}
    &
    smooth:
    &
    $\smooth(X) \simeq X$
    &
    (Def. \ref{CohesiveTopos})\;,
  \end{tabular}

 \vspace{-1mm}
   \item
  and is equipped with a $G$-action (Prop. \ref{InfinityAction})
  such that $\mathcal{X}$ is
  the orbi-singularization \eqref{SingularityModalities} of
  the corresponding homotopy quotient:
  \begin{equation}
    \label{GOrbiSingularSpaceAsOrbSnglrOfHomotopyQuotient}
    \mathcal{X}
    \;\simeq\;
    \orbisingular
    \big(
      X \!\sslash\! G
    \big)
    \,.
  \end{equation}
  \end{enumerate}
\end{prop}
\begin{proof}
 \noindent {\bf (i)}
   That $X$ is {\bf (a)} 0-truncated follows from
  the condition that $p$ is 0-truncated and using
  Lemma \ref{nTruncatedMorphismViannTruncatedHomotopyFiber}.
   To see that $X$ is \noindent {\bf (b)} smooth,
  observe that by the  defining assumption
  \eqref{ConditionsDefiningGOrbiSingularSpaces}
  that $p$ is orbi-singular,
  it is the image under $\mathrm{OrbSnglr}$ \eqref{CohesionSingular}
  of a morphism $p_{\tiny\smooth}$ in $\mathbf{H}_{\tiny\smooth}$:
  \vspace{-2mm}
  \begin{equation}
    \label{CohesiveCoveringSpaceAsSmoothHomotopyFiber}
    \raisebox{13pt}{
    \xymatrix@R=1em{
      X \ar[r]^-{\mathrm{fib}(p)}
      &
      \mathcal{X}
      \ar[d]^-p
      \\
      & \orbisingularG
    }
    }
    \;\;\;\;
    \simeq
    \;\;\;\;
    \mathrm{OrbSnglr}
    \left(\!\!\!\!
    \raisebox{15pt}{
    \xymatrix@R=1em{
      X_{\tiny\smooth} \ar[r]^-{\mathrm{fib}(p_{\tiny\smooth})}
      &
      \mathcal{X}_{\tiny\smooth}
      \ar[d]^-{p_{\tiny \smooth}}
      \\
      &
      \ast \!\sslash\! G
    }
    }
    \!\!\! \right).
  \end{equation}

\vspace{-1mm}
  \noindent
  We claim that in fact $X \simeq \mathrm{NnOrbSinglr}(X_{\tiny\smooth})$,
  whence $X \simeq \smooth(X)$:
  First, since $\mathrm{OrbSnglr}$ is a right adjoint it preserves
  homotopy fibers (Prop. \ref{AdjointsPreserveCoLimits}),
  $\mathrm{fib}(p) \simeq \mathrm{OrbSnglr}
  \big(
    \mathrm{fib}(p_{\tiny\smooth})
  \big)$,
  hence we have $X \simeq \mathrm{OrbSnglr}(X_{\tiny\smooth})$.
  It follows, in particular,
  that $X_{\tiny\smooth}$ is 0-truncated,
  since $X \simeq \mathrm{OrbSnglr}(X_{\tiny\smooth})$ is 0-truncated
  by part {\bf(a)}, and using that $\mathrm{OrbSnglr}$ is fully faithful.
  From this it follows that $\mathrm{OrbSnglr}(X_{\tiny\smooth})
  \simeq \mathrm{NnOrbSinglr}(X_{\tiny\smooth})$,
  by Lemma \ref{Smooth0TruncatedObjectsAreOrbiSingular}.
  Together this gives the claim {\bf (b)}.

  \noindent With this, part {\bf (ii)} now follows by comparison with \eqref{InfinityActionHomotopyFiberSequence}.
\hfill \end{proof}

\medskip
\noindent {\bf Shape of Cohesive $G$-orbispaces.}
We derive the following formula \eqref{PlotsOfShapeOfConcreteGroupoidalGOrbiSpace}
in Prop. \ref{SingularitiesFaithfulSliceOverGSingIsGEquivariantHomotopyTheory}
which generalizes the embedding of $G$-spaces into
global equivariant homotopy theory, discussed in
\cite[p. 7]{Rezk14}\cite[3.2.17]{Lurie19}, from topological $G$-spaces
to general cohesive $G$-spaces. Below in \cref{EquivariantCohomology}
this serves to prove that the intrinsic cohomology of good cohesive
orbispaces subsumes proper equivariant cohomology
(Theorem \ref{OrbifoldCohomologyEquivariant}).

\begin{lemma}[Shape of Cohesive $G$-orbispaces]
  \label{ShapeOfOrbiSingularSpacesAsPresheafOnSingularities}
  Let $\mathbf{H}$ be a singular-cohesive $\infty$-topos
  (Def. \ref{SingularCohesiveInfinityTopos}).
  \eqref{InfinityCategoryOfInfinityGroupoids},
  $G \in \mathrm{Groups}\big(\mathbf{H}\big)$ (Prop. \ref{LoopingAndDelooping})
  be a 0-truncated $G \simeq \tau_0 G$ and discrete
  $G \simeq \flat G$ and let
  $X \in \mathbf{H}$
  be smooth $X \simeq \smooth X$ and 0-truncated $X \simeq \tau_0 X$
  and
  equipped with a $G$-action
  $(X,\rho) \in G\mathrm{Actions}(\mathbf{H})$ (Prop. \ref{InfinityAction}).

 \item {\bf (i)} Then the orbi-singularization \eqref{CohesionSingular}
  of the corresponding homotopy quotient \eqref{InfinityActionHomotopyFiberSequence}
  $$
    \mathcal{X}
    \;:=\;
    \orbisingular
    \big(
      X \!\sslash\! G
    \big)
    \;\in\;
    \mathbf{H}
    \;:=\;
    \mathrm{Sheaves}_\infty
    \big(
      \mathrm{Singularities},
      \,
      \mathbf{H}_{\tiny\smooth}
    \big)
    \,,
  $$
  when regarded as an $\mathbf{H}_{\tiny\smooth}$-valued
  $\infty$-presheaf on $\mathrm{Singularities}$ \eqref{InfinitToposOverSingularities},
  assigns to a singularity $\orbisingularK$
  \eqref{SingularitiesUnderYoneda}
  the disjoint union of
  fixed loci $X^{\phi(K)}$ (Def. \ref{FixedPoints})
  of the smooth covering space $X$ (Def. \ref{CoveringSpaceOfGOrbiSingularSpace})
  for all group homomorphisms $\phi : K \to G$
  homotopy-quotiented \eqref{HomotopyQuotientAsColimit}
  by
  the residual $G$-action (Prop. \ref{PropertiesOfUniversalCoveringSpaces}):
 \begin{equation}
    \label{PlotsOfShapeOfConcreteGroupoidalGOrbiSpace}
  \mathcal{X}
  \;:\;
  \orbisingularK
  \;\longmapsto
    \left(
      \;\;\;\;\;\;
    \raisebox{5pt}{$\underset{{}_{
        \mathclap{
          \phi \in \mathrm{Groups}(K,G)
        }}
      }{
        \bigsqcup
      }
      \;\;\;
      X^{\phi(K)}
     $ }
    \right)
    \sslash G
    \,.
  \end{equation}

\vspace{-2mm}
  \item {\bf (ii)}  Moreover, its shape \eqref{CohesionOfSingularCohesive}
 \vspace{-2mm}
  $$
    \mathrm{Shp}
    \Big(
      \orbisingular
      \big(
        X \!\sslash\! G
      \big)
    \Big)
    \;\in\;
    \mathrm{SingularGroupoids}
    \;:=\;
    \mathrm{Sheaves}_\infty
    \big(
      \mathrm{Singularities}
    \big)
  $$

   \vspace{-2mm}
\noindent
  assigns to a singularity $\orbisingularK$
  \eqref{SingularitiesUnderYoneda}
  the cohesive shape \eqref{AdjunctionCohesion}
  of these disjoint unions of fixed loci (Def. \ref{FixedPoints})
  of the smooth covering space $X$ (Def. \ref{CoveringSpaceOfGOrbiSingularSpace})
  homotopy-quotiented by its $G$-action (Prop. \ref{PropertiesOfUniversalCoveringSpaces}):
 \begin{equation}
    \label{PlotsOfShapeOfConcreteGroupoidalGOrbiSpace2}
  \mathrm{Shp}\big(\mathcal{X}\big)
  \;:\;
  \orbisingularK
  \;\longmapsto
    \mathrm{Shp}
    \left(
      \;\;\;\;\;\;
     \raisebox{5pt}{$\underset{{}_{
        \mathclap{
          \phi \in \mathrm{Groups}(K,G)
        }}
      }{
        \bigsqcup
      }
      \;\;\;\;\;
      X^{\phi(K)}
     $ }
    \right)
    \sslash G
    \,.
  \end{equation}
\end{lemma}
\begin{proof}
We claim that
for $U \in \mathrm{Charts}$ (Def. \ref{ChartsForCohesion})
and $\orbisingularK \in \mathrm{Singularities}$ (Def. \ref{CategoryOfSingularities}),
hence $U \times \orbisingularK \in
\mathrm{Charts} \times \mathrm{Singularities}$ (Lemma \ref{YonedaOnProductSite}),
we have the following sequence of natural equivalences:
\begin{equation}
  \label{ValueOfGOrbiSingularSpacesOnRepresentables}
  \hspace{-1cm}
  \begin{aligned}
    \mathbf{H}
    \big(
      U \times \orbisingularK,
      \,
      \mathcal{X}
    \big)
    & =
    \mathbf{H}
    \big(
      U \times \orbisingularK,
      \,
      \mathrm{ OrbSnglr }
      \big(
        X \!\sslash\! G
      \big)
    \big)
    \\
    & \simeq
    \mathbf{H}_{\tiny\smooth}
    \big(
      \underset{
        \simeq \, U \times (\ast \sslash K)
      }{\;
      \underbrace{
        \mathrm{Smth}
        \big(
          U \times \orbisingularK
        \big)
      }
      }
      \;,
      \,
      X \!\sslash\! G
    \big)
    \\
    & \simeq
    \mathrm{Groupoids}_\infty
    \big(
      (\ast \!\sslash\! K)
      ,
      \,
      \mathbf{H}_{\tiny\smooth}
      \big(
        U,
        X \!\sslash\! G
      \big)
    \big)
    \\
    & \simeq
    \mathrm{Groupoids}_1
    \big(
      (\ast \!\sslash\! K)
      ,
      \,
      \mathbf{H}_{\tiny\smooth}
      \big(
        U, X
      \big)
      \!\sslash\! G
    \big)
    \\
    & \simeq
    \left(
      \;\;\;\;\;\;\;
     \raisebox{4pt}{$ \underset{{}_{
        \mathclap{
          \phi \in \mathrm{Groups}(K,G)
        }}
      }{
        \bigsqcup
      }
      \;\;\;\;\;
      \mathbf{H}_{\tiny\smooth}
      \big(
        U,
        X
      \big)^{ \phi(K) }
    $}
    \right)
    \!\sslash\! G
    \;\;
    \simeq
    \;
    \left(
      \;\;\;\;\;\;\;
    \raisebox{4pt}{$   \underset{{}_{
        \mathclap{
          \phi \in \mathrm{Groups}(K,G)
        }}
      }{
        \bigsqcup
      }
      \;\;\;\;\;
    \mathbf{H}_{\tiny\smooth}
    \big(
      U, X^{ \phi(K) }
    \big)
    $}
    \right)
    \!\sslash\! G
    \\
    & \simeq
    \left(
     \mathbf{H}_{\tiny\smooth}
     \left(
      U,
      \;\;\;\;\;\;\;
   \raisebox{4pt}{$    \underset{{}_{
        \mathclap{
          \phi \in \mathrm{Groups}(K,G)
        }}
      }{
        \bigsqcup
      }
      \;\;\;\;\;
      X^{ \phi(K) }
      $}
    \right)
    \right)
    \!\sslash\! G
    \;\;
    \simeq
    \;
    \mathbf{H}_{\tiny\smooth}
    \left(
      U,
      \left(
    \raisebox{5pt}{$   \underset{{}_{
        {
          \phi \in \mathrm{Groups}(K,G)
        }}
      }{
        \bigsqcup
      }
      \!\!\!\!
      X^{ \phi(K) }
      $}
      \right)
      \!\sslash\! G
    \right).
  \end{aligned}
\end{equation}
Here the first step is \eqref{GOrbiSingularSpaceAsOrbSnglrOfHomotopyQuotient},
the
second is the adjunction $\mathrm{Smth} \dashv \mathrm{OrbSnglr}$
\eqref{CohesionSingular} and using under the brace that
$\mathrm{Smth}$ preserves products (by Prop. \ref{AdjointsPreserveCoLimits}), that $U$ is already smooth by
assumption, and that
$\mathrm{Smth}\big( \orbisingularK \big)
\simeq (\ast \!\sslash\! K )$ by \eqref{ImagedAndPreimagesofBasicOrbiSngularities}.
The third step is the tensoring of $\mathbf{H}$
over $\infty$-groupoids (Prop. \ref{TensoringOfInfinityToposesOverInfinityGroupoids})
(using the geometric discreteness $(\ast \!\sslash\! K)
\simeq \mathrm{Disc}(\ast \!\sslash\! K)$ by Remark \ref{DiscreteObjectsInSingularCohesion})
The fourth step uses the geometric contractibility of $U$
and the discreteness of $G$ to identify
$\mathbf{H}_{\tiny\smooth}( U, X \!\sslash\! G ) \simeq \mathbf{H}_{\tiny\smooth}( U, X )\!\sslash\!G$
(Lemma \ref{HommingChartsIntoHomotopyFiberSequences}).
The fifth is the general observation of Example \ref{HomGroupoidFromBGIntoActionGroupoids}
about hom-groupoids between quotient groupoids of sets.
The sixth step uses Prop. \ref{FixedLocusIn0TruncatedObjectsForDiscreteGroups}
to find that the fixed points in the set of maps are
the maps into the fixed point locus.
After this key step, we
just re-organize term:
The seventh step uses the connectedness of $U$ (Lemma \ref{ChartsAreCohesivelyConnected})
to find that a coproduct of homs out of $U$ is a hom into
the coproduct. Finally, the eighth step
uses again Lemma \ref{HommingChartsIntoHomotopyFiberSequences}.
\begin{enumerate}[{\bf (i)}]
 \vspace{-3mm}
\item The composite equivalence
\eqref{ValueOfGOrbiSingularSpacesOnRepresentables} implies
the first claim \eqref{PlotsOfShapeOfConcreteGroupoidalGOrbiSpace} by
the $\infty$-Yoneda embedding (Prop. \ref{InfinityYonedaEmbedding}),
using Lemma \ref{YonedaOnProductSite}.
 \vspace{-3mm}
\item From this, the second claim
\eqref{PlotsOfShapeOfConcreteGroupoidalGOrbiSpace2}
follows, using that $\mathrm{Shp}$ acts objectwise over
$\mathrm{Singularities}$ \eqref{CohesionOfSingularCohesive},
and preserves homotopy quotients by discrete groups
(Prop. \ref{ShapePreservesHomotopyColimitsByDiscreteGroups}).
\end{enumerate}
\vspace{-7mm}
\hfill \end{proof}

\begin{remark}[Relevance of 0-truncated orbi-singular spaces]
  \label{RelevanceOfConcrete0TruncatedOrbiSingularSpacesAndDiscreteIsotropyGroups}
 $\phantom{A}$
  \vspace{-1mm}
 \item {\bf (i)}  The crucial assumption that makes the proof of Lemma
  \ref{ShapeOfOrbiSingularSpacesAsPresheafOnSingularities}
 work is, {\bf (a)} that  $G$ is discrete
 and {\bf (b)}
 that $X$ is 0-truncated.
    This is what yields 1-groupoidal homs
    in the middle step of \eqref{ValueOfGOrbiSingularSpacesOnRepresentables}
    and thus the form of the expression
    in the next step, as on the right hand side of
  \eqref{GroupoidHomsBetweenQuotientGroupoidsOfSets}.

   \vspace{-1mm}
 \item {\bf (ii)}   Without the assumption of $\mathcal{X}$ being
  0-truncated over $\ast \!\sslash\! G$, the proof of
  Lemma \ref{ShapeOfOrbiSingularSpacesAsPresheafOnSingularities}
  would proceed verbatim up to that middle step, but then
  would break
  as the nontrivial morphisms present in $\mathcal{X}$ would then
  mix with those of the action by $G$.

   \vspace{-1mm}
\item {\bf (iii)}    Lemma \ref{ShapeOfOrbiSingularSpacesAsPresheafOnSingularities}
  shows that this subltety is closely related
  to the cohesive nature of the problem: We either have
  a space which is 0-truncated but carries cohesive (i.e. geometric) structure,
  or we turn it into its cohesive shape which is
  un-truncated but geometrically discrete.
\end{remark}

\medskip

\noindent {\bf Singular quotient of Cohesive $G$-orbispaces.}
\begin{prop}[Singular quotient of $G$-orbi-singular space]
  \label{SingularQuotientOfGOrbiSingularSpaces}
  Let $\mathbf{H}$ be a singular-cohesive $\infty$-topos
  (Def. \ref{SingularCohesiveInfinityTopos}),
  $G \in \mathrm{Groups}(\mathbf{H})$ being discrete $G \simeq \flat G$
  and 0-truncated $G \simeq \tau_0 G$.
  For $\mathcal{X}$ be a $G$-orbi-singular space (Def. \ref{GOrbiSpace})
  with universal covering space
   $
    X
    \;\in\;
    \mathbf{H}_{\tiny\smooth,0}
    \hookrightarrow
    \mathbf{H}
  $
  equipped with its induced $G$-action
  (Def. \ref{CoveringSpaceOfGOrbiSingularSpace}, Prop. \ref{PropertiesOfUniversalCoveringSpaces}).
  Then the singularization \eqref{CohesionSingular}
  of $\mathcal{X}$ is the plain $G$-quotient of $X$
     $$
    \mathrm{Snglr}
    \big(
      \mathcal{X}
    \big)
    \;\simeq\;
    X/G
    \;\;\;\in\;
    \mathbf{H}_{\tiny\smooth,0}
    \longhookrightarrow
    \mathbf{H}_{\tiny\smooth}
  $$
  (i.e., the quotient of the $G$-action formed in the
  1-topos $\mathbf{H}_{\tiny\smooth,0}$ of 0-truncated objects).
\end{prop}
\begin{proof}
  For $U \in \mathrm{Charts}$, write
   \vspace{-2mm}
  \begin{equation}
    \label{DecompositionOfGroupoidIntoConnectedComponents}
    \mathbf{H}_{\tiny\smooth}
    \big(
      U, X
    \big)
    \!\sslash\!
    G
    \;\;\simeq\;\;
    \underset{
      c
    }{\bigsqcup}
    \big(
      \ast \!\sslash\! H_c
    \big)
    \;\;\;\;
    \in\;
    \mathrm{Groupoids}_1
  \end{equation}

   \vspace{-2mm}
\noindent
  for the essentially unique decomposition
  of the groupoid on the left into its
  connected components
  \begin{equation}
    \label{ConnectedComponentOfHomotopyQuotientGroupoidIsQuotientSet}
    c
    \;\in\;
    \pi_0
    \Big(
      \mathbf{H}_{\tiny\smooth}
      \big(
        U, X
     \big)
     \!\sslash\!
     G
    \Big)
    \;\simeq\;
      \mathbf{H}_{\tiny\smooth}
      \big(
        U, X
     \big)
     /
     G
     \;\;
     \in
     \mathrm{Set}
  \end{equation}
  each of which is equivalent to the delooping groupoid
  (Example \ref{DeloopingGroupoids}) of its fundamental group
  $$
    H_c
    \;:=\;
    \pi_1
    \Big(
      \mathbf{H}_{\tiny\smooth}
      \big(
        U, X
      \big)
      \!\sslash\!
      G
      \,,\,
      c
    \Big)
    \;\;
    \in
    \mathrm{Groups}
    \,.
  $$
  Now, by Lemma \ref{ShapeOfOrbiSingularSpacesAsPresheafOnSingularities}
  and re-instantiating the last few manipulations in
  \eqref{ValueOfGOrbiSingularSpacesOnRepresentables},
  we have that over each $U \in \mathrm{Charts}$
  the incarnation of
  the $G$-orbi-singular space $\mathcal{X}$ as an $\infty$-presheaf
  on $\mathrm{Singularities}$ is given by:
  \vspace{-1mm}
  \begin{equation}
    \begin{aligned}
    \mathcal{X}(U)
    \;\;:\;
    \orbisingularK
    \;\longmapsto\;
    &
    \phantom{\simeq}\;
    \mathrm{Groupoids}_1
    \Big(
      \ast \!\sslash\! K
      \,,\,
      \mathbf{H}_{\tiny\smooth}
      \big(
        U, X
      \big)
      \!\sslash\!
      G
    \Big)
    \\
    & \simeq
    \mathrm{Groupoids}_1
    \Big(
      \ast \!\sslash\! K
      \,,\,
      \underset{
        c
      }{\bigsqcup}
      \,
      \big(
        \ast \!\sslash\! H_c
      \big)
    \Big)
    \\
    & \simeq
    \underset{
      c
    }{\bigsqcup}
    \;
    \mathrm{Groupoids}_1
    \Big(
      \ast \!\sslash\! K
      \,,\,
        \ast \!\sslash\! H_c
    \Big)
    \\
    & \simeq
    \underset{
      c
    }{\bigsqcup}
    \;
    \mathrm{Singularities}
    \Big(
      \orbisingularK
      \,,\,
      \orbisingularHc
    \Big).
    \end{aligned}
  \end{equation}

    \vspace{-2mm}
\noindent Here the first step is \eqref{DecompositionOfGroupoidIntoConnectedComponents},
the second step uses that the delooping groupoids
$\ast \!\sslash\! K$ are connected and the last step
observes the definition of $\mathrm{Singularities}$
(Def. \ref{CategoryOfSingularities}).
By the $\infty$-Yoneda embedding (Prop. \ref{InfinityYonedaEmbedding})
over the site of $\mathrm{Singularities}$ \eqref{Singularities}
this means that
  \vspace{-1mm}
\begin{equation}
  \label{OrbiSpaceValueOnUAsCoproductOfSingularities}
  \mathcal{X}(U)
  \;\simeq\;
  \underset{c}{\bigsqcup}
  \;
  \orbisingularHc
  \;\;\;
  \in
  \;
  \mathrm{Sheaves}_\infty
  \big(
    \mathrm{Singularities}
  \big)
  \,.
\end{equation}

  \vspace{-2mm}
\noindent With this, we find that
$\mathrm{Snglr}(\mathcal{X}) \in \mathrm{PreSheaves}_\infty(\mathrm{Charts})$
is given by

\newpage
$\phantom{A}$
\vspace{-6mm}
$$
  \begin{aligned}
    \mathrm{Snglr}
    \big(
      \mathcal{X}
    \big)
    \;\;:\;
    U
    \longmapsto
    &
    \phantom{\simeq}\;
    \mathrm{Snglr}
    \big(
      \mathcal{X}
      (U)
    \big)
    \\
    &
    \;\simeq\;
    \mathrm{Snglr}
    \left(
    \raisebox{2pt}{$
      \underset{c}{\bigsqcup}
      \,
      \orbisingularHc
      $}
    \right)
    \\
    &
    \;\simeq\;
    \underset{c}{\bigsqcup}
    \,
    \mathrm{Snglr}
    \left(
      \orbisingularHc
    \right)
    \\
    &
    \;\simeq\;
    \underset{c}{\bigsqcup}
    \,
    \ast
    \\
    &
    \;\simeq\;
    \pi_0
    \Big(
      \mathbf{H}_{\tiny\smooth}
      \big(
        U, X
     \big)
     \!\sslash\!
     G
    \Big)
    \\
    &
    \;\simeq\;
      \mathbf{H}_{\tiny\smooth}
      \big(
        U, X
     \big)
     /
     G\;.
  \end{aligned}
$$

  \vspace{0mm}
\noindent
Here the first line is the object-wise application of
$\mathrm{Snglr}$ (Remark \ref{ObjectwiseApplicationOfSingularityModalities}),
while the next line is \eqref{OrbiSpaceValueOnUAsCoproductOfSingularities}.
From there we use that $\mathrm{Snglr}$, being a left adjoint,
preserves coproducts (Prop. \ref{AdjointsPreserveCoLimits}) and then that it takes the
elementary singularies to points, by Lemma \ref{ImagesAndPreimagesOfOrbisingularities}.
Finally, we identify \eqref{ConnectedComponentOfHomotopyQuotientGroupoidIsQuotientSet}.
But this resulting assignment is just that of
$X/G \in \mathrm{PreSheaves}(\mathrm{Charts})$:
  \vspace{0mm}
$$
  X/G
  \;\;:\;
  U \;\longmapsto\;
  \mathbf{H}(U,X)/G
$$

  \vspace{-2mm}
\noindent and hence the claim follows.
\hfill \end{proof}

\medskip

\noindent {\bf Examples of Cohesive $G$-orbispaces.}
We make explicit two classes of examples of
cohesive $G$-orbispaces (Def. \ref{GOrbiSpace}):
Fr{\'e}chet-smooth orbispaces and topological orbispaces.

\begin{example}[Fr{\'e}chet smooth $G$-orbispaces]
  \label{FrechetSmoothGOrbifolds}
  Consider
  \vspace{-2mm}
  $$
    X
      \;\in\;
    \mbox{Fr{\'e}chetManifolds}
    \xymatrix{
        \;   \ar@{^{(}->}[r]
      &
}
    \mathrm{SmoothGroupoids}_\infty
  $$

  \vspace{-2mm}
\noindent
  a (possibly infinite-dimensional Fr{\'e}chet-)smooth manifold
  regarded as a 0-truncated concrete smooth $\infty$-groupoid
  \eqref{ConcreteSmoothInfinityGroupoidsAreDiffeologicalSpaces}.
  Given a
  $G \in \mathrm{Groups}(\mathbf{H})$
  \eqref{InclusionOfFiniteGroups} being discrete $G \simeq \flat G$,
  a smooth action $\rho$ of $G$ on $X$ is equivalently
  a homotopy fiber sequence in $\mathrm{SmoothGroupoids}_\infty$
  of this form (Prop. \ref{InfinityAction}):
  \vspace{-2mm}
  $$
    \xymatrix@R=1em@C=3em{
      X
      \ar[r]^-{ \mathrm{fib}(\rho) }
      &
      X \!\sslash\! G
      \ar[d]^-{ \rho }
      \\
      &
      \ast \!\sslash\! G
    }
    \,.
  $$

  \vspace{-2mm}
\noindent
  Here the homotopy quotient \eqref{InfinityActionHomotopyFiberSequence}
  \vspace{-2mm}
  $$
    X
      \!\sslash\!
    G
    \;\in\;
    \xymatrix{
      \mathrm{LieGroupoids}
   \;   \ar@{^{(}->}[r]
      &
      \mathrm{SmoothGroupoids}_\infty
    }
  $$
  is the corresponding
  (possibly infinite-dimensional Fr{\'e}chet-)Lie groupoid,
  regarded as a smooth $\infty$-groupoid
  via the embedding \eqref{1TrunConcreteSmoothInfinityGroupoidsAreDiffeologicalSpaces}.
  Its orbi-singularization \eqref{CohesionSingular}
  is a $G$-orbi-singular space, in the sense of Def. \ref{GOrbiSpace},
  in the $\infty$-topos
  $\mathrm{SingularSmoothGroupoids}_\infty$ \eqref{SingularSmoothGroupoids}:
   \vspace{-3mm}
  \begin{equation}
    \label{SmoothOrbifoldAsOrbiSingularGroupoid}
    \raisebox{16pt}{
    \xymatrix@R=.7em{
      \mathcal{X}
      \ar[d]
      \\
      \orbisingularG
    }
    }
    \;\;:=\;\;
    \mathrm{OrbSnglr}
    \left(\!\!\!\!\!
      \raisebox{20pt}{
      \xymatrix@R=1.5em{
        X
        \!\sslash\! G
        \ar[d]
        \\
        \ast \!\sslash\! G
      }
      }
    \!\!\right).
  \end{equation}

  \vspace{-2mm}
\noindent
  This orbi-singular smooth groupoid \eqref{SmoothOrbifoldAsOrbiSingularGroupoid}
  what we suggest is the proper incarnation of the
  quotient orbifold that is presented by the smooth manifold $X$
  with its $G$-action. Notice that (see \hyperlink{FigureG}{\it Figure G}):
  \begin{enumerate}[{\bf (i)}]
  \vspace{-2mm}
    \item
     its {\it purely smooth aspect} is the Lie groupoid
     $$
       \smooth
       \big(
         \mathcal{X}
       \big)
       \;\simeq\;
       X \!\sslash\! G
       \;\in\;
       \mathrm{LieGroupoids}
           \xymatrix{
        \;   \ar@{^{(}->}[r]
      &
}
       \mathrm{SingularSmoothGroupoids}_{\infty}
       \,,
     $$

     \vspace{-2mm}
\noindent     (by Prop. \ref{PropertiesOfUniversalCoveringSpaces})
     which is the incarnation of this
     orbifold, according to
     \cite{MoerdijkPronk97}\cite{PronkScull10}
     \vspace{-2mm}
    \item
    its {\it purely singular aspect} is the diffeological space
    $$
      \singular
      \big(
        \mathcal{X}
      \big)
      \;\simeq\;
      X/G
      \;\in\;
      \mathrm{DiffeologicalSpaces}
          \xymatrix{
        \;   \ar@{^{(}->}[r]
      &
}
      \mathrm{SingularSmoothGroupoid}_\infty
    $$

        \vspace{-1.5mm}
\noindent
    (by Prop. \ref{SingularQuotientOfGOrbiSingularSpaces})
    which is the incarnation of this orbifold,
    according to \cite{IKZ10}.
  \end{enumerate}

  \vspace{-2mm}
\noindent  However, it is only
  the full orbi-singular object
  $\mathcal{X}$  which is structured enough to
  have proper (Bredon-)equivariant cohomology.
  This is the content of Theorem \ref{OrbifoldCohomologyEquivariant} below.
\end{example}

\newpage

\begin{example}[Topological $G$-orbispaces]
  \label{TopologicalGOrbiSingularSpaces}
  For $G$ a finite group,
  let $G \acts X_{\mathrm{top}}$ be a topological $G$-space
  (Def \ref{GSpaces}) with Borel construction
 \vspace{-2mm}
  $$
    \xymatrix@R=1.5em{
      X_{\mathrm{top}}
      \ar[r]
      &
      X \underset{G}{\times} E G
      \ar[d]
      \\
      &
      B G
    }
  $$
  Via its continuous diffeology
  \eqref{AdjunctionBetweenTopologicalAndDiffeologicalSpaces},
  this is equivalently a 0-truncated (and concrete) object in
  $\mathbf{H}_{\tiny\smooth} := \mathrm{SmoothGroupoids}_\infty$
  (Example \ref{SmoothInfinityGroupoids})
  $$
    X
    \;:=\;
    \mathrm{Cdfflg}(X_{\mathrm{top}})
    \;\in\;
    \mathbf{H}_{\tiny\smooth,0}
  $$
  equipped with a smooth $G$-action (Prop. \ref{InfinityAction})
   \vspace{-3mm}
  $$
    \xymatrix@R=1em{
      X \ar[r]
      &
      X \!\sslash\! G
      \ar[d]
      \\
      &
      \ast \!\sslash\! G
      \,.
    }
  $$

   \vspace{-2mm}
\noindent
  The orbi-singularization \eqref{CohesionSingular}
  of the corresponding homotopy quotient
  is a $G$-orbi-singular space (Def. \ref{GOrbiSpace})
  $$
    \raisebox{16pt}{
    \xymatrix@R=1.5em{
      \mathcal{X}
      \ar[d]
      \\
      \orbisingularG
    }
    }
    \;\;:=\;\;
    \mathrm{OrbSnglr}
    \left(\!\!\!\!
      \raisebox{20pt}{
      \xymatrix@R=1.5em{
        \mathrm{Cdfflg}(X_{\mathrm{top}})
        \!\sslash\! G
        \ar[d]
        \\
        \ast \!\sslash\! G
      }
      }
    \!\!\right).
  $$
\end{example}

\begin{prop}[Shape of good orbifolds]
  \label{ShapeOfFrechetSmoothGOrbifold}
  Consider a finite-dimensional smooth $G$-orbifold, as in
  Example \ref{FrechetSmoothGOrbifolds} (a good orbifold, Remark \ref{GoodOrbifolds})
  $$
    \mathcal{X}
    \;:=\;
    \mathrm{OrbSnglr}
    \big(
      X \!\sslash\! G
    \big).
  $$
    Then its cohesive shape \eqref{CohesionSingular}
  $
    \mathrm{Shp}
    \big(
      \mathcal{X}
    \big)
    \;\in\;
    \mathrm{Sheaves}_\infty
    \big(
      \mathrm{Singularities}
    \big)
  $
  is, over any singularity $\orbisingularK$
  \eqref{AnOrbifoldSingularity}, the
  topological shape \eqref{ShapeOfTopologicalSpaces}
  of the $G$-Borel construction on the
  disjoint union of all $K$-fixed subspaces
  $X_{\mathrm{top}}^{\phi(K)} \subset X_{\mathrm{top}}$ \eqref{HFixedPointSpaces}
  in the underlying \eqref{AdjunctionBetweenTopologicalAndDiffeologicalSpaces}
  D-topological $G$-space (Def. \ref{GSpaces}):
  \begin{equation}
    \label{KValueOfShapeOfGOrbiSingularFrechetManifold}
    \mathrm{Shp}
    \big(
      \mathcal{X}
    \big)
    \;\;:\;\;
    \orbisingularK
    \;\longmapsto
    \;\;\;\;\;\;
    \mathrm{Shp}_{\mathrm{Top}}
    \left(
    \left(
    \;\;\;\;\;\;\;
   \raisebox{5pt}{$ \underset{{}_{
      \mathclap{
        \phi \in \mathrm{Groups}(K,G)
      }}
    }{
      \bigsqcup
    }
    \;\;\;\;
    (\mathrm{Dtplg}(X))^{\phi(K)}
    $}
    \right)
    \underset{G}{\times}
    E G
    \right).
  \end{equation}
\end{prop}

\begin{proof}
  With Lemma \ref{ShapeOfOrbiSingularSpacesAsPresheafOnSingularities},
  the task is reduced to showing that, for $\phi(K) \subset G$ any specified
  subgroup, we have an equivalence
  $$
    \mathrm{Shp}
     \big(
       X^{\phi(K)}
     \big)
     \;\simeq\;
     \mathrm{Shp}_{\mathrm{Top}}
     \big(
       (\mathrm{Dtplg}(X))^{\phi(K)}
     \big)
     \;\;\;\;
     \in \
     \mathrm{Groupoids}_\infty
  $$
  between the cohesive shape \eqref{AdjunctionCohesion}
  of the orbi-singular homotopy quotient of $X$ by $G$
  and the ordinary topological shape \eqref{ShapeOfTopologicalSpaces}
  of the D-topological space underlying $X$.
  But this is \eqref{CohesiveShapeOfSmoothManifolds}
  in Example \ref{SmoothInfinityGroupoids}, given by \cite[4.3.29]{dcct}.
\hfill \end{proof}

\begin{prop}[Shape of topological $G$-orbi spaces]
  \label{ShapeOfGOrbiSingularTopologicalSpaces}
  Consider the topological $G$-orbi-singular space, as in
  Example \ref{TopologicalGOrbiSingularSpaces},
  $$
    \mathcal{X}
    \;:=\;
    \mathrm{OrbSnglr}
    \big(
      \mathrm{Cdfflg}(X_{\mathrm{top}}) \!\sslash\! G
    \big).
  $$
    Then its cohesive shape \eqref{CohesionSingular}
  $
    \mathrm{Shp}
    \big(
      \mathcal{X}
    \big)
    \;\in\;
    \mathrm{Sheaves}_\infty
    \big(
      \mathrm{Singularities}
    \big)
  $
  is, over any singularity $\orbisingularK$
  \eqref{AnOrbifoldSingularity}, the
  topological space \eqref{ShapeOfTopologicalSpaces}
  of the $G$-Borel construction on the
  disjoint union of all $K$-fixed subspaces
  $X_{\mathrm{top}}^{\phi(K)} \subset X_{\mathrm{top}}$ \eqref{HFixedPointSpaces}:
  \vspace{-1mm}
  \begin{equation}
    \label{KValueOfShapeOfGOrbiSingularTopologicalSpace}
    \mathrm{Shp}
    \big(
      \mathcal{X}
    \big)
    \;\;:\;\;
    \orbisingularK
    \;\longmapsto
    \;\;\;\;\;\;\;
    \mathrm{Shp}_{\mathrm{Top}}
    \left(
    \left(
    \;\;\;\;\;\;\;
   \raisebox{4pt}{$ \underset{{}_{
      \mathclap{
        \phi \in \mathrm{Groups}(K,G)
      }}
    }{
      \bigsqcup
    }
    \;\;\;\;\;
    X_{\mathrm{top}}^{\phi(K)}
  $}
    \right)
    \underset{G}{\times}
    E G
    \right).
  \end{equation}
\end{prop}

\vspace{0mm}
\begin{proof}
  With Lemma \ref{ShapeOfOrbiSingularSpacesAsPresheafOnSingularities},
  the task is reduced to showing that, for $\phi(K) \subset G$ any specified
  subgroup, we have an equivalence
  $$
    \mathrm{Shp}
     \big(
       \mathrm{Cdfflg}(X_{\mathrm{top}})^{\phi(K)}
     \big)
     \;\simeq\;
     \mathrm{Shp}_{\mathrm{Top}}
     \big(
       X_{\mathrm{top}}^{\phi(K)}
     \big)
     \;\;\;\;
     \in \
     \mathrm{Groupoids}_\infty
  $$
  between the cohesive shape \eqref{AdjunctionCohesion}
  of the orbi-singular homotopy quotient by $G$
  of the continuous-diffeological space
  and the ordinary topological shape \eqref{ShapeOfTopologicalSpaces}
  But this is item \eqref{CohesiveShapeOfTopologicalSpaces}
  in Example \ref{SmoothInfinityGroupoids}, given by
  combining the result \eqref{SmoothShapeViaPathInfinityGroupoid} of \cite{Pavlov19} with Prop. \ref{DiffeologicalSingularSimpliciaSetOfContinuousDiffeology}
  from \cite{ChristensenWu14}.
\hfill \end{proof}

\medskip

\medskip

\subsection{Orbifolds}
\label{VFolds}

We introduce a general theory of orbi-singular spaces, whose
underlying smooth cohesive groupoid is
locally diffeomorphic to a fixed local model space $V$.
Since, for $V = \mathbb{R}^n \in \mathrm{JetsOfSmoothGroupoids}_\infty$,
these are ordinary $n$-folds
(i.e., ordinary $n$-dimensional manifolds for any $n$, see Example \ref{OrdinaryManifolds}),
or, more generally, {\'e}tale $\infty$-groupoids
with atlases by $n$-folds (Example \ref{EtaleLieGroupoidAsRnFold}),
including ordinary orbifolds,
we generally speak of \emph{$V$-folds}, with a hat tip to \cite{Satake56}.
Externally these are $V$-{\'e}tale $\infty$-stacks (Remark \ref{VFoldsAndVEtaleGroupoids})
but their theory internal to the ambient elastic $\infty$-topos
(such as the construction of their frame bundles in Prop. \ref{TangentBundleOfVFoldIsFiverBundle})
is elegant and finitary and lends itself to full formalization
in homotopy type theory \cite{Wellen18} (see p. \pageref{HoTTFormalizations}).
The \emph{proper} incarnation (see Remark \ref{Proper})
of these $V$-folds as orbifolds is via their orbi-singularization
(Def. \ref{OrbiVFolds}, Remark \ref{Proper}).

\medskip

\noindent {\bf $V$-folds and $V$-{\'e}tale groupoids.}
\begin{defn}[$V$-folds]
  \label{VManifold}
  Let $\mathbf{H}$ be an elastic $\infty$-topos $\mathbf{H}$ (Def. \ref{ElasticInfinityTopos}).
  \vspace{-1mm}
  \item {\bf (i)}  Given
  $V \in \mathrm{Groups}(\mathbf{H})$  (Prop. \ref{LoopingAndDelooping}),
  we say that an object
  $X \in \mathbf{H}$ is a \emph{$V$-fold}
  if there exists a correspondence between $V$ and $X$
    \vspace{-2mm}
  \begin{equation}
    \label{VAtlas}
    \xymatrix@R=.5em@C=3em{
      &
      U
      \ar[dl]_-{\mbox{\tiny \'et}}
      \ar@{->>}[dr]^-{\mbox{\tiny \'et}}
      \\
      V && X
    }
  \end{equation}

   \vspace{-2mm}
\noindent
  such that
  \begin{itemize}
    \vspace{-2mm}
\item[{\bf (a)}] both morphisms are local diffeomorphisms (Def. \ref{FormallyEtaleMorphism})
  and
    \vspace{-3mm}
  \item[{\bf (b)}] the right one is, in addition, an effective epimorphism
  (Def. \ref{EffectiveEpimorphisms}),
  then called a {\it $V$-atlas of $X$} \eqref{StacksAtlasesAndGroupoids}.
\end{itemize}
 \vspace{-3mm}
 \item {\bf (ii)}
 We write
  \begin{equation}
    \label{VFoldsCategory}
    V \mathrm{Folds}(\mathbf{H})
    \;\subset\;
    \mathbf{H}
  \end{equation}
  for the full sub-$\infty$-category
  of $V$-folds in $\mathbf{H}$
  and we write
  \vspace{-2mm}
  \begin{equation}
    \label{VFoldsAndLocalDiffeomorphismsCategory}
    V \mathrm{Folds}(\mathbf{H})^{\mbox{\tiny{\'et}}}
    \;\subset\;
    \mathbf{H}
  \end{equation}

  \vspace{-2mm}
\noindent
  for its wide subcategory on those morphisms
  which are local diffeomorphisms (Def. \ref{FormallyEtaleMorphism}).
\end{defn}

\begin{remark}[$V$-folds and $V$-{\'e}tale groupoids]
  \label{VFoldsAndVEtaleGroupoids}
  By Prop. \ref{EtaleGroupoidsAndEtaleAtlases},
  a $V$-fold (Def. \ref{VManifold})
  is a stack \eqref{StacksAtlasesAndGroupoids}
  whose choice of $V$-atlas \eqref{VAtlas} realizes
  it as an {\'e}tale groupoid (Def. \ref{EtaleGroupoids})
  with space of objects locally diffeomorphic over $V$:
 \vspace{-2mm}
  \begin{equation}
  \label{VFoldsVAtlasesAndGroupoids}
  \raisebox{80pt}{
  \xymatrix@C=10pt{
    &
    \ar@<-12pt>@{..>}[d]
    \ar@<-6pt>@{<..}[d]
    \ar@{..>}[d]
    \ar@<+6pt>@{<..}[d]
    \ar@<+12pt>@{..>}[d]
    &
    \ar@<-12pt>@{..>}[d]
    \ar@<-6pt>@{<..}[d]
    \ar@{..>}[d]
    \ar@<+6pt>@{<..}[d]
    \ar@<+12pt>@{..>}[d]
    \\
    &
    U \times_X U
    \ar@{}[r]|-{\simeq}
    \ar@<-6pt>[d]_-{\mathrm{pr}_1}
    \ar@{<-}[d]|-{\Delta}
    \ar@<+6pt>[d]^{\mathrm{pr}_2}
    &
    U_1
    \ar@<-6pt>[d]_-{s}
    \ar@{<-}[d]|-{e}
    \ar@<+6pt>[d]^-{t}
    &
    \mbox{
      \footnotesize
      \color{darkblue}
      \bf
      ``$V$-{\'e}tale groupoid''
    }
    \\
    V
    \ar@{<-}[r]^-{ \mbox{\tiny{\'e}t} }
    &
    U
    \ar@{->>}[d]_-{a}^-{ \mbox{\tiny{\'e}t} }
    \ar@{=}[r]
    &
    U_0
    \ar@{->>}[d]
    &
    \ar@{}[d]|-{
      \mbox{
        \footnotesize
        \color{darkblue}
        \bf
        ``$V$-atlas''
      }
    }
    \\
    &
    X
    \ar@{}[r]|-{\simeq}
    &
    \underset{\longrightarrow}{\mathrm{lim}}\, U_\bullet
    &
    \mbox{
      \footnotesize
      \color{darkblue}
      \bf
      ``$V$-fold''
    }
  }
  }
  \end{equation}
\end{remark}

\begin{example}[$V$ is a $V$-fold]
  \label{VAsAVFold}
  Let $\mathbf{H}$ be an elastic $\infty$-topos $\mathbf{H}$ (Def. \ref{ElasticInfinityTopos}) and
  $V \in \mathrm{Groups}(\mathbf{H})$  (Prop. \ref{LoopingAndDelooping}).
  Then the underlying object $V \in \mathbf{H}$
  itself is a $V$-fold (Def. \ref{VManifold}):
  A $V$-atlas \eqref{VAtlas} is given by the identity morphisms
  \vspace{-2mm}
  \begin{equation}
    \xymatrix@R=.5em@C=3em{
      &
      V
      \ar[dl]_-{\mathrm{id}} ^-{\mbox{\tiny \'et}}
      \ar@{->>}[dr]^-{ \mathrm{id} }_-{\mbox{\tiny \'et}}
      \\
      V && V \;.
    }
  \end{equation}
\end{example}
\begin{example}[Smooth manifolds are $\mathbb{R}^n$-folds]
  \label{OrdinaryManifolds}
  For $k \in \mathbb{N}$ with $k \geq 1$,
  let $\mathbf{H} = k\mathrm{JetsOfSmoothGroupoids}_\infty$
  (Example \ref{FormalSmoothInfinityGroupoids}).
  Then, for every $n \in \mathbb{N}$, the object
  \vspace{-1mm}
  \begin{equation}
    \label{RnasV}
    V \;:=\;
    \mathbb{R}^n
    \;\in\;
    \mathrm{CartesianSpaces}
    \xymatrix{
        \;   \ar@{^{(}->}[r]
      &
}
    k\mathrm{JetsOfSmoothGroupoids}_\infty
  \end{equation}

  \vspace{-1mm}
\noindent
  canonically carries the structure of a group object
  $(\mathbb{R}^n,+) \in \mathrm{Groups}(\mathbf{H})$,
  via addition in $\mathbb{R}^n$ regarded as a vector space.
  Now every smooth manifold
  \vspace{-1mm}
  $$
    X
    \;\in\;
    \mathrm{SmoothManifolds}
    \xymatrix{
        \;   \ar@{^{(}->}[r]
      &
}
    k\mathrm{JetsOfSmoothGroupoids}_\infty
  $$

  \vspace{-1mm}
\noindent
  of dimension $n$
  is a $V$-fold, hence an $\mathbb{R}^n$-fold
  in the sense of Def. \ref{VManifold}:
  For any choice of atlas in the traditional  sense of
  manifold theory, namely an open cover
   \vspace{-3mm}
  $$
    \big\{ \xymatrix{ U_j \ar[r]^-{\phi_i} & X }  \big\}_{j \in J}
  $$

   \vspace{-2mm}
\noindent
  by local diffeomorphisms $\phi_i$ from open subsets of Cartesian space
 \vspace{-2mm}
  $$
    \xymatrix{
      U_j
    \;  \ar@{^{(}->}[r]^{ \iota_j }
      &
      \mathbb{R}^n
    }
    \,,
  $$

 \vspace{-3mm}
\noindent
  a $V$-atlas \eqref{VAtlas} is obtained by setting:
   \vspace{-3mm}
  \begin{equation}
    \label{TraditionalAtlas}
    \xymatrix@R=.5em@C=3em{
      &
      \underset{j \in J}{\sqcup} U_j
      \ar[dl]^-{\mbox{\tiny \'et}}_-{ (\iota_j)_{j \in J} }
      \ar@{->>}[dr]_-{\mbox{\tiny \'et}}^-{ (\phi_j)_{j \in J} }
      \\
      \mathbb{R}^n && X
    }
  \end{equation}
  \end{example}

\begin{example}[Differentiable {\'e}tale stacks are $\mathbb{R}^n$-folds {\cite[Prop. 4.5.56]{dcct}}]
  \label{EtaleLieGroupoidAsRnFold}
  Let $\mathbf{H} = \mathrm{JetsOfSmoothGroupoids}_\infty$
  (Example \ref{FormalSmoothInfinityGroupoids})
  and take $V = (\mathbb{R}^n, +)$ as in Example \ref{OrdinaryManifolds}.
  Then a diffeological groupoid $X \in \mathbf{H}$
  \eqref{ConcreteSmoothFormalInfinityGroupoidsAreDiffeologicalSpaces}
  is a $V$-fold (Def. \ref{VManifold})
  for $V = \mathbb{R}^n$ \eqref{RnasV} if it is an $n$-dimensional
  {\it differentiable {\'e}tale stack} in that:
  \vspace{-.2cm}
  \begin{itemize}
    \item [\bf (i)] it admits an atlas (effective epimorphism)
   $\xymatrix{X_0 \ar@{->>}[r] & X}$ from a smooth $n$-manifold
   $X_0$ (via \eqref{ConcreteSmoothInfinityGroupoidsAreDiffeologicalSpaces}
   and \eqref{SmoothGroupoidsAsReduced})
  \vspace{-.2cm}
   \item [\bf (ii)] its source and target morphisms
   with respect to this atlas are local diffeomorphisms.
   \end{itemize}

   Generally, a smooth $\infty$-groupoid presented by
   a Kan simplicial smooth manifold is an $\mathbb{R}^n$-fold
   in the sense of Def. \ref{VManifold} if
   it presents an {\it {\'e}tale $\infty$-groupoid}
   in that all its simplicial face maps are local diffeomorphisms.
\end{example}

\begin{examples}[Super-manifolds are $\mathbb{R}^{n\vert q}$-folds]
  \label{SuperfoldsAsVfolds}
  Let $\mathbf{H} = \infty\mathrm{JetsOfSupergeometricGroupoids}_\infty$
  (Example \ref{JetsOfSupergeometricGroupoids}).
  Then, for every $n, q \in \mathbb{N}$, the super-Cartesian space
  (Def. \ref{SuperFormalCartesianSpaces})
  \vspace{-1mm}
  \begin{equation}
    \label{RnasV}
    V \;:=\;
    \mathbb{R}^{n \vert q}
    \;\in\;
    \infty\mathrm{JetsOfSuperCartesianSpaces}
    \xymatrix{
        \;   \ar@{^{(}->}[r]
      &
}
    \infty\mathrm{JetsOfSupergeometricGroupoids}_\infty
  \end{equation}

  \vspace{-1mm}
\noindent
  carries the structure of a group object, whose
  bosonic aspect \eqref{SolidityAdjunctions} is \eqref{RnasV}.
  The corresponding $V$-folds (Def. \ref{VManifold})
  are the $(n\vert q)$-dimensional supermanifolds
  \eqref{EmbeddingOfSupermanifoldsInSuperGroupoids}.

\end{examples}
\begin{example}[General super {\'e}tale $\infty$-stacks]
  \label{GeneralSuperEtaleStacks}
  Let $\mathbf{H} = \infty\mathrm{JetsOfSupergeometricGroupoids}_\infty$
  (Example \ref{JetsOfSupergeometricGroupoids}).
  Then for any $V \in \mathrm{Groups}(\mathbf{H})$ the
  corresponding $V$-{\'e}tale $\infty$-stacks (Remark \ref{VFoldsAndVEtaleGroupoids})
  realize a flavor of \emph{super {\'e}tale $\infty$-stacks},
  locally modeled on $V$.
  Lemma \ref{InSuperGeometricGroupoidsImModalObjectsAreBosonic}
  implies that, generally, the bosonic part $\overset{\rightsquigarrow}{X}$
  of a super {\'e}tale $\infty$-stack is a bosonic {\'e}tale $\infty$-stack
  locally modeled on the bosonic part $\overset{\rightsquigarrow}{V}$
  of $V$:

  \vspace{-.4cm}

  $$
    \xymatrix@R=-5pt{
      V \mathrm{Folds}(\mathbf{H})
      \ar[rr]^-{ \rightsquigarrow }
      &&
      \overset{\rightsquigarrow}{V}
      \mathrm{Folds}(\mathbf{H})
      \\
      \mathllap{
        \mbox{
          \tiny
          \color{darkblue}
          \bf
          \begin{tabular}{c}
            supergeometric
            \\
            {\'e}tale $\infty$-stack
          \end{tabular}
        }
        \;\;
      }
      X
        \ar@{}[rr]|-{ \longmapsto }
        &&
        \overset{\rightsquigarrow}{X}
        \mathrlap{
          \;\;
          \mbox{
            \tiny
            \color{darkblue}
            \bf
            \begin{tabular}{c}
              underlying bosonic
              \\
              {\'e}tale $\infty$-stack
            \end{tabular}
          }
        }
    }
  $$
\end{example}

\medskip

\noindent {\bf Quotients of $V$-folds.}
\begin{prop}[Orbifolding of a $V$-fold is a $V$-fold]
  \label{OrbifoldingOfVFoldIsAVFold}
  Let $\mathbf{H}$ be an elastic $\infty$-topos $\mathbf{H}$ (Def. \ref{ElasticInfinityTopos}),
  $V, G \in \mathrm{Groups}(\mathbf{H})$  (Prop. \ref{LoopingAndDelooping})
  with $G \simeq \flat G $ discrete,
  and $(X,\rho) \in G\mathrm{Actions}(\mathbf{H})$
  (Prop. \ref{InfinityAction}).
  Then if $X$ is a $V$-fold (Def. \ref{VManifold})
  so is its homotopy quotient $X \!\sslash\! G$ \eqref{HomotopyQuotientAsColimit}.
  Specifically,
  if $\xymatrix{ U \ar@{->>}[r]^-{\mbox{\tiny{\'e}t}} & X}$
  is a $V$-atlas for $X$ \eqref{VAtlas},
  then a $V$-atlas for $V \!\sslash\! G$ is given by composition with
  the homotopy fiber inclusion map $\mathrm{fib}(\rho)$ \eqref{InfinityActionHomotopyFiberSequence}:
   \vspace{-2mm}
  \begin{equation}
    \label{VAtlasForVModG}
    \xymatrix@R=.5em@C=3em{
      &
      U
      \ar[dl]_-{\mbox{\tiny \'et}}
      \ar@{->>}[dr]^-{\mbox{\tiny \'et}}
      \\
      V
        &&
      X
      \ar[r]^-{ \mathrm{fib}(\rho) }
      &
      X
      \!\sslash\!
      G \;.
    }
  \end{equation}
\end{prop}
\begin{proof}
  We need to show that the composite morphism on the right of
  \eqref{VAtlasForVModG} is {\bf (a)} an effective epimorphism
  and {\bf (b)} a local diffeomorphism. Since
  both of these classes of morphisms are closed under composition
  (Lemma \ref{EffectiveEpimorphismsArePreservedByPullback} and
  Lemma \ref{ClosureOfLocalDiffeomorphisms}),
  it is sufficient to show that $\mathrm{fib}(\rho)$ itself
  has these two properties.

  \newpage

  For {\bf (a)} observe that, by definition of homotopy fibers
  \eqref{InfinityActionHomotopyFiberSequence},
  we have a Cartesian square
  \vspace{-2mm}
  \begin{equation}
    \label{FirstSquareForVModGIsVFold}
    \raisebox{20pt}{
    \xymatrix@C=3em@R=1.5em{
      X
      \ar[d]
      \ar[rr]^-{ \mathrm{fib}(\rho) }
      \ar@{}[drr]|-{
        \mbox{
          \tiny (pb)
        }
      }
      &&
      X \!\sslash\! G
      \ar[d]^-{ \mathrm{\rho} }
      \\
      \ast
      \ar@{->>}[rr]
      &&
      \mathbf{B}G
    }
    }
  \end{equation}

   \vspace{-2mm}
\noindent
  Here the bottom morphism is an effective epimorphism
  (Example \ref{PointInDeloopingIsEffectiveEpi}).
  Since these are preserved by homotopy pullback, also
  $\mathrm{fib}(\rho)$ is an effective epimorphism.

  For {\bf (b)} consider the image of this square \eqref{FirstSquareForVModGIsVFold}
  under $\Im$. Since $\Im$ is both a right and a left adjoint it preserves
  Cartesian squares and homotopy quotients
  (by Prop. \ref{AdjointsPreserveCoLimits}),
  while it preserves discrete objects by elasticity
  \eqref{ElasticityAdjunctions}
  and idempotency
  (Prop. \ref{CharacterizationOfFullyFaithfulAdjoints},
  Prop. \ref{IdempotentMonads}).
  Therefore
   \vspace{-2mm}
  \begin{equation}
    \label{SecondSquareForVModGIsVFold}
    \raisebox{20pt}{
    \xymatrix@C=3em@R=2em{
      \Im X
      \ar[d]
      \ar[rr]^-{
        \Im \mathrm{fib}(\rho)
        \;\simeq\;
        \mathrm{fib}
        (\Im \rho)
      }
      \ar@{}[drr]|-{ \mbox{\tiny(pb)} }
      &&
      \big(
        \Im X
      \big)
      \!\sslash\! G
      \ar[d]^-{ \Im \mathrm{\rho} }
      \\
      \ast
      \ar@{->>}[rr]
      &&
      \mathbf{B}G
    }
    }
  \end{equation}

   \vspace{-2mm}
\noindent
  is Cartesian.
  Consider finally the pasting composite of this second square
  \eqref{SecondSquareForVModGIsVFold}
  with the naturality square of $\eta^\Im$ on $\mathrm{fib}(\rho)$:
  \vspace{-2mm}
  \begin{equation}
    \label{ThirdSquareForVModGIsVFold}
    \raisebox{40pt}{
    \xymatrix@R=1.8em{
      X
      \ar[d]_-{ \eta^\Im_X }
      \ar[rr]^{ \mathrm{fib}(\rho) }
      \ar@{}[drr]|-{  }
      &&
      X \!\sslash\! G
      \ar[d]_-{ \eta_{X \!\sslash\! G}^\Im }
      \ar@/^3pc/[dd]^-{ \rho }
      \\
      \Im X
      \ar[d]
      \ar[rr]
      \ar@{}[drr]|-{\mbox{\tiny(pb)}}
      &&
      \big(
        \Im X
      \big)
      \!\sslash\!
      G
      \ar[d]_-{ \Im \rho }
      \\
      \ast
      \ar[rr]
      &&
      \mathbf{B}G
    }
    }
  \end{equation}

   \vspace{-1mm}
\noindent
  Here the composite morphism on the right is equivalent to
  $\rho$, as shown, by the naturality of $\eta^\Im$ and
  using that the object $\mathbf{B}G$, being discrete, is $\Im$-modal.
  Therefore, the total outer rectangle of \eqref{ThirdSquareForVModGIsVFold}
  is Cartesian by \eqref{FirstSquareForVModGIsVFold}.
  Moreover, the bottom square of \eqref{ThirdSquareForVModGIsVFold}
  is Cartesian by \eqref{SecondSquareForVModGIsVFold}.
  Therefore the pasting law (Prop. \ref{PastingLaw})
  implies that the top square of \eqref{ThirdSquareForVModGIsVFold}
  is Cartesian. But this means \eqref{PullbackSquareForLocalDiffeomorphisms} that
  $\mathrm{fib}(\rho)$ is a local diffeomorphism.
\hfill \end{proof}

\begin{prop}[Induced $G$-action on the tangent bundle]
  \label{OnTXInducedGAction}
  Let $\mathbf{H}$ be an elastic $\infty$-topos $\mathbf{H}$ (Def. \ref{ElasticInfinityTopos}),
  $V, G \in \mathrm{Groups}(\mathbf{H})$  (Prop. \ref{LoopingAndDelooping}),
  with $G \simeq \flat G $ discrete,
  $(X,\rho) \in G\mathrm{Actions}(\mathbf{H})$
  (Prop. \ref{InfinityAction}) and
  $X \in V \mathrm{Folds}(X)$ (Def. \ref{VManifold}).
  Then the tangent bundle $T X$ (Def. \ref{InfinitesimalTangentBundle})
  carries an essentially unique
  $G$-action $T \rho$ such that:

  {\bf (i)} the defining projection $T X \to X$ is $G$-equivariant (Def. \ref{Equivariance});

  {\bf (ii)} the homotopy quotient of $T X$ is the
   tangent bundle
   of the orbifolded $V$-fold $X \!\sslash\! G$
  (Prop. \ref{OrbifoldingOfVFoldIsAVFold}):
   \vspace{-2mm}
  \begin{equation}
    \label{QuotientTangentBundle}
    (T X)
    \!\sslash\!
    G
    \;\simeq\;
    T(X \!\sslash\! G)
    \;\;\;
    \in
    \mathbf{H}_{\big/ X \!\sslash\! G}
    \,.
  \end{equation}
\end{prop}
\begin{proof}
Consider the following diagram:
\vspace{-2mm}
\begin{equation}
  \label{TangentAction}
  \raisebox{40pt}{
  \xymatrix@C=4em@R=1em{
    T X
    \ar@{-->}[dr]
    \ar[dd]
    \ar[rr]
    &&
    X
    \ar[dd]|-{
      {\phantom{AA}}
    }    \ar[dr]^-{
      \mathrm{fib}(\rho)
    }
    \\
    &
    T (X \!\sslash\! G)
    \ar[rr]
    \ar[dd]
    \ar@/^1.5pc/[dddd]|-{ \phantom{AA} }^>>>>>>>>>>{ T \rho }
    &&
    X \!\sslash\! G
    \ar[dd]^-{\scalebox{.7}{$
      \eta^\Im_X \!\sslash\! G
      \mathrlap{
        \;
        \simeq
        \;
        \eta^\Im_{X \!\sslash\! G}
      }
      $}
          }
    \\
    X
    \ar[dd]
    \ar[dr]_-{ \mathrm{fib}(\rho) }
    \ar[rr]^-{
      \phantom{AA}
    }|>>>>>>>>>>>>>>>>>>>{
      \phantom{AA}
    }^>>>>>>>>>{ \eta^\Im_X }
    &&
    \Im X
    \ar[dr]^-{ \Im \mathrm{fib}(\rho) }
    \\
    &
    X \!\sslash\! G
    \ar[dd]_-{\rho}
    \ar[rr]_-{ \,\eta^\Im_{X \!\sslash\! G} \, }
    &&
    (\Im X)\!\sslash\! G
    \\
    \mathclap{\phantom{A}}
    \;\ast\;
    \ar[dr]
    \\
    &
    \mathbf{B}G
  }
  }
\end{equation}

 \vspace{-2mm}
\noindent
Here the bottom left square is that characterizing the
$G$-action on $X$, by \eqref{InfinityActionHomotopyFiberSequence};
while the bottom and right squares are both the naturality
square of $\eta^\Im$ on the morphism $\mathrm{fib}(\rho)$
(where we use that $\Im$ commutes with taking the
homotopy quotient by the discrete group $G$).  Now observe that:

\vspace{-.3cm}
\begin{itemize}
\item[\bf (a)]
The bottom and right squares are pullback squares since $\mathrm{fib}(\rho)$
is a local diffeomorphism (Def. \ref{FormallyEtaleMorphism})
by Prop. \ref{OrbifoldingOfVFoldIsAVFold}.

 \vspace{-3mm}
\item[\bf (b)] The front and back squares are pullback squares
by the definition of tangent bundles (Def. \ref{InfinitesimalTangentBundle}).
\end{itemize}

 \vspace{-3mm}
\noindent
In particular, the solid part of the diagram is homotopy-commutative,
so that, by the universal property of the front pullback square,
the dashed morphism exists, essentially uniquely,
such as to make the top and the top left square homotopy-commutative.
Further observe, by repeatedly applying
the pasting law (Prop. \ref{PastingLaw}), that:
\begin{itemize}
\vspace{-2mm}
  \item[\bf (c)] The top left square is a homotopy pullback
  since the back, right and front squares
  are pullbacks by (a) and (b).
  \vspace{-3mm}
  \item[\bf (d)] The total left rectangle is a pullback,
  since the top one is so, by (c), and the bottom one is so, by the
  action property \eqref{InfinityActionHomotopyFiberSequence}.
\end{itemize}

 \vspace{-2mm}
\noindent
Thus, again by the action property \eqref{InfinityActionHomotopyFiberSequence},
the total left rectangle exhibits a $G$-action on $T X$ whose homotopy quotient is
as claimed \eqref{QuotientTangentBundle}, and its factorization
into two pullback squares as shown
exhibits the projection $T X \to X$ as a homomorphism
of $G$-actions, hence as being $G$-equivariant (Def. \ref{Equivariance}).
\hfill \end{proof}

\begin{prop}[Induced $G$-action on local neighborhood of fixed point]
  \label{InducedGActionOnLocalNeighbourhoodOfFixedPoint}
  Let $\mathbf{H}$ be an elastic $\infty$-topos $\mathbf{H}$ (Def. \ref{ElasticInfinityTopos}),
  $V, G \in \mathrm{Groups}(\mathbf{H})$  (Prop. \ref{LoopingAndDelooping}),
  with $G \simeq \flat G $ discrete,
  $(X,\rho) \in G\mathrm{Actions}(\mathbf{H})$
  (Prop. \ref{InfinityAction}) with
  $X \in V \mathrm{Folds}(X)$ (Def. \ref{VManifold})
  and $\xymatrix{ \ast \ar[r]^x & X}$
  a homotopy fixed point (Def. \ref{FixedPoints}).
  Then the induced $G$-action $T \rho$ on the tangent bundle $T X$,
  from Prop. \ref{OnTXInducedGAction}, restricts to a $G$-action $T_x \rho$
  on the local neighborhood $T_x X$
  (Example \ref{LocalNeighbourhoodOfAPoint})
      of the homotopy fixed point $x$.
\end{prop}
\begin{proof}
  Consider the following diagram:
   \vspace{-2mm}
  \begin{equation}
  \label{TowardsInducedActionOnLocalNeighbourhoodOfPoint}
  \xymatrix@C=4em@R=1em{
    T_x X
    \ar[rr]
    \ar[dd]
    \ar[dr]
    &&
    T X
    \ar[dr]
    \ar[dd]|-{\phantom{AA}}
    \\
    &
    (T_x X) \!\sslash\! G
    \ar[rr]
    \ar[dd]|<<<<{ T_x \rho }
    &
    &
    (T X) \!\sslash\! G
    \ar[dd]
    \ar@/^1.5pc/[dddd]^-{ T \rho }
    \\
    \ast
    \ar[rr]^<<<<<<<<<<<<{\,x\,}|-{\phantom{AA}}
    \ar[dr]
    &&
    X
    \ar[dd]|-{\phantom{AA}}
    \ar[dr]^-{ \mathrm{fib}(\rho) }
    \\
    &
    \mathbf{B}G
    \ar[rr]_<<<<<<<<<<<<<{\, x \!\sslash\! G  \,}
    \ar@{=}@/_1.3pc/[ddrr]
    &&
    X \!\sslash\! G
    \ar[dd]_-{\rho}
    \\
    &&
    \mathclap{\phantom{A}}
    \;\ast\;
    \ar[dr]
    \\
    &
    &
    &
    \mathbf{B}G
  }
  \end{equation}

   \vspace{-2mm}
\noindent
  Here the squares on the right are from \eqref{TangentAction}
  and are thus both homtopy Cartesian.
  The rear square is the homotopy pullback square defining
  the tangent fiber, and we define the front square to be
  a homotopy pullback, giving us the object denoted
  $(T_x X) \!\sslash\! G$. We need to show
  that this object really is the homotopy quotient of the restricted action.
  But
  the bottom horizontal square homotopy-commutes,
  exhibiting the homotopy fixed point by \eqref{HomotopyFixedPoints},
  so that, by applying the pasting law (Prop. \ref{PastingLaw})
  to the top vertical squares,
  it follows that also the top left square is Cartesian.
  This already identifies $(T_x X) \!\sslash\! G$ as the
  homotopy quotient of some $G$-action on $T_x X$, by
  Prop. \ref{InfinityAction}.
  To see that this is indeed the restricted action,
  observe that the front triangle commutes,
  again by \eqref{HomotopyFixedPoints}, so that the
  total diagram exhibits the fiber inclusion $T_x X \to T X$
  as being a homomorphism $G$-actions
  $T_x \rho \to T_\rho$ (by Prop. \ref{InfinityAction}).
\hfill \end{proof}

\medskip
\noindent {\bf Frame bundles.}
\begin{defn}[Structure group of $V$-folds]
  \label{StructureGroupOfVFolds}
  Let $\mathbf{H}$ be an elastic $\infty$-topos (Def. \ref{ElasticInfinityTopos})
  and
  $V \in \mathrm{Groups}(\mathbf{H})$ (Prop. \ref{LoopingAndDelooping}),
  to be regarded as the local model space of
  $V$-folds (Def. \ref{VManifold}).

\noindent  {\bf (i)} Then we say that the automorphism group (Def. \ref{AutomorphismGroup})
  of the local neighborhood (Example \ref{LocalNeighbourhoodOfAPoint})
  of the neutral element $\xymatrix{\ast \ar[r]^-{e} & V}$ (Example \ref{NeutralElement})
  \begin{equation}
    \label{StructureGroup}
    \mathrm{Aut}(T_e V)
    \;\;
    \in
    \mathrm{Groups}(T_e V)
  \end{equation}
  is the \emph{structure group of $V$-folds}.

 \noindent  {\bf (ii)} We write

  \vspace{-.5cm}

  $$
    (T_e V, \rho_{\mathrm{Aut}}) \in \mathrm{Aut}(T_e V)\mathrm{Actions}(\mathbf{H})
  $$

  \vspace{-.1cm}

 \noindent  for its canonical action \eqref{AutomorphismAction}.
\end{defn}

\begin{example}[Ordinary general linear group]
  \label{OrdinaryGeneralLinearGroup}
  Let $\mathbf{H} = \mathrm{JetsOfSmoothGroupoids}_\infty$
  (Example \ref{FormalSmoothInfinityGroupoids})
  and let
   \vspace{-1mm}
  $$
    V
      :=
    (\mathbb{R}^n,+)
    \;\in\;
    \xymatrix{
      \mathrm{Groups}(\mathrm{SmoothManifolds})
    \;  \ar@{^{(}->}[r]
      &
      \mathrm{Groups}(\mathbf{H})
    }
  $$

   \vspace{-1mm}
\noindent
  via the full inclusion \eqref{ConcreteSmoothFormalInfinityGroupoidsAreDiffeologicalSpaces},
  with $\mathbb{R}^n$ regarded as a group under addition
  of tuples of real numbers. Then the structure group
  of $\mathbb{R}^n$-folds, according to Def. \ref{StructureGroupOfVFolds},
  is the traditional general linear group, regarded as a Lie group:
   \vspace{-1mm}
  $$
    \mathrm{Aut}(T_0 \mathbb{R}^n)
    \;\simeq\;
    \mathrm{GL}(n)
    \,.
  $$
\end{example}

\begin{prop}[Frame bundle]
  \label{TangentBundleOfVFoldIsFiverBundle}
  Let $\mathbf{H}$ be an elastic $\infty$-topos $\mathbf{H}$ (Def. \ref{ElasticInfinityTopos}),
  $V \in \mathrm{Groups}(\mathbf{H})$ (Prop. \ref{LoopingAndDelooping})
  and
  $X \in \mathbf{H}$ a $V$-fold (Def. \ref{VManifold}).
  Then the tangent bundle of $X$ (Def. \ref{InfinitesimalTangentBundle})
  is a fiber bundle (Def. \ref{FiberBundle})
  with typical fiber
  the local neighborhood $T_e V$ (Def. \ref{InfinitesimalNeighborhood})
  of the neutral element $\ast \overset{e}{\longrightarrow} V$,
  hence is the associated bundle of an
  $\mathrm{Aut}(T_e V)$-principal \eqref{StructureGroup}
  bundle (Prop. \ref{ClassificationPrincipalInfinityBundles}),
  to be called the \emph{frame bundle} of $X$:
   \vspace{-2mm}
  \begin{equation}
    \label{FrameBundleClassifyingMap}
    \raisebox{20pt}{
    \xymatrix@C=3em@R=1.5em{
      \overset{
        \raisebox{5pt}{
          \tiny
          \color{darkblue}
          \bf
          tangent bundle
        }
      }{
        T X
      }
      \ar[rr]
      \ar[d]
      \ar@{}[drr]|-{
        \mbox{\tiny\rm(pb)}
      }
      &&
      \big(
        T_e V
      \big)
      \!\sslash\! \mathrm{Aut}\big( T_e V \big)
      \ar[d]
      \\
      X
      \ar[rr]_-{ \scalebox{.6}{$\vdash\,\mathrm{Frames}(X) $}}
      &&
      \mathbf{B} \mathrm{Aut}(T_e V)
    }
    }
    \phantom{AAAAA}
    \raisebox{20pt}{
    \xymatrix@C=3em@R=.9em{
      \overset{
        \raisebox{4pt}{
          \tiny
          \color{darkblue}
          \bf
          frame bundle
        }
      }{
        \mathrm{Frames}(X)
      }
      \ar[d]
      \ar[rr]
      \ar@{}[drr]|-{ \mbox{\tiny\rm(pb)} }
      &&
      \ast
      \ar[d]
      \\
      X
      \ar[rr]_-{\scalebox{.6}{$ \vdash\mathrm{Frames}(X)$} }
      &&
      \mathbf{B}
      \underset{
        \mathclap{
        \raisebox{-3pt}{
          \tiny
          \color{darkblue} \bf
          \begin{tabular}{c}
            structure
            \\
            group
          \end{tabular}
        }
        }
      }{
        \mathrm{Aut}(T_e V)
      }
    }
    }
  \end{equation}
\end{prop}
\begin{proof}
  By Prop. \ref{PullbackAlongLocalDiffeomorphismsPreservesTangentBundles}
  the tangent bundles over any $V$-atlas \eqref{VAtlas} for $X$
  form two Cartesian squares as follows:
  \begin{equation}
    \label{FrameBundleFromVAtlas}
    \raisebox{30pt}{
    \xymatrix@R=.5em@C=3em{
      &
      T U
      \ar[dl]_-{}
      \ar@{->>}[dr]^-{}
      \ar[dd]
      \ar@{}[dddr]|-{ \mbox{\tiny(pb)} }
      \\
      \mathllap{
        V \times T_e V
        \simeq
        \;
      }
      T V
      \ar@{}[dr]|-{ \mbox{\tiny(pb)} }
      \ar[dd]
        &&
      T X
      \ar[dd]
      \\
      &
      U
      \ar[dl]^-{\mbox{\tiny \'et}}
      \ar@{->>}[dr]_-{\mbox{\tiny \'et}}
      \\
      V && X
    }
    }
  \end{equation}
  Moreover, by Prop. \ref{GroupsAdmitFraming}
  the tangent bundle of $V$ is trivial, as shown on the left.
  Since Cartesian products are preserved by homotopy pullback,
  the left square implies that also $T U \simeq U \times T_e V$
  is trivial. But with this the
  existence of the right square is the defining characterization for
  $T X$ being a $T_e V$-fiber bundle.
\hfill \end{proof}
\begin{remark}[Frame bundles are well-defined]
  The frame bundle (Def. \ref{TangentBundleOfVFoldIsFiverBundle})
  of a $V$-fold (Def. \ref{VManifold}) is
  independent, up to a contractible space of equivalences,
  of the choice of $V$-atlas \eqref{VAtlas}
  in the construction \eqref{FrameBundleFromVAtlas}:
  This follows as a special case of the essential independence
  of classifying maps
  of fiber bundles from the choice of trivializing cover,
  as in Prop. \ref{FiberBundlesClassified}, using that
  not only the class of effective epimorphisms
  but also that of local diffeomorphisms is closed
  under pullback and composition (Lemma \ref{ClosureOfLocalDiffeomorphisms}).
\end{remark}

\begin{prop}[$V$-fold is $\mathrm{Aut}(T_e V)$-quotient of its frame bundle]
  \label{VFoldAsQuotientOfItsFrameBundle}
  Let $\mathbf{H}$ be an elastic $\infty$-topos $\mathbf{H}$ (Def. \ref{ElasticInfinityTopos}),
  $V \in \mathrm{Groups}(\mathbf{H})$ (Prop. \ref{LoopingAndDelooping})
  and $X \in V\mathrm{Folds}(\mathbf{H})$  (Def. \ref{VManifold}).
Then $X$ is equivalent to the homotopy quotient \eqref{HomotopyQuotientAsColimit}
of its
  own frame bundle (Prop. \ref{TangentBundleOfVFoldIsFiverBundle})
  by $\mathrm{Aut}(T_e V)$:
   \vspace{-2mm}
  $$
    X \;\simeq\;
    \mathrm{Frames}(X) \!\sslash\! \mathrm{Aut}(T_e V)
    \,.
  $$
\end{prop}
\begin{proof}
  This is immediate from the equivalence between principal bundles
  and homotopy qotient projections
  (Remark \ref{PrincipalBaseSpacesAreHomotopyQuotients}) applied to
  the frame bundle \eqref{FrameBundleClassifyingMap}.
\hfill \end{proof}


\begin{example}[Frame bundles on smooth manifolds]
  \label{FrameBundleOnSmoothManifold}
  Let $\mathbf{H} = \mathrm{JetsOfSmoothGroupoids}_\infty$
  (Example \ref{FormalSmoothInfinityGroupoids})
  and $X \in \mathrm{SmoothManifolds} \hookrightarrow \mathbf{H}$
  a smooth manifold \eqref{ConcreteSmoothFormalInfinityGroupoidsAreDiffeologicalSpaces}
  regarded as an $\mathbb{R}^n$-fold according to Example \ref{OrdinaryManifolds}.

 \vspace{-2mm}
\item {\bf (i)} Then its frame bundle, according to Prop.
\ref{TangentBundleOfVFoldIsFiverBundle},
  is the $\mathrm{GL}(n)$-principal bundle on $X$ which is the
  frame bundle in the traditional sense of differential geometry.
 \vspace{-2mm}
\item {\bf (ii)}  For the same manifold but regarded in
  $\mathbf{H} = k\mathrm{JetsOfSmoothGroupoids}_\infty$
  with $k \geq 1$ we instead get the corresponding jet version of the
  frame bundle (see e.g.  \cite[12.12]{KMS93}).
\end{example}

\medskip
\medskip

\noindent {\bf Framed $V$-folds.}
\begin{defn}[Framing]
  \label{Framing}
  Let $\mathbf{H}$ be an elastic $\infty$-topos $\mathbf{H}$
  (Def. \ref{ElasticInfinityTopos}).
  A \emph{framing} of an objext $X \in \mathbf{H}$ is
  a trivialization of its tangent bundle Def. \ref{InfinitesimalTangentBundle}, hence
  an equivalence
   \vspace{-2mm}
  $$
    T X \;\simeq\; X \times T_x X
    \;\;
    \in
    \mathbf{H}_{/X}
  $$

   \vspace{-2mm}
\noindent
  for $\xymatrix{ \ast \ar[r]^-{x} & X }$ any point.
\end{defn}

\begin{remark}[Framing on a $V$-fold]
  \label{FramingVFold}
  If $X$ is a $V$-fold (Def. \ref{VManifold})
  then a \emph{framing} on $V$ in the sense of Def. \ref{Framing}
  is equivalent, by Prop. \ref{TangentBundleOfVFoldIsFiverBundle},
  to a trivialization of the frame bundle, hence to a trivialization of its
  classifying map \eqref{FrameBundleClassifyingMap}:
  \vspace{-2mm}
  \begin{equation}
    \label{FramingAsTrivializationOfFrameBundle}
    \raisebox{20pt}{
    \xymatrix{
      T X
      \ar[dr]
      \ar[rr]^-{ \mathrm{fr} }_-{ \simeq }
      &&
      X \times T_e V
      \ar[dl]
      \\
      & X
    }
    }
    \phantom{AAA}
    \Leftrightarrow
    \phantom{AAA}
    \raisebox{20pt}{
    \xymatrix@R=5pt@C=3em{
      &
      \ast
      \ar[dr]_<<<{\ }="s"
      \\
      X
      \ar[ur]
      \ar@/_1.4pc/[rr]_-{ \scalebox{.6}{$\vdash\,\mathrm{Frames}(X)$} }^{\ }="t"
      &
      &
      \mathbf{B} \mathrm{Aut}(T_e V)\;.
      \ar@{<=}_-{
        \simeq
      }^-{\vdash\,\mathrm{fr}} "s"; "t"
    }
    }
  \end{equation}
\end{remark}

\begin{prop}[Groups carry canonical framings by left-translation]
  \label{GroupsAdmitFraming}
  In an elastic $\infty$-topos $\mathbf{H}$
  (Def. \ref{ElasticInfinityTopos})
  every group object
  $V \in \mathrm{Groups}(\mathbf{H})$  (Prop. \ref{LoopingAndDelooping})
  carries a canonical framing
  (Def. \ref{Framing}), which we call
  the \emph{framing by left translation}:
   \vspace{-2mm}
  \begin{equation}
    \label{LeftTranslationFramingOnGroup}
    \xymatrix{
      T V
      \ar[rr]^-{ \mathrm{fr}_\ell }_-{ \simeq }
      &&
      V \times T_e V
    }
    \mathrlap{
    \;\;\;\;\;
    \in
    \mathbf{H}_{/V}
    \,.
    }
  \end{equation}
\end{prop}
\begin{proof}
  Since $\Im$ preserves group structure
  (as in Prop. \ref{ShapePreservesHomotopyColimitsByDiscreteGroups}),
  the defining homotopy fiber product of the tangent bundle of $V$
  \eqref{TangentBundlePullbackDefinition} sits in a
  Mayer-Vietoris sequence (Prop. \ref{MayerVietorisSequence})
  as shown on the left of the following:
   \vspace{-2mm}
  \begin{equation}
    \label{TowardsFramingByLeftTranslation}
   \raisebox{20pt}{
    \xymatrix@C=3em{
      T V
      \ar[rrrr]
      \ar[d]
      \ar@{}[drrrr]|-{ \mbox{\tiny(pb)} }
      && &&
      \ast
      \ar[d]^-{ \vdash e }
      \\
      V \times V
      \ar[rr]_-{
        (\eta^\Im_V, \eta^\Im_V) = \eta^{\Im}_{V \times V}
      }
      &&
      \Im V \times \Im V
      \ar[rr]_-{
        (-) \cdot (-)^{-1}
      }
      &&
      \Im V
    }
    }
    \;\;\;
    \simeq
    \;\;\;
    \raisebox{20pt}{
    \xymatrix{
      T V
      \ar@{}[drr]|-{\mbox{\tiny(pb)}}
      \ar[d]
      \ar[rr]
      &&
      T_e V
      \ar[d]
      \ar[r]
      \ar@{}[dr]|-{\mbox{\tiny(pb)}}
      &
      \ast
      \ar[d]^-{ \vdash e }
      \\
      V \times V
      \ar[rr]_-{ (-)\cdot (-)^{-1} }
      &&
      V
      \ar[r]_-{ \eta^\Im_V }
      &
      \Im V
    }
    }
  \end{equation}

   \vspace{-2mm}
\noindent
  Using that $\Im$ preserves products (by Prop. \ref{AdjointsPreserveCoLimits})
  and using the naturality of
  its unit transformation $\eta^\Im$ \eqref{AdjunctionUnit}, this Cartesian
  square on the left is equivalent to the total rectangle shown on the
  right. By the pasting law (Prop. \ref{PastingLaw}),
  this is the pasting of two Cartesian squares, the
  right one of which exhibits the local neighborhood $T_e V$
  (Def. \ref{InfinitesimalNeighborhood}) as shown.
  To see what the Cartesian property of the
  left square on the right says, consider pasting to it
  the top square appearing in
  the diagram \eqref{GroupOperationFromLooping} which
  exhibits the group division $(-)\cdot (-)^{-1}$ in
  Example \ref{GroupOperation}:
   \vspace{-1mm}
  \begin{equation}
    \label{FurtherTowardsFramingByLeftTranslation}
    \xymatrix@C=3em@R=1.5em{
      T V
      \ar[rr]
      \ar@{}[drr]|-{\mbox{\tiny(pb)}}
      \ar[d]
      &&
      T_e V
      \ar[d]
      \\
      V \times V
      \ar[d]_-{ \mathrm{pr}_1 }
      \ar[rr]|-{\; (-) \cdot (-)^{-1} \; }
      \ar@{}[drr]|-{\mbox{\tiny(pb)}}
      &&
      V
      \ar[d]
      \\
      V \ar[rr]
      &&
      \ast
    }
  \end{equation}

   \vspace{-2mm}
\noindent
  Since both squares are Cartesian, the pasting law
  (Prop. \ref{PastingLaw}) says that the total rectangle is
  Cartesian. This is the equivalence \eqref{LeftTranslationFramingOnGroup}.
\hfill \end{proof}

\begin{prop}[Canonical framing on group is equivariant under group automorphisms]
  \label{CanonicalFramingIsEquivariantUnderGroupAutomorphisms}
  Consider an elastic $\infty$-topos $\mathbf{H}$
  (Def. \ref{ElasticInfinityTopos}),
  $V, G \in \mathrm{Groups}(\mathbf{H})$  (Prop. \ref{LoopingAndDelooping}).
  with 0-truncated  $V \simeq \tau_0 V$
  and
  $(V,\rho_G) \in G\mathrm{Actions}(\mathbf{H})$
  (Prop. \ref{InfinityAction})
  acting by
  group-automorphisms
  (Prop. \ref{CanonicalActionOfGroupAutomomorphisms})
  hence by restriction
  $\rho_G = \mathbf{B}i^\ast \rho_{\mathrm{Aut}_{\mathrm{Grp}}} $ (Prop. \ref{PullbackAction})
  along a group homomorphism
  $
    \xymatrix{
      G
      \ar[r]|-{\,i\,}
      &
      \mathrm{Aut}_{e}(V)
    }
  $,
  to the
  group-automorphism group $\mathrm{Aut}_{\mathrm{Grp}}(V)$
  (Def. \ref{GroupOFGroupAutomorphisms}).
Then  the canonical framing
  $\mathrm{fr}_\ell$ on $V$ from Prop. \ref{GroupsAdmitFraming}
  is $G$-equivariant (Def. \ref{Equivariance}),
  in that it lifts to a morphism of $G$-actions
  (Prop. \ref{InfinityAction}) of the form
  \vspace{-2mm}
  $$
  \hspace{-1cm}
    \xymatrix{
      (T V, T \rho)
      \ar[rr]^-{ \mathrm{fr}_\ell }
      &&
      (V, \rho)
        \times
      (T_e V, T_e \rho)
      \mathrlap{
        \;\;\;\;
        \in
        G\mathrm{Actions}(\mathbf{H})
        \,,
      }
    }
  $$

  \vspace{-0mm}
\noindent
  where $T \rho$ is the induced action on $T V$
  from Prop. \ref{OnTXInducedGAction},
  and $T_e \rho$ is the induced action on $T_e V$
  from Prop. \ref{InducedGActionOnLocalNeighbourhoodOfFixedPoint}
  (which exists since group-automorphisms of $V$ are in particular
  pointed automorphisms of $V$ (Def. \ref{PointedAutomorphismGroup}).
\end{prop}
\begin{proof}
  Consider the following diagram:
  \vspace{-2mm}
  \begin{equation}
    \label{TowardsEquivarianceOfTheCanonicalFraming}
    \raisebox{80pt}{
    \xymatrix@R=10pt{
      T V
      \ar@{=}[dr]
      \ar@{->>}[rr]
      \ar[dddr]
      &&
      (T V) \!\sslash\! G
      \ar@{-->}[dr]^-{\phi_3}
      \ar@{-->}@/^2.4pc/[ddrr]^-{\phi_2}
      \ar[dddr]|<<<<<<<<<<<<{\phantom{AA}}|>>>>>>>>>>>>{\phantom{AA}}
      \\
      &
      T V
      \ar[dr]
      \ar[dd]
      \ar@{-->}[rr]^-{\phi_1}
      &&
      (V \!\sslash\! G)
        \underset{
          \ast \!\sslash\! G
        }{\times}
      \big(
        (T_e V) \!\sslash\! G
      \big)
      \ar[dr]
      \ar[dd]|<<<<<<{\phantom{AA}}
      \\
      &
      &
      T_e V
      \ar@{^{(}->}[dd]
      \ar[rr]
      &&
      (T_e V) \!\sslash\! G
      \ar@{^{(}->}[dd]
      \\
      &
      V \times V
      \ar[dr]|-{ (-)\cdot(-)^{-1} }
      \ar[dd]_-{ \mathrm{pr}_1 }
      \ar[rr]|<<<<<<<<<<<<<<<{\phantom{AA}}
      &&
      (V \!\sslash\! G)
      \underset{
        \ast \!\sslash\! G
      }{\times}
      (V \!\sslash\! G)
      \ar[dr]
      \ar[dd]|<<<<<<{\phantom{AA}}
      \\
      &
      &
      V
      \ar[rr]
      \ar[dd]
      &&
      V \!\sslash\! G
      \ar[dd]
      \\
      &
      V
      \ar[dr]
      \ar[rr]|<<<<<<<<<<<<<<<<<<<<{ \phantom{AA} }|>>>>>>>>>>>>>>{\, \mathrm{fib}(\rho)\, }
      &&
      V \!\sslash\! G
      \ar[dr]
      \\
      &
      &
      \ast
      \ar[rr]
      &&
      \ast \!\sslash\! G
    }
    }
  \end{equation}

  \vspace{-2mm}
\noindent
  Here
  \begin{itemize}
    \vspace{-.25cm}
    \item the bottom square is
      the Cartesian square \eqref{InfinityActionHomotopyFiberSequence}
      which exhibits the action on $V$,
    \vspace{-.25cm}
    \item the middle horizontal square is the
    Cartesian square
    which exhibits the equivariance under group-automorphisms of
    the group division operatoin
    (Prop. \ref{GroupDivisionIsEquivariantUnderGroupAutomorphisms}),
    \vspace{-.25cm}
  \item
    the total left rectangle is
    the Cartesian square from \eqref{FurtherTowardsFramingByLeftTranslation}
    which exhibits the canonical framing,
    \vspace{-.25cm}
  \item
    the total front face is the pasting of
    \begin{itemize}
      \vspace{-.25cm}
      \item on the bottom: the Cartesian square \eqref{InfinityActionHomotopyFiberSequence}
      which exhibits the action on $V$,
      \vspace{-.1cm}
      \item on the top: the Cartesian square
      which is the pasting of the top and the top-right squares in
      \eqref{TowardsInducedActionOnLocalNeighbourhoodOfPoint}
      equibiting the action on $T_e V$
    \end{itemize}
    \vspace{-.2cm}
    and hence is itself Cartesian,
    \vspace{-.25cm}
  \item
    the bottom and the total right squares are the
    defining Cartesian squares of the fiber products,
    and hence, by the pasting law,
    also their pasting to the total right square
    is Cartesian,
    \vspace{-.25cm}
  \item
    the total vertical rear square, with the dashed morphism $\phi_1$ on top,
    is the one thus induced from the
    universal property of the fiber product,
    and is itself Cartesian, by the pasting law (Prop. \ref{PastingLaw}),
    (using, by the above items, that the left, right and front squares are Cartesian
    and that the diagram of squares commutes)
    \vspace{-.25cm}
  \item
    the slanted square in the rear is the pasting of the
    Cartesian square on the left
    of \eqref{TangentAction},
    that exhibits the induced $G$-action on $T V$, with the
    diagonal square on $\mathrm{fib}(\rho)$.
  \end{itemize}

  \vspace{-2mm}
\noindent
  Now observe that inside this big diagram \eqref{TowardsEquivarianceOfTheCanonicalFraming}
  we find the following
  solid homotopy-commutative sub-diagram
  \vspace{-2mm}
  $$
    \xymatrix@C=3em@R=1.5em{
      T V
      \ar@{->>}[d]
      \ar[r]
      &
      (T_e V) \!\sslash\! G
      \ar@{^{(}->}[d]
      \\
      (T V) \!\sslash\! G
      \ar[r]
      \ar@{-->}[ur]^{\phi_2}
      &
      V \!\sslash\! G
      \,.
    }
  $$

  \vspace{-2mm}
\noindent
  Here the left morphism is an effective epimorphism
  (by Lemma \ref{HomotopyQuotientMapsAreEffectiveEpimorphisms})
  and the right morphism is (-1)-truncated by the assumption that
  $V$ is 0-truncated (Lemma xyz).
  Therefore,
  the connected/truncated factorization system (Prop. \ref{nConnectednTruncatedFactorizationSystem})
  implies an essentially unique lift $\phi_2$, as shown.
  This, in turn, implies the morphism $\phi_3$ in
  \eqref{TowardsEquivarianceOfTheCanonicalFraming},
  again by the universal property of the homotopy fiber product.

  Now, since both the slanted as well as the
  vertical total rear squares are Cartesian,
  the diagram \eqref{TowardsEquivarianceOfTheCanonicalFraming}
  shows that the contravariant base change (Prop. \ref{BaseChange})
  of $\phi_3$ along
  $\mathrm{fib}(\rho)$ is an equivalence.
  But since  $\mathrm{fib}(\rho)$ is an effective epimorphism
  (Lemma \ref{BaseChangeAlongEffectiveEpimorphismsIsConservative}) ,
  base change along it is conservative (Prop. \ref{BaseChangeAlongEffectiveEpimorphismsIsConservative}),
  and hence it follows
  that $\phi_3$ itself is already an equivalence.

  With that identification, the total cube in \eqref{TowardsEquivarianceOfTheCanonicalFraming}
  exhibits the $G$-equivariance of the framing.
\hfill \end{proof}

\begin{prop}[Orbifolding of framed $V$-folds]
  \label{FactoringVramesOfQuotientOfFramed}
  Let $\mathbf{H}$ be an elastic $\infty$-topos $\mathbf{H}$ (Def. \ref{ElasticInfinityTopos}),
  $V, G \in \mathrm{Groups}(\mathbf{H})$  (Prop. \ref{LoopingAndDelooping})
  with $G \simeq \flat G$ discrete,
  and
  $(X,\rho_X), (T_e V, \rho_{T_e V}) \in G\mathrm{Actions}(\mathbf{H})$
  (Prop. \ref{InfinityAction})
  for $X$ a $V$-fold (Def. \ref{VManifold})
  equipped with a framing $\mathrm{fr}$
  (Def. \ref{Framing}).
 Then the following are equivalent:

  \begin{itemize}
   \vspace{-3mm}
   \item[\bf (i)] The framing is $G$-equivariant (Def. \ref{Equivariance})
    with respect to the induced action on $T X$ (from Prop. \ref{OnTXInducedGAction}) and the product action
    $\rho_X \times \rho_{T_e V}$ on
    $X \times T_e V$, hence lifts to a morphism
     \vspace{-2mm}
    \begin{equation}
      \label{FramingGEquivariant}
      \xymatrix{
        (T X, \rho_{T X})
        \ar[rr]^-{ \mathrm{fr} }_-{ \simeq }
        &&
        (X,\rho_X) \times (T_e V, \rho_{T_e V})
      }
      \mathrlap{
        \;\;\;\;
        \in
        G\mathrm{Actions}(\mathbf{H})
      }
    \end{equation}
 \vspace{-7mm}
  \item[\bf (ii)] The classifying map \eqref{FrameBundleClassifyingMap}
  of the frame bundle (Def. \ref{TangentBundleOfVFoldIsFiverBundle})
  of the orbifolded $V$-fold $X \!\sslash\! G$ (Prop. \ref{OrbifoldingOfVFoldIsAVFold})
  factors through $\mathbf{B}G$ as
   \vspace{-1mm}
  \begin{equation}
    \label{FractoringFrameBundleOnQuotient}
    \xymatrix{
      X \!\sslash\! G
      \ar@/^2pc/[rrrr]^{\scalebox{.6}{$
        \vdash\, \mathrm{Frames}(X \!\sslash\! G)
        $}
      }_-{\ }="s"
      \ar[rr]_-{
        \rho_X
      }
      &&
      \mathbf{B}G
      \ar[rr]_-{
        \vdash\,\rho_{T_e V}
      }
      &&
      \mathbf{B}
      \mathrm{Aut}(T_e V)
      \ar@{=>} "s"; "s"+(0,-5)
    }
  \end{equation}
  \end{itemize}
\end{prop}
\begin{proof}
Consider the following diagram:
$$
  \xymatrix@R=10pt@C=12pt{
    T X
    \ar[dddr]
    \ar[dr]^-{\mathrm{fr}}
    \ar[rr]
    &&
    T (X \!\sslash\! G)
    \ar[dddr]|<<<<<<<<<<<<{ \phantom{ {AA} \atop {AA} } }
    \ar[dr]^-{ \mathrm{fr} \!\sslash\! G }
    \ar[rrrr]
    && &&
    (T_e V) \!\sslash\! \mathrm{Aut}(T_e V)
    \ar@{=}[dr]
    \ar[dddr]|<<<<<<<<<<<<<{ \phantom{ {AA} \atop {AA}} }
    \\
    &
    X \times T_e V
    \ar[dd]
    \ar[rr]
    &&
    X \!\sslash\! G
      \underset{\mathbf{B}G}{\times}
    (T_e V) \!\sslash\! G
    \ar[dd]
    \ar[rr]
    &&
    (T_e V) \!\sslash\! G
    \ar[rr]
    \ar[dd]|-{ \rho_{T_e V} }
    &&
    (T_e V) \!\sslash\! \mathrm{Aut}\big( T_e V \big)
    \ar[dd]
    \\
    \\
    &
    X
    \ar[rr]|-{ \, \mathrm{fib}(\rho_X) \,}
    \ar[drrr]
    \ar@/_3pc/[rrrrrr]_{\scalebox{.6}{$ \vdash\, \mathrm{Frames}(X)$} }^>>>>>>>>>>>>>>>>>>>>>>>>>>>{\ }="s1"
    &&
    X \!\sslash\! G
    \ar[rr]|-{\, \rho_X }
    \ar@{}[rrrr]|-{\phantom{\vert}}="t2"
    \ar@/^1.43pc/[rrrr]|<<<<<<<<<<<<<<<{\phantom{AAA}}^>>>>>>>>>>>{\scalebox{.6}{$\vdash\,\mathrm{Frames}(X \!\sslash\! G)$}}|-{\phantom{\vert}}="s2"
    &&
    \mathbf{B}G
    \ar[rr]_-{
      \vdash
      \,
      \rho_{T_e V}
    }
    &&
    \mathbf{B}
    \mathrm{Aut}
    \big(
      T_e V
    \big)
    \\
    &
    && &
    \ast
    \ar[ur]_<<<<<{\ }="t1"
    \ar@{=>}_-{\, \vdash\,\mathrm{fr}} "s1"; "t1"
    \ar@{=>}^-{\; \vdash\,\mathrm{fr}\!\sslash\! G} "s2"; "t2"
  }
$$

\vspace{-.5cm}
\noindent
Note that here:
\begin{itemize}
 \vspace{-2mm}
\item[{\bf (a)}] The total outer part of the diagram exhibits the
given framing $\mathrm{fr}$ via its classifying homotopy
$\vdash\mathrm{fr}$, according to Remark \ref{FramingVFold}.
 \vspace{-2mm}
\item[{\bf (b)}] The front squares in the middle and on the right
are the pullback squares that defines the diagonal $G$-action
and the classification of the $\rho_{T_e V}$-action respectively.
Hence also their pasting composite is a pullback, by the pasting law
(Prop. \ref{PastingLaw}).
\end{itemize}

 \vspace{-2mm}
\noindent
\noindent {\bf (i)}
First to see that $G$-equivariance of $\mathrm{fr}$
implies the factorization \eqref{FractoringFrameBundleOnQuotient}:
By the characterization of $G$-actions \eqref{InfinityActionHomotopyFiberSequence}
$G$-equivariance of $\mathrm{fr}$ means, equivalently,
that $\mathrm{fr}$ is the morphism
on homotopy fibers over $\mathbf{B}G$ induced from
an equivalence $\mathrm{fr} \!\sslash\! G$ on
homotopy quotients.
But, by {\bf (b)} and Prop. \ref{OnTXInducedGAction},
such an equivalence is classified by a homotopy of the form \eqref{FractoringFrameBundleOnQuotient}.

\vspace{1mm}
\noindent {\bf (ii) } Now to see that, conversely, the existence of a homotopy
`` $\vdash\mathrm{fr}\!\sslash\!G$ ''
of the form
\eqref{FractoringFrameBundleOnQuotient} implies the
existence of a $G$-equivariant framing $\mathrm{fr}$
(quotation marks now since we yet have to show that the two are related
in this way).
For this, we have to show that the morphism on
homotopy fibers induced by $\mathrm{fr}\!\sslash\!G$
is a framing $\mathrm{fr}$.
But, by the nature of the $G$-action on $T X$ from Prop. \ref{OnTXInducedGAction},
the nature of the diagonal $G$-action
exhibited by the middle front square, and using the
pasting law (Prop. \ref{PastingLaw}), this means to show
that the left front and rear squares are homotopy pullbacks.
For the front left square this follows by the factorization of
$\rho_X \circ \mathrm{fib}(\rho_X)$ through the point,
using {\bf (a)}, {\bf (b)} and  the pasting law (Prop. \ref{PastingLaw}).
For the rear left square,
this follows
by Prop. \ref{PullbackAlongLocalDiffeomorphismsPreservesTangentBundles},
since $\mathrm{fib}(\rho)$ is a local diffeomorphism by
Prop. \ref{OrbifoldingOfVFoldIsAVFold}.
\hfill \end{proof}

\medskip

\noindent {\bf $G$-Structures.}
\begin{defn}[$G$-Structure coefficients]
  \label{GStructureCoefficients}
  Let $\mathbf{H}$ be an elastic $\infty$-topos  (Def. \ref{ElasticInfinityTopos})
  and
  $V \in \mathrm{Groups}(\mathbf{H})$ (Prop. \ref{LoopingAndDelooping}).
  Then a \emph{coefficient for $G$-structure}
  \vspace{-2mm}
  $$
    (G,\phi) \;\;\in\; \mathrm{Groups}(\mathbf{H})_{/\mathrm{Aut}(T_e V)}
  $$

  \vspace{-2mm}
\noindent
  is a group $G$ equipped with a homomorphism of groups
  $\!\xymatrix@C=20pt{ G \ar[r]|-{\,\phi \,} & \mathrm{Aut}(T_e V) }\!$
  to the structure group (Def. \ref{StructureGroupOfVFolds}) of
  $V$-folds.
  Under delooping \eqref{LoopingDeloopingEquivalence} this is equivalently a
  morphism in $\mathbf{H}$ of the form
  $\xymatrix@C=28pt{
    \mathbf{B}G \ar[r]|-{ \;\mathbf{B}\phi \,}
    &
    \mathbf{B}\mathrm{Aut}(T_e V)
  }$.
\end{defn}

\begin{defn}[G-structures on $V$-folds]
  \label{GStructures}
  Let $\mathbf{H}$ be an elastic $\infty$-topos  (Def. \ref{ElasticInfinityTopos}),
  $V \in \mathrm{Groups}(\mathbf{H})$ (Prop. \ref{LoopingAndDelooping}),
  $(G,\phi) \in \mathrm{Groups}(\mathbf{H})_{/\mathrm{Aut}(T_e V)}$
  (Def. \ref{GStructureCoefficients})
  and $X \in V\mathrm{Folds}(\mathbf{H})$ (Def. \ref{VManifold}).

 \vspace{-1mm}
 \noindent \item {\bf (i)} We say that
  \begin{itemize}
   \vspace{-.2cm}
   \item
   \emph{a \emph{$(G,\phi)$-structure} on $X$}
   (often just \emph{$G$-structure} if $\phi$ is understood),
   \vspace{-.2cm}
   \item
   or
   \emph{$(G,\phi)$-structure on its
  frame bundle} (Def. \ref{TangentBundleOfVFoldIsFiverBundle}),
   \vspace{-.2cm}
  \item
  or \emph{reduction of the structure group \eqref{StructureGroupOfVFolds}
  along $\phi$}
  \end{itemize}
  \vspace{-.2cm}
  is a lift $(\tau,g)$ of the frame bundle classifying
  map \eqref{FrameBundleClassifyingMap} through $\mathbf{B}\phi$:
  \begin{equation}
    \label{GphiStructure}
    \raisebox{20pt}{
    \xymatrix@C=3em@R=1.5em{
      &&
      \mathbf{B}G
      \ar[d]^-{ \mathbf{B}\phi }
      \\
      \mathllap{
        \mbox{
          \tiny
          \color{darkblue}
          \bf
          $V$-fold
        }
        \;\;\;
      }
      X
      \ar[rr]_-{ \scalebox{.6}{$\vdash\, \mathrm{Frames}(X)$} }^>>>>>{\ }="t"
      \ar@{-->}[urr]^-{
        \overset{
          \mathclap{
          \mbox{
            \tiny
            \color{darkblue}
            \bf
            \begin{tabular}{c}
              $G$-structure
              \\
              \phantom{A}
            \end{tabular}
          }
          }
        }{
          \tau
        }
      }_>>>>>>>>>{\ }="s"
      &&
      \mathbf{B}
      \mathrm{Aut}(T_e V)
      \mathrlap{
        \mbox{
          \tiny
          \color{darkblue}
          \bf
          \begin{tabular}{c}
            structure group
            \\
            of frame bundle
          \end{tabular}
        }
      }
      \ar@{<=}^-{ g } "s"; "t"
    }
    }
  \end{equation}

  \vspace{-2mm}
\noindent \item{\bf (ii)} We say that the \emph{$G$-frame bundle}
  $G \mathrm{Frames}(X)$
  of a $V$-fold $X$ equipped with such a
  $(G,\phi)$-structure  is the
  $G$-principal bundle which is classified
  (via Prop. \ref{ClassificationPrincipalInfinityBundles}):
  by $\tau$, hence the object in the following diagram:
  \vspace{-1mm}
  \begin{equation}
    \label{GFrameBundle}
    \xymatrix@C=3em@R=1.2em{
      G\mathrm{Frames}(X,\tau)
      \ar[d]
      \ar[rr]
      \ar@{}[drr]|-{\mbox{\tiny\rm(pb)}}
      &&
      \ast
      \ar[dd]|-{\phantom{AA}}
      \\
      \mathrm{Frames}(X)
      \ar[d]
      \ar[rrrr]
       \ar@{}[drrr]|-{\mbox{\tiny\rm(pb)}}
     &&&&
      \ast
      \ar[d]
      \\
      X
      \ar[rr]_-{ \tau }
      \ar@/_1.7pc/[rrrr]_-{\scalebox{.6}{$ \vdash\mathrm{Frames}(X)$} }
      &&
      \mathbf{B}G
      \ar[rr]_-{ \mathbf{B}\phi }
      &&
      \mathbf{B} \mathrm{Aut}(T_e V)
    }
  \end{equation}

  \vspace{-4mm}
\item {\bf (iii)}  We write

  \vspace{-3mm}
  \begin{equation}
    \label{GroupoidOfGStructures}
    (G,\phi)\mathrm{Structures}_{X}(\mathbf{H})
    \;:=\;
    \mathbf{H}_{\!/_{\mathbf{B}\mathrm{Aut}(T_e V)}}
    \big(
      \vdash\mathrm{Frames}(X)
      \,,\,
      \mathbf{B}\phi
    \big)
    \;\in\;
    \mathrm{Groupoids}_\infty
  \end{equation}
  for the $\infty$-groupoid of $(G,\phi)$-structures
  on the $V$-fold $X$.
\end{defn}

In direct generalization of Prop. \ref{VFoldAsQuotientOfItsFrameBundle} we have:

\begin{prop}[$G$-structured $V$-fold is $G$-quotient of its $G$-frame bundle]
  \label{GStructuredVFoldAsQuotientOfGFrameBundle}
  Let $\mathbf{H}$ be an elastic $\infty$-topos  (Def. \ref{ElasticInfinityTopos}),
  $V \in \mathrm{Groups}(\mathbf{H})$ (Prop. \ref{LoopingAndDelooping}),
  $(G,\phi) \in \mathrm{Groups}(\mathbf{H})_{/\mathrm{Aut}(T_e V)}$
  (Def. \ref{GStructureCoefficients}),
  $X \in V\mathrm{Folds}(\mathbf{H})$ (Def. \ref{VManifold})
  and $(\tau,g) \in (G,\phi)\mathrm{Structures}_X(\mathbf{H})$
  (Def. \ref{GStructures}).
    Then
\vspace{-1mm}
  \item {\bf (i)}
  $X$  is equivalently the homotopy quotient \eqref{HomotopyQuotientAsColimit}
  of its
  $G$-frame bundle \eqref{GFrameBundle} by $G$:
  \vspace{-1mm}
  $$
    X \;\simeq\; G \mathrm{Frames}(X,\tau) \!\sslash\! G
    \,.
  $$

   \vspace{-3mm}
 \item {\bf (ii)} the classifying map of the $G$-frame bundle on $X$
  exhibits the action of $G$ on $G \mathrm{Frames}(X,\tau)$
  according to \eqref{InfinityActionHomotopyFiberSequence}.
\end{prop}
\begin{proof}
  This is immediate from the equivalence between principal bundles
  and homotopy quotient projections
  (Remark \ref{PrincipalBaseSpacesAreHomotopyQuotients}) applied to
  the $G$-frame bundle \eqref{GFrameBundle}:
  \vspace{-2mm}
  $$
  \hspace{1cm}
    \xymatrix@R=1.2em{
      G \mathrm{Frames}(X,\tau)
      \ar[d]_-{ \mathrm{fib}(\rho) \simeq \mathrm{fib}(\tau) }
      \\
      G\mathrm{Frames}(X,\tau) \!\sslash\! G
      \ar@/_1pc/[rr]_-{\rho_G}
      \ar@{}[r]|-{\simeq}
      &
      X
      \ar[r]^-{\tau}
      &
      \mathbf{B}G
    }
  $$

  \vspace{-7mm}
\hfill \end{proof}

\begin{example}[$G$-structure induced from framing]
  \label{GStructuresInducedFromFramings}
  Let $\mathbf{H}$ be an elastic $\infty$-topos (Def. \ref{ElasticInfinityTopos}),
  $V \in \mathrm{Groups}(\mathbf{H})$ (Prop. \ref{LoopingAndDelooping})
  and
  $X \in V \mathrm{Folds}(\mathbf{H})$ (Def. \ref{VManifold}).
  Then a framing on $X$ (Def. \ref{Framing})
  induces a $(G,\phi)$-structure (Def. \ref{GStructures}) for any
  $(G,\phi) \in \mathrm{Groups}(\mathbf{H})_{/\mathrm{Aut}(T_e V)}$,
  given by the pasting
  \vspace{-2mm}
  \begin{equation}
    \label{GStructureInducedFromFraming}
    \raisebox{20pt}{
    \xymatrix@C=3em{
      &
      \ast
      \ar[r]
      \ar[dr]_<<{\ }="s"^-{\ }="t2"
      &
      \mathbf{B}G
      \ar[d]^-{ \mathbf{B}\phi }_<<{\ }="s2"
      \\
      X
      \ar[ur]
      \ar[rr]_-{ \scalebox{.6}{$\vdash \mathrm{Frames}(X)$} }^-{\ }="t"
      &&
      \mathbf{B}\mathrm{Aut}(T_e V)
      \ar@{<=}|-{
        \mathclap{\phantom{\vert^\vert}}
        \vdash\mathrm{fr}
        \mathclap{\phantom{\vert_\vert}}
      } "s"; "t"
      \ar@{<=} "s2"; "t2"
    }
    }
  \end{equation}

  \vspace{-2mm}
\noindent
  of the homotopy $\vdash\mathrm{fr}$ \eqref{FramingAsTrivializationOfFrameBundle}
  which classifies the framing
  (Remark \ref{FramingVFold})  with the homotopy that
  exhibits the group homomorphism $\phi$ as a morphism
  of pointed objects (Prop. \ref{LoopingAndDelooping}).
\end{example}

\newpage

\begin{example}[Canonical $G$-structure]
  \label{CanonicalGStructure}
  Let $\mathbf{H}$ be an elastic $\infty$-topos $\mathbf{H}$ (Def. \ref{ElasticInfinityTopos}),
  and $V \in \mathrm{Groups}(\mathbf{H})$ (Prop. \ref{LoopingAndDelooping}).
  Then $V$ itself, regarded as a $V$-fold by Example \ref{VAsAVFold},
  carries a $(G,\phi)$-structure (Def. \ref{GStructures}) for any
  $(G,\phi) \in \mathrm{Groups}(\mathbf{H})_{/\mathrm{Aut}(T_e V)}$,
  induced via Example \ref{GStructuresInducedFromFramings}
  from its canonical framing $\mathrm{fr}_\ell$ \eqref{LeftTranslationFramingOnGroup}
  via left-translation
  (Prop. \ref{GroupsAdmitFraming}). We call this the
  \emph{canonical $(G,\phi)$-structure} on $V$:
   \vspace{-2mm}
  \begin{equation}
    \label{GStructureInducedFromFraming}
    \raisebox{20pt}{
    \xymatrix{
      &&
      \mathbf{B}G
      \ar[d]^-{ \mathbf{B}\phi }
      \\
      V
      \ar@{-->}[urr]^-{ \tau_V }_>>>>>>>>{\ }="t"
      \ar[rr]^{\ }="s"
      &&
      \mathbf{B}\mathrm{Aut}(T_e V)
      \ar@{=>}_-{ g_V } "s"; "t"
    }
    }
    \phantom{AAA}
    :=
    \phantom{AAA}
    \raisebox{20pt}{
    \xymatrix{
      &
      \ast
      \ar[r]
      \ar[dr]_<<{\ }="s"^-{\ }="t2"
      &
      \mathbf{B}G
      \ar[d]^-{ \mathbf{B}\phi }_<<{\ }="s2"
      \\
      V
      \ar[ur]
      \ar[rr]_-{\scalebox{.6}{$ \vdash \mathrm{Frames}(V)$} }^-{\ }="t"
      &&
      \mathbf{B}\mathrm{Aut}(T_e V)
      \ar@{<=}|-{
        \mathclap{\phantom{\vert^\vert}}
        \vdash\mathrm{fr}_\ell
        \mathclap{\phantom{\vert_\vert}}
      } "s"; "t"
      \ar@{<=} "s2"; "t2"
    }
    }
  \end{equation}
\end{example}

\medskip

\noindent {\bf Local isometries}
\begin{lemma}[$G$-structures pull back along local diffeomorphisms]
  \label{GStructuresPullBackAlongLocalDiffeomorphisms}
  Let $\mathbf{H}$ be an elastic $\infty$-topos  (Def. \ref{ElasticInfinityTopos}),
  $V \in \mathrm{Groups}(\mathbf{H})$ (Prop. \ref{LoopingAndDelooping})
  and
  $(G,\phi) \in \mathrm{Groups}(\mathbf{H})_{/\mathrm{Aut}(T_e V)}$
    (Prop. \ref{LoopingAndDelooping},  Def. \ref{AutomorphismGroup},
    Example \ref{LocalNeighbourhoodOfAPoint}).
  Then pre-composition
  constitutes a contravariant $\infty$-functor
  (``pullback of $(G,\phi)$-structures'')
   \vspace{-3mm}
  \begin{equation}
    \label{PullbackOfGStructures}
    \xymatrix@R=3pt{
      \big(
        V\mathrm{Folds}(\mathbf{H})^{\mbox{\tiny\rm{\'e}t}}
      \big)^{\mathrm{op}}
      \ar[rr]
      &&
      \mathrm{Groupoids}_\infty
      \\
      X_1
      \ar[dd]^-{ f }_-{\mbox{\tiny\rm{\'e}t}}
      \ar@{}[rr]|-{\longmapsto}
      &&
      (G,\phi)\mathrm{Structures}_{X_1}(\mathbf{H})
      \ar@{<-}[dd]^-{f^\ast}
      \ar@{}[r]|-{
        \begin{rotate}{180}
          $\mathclap{\in}$
        \end{rotate}
      }
      &
      \tau \circ f
      \ar@{<-|}[dd]
      \\
      \\
      X_2
      \ar@{}[rr]|-{\longmapsto}
      &&
      (G,\phi)\mathrm{Structures}_{X_2}(\mathbf{H})
      \ar@{}[r]|-{
        \begin{rotate}{180}
          $\mathclap{\in}$
        \end{rotate}
      }
      &
      \tau
    }
  \end{equation}

   \vspace{-2mm}
\noindent
  from the $\infty$-category \eqref{VFoldsAndLocalDiffeomorphismsCategory}
  of $V$-folds and local diffeomorphisms,
  which assigns to any $V$-fold its $\infty$-groupoid \eqref{GroupoidOfGStructures}
  of $(G,\phi)$-structres (Def. \ref{GStructures}).
\end{lemma}
\begin{proof}
  We need to show that for $(\tau,g)$ a $(G,\phi)$-structure on $X_2$,
  the composite
   \vspace{-3mm}
  \begin{equation}
    \label{PullbackOfGStructures}
    \xymatrix@C=3em@R=1em{
      & && \mathbf{B} G
      \ar[d]^-{ \mathbf{B} \phi }
      \\
      X_1
      \ar[r]_-{ f }^-{\mbox{\tiny\rm{\'e}t}}
      &
      X_2
      \ar[urr]^-{ \tau }_>>>>>>>{\ }="t"
      \ar[rr]_-{  \scalebox{.6}{$\vdash\mathrm{Frames}(X_2)$} }^-{\ }="s"
      &&
      \mathbf{B} \mathrm{Aut}(T_e V)
      \ar@{=>}|{\; g\,} "s"; "t"
    }
  \end{equation}

   \vspace{-2mm}
\noindent
  is a $(G,\phi)$-structure on $X_1$.
  For this we need to exhibit a natural equivalence
   \vspace{-3mm}
  $$
  \hspace{-1cm}
    \big( \vdash\mathrm{Frames}(X_2) \big) \circ f_1
    \;\;\simeq\;\;
    \vdash\mathrm{Frames}(X_1)
    \phantom{AA}
    \mbox{so that}
    \phantom{AA}
    \raisebox{20pt}{
    \xymatrix@C=3em{
      X_1
        \ar[rr]^-{ f }|-{\; \mbox{\tiny{\'e}t}\; }_>>>{\ }="t1"
        \ar[drr]_{ \scalebox{.6}{$
          \vdash \mathrm{Frames}(X_1)
          $}
        }^-{\ }="s1"
      &&
      X_2
      \ar[d]|>>>>{ \scalebox{.6}{$\vdash \mathrm{Frames}(X_2)$}}^<<<{\ }="s2"
      \ar[rr]^-{\tau}_-{\ }="t2"
      &&
      \mathbf{B}G
      \ar[dll]^-{ \mathbf{B}\phi }
      \\
      &&
      \mathbf{B}\mathrm{Aut}(T_e V)
      \ar@{=>}|-{\;\simeq\;} "s1"; "t1"
      \ar@{=>}|-{\;g\;} "s2"; "t2"
    }
    }
  $$

   \vspace{-2mm}
\noindent
  But this exists by Prop. \ref{PullbackAlongLocalDiffeomorphismsPreservesTangentBundles}.
\hfill \end{proof}

\begin{defn}[Local isometries between $G$-structured $V$-folds]
  \label{Isometries}
  Let $\mathbf{H}$ be an elastic $\infty$-topos  (Def. \ref{ElasticInfinityTopos}),
  $V \in \mathrm{Groups}(\mathbf{H})$ (Prop. \ref{LoopingAndDelooping})
  and
  $(G,\phi) \in \mathrm{Groups}(\mathbf{H})_{/\mathrm{Aut}(T_e V)}$
    (Prop. \ref{LoopingAndDelooping},  Def. \ref{AutomorphismGroup},
    Example \ref{LocalNeighbourhoodOfAPoint}).

  \noindent {\bf (i)}
  For $X_1, X_2 \in V \mathrm{Folds}$ (Def. \ref{VManifold})
  and $(\tau_i, g_i) \in (G,\phi)\mathrm{Structures}_{X_i}(\mathbf{H})$
  \eqref{GroupoidOfGStructures},
  we say a \emph{local isometry}, to be denoted
  \vspace{-2mm}
  $$
    \xymatrix{
      \big(
        X_1, (\tau_1, g_1)
      \big)
      \ar[rr]_-{ \mbox{\tiny met} }^{ (f,\sigma) }
      &&
      \big(
        X_2, (\tau_2, g_2)
      \big)
         }
  $$

  \vspace{-2mm}
\noindent
  is a pair
   \vspace{-2mm}
  \begin{equation}
    \label{IsometryExplicitly}
    \hspace{1cm}
    \xymatrix@R6pt{
      X_1
        \ar[r]^-{f}_-{ \mbox{\tiny{\'e}t} }
      &
      X_2
    }
    \;\;,
    \phantom{AAA}
    \xymatrix{
      f^\ast (\tau_2, g_2)
      \ar[r]_-{\simeq}^-{\sigma}
      &
      (\tau_1, g_1)\;,
    }
  \end{equation}
  consisting of a local diffeomorphism (Def. \ref{FormallyEtaleMorphism})
  and an equivalence of $(G,\phi)$-structures \eqref{GroupoidOfGStructures}
  between that on its domain $V$-fold
  and the pullback \eqref{PullbackOfGStructures}
  of the $(G,\phi)$-structure on its codomain $V$-fold.

  \noindent {\bf (ii)} Equivalently, by \eqref{PullbackOfGStructures},
  a local isometry \eqref{IsometryExplicitly}
  is a morphism between $(G,\phi)$-structured $V$-folds
  regarded as objects in the iterated slice
  $\infty$-topos (Example \ref{IteratedSliceTopos})

  {\bf (a)} over $\mathbf{B} \mathrm{Aut}(T_e V)$ via their classifying maps
  of their frame bundles \eqref{FrameBundleClassifyingMap}

  {\bf (b)} over $\big(\mathbf{B}G , \mathbf{B}\phi\big)$ via their $(G,\phi)$-structure
  \eqref{GphiStructure}

  \noindent of this form:
  \vspace{-2mm}
  \begin{equation}
    \label{IsometriesAsMorphismsInSlice}
    \hspace{-.5cm}
    \raisebox{30pt}{
    \xymatrix@R=20pt@C=5em{
      X_1
      \ar[ddr]^>>>>>{\ }="s1"
      \ar@{}[ddr]|-{
        \rotatebox[origin=t]{-38}{
          \scalebox{.63}{
            $\mathclap{
              \begin{array}{c}
              \phantom{-}
              \\
              {\!\!\!\vdash \mathrm{Frames}(X_1)}
              \end{array}
            }$
          }
        }
      }^>>>>>>>{\ }="s4"
      \ar@/_.46pc/[drr]|>>>>>>>{\,\tau_1\,}_->>>>>{\ }="t1"^-<<<<<<{\ }="t2"
      \ar[r]^-f|-{\;\mbox{\tiny{\'et}}\;}
      &
      X_2
      \ar@/^.55pc/[dr]^>>>>>>{\,\tau_2\,}_>>>>>>{\ }="s2"
      \ar[dd]|<<<{\phantom{A\vert}}|<<<<<<{\phantom{AA}}|>>>>>{\phantom{AA} \atop \phantom{a}}^<<<{\ }="s3"_-{\ }="t4"
      \\
      &&
      \mathbf{B}G
      \\
      &
      \mathbf{B}\mathrm{Aut}(T_e V)
      \ar@{<-}[ur]_-{ \mathbf{B} \phi }^-{\ }="t3"
      \ar@/_.23pc/@{=>} "s1"; "t1"|<<<<<<<<<<<<<<<<<{\,g_1\,}
      \ar@{=>}|<<<<<<<<<<<<<{\,\sigma\,} "s2"; "t2"
      \ar@<+2pt>@{=>} "t3"+(-16,4); "t3"|<<<<{\,g_2\,}|>>>>>>{\phantom{AAA}}
      \ar@{=>} "s4"; "t4"
    }
    }
    \;\;\;
    \in
    \;\;
    \big(
      \mathbf{H}_{/\mathbf{B}\mathrm{Aut}(T_e V)}
    \big)_{\!\!/(\mathbf{B}G , \mathbf{B}\phi)}
    \Big(
      \big(
        X_1, (\tau_1,g_1)
      \big)
      \,,\,
      \big(
        X_2, (\tau_2,g_2)
      \big)
    \Big).
  \end{equation}

\vspace{-2mm}
  \noindent {\bf (iii)} Hence we write

  \vspace{-1cm}

  \begin{equation}
    \label{CategoryOfGStructuredVFolds}
    (G,\phi)\mathrm{Structured}V\mathrm{Folds}(\mathbf{H})
    \;\longrightarrow\;
    \big(
      \mathbf{H}_{/\mathbf{B}\mathrm{Aut}(T_e V)}
    \big)_{\!\!/\mathbf{B}G}
    \;\;
    \in
    \mathrm{Categories}_\infty
  \end{equation}

  \vspace{-2mm}
\noindent
  for the sub-$\infty$-category of this iterated slice
  on 1-morphisms of the form \eqref{IsometriesAsMorphismsInSlice}.
\end{defn}

\medskip

\noindent {\bf Integrability of  $G$-structures.}
\begin{defn}[Integrable $G$-structure]
  \label{IntegrableGStructure}
  Let $\mathbf{H}$ be an elastic $\infty$-topos (Def. \ref{ElasticInfinityTopos}),
  $V \in \mathrm{Groups}(\mathbf{H})$ (Prop. \ref{LoopingAndDelooping}),
  $(G,\phi) \in \mathrm{Groups}(\mathbf{H})_{/\mathrm{Aut}(T_e V)}$
  (Def. \ref{GStructureCoefficients}).

  \noindent {\bf (i)}
  Given
  $(X, (\tau_X,g_X)) \in
    (G,\phi)\mathrm{Structured}V\mathrm{Folds}(\mathbf{H})$
    (Def. \ref{Isometries}),
  we say that $(\tau,g)$ is an \emph{integrable} $(G,\phi)$-structure
  on the $V$-fold $X$ if there exists a correspondence
  of local isometries \eqref{IsometryExplicitly}
  between $V$ equipped with its canonical $(G,\phi)$-structure
  $(\tau_V, g_V)$
  (Def. \ref{CanonicalGStructure}) to $(X,(\tau_X,g_X))$:
    \vspace{-2mm}
  \begin{equation}
    \label{IsometryAtlas}
    \raisebox{10pt}{
    \xymatrix@R=-2pt@C=3em{
      &
      \big(
        U,
        (\tau_U, g_U)
      \big)
      \ar[dl]_-{\mbox{\tiny met}}
      \ar@{->>}[dr]^-{\mbox{\tiny met}}
      \\
      \big(
        V,
        (\tau_V, g_V)
      \big)
      &&
      \big(
        X,
        (\tau_X, g_X)
      \big)
    }
    }
  \end{equation}

  \vspace{-2mm}
  \noindent
  such that
  the right left is, in addition, an effective epimorphism
  (Def. \ref{EffectiveEpimorphisms}),
  then called a {\it $(V, (\tau_V, g_V))$-atlas of $(X,(\tau_X, g_X))$}
  \eqref{StacksAtlasesAndGroupoids}.
  (Underlying this, forgetting the $(G,\phi)$-structures,
  is a $V$-atlas \eqref{VAtlas}.)

  \noindent {\bf (ii)} We write
  \vspace{-2mm}
  \begin{equation}
    \label{IntegrablyGStructuredVFolds}
    \mathrm{Integrably}(G,\phi)\mathrm{Structured}V\mathrm{Folds}(\mathbf{H})
      \xymatrix{
        \;   \ar@{^{(}->}[r]
      &
}
    (G,\phi)\mathrm{Structured}V\mathrm{Folds}(\mathbf{H})
    \;\;
    \in
    \mathrm{Categories}_\infty
  \end{equation}

  \vspace{-2mm}
\noindent
  for the full sub-$\infty$-category of that of
  $(G,\phi)$-structured $V$-folds \eqref{CategoryOfGStructuredVFolds}
  on those that are integrable.
\end{defn}
\begin{defn}[Locally integrable $G$-structure]
  \label{LocallyIntegrableGStructures}
  Let $\mathbf{H}$ be an elastic $\infty$-topos (Def. \ref{ElasticInfinityTopos}),
  $V \in \mathrm{Groups}(\mathbf{H})$ (Prop. \ref{LoopingAndDelooping}),
  $(G,\phi) \in \mathrm{Groups}(\mathbf{H})_{/\mathrm{Aut}(T_e V)}$
  (Def. \ref{GStructureCoefficients}),
  $X \in V \mathrm{Folds}(\mathbf{H})$ (Def. \ref{VManifold})
  and
  $(\tau,g) \in (G,\phi)\mathrm{Structures}_X(\mathbf{H})$
  (Def. \ref{GStructures}).
  We say that $(\tau,g)$
  is a \emph{locally integrable} $(G,\phi)$-structure if, for each point
  $\!\xymatrix{\ast \ar[r]^-{x} & X}\!\!\!$, there is a local diffeomorphism
  $\phi_x$
  of the local neighborhood (Def. \ref{InfinitesimalNeighborhood})
  of $\!\xymatrix{\ast \ar[r]^-{e} & V}\!\!$ onto a local neighborhood of $x$
  such that the restriction of $(\tau,g)$ along $\phi$
  is equivalent to the canonical $(G,\phi)$-structure
  (Def. \ref{CanonicalGStructure}) on $T_e V$:
\vspace{-2mm}
  $$
    \underset{
      \scalebox{.75}{
      \raisebox{0pt}{
      \xymatrix{
        \ast \ar[r]^-x & X
      }
      }}
    }{
      \forall
    }
    \;\;
    \underset{
    \scalebox{.75}{
    \raisebox{0pt}{
    \xymatrix@R=12.5pt@C=10pt{
      T_e V
      \ar[rr]^-{\phi_x}_-{\mbox{\tiny{\'e}t}}
      &&
      X
      \\
      &
      \ast
      \ar[ur]_-{x}
      \ar[ul]^-{e}
    }
    \;\;
    }
    }
    }{
      \exists
    }
    \;:\;
    \phantom{AA}
        \phi_x^\ast (\tau, g)
    \;\simeq\;
    (\tau_{T_e V}, g_{T_e V})
    \,.
     $$

    \vspace{-.3cm}

    Another way to say this: We have a correspondence of
    local isometries as in \eqref{IsometryAtlas}, but with the
    right leg required to be an effective epimorphism only
    under $\flat$.
\end{defn}

\begin{example}[$G$-Structures on smooth manifolds and orbifolds]
  \label{TorsionFreeGStructureOnSmoothManifolds}
$\phantom{A}$
\vspace{-2mm}
  \item {\bf (i)}
  Let $\mathbf{H} = \mathrm{JetsOfSmoothGroupoids}_\infty$
  (Example \ref{FormalSmoothInfinityGroupoids})
  $G \in \!
    \xymatrix{
      \mathrm{LieGroups}
      \longhookrightarrow
      \mathrm{Groups}
      \big(
        \mathbf{H}
      \big)
    }
    \!\!
  $ (see \eqref{ConcreteSmoothFormalInfinityGroupoidsAreDiffeologicalSpaces})
  and $X \in$  \newline $\mathrm{SmoothManifolds} \longhookrightarrow \mathbf{H}$
  regarded as an $\mathbb{R}^n$-fold according to Example \ref{OrdinaryManifolds}.
  In this case, the structure group of $X$ (Def. \ref{StructureGroupOfVFolds})
  is the ordinary general linear group $\mathrm{GL}_{\mathbb{R}}(n)$
  (Example \ref{OrdinaryGeneralLinearGroup}).
  Therefore, a $G$-structure on $X$ in the sense of Def. \ref{GStructures}
  is
  (by Example \ref{FrameBundleOnSmoothManifold})
  a $G$-structure  in the traditional sense of differential geometry
  \cite[VII]{Sternberg64}\cite{Kobayashi72}\cite{Molino77};
  and it is
  integrable according to Def. \ref{IntegrableGStructure}
  if it is ``flat'' in the traditional sense of \cite{Guillemin65}
  and
  locally integrable according to Def. \ref{LocallyIntegrableGStructures}
  precisely if it is ``uniformly 1-flat''
  in the traditional sense of \cite{Guillemin65}, namely
  if it is torsion-free (review in \cite{Lott01}).
  Examples include:

\vspace{2mm}
{\small
  \hspace{-.7cm}
  \begin{tabular}{|r||c|c|c||c|}
    \hline
    $
      \xymatrix@C=10pt{
        G
        \ar[r]^-{\phi}
        &
        \mathrm{GL}_{{}_{\mathbb{R}}}(n)
      }
    $
    &
    \bf
    $(G,\phi)$-structure
    &
    \bf
    Locally integrable
    &
    \bf
    Integrable
    &
    see
    \\
    \hline
    \hline
    $
      \xymatrix@C=10pt{
        \mathrm{Sp}_{{}_{\mathbb{R}}}(n)
        \;
        \ar@{^{(}->}[r]
        &
        \mathrm{GL}_{{}_{\mathbb{R}}}(n)
      }
    $
    &
    \begin{tabular}{c}
      almost
      \\
      symplectic
    \end{tabular}
    &
    \begin{tabular}{c}
      \phantom{almost}
      \\
      symplectic
    \end{tabular}
    &
    \begin{tabular}{c}
      \phantom{almost}
      \\
      symplectic
    \end{tabular}
    &
    \multirow{4}{*}{
      \scalebox{.9}{
        \cite[VII.2]{Sternberg64}
      }
    }
    \\
    \cline{1-4}
    $
      \xymatrix@C=10pt{
        \mathrm{GL}_{{}_{\mathbb{C}}}(n/2)
        \;
        \ar@{^{(}->}[r]
        &
        \mathrm{GL}_{{}_{\mathbb{R}}}(n)
      }
    $
    &
    \begin{tabular}{c}
      almost
      \\
      complex
    \end{tabular}
    &
    \begin{tabular}{c}
      \phantom{almost}
      \\
      complex
    \end{tabular}
    &
    \begin{tabular}{c}
      \phantom{almost}
      \\
      complex
    \end{tabular}
    &
    \\
    \cline{1-4}
    $
      \xymatrix@C=10pt{
        \mathrm{O}(n)
        \;
        \ar@{^{(}->}[r]
        &
        \mathrm{GL}_{{}_{\mathbb{R}}}(n)
      }
    $
    &
    \begin{tabular}{c}
      \phantom{flat}
      \\
      Riemannian
    \end{tabular}
    &
    \begin{tabular}{c}
      torsion-free
      \\
      Riemannian
    \end{tabular}
    &
    \begin{tabular}{c}
      flat
      \\
      Riemannian
    \end{tabular}
    &
    \\
    \hline
    $
      \xymatrix@C=10pt{
        \mathrm{O}(n-1,1)
        \;
        \ar@{^{(}->}[r]
        &
        \mathrm{GL}_{{}_{\mathbb{R}}}(n)
      }
    $
    &
    \begin{tabular}{c}
      \phantom{flat}
      \\
      Lorentzian
    \end{tabular}
    &
    \begin{tabular}{c}
      torsion-free
      \\
      Lorentzian
    \end{tabular}
    &
    \begin{tabular}{c}
      flat
      \\
      Lorentzian
    \end{tabular}
    &
    \scalebox{.9}{
      \cite{LPZ13}
    }
    \\
    \hline
    $
      \xymatrix@C=10pt{
        \mathrm{O}(n)\times \mathbb{R}
        \;
        \ar@{^{(}->}[r]
        &
        \mathrm{GL}_{{}_{\mathbb{R}}}(n)
      }
    $
    &
    \begin{tabular}{c}
      \phantom{flat}
      \\
      $\mathrm{CO}(n)$-structure
    \end{tabular}
    &
    \begin{tabular}{c}
      \phantom{flat}
      \\
      conformal
    \end{tabular}
    &
    \begin{tabular}{c}
      flat
      \\
      conformal
    \end{tabular}
    &
    \scalebox{.9}{
      \cite{AkivisGoldberg98}
    }
    \\
    \hline
    $
      \xymatrix@C=10pt{
        \mathrm{CR}(n/2-1)
        \;
        \ar@{^{(}->}[r]
        &
        \mathrm{GL}_{{}_{\mathbb{R}}}(n)
      }
    $
    &
    \begin{tabular}{c}
      \phantom{flat}
      \\
      $\mathrm{CR}(n)$-structure
    \end{tabular}
    &
    \begin{tabular}{c}
      \phantom{flat}
      \\
      Cauchy-Riemann
    \end{tabular}
    &
    \begin{tabular}{c}
      flat
      \\
      Cauchy-Riemann
    \end{tabular}
    &
    \scalebox{.9}{
      \cite{DragomiTomassini06}
    }
    \\
    \hline
    $
      \xymatrix@C=10pt{
        \mathrm{GL}_{{}_{\mathbb{H}}}(n/4)
        \;
        \ar@{^{(}->}[r]
        &
        \mathrm{GL}_{{}_{\mathbb{R}}}(n)
      }
    $
    &
    \begin{tabular}{c}
      \phantom{flat}
      \\
      $\mathrm{GL}_{{}_{\mathbb{H}}}(n/4)$-structure
    \end{tabular}
    &
    \begin{tabular}{c}
      \phantom{flat}
      \\
      hypercomplex
    \end{tabular}
    &
    \begin{tabular}{c}
      flat
      \\
      hypercomplex
    \end{tabular}
    &
    \scalebox{.9}{
      \cite{Joyce95}
    }
    \\
    \hline
    \hline
    $
      \xymatrix@C=10pt{
        \mathrm{U}(n/2)
        \;
        \ar@{^{(}->}[r]
        &
        \mathrm{GL}_{{}_{\mathbb{R}}}(n)
      }
    $
    &
    \begin{tabular}{c}
      hermitian
      \\
      almost complex
    \end{tabular}
    &
    K{\"a}hler
    &
    K{\"a}hler
    &
    \scalebox{.9}{
      \cite[11.1]{Moroianu07}
    }
    \\
    \hline
    $
      \xymatrix@C=10pt{
        \mathrm{SU}(n/2)
        \;
        \ar@{^{(}->}[r]
        &
        \mathrm{GL}_{{}_{\mathbb{R}}}(n)
      }
    $
    &
    $\mathrm{SU}(n)$-structure
    &
    Calabi-Yau
    &
    Calabi-Yau
    &
    \scalebox{.9}{
      \cite[1.3]{Prins15}
    }
    \\
    \hline
    $
      \!\!\!\!\!\!
      \xymatrix@C=10pt{
        \mathrm{Sp}(n/4) \!{\boldsymbol\cdot}\! \mathrm{Sp}(1)
        \;
        \ar@{^{(}->}[r]
        &
        \mathrm{GL}_{{}_{\mathbb{R}}}(n)
      }
    $
    &
  \!\!\!\!\!\!\!  \begin{tabular}{c}
      almost unimodular
      \\
      quaternionic
    \end{tabular}
    &
  \!\!\!\!\!  \begin{tabular}{c}
      \phantom{almost}
      \\
      quaternionic K{\"a}hler
    \end{tabular}
    \!\!\!
    &
  \!\!\!\!  \begin{tabular}{c}
      flat
      \\
      quaternionic K{\"a}hler
    \end{tabular}
    \!\!
    &
    \multirow{2}{*}{\!\!\!\!\!
      \scalebox{.9}{
        \cite{AlekseevskyMarchiafava93a}\cite{AlekseevskyMarchiafava93b}
      }
   \!\!\!\!\!\! }
    \\
    \cline{1-4}
    $
      \xymatrix@C=10pt{
        \mathrm{Sp}(n/4)
        \;
        \ar@{^{(}->}[r]
        &
        \mathrm{GL}_{{}_{\mathbb{R}}}(n)
      }
    $
    &
    \begin{tabular}{c}
      almost
      \\
      Hyperk{\"a}hler
    \end{tabular}
    &
    \begin{tabular}{c}
      \phantom{almost}
      \\
      Hyper{\"a}hler
    \end{tabular}
    &
    \begin{tabular}{c}
      flat
      \\
      Hyperk{\"a}hler
    \end{tabular}
    &
    \\
    \hline
    $
      \xymatrix@C=10pt{
        G_2
        \;
        \ar@{^{(}->}[r]
        &
        \mathrm{GL}_{{}_{\mathbb{R}}}(7)
      }
    $
    &
    \begin{tabular}{c}
      $G_2$-structure
    \end{tabular}
    &
    \begin{tabular}{c}
      torsion-free
      \\
      $G_2$-structure
    \end{tabular}
    &
    \begin{tabular}{c}
      flat/interable
      \\
      $G_2$-structure
    \end{tabular}
    &
    \scalebox{.9}{
      \cite{Bryant05}
    }
    \\
    \hline
    $
      \xymatrix@C=10pt{
        \mathrm{Spin}(7)
        \;
        \ar@{^{(}->}[r]
        &
        \mathrm{GL}_{{}_{\mathbb{R}}}(8)
      }
    $
    &
    \begin{tabular}{c}
      $\mathrm{Spin}(7)$-structure
    \end{tabular}
    &
    \begin{tabular}{c}
      torsion-free
      \\
      $\mathrm{Spin}(7)$-structure
    \end{tabular}
    &
    \begin{tabular}{c}
      flat
      \\
      $\mathrm{Spin}(7)$-structure
    \end{tabular}
    &
    \scalebox{.9}{
      \begin{tabular}{c}
        \cite{Bryant87}
        \\
        \cite{Joyce01}
      \end{tabular}
    }
    \\
    \hline
  \end{tabular}
  }
\vspace{-1mm}
 \item {\bf (ii)}  For $k \gt 1$ and $\mathbf{H} = k\mathrm{JetsOfSmoothGroupoids}_\infty$
  (Def. \ref{FormalSmoothInfinityGroupoids})
  the local integrability condition of Def. \ref{LocallyIntegrableGStructures}
  is of the form of the ``uniformly $k$-flatness''-condition of
  \cite{Guillemin65}.  But beware that, according to Def. \ref{GStructures} but in
  contrast to \cite{Guillemin65}, in this case the $G$-structure itself
  is not on the plain frame bundle but on the order-$k$ jet frame bundle
  (by Example \ref{FrameBundleOnSmoothManifold}).
\end{example}

\medskip

\noindent {\bf Haefliger groupoids.}
\begin{defn}[Haefliger groupoid]
  \label{HaefligerGroupoid}
  Let $\mathbf{H}$ be an elastic $\infty$-topos (Def. \ref{ElasticInfinityTopos})
  and $V \in \mathrm{Groups}(\mathbf{H})$ (Prop. \ref{LoopingAndDelooping}).

  \noindent {\bf (i)} With no further structure,
 \begin{itemize}
\vspace{-1.5mm}
\item[{\bf (a)}]
  The $V$-\emph{Haefliger groupoid} is the
  {\'e}tale groupoid (Def. \ref{EtaleGroupoids})
  \vspace{-1.5mm}
  $$
    \mathrm{Haef}_\bullet(V)
    \;\in\;
    \mbox{\'EtaleGroupoids}(\mathbf{H})
  $$

\vspace{-2.5mm}
\noindent
  which is the
  {\'e}talification (Def. \ref{EtalificationOfGroupoids})
  of the Atiyah groupoid (Def. \ref{TheAtiyahGroupoid})
  of the frame bundle (Def. \ref{StructureGroupOfVFolds})
  of $V$ regarded as a $V$-fold (Example \ref{VAsAVFold}):
  \vspace{-1.5mm}
  \begin{equation}
    \label{HaefligerGroupoidFromAtiyah}
    \mathrm{Haef}_\bullet(V)
    \;:=\;
    \mathrm{At}^{\mbox{\tiny{\'et}}}_\bullet
    \big(
      \mathrm{Frames}(V)
    \big).
  \end{equation}

\vspace{-4mm}
 \item[{\bf (b)}] The {\it $V$-Haefliger stack} of $V$ is
  the corresponding $V$-fold
  (according to Remark \ref{VFoldsAndVEtaleGroupoids}):
  \vspace{-1mm}
  \begin{equation}
    \label{HaefligerStack}
    \mathcal{H}\mathrm{aef}(V)
    \;:=\;
    \mathcal{A}\mathrm{t}^{\mbox{\tiny{\'et}}}
    \big(
      \mathrm{Frames}(V)
    \big)
    \;\in\;
    V\mathrm{Folds}
    \,.
  \end{equation}
  \end{itemize}

\vspace{-2.5mm}
  \noindent {\bf (ii)}
  Given, in addition,
  $(G,\phi) \in \mathrm{Groups}(\mathbf{H})_{/\mathrm{Aut}(T_e V)}$
  (Def. \ref{GStructureCoefficients}),
  with $G\mathrm{Frames}(V) \to V$ denoting the $G$-frame bundle
  \eqref{GFrameBundle} corresponding to the canonical $(G,\phi)$-structure
  on $V$ (Example \ref{CanonicalGStructure}), we say

 \begin{itemize}
  \vspace{-.15cm}
\item[{\bf (a)}] the
    $\big(V, (G,\phi) \big)$-\emph{Haefliger groupoid}
  is the
  {\'e}tale groupoid (Def. \ref{EtaleGroupoids})
  \vspace{-1.5mm}
  $$
    \mathrm{Haef}_\bullet
    \big(
      V, (G,\phi)
    \big)
    \;\in\;
    \mbox{\'EtaleGroupoids}(\mathbf{H})
  $$

  \vspace{-2mm}
\noindent
  which is the
  {\'e}talification (Def. \ref{EtalificationOfGroupoids})
  of the Atiyah groupoid (Def. \ref{TheAtiyahGroupoid})
  of the $G$-frame bundle \eqref{GFrameBundle}:
  \vspace{-1mm}
  \begin{equation}
    \label{GHaefligerGroupoidFromAtiyah}
    \mathrm{Haef}_\bullet
    \big(
      V, (G,\phi)
    \big)
    \;:=\;
    \mathrm{At}^{\mbox{\tiny{\'et}}}_\bullet
    \big(
      G\mathrm{Frames}(V)
    \big)
    \,.
  \end{equation}

\vspace{-4mm}
  \item[{\bf (b)}] The {\it $\big(V, (G,\phi) \big)$-Haefliger stack} of $V$ is
  the corresponding $V$-fold
  (according to Remark \ref{VFoldsAndVEtaleGroupoids}):
  \vspace{-1mm}
  \begin{equation}
    \label{GHaefligerStack}
    \mathcal{H}\mathrm{aef}\big(V, (G,\phi) \big)
    \;:=\;
    \mathcal{A}\mathrm{t}^{\mbox{\tiny{\'et}}}
    \big(
      G\mathrm{Frames}(V)
    \big)
    \;\in\;
    V\mathrm{Folds}
    \,.
  \end{equation}
   \end{itemize}
\end{defn}

\newpage

\begin{prop}[Haefliger stack represents $V$-fold structure]
  \label{VFoldStructureRepresentedByHaefliger}
  Let $\mathbf{H}$ be an elastic $\infty$-topos (Def. \ref{ElasticInfinityTopos})
  $V \in \mathrm{Groups}(\mathbf{H})$ (Prop. \ref{LoopingAndDelooping})
  and $X \in \mathbf{H}$.  Then the following are equivalent:

  {\bf (i)} $X$ is a $V$-fold (Def. \ref{VManifold});

  {\bf (ii)} $X$ admits a local diffeomorphism to the
    $V$-Haefliger stack (Def. \ref{HaefligerGroupoid}).
\end{prop}

\begin{proof}
  First consider the implication {\bf (i)} $\Rightarrow$ {\bf (ii)}:
Assuming $X$ is a $V$-fold, consider a $V$-atlas \eqref{VAtlas}
  $\xymatrix@C=15pt{V \ar@{<-}[r]^-{\mbox{\tiny{\'e}t}} & U \ar@{->>}[r]^-{\mbox{\tiny{\'e}t}} & X\;. }$
  By Prop. \ref{PullbackAlongLocalDiffeomorphismsPreservesTangentBundles}
  (and as in the proof of Prop. \ref{TangentBundleOfVFoldIsFiverBundle})
  the pullbacks of the frame bundles of $V$ and of $X$
  along this $V$-atlas to $U$ coincide there,
  which means that we have a homotopy-commutative square
  of their classifying maps \eqref{FrameBundleClassifyingMap}
  as shown on the bottom left of the following diagram:
  \vspace{-3mm}
  \begin{equation}
    \label{TowardsUniversalMapsToHaefliger}
    \raisebox{80pt}{
    \xymatrix{
      \ar@<-12pt>@{..>}[d]
      \ar@<-6pt>@{<..}[d]
      \ar@{..>}[d]
      \ar@<+6pt>@{<..}[d]
      \ar@<+12pt>@{..>}[d]
      &&
      \ar@<-12pt>@{..>}[d]
      \ar@<-6pt>@{<..}[d]
      \ar@{..>}[d]
      \ar@<+6pt>@{<..}[d]
      \ar@<+12pt>@{..>}[d]
      \\
      U \times_X U
      \ar@<-6pt>[d]_-{ \mbox{\tiny{\'e}t} }
      \ar@{<-}[d]|-{ \mbox{\tiny{\'e}t} }
      \ar@<+6pt>[d]^-{ \mbox{\tiny{\'e}t} }
      \ar@{-->}[rr]
      &&
      \mathrm{At}_1
      \big(
        \mathrm{Frames}(V)
      \big)
      \ar@<-6pt>[d]
      \ar@{<-}[d]
      \ar@<+6pt>[d]
      \\
      U
      \ar@{->>}[d]_{
        \mbox{\tiny{\'e}t}
      }
      \ar[rr]|-{\; \mbox{\tiny{\'e}t} \;}
      &&
      V
      \ar@{->>}[d]^-{\scalebox{.6}{$
        \vdash \mathrm{Frames}(V)
        $}
      }
            \\
      X
      \ar[rr]_-{\scalebox{.6}{$ \vdash \mathrm{Frames}(X)$} }
      &&
      \mathbf{B}\mathrm{Aut}(T_e V)
    }
    }
    \;\;\;\;\;\;\;\;\;\;\;\;\;
    \Leftrightarrow
    \;\;\;\;\;\;\;\;\;\;\;\;\;
    \raisebox{80pt}{
    \xymatrix{
      \ar@<-12pt>@{..>}[d]
      \ar@<-6pt>@{<..}[d]
      \ar@{..>}[d]
      \ar@<+6pt>@{<..}[d]
      \ar@<+12pt>@{..>}[d]
      &&
      \ar@<-12pt>@{..>}[d]
      \ar@<-6pt>@{<..}[d]
      \ar@{..>}[d]
      \ar@<+6pt>@{<..}[d]
      \ar@<+12pt>@{..>}[d]
      \\
      U \times_X U
      \ar@<-6pt>[d]_-{ \mbox{\tiny{\'e}t} }
      \ar@{<-}[d]|-{ \mbox{\tiny{\'e}t} }
      \ar@<+6pt>[d]^-{ \mbox{\tiny{\'e}t} }
      \ar[rr]|-{ \;\mbox{\tiny{\'e}t} \;}
      &&
      \mathrm{At}^{\mbox{\tiny{\'et}}}_1
      \big(
        \mathrm{Frames}(V)
      \big)
      \ar@<-6pt>[d]_-{ \mbox{\tiny{\'e}t} }
      \ar@{<-}[d]|-{ \mbox{\tiny{\'e}t} }
      \ar@<+6pt>[d]^-{ \mbox{\tiny{\'e}t} }
      \\
      U
      \ar@{->>}[d]_{
        \mbox{\tiny{\'e}t}
      }
      \ar[rr]|-{\; \mbox{\tiny{\'e}t}\; }
      &&
      V
      \ar@{->>}[d]^-{\mbox{\tiny{\'e}t }}
           \\
      X
      \ar@{-->}[rr]
      &&
      \mathcal{H}\!\mathrm{aef}(V)
    }
    }
  \end{equation}

  \vspace{-2mm}
\noindent  By passing to nerves (Example \ref{Nerve}) of the vertical morphisms,
  this induces a morphism of groupoids
  as shown on the top left.
    But $U_\bullet$ is an {\'e}tale groupoid (by Prop. \ref{EtaleGroupoidsAndEtaleAtlases}),
  and $U \longrightarrow V$ is a local diffeomorphism by
  definition of $V$-atlases,
  so that the top left part of the left diagram in \eqref{TowardsUniversalMapsToHaefliger}
  is in the {\'e}tale slice over $V$ (Def. \ref{EtaleTopos}).
  Therefore, the adjunction \eqref{AdjointsToInclusionOfEtaleSlice} of Prop. \ref{PropertiesOfEtaleToposes}
  implies that the top part of the diagram on the left of
  \eqref{TowardsUniversalMapsToHaefliger} factors
  through the {\'e}talification (Def. \ref{EtalificationOfGroupoids})
  as shown in the top part on the right.
  With this we get the dashed morphism on the right
  by passing to colimits over the vertical simplicial diagrams
  (as in Prop. \ref{EtaleGroupoidsAndEtaleAtlases}).

  It only remains to see that the dashed morphism on the right
  is itself a local diffeomorphism.
  For this observe that al the horizontal morphisms are
  local diffeomorphisms, using the assumptions and then left-cancellability
  (Lemma \ref{ClosureOfLocalDiffeomorphisms}). Therefore
  the statement follows with Lemma \ref{DegreewiseLocalDiffeomorphismsOfEtaleGroupoids}.

\vspace{1mm}
  For the converse implication {\bf (ii)} $\Rightarrow$ {\bf (i)}:
  Given a local diffeomorphism as shown dashed on the right of
  \eqref{TowardsUniversalMapsToHaefliger},
  we need to produce a $V$-atlas for $X$.
  So now define the bottom square on the right of \eqref{TowardsUniversalMapsToHaefliger}
  to be the pullback of the {\'e}tale atlas of the Haefliger stack
  along the griven morphism. This does
  make the top left span of the square a $V$-atlas by the fact that
  the classes of local diffeomorphisms and of effective epimorphisms
  are both closed under pullback
  (by Lemma \ref{EffectiveEpimorphismsArePreservedByPullback} and
  Lemma \ref{ClosureOfLocalDiffeomorphisms}).
\hfill \end{proof}

\begin{prop}[$G$-Structured Haefliger stack represents integrable $G$-structure]
  \label{GStructuredHaefligerGroupoidClassifiesIntegrability}
  Let $\mathbf{H}$ be an elastic $\infty$-topos  (Def. \ref{ElasticInfinityTopos}),
  $V \in \mathrm{Groups}(\mathbf{H})$ (Prop. \ref{LoopingAndDelooping}),
  $(G,\phi) \in \mathrm{Groups}(\mathbf{H})_{/\mathrm{Aut}(T_e V)}$
  (Def. \ref{GStructureCoefficients}).
  The
  $\big( V, (G,\phi) \big)$-Haefliger groupoid (Def. \ref{HaefligerGroupoid}),
  carries a canonical integrable $(G,\phi)$-structure
  (Def. \ref{IntegrableGStructure})
  \vspace{-2mm}
  \begin{equation}
    \label{CanonicalGStructureOnHaefliger}
    (\tau_{\mathcal{H}}, g_{\mathcal{H}})
    \;\in\;
    (G,\phi)\mathrm{Structures}_{\mathcal{H}\!\mathrm{Haef}(V)}(\mathbf{H})
  \end{equation}

  \vspace{-2mm}
\noindent
  such that the operation of pullback of \eqref{PullbackOfGStructures}
  along local diffeomorphism (Lemma \ref{GStructuresPullBackAlongLocalDiffeomorphisms})
  constitutes a natural bijection
  \begin{equation}
    \xymatrix@R=-4pt{
      \pi_0
      \;
      \mathrm{Integrably}(G,\phi)\mathrm{Structured}V\mathrm{Folds}(\mathbf{H})
      \ar@{}[rr]|-{\simeq}
      &&
      \pi_0
      \;
      \mbox{\rm\'Et}_{
        \mathcal{H}\!\mathrm{aef}
        (
          V, (G,\phi)
        )
      }
      \\
      \big(
        X, (\tau,g)
      \big)
      \ar@{}[rr]|-{ \longmapsto }
      &&
      \Big(
        X
        \xrightarrow{
          \vdash (\tau,g)
        }
        \mathcal{H}\!\mathrm{aef}
        \big(
          V, (G,\phi)
        \big)
      \Big)
    }
  \end{equation}

  \vspace{-2mm}
\noindent
  between the sets of equivalence classes of:

\vspace{1mm}
  {\bf( i)} integrably $(G,\phi)$-structured $V$-folds (Def. \ref{IntegrableGStructure}),

  {\bf (ii)} local diffeomorphisms into the
  $\big( V, (G,\phi)\big)$-Haefliger stack,
  hence objects in its {\'e}tale topos (Def. \ref{EtaleTopos}).
  \end{prop}
\begin{proof}
  We proceed as in the proof of Prop. \ref{VFoldStructureRepresentedByHaefliger},
  but lifting the diagram there from $\mathbf{H}$ to the
  iterated slice
  $\left(
    \mathbf{H}_{/\mathbf{B}\mathrm{Aut}(T_e V)}
   \right)_{/\mathbf{B}G}$ \eqref{IsometriesAsMorphismsInSlice}.

   \noindent {\bf (i)} First consider an integrably
   $G$-structured $V$-fold $\big( X, (\tau,g) \big)$.
   We describe the construction of a local diffeomorphism into the
   Haefliger stack from this:
      Pick any $\big( V, (\tau_V, g_V) \big)$-atlas
   $\!\xymatrix@C=11pt{
     \big(
       V, (\tau_V, g_V)
     \big)
     \ar@{<-}[r]^-{\mbox{\tiny met}}
     &
     \big(
       U, (\tau_U, g_U)
     \big)
     \ar@{->>}[r]^-{\mbox{\tiny met}}
     &
     \big(
       X, (\tau_X, g_X)
     \big)
   }$ \eqref{IsometryAtlas}.
   By Def. \ref{GStructures}, this is equivalently a
   choice of equivalence between the pullbacks to $U$ of the
   $G$-structures on $V$ and on $X$.
   Regarded in the iterated slice \eqref{IsometriesAsMorphismsInSlice},
   this equivalently means that we have a square in
    $\left(
    \mathbf{H}_{/\mathbf{B}\mathrm{Aut}(T_e V)}
   \right)_{/\mathbf{B}G}$ \eqref{IsometriesAsMorphismsInSlice},
  as shown on the left of the following:
  \vspace{-2mm}
  \begin{equation}
    \label{TowardsGStructuredHaefliger}
    \hspace{-1cm}
    \raisebox{80pt}{
    \xymatrix@R=8pt@C=-32pt{
      \ar@<-12pt>@{..>}[dd]
      \ar@<-6pt>@{<..}[dd]
      \ar@{..>}[dd]
      \ar@<+6pt>@{<..}[dd]
      \ar@<+12pt>@{..>}[dd]
      && &&
      \ar@<-12pt>@{..>}[dd]
      \ar@<-6pt>@{<..}[dd]
      \ar@{..>}[dd]
      \ar@<+6pt>@{<..}[dd]
      \ar@<+12pt>@{..>}[dd]
      \\
      &&
      {\phantom{AAAAAAAAAAAAAAAAAAAA}}
      \\
      U \times_X U
      \ar@<-6pt>[dd]_-{ \mbox{\tiny{\'e}t} }
      \ar@{<-}[dd]|-{ \mbox{\tiny{\'e}t} }
      \ar@<+6pt>[dd]^-{ \mbox{\tiny{\'e}t} }
      \ar@{-->}[rrrr]
      && &&
      \mathrm{At}_1
      \big(
        G\mathrm{Frames}(V)
      \big)
      \ar@<-6pt>[dd]
      \ar@{<-}[dd]
      \ar@<+6pt>[dd]
      \\
      \\
      U
      \ar@{->>}[dddd]_{
        \mbox{\tiny{\'e}t}
      }
      \ar[rrrr]|-{\; \mbox{\tiny{\'e}t} \;}
      && &&
      V
      \ar@{->>}[dddd]|-{ \tau_V }^>>>{\ }="t2"
      \ar@/^1pc/[dddddr]^-{\scalebox{.6}{$
        \vdash\mathrm{Frames}(V)
        $}
      }_-{\ }="s2"
      \\
      \\
      \\
      \\
      X
      \ar[rrrr]|-{ \tau_X }_>>>{\ }="t1"
      \ar@/_1pc/[drrrrr]_-{\scalebox{.6}{$
        \vdash \mathrm{Frames}(X)
        $}
      }^>>>>>>>>>{\ }="s1"
      && &&
      \mathbf{B}G
      \ar[dr]
      \\
      &&&& &
      \scalebox{.85}{$
        \mathbf{B}\mathrm{Aut}(T_e V)
      $}
      \ar@{=>}_{ g_X } "s1"; "t1"
      \ar@{=>}^{ g_V } "s2"; "t2"
    }
    }
        \Leftrightarrow
       \raisebox{80pt}{
    \xymatrix@R=8pt@C=3pt{
      \ar@<-12pt>@{..>}[dd]
      \ar@<-6pt>@{<..}[dd]
      \ar@{..>}[dd]
      \ar@<+6pt>@{<..}[dd]
      \ar@<+12pt>@{..>}[dd]
      && &&
      \ar@<-12pt>@{..>}[dd]
      \ar@<-6pt>@{<..}[dd]
      \ar@{..>}[dd]
      \ar@<+6pt>@{<..}[dd]
      \ar@<+12pt>@{..>}[dd]
      \\
      &&
      {\phantom{AAAAAAAAA}}
      \\
      U \times_X U
      \ar@<-6pt>[dd]_-{ \mbox{\tiny{\'e}t} }
      \ar@{<-}[dd]|-{ \mbox{\tiny{\'e}t} }
      \ar@<+6pt>[dd]^-{ \mbox{\tiny{\'e}t} }
      \ar[rrrr]|-{ \;\mbox{\tiny{\'e}t} \;}
      && &&
      \mathrm{At}^{\mbox{\tiny{\'et}}}_1
      \big(
        G\mathrm{Frames}(V)
      \big)
      \ar@<-6pt>[dd]_-{ \mbox{\tiny{\'e}t} }
      \ar@{<-}[dd]|-{ \mbox{\tiny{\'e}t} }
      \ar@<+6pt>[dd]^-{ \mbox{\tiny{\'e}t} }
      \\
      \\
      U
      \ar@{->>}[dddd]_{
        \mbox{\tiny{\'e}t}
      }
      \ar[rrrr]|-{\; \mbox{\tiny{\'e}t}\; }
      && &&
      V
      \ar@{->>}[dddd]_-{\mbox{\tiny{\'e}t }}
      \ar@/^1.4pc/[ddddddrr]^-{ \scalebox{.6}{$\vdash\mathrm{Frames}(V)$} }_-{\ }="s2"
      \ar[dddddr]|-{ \tau_V }^>>>>>>>>{\ }="t2"
      \\
      \\
      \\
      \\
      X
      \ar@{-->}[rrrr]^-{
        \vdash (\tau_X, g_X)
      }
      \ar@/_1.4pc/[ddrrrrrr]_-{\scalebox{.6}{$ \vdash \mathrm{Frames}(X)$} }^-{\ }="s1"
      \ar[drrrrr]|-{\tau_X}_>>>>>>>>>>>>>>>>>>{\ }="t1"
      && &&
      \scalebox{.9}{$
        \mathcal{H}\!\mathrm{aef}
        \big(
          V,(G,\phi)
        \big)
      $}
      \ar[dr]|<<<<{ \tau_{\mathcal{H}} }
      \ar@/_1.2pc/[ddrr]|<<<<{\phantom{AA}}|>>>>>{
        \scalebox{.5}{$
          \vdash\mathrm{Frames}(\mathcal{H})
        $}
      }
      \\
      && && &
      \scalebox{.85}{$
        \mathbf{B}G
        $}
      \ar[dr]^{\!\!\!\!\!
        \scalebox{.6}{$
          \mathbf{B}\phi
        $}
      }
      \\
      && && &&
      \scalebox{.85}{$
        \mathbf{B} \mathrm{Aut}(T_e V)
      $}
      \ar@{=>}_{g_X} "s1"; "t1"
      \ar@{=>}^{g_X} "s2"; "t2"
    }
    }
  \end{equation}

  \vspace{2mm}
 \noindent  Now we proceed as follows:

\noindent  {\bf (a)}
  Observing (with Prop. \ref{HomotopyFiberProductsInOverToposes})
  that
  fiber products in the iterated slice are actually given
  by the plain fiber products in $\mathbf{H}$
  equipped with canonical morphisms to the slicing objects,
  we find that
  passing to nerves (Example \ref{Nerve})
  of the vertical morphisms on the left of \eqref{TowardsGStructuredHaefliger}
  yields a morphism from the {\'e}tale groupoid induced by
  the given $V$-cover of $X$ to the Atiyah groupoid of
  $G\mathrm{Frames}(X)$ (Def. \ref{TheAtiyahGroupoid}) --
  just as in \eqref{TowardsUniversalMapsToHaefliger},
  but now equipped with coherent maps
  to $\mathbf{B}\phi$.

 \noindent  {\bf (b)} Therefore, we obtain the factorization through the
  $\big(V, (G,\phi)\big)$-Haefliger groupoid
  (the {\'e}talification of the Atiyah groupoid of the $G$-frame bundle
  shown on the top right of \eqref{TowardsGStructuredHaefliger})
  just as in \eqref{TowardsUniversalMapsToHaefliger},
  but now, in addition, coherently equipped with maps to $\mathbf{B}\phi$.

 \noindent  {\bf (c)} After this {\'e}talification we may identify these maps:
  Since those on $V$ remain unchanged by {\'e}talification over $V$,
  these still give the canonical $(G,\phi)$-structure $(\tau_V, g_V)$,
  as shown on the far right of \eqref{TowardsGStructuredHaefliger}.
  But since now the vertical simplicial morphisms are all local diffeomorphisms,
  pullback along which preserves $(G,\phi)$-structure
  (by Lemma \ref{GStructuresPullBackAlongLocalDiffeomorphisms}) and
  in particular preserves tangent- and frame bundles
  (by Prop. \ref{PullbackAlongLocalDiffeomorphismsPreservesTangentBundles})
  it follows that all stages of the $\big( V, (G,\phi) \big)$-Haefliger
  groupoid in the top right are now equipped with the classifying map
  of their frame bundles.

 \noindent  {\bf (d)} Since colimits in the slice are given by colimits
  in the underlying topos
  (by Example \ref{BaseChangeAlongTerminalMorphism}), the
  colimit over the simplicial sub-diagram on the far right of
  \eqref{TowardsGStructuredHaefliger} still yields
  the $\big( V, (\tau,g)\big)$-Haefliger stack \eqref{GHaefligerStack},
  as shown, now equipped with canonical maps to $\mathbf{B}\phi$.

\noindent  {\bf (e)} We claim that the induced map
  from the Haefliger stack to
  $\mathbf{B}\mathrm{Aut}(T_e V)$,
  denoted $\vdash \mathrm{Frames}(\mathcal{H})$ in
  \eqref{TowardsGStructuredHaefliger},
  is indeed the classifying map
  of the frame bundle of the Haefliger stack:
  \vspace{-2mm}
  \begin{equation}
    \vdash \mathrm{Frames}(\mathcal{H})
    \;\;\simeq\;\;\;
    \vdash \mathrm{Frames}\Big(
      \mathcal{H}\!\mathrm{aef}
      \big(
        V, (G,\phi)
      \big)
    \big).
  \end{equation}

  \vspace{-3mm}
\noindent
  This follows because:
  \begin{itemize}

  \vspace{-.3cm}
  \item
  by {\bf (c)} above, the component maps of the
  colimiting map classify the frame bundles
  of the stages of the simplicial nerve;

  \vspace{-.3cm}
  \item  therefore, the colimiting map classifies the
  colimit of the frame bundles of the simplicial nerve,
  by Prop. \ref{ColimitsOfClassifyingMapsAreClassifyingMapsOfColimits},

  \vspace{-.3cm}
  \item but the colimit of the tangent bundles of the
  {\'e}tale cover is the
  tangent bundle of the corresponding {\'e}tale stack, by
  Prop. \ref{TangentGroupoids}.
  \end{itemize}

\vspace{-3mm}
\noindent  {\bf (f)} In particular, this implies that the induced homotopy
  which fills the bottom right part of \eqref{TowardsGStructuredHaefliger}
  \vspace{-2mm}
  \begin{equation}
    \label{UniversalGStructureOnHaefligerStack}
    \xymatrix{
      \vdash \mathrm{Frames}
      \big(
        \mathcal{H}
      \big)
      \;
      \ar@{=>}[r]^-{ g_{\mathcal{H}} }
      &
      \;
      \mathbf{B} \phi\circ \tau_{\mathcal{H}}
      \,,
    }
  \end{equation}

  \vspace{-3mm}
\noindent
  canonically given by the colimit construction in the iterated slice,
  constitutes a $(G,\phi)$-structure on the
  $\big( V, (G,\phi))\big)$-Haefliger stack.

\noindent  {\bf (g)} In conclusion, the dashed morphism on the right of
  \eqref{TowardsGStructuredHaefliger} exists and is a local
  diffeomorphism, as in the proof of Prop. \ref{VFoldStructureRepresentedByHaefliger};
  but, by construction in the iterated slice, it is now
  exhibited as a local isometry
  to the Haefliger stack equipped with the induced $(G,\phi)$-structure
  \eqref{UniversalGStructureOnHaefligerStack}.

\vspace{1mm}
  \noindent {\bf (ii)}
  The converse construction is now immediate:
  Given a local diffeomorphism of the form shown dashed
  on the right of \eqref{TowardsGStructuredHaefliger},
  pulling back the {\'e}tale atlas of the Haefliger stack
  along it yields a $V$-atlas for $X$
  (just as in the proof of this converse step in Prop. \ref{VFoldStructureRepresentedByHaefliger})
  and
  pulling (via Lemma \ref{GStructuresPullBackAlongLocalDiffeomorphisms})
  the $(G,\phi)$-structure \eqref{UniversalGStructureOnHaefligerStack}
  around the resulting Cartesian square makes this
  a $\big(V, (G,\phi) \big)$-atlas that exhibits
  $X$ as equipped with an integrable $(G,\phi)$-structure.
  This construction is clearly injective on equivalence classes,
  by $\infty$-functoriality of the pullback construction \eqref{PullbackOfGStructures}
  of $(G,\phi)$-structures;
  and it is surjective on equivalence classes
  by item {(i)} above. Hence this is a bijection on equivalence classes, as claimed.
\hfill \end{proof}

\medskip

\noindent {\bf Tangential structures.}
Closely akin to $G$-structures (Def. \ref{GStructures})
are \emph{tangential structures} (Def. \ref{TangentialStructures} below)
where not the structure group itself is lifted, but only its shape:
\begin{defn}[Tangential structure]
  \label{TangentialStructures}
  Let $\mathbf{H}$ be an elastic $\infty$-topos  (Def. \ref{ElasticInfinityTopos}),
  $V \in \mathrm{Groups}(\mathbf{H})$ (Prop. \ref{LoopingAndDelooping}),
  $(G,\phi) \in \mathrm{Groups}(\mathbf{H})_{
    /_{
      \scalebox{.7}{$
        \raisebox{1pt}{\rm\textesh}
        \mathrm{Aut}(T_e V)
      $}
    }
  }$
  (Def. \ref{GStructureCoefficients})
  and $X \in V\mathrm{Folds}(\mathbf{H})$ (Def. \ref{VManifold}).

  \noindent {\bf (i)}
  We say that a \emph{tangential $(G,\tau)$-structure} on $X$
  is a lift $(\tau, g)$
  through $\mathbf{B}\phi$ of the
  composite of the frame bundle classifying
  map \eqref{FrameBundleClassifyingMap}
  with the shape-unit \eqref{AdjunctionUnit}:

   \vspace{-5mm}
  \begin{equation}
    \label{TangentialGphiStructure}
    \raisebox{20pt}{
    \xymatrix@C=3em@R=1.5em{
      &&&
      \mathbf{B}G
      \ar[d]^-{ \mathbf{B}\phi }
      \\
      \mathllap{
        \mbox{
          \tiny
          \color{darkblue}
          \bf
          $V$-fold
        }
        \;\;\;
      }
      X
      \ar[rr]_-{ \vdash\, \mathrm{Frames}(X) }
      \ar@{-->}[urrr]^-{
        \overset{
          \mathclap{
          \mbox{
            \tiny
            \color{darkblue}
            \bf
            \begin{tabular}{c}
              tangential
              \\
              structure
              \\
              \phantom{A}
            \end{tabular}
          }
          }
        }{
          \tau
        }
      }_>>>>>>>>>{\ }="s"
      &&
      \mathbf{B}
      \mathrm{Aut}(T_e V)
      \ar[r]_-{ \eta^{\scalebox{.6}{\textesh}} }^<{\ }="t"
      &
      \mathbf{B}
      \raisebox{1pt}{\rm\textesh}
      \mathrm{Aut}(T_e V)
      \mathrlap{
        \mbox{
          \tiny
          \color{darkblue}
          \bf
          \begin{tabular}{c}
            shape of
            \\
            structure group
            \\
            of frame bundle
          \end{tabular}
        }
      }
      \ar@{<=}^-{ g } "s"; "t"
    }
    }
  \end{equation}

 \noindent  {\bf (ii)}
  We write
  \begin{equation}
    \label{GroupoidOfTangentialGphiStructures}
    \mathrm{Tangential}(G,\phi)\mathrm{Structures}_X(\mathbf{H})
    \;:=\;
    \mathbf{H}_{/_{
      \scalebox{.5}{$
        \mathbf{B}
        \raisebox{1pt}{\textesh}
        \mathrm{Aut}(T_e V)
      $}
    }}
    \big(
      \eta^{\scalebox{.6}{\textesh}}
      \circ
      \vdash\mathrm{Frames}(X)
      \,,\,
      \mathbf{B}\phi
    \big)
  \end{equation}
  for the {\it $\infty$-groupoid of $(G,\phi)$-tangential structures}
  on the $V$-fold $X$.
\end{defn}

\begin{example}[Tangential structures on smooth manifolds]
  \label{TangentialStructureOnSmoothManifolds}
  Let $\mathbf{H} = \mathrm{JetsOfSmoothGroupoids}_\infty$
  (Example \ref{FormalSmoothInfinityGroupoids})
  $G \in \!
    \xymatrix{
      \mathrm{LieGroups}
      \longhookrightarrow
      \mathrm{Groups}
      \big(
        \mathbf{H}
      \big)
    }
    \!\!
  $ (see \eqref{ConcreteSmoothFormalInfinityGroupoidsAreDiffeologicalSpaces})
  and $X \in \mathrm{SmoothManifolds} \longhookrightarrow \mathbf{H}$
  regarded as an $\mathbb{R}^n$-fold according to Example \ref{OrdinaryManifolds}.
  In this case, the structure group of $X$ (Def. \ref{StructureGroupOfVFolds})
  is the ordinary general linear group $\mathrm{GL}_{\mathbb{R}}(n)$
  (Example \ref{OrdinaryGeneralLinearGroup}).
  Hence here tangential structure in the general sense of
  Def. \ref{TangentialStructures} is tangential structure in the
  traditional sense of differential topology
  (popularized under this name in \cite[5]{GMTW06},
  originally introduced as ``$(B,f)$-structure''
  \cite{Lashof63}\cite[II]{Stong68}, review in \cite[1.4]{Kochman96}).
\end{example}

\begin{example}[Cohesive refinement of tangential structure]
  \label{CohesiveRefinementOfTangentiaStructure}
  Every $(G,\phi)$-structure (Def. \ref{GStructures})
  induces tangential
  $(\raisebox{1pt}{\textesh}G, \raisebox{1pt}{\textesh}\phi)$-structure
  (Def. \ref{TangentialStructures}) by composition with the
  naturality square of $\eta^{\scalebox{.6}{\textesh}}$ on
  $\mathbf{B}\phi$:
  \vspace{-2mm}
  \begin{equation}
    \label{TangentialGphiStructureFromGphiStructure}
    \raisebox{20pt}{
    \xymatrix@C=4.5em@R=1.5em{
      &
      \mathbf{B}G
      \ar[d]^-{ \mathbf{B}\phi }
      \ar[r]^-{
        \eta^{\scalebox{.5}{\textesh}}_{ \mathbf{B}G }
      }
      &
      \mathbf{B}\raisebox{1pt}{\textesh}G
      \ar[d]^-{
        \mathbf{B}\raisebox{0pt}{\textesh}\phi
      }
      \\
      \mathllap{
        \mbox{
          \tiny
          \color{darkblue}
          \bf
          $V$-fold
        }
        \;\;\;
      }
      X
      \ar[r]_-{ \vdash\, \mathrm{Frames}(X) }^>>>{\ }="t"
      \ar@{-->}[ur]^-{
        \overset{
          \mathclap{
          \mbox{
            \tiny
            \color{darkblue}
            \bf
            \begin{tabular}{c}
              $(G,\phi)$-structure
              \\
              \phantom{A}
            \end{tabular}
          }
          }
        }{
          \tau
        }
      }_>>>>>>>>>{\ }="s"
      &
      \mathbf{B}
      \mathrm{Aut}(T_e V)
      \ar[r]_-{
        \eta^{\scalebox{.6}{\textesh}}_{
          \mathbf{B} \mathrm{Aut}(T_e V)
        }
      }
      &
      \mathbf{B}
      \raisebox{1pt}{\rm\textesh}
      \mathrm{Aut}(T_e V)
      \mathrlap{
        \mbox{
          \tiny
          \color{darkblue}
          \bf
          \begin{tabular}{c}
            shape of
            \\
            structure group
            \\
            of frame bundle
          \end{tabular}
        }
      }
      \ar@{<=}^-{ g } "s"; "t"
    }
    }
  \end{equation}

  \vspace{-2mm}
\noindent
Conversely, realizing a tangent structure as obtained
  from a $G$-structure this way means to find a
  geometric (differential) refinement.
\end{example}
\begin{example}[Orientation structure]
  \label{OrientationStructure}
  Let $\mathbf{H} = \mathrm{JetsOfSmoothGroupoids}_\infty$
  (Example \ref{FormalSmoothInfinityGroupoids})
  and $X \in \mathbf{H}$ an $\mathbb{R}^n$-fold (Def. \ref{VManifold})
  hence an ordinary manifold (Example \ref{OrdinaryManifolds})
  or, more generally, an ordinary {\'e}tale Lie groupoid (Example \ref{EtaleLieGroupoidAsRnFold}).
  With the general linear and the
  (special) orthogonal group regarded as smooth groups via \eqref{ConcreteSmoothFormalInfinityGroupoidsAreDiffeologicalSpaces}
   \vspace{-3mm}
  \begin{equation}
    \label{OrthogonalGroupsInGeneralLinearGroup}
    \xymatrix{
      \mathrm{SO}(n)
        \ar[r]^-{ i_{\mathrm{SO}} }
        &
      \mathrm{O}(n)
        \ar[r]^-{ i_O  }
        &
      \mathrm{GL}(n)
    }
    \;\in\;
    \xymatrix{
      \mathrm{Groups}(\mathrm{SmoothManifolds})
      \ar[r]
      &
      \mathrm{Groups}(\mathbf{H})
    }
  \end{equation}

   \vspace{-3mm}
\noindent
  we have:
  \begin{itemize}
    \vspace{-.2cm}
    \item[{\bf (i)}] an $\mathrm{O}(n)$-structure (Def. \ref{GStructures})
    on $X$
    is equivalently a Riemannian structure (Example \ref{TorsionFreeGStructureOnSmoothManifolds});
    \vspace{-.2cm}
    \item [{\bf (ii)}] but a
    tangential $\raisebox{1pt}{\textesh} \mathrm{O}(n)$-structure
    (Def. \ref{TangentialGphiStructure}) is equivalently
    \emph{no structure}, since
    $
      \xymatrix{\!
        \raisebox{1pt}{\textesh} \mathrm{O}(n)
        \ar[r]^-{ \scalebox{.7}{$
          \raisebox{1pt}{\textesh}
          i_{\mathrm{O}}
        $} }_-{ \simeq }
        &
        \raisebox{1pt}{\textesh} \mathrm{GL}(n)
      \!}
    $
    is an equivalence of underlying shapes
    (since $\mathrm{O}(n)$ is the maximal compact subgroup
    of $\mathrm{GL}(n)$),
    \vspace{-.2cm}
    \item[{\bf (iii)}] while a tangential
      $\raisebox{1pt}{\textesh} \mathrm{SO}(n)$-structure
      (Def. \ref{TangentialGphiStructure})
      is an \emph{orientation} of $X$.

    \vspace{-.2cm}
    \item[{\bf (iv)}]
      A differential refinement, in the sense of Example
      \ref{CohesiveRefinementOfTangentiaStructure}, of
      such an orientation structure is an oriented Riemannian
      structure (via its induced volume form).
\end{itemize}
\end{example}

\begin{example}[Higher Spin structure {\cite{SSS08}\cite{SSS09}}]
  \label{HigherSpinStructure}
  Let $\mathbf{H} = \mathrm{JetsOfSmoothGroupoids}_\infty$
  (Example \ref{FormalSmoothInfinityGroupoids})
  and $X \in \mathbf{H}$ an $\mathbb{R}^n$-fold (Def. \ref{VManifold})
  hence an ordinary manifold (Example \ref{OrdinaryManifolds})
  or, more generally, an ordinary {\'e}tale Lie groupoid (Example \ref{EtaleLieGroupoidAsRnFold}).
  The sequence of groups \eqref{OrthogonalGroupsInGeneralLinearGroup}
  in Example \ref{OrientationStructure} is, under shape,
  the beginning of the \emph{Whitehead tower}
  of $\raisebox{1pt}{\rm\textesh} O(n) \simeq \raisebox{1pt}{\rm\textesh} \mathrm{GL}(n)$.
  The tangential structures (Def. \ref{TangentialStructures},
  Example \ref{TangentialStructureOnSmoothManifolds})
  corresponding to the stages in this tower are the \emph{Spin structure}
  and its higher analogues:
   \vspace{-4mm}
  \begin{equation}
    \label{HigherSpinStructures}
    \hspace{-1cm}
    \raisebox{80pt}{
    \xymatrix@C=5em@R=1.8em{
      &&& &
      \ar@{..>}[d]
      \\
      &&& &
      \mathbf{B} \raisebox{1pt}{\textesh} \mathrm{Fivebrane}(n)
      \ar[d]
      \\
      &&& &
      \mathbf{B} \raisebox{1pt}{\textesh} \mathrm{String}(n)
      \ar[d]
      \\
      &&& &
      \mathbf{B} \raisebox{1pt}{\textesh} \mathrm{Spin}(n)
      \ar[d]
      \\
      &&& &
      \mathbf{B} \raisebox{1pt}{\textesh} \mathrm{SO}(n)
      \ar[d]
      \\
      X
      \ar[drrr]_-{
      \vdash\mathrm{Frames}(X)  \phantom{AAAAA}
       }
      \ar[rrr]|-{
        \mbox{
          \tiny
          \color{darkblue}
          \bf
          \begin{tabular}{c}
            Riemannian
            \\
            structure
          \end{tabular}
        }
      }
      \ar@/^.2pc/[urrrr]|-{
        \mbox{
          \tiny
          \color{darkblue}
          \bf
          \begin{tabular}{c}
          Orientation
          \\ structure
          \end{tabular}
        }
      }
      \ar@/^.4pc/[uurrrr]|-{
        \mbox{
          \tiny
          \color{darkblue}
          \bf
          \begin{tabular}{c}
          Spin \\ structure
          \end{tabular}
        }
      }
      \ar@/^.6pc/[uuurrrr]|-{
        \mbox{
          \tiny
          \color{darkblue}
          \bf
          \begin{tabular}{c}
          String \\ structure
          \end{tabular}
        }
      }
      \ar@/^.8pc/[uuuurrrr]|-{
        \mbox{
          \tiny
          \color{darkblue}
          \bf
          \begin{tabular}{c}
          Fivebrane \\ structure
          \end{tabular}
        }
      }
      &&&
      \mathbf{B} \mathrm{O}(n)
      \ar[r]|-{\;
        \eta^{\scalebox{.6}{\textesh}}_{\mathbf{B}\mathrm{O}(n)}
      }
      \ar[d]
      &
      \mathbf{B} \raisebox{1pt}{\textesh} \mathrm{O}(n)
      \ar[d]^-{\simeq}
      \\
      &&&
      \mathbf{B}\mathrm{GL}(n)
      \ar[r]|-{\;
        \eta^{\scalebox{.6}{\textesh}}_{\mathbf{B}\mathrm{GL}(n)}
      }
      &
      \mathbf{B} \raisebox{1pt}{\textesh} \mathrm{GL}(n)
    }
    }
  \end{equation}
\end{example}

\medskip
\medskip

\noindent {\bf Flat $V$-folds.}
\begin{defn}[Flat $V$-folds]
\label{FlatVFolds}
Let $\mathbf{H}$ be an elastic $\infty$-topos (Def. \ref{ElasticInfinityTopos}),
$V \in \mathrm{Groups}(\mathbf{H})$ (Prop. \ref{LoopingAndDelooping})
and $X \in V\mathrm{Folds}(\mathbf{H})$ (Def. \ref{VManifold}).
We say that $X$ is \emph{flat} if
the classifying map \eqref{FrameBundleClassifyingMap}
of its frame bundle (Prop. \ref{TangentBundleOfVFoldIsFiverBundle})
factors through the $\flat$-counit \eqref{CounitOfAdjunction}, hence if it carries
$(G,\phi)$-structure (Def. \ref{GStructures})
for $(G,\phi) = (\flat \mathrm{Aut}(T_e V), \epsilon^\flat_{\mathrm{Aut}(T_e V)})$:
 \vspace{-2mm}
\begin{equation}
  \label{FlatFrameBundleClassifyingMap}
  \xymatrix@C=30pt{
    &&
    \flat
    \mathbf{B} \mathrm{Aut}(T_e V)
    \ar[d]^-{ \epsilon^\flat_{\mathbf{B}\mathrm{Aut}(T_e V)} }
    \\
    X
    \ar[rr]_-{
      \vdash\mathrm{Frames}(X)
    }^-{\ }="s"
    \ar[urr]^-{ \tau }_>>>>>>>{\ }="t"
    &
    &
    \mathbf{B} \mathrm{Aut}(T_e V)
    \ar@{=>} "s"; "t"
  }
\end{equation}

 \vspace{-3mm}
\noindent
By the universal property of $\epsilon^\flat$ and since $\flat$
commutes with $\mathbf{B}$,
this means equivalently that $X$ carries $G$-structure for
any discrete group $G \simeq \flat \flat G$.
\end{defn}

\newpage

\begin{prop}[Flat frame bundles are $V$-folds]
 \label{FlatFrameBundlesAreVFolds}
Let $\mathbf{H}$ be an elastic $\infty$-topos (Def. \ref{ElasticInfinityTopos}),
$V \in \mathrm{Groups}(\mathbf{H})$ (Prop. \ref{LoopingAndDelooping})
and $X \in V\mathrm{Folds}(\mathbf{H})$ (Def. \ref{VManifold}).
If $X$ is flat (Def. \ref{FlatVFolds}), then
 \vspace{-1mm}
\item {\bf (i)}
its flat frame bundle $(\flat \mathrm{Aut}(T_eV))\mathrm{Frames}(X)$
\eqref{GFrameBundle} is itself a $V$-fold (Def. \ref{VManifold})
and
 \vspace{-2mm}
\item {\bf (ii)} the bundle morphism
is a local diffeomorphism (Def. \ref{FormallyEtaleMorphism}):
$\xymatrix{
  (\flat \mathrm{Aut}(T_eV))\mathrm{Frames}(X) \ar[r]^-{\mbox{\tiny\rm{\'e}t}} &  X
  }$.
\end{prop}
\begin{proof}
  First consider {\bf (ii)}:
  We need to show that the left square in the following pasting
  diagram is Cartesian:
   \vspace{-2mm}
  $$
    \xymatrix@C=3em@R=1.8em{
      (\flat \mathrm{Aut}(T_eV))\mathrm{Frames}(X)
      \ar[rr]^{
        \!\!\!\!\!\!\!\!\!\!
        \eta^\Im_{(\flat \mathrm{Aut}(T_eV))\mathrm{Frames}(X)}
      }
      \ar[d]_-p
      &&
      \Im
      \big(
        (\flat \mathrm{Aut}(T_eV))\mathrm{Frames}(X)
      \big)
      \ar[r]
      \ar[d]_-{ \Im p }
      \ar@{}[dr]|-{ \mbox{\tiny(pb)} }
      &
      \Im \ast
      \ar[d]
      \\
      X
      \ar[rr]_-{ \eta^\Im_X }
      &&
      \Im X
      \ar[r]_{ \Im \tau }
      &
      \Im \flat \mathrm{Aut}(T_e V)
    }
  $$

   \vspace{-2mm}
\noindent
  Here the right square is Cartesian, by
  definition \eqref{GFrameBundle} and since $\Im$, being
  a right adjoint, preserves Cartesian squares (by Prop. \ref{AdjointsPreserveCoLimits}).
  Hence, by the pasting law (Prop. \ref{PastingLaw})
  it is sufficient to show that the total rectangle is Cartesian.
  But, by the naturality of $\eta^\Im$, the total rectangle is equivalent to
  that of  the following pasting diagram:
   \vspace{-2mm}
  $$
    \xymatrix{
      (\flat \mathrm{Aut}(T_eV))\mathrm{Frames}(X)
      \ar[d]
      \ar[r]
      \ar@{}[dr]|-{ \mbox{\tiny(pb)} }
      &
      \ast
      \ar[d]
      \ar[rr]^-{ \eta^\Im_\ast }
      \ar@{}[drr]|-{ \mbox{\tiny(pb)} }
      &&
      \Im \ast
      \ar[d]
      \\
      X
      \ar[r]_-{\tau}
      &
      \flat \mathbf{B} \mathrm{Aut}(T_e V)
      \ar[rr]_{ \eta^\Im_{\flat \mathbf{B} \mathrm{Aut}(T_e V)} }
      &&
      \Im \flat \mathbf{B} \mathrm{Aut}(T_e V)
    }
  $$

   \vspace{-2mm}
\noindent
  Here the left square is Cartesian by the definition
  \eqref{GFrameBundle}, while the right square is Cartesian
  since its two horizontal morphisms are equivalences, by
  elasticity. Hence the total rectangle is Cartesian
  by the pasting law (Prop. \ref{PastingLaw}).

\noindent   Regarding {\bf (i)}: We need to exhibit a $V$-atlas \eqref{VAtlas}  for
  the flat frame bundle. So let
    $
    \xymatrix@R=1.5em{
      V
      \ar@{<-}[r]^-{ \mbox{\tiny{\'e}t} }
      &
      U
      \ar@{->>}[r]^-{ \mbox{\tiny{\'e}t} }
      &
      X
    }
  $
  be a $V$-atlas for $X$, and consider the following
  pullback diagram:
 \vspace{-2mm}
  $$
    \xymatrix@C=4em{
      U
        \underset{X}{\times}
      (\flat\mathrm{Aut}(T_e V))\mathrm{Frames}(X)
      \ar@{->>}[d]_{ \mbox{\tiny{\'e}t} }
      \ar@{->>}[r]^{ \mbox{\tiny{\'e}t} }
      \ar@{}[dr]|<<<<<<{\mbox{\tiny(pb)}}
      &
      (\flat\mathrm{Aut}(T_e V))\mathrm{Frames}(X)
      \ar@{->>}[d]^-{\mbox{\tiny{\'e}t}}
      \\
      U
      \ar[d]_-{ \mbox{\tiny{\'e}t} }
      \ar@{->>}[r]_-{ \mbox{\tiny{\'e}t} }
      &
      X
      \\
      V
    }
  $$

   \vspace{-2mm}
\noindent
  Observe that all four morphisms in the square
  are effective epimorphisms (Def. \ref{EffectiveEpimorphisms})
  and local diffeomorphisms (Def. \ref{FormallyEtaleMorphism}):
  The bottom one by definition,
  the right one by {\bf (ii)} and hence the other two since
  both classes of morphisms are closed under pullback
  (Lemma \ref{EffectiveEpimorphismsArePreservedByPullback} and
  Lemma \ref{ClosureOfLocalDiffeomorphisms}).
  Finally, since the class of local diffeomorphisms is also closed under
  composition (Lemma \ref{ClosureOfLocalDiffeomorphisms}),
  the total vertical morphisms is a local
  diffeomorphism, and hence the total outer diagram is a $V$-atlas
  of the flat frame bundle.
\hfill \end{proof}

\begin{prop}[$\flat G$-frame bundles are $V$-folds]
\label{FlatGFrameBundlesAreVFolds}
Let $\mathbf{H}$ be an elastic $\infty$-topos (Def. \ref{ElasticInfinityTopos}),
$V \in \mathrm{Groups}(\mathbf{H})$ (Prop. \ref{LoopingAndDelooping})
$X \in V\mathrm{Folds}(\mathbf{H})$ (Def. \ref{VManifold}),
$(G,\phi) \in \mathrm{Groups}(\mathbf{H})_{/\mathrm{Aut}(T_e V)}$
    (Prop. \ref{LoopingAndDelooping},  Def. \ref{AutomorphismGroup},
    Example \ref{LocalNeighbourhoodOfAPoint})
with $G \simeq \flat G$ discrete,
and $(\tau,g) \in (G,\phi)\mathrm{Structures}_X(\mathbf{H})$.
Then the corresponding $G$-frame bundle \eqref{GFrameBundle} is itself a $V$-fold:
 \vspace{-2mm}
$$
  G \;\simeq\; \flat G
  \phantom{AAA}
  \Leftrightarrow
  \phantom{AAA}
  G \mathrm{Frames}(X)
  \;\in\;
  V\mathrm{Folds}(\mathbf{H})
  \,.
$$
\end{prop}
\begin{proof}
  The proof  proceeds verbatim as that for
  Prop. \ref{FlatFrameBundlesAreVFolds},
  just with the structure group restricted along
  $\flat G \to \flat \mathrm{Aut}(T_e V)$.
\hfill \end{proof}

In summary, we have found the general abstract version of the local
model spaces of orbifolds:

\begin{prop}[Local orbifold model spaces]
  \label{LocalOrbifoldModelSpaces}
 Let $\mathbf{H}$ be an elastic $\infty$-topos (Def. \ref{ElasticInfinityTopos}),
 $G, V \in \mathrm{Groups}(\mathbf{H})$ (Prop. \ref{LoopingAndDelooping}),
 with $G \simeq \flat G$ discrete, and
 $(V, \rho) \in G\mathrm{Actions}(\mathbf{H})$ (Prop. \ref{InfinityAction})
 a restriction (Prop. \ref{PullbackAction}) of the action
 $(V,\rho_{\mathrm{Aut}})$ by group-automorphisms
  (Prop. \ref{CanonicalActionOfGroupAutomomorphisms}).
Then the homotopy quotient
 \eqref{HomotopyQuotientAsColimit}
 \vspace{-2mm}
 $$
   V \!\sslash\! G \;\in \; \mathbf{H}
 $$

 \vspace{-2mm}
\noindent
 of $V$ regarded with its canonical framing (Prop. \ref{GroupsAdmitFraming})
 \begin{itemize}
   \vspace{-.2cm}
   \item[\bf (i)]
     is a flat $V$-fold (Def. \ref{FlatVFolds});
   \vspace{-.2cm}
   \item[\bf (ii)]
     with $G$-structure (Def. \ref{GStructures})\;
   \vspace{-.2cm}
   \item[\bf (iii)]
     whose $G$-frame bundle \eqref{GFrameBundle}
     is $G$-equivariantly (Def. \ref{Equivariance}) equivalent to $V$ itself:
     \vspace{-2mm}
     $$
       G \mathrm{Frames}
       \big(
         V \!\sslash\! G
       \big)
       \;\simeq\;
       V\;.
     $$
 \end{itemize}
\end{prop}
\begin{proof}
 First observe that $V \!\sslash\! G$
  is a $V$-fold, by Prop. \ref{OrbifoldingOfVFoldIsAVFold}
  applied to Example \ref{VAsAVFold}.
  That this is flat {\bf (i)} is implied by {\bf (ii)}, since $G$ is assumed
  to be discrete.
    For {\bf (ii)} and {\bf (iii)} observe that
  the canonical framing on $V$
  is $G$-equivariant, by Prop. \ref{CanonicalFramingIsEquivariantUnderGroupAutomorphisms},
  so that Prop. \ref{FactoringVramesOfQuotientOfFramed} implies
  $G$-structure on $V \!\sslash\! G$
  classified by the action morphism $\rho$ itself.
  But this means that its homotopy fiber, hence the
  corresponding $G$-frame bundle (Def. \ref{GFrameBundle})
  is $V$ itself, by \eqref{InfinityActionHomotopyFiberSequence}
  (and in accord with Prop. \ref{FlatGFrameBundlesAreVFolds}).
\hfill \end{proof}

\begin{example}[Ordinary orbifold singularities]
  \label{OrdinaryOrbifoldSingularities}
  Let $\mathbf{H} := \mathrm{JetsOfSmoothGroupoids}_\infty$
  (Example \ref{FormalSmoothInfinityGroupoids})
  and $V := (\mathbb{R}^n, +)$ as in Example \ref{OrdinaryManifolds}.
  Then a group automorphism of $V$ is a \emph{linear} isomorphism,
  hence $\mathrm{Aut}_{\mathrm{Grp}}(\mathbb{R}^n,+) \simeq
  \mathrm{GL}(n)$.
  Therefore, in this case
  the assumptions of Prop. \ref{LocalOrbifoldModelSpaces}
  hold precisely for $V$ a \emph{linear representation} of
  the discrete group $G$, and thus we recover the
  traditional local orbifold models $V \!\sslash\! G$
  from \cite{Satake56} (in their incarnation as {\'e}tale groupoids).
\end{example}

\medskip

\noindent{\bf Orbi-$V$-folds.}
Finally, we promote $V$-folds to orbifolds proper,
in that we promote the $\infty$-category of {\'e}tale stacks to a
proper $\infty$-category of higher orbifolds:
\begin{defn}[Orbi-$V$-folds]
  \label{OrbiVFolds}
  Let $\mathbf{H}$ be a singular-elastic $\infty$-topos
  (Def. \ref{SingularSolidInfinityTopos})
  and $V \in \mathrm{Groups}(\mathbf{H}_{\tiny\smooth})$.
  We say that an \emph{orbi-$V$-fold} is an object
  $\mathcal{X} \in \mathbf{H}$ whose purely smooth aspect
  \eqref{CohesionSingular} is a $V$-fold (Def. \ref{VManifold}).

  \noindent {\bf (i)} We
  write $V\mathrm{Orbifolds}(\mathbf{H}) \subset \mathbf{H}$
  for the full sub-$\infty$-category on orbi-$V$-folds:
  \vspace{-2mm}
  $$
    \mathcal{X}
    \;\in\;
    V\mathrm{Orbifolds}(\mathbf{H})
    \phantom{AAA}
    \Leftrightarrow
    \phantom{AAA}
    \smooth \mathcal{X}
    \;\in\;
    V\mathrm{Folds}(\mathbf{H})
    \,.
  $$

   \vspace{-2mm}
\noindent
 This means, equivalently,
  that the orbi-$V$-folds in $\mathbf{H}$
  are the orbi-singularizations \eqref{CohesionSingular}
  of the $V$-folds in $\mathbf{H}_{\tiny\smooth}$:
  \begin{equation}
    \label{VFoldsInHSmoothAreVOrbifoldsInH}
    \xymatrix@R=7pt{
      V \mathrm{Folds}(\mathbf{H}_{\tiny\smooth})
      \;
      \ar@<+10pt>@{<-}[rr]|-{\;\mathrm{Smth}\;}
      \ar@<-10pt>@{->}[rr]|-{\;
        \mathrm{OrbSnglr}
      \;}^-{ \raisebox{5pt}{$\simeq$} }
      &&
      \;
      V\mathrm{Orbifolds}(\mathbf{H})
      \\
      \\
      \mathrm{Smth}(\mathcal{X})
      \ar@{<-|}[rr]|-{}
      \ar@{}[d]|-{
        \;\;\;
        \raisebox{1pt}{
        \begin{rotate}{+90}
          $
          \mathclap{
            :=
          }
          $
        \end{rotate}
        }
      }
      &&
      \mathcal{X}
      \ar@{}[d]|-{
        \!\!\!\!\!
        \raisebox{-1pt}{
        \begin{rotate}{-90}
          $
          \mathclap{
            :=
          }
          $
        \end{rotate}
        }
      }
      \\
      X
      \ar@{|->}[rr]
      && \mathrm{OrbSnglr}(X)
    }
  \end{equation}

  \noindent
  {\bf (ii)} Similarly, given, in addition,
  $(G,\phi) \in \mathrm{Groups}(\mathbf{H})_{/\mathrm{Aut}(T_e V)}$
  (Def. \ref{GStructureCoefficients}),  we write
  $(G,\phi)\mathrm{Structured}V\mathrm{Orbifolds}(\mathbf{H}) \subset \mathbf{H}$
  for the full sub-$\infty$-category on
  $(G,\phi)$-structured orbi-$V$-folds (Def. \ref{Isometries}):
  \vspace{-1mm}
  \begin{equation}
    \label{VFoldsInHSmoothAreVOrbifoldsInH}
    \xymatrix@R=5pt{
      (G,\phi)\mathrm{Structured}V\mathrm{Folds}(\mathbf{H}_{\tiny\smooth})
      \;
      \ar@<+10pt>@{<-}[rr]|-{\;\mathrm{Smth}\;}
      \ar@<-10pt>@{->}[rr]|-{\;
        \mathrm{OrbSnglr}
      \;}^-{ \raisebox{5pt}{$\simeq$} }
      &&
      \;
      (G,\phi)\mathrm{Structured}V\mathrm{Orbifolds}(\mathbf{H})
      \\
      \\
      (\mathrm{Smth}(\mathcal{X}),(\tau,g))
      \ar@{<-|}[rr]|-{}
      \ar@{}[d]|-{
        \;\;\;
        \raisebox{1pt}{
        \begin{rotate}{+90}
          $
          \mathclap{
            :=
          }
          $
        \end{rotate}
        }
      }
      &&
      (\mathcal{X}, (\tau,g))
      \ar@{}[d]|-{
        \!\!\!\!\!
        \raisebox{-1pt}{
        \begin{rotate}{-90}
          $
          \mathclap{
            :=
          }
          $
        \end{rotate}
        }
      }
      \\
      (X,(\tau,g))
      \ar@{|->}[rr]
      &&
      (\mathrm{OrbSnglr}(X),(\tau,g))
    }
  \end{equation}
\end{defn}

\begin{remark}[Coefficients for orbifold cohomology]
  \label{MoreCoefficients}
  The point of Def. \ref{OrbiVFolds} is that, by regarding a $V$-fold
  in the elastic $\infty$-topos $\mathbf{H}_{\tiny\smooth}$
  equivalently as an orbi-$V$-fold in the larger
  singular-elastic $\infty$-topos $\mathbf{H}$,
  a larger class of coefficients
  for intrinsic cohomology theories \eqref{IntrinsicCohomologyOfAnInfinityTopos}
  becomes available, notably coefficients
  of the form $\raisebox{1pt}{\textesh} \orbisingular (A \!\sslash\! G)$
  (see Lemma \ref{ShapeOfOrbiSingularSpacesAsPresheafOnSingularities} below).
  This is what gives rise, in \cref{OrbifoldCohomology},
  to proper orbifold cohomology (Def. \ref{OrbifoldCohomologyTangentiallyTwisted} below)
  in contrast to the coarser cohomology of
  underlying {\'e}tale groupoids (Def. \ref{EtaleCohomologyOfVFolds} below).
\end{remark}

\newpage

\begin{remark}[The proper $\infty$-category of higher orbifolds]
  \label{Proper}
  While \eqref{VFoldsInHSmoothAreVOrbifoldsInH} is
  an equivalence of abstract $\infty$-categories,

\vspace{-1mm}
\item {\bf (i)}   it is not an equivalence of sub-$\infty$-categories
  of the ambient singular-elastic $\infty$-topos $\mathbf{H}$:
  \vspace{-3mm}
  $$
    \xymatrix@R=8pt{
      \overset{
        \mathclap{
        \mbox{
          \tiny
          \color{darkblue}
          \bf
          \begin{tabular}{c}
            $\infty$-category of
            \\
            of {\'e}tale groupoids
            \\
            \phantom{a}
          \end{tabular}
        }
        }
      }{
        V \mathrm{Folds}(\mathbf{H}_{\tiny\smooth})
      }
      \ar@{_{(}->}[dr]_-{ \mathrm{Smth} }
      &
      \not\simeq
      &
      \overset{
        \mbox{
          \tiny
          \color{darkblue}
          \bf
          \begin{tabular}{c}
            proper
            \\
            $\infty$-category
            \\
            of orbifolds
            \\
            \phantom{a}
          \end{tabular}
        }
      }{
        V\mathrm{Orbifolds}(\mathbf{H})
      }
      \ar@{^{(}->}[dl]^-{ \mathrm{OrbSnglr} }
      &
      \mathrlap{
        \;\;
        \in
        \;
        (\mathrm{Categories}_\infty)_{/\mathbf{H}}
      }
      \\
      &
      \mathbf{H}
    }
  $$

   \vspace{-4mm}
\noindent
 \item {\bf (ii)}  To bring out this distinction, also in view of
  Remark \ref{MoreCoefficients}, we call
  $V\mathrm{Orbifolds}(\mathbf{H})$ (Def. \ref{OrbiVFolds})
  the
  \emph{proper $\infty$-category of orbifolds},
  in contrast to the  $\infty$-category
  $V\mathrm{Folds}(\mathbf{H}_{\tiny\smooth})$ \eqref{VFoldsCategory}
  of {\'e}tale $\infty$-groupoids.

   \vspace{-1mm}
 \item {\bf (iii)} It is a happy coincidence that \emph{proper}
  is also the technical adjective chosen in \cite{DHLPS19}
  for equivariant homotopy theories presented
  by $\infty$-presheaves over categories of
  orbits with compact -- hence finite if discrete --
  isotropy groups: In this terminology the singular-cohesive
  $\infty$-topos $\mathbf{H}$
  is, according to Def. \ref{SingularCohesiveInfinityTopos}, indeed
  a \emph{proper} global equivariant homotopy theory.
\end{remark}

\begin{example}[$\mathcal{V}$-folds]
  Let $\mathbf{H}$ be a singular-elastic $\infty$-topos (Def. \ref{SingularSolidInfinityTopos})
  and $\mathcal{V} \in \mathrm{Groups}(\mathbf{H})$ (Prop. \ref{LoopingAndDelooping})
  any group object, not necessarily smooth.
  Then a $\mathcal{V}$-fold according to Def. \ref{VManifold}
  is, in particular, an orbi-$\smooth \mathcal{V}$-fold according to
  Def. \ref{OrbiVFolds}, hence a $V$-fold for
  $V := \smooth \mathcal{V}$ the purely smooth aspect of
  $\mathcal{V}$:
       \vspace{-2mm}
  $$
    \xymatrix{
      \mathcal{V}\mathrm{Folds}(\mathbf{H})
     \; \ar@{^{(}->}[r]
      &
      (\smooth \mathcal{V})\mathrm{Orbifolds}(\mathbf{H})
    }
    \,.
  $$

     \vspace{-2mm}
\noindent
But, in general, being a $\mathcal{V}$-fold is a much stronger
  condition than being an $(\smooth \mathcal{V})$-orbifold,
  even (and in particular) if $\mathcal{V}$ is already smooth:
  For a $\mathcal{V}$-fold $\mathcal{X}$ not only the full
  underlying $\smooth \mathcal{X}$ is required to be locally
  equivalent to
  $\smooth \mathcal{V}$, but moreover, for each
  $K \in \mathrm{Groups}^{\mathrm{fin}}$, the geometric $K$-fixed
  locus of $\smooth \mathcal{X}$ is required to be locally
  equivalent to the geometric $K$-fixed locus of $\smooth \mathcal{V}$.
\end{example}

\begin{example}[Subcategories of smooth and of flat orbifolds]
  \label{SubcategorySmoothAndFlatOrbifolds}
  Let $\mathbf{H}$ be an elastic $\infty$-topos (Def. \ref{ElasticInfinityTopos}),
  $V \in \mathrm{Groups}(\mathbf{H})$ (Prop. \ref{LoopingAndDelooping})
  and
  $(G,\phi) \in \mathrm{Groups}(\mathbf{H})_{/\mathrm{Aut}(T_e V)}$
    (Prop. \ref{LoopingAndDelooping},  Def. \ref{AutomorphismGroup},
    Example \ref{LocalNeighbourhoodOfAPoint}).
    We have fully faithful inclusions into the $\infty$-category
  of $(G,\phi)$-structured orbi-$V$-folds (Def. \ref{OrbiVFolds})
   \vspace{-2mm}
  \begin{equation}
    \label{SmoothAndFlatOrbifolds}
    \raisebox{20pt}{
    \xymatrix@C=-6pt{
      &
      (G,\phi)\mathrm{Structured}V\mathrm{Orbifolds}(\mathbf{H})
      \\
      (G,\phi)\mathrm{Structured}V\mathrm{Folds}(\mathbf{H}_0)
      \ar@{^{(}->}[ur]_-{ i_{\scalebox{.6}{$\smooth$}} }^-{
        \mbox{
          \tiny
          \color{darkblue}
          \bf
          \begin{rotate}{16}
  \hspace{-1.2cm} smooth orbifolds
  \end{rotate}
                   }
      }
      &&
   \;\;\;\;\;\;  (\flat G, \phi \circ \epsilon^{\flat})\mathrm{Structured}V\mathrm{Orbifolds}(\mathbf{H})
      \ar@{_{(}->}[ul]^-{ i_{\scalebox{.6}{$\flat$}}}_-{
        \mbox{
          \tiny
          \color{darkblue}
          \bf
            \begin{rotate}{-14}
  \hspace{-.2cm} flat orbifolds
  \end{rotate}
                   }
      }
      &&
    }
    }
  \end{equation}
  of

  \noindent {\bf (i)} smooth $(G,\phi)$-structured $V$-folds,
  via Lemma \ref{Smooth0TruncatedObjectsAreOrbiSingular};

  \noindent {\bf (ii)} flat
  $(\flat G, \phi \circ \epsilon^\flat )$-structured $V$-folds
  (Def. \ref{FlatVFolds}).
\end{example}

\newpage

\section{Orbifold cohomology}
\label{OrbifoldCohomology}

With an internal $\infty$-topos-theoretic characterization of orbifolds
in hand (from \cref{OrbifoldGeometry}),
we immediately obtain an induced notion of
(differential, geometric, {\'e}tale) \emph{orbifold cohomology},
given by the intrinsic cohomology \eqref{IntrinsicCohomologyOfAnInfinityTopos}
of the ambient singular-cohesive $\infty$-topos.
Here we discuss how this new intrinsic notion of orbifold cohomology

\noindent - subsumes proper equivariant cohomology theory (\cref{EquivariantCohomology})

\noindent - and unifies it with tangentially twisted cohomology (\cref{JTwistedOrbifoldCohomotopy}).

\subsection{Proper equivariant cohomology}
\label{EquivariantCohomology}

\noindent {\bf Proper equivariant cohomology.}

\begin{defn}[Borel equivariant cohomology]
\label{BorelEquivariantCohomologyInSingularCohsion}
Let $\mathbf{H}_{\tiny\smooth}$ be a cohesive $\infty$-topos (Def. \ref{CohesiveTopos})
$G \in \mathrm{Groups}(\mathbf{H}_{\tiny\smooth})$ (Prop. \ref{LoopingAndDelooping})
and $(X,\tau), (A,\rho) \in G\mathrm{Actions}(\mathbf{H}_{\tiny\smooth})$
(Prop. \ref{InfinityAction}).
Then the \emph{Borel equivariant cohomology} of $X$ with
coefficients in $A$ is the intrinsic cohomology
\eqref{IntrinsicCohomologyOfAnInfinityTopos} in
the slice
$\mathbf{H}_{/_{\!\scalebox{.6}{$\mathbf{B}G$}}}$ (Prop. \ref{SliceInfinityTopos})
of the homotopy quotient \eqref{HomotopyQuotientAsColimit}
of $X$ with coefficients in the shape \eqref{CohesiveModalitiesFromAdjointQuadruple}
of the homotopy quotient of $A$:
\vspace{-2mm}
\begin{equation}
  \label{BorelEquivariantCohomology}
  \overset{
    \mathclap{
    \mbox{
      \tiny
      \color{darkblue}
      \bf
      \begin{tabular}{c}
        Borel equivariant
        \\
        cohomology
      \end{tabular}
    }
    }
  }{
    H_{\mathrm{Borel}}(X,A)
  }
  \;:=\;
  \pi_0 \;
  \mathbf{H}_{/_{\!\scalebox{.6}{$\mathbf{B}G$}}}
  \big(
    ( X \!\sslash\! G )
    \,,\,
    \raisebox{1pt}{\textesh} ( A \!\sslash\! G)
  \big)
  \;=\;
  \left(\!\!\!\!\!\!
  \raisebox{16pt}{
  \xymatrix{
    (X \!\sslash\! G)
    \ar[rr]^-{
      \overset{
        \raisebox{3pt}{
          \tiny
          \color{darkblue}
          \bf
          cocycle
        }
      }{c}
    }_-{\ }="s"
    \ar[dr]_-{
      \tau
    }^{\ }="t"
    &&
    A \!\sslash\! G
    \ar[dl]^-{
      \scalebox{.7}{$
        \rho
      $}
    }
    \\
    &
    \mathbf{B}G
    \ar@{=>} "s"; "t"
  }
  }
  \!\!\! \right)
\end{equation}
\end{defn}

\begin{defn}[Proper equivariant cohomology]
\label{ProperEquivariantCohomologyInSingularCohsion}
Let $\mathbf{H}$ be a singular-cohesive $\infty$-topos (Def. \ref{SingularCohesiveInfinityTopos})
$G \in \mathrm{Groups}(\mathbf{H}_\flat)$ (Prop. \ref{LoopingAndDelooping})
a discrete $\infty$-group,
and $(X,\tau), (A,\rho) \in G\mathrm{Actions}(\mathbf{H})$
(Prop. \ref{InfinityAction}).
Then we say that the \emph{proper equivariant cohomology} of $X$ with
coefficients in $A$ is the intrinsic cohomology
\eqref{IntrinsicCohomologyOfAnInfinityTopos} in
the slice
$\mathbf{H}_{/_{\!\scalebox{.6}{$\orbisingular \mathbf{B}G$}}}$ (Prop. \ref{SliceInfinityTopos})
of the orbi-singularization \eqref{SingularityModalities}
of the homotopy quotient \eqref{HomotopyQuotientAsColimit}
of $X$ with coefficients in the shape \eqref{CohesiveModalitiesFromAdjointQuadruple}
of the orb-singularization of the homotopy quotient of $A$:
\vspace{-2mm}
\begin{equation}
  \label{ProperEquivariantCohomology}
  \overset{
    \mathclap{
    \mbox{
      \tiny
      \color{darkblue}
      \bf
      \begin{tabular}{c}
        proper equivariant
        \\
        cohomology
      \end{tabular}
    }
    }
  }{
    H_G(X,A)
  }
  \;:=\;
  \pi_0 \;
  \mathbf{H}_{/_{\!\scalebox{.6}{$\orbisingular \mathbf{B}G$}}}
  \big(
    \orbisingular ( X \!\sslash\! G )
    \,,\,
    \raisebox{1pt}{\textesh} \orbisingular ( A \!\sslash\! G)
  \big)
  \;=\;
  \left(\!\!\!\!\!\!
  \raisebox{19pt}{
  \xymatrix{
    \orbisingular (X \!\sslash\! G)
    \ar[rr]^-{
      \overset{
        \mbox{
          \tiny
          \color{darkblue}
          \bf
          cocycle
        }
      }{c}
    }_-{\ }="s"
    \ar[dr]_-{
      \scalebox{.7}{$
        \orbisingular(\tau)
      $}
    }^{\ }="t"
    &&
    \raisebox{1pt}{\textesh} \orbisingular ( A \!\sslash\! G)
    \ar[dl]^-{\;\;
      \scalebox{.7}{$
        \big(
          \eta^{\scalebox{.6}{\textesh}}_{ \orbisingular \mathbf{B}G }
        \big)^{-1}
        \circ
        \raisebox{1pt}{\textesh}\orbisingular(\rho)
      $}
    }
    \\
    &
    \orbisingular \mathbf{B}G
    \ar@{=>} "s"; "t"
  }
  }
 \!\!\! \right)
\end{equation}

\end{defn}

\noindent {\bf Recovering traditional $G$-equivariant cohomology.}
We discuss how in the case of a finite group $G$, traditional $G$-equivariant
cohomology (see \cref{GEquivariantHomotopyTheory}) is a special case of
proper equivariant cohomology (Def. \ref{ProperEquivariantCohomologyInSingularCohsion}).
We take the key observation from \cite{Rezk14}
(Prop. \ref{SingularitiesFaithfulSliceOverGSingIsGEquivariantHomotopyTheory} below).

\begin{defn}[$G$-equivariant cohesive $\infty$-topos]
  \label{GEquivariantCohesiveInfinityTopos}
  Let $\mathbf{H}_{\tiny\smooth}$ be a cohesive $\infty$-topos
  (Def. \ref{CohesiveTopos})
  and $G \in \mathrm{Groups}^{\mathrm{fin}}$ a
  finite group \eqref{InclusionOfFiniteGroups}.
  We write
   \vspace{-2mm}
  \begin{equation}
    \label{TheGEquivariantCohesiveInfinityTopos}
    G\mathbf{H}_{\tiny\smooth}
    \;:=\;
    \mathrm{Sheaves}_\infty
    \big(
      G\mathrm{Orbits},
      \,
      \mathbf{H}_{\tiny\smooth}
    \big)
    \;=\;
    \mathrm{Func}_\infty
    \big(
      G \mathrm{Orbits}^{\mathrm{op}},
      \,
      \mathbf{H}_{\tiny\smooth}
    \big)
  \end{equation}

  \vspace{-2mm}
\noindent
  for the $\infty$-topos of $\mathbf{H}_{\tiny\smooth}$-valued
  $\infty$-sheaves on the $G$-orbit category
  (Def. \ref{OrdinaryOrbitCategory}),
  to be called the
  corresponding \emph{$G$-equivariant cohesive $\infty$-topos}.
\end{defn}

\begin{remark}[Proper equivariant cohomology theory in singular $\infty$-toposes]
  \label{ProperEquivariantCohomologyTheory}
  In the case $\mathbf{H}_{\tiny\smooth} \simeq \mathrm{Groupoids}_{\infty}$
  \eqref{InfinityCategoryOfInfinityGroupoids},
  Def. \ref{GEquivariantCohesiveInfinityTopos} reduces to the
  $\infty$-category
  $G \mathrm{Groupoids}_\infty$ (Def. \ref{HomotopyTheoryOfGSpaces})
  of traditional $G$-equivariant homotopy theory
  (recalled in \cref{GEquivariantHomotopyTheory}).
  The intrinsic cohomology \eqref{IntrinsicCohomologyOfAnInfinityTopos}
  of the $\infty$-topos $G \mathrm{Groupoids}_\infty$
  --
  or of its tangent $\infty$-topos
  $T\big( G \mathrm{Groupoids}_\infty \big)$
  (Example \ref{TangentInfinityTopos})
  in the twisted abelian case (Remark \ref{AbelianTwistedCohomology})
  --
  is
  \emph{proper equivariant cohomology}
  (following terminology in \cite{DHLPS19}),
  including
  $G$-Bredon cohomology
  \cite{Bredon67a}\cite{Bredon67b} (review in \cite[\S 1.4]{Blu17}\cite[\S 7]{tomDieck79}),
  $G$-equivariant K-theory
  \cite{Segal68}\cite{AtiyahSegal69}
  (which is proper equivariant by \cite[A3.2]{AtiyahSegal04}\cite[A.5]{FHT07}\cite{DavisLuck98}),
  $G$-equivariant Cohomotopy theory
  \cite{Segal71}\cite[\S 8]{tomDieck79}\cite{SS19a}\cite{SS19b},
  etc.

  Hence,
  by Remark \ref{DifferentialCohomologyTheory},
  to the extent that the objects of
  the cohesive $\infty$-topos $\mathbf{H}_{\tiny\smooth}$
  in
  Def. \ref{GEquivariantCohesiveInfinityTopos} are $\infty$-groupoids
  equipped with further geometric or differential-geometric structure,
  the intrinsic cohomology theory \eqref{IntrinsicCohomologyOfAnInfinityTopos}
  in
  $G\mathbf{H}_{\tiny\smooth}$
  \eqref{TheGEquivariantCohesiveInfinityTopos}
  is an enhancement of plain $G$-equivariant cohomology to
  a flavor of
  \emph{proper $G$-equivariant differential cohomology} theory
  (by Remark \ref{DifferentialCohomologyTheory}).
\end{remark}

\begin{prop}[Cohesive Elmendorf theorem]
  \label{CohesiveElmendorfEquivalence}
  Consider a cohesive $\infty$-topos $\mathbf{H}_{\tiny\smooth}$
  (Def. \ref{CohesiveTopos})
  with an $\infty$-site $\mathrm{Charts}$ of charts
  (Def. \ref{ChartsForCohesion}).
  Then for $G \in \mathrm{Groups}^{\mathrm{fin}}$ a finite group,
  we have an equivalence of $\infty$-categories
\vspace{-1mm}
  \begin{equation}
    G\mathbf{H}_{\tiny\smooth}
    \;\simeq\;
    \mathrm{Sheaves}_\infty
    \big(
      \mathrm{Charts},
      G \mathrm{Groupoids}_\infty
    \big)\,,
  \end{equation}

   \vspace{-1mm}
\noindent
  where $G \mathrm{Groupoids}_\infty$ is the
  $\infty$-category of D-topological $G$-spaces (Def. \ref{HomotopyTheoryOfGSpaces}).
\end{prop}
\begin{proof}
  Consider the following sequence of $\infty$-functors
   \vspace{-2mm}
  $$
    \begin{aligned}
      G\mathbf{H}_{\tiny\smooth}
      &  := \;
      \mathrm{Sheaves}_\infty
      \big(
        G\mathrm{Orbits},
        \mathbf{H}_{\tiny\smooth}
      \big)
      \\
      & = \;
      \mathrm{Sheaves}_\infty
      \big(
        G\mathrm{Orbits},
        \mathrm{Sheaves}_\infty(\mathrm{Charts})
      \big)
      \\
&      \overset{\simeq}{\to}\;
      \mathrm{Sheaves}_\infty
      \big(
        G\mathrm{Orbits}
        \times
        \mathrm{Charts}
      \big)
      \\
&      \overset{\simeq}{\to}\;
      \mathrm{Sheaves}_\infty
      \big(
        \mathrm{Charts},
        \mathrm{Sheaves}_\infty(G \mathrm{Orbits})
      \big)
      \\
&      \overset{\simeq}{\to}\;
      \mathrm{Sheaves}_\infty
      \big(
        \mathrm{Charts},
        G \mathrm{Groupoids}
      \big).
    \end{aligned}
  $$

   \vspace{-2mm}
\noindent
  That the first and second of these
  $\infty$-functors are equivalences
  follows by the product/hom-adjunction for $\infty$-functors.
  With that, the last equivalence follows,
  objectwise, by Elmendorf's theorem
  (Prop. \ref{ElmendorfTheorem}).
\hfill \end{proof}

\begin{prop}[G-equivariant homotopy theory embeds into $G$-singular cohesion]
  \label{SingularitiesFaithfulSliceOverGSingIsGEquivariantHomotopyTheory}
  Let $\mathbf{H}$ be a singular-cohesive $\infty$-topos
  (Def. \ref{SingularCohesiveInfinityTopos})
  over $\mathrm{Groupoids}_\infty$ \eqref{InfinityCategoryOfInfinityGroupoids}
  and let $G \in \mathrm{Grp}_{\mathrm{fin}}$ be a finite group
  \eqref{InclusionOfFiniteGroups}.
   \vspace{-2mm}
 \item {\bf (i)}  Then there is
  a full sub-$\infty$-category inclusion
  \vspace{-3mm}
  \begin{equation}
    \label{EquivariantCohomologyBySlicing}
    \xymatrix@R=9pt{
      G\mathbf{H}_{\tiny\smooth}
      \ar@{^{(}->}[rr]^-{\Delta_G}_-{\simeq}
      &&
      \mathbf{H}_{ /\scalebox{.7}{$\orbisingularG$}}
    }
  \end{equation}

 \vspace{-3mm}
\noindent   of
  the $G$-equivariant non-singular cohesive $\infty$-topos
  (Def. \ref{GEquivariantCohesiveInfinityTopos})
  into
  the slice of $\mathbf{H}$ (Prop. \ref{SliceInfinityTopos})
  over the generic $G$-orbi singularity \eqref{SingularitiesUnderYoneda}.

 \vspace{-2mm}
  \item {\bf (ii)} This is such that,
  when pre-composed with the cohesive Elmendorf equivalence
  (Prop. \ref{CohesiveElmendorfEquivalence}),
  a cohesive sheaf (on $\mathrm{Charts}$)
  of $G \mathrm{Groupoids}$
  \eqref{InfinityCategoryOfGInfinityGroupoids}
  presented \eqref{ShapeOfGTopologicalSpaces} by
  D-topological $G$-spaces $X_U$ (Def. \ref{GSpaces})
  is sent to the presheaf on $\mathrm{Singularities}$ that
  is given as follows:
  \vspace{-2mm}
  \begin{equation}
    \label{DeltaGInComponents}
    \hspace{-4mm}
    \raisebox{30pt}{
    \xymatrix@R=-2pt{
      \mathrm{Sheaves}_\infty
      \big(
        \mathrm{Charts},
        \,
        G\mathrm{Groupoids}_\infty
      \big)
      \ar[r]^-{\simeq}
      &
      G\mathbf{H}_{\tiny\smooth}
      \ar[r]|-{\;\;  \small \underset{\scalebox{.7}{$\orbisingularG$} }{\sum} \Delta_G \; }
      &
      \mathrm{Sheaves}_{\infty}
      \big(
        \mathrm{Charts}
        \times
        \mathrm{Singularities}
      \big)
      \\
      \left(
        U
          \;\mapsto\;
        \mathrm{Shp}_{G\mathrm{Top}}
        \big(
          X_U
        \big)
      \right)
   \;\;\; \;\;\;\;\;\;  \ar@{|->}[rr]
      &&
      \;\;\;\;\;\;\;\;\;\left(\!\!
        \big(
          U,
          \orbisingularK
        \big)
        \mapsto
        \mathrm{Shp}_{\mathrm{Top}}
        \left(\!\!
          \left(
          \;\;\;\;\;\;\;
      \raisebox{4pt}{$    \underset{{}_{
            \mathclap{
              \phi \in \mathrm{Groups}(K,G)
            }}
          }{
            \bigsqcup
          }
          \;\;\;\;\;
          X_U^{\phi(K)}
          $}
          \right)
        \!\underset{G}{\times}
        E G
      \right)
    \!\!\right)
        }
    }
  \end{equation}

  \vspace{-2mm}
 \noindent where on the right we have the topological shape \eqref{ShapeOfTopologicalSpaces}
  of
  the Borel construction by the residual $G$-action
  on the fixed point subspaces $ X_U^{\phi(K)} \subset X_U$
  \eqref{HFixedPointSpaces}.
\end{prop}
\begin{proof}
  For $\mathbf{H}_{\tiny\smooth} \simeq \mathrm{Groupoids}_{\infty}$
  this is \cite[Prop. 3.5.1]{Rezk14};
  our expression
  $
    \mathrm{Shp}_{\mathrm{Top}}
    \big(
      X_U^{\phi(K)} \underset{G}{\times} E G
    \big)
  $
  is, up to convention of notation, the expression for
  $B \mathrm{Fun}(
    H,
    G \curvearrowright X_U
  )
    $ that is spelled out in \cite[p. 7]{Rezk14}\cite[3.2.17]{Lurie19}
    (using that our $G$ is discrete).
    The generalization here follows immediately by applying this equivalence
   objectwise in the $\infty$-site $\mathrm{Charts}$.
\hfill \end{proof}

The following is our key class of examples:

\begin{example}[Cohesive shape of $G$-orbi-singular space is $G$-homotopy type]
  \label{CohesiveShapeOfGOrbiSungularSpacesIsGHomotopyType}

  In the cohesive $\infty$-topos
  $\mathbf{H}_{\tiny\smooth} \coloneqq \mathrm{SmoothGroupoids}_\infty$
  (Example \ref{SmoothInfinityGroupoids})
  consider a 0-truncated object
  $X \in \mathbf{H}_{{\tiny\smooth},0}$
  equipped with a
  $G$-action (Def. \ref{InfinityActionHomotopyFiberSequence})
  of a discrete group $G$,
  and with corresponding Cohesive $G$-orbispace
  (Prop. \ref{PropertiesOfUniversalCoveringSpaces})
  \vspace{-2mm}
  $$
    \mathcal{X} \;:=\; \mathrm{OrbSnglr}(X\!\sslash\!G)
  $$

  \vspace{-2mm}
\noindent
  in $\mathbf{H} \coloneqq \mathrm{SingularSmoothGroupoids}_\infty$
  (Example \ref{OrbiSingularSmoothInfinityGroupoids}),
  which is either of:

  \vspace{.2cm}

  \hspace{-.9cm}
  \begin{tabular}{lll}
    {\bf (i)}
       a smooth $G$-orbifold:
    &
    $X \;\in\;
     \xymatrix{
       \mbox{SmoothManifolds}
       \;\ar@{^{(}->}[r]
       &
       \mathrm{DiffeologicalSpaces}
      \; \ar@{^{(}->}[r]
       &
       \mathbf{H}_{\tiny\smooth}
     }$
    &
    (Example \ref{FrechetSmoothGOrbifolds})
    \\
    {\bf (ii)}
        a topological $G$-orbi space:
    &
    $X \in
      \xymatrix{
        \mathrm{TopologicalSpaces}
       \ar[r]^-{\scalebox{.6}{$\mathrm{Cdfflg}$}}
        &
        \mathrm{DTopologicalSpaces}
     \;   \ar@{^{(}->}[r]
        &
        \mathbf{H}_{\tiny\smooth}
      }
    $
    &
    (Example \ref{TopologicalGOrbiSingularSpaces})
  \end{tabular}

  \vspace{2mm}
\noindent   Then the cohesive shape \eqref{CohesionSingular}
  of the $G$-orbi-singular space $\mathcal{X} \in \mathbf{H}$
  is equivalent, under the identification of Prop.
  \ref{SingularitiesFaithfulSliceOverGSingIsGEquivariantHomotopyTheory},
  to the $G$-topological shape \eqref{ShapeOfGTopologicalSpaces}
  of the underlying topological $G$-space of $X$:

  \vspace{.1cm}

  \noindent {\bf (i)} By Prop. \ref{ShapeOfFrechetSmoothGOrbifold},
  comparing \eqref{KValueOfShapeOfGOrbiSingularFrechetManifold} with
  \eqref{DeltaGInComponents} we have:

  \vspace{-.2cm}

  \begin{equation}
    \raisebox{40pt}{
    \xymatrix@R=20pt@C45pt{
      G \mbox{SmoothManifolds}
      \ar[dd]^-{
        \mathllap{
          \mbox{
            \tiny
            \color{darkblue}
            \bf
            \begin{tabular}{c}
              form
              \\
              $G$-topological shape
            \end{tabular}
          }
          \;\;\;
        }
     \scalebox{.6}{$    \mathrm{Shp}_{G\mathrm{Top}}
        \left(
          \mathrm{Dtplg}(-)
        \right)
        $}
      }
      \ar[rr]^-{
            \tiny
            \color{darkblue}
            \bf
            \scalebox{1}{\bf form Fr{\'e}chet-smooth orbifold}}_-{
              \scalebox{.6}{$
          \mathrm{OrbSnglr}
          \left(
            (-) \sslash G
          \right)
          $}
        }
      &&
      {\mathrm{SingularSmoothGroupoids}_\infty}_{/\scalebox{.7}{$\orbisingularG$}}
      \ar[dd]_-{ \scalebox{.6}{$
        \mathrm{Shp}
        $}
        \mathrlap{
          \;\;\;
          \mbox{
            \tiny
            \color{darkblue}
            \bf
            \begin{tabular}{c}
              form
              \\
              cohesive shape
            \end{tabular}
          }
        }
      }
      \\
      \\
      G \mathrm{Groupoids}_\infty
      \ar@{^{(}->}[rr]_-{
                \tiny
            \color{darkblue}
              \begin{tabular}{c}
        \bf      include $G$-equivariant homotopy theory
            \end{tabular}
          }^-{ \scalebox{.7}{$
          \Delta_G$}
        }
            &&
      {\mathrm{SingularGroupoids}_\infty}_{/\scalebox{.7}{$\orbisingularG$}}
    }
    }
  \end{equation}

  \vspace{.1cm}

  \noindent {\bf (ii)} By Prop. \ref{ShapeOfGOrbiSingularTopologicalSpaces},
  comparing \eqref{KValueOfShapeOfGOrbiSingularTopologicalSpace} with
  \eqref{DeltaGInComponents} we have:
  \vspace{-.2cm}

  \begin{equation}
    \raisebox{40pt}{
    \xymatrix@R=20pt@C55pt{
      G \mathrm{TopologicalSpaces}
      \ar[dd]^-{
        \mathllap{
          \mbox{
            \tiny
            \color{darkblue}
            \bf
            \begin{tabular}{c}
              form
              \\
              $G$-topological shape
            \end{tabular}
          }
          \;\;\;
        }
      \scalebox{.6}{$   \mathrm{Shp}_{G\mathrm{Top}}
      $}
      }
      \ar[rr]^-{
                \tiny
            \color{darkblue}
            \bf
            \begin{tabular}{c}
       \bf       form topological $G$-orbi space
            \end{tabular}
          }_-{ \scalebox{.6}{$
          \mathrm{OrbSnglr}
          \left(
            \mathrm{Cdfflg}(-) \sslash G
          \right)
          $}
        }
      &&
      {\mathrm{SingularSmoothGroupoids}_\infty}_{/\scalebox{.7}{$\orbisingularG$}}
      \ar[dd]_-{ \scalebox{.6}{$
        \mathrm{Shp}
        $}
        \mathrlap{
          \;\;\;
          \mbox{
            \tiny
            \color{darkblue}
            \bf
            \begin{tabular}{c}
              form
              \\
              cohesive shape
            \end{tabular}
          }
        }
      }
      \\
      \\
      G \mathrm{Groupoids}_\infty
      \ar@{^{(}->}[rr]_-{
             \tiny
            \color{darkblue}
            \begin{tabular}{c}
          \bf    include $G$-equivariant homotopy theory
            \end{tabular}
          }^-{ \scalebox{.7}{$
          \Delta_G$}
        }
          &&
      {\mathrm{SingularGroupoids}_\infty}_{/\scalebox{.7}{$\orbisingularG$}}
    }
    }
  \end{equation}

\end{example}

\begin{lemma}[$\Delta_G$ commutes with Disc]
  \label{DeltaGCommutesWithDisc}
  The construction $\Delta_G$ from Prop.
  \ref{SingularitiesFaithfulSliceOverGSingIsGEquivariantHomotopyTheory}
  commutes with embedding of discrete cohesive structure
  \eqref{CohesionOfSingularCohesive}:
 \vspace{-2mm}
  $$
    \xymatrix@R=2pt@C=4em{
      &
      \mathrm{Sheaves}_\infty
      \big(
        \mathrm{Singularities},
        \mathrm{Groupoids}_\infty
      \big)_{/\scalebox{.7}{$\orbisingularG$}}
      \ar[dr]^-{  \scalebox{.6}{$ \mathrm{Disc}$} }
      \\
      G \mathrm{Groupoids}
      \ar[dr]_-{ \scalebox{.6}{$ \mathrm{Disc}$} }
      \ar[ur]^-{ \Delta_G }
      &&
      \mathrm{Sheaves}_\infty
      \big(
        \mathrm{Singularities},
        \,
        \mathbf{H}_{\tiny\smooth}
      \big)_{/\scalebox{.7}{$\orbisingularG$}}
      \\
      &
      G \mathbf{H}_{\tiny\smooth}
      \ar[ur]_-{ \Delta_G }
    }
  $$
\end{lemma}

\begin{theorem}[Cohomology of good  orbispaces is proper equivariant cohomology]
  \label{OrbifoldCohomologyEquivariant}
  Consider the singular-cohesive $\infty$-topos
  $
    \mathbf{H}
    \;:=\;
    \mathrm{SingularCohesiveGroupoids}_\infty
  $
  (Example \ref{OrbiSingularSmoothInfinityGroupoids})
  and let $G \in \mathrm{Groups}^{\mathrm{fin}}$ be
  a discrete group \eqref{InclusionOfFiniteGroups}.
  Then the intrinsic cohomology \eqref{IntrinsicCohomologyOfAnInfinityTopos}

\item{\bf (i)}
  of a $G$-orbi-singular space $\mathcal{X} \in \mathbf{H}_{/\scalebox{.7}{\orbisingularG}}$
  (Def. \ref{GOrbiSpace})
  which is either
  \begin{itemize}
   \vspace{-2mm}
  \item[{\bf (a)}] a topological $G$-orbi-space (Example \ref{TopologicalGOrbiSingularSpaces}) with universal covering
  space (Def. \ref{CoveringSpaceOfGOrbiSingularSpace})
  $X_{G\mathrm{top}} \in G \mathrm{TopologicalSpaces}$ \eqref{GTopologicalSpaces};
 \vspace{-2mm}
  \item[{\bf (b)}] a Fr{\'e}chet-smooth $G$-orbifold (Example \ref{FrechetSmoothGOrbifolds}) with universal covering
  space (Def. \ref{CoveringSpaceOfGOrbiSingularSpace})
  $X \in \mbox{Fr{\'e}chetManifolds}$ and underlying $G$-topological space
  $X_{G\mathrm{top}} := \mathrm{Dtplg(X)}$ \eqref{AdjunctionBetweenTopologicalAndDiffeologicalSpaces};
  \end{itemize}
 \vspace{-3mm}
  \item {\bf (ii)} with
  coefficients in a cohesively discrete
  $G$-$\infty$-groupoid $A$ \eqref{InfinityCategoryOfGInfinityGroupoids}
  (hence the $G$-topological shape \eqref{ShapeOfGTopologicalSpaces}
  of some topological $G$-space $A_{G\mathrm{top}}$)
  regarded as a geometrically discrete
  orbi-singular $\infty$-groupoid $\mathcal{A}$
  via \eqref{EquivariantCohomologyBySlicing}:
  \vspace{-1mm}
  $$
    \xymatrix@R=-1pt{
      G \mathrm{TopologicalSpaces}
      \ar[rr]^-{ \scalebox{.6}{$ \mathrm{Shp}_{G\mathrm{Top}}$} }
      &&
      G \mathrm{Groupoids}_\infty
      \ar[rr]^-{ \scalebox{.6}{$
        \mathrm{Disc}
        $}
      }
      &&
      G\mathbf{H}_{\tiny\smooth}
      \ar[rr]^-{ \Delta_G }
      &&
      \mathbf{H}_{/\scalebox{.7}{$\orbisingularG$}}
      \\
      A_{\mathrm{top}}
      \ar[rr]
      &&
      A
      \ar@{|->}[rrrr]
      &&&&
      \mathcal{A}
    }
  $$
  equals the proper $G$-equivariant cohomology (Def. \ref{GenuineGEquivariantCohomology})
  of $X_{G\mathrm{top}}$
  with coefficients in $A$:
 \vspace{-1mm}
  $$
    \begin{array}{lccc}
    &
    \mathbf{H}_{/\scalebox{.7}{$\orbisingularG$}}
    \big(
      \mathcal{X},
      \,
      \mathcal{A}
    \big)
    &\;\simeq\;&
    G \mathrm{Groupoids}_\infty
    \big(
      \mathrm{Shp}_{G\mathrm{Top}}(X_{G\mathrm{top}}),
      \,
      A
    \big)
    \\
    \phantom{\mathclap{\vert^{\vert^{\vert^{\vert^{\vert}}}}}}
    \mbox{hence:}
    \phantom{AAA}
    \phantom{\mathclap{\vert_{\vert_{\vert_{\vert_{\vert}}}}}}
    &
    \pi_{n}
    \mathbf{H}_{/\scalebox{.7}{$\orbisingularG$}}
    \big(
      \mathcal{X},
      \,
      \mathcal{A}
    \big)
    &\;\simeq\;&
    H^{-n}_G(X_{G\mathrm{top}},A)
    \\
    &
    \mathclap{
    \mbox{
      \tiny
      \color{darkblue}
      \bf
      \begin{tabular}{c}
        intrinsic
        \\
        equivariant differential cohomology
        \\
        in $\infty$-topos of
        \\
        singular smooth $\infty$-groupoids
      \end{tabular}
    }
    }
    &&
    \mbox{
      \tiny
      \color{darkblue}
      \bf
      \begin{tabular}{c}
        proper
        \\
        $G$-equivariant cohomology
      \end{tabular}
    }
        \end{array}
  $$
\end{theorem}

\newpage

\begin{proof}
  {\bf (i)}
  By Example \ref{TopologicalGOrbiSingularSpaces}
  the topological $G$-orbi space $\mathcal{X}$  is given by
  \vspace{-2mm}
  $$
    \mathcal{X}
    \;\simeq\;
    \mathrm{OrbSnglr}
    \big(
      \mathrm{Cdfflg}(X) \!\sslash \! G
    \big)
    \,.
  $$

  \vspace{-2mm}
\noindent
  With this, we compute as follows:

  \vspace{-3mm}
  \begin{equation}
    \label{ComputingGOrbiSpaceCohomology}
    \begin{aligned}
    \xymatrix{
      \mathbf{H}_{/\scalebox{.7}{$\orbisingularG$}}
      \big(
        \mathcal{X},
        \mathcal{A}
      \big)
    }
    & =
    \mathbf{H}_{/\scalebox{.7}{$\orbisingularG$}}
    \Big(
      \mathrm{OrbSnglr}
      \big(
        \mathrm{Cdfflg}(X_{\mathrm{top}}) \!\sslash\! G
      \big)
      ,
      \Delta_G \mathrm{Disc}(A)
    \Big)
    \\
    &
    \simeq
    \mathbf{H}_{/\scalebox{.7}{$\orbisingularG$}}
    \Big(
      \mathrm{OrbSnglr}
      \big(
        \mathrm{Cdfflg}(X_{\mathrm{top}}) \!\sslash\! G
      \big)
      ,
      \mathrm{Disc}(\Delta_G A)
    \Big)
    \\
    & \simeq
    (\mathrm{Groupoids}_\infty)_{/\scalebox{.7}{$\orbisingularG$}}
    \Big(
      \mathrm{Shp}
      \big(
      \mathrm{OrbSnglr}
      (
        \mathrm{Cdfflg}(X_{\mathrm{top}}) \!\sslash\! G
      )
      \big)
      ,
      \Delta_G A
    \Big)
    \\
    & \simeq
    (\mathrm{Groupoids}_\infty)_{/\scalebox{.7}{$\orbisingularG$}}
    \big(
      \Delta_G X,
      \Delta_G A
    \big)
    \\
    & \simeq
    G \mathrm{Groupoids}
    \big(
      \mathrm{Shp}_{G\mathrm{Top}}(X_{\mathrm{top}}),
      A
    \big)
    \,.
    \end{aligned}
  \end{equation}

  \vspace{-1mm}
\noindent
  Here the first step, after unwinding the definitions, is Lemma
  \ref{DeltaGCommutesWithDisc}. The second step is the
  $\mathrm{Shp} \dashv \mathrm{Disc}fc$-adjunction \eqref{CohesionOfSingularCohesive}.
  The third step is Prop. \ref{ShapeOfGOrbiSingularTopologicalSpaces}.
  The last step is Prop.
  \ref{SingularitiesFaithfulSliceOverGSingIsGEquivariantHomotopyTheory}

  \vspace{.1cm}

  \noindent {\bf (ii)}
  By Example \ref{FrechetSmoothGOrbifolds}
  the Fr{\'e}chet-smooth $G$-orbifold $\mathcal{X}$  is given by
  \vspace{-2mm}
  $$
    \mathcal{X}
    \;\simeq\;
    \mathrm{OrbSnglr}
    \big(
      X \!\sslash\! G
    \big)
    \,.
  $$

  \vspace{-2mm}
\noindent
  With this, we compute just as in \eqref{ComputingGOrbiSpaceCohomology}
  only that now in the third step we use Prop. \ref{ShapeOfFrechetSmoothGOrbifold}.
\hfill \end{proof}

\begin{example}[Orientifold cohomology]
  \label{OrientifoldCohomology}
  Take the singular elastic $\infty$-topos
  $\mathbf{H} = \mathrm{SingularJetsOfSmoothGroupoids}_\infty$
  (Example \ref{OrbiSingularSmoothInfinityGroupoids})
  and $V = (\mathbb{R}^n, +) \in \mathbf{H}_{\tiny\smooth}$ \eqref{RnasV}.
  Then a $\mathcal{X}_{\tiny\smooth} \in V\mathrm{Folds}(\mathbf{H}_{\tiny\smooth})$
  (Def. \ref{VManifold}) is an
  ordinary $n$-dimensional orbifold or, more generally,
  an $n$-dimensional {\'e}tale $\infty$-stack
  (by Example \ref{EtaleLieGroupoidAsRnFold}) with
  structure group (Def. \ref{StructureGroupOfVFolds}) the
  ordinary general linear group
  $\mathbf{Aut}(T_e V) \simeq \mathrm{GL}(n)$
  (by Example \ref{OrdinaryGeneralLinearGroup}).
  Hence, the composition of the
  delooping \eqref{LoopingDeloopingEquivalence}
  of the ordinary determinant group homomorphism
  $\mathrm{GL}(n) \overset{\mathrm{det}}{\longrightarrow} \mathbb{Z}_2$
  with the
  classifying map $\vdash \mathrm{Frames}(\mathcal{X}_{\tiny\smooth})$
  \eqref{FrameBundleClassifyingMap}
  of the frame bundle of $X$ (Def. \ref{TangentBundleOfVFoldIsFiverBundle})
  realizes $\mathcal{X}_{\tiny\smooth}$ as an object in the
  slice $\infty$-topos (Prop. \ref{SliceInfinityTopos})
  over $\mathbf{B}\mathbb{Z}_2$. Consequently, it realizes its
  orbi-singualrization
  $\mathcal{X} := \orbisingular \mathcal{X}_{\tiny\smooth} \in \mathbf{H}$ \eqref{SincularCohesionModalies}
  as an object in the slice over
  $\orbisingular^{\hspace{-4pt}\scalebox{.61}{\raisebox{2pt}{$\mathbb{Z}_2$}}}$
  \eqref{AnOrbifoldSingularity}:
  \vspace{-2mm}
  \begin{equation}
    \label{OrbiOrientifold}
    \raisebox{20pt}{
    \xymatrix{
      \mathcal{X}_{\tiny\smooth}
      \ar[d]_-{\scalebox{0.7}{$
        \mathbf{B}\mathrm{det}
          \,\circ\,
        \vdash \mathrm{Frames}(\mathcal{X}_{\tiny\smooth})
        $}
      }
      \\
      \mathbf{B}\mathbb{Z}_2
    }
    }
    \;\;
    \in
    \big(\mathbf{H}_{\tiny\smooth}\big)_{/\mathbf{B}\mathbb{Z}_2}
    \;\;\;\;\;\;\;\;\;\;
    \Leftrightarrow
    \;\;\;\;\;\;\;\;\;\;
    \raisebox{20pt}{
    \xymatrix{
      \mathcal{X}
      \ar[d]_-{\scalebox{0.7}{$
        \orbisingular
        \big(
          \mathbf{B}\mathrm{det}
            \,\circ\,
          \vdash \mathrm{Frames}(\mathcal{X}_{\tiny\smooth})
        \big)
        $}
      }
      \\
      \orbisingular^{\hspace{-4pt}\scalebox{.5}{$\mathbb{Z}_2$}}
    }
    }
    \;\;
    \in
    \big(\mathbf{H}_{\tiny\smooth}\big)_{\big/
      \scalebox{.8}{$\orbisingular^{\hspace{-1.5mm}\raisebox{2pt}{\scalebox{.6}{$\mathbb{Z}_2$}}}$}
    }
    \,.
  \end{equation}

  \vspace{-2mm}
\noindent  This is the incarnation of the orbifold as an
  \emph{orbi-orientifold} \cite{DFM11}\cite[4.4]{FSS15}\cite{SS19a}.
  In particular, if the covering space (Def. \ref{CoveringSpaceOfGOrbiSingularSpace})
  $$
    X := \mathrm{fib}
    \big(
      \mathbf{B}\mathrm{det}
        \,\circ\,
      \vdash \mathrm{Frames}(\mathcal{X}_{\tiny\smooth})
    \big)
  $$
  happens to be an $\mathbb{R}^n$-fold (Example \ref{OrdinaryManifolds}),
  we have just a plain \emph{orientifold} (without further orbifolding)
  and then
  the intrinsic cohomology \eqref{IntrinsicCohomologyOfAnInfinityTopos}
  of $\mathcal{X}$ regarded in the slice over
  $\orbisingular^{\hspace{-4pt}\scalebox{.5}{$\mathbb{Z}_2$}}$\eqref{OrbiOrientifold}
  is, by Theorem \ref{OrbifoldCohomologyEquivariant}
  the proper $\mathbb{Z}_2$-equivariant cohomology of $X$,
  such as, for instance, Real K-theory \cite{Atiyah66}
  (see \cite{Masulli11} for the perspective in proper equivariant cohomology)
  or $\mathbb{Z}_2$-Equivariant Cohomotopy \cite[8.4]{tomDieck79}\cite{SS19a}.
\end{example}

\subsection{Proper orbifold cohomology}
\label{JTwistedOrbifoldCohomotopy}

We introduce general \emph{\'etale cohomology} of
{\'e}tale $\infty$-stacks (Def. \ref{EtaleCohomologyOfVFolds}),
which is sensitive to geometric $G$-structure
and to tangential structure (Def. \ref{TangentiallyTwistedCohomology}).
Promoting this to the \emph{proper} incarnation of orbifolds
(Remark \ref{Proper}),
we finally obtain \emph{tangentially twisted proper orbifold cohomology}
(Def. \ref{OrbifoldCohomologyTangentiallyTwisted}) which
we prove unifies tangentially twisted topological cohomology
away from orbifold singularities
with proper equivariant cohomology at the singularities
(Theorem \ref{OrbifoldCohomologyTangentiallyTwistedReduces}).
As a fundamental class of examples, we construct
J-twisted proper orbifold Cohomotopy theories (Def. \ref{JTwistedOrbifoldCohomotopyTheory})
and observe, as an application, that these subsume the relevant cohomology theories
for non-perturbative string theory, according to ``Hypothesis H''
(Remark \ref{HypothesisH}).

\newpage

\noindent {\bf Cohomology of $V$-{\'e}tale $\infty$-stacks.}
\begin{defn}[{\'E}tale cohomology]
  \label{EtaleCohomologyOfVFolds}
  Let $\mathbf{H}$ be an elastic $\infty$-topos  (Def. \ref{ElasticInfinityTopos}),
  $V \in \mathrm{Groups}(\mathbf{H})$
  (Prop. \ref{LoopingAndDelooping}),
  $(G,\phi) \in \mathrm{Groups}(\mathbf{H})_{/_{\mathrm{Aut}(T_e V)}}$
  (Def. \ref{GStructureCoefficients}),
  and $X \in \mathrm{Integrably}(G,\phi)\mathrm{Structured}V\mathrm{Folds}(\mathbf{H})$
  (Def. \ref{IntegrableGStructure}).
      The \emph{{\'e}tale cohomology}
  of $\big( X, (\tau,g)\big)$
  is its intrinsic cohomology \eqref{IntrinsicCohomologyOfAnInfinityTopos}
  when regarded (via Prop. \ref{VFoldStructureRepresentedByHaefliger})
  \vspace{-2mm}
  $$
    \xymatrix@R=-2pt{
      \mathrm{Integrably}(G,\phi)\mathrm{Structured}V\mathrm{Folds}(\mathbf{H})
      \ar[rr]
      &&
      \left(
        \left(
          \mathbf{H}_{/\mathbf{B}\mathrm{Aut}(T_e V)}
        \right)_{/(\mathbf{B}G , \mathbf{B}\phi )}
      \right)_{
        \!\!\!\big/\left(
          \mathcal{H}\!\mathrm{aef}
          \big(
            V, (G,\phi)
          \big)
          ,
          (\tau_{\mathcal{H}}, g_{\mathcal{H}})
        \right)
      }
      \\
      \big(
        X, (\tau, g)
      \big)
      \ar@{}[rr]|-{\longmapsto}
      &&
      \left(
      \big(
        X, (\tau, g)
      \big)
      \underset{
        \mbox{\tiny{met}}
      }{
        \xrightarrow{ \vdash(\tau,g) }
      }
      \mathcal{H}\!\mathrm{aef}
      \big(
        V, (G,\phi)
      \big)
      \right)
    }
  $$

  \vspace{-2mm}
\noindent
  in the  iterated slice
  of \eqref{IsometriesAsMorphismsInSlice}
  over the $\big( V, (G,\phi)\big)$-Haefliger stack
  (Def. \ref{HaefligerGroupoid}) equipped with its canonical
  $(G,\phi)$-structure $(\tau_{\mathcal{H}}, g_{\mathcal{H}})$
  (Prop. \ref{VFoldStructureRepresentedByHaefliger}),
  hence is \emph{$G$-structure-twisted cohomology} (Remark \ref{TwistedCohomology}):
     \vspace{-2mm}
  \begin{equation}
    \label{EtaleCohomologyWithLocalCoefficients}
    \begin{aligned}
    \overset{
      \mathclap{
      \raisebox{6pt}{
        \tiny
        \color{darkblue}
        \bf
        \begin{tabular}{c}
          {\'e}tale cohomology
        \end{tabular}
      }
      }
    }{
      H^{(\tau,g)}
      \big(
        X, A
      \big)
    }
    & \;:=\;
      \left(
        \left(
          \mathbf{H}_{/\mathbf{B}\mathrm{Aut}(T_e V)}
        \right)_{/(\mathbf{B}G, \mathbf{B}\phi )}
      \right)_{
        \!\!\!\big/
        \left(
          \mathcal{H}\!\mathrm{aef}
          \big(
            V, (G,\phi)
          \big)
          ,
          (\tau_{\mathcal{H}}, g_{\mathcal{H}})
        \right)
      }
    \Big(
      \big(
        X, (\tau,g)
      \big)
      \,,\,
      \big(
        A, p
      \big)
    \Big)
    \\
    & \;=\;
    \left\{\!\!\!\!
    \raisebox{18pt}{
    \xymatrix@C=13pt{
      X
      \ar[rr]^-{
        \overset{
          \mathclap{
          \raisebox{4pt}{
            \tiny
            \color{darkblue}
            \bf
            cocycle
          }
          }
        }{
          c
        }
      }_-{\ }="s"
      \ar[dr]_-{
        \vdash (\tau,g)
      }^-{\ }="t"
      &&
      A
      \ar[dl]^-{p}
      \\
      &
      \mathcal{H}\!\mathrm{aef}(V, (G,\phi))
      \ar@{=>} "s"; "t"
    }
    }
   \!\! \right\}
  \end{aligned}
  \end{equation}
\end{defn}
\begin{remark}[{\'E}tale cohomology is geometric]
  As the notation in Def. \ref{EtaleCohomologyOfVFolds} indicates,
  {\'e}tale cohomology is a ``geometric cohomology theory''
  in that it does depend (in general) on the $G$-structure $g$ on the $V$-fold $X$,
  (for instance its complex- or symplectic- or Riemannian-
  or Lorentzian structure structure, by Example \ref{TorsionFreeGStructureOnSmoothManifolds}).
\end{remark}

Next we consider cohomology theories that
are not sensitive to the metric part $g$ of a $G$-structure $(\tau,g)$,
but just to its tangential structure $\tau$.

\begin{defn}[Tangentially twisted cohomology]
  \label{TangentiallyTwistedCohomology}
  Let $\mathbf{H}$ be an elastic $\infty$-topos  (Def. \ref{ElasticInfinityTopos}),
  $V \in \mathrm{Groups}(\mathbf{H})$
  (Prop. \ref{LoopingAndDelooping}),
  $(G,\phi) \in \mathrm{Groups}(\mathbf{H})_{/_{\mathrm{Aut}(T_e V)}}$
  (Def. \ref{GStructureCoefficients}),
  $(A,\rho) \in G\mathrm{Actions}(\mathbf{H})$
  and $X \in (G,\phi)\mathrm{Structured}V\mathrm{Folds}(\mathbf{H})$
  \eqref{CategoryOfGStructuredVFolds}.
    Then, for $A \in \mathbf{H}_{/_{
    \scalebox{.7}{$
      \mathbf{B}\raisebox{1pt}{\textesh}G
    $}
  }}$, the \emph{tangentially twisted cohomology}
  of $V$ with coefficients in $A$ is (see Remark \ref{TwistedCohomology})
  \vspace{-2mm}
  \begin{equation}
    \label{CohomologyTangentiallyTwisted}
    \overset{
      \mathclap{
      \raisebox{6pt}{
        \tiny
        \color{darkblue}
        \bf
        \begin{tabular}{c}
          tangentially twisted
          \\
          cohomology
        \end{tabular}
      }
      }
    }{
      H^{
        \scalebox{.7}{$
          \raisebox{1pt}{\textesh}
          \tau
        $}
      }
      \big(
        X, A
      \big)
    }
    \;:=\;
    \mathbf{H}_{/_{
      \scalebox{.7}{$
        \raisebox{1pt}{\textesh}\mathrm{Aut}(T_e V)
      $}
    }}
    \big(
      (X, \eta^{\scalebox{.6}{\textesh}} \circ  \tau),
      (A \!\sslash\!G, \rho )
    \big)
    \;=\;
    \left\{\!\!\!\!
    \raisebox{20pt}{
    \xymatrix{
      X
      \ar[rr]^-{
        \overset{
          \mathclap{
          \mbox{
            \tiny
            \color{darkblue}
            \bf
            cocycle
          }
          }
        }{
          c
        }
      }_-{\ }="s"
      \ar[dr]_-{
        \eta^{\scalebox{.5}{\textesh}} \circ \tau
      }^-{\ }="t"
      &&
      (\raisebox{1pt}{\textesh}{A})
        \!\sslash\!
      (\raisebox{1pt}{\textesh}G)
      \ar[dl]^-{
        \scalebox{.7}{$
          \raisebox{1pt}{\textesh}\rho
        $}
      }
      \\
      &
      \mathbf{B}
      \raisebox{1pt}{\textesh}G
      \ar@{=>} "s"; "t"
    }
    }
    \!\! \right\}
  \end{equation}
\end{defn}

\begin{remark}[Need for $G$-Structure vs. tangential structure]
  $\phantom{A}$

  \vspace{-2mm}
\item {\bf (i)}   The notion of tangentially twisted cohomology
  in Definition \ref{TangentiallyTwistedCohomology}
  make sense more generally for $V$-folds
  equipped only with tangential structure (Def. \ref{TangentialStructures})
  instead of full $G$-structure (Def. \ref{GStructures})
  (hence
  only with a reduction of the shape of their structure group,
  instead of the actual structure group (Def. \ref{StructureGroupOfVFolds}))
  and it only need $A$ to be equipped with a $\raisebox{1pt}{\textesh}G$-action.

  \vspace{-1mm}
 \item {\bf (ii)} We state the definition in the more restrictive form above
  just
  in order to bring out the following promotion of this notion to
  its proper orbifold version (Remark \ref{Proper}), in
  Def. \ref{OrbifoldCohomologyTangentiallyTwisted} below.
  The process of orbi-singularization is
  in fact sensitive to the full $G$-structure, and not just to its
  tangential shape. More precisely, it is sensitive to the
  \emph{geometric fixed point spaces} of the $G$-structure
  and not just its homotopy fixed point spaces
  (as per Remark \ref{RootOfProperAsComparedToBorelEquivariantCohomology}
  Example \ref{GeometricFixedPointSpacesDiffer}).
\end{remark}

\newpage

\noindent {\bf Tangentially twisted proper orbifold cohomology.}
We now promote tangentially twisted cohomology of $V$-folds
(Def. \ref{TangentiallyTwistedCohomology}) to a
\emph{proper orbifold cohomology} theory in the sense of
Def. \ref{OrbiVFolds}.
\begin{defn}[Tangentially twisted proper orbifold cohomology]
 \label{OrbifoldCohomologyTangentiallyTwisted}
Let
\begin{itemize}
  \vspace{-2mm}
  \item[$\circ$]
$\mathbf{H}$ be a singular-elastic $\infty$-topos (Def. \ref{SingularSolidInfinityTopos}).

\vspace{-2.5mm}
  \item[$\circ$]
$V \in \mathrm{Groups}(\mathbf{H}_{\tiny\smooth})$ (Prop. \ref{LoopingAndDelooping}).

\vspace{-2.5mm}
  \item[$\circ$]
 $(G,\phi) \in \mathrm{Groups}(\mathbf{H}_{\tiny\smooth})_{/\mathrm{Aut}(T_e V)}$
  (Prop. \ref{LoopingAndDelooping},  Def. \ref{AutomorphismGroup},
    Example \ref{LocalNeighbourhoodOfAPoint}).

\vspace{-2.5mm}
  \item[$\circ$]
$\mathcal{X}_{\tiny\smooth} \in V\mathrm{Folds}(\mathbf{H}_{\tiny\smooth})$
  (Def. \ref{VManifold}).

\vspace{-2.5mm}
  \item[$\circ$]
 $(\tau,g) \in (G,\phi)\mathrm{Structures}_{\mathcal{X}_{\scalebox{.4}{$\smooth$}}}(\mathbf{H}_{\tiny\smooth})$
  (Def. \ref{GStructures}).

\vspace{-2.5mm}
  \item[$\circ$]
$(A, \rho) \in G\mathrm{Actions}(\mathbf{H}_{\tiny\smooth})$.
\end{itemize}

\vspace{-2.5mm}
\noindent and set
 $\mathcal{A} := \orbisingular( A \!\sslash\! G )
  \;\;\;
  \mbox{and}
  \;\;\;
  \mathcal{X} \:= \orbisingular \mathcal{X}_{\tiny\smooth}$.

\noindent
The
\emph{tangentially twisted proper orbifold cohomology}
of $\mathcal{X}$ with coefficients in
$\raisebox{1pt}{\rm\textesh}\mathcal{A}$ is (see Remark \ref{TwistedCohomology})
\vspace{-2mm}
$$
  H^{
    \scalebox{.7}{$ \raisebox{1pt}{\rm\textesh}\orbisingular\tau $}
  }
  \big(
    \mathcal{X}
    \!,\,
    \mathcal{A}
  \big)
  \;:=\;
  \pi_0
  \,
  \mathbf{H}_{
    \!
    /_{
      \!\!
      \scalebox{.7}{$
        \,
        \raisebox{1pt}{\textesh}
        \!
        \orbisingular
        \mathbf{B} G
      $}
    }
  }
  \!
  \Big(
    (\mathcal{X},
    \eta^{\scalebox{.6}{\textesh}}
     \circ
     \orbisingular(\tau)
    )
    ,\,
    (
      \raisebox{1pt}{\textesh} \mathcal{A},
      \raisebox{1pt}{\rm\textesh}\orbisingular \rho
    )
  \Big)
  \;=\;
  \left\{\!\!\!\!
  \raisebox{18pt}{
  \xymatrix{
    \mathcal{X}
    \ar@{-->}[rr]^-{
      \overset{
        \mbox{
          \tiny
          \color{darkblue}
          \bf
          \begin{tabular}{c}
            cocycle
          \end{tabular}
        }
      }{
        c
      }
    }
    \ar[dr]_-{
      \eta^{\scalebox{.6}{\textesh}}
      \circ
      \orbisingular(\tau)
      \!\!\!\!
    }
    &&
    \mbox{\textesh}
    \orbisingular
    \big(
      A \!\sslash\! G
    \big)
    \ar[dl]^{
      \scalebox{.7}{$
        \raisebox{1pt}{\textesh} \orbisingular(\rho)
        \mathrlap{
          \;\;
          \mbox{
            \tiny
            \color{darkblue}
            \bf
          }
        }
      $}
    }
    \\
    &
    \raisebox{1pt}{\textesh}\orbisingular \mathbf{B} G
  }
  }
 \!\! \right\}_{\!\!\big/\sim}
$$

\vspace{-.3cm}

\end{defn}

\begin{theorem}[Tangentially twisted orbifold cohomology at and away from singularities]
  \label{OrbifoldCohomologyTangentiallyTwistedReduces}
Consider the tangentially twisted orbifold cohomology of Def. \ref{OrbifoldCohomologyTangentiallyTwisted}
restriction to {\bf (1)} smooth and {\bf (2)} flat orbifolds,
according to Example \ref{SubcategorySmoothAndFlatOrbifolds}.
Then (see the first diagram on p. \pageref{JTwistedOrbifoldCohomotopyDiagram}):

 \noindent {\bf (i)} The tangentially twisted orbifold cohomology
  of flat orbifolds for 0-truncated coefficients $A$
  is naturally equivalent to the
  proper equivariant cohomology (Def. \ref{ProperEquivariantCohomologyInSingularCohsion})
  of their $\flat G$-frame bundle \eqref{GFrameBundle}:

 \vspace{-2mm}
  $$
\hspace{2cm}
    \overset{
      \mathclap{
      \raisebox{10pt}{
        \tiny
        \color{darkblue}
        \bf
        \begin{tabular}{c}
          tangentially twisted
          \\
          orbifold cohomology
        \end{tabular}
      }
      }
    }{
    H^{\scalebox{.6}{$\raisebox{1pt}{\rm\textesh}\orbisingular \tau$}}
    \big(
      \underset{
        \mathclap{
        \mbox{
          \tiny
          \color{darkblue}
          \bf
          \begin{tabular}{c}
            flat
            \\
            orbifold
          \end{tabular}
        }
        }
      }{
        i_\flat \mathcal{X}
      }
      ,
      \mathcal{A}
    \big)
    }
    \;\;\simeq\;\;
    \overset{
      \mathclap{
      \raisebox{10pt}{
        \tiny
        \color{darkblue}
        \bf
        \begin{tabular}{c}
          proper
          \\
          equivariant cohomology
        \end{tabular}
      }
      }
    }{
    H_{\flat G}
    \big(
      \underset{
        \mathclap{
        \mbox{
          \tiny
          \color{darkblue}
          \bf
          \begin{tabular}{c}
            $\flat G$-frame
            bundle
          \end{tabular}
        }
        }
      }{
        {(\flat G)}\mathrm{Frames}(\mathcal{X}_{\tiny\smooth})
      }
      \,,\,
      A
    \big).
    }
  $$

  \noindent {\bf (ii)}
  The tangentially twisted orbifold cohomology of
  smooth
  (non-orbi-singular) orbifolds is equivalently
  the tangentially twisted cohomology (Def. \ref{EtaleCohomologyOfVFolds})
  of the underlying $V$-folds:
  \vspace{-2mm}
  $$
    \overset{
      \mathclap{
      \raisebox{10pt}{
        \tiny
        \color{darkblue}
        \bf
        \begin{tabular}{c}
          tangentially twisted
          \\
          orbifold cohomology
        \end{tabular}
      }
      }
    }{
    H^{\scalebox{.6}{$\raisebox{1pt}{\rm\textesh}\orbisingular \tau$}}
    \big(
      \underset{
        \mathclap{
        \mbox{
          \tiny
          \color{darkblue}
          \bf
          \begin{tabular}{c}
            smooth
            \\
            orbifold
          \end{tabular}
        }
        }
      }{
        i_{\scalebox{.5}{$\smooth$}} \mathcal{X}
      }
      ,
      \mathcal{A}
    \big)
    }
    \;\;\simeq\;\;
    \overset{
      \mathclap{
      \raisebox{10pt}{
        \tiny
        \color{darkblue}
        \bf
        \begin{tabular}{c}
          tangentially twisted
          \\
          $V$-fold cohomology
        \end{tabular}
      }
      }
    }{
    H^{ \scalebox{.6}{$\raisebox{1pt}{\rm\textesh} \tau$} }
    \big(
      \underset{
        \mathclap{
        \raisebox{-5.5pt}{
          \tiny
          \color{darkblue}
          \bf
          \begin{tabular}{c}
            0-truncated
            \\
            $V$-fold
          \end{tabular}
        }
        }
      }{
        \mathcal{X}_{\tiny\smooth}
      }
      \,,\,
      A
    \big).
    }
  $$

\end{theorem}

\begin{proof}

The case {\bf (i)} means that the classifying map of the $G$-structure
 in question factors as follows, where we use
 Prop. \ref{GStructuredVFoldAsQuotientOfGFrameBundle}
 to identify the leftmost morphism $\rho$ as exhibiting the
 action \eqref{InfinityActionHomotopyFiberSequence} of
 $\flat G$ on $(\flat G)\mathrm{Frames}(X)$:
\vspace{-2mm}
$$
  \xymatrix@C=40pt{
    (\flat G)\mathrm{Frames}(X)
    \!\sslash\!
    \flat G
    \ar[r]^-{\rho}
    \ar@/_.8pc/[rr]_>>>>>>>{
      \vdash \tau
    }
    \ar@/_1.9pc/[rrr]_{
      \vdash \mathrm{Frames}(X)
    }
    &
    \mathbf{B} \flat G
    \ar[r]^{ \epsilon^\flat_{\mathbf{B}G} }
    &
    \mathbf{B}G
    \ar[r]
    &
    \mathbf{B} \mathrm{Aut}(T_e V)
    \,.
  }
$$

 \vspace{-2mm}
\noindent
Now we observe
\begin{itemize}
 \vspace{-.2cm}
 \item[{\bf (a)}]
 with Def. \ref{SingularCohesiveInfinityTopos}
 that \raisebox{1pt}{\textesh} acts objectwise over $\mathrm{Singularities}$,
 \vspace{-.2cm}
 \item[{\bf (b)}]
 with Prop. \ref{LimitsAndColimitsOfPresheavesComputedObjectwise} that
 the pullback of presheaves over $\mathrm{Singularities}$ is computed
 objectwise,
 \vspace{-.2cm}
 \item[{\bf (c)}]
and with Lemma \ref{ShapeOfOrbiSingularSpacesAsPresheafOnSingularities}
that $\orbisingular (A \!\sslash\! G)$ is objectwise over
$\mathrm{Singularities}$ a homotopy quotient by $G$,
\end{itemize}

 \vspace{-2mm}
\noindent
so that Lemma \ref{ShapeOfAnInducedAction} applies objectwise
over $\mathrm{Singularities}$ to give the pullback square
shown on the right here:

\newpage
\vspace{-2mm}
$$
  \xymatrix@C=5em@R=1.5em{
    &
    \raisebox{1pt}{\textesh}
    \orbisingular
    \big(
      A \!\sslash\! \flat G
    \big)
    \ar@{}[dr]|-{
      \mbox{
        \tiny
        (pb)
      }
    }
    \ar[d]
    \ar[r]
    &
    \raisebox{1pt}{\textesh}
    \orbisingular
    \big(
      A \!\sslash\! G
    \big)
    \ar[d]
    \\
    \orbisingular
    \big(
      (\flat G)\mathrm{Frames}(X) \!\sslash\! (\flat G)
    \big)
    \ar@{-->}[ur]
    \ar[r]|-{ \orbisingular \rho }
    \ar@/_1.5pc/[rr]_-{
      \eta^{\scalebox{.6}{\textesh}}
      \circ
      \orbisingular\tau
    }
    &
    \orbisingular \mathbf{B}\flat G
    \ar[r]
    &
    \raisebox{1pt}{\textesh}
    \orbisingular
    \mathbf{B} G
  }
$$

 \vspace{-1mm}
\noindent
By the universal property of the pullback, this means that
every cocycle factors naturally as shown by the dashed morphism.
But by Def. \ref{ProperEquivariantCohomologyInSingularCohsion} this
dashed morphism is equivalently a cocycle in proper equivariant cohomology,
as claimed.

\noindent The case {\bf (ii)} means
(using Lemma \ref{Smooth0TruncatedObjectsAreOrbiSingular}) that
the orbi-singular space $\mathcal{X}$ is in fact smooth
 \vspace{-2mm}
$$
  \mathcal{X} \;\simeq\; \smooth \mathcal{X}
  \;\simeq\;
  \mathrm{NnOrbSnglr}
  \big(
    \mathcal{X}_{\smooth}
  \big)
  \,.
$$

 \vspace{-2mm}
\noindent
Therefore, we have the following natural equivalences of
spaces of dashed morphisms:
 \vspace{-2mm}
\begin{equation}
  \label{JTwistedOrbifoldCocycleOnSmoothSpace}
  \begin{array}{ccccc}
  \raisebox{20pt}{
  \xymatrix{
    &
    \raisebox{1pt}{\textesh}\orbisingular
    \big(
      A \!\!\sslash\! G
    \big)
    \ar[d]
    \\
    \mathllap{
      \smooth \mathcal{X}
      \simeq
      \;
    }
    \mathcal{X}
    \ar@{-->}[ur]
    \ar[r]_{  }
    &
    \raisebox{1pt}{\textesh} \orbisingular \mathbf{B} G
  }
  }
  &
  \;\;
  \Leftrightarrow
  \;\;
  &
  \raisebox{20pt}{
  \xymatrix{
    &
    \smooth \raisebox{1pt}{\textesh}\orbisingular
    \big(
      A \!\!\sslash\! G
    \big)
    \ar[d]
    \\
    \smooth \mathcal{X}
    \ar@{-->}[ur]
    \ar[r]_{  }
    &
    \smooth \raisebox{1pt}{\textesh} \orbisingular \mathbf{B} G
  }
  }
  &
  \;\;
  \Leftrightarrow
  \;\;
  &
  \raisebox{20pt}{
  \xymatrix{
    &
    \raisebox{1pt}{\textesh}
    \big(
      A \!\!\sslash\! G
    \big)
    \ar[d]
    \\
    \mathcal{X}_{\tiny\smooth}
    \ar@{-->}[ur]
    \ar[r]
    &
    \raisebox{1pt}{\textesh} \mathbf{B} G
  }
  }
  \\
  \begin{rotate}{-90}
    $\mathclap{\in}$
  \end{rotate}
  && &&
  \begin{rotate}{-90}
    $\mathclap{\in}$
  \end{rotate}
  \\
  \mathclap{
  \mathbf{H}_{
    \big/
    \scalebox{.7}{
      $
      \raisebox{1pt}{\textesh}\orbisingular
      \mathbf{B} G
      $
    }
  }
  \Big(
    \mathcal{X},
    \raisebox{1pt}{\textesh}
    \orbisingular
    \big(
      A \!\!\sslash\! G
    \big)
  \Big)
  }
  && \simeq &&
  \mathclap{
  \mathbf{H}_{\big/
    \scalebox{.7}{
      $
      \raisebox{1pt}{\textesh}
      \mathbf{B} G
      $
    }
  }
  \Big(
    \mathcal{X}_{\tiny\smooth}
    \,,\,
    \raisebox{1pt}{\textesh}
    \big(
      A \!\!\sslash\! G
    \big)
  \Big)
  }
\end{array}
\end{equation}

 \vspace{-2mm}
\noindent
Here the first equivalence is by the adjunction
$\mathrm{NnOrbSnglr} \dashv \mathrm{Smth}$
and the fully faithfulness of $\mathrm{NnOrbSnglr}$ \eqref{CohesionOfSingularCohesive}.
The second step uses
$\smooth \circ \raisebox{1pt}{\textesh} \simeq \mbox{\textesh} \circ \smooth$
(Lemma \ref{SmoothCommutesWithShape}) and
$\smooth \circ \orbisingular \simeq \smooth$
(Remark \ref{SmoothOrbiSingularIsSmooth})
But on the right of \eqref{JTwistedOrbifoldCocycleOnSmoothSpace}
we see the tangentially twisted cohomology of $\mathcal{X}_{\tiny\smooth}$,
as claimed.
\hfill \end{proof}

\medskip

\noindent {\bf J-Twisted orbifold Cohomotopy theory.}
We discuss now the example of tangentially twisted proper orbifold
cohomology (Def. \ref{OrbifoldCohomologyTangentiallyTwisted})
where the coefficients are (shapes of) spheres, specifically
of \emph{Tate $V$-spheres} (Def. \ref{TheVSphere}),
In this case the tangential twist
is the \emph{J-homomorphism} (Def. \ref{TateJHomomorphism}) whence we speak of
\emph{J-twisted Cohomotopy theory} (Def. \ref{JTwistedOrbifoldCohomotopyTheory}).

\begin{defn}[Complement of neutral element]
  \label{ComplemenetOfNeutralElement}
  Let $\mathbf{H}$ be an $\infty$-topos (Def. \ref{InfinityTopos})
  and $V \in \mathrm{Groups}(\mathbf{H})$ (Prop. \ref{LoopingAndDelooping}).
  Let $(V, \rho_{\mathrm{Aut}_{\mathrm{Grp}}}) \in
  \mathrm{Aut}_{\mathrm{Grp}}(V)\mathrm{Actions}(\mathbf{H})$
  denote
  the group-automorphism action on $V$ (Prop. \ref{CanonicalActionOfGroupAutomomorphisms}).

\noindent {\bf (i)}
 Consider those subobjects (Def. \ref{PosetOfSubobjects}) of the
  homotopy quotient $V \sslash \mathrm{Aut}_{\mathrm{Grp}}$
  \eqref{HomotopyQuotientOfGByGroupAutomorphisms}
  whose pullback along the morphism
   \vspace{-2mm}
  $$
    \xymatrix{
      \ast \!\sslash\! \mathrm{Aut}_{\mathrm{Grp}}(V)
      \ar[rr]^-{\scalebox{.7}{$
        e \!\sslash\! \mathrm{Aut}_{\mathrm{Grp}}(V)
        $}
      }
      &&
      V \!\sslash\! \mathrm{Aut}_{\mathrm{Grp}}(V)
    },
  $$

   \vspace{-1mm}
\noindent
  which exhibits the neutral element as a fixed point of the
  group-automorphism action (Prop. \ref{CanonicalActionOfGroupAutomomorphisms}),
  is empty. These are the subobjects forming the poset
  in the top left of the following
  Cartesian square (of $\infty$-categories):

 \vspace{-3mm}
  \begin{equation}
    \label{SubobjectsNotIncludingTheNeutralElement}
    \xymatrix@R=1.5em{
      \mathrm{SubObjects}_{e\!\!\!\!/}
      \big(
        V \!\sslash\! \mathrm{Aut}_{\mathrm{Grp}}(V)
      \big)
      \ar@{}[rrd]|-{ \mbox{\tiny\rm(pb)} }
      \ar[rr]
      \ar@{^{(}->}[d]
      && \ast
      \ar@{^{(}->}[d]^-{ \varnothing }
      \\
      \mathrm{SubObjects}
      \big(
        V \!\sslash\! \mathrm{Aut}_{\mathrm{Grp}}(V)
      \big)
      \ar[rr]_-{\scalebox{.7}{$
        (e \!\sslash\! \mathrm{Aut}_{\mathrm{Grp}}(V))^\ast
        $}
      }
      &&
      \mathrm{SubObjects}(\ast)
    }
  \end{equation}
\noindent   {\bf (ii)} Consider next
  the union of these subobjects, hence the colimit over the
  left vertical functor in \eqref{SubobjectsNotIncludingTheNeutralElement},
  which we denote as follows:
   \vspace{-2mm}
  \begin{equation}
    \label{HomotopyQuotientOfComplementOfNeutralElement}
    \big(
      V \setminus \{e\}
    \big)
    \!\sslash\! \mathrm{Aut}_{\mathrm{Grp}}(V)
    \;:=\;
    \underset{
      \longrightarrow
    }{\mathrm{lim}}
    \Big(
      \mathrm{SubObjects}_{e\!\!\!\!/}
      \big(
        V \!\sslash\! \mathrm{Aut}_{\mathrm{Grp}}(V)
      \big)
      \longhookrightarrow
      \mathrm{SubObjects}
      \big(
        V \!\sslash\! \mathrm{Aut}_{\mathrm{Grp}}(V)
      \big)
    \Big).
  \end{equation}

 \vspace{-2mm}
  \noindent {\bf (iii)} We call  the homotopy fiber $V \setminus \{e\}$
  of the canonical morphism from this
  object \eqref{HomotopyQuotientOfComplementOfNeutralElement} to
  $\mathbf{B}\mathrm{Aut}_{\mathrm{Grp}}(G)$
  the \emph{complement of the neutral element of $V$}
  \vspace{-.3cm}
  \begin{equation}
    \label{ComplementOfNeutralElementAsHomotopyFiber}
    \raisebox{40pt}{
    \xymatrix@C=3em@R=1.5em{
      V \setminus \{e\}
      \ar[rr]^-{\scalebox{.7}{$
        \mathrm{fib}
        \big(
          \rho_{\mathrm{Aut}_{\mathrm{Grp}}} \setminus \{e\}
        \big)
        $}
      }
      &&
      (V \setminus \{e\})
      \!\sslash\!
      \mathrm{Aut}_{\mathrm{Grp}}(V)
      \ar@/^3.2pc/[dd]^-{\scalebox{.7}{$
        \rho_{ \mathrm{Aut}_{\mathrm{Grp}} } \setminus \{e\}
        $}
      }
      \ar@{^{(}->}[d]
      \\
      &&
      V \!\sslash\! \mathrm{Aut}_{\mathrm{Grp}}(V)
      \ar[d]_-{ \scalebox{.7}{$\rho_{\mathrm{Aut}_{\mathrm{Grp}}} $}}
      \\
      &&
      \mathbf{B}\mathrm{Aut}_{\mathrm{Grp}}(V)
    }
    }
  \end{equation}

\vspace{-2mm}
  \noindent
  {\bf (iv)}
  We regard the complement of the neutral element as equipped with
  the $\mathrm{Aut}_{\mathrm{Grp}}(V)$-action
  which is exhibited by
  the homotopy fiber sequence \eqref{ComplementOfNeutralElementAsHomotopyFiber}
  (by Prop. \ref{ComplementOfNeutralElementAsHomotopyFiber}):
  $$
    \big(
      V \setminus \{e\}
      \,,\,
      \rho_{\mathrm{Aut}_{\mathrm{Grp}}} \setminus \{E\}
    \big))
    \;\in\;
    \mathrm{Aut}_{\mathrm{Grp}}(V)\mathrm{Actions}(\mathbf{H})
    \,.
  $$
\end{defn}

\begin{prop}[Basic properties of complement of neutral element]
  \label{BasicPropertiesOfComplementOfNeutralElement}
  Let $\mathbf{H}$ be an $\infty$-topos (Def. \ref{InfinityTopos})
  and $V \in \mathrm{Groups}(\mathbf{H})$ (Prop. \ref{LoopingAndDelooping}).
  Then the complement $V \setminus \{e\}$ of the neutral element
  (Def. \ref{ComplemenetOfNeutralElement})

  \noindent {\bf (i)} is a subobject (Def. \ref{PosetOfSubobjects}) of $V$
  \vspace{-1mm}
  \begin{equation}
    \label{ComplementOfNeutralElementAsSubobject}
    \xymatrix@C=3em@R=1em{
      V \setminus \{e\}
     \; \ar@{^{(}->}[r]
      &
      V
    }
  \end{equation}
  \vspace{-.2cm}
  \noindent {\bf (ii)} which is disjoint from the neutral element:
 \vspace{-2mm}
$$
  \xymatrix@C=4.5em@R=1.2em{
    \varnothing
    \ar[r]
    \ar[d]
    \ar@{}[dr]|-{ \mbox{\tiny\rm(pb)} }
    &
    V \setminus \{ e \}
    \ar@{^{(}->}[d]
    \\
    \ast
    \ar[r]_-{ e }
    &
    V
  }
$$
\end{prop}
\begin{proof}
  For {\bf (i)} we use the pasting law (Prop. \ref{PastingLaw})
  and the homotopy fiber characterization of the
  group-automorphism action \eqref{TowardsUnderstandingGroupAutomorphismActionOnG}
  to decompose \eqref{ComplementOfNeutralElementAsHomotopyFiber} as the
  pasting of two Cartesian squares, as follows:
 \vspace{-1mm}
 $$
    \xymatrix@R=1.8em{
      V \setminus \{e\}
      \ar[rr]
      \ar@{^{(}->}[d]
      \ar@{}[rrd]|-{ \mbox{\tiny\rm(pb)} }
      &&
      (V \setminus \{e\})
      \!\sslash\!
      \mathrm{Aut}_{\mathrm{Grp}}(V)
      \ar@/^3.2pc/[dd]^-{
        \rho_{ \mathrm{Aut}_{\mathrm{Grp}} } \setminus \{e\}
      }
      \ar@{^{(}->}[d]
      \\
      V \ar[rr]
      \ar[d]
      \ar@{}[rrd]|-{ \mbox{\tiny\rm(pb)} }
      &&
      V \!\sslash\! \mathrm{Aut}_{\mathrm{Grp}}(V)
      \ar[d]_-{ \rho_{\mathrm{Aut}_{\mathrm{Grp}}} }
      \\
      \ast
      \ar[rr]
      &&
      \mathbf{B}\mathrm{Aut}_{\mathrm{Grp}}(V)
          }
  $$

 \vspace{-2mm}
\noindent
  Since monomorphisms are preserved by pullback
  (by Prop. \ref{nConnectednTruncatedFactorizationSystem}), this
  shows the first claim from the construction \eqref{HomotopyQuotientOfComplementOfNeutralElement}.

  For {\bf (ii)} we paste to the middle horizontal morphism in
  this diagram the square \eqref{HomotopyFixedPoints}
  which exhibits the neutral element as a fixed point of the
  group-automorphims action (Prop. \ref{CanonicalActionOfGroupAutomomorphisms})
  and then we pull back the right vertical morphism along the boundary
  of that square, as shown in the following:
\vspace{-2mm}
  $$
    \xymatrix@R=6pt{
      \varnothing
      \ar[dd]
      \ar[dr]
      \ar[rr]
      &&
      \varnothing
      \mathrlap{
        \;\;\;\;\;\;
        \simeq\;\;\;\;\;\;
        \underset{\underset{i}{\longrightarrow}}{\lim} \, \varnothing
      }
      \ar[dd]|<<<<<<{\phantom{A\vert^{\vert}}}
      \ar[dr]
      \\
      &
      V \setminus \{e\}
      \ar[rr]
      \ar@{^{(}->}[dd]
      &&
      (V \setminus \{e\})
      \!\sslash\!
      \mathrm{Aut}_{\mathrm{Grp}}(V)
      \mathrlap{
        \;\;
        :=
        \;
        \underset{\underset{i}{\longrightarrow}}{\lim} \, U_i
      }
      \ar@{^{(}->}[dd]
      \\
      \ast
      \ar[dr]_-{ e }
      \ar[rr]|-{\phantom{AA}}
      &&
      \ast \!\sslash\! \mathrm{Aut}_{\mathrm{Grp}}(V)
      \ar[dr]^-{\scalebox{.6}{$ e \sslash \mathrm{Aut}_{\mathrm{Grp}}(V) $}}
      \\
      &
      V \ar[rr]
      &&
      V \!\sslash\! \mathrm{Aut}_{\mathrm{Grp}}(V)
    }
  $$

  \vspace{-1mm}
\noindent
  Here the right square is Cartesian
  since colimits in an $\infty$-topos are preserved by pullback
  \eqref{PreserveColimitsByPullback}
  and using the definition \eqref{SubobjectsNotIncludingTheNeutralElement},
  as indicated in the top right. Similarly the rear square is Cartesian,
  since pullback preserves the initial object
  (this being the empty colimit, Example \ref{InitialObjectInInfinityToposIsEmpty}).
  With this, and since the front square is Cartesian by {\bf (i)}, the
  pasting law (Prop. \ref{PastingLaw}) implies that also the
  left square is Cartesian, which was to be shown.
\hfill \end{proof}

\newpage

\begin{defn}[Tate $V$-sphere]
  \label{TheVSphere}
  Let $\mathbf{H}$ be an $\infty$-topos (Def. \ref{InfinityTopos})
  and $V \in \mathrm{Groups}(\mathbf{H})$ (Prop. \ref{LoopingAndDelooping}).
  Then we say that the \emph{Tate $V$-sphere} is the
  homotopy cofiber
   \vspace{-2mm}
  $$
    \mathbf{S}^V \;:=\; V / (V \setminus \{e\})
  $$

   \vspace{-2mm}
\noindent
  of the inclusion
  \eqref{ComplementOfNeutralElementAsSubobject}
  of the complement of the neutral element
  into $V$ (Def. \ref{ComplemenetOfNeutralElement}), hence
  the object in this homotopy pushout square:
  \begin{equation}
    \label{VSphere}
    \raisebox{10pt}{
    \xymatrix@C=3em@R=1.2em{
      V \setminus \{e\}
      \;\ar@{^{(}->}[rr]
      \ar[d]
      \ar@{}[drr]|-{ \mbox{\tiny(po)} }
      &&
      V
      \ar[d]
      \\
      \ast
      \ar[rr]
      &&
      \mathbf{S}^V
    }
    }
  \end{equation}
\end{defn}
\begin{example}[Tate sphere in unstable motivic homotopy theory]
  \label{TateSphereInMotivicHomotopyTheory}
  For $\mathbf{H} \;:=\; \mathrm{Sheaves}_\infty\big(\mathrm{Schemes}_{{}_{\mathrm{Nis}}}\big)$
  and $V \;:=\; \mathbb{A}^1$ the Tate $V$-sphere
  of Def. \ref{TheVSphere} is the Tate sphere in the traditional
  sense of (unstable) motivic homotopy theory,
  see \cite[2.22]{VRO07}.
\end{example}

\begin{example}[Tate spheres with shape of ordinary spheres]
  \label{TateSpheresWithShapeOfOrdinarySpheres}
  Let $\mathbf{H} = \mathrm{JetsOfSmoothGroupoids}_\infty$
  (Def. \ref{FormalSmoothInfinityGroupoids})
  and $V := (\mathbb{R}^n,+)$ as in Example \ref{OrdinaryManifolds}.
  Then $\mathrm{Aut}_{\mathrm{Grp}}(\mathbb{R}^n,+) = \mathrm{GL}(n)$
  (as in Example \ref{OrdinaryGeneralLinearGroup}) and
  the complement of the neutral element (Def. \ref{ComplemenetOfNeutralElement})
  is the ordinary complement $\mathbb{R}^n \setminus \{0\}$,
  whose shape is that or the ordinary $n-1$-sphere:
  \begin{equation}
    \label{ShapeOfComplementOfNeutralElementInRn}
    \raisebox{1pt}{\textesh}
    \,
    \big(
      \mathbb{R}^n \setminus \{0\}
    \big)
    \;\simeq\;
    \raisebox{1pt}{\textesh}\, S^{n-1}.
  \end{equation}

  \vspace{-2mm}
\noindent
  Hence the Tate $\mathbb{R}^n$-sphere (Def. \ref{TheVSphere})
  is the homotopy pushout shown on the left here:
  \vspace{-2mm}
  $$
    \raisebox{20pt}{
    \xymatrix@R=1.5em{
      \mathbb{R}^n \setminus \{e\}
      \ar[d]
      \ar@{^{(}->}[rr]
      \ar@{}[drr]|-{ \mbox{\tiny\rm(pb)} }
      &&
      \mathbb{R}^n
      \ar[d]
      \\
      \ast
      \ar[rr]
      &&
      S^{(\mathbb{R}^n)}
    }
    }
    \;\;\;\;
    \overset{
      \raisebox{1pt}{\textesh}
    }{\longmapsto}
    \;\;\;\;
    \raisebox{20pt}{
    \xymatrix@R=1.5em{
      \raisebox{1pt}{\textesh} \, S^{n-1}
      \ar[d]
      \ar@{^{(}->}[rr]
      \ar@{}[drr]|-{ \mbox{\tiny\rm(pb)} }
      &&
      \ast
      \ar[d]
      \\
      \ast
      \ar[rr]
      &&
      \raisebox{1pt}{\textesh}\,\mathbf{S}^{(\mathbb{R}^n)}
    }
    }
  $$
  Since the shape modality \eqref{CohesiveModalitiesFromAdjointQuadruple}
  is left adjoint it preserves homotopy pushouts (Prop. \ref{AdjointsPreserveCoLimits}),
  so that the shape of the
  Tate $\mathbb{R}^n$-sphere is that of the ordinary $n$-sphere:
  \begin{equation}
    \label{ShapeOfTatenSphereIsShapeOfNSphere}
    \raisebox{1pt}{\textesh}
    \,\mathbf{S}^{\mathbb{R}^n}
    \;\simeq\;
    \raisebox{1pt}{\textesh}\,S^n.
  \end{equation}
  In contrast, the Tate $\mathbb{R}^n$-sphere itself is the
  ``germ of a smooth sphere''.
\end{example}

\medskip
\begin{prop}[Canonical action on Tate $V$-sphere]
  \label{TateVSphereCanonicalAction}
  Let $\mathbf{H}$ be an $\infty$-topos (Def. \ref{InfinityTopos})
  and $V \in \mathrm{Groups}(\mathbf{H})$ (Prop. \ref{LoopingAndDelooping}).
  The Tate $V$-sphere (Def. \ref{TheVSphere})
  inherits a canonical
  action (Prop. \ref{InfinityAction})
  of the group-automorphism group
  $\mathrm{Aut}_{\mathrm{Grp}}(V)$ (Def. \ref{GroupOFGroupAutomorphisms}),
  associated (via Prop. \ref{AutomorphimsGroupIsUniversal}) to a
  group homomorphism
   \vspace{-2mm}
  \begin{equation}
    \label{ActionOnVSphereByGroupHomomorphism}
    \xymatrix{
      \mathrm{Aut}_{\mathrm{Grp}}(V)
      \ar[r]
      &
      \mathrm{Aut}(\mathbf{S}^V)
    }
  \end{equation}

   \vspace{-2mm}
\noindent
  whose homotopy quotient \eqref{HomotopyQuotientAsColimit}
  is given by
  the following homotopy pushout
  \begin{equation}
    \label{VSphere}
    \xymatrix@C=3em@R=1.2em{
      (V \setminus \{e\})
        \!\sslash\! \mathrm{Aut}_{\mathrm{Grp}}(V)
      \; \ar@{^{(}->}[rr]
      \ar[d]
      \ar@{}[drr]|-{ \mbox{\tiny\rm(po)} }
      &&
      V \!\sslash\! \mathrm{Aut}_{\mathrm{Grp}}(V)
      \ar[d]
      \\
      \ast
      \!\sslash\! \mathrm{Aut}_{\mathrm{Grp}}(V)
      \ar[rr]
      &&
      \mathbf{S}^V \!\sslash\! \mathrm{Aut}_{\mathrm{Grp}}(V)
    }
  \end{equation}
  of the defining morphisms in \eqref{ComplementOfNeutralElementAsHomotopyFiber}.
\end{prop}
\begin{proof}
  Since the forgetful $\infty$-functor
  $\mathbf{H}_{/\mathbf{B}\mathrm{Aut}_{\mathrm{Grp}}(V)}
  \longrightarrow \mathbf{H}$ preserbes colimits (Example \ref{BaseChangeAlongTerminalMorphism}),
  the diagram \eqref{VSphere} extends to a diagram over
  $\mathbf{B}\mathrm{Aut}_{\mathrm{Grp}}(V)$. Pulling this
  back along the point inclusion \eqref{PointInclusionIntoDelooping}
  and using that colimits in an $\infty$-topos are preserved by pullback
  \eqref{PreserveColimitsByPullback}, we find that the homotopy fiber of
  $\mathbf{S}^V \!\sslash\! \mathrm{Aut}_{\mathrm{Grp}}(V)
  \to \mathbf{B}\mathrm{Aut}_{\mathrm{Grp}}(V)$
  is given by the defining homotopy pushout \eqref{VSphere}
  of the Tate $V$-sphere.
\hfill \end{proof}

\begin{defn}[Linear group]
  \label{LinearGroup}
  Let $\mathbf{H}$ be an
  elastic $\infty$-topos (Def. \ref{ElasticInfinityTopos})
  and $V \in \mathrm{Groups}(\mathbf{H})$ (Prop. \ref{LoopingAndDelooping}).

  \noindent {\bf (i)}
    We say that $V$ is a \emph{linear group}
  if it is equipped with an equivalence
  \begin{equation}
    \label{DefiningEquivalenceForLinearGroup}
    \xymatrix{
      \mathrm{Aut}(T_e V)
      \ar[rr]_-{\simeq}^{ \mathrm{exp} }
      &&
      \mathrm{Aut}_{\mathrm{Grp}}(V)
    }
    \;\;
    \in\
    \mathrm{Groups}(\mathbf{H})
  \end{equation}
  between {\bf (a)} the plain automorphism group of the
  local neighborhood of the neutral element (Def. \ref{StructureGroupOfVFolds})
  and {\bf (b)} the group-automorphism group of $V$ (Def. \ref{GroupOFGroupAutomorphisms})

 \noindent {\bf (ii)}   We write
  $$
    \mathrm{LinearGroups}(\mathbf{H})
    \in
    \mathrm{Categories}_\infty
  $$
  for the $\infty$-category of linear groups in $\mathbf{H}$.
\end{defn}

\begin{defn}[Tate J-homomorphism]
  \label{TateJHomomorphism}
  Let $\mathbf{H}$ be an
  elastic $\infty$-topos (Def. \ref{ElasticInfinityTopos})
  and $V \in \mathrm{LinearGroups}(\mathbf{H})$ (Prop. \ref{LinearGroup}).

\noindent {\bf (i)}  The \emph{Tate J-homomorphism} is the composite
 \vspace{-2mm}
  \begin{equation}
    \label{TheTateJHomomorphism}
    \mathbf{J}_V
    \;:\;
    \xymatrix{
      \mathrm{Aut}(T_e V)
      \ar[r]_-{\simeq}^{ \mathrm{exp} }
      &
      \mathrm{Aut}_{\mathrm{Grp}}(V)
      \ar[r]
      &
      \mathrm{Aut}(\mathbf{S}^V)
    }
  \end{equation}

   \vspace{-2mm}
\noindent
  of {\bf (a)} the defining equivalence \eqref{DefiningEquivalenceForLinearGroup}
  with {\bf (b)}
  the homomorphism \eqref{ActionOnVSphereByGroupHomomorphism}
  which reflects the canonical $\mathrm{Aut}_{\mathrm{Grp}}(V)$-action
  on the Tate $V$-sphere (Def. \ref{TateVSphereCanonicalAction}).

\noindent {\bf (ii)}
 The corresponding $\mathrm{Aut}(T_e V)$-actions on $\mathbf{S}^V$
 and on $\raisebox{1pt}{\textesh}(\mathbf{S}^V)$,
 by restriction
 along \eqref{TheTateJHomomorphism} and \eqref{JHomomorphism}
 of the canonical automorphism actions (Prop. \ref{AutomorphimsGroupIsUniversal}),
 we denote, respectively, by
  \begin{equation}
    \label{ActionOfJHomomorphism}
    (\mathbf{S}^V, \rho_{\mathbf{J}})
    \;
    \;\in\;
    \mathrm{Aut}(T_e V)\mathrm{Actions}(\mathbf{H})\;.
  \end{equation}

\vspace{-2mm}
 \noindent {\bf (iii)} The actual \emph{J-homomorphism}
  is the shape of the further composite with the homomorphism
  $\mathrm{Aut}(\eta^{\scalebox{.7}{\textesh}}_{\mathbf{S}^V})$
  from Prop. \ref{AutomorphismsAlongShapeUnit}:
  \vspace{-2mm}
  \begin{equation}
    \label{JHomomorphism}
    J_V
    \;:\;
    \xymatrix{
      \raisebox{1pt}{\textesh}\mathrm{Aut}(T_e V)
      \ar[r]_-{\simeq}^{
        \scalebox{.6}{$
          \raisebox{1pt}{\textesh}\mathrm{exp}
        $}
      }
      \ar@/_1.5pc/[rr]_-{
        \scalebox{.6}{$
          \raisebox{1pt}{\textesh}\mathbf{J}_V
        $}
      }
      &
      \raisebox{1pt}{\textesh}
      \mathrm{Aut}_{\mathrm{Grp}}(V)
      \ar[r]
      &
      \raisebox{1pt}{\textesh}
      \mathrm{Aut}(\mathbf{S}^V)
      \ar[rr]^-{
        \scalebox{.6}{$
          \raisebox{1pt}{\textesh}
          \mathrm{Aut}
          \big(
            \eta^{\scalebox{.7}{\textesh}}_{\mathrm{Aut}(\mathbf{S}^V)}
          \big)
        $}
      }
      &&
      \raisebox{1pt}{\textesh}
      \mathrm{Aut}
      \big(
        \raisebox{1pt}{\textesh} \mathbf{S}^V
      \big).
    }
  \end{equation}

\end{defn}

\begin{example}[Ordinary J-homomorphism]
  \label{OrdinaryJHomomorphism}
  Let $\mathbf{H} = \mathrm{SingularJetsOfSmoothGroupoids}_\infty$
  (Example \ref{OrbiSingularSmoothInfinityGroupoids}) and
  $V := (\mathbb{R}^n, +)$ as in Example \ref{OrdinaryManifolds}.
  This is a linear group in the sense of Def. \ref{LinearGroup},
  with $\mathrm{Aut}(T_0 \mathbb{R}^n) \simeq \mathrm{GL}(n)$
  (Example \ref{OrdinaryGeneralLinearGroup}).
  Via Example \ref{TateSpheresWithShapeOfOrdinarySpheres} the
  induced action
  on the shape of the Tate $\mathbb{R}^n$-sphere (Def. \ref{TateJHomomorphism})
  is the classical J-homomorphism
  (going back to \cite{Whitehead42},
  reviewed in \cite[p. 4]{Ravenel86}):
  \vspace{-2mm}
  \begin{equation}
    \label{OrdinaryJHomomorphism}
    J
    \;:\;
    \xymatrix{
      \raisebox{1pt}{\textesh}\, \mathrm{O}(n)
      \;\simeq\;
      \raisebox{1pt}{\textesh}\, \mathrm{GL}(n)
      \ar[r]
      &
      \mathrm{Aut}
      \big(
        \raisebox{1pt}{\textesh} S^n
      \big)
    }
  \end{equation}

   \vspace{-2mm}
\noindent
  being the image
  under topological shape (Def. \ref{ShapeOfTopologicalSpaces})
  of the defining action of $\mathrm{GL}(n)$
  on $\mathbb{R}^n$ and hence on its one-point compactification
  $S^n$.
\end{example}

\begin{defn}[Representation spheres]
  \label{RepresentationSphere}
  Let $\mathbf{H}$ be a singular-elastic $\infty$-topos (Def. \ref{SingularSolidInfinityTopos}),
  $V \mathrm{Groups}(\mathbf{H}_{\tiny\smooth})$ (Prop. \ref{LoopingAndDelooping}),
  and
  $(G, \phi) \in \mathrm{Groups}(\mathbf{H}_{\tiny\smooth})_{/ \mathrm{Aut}(T_e V) }$
  (Prop. \ref{LoopingAndDelooping},  Def. \ref{AutomorphismGroup},
    Example \ref{LocalNeighbourhoodOfAPoint}).
     Then we say that the \emph{representation sphere}
  $S^{V_\phi}$
  of the $G$-action $\phi$  on $V$ (via Prop. \ref{AutomorphimsGroupIsUniversal})
  is the shape (Def. \ref{CohesiveTopos}) of the
  orbi-singularization (Def. \ref{SincularCohesionModalies}) of the
  homotopy quotient \eqref{HomotopyQuotientAsColimit}
  of the Tate $V$-sphere (Def. \ref{TheVSphere})
  by the restricted action (Prop. \ref{PullbackAction}) along $\phi$
  of the action $\rho_{\mathbf{J}}$
  \eqref{ActionOfJHomomorphism}
  induced by the J-homomorphism (Def. \ref{TateJHomomorphism}):
  \vspace{-2mm}
  $$
    S^{V_\phi}
    \;:=\;
    \raisebox{1pt}{\textesh}
    \orbisingular
    \big(
      \mathbf{S}^V \!\sslash_{{}_{\phi}}\! G
    \big)
    \;\;\;\;\;
    \in
    \mathbf{H}_{/_{\!\scalebox{.7}{$\orbisingular G$}}}\;.
  $$
\end{defn}

\begin{example}[Ordinary representation spheres]
  \label{OrdinaryRepresentationSpheres}
  Let $\mathbf{H} = \mathrm{SingularJetsOfSmoothGroupoids}_\infty$
  (Example \ref{SingularInfinityGroupoids}) and
  $V := (\mathbb{R}^n, +)$ as in Example \ref{OrdinaryManifolds},
  whence $\mathrm{Aut}(T_e V) \simeq \mathrm{GL}(n)$
  (Example \ref{OrdinaryGeneralLinearGroup}).
    For
    \vspace{-1mm}
  $$
    G
      \xymatrix{
        \;   \ar@{^{(}->}[r]^-{\phi}
      &
}
  \mathrm{GL}(n) \subset \mathrm{Aut}(T_e V)
  $$

  \vspace{-2mm}
\noindent
  a finite subgroup, hence  a linear $G$-representation,
  we have that the representation sphere
  $S^{\mathbb{R}^n_\phi}$ according to
  Def. \ref{RepresentationSphere} is the ordinary representation sphere,
  as an object in $G$-equivariant homotopy theory.
\end{example}

\begin{defn}[$J$-twisted proper orbifold Cohomotopy theory]
  \label{JTwistedOrbifoldCohomotopyTheory}
  Let $\mathbf{H}$ be a
  singular-elastic $\infty$-topos (Def. \ref{SingularSolidInfinityTopos})
  $V \in \mathrm{Groups}(\mathbf{H})$ (Prop. \ref{LoopingAndDelooping}),
  $W \in \mathrm{LinearGroups}(\mathbf{H})$ (Def. \ref{LinearGroup})
  and
  $\phi : \!\xymatrix@C=15pt{
    \mathrm{Aut}(T_e W)
    \ar[r]
    &
    \mathrm{Aut}(T_E V)
    \,.
  }$
  Then \emph{J-twisted proper orbifold Cohomotopy}
  is the tangentially twisted proper orbifold cohomotopy (Def. \ref{OrbifoldCohomologyTangentiallyTwisted})
  with coefficients
   \vspace{-2mm}
  $$
    (A,\rho) \;:=\; (\mathbf{S}^V, \rho_{\mathbf{J}})
  $$

  \vspace{-2mm}
\noindent
  the Tate $W$-sphere (Def. \ref{TheVSphere})
  with its Tate J-homomorphism action (Def. \ref{TateJHomomorphism}):
   \vspace{-2mm}
  $$
    \overset{
      \mathclap{
      \raisebox{7pt}{
        \tiny
        \color{darkblue}
        \bf
        \begin{tabular}{c}
          J-twisted
          \\
          orbifold Cohomotopy
        \end{tabular}
      }
      }
    }{
      \pi^{\scalebox{.8}{$
        \raisebox{1pt}{\textesh} \orbisingular \tau
      $}}
      (-)
    }
    \;\;\; :=\;\;\;
    \overset{
      \mathclap{
      \raisebox{7pt}{
        \tiny
        \color{darkblue}
        \bf
        \begin{tabular}{c}
          tangentially twisted
          orbifold cohomology
        \end{tabular}
      }
      }
    }{
      H^{\scalebox{.8}{$
        \raisebox{1pt}{\textesh} \orbisingular \tau
      $}}
      \big(
        -,
        \underset{
          \mathclap{
          \raisebox{-7pt}{
            \tiny
            \color{darkblue}
            \bf
            \begin{tabular}{c}
              Tate $V$-sphere with
              \\
              $J$-homomorphism action
            \end{tabular}
          }
          }
        }{
          (\mathbf{S}^V, \rho_{\mathbf{J}})
        }
      \big).
    }
  $$

   \vspace{-2mm}
\noindent
  Hence for a structured orbifold (Def. \ref{OrbiVFolds})
   \vspace{-2mm}
  $$
    \big(
      \mathcal{X}, (\tau,g)
    \big)
    \;
    \in
    \;
    \big(
      \mathrm{Aut}(T_e W), \phi
    \big)
    \mathrm{Structured}
    V\mathrm{Orbifolds}(\mathbf{H})\;,
  $$

  \vspace{-2mm}
\noindent
  we have:
   \vspace{-2mm}
  $$
    \overset{
      \raisebox{7pt}{
        \tiny
        \color{darkblue}
        \bf
        \begin{tabular}{c}
          J-twisted
          \\
          orbifold Cohomotopy
        \end{tabular}
      }
    }{
      \pi^{\scalebox{.8}{$
        \raisebox{1pt}{\textesh} \orbisingular \tau
      $}}
      (\mathcal{X})
    }
      \;\;=\;\;
  \left\{\!\!\!\!
  \raisebox{19pt}{
  \xymatrix@C=4em{
    \overset{
      \mbox{
        \tiny
        \color{darkblue}
        \bf
        orbifold
      }
    }{
      \mathcal{X}
    }
    \ar@{-->}[rr]^-{
      \overset{
        \raisebox{5pt}{
          \tiny
          \color{darkblue}
          \bf
          \begin{tabular}{c}
            cocycle
          \end{tabular}
        }
      }{
        c
      }
    }_-{\ }="s"
    \ar[dr]_-{
      \underset{
        \mathllap{
        \mbox{
          \tiny
          \color{darkblue}
          \bf
          \begin{tabular}{c}
            tangential
            \\
            twist
          \end{tabular}
        }
        }
      }{
        \eta^{\scalebox{.6}{\textesh}}
        \circ
        \orbisingular(\tau)
      }
      \!\!\!\!
    }^-{\ }="t"
    &&
    \mbox{\textesh}
    \orbisingular
    \big(
      \overset{
        \raisebox{5pt}{
          \tiny
          \color{darkblue}
          \bf
          \begin{tabular}{c}
            orbi-singularized
            \\
            Tate $W$-sphere
          \end{tabular}
        }
      }{
        \mathbf{S}^W \!\!\sslash\! \mathrm{Aut}(T_e W)
      }
    \big)
    \ar[dl]^{\;\;
      \underset{
        \mathrlap{
        \mbox{
          \tiny
          \color{darkblue}
          \bf
          \begin{tabular}{c}
            twisting via
            \\
            orbi-singularized
            \\
            J-homomorphism
          \end{tabular}
        }
        }
      }{
      \scalebox{.7}{$
        \raisebox{1pt}{\textesh} \orbisingular
        \big(
          \rho_{\mathbf{J}}
        \big)
        \mathrlap{
          \;\;
          \mbox{
            \tiny
            \color{darkblue}
            \bf
          }
        }
      $}
      }
    }
    \\
    &
    \raisebox{1pt}{\textesh}\orbisingular
    \mathbf{B} \mathrm{Aut}(T_e W)
    \ar@{=>} "s"; "t"
  }
  }
  \!\! \right\}_{\!\!\big/\sim}
$$
\end{defn}

\begin{example}[$J$-Twisted proper orbifold Cohomotopy of ordinary orbifolds]
  \label{JTwistedOrbifoldCohomotopyOnOrdinaryOrbifolds}
  $\,$
  \\
  \noindent
  Let $\mathbf{H} = \mathrm{SingularJetsOfSmoothGroupoids}_\infty$
  (Example \ref{SingularInfinityGroupoids}) and
  $V := (\mathbb{R}^n, +)$, $W := (\mathbb{R}^p, +)$
  as in Example \ref{OrdinaryManifolds}, with $p \leq n$,
  and $\phi : (\mathbb{R}^p, + ) \hookrightarrow (\mathbb{R}^n, +)$
  be the canonical inclusion.
  Then
  the corresponding J-twisted proper orbifold Cohomotopy theory
  $
      \pi^{\scalebox{.8}{$
        \raisebox{1pt}{\textesh} \orbisingular \tau
      $}}$
  (Def. \ref{JTwistedOrbifoldCohomotopyTheory})
  is defined on ordinary $n$-dimensional orbifolds
  (by Example \ref{EtaleLieGroupoidAsRnFold})
  with
  $\mathrm{GL}(p)$-structure (by Example \ref{OrdinaryGeneralLinearGroup})
  and
  it unifies the following two special cases
  (by Theorem \ref{OrbifoldCohomologyTangentiallyTwistedReduces},
  see the second diagram on p. \pageref{JTwistedOrbifoldCohomotopyDiagram})):

 \noindent  {\bf (i)} On smooth orbifolds, i.e., on ordinary manifolds
  (Example \ref{OrdinaryManifolds}) it reduces to
  non-abelian cohomology with coefficients the shape of the
  ordinary
  $p$-sphere (by Example \ref{TateSpheresWithShapeOfOrdinarySpheres})
  and tangentially twisted
  via the traditional J-homomorphism (by Example \ref{OrdinaryJHomomorphism}).
  This is the \emph{J-twisted Cohomotopy theory}
  considered in \cite{FSS19b}\cite{FSS19c} \cite{SS19b}.

\vspace{-1mm}
\noindent  {\bf (ii)} On flat orbifolds, such as the vicinity
  of ordinary
  orbifold singularities $\mathbb{R}^p \!\sslash\! G$
  for finite subgroups $G \overset{\phi}{\hookrightarrow} \mathrm{GL}(p)$
  (by Example \ref{OrdinaryOrbifoldSingularities}),
  hence for linear $G$-representations $\phi$,
  it reduces to proper equivariant cohomology in
  $\mathrm{RO}$-degree $\phi$
  and with
  coefficients the representation sphere $S^{\mathbb{R}^n_\phi}$
  (by Example \ref{OrdinaryRepresentationSpheres}).
  This is the \emph{tangentially RO-graded equivariant Cohomotopy theory}
  considered in \cite{SS19a}\cite{SS19b}.
\end{example}

By way of conclusion and outlook, we highlight the following:

\begin{remark}[Orbifold cohomology in non-perturbative string theory and
{\it Hypothesis H}]
  \label{HypothesisH}
  Traditional discussion of orbifold cohomology has been strongly
  motivated by its application to \emph{perturbative string theory}
  (e.g. \cite{AMR02}\cite{ARZ06}\cite{ALR07}\cite{BecerraUribe09}\cite{DFM11}).
  However, perturbative string theory is famously in need of
  a non-perturbative completion (``M-theory'', see \cite[2]{HSS18}\cite{FSS19} for
  review and pointers) whose mathematical formulation has remained
  an open problem. Therefore, it is to be expected that the
  historically rich interaction between orbifold cohomology theory and string theory
  is just the tip of an iceberg, whose full scope is
  a cohomology theory of M-theoretic orbifolds.

  Elsewhere we have put forward a precise hypothesis
  as to what this mathematical theory should be.
  This \emph{Hypothesis H} says that:
  \\
  \noindent {\bf (i)} far from singularities, M-theory is controlled by
  twisted Cohomotopy theory \cite{FSS19b}\cite{FSS19c}\cite{SS19b}\cite{FSS20a};

  \noindent {\bf (ii)} at singularities, M-theory is controlled by
  RO-graded equivariant Cohomotopy theory \cite{HSS18}\cite{SS19a}\cite{SS19b}.

  \noindent
  (See these references for various consistency checks of this
  hypothesis.)

  The impact of Theorem \ref{OrbifoldCohomologyTangentiallyTwistedReduces},
  in its specialization to Example \ref{JTwistedOrbifoldCohomotopyOnOrdinaryOrbifolds},
  is to show that these two cases are indeed two aspects of
  a single unified cohomology theory: J-twisted proper orbifold Cohomotopy theory.

\end{remark}

\medskip

\newpage


\begin{appendices}

\section{Model category presentations}
\label{ModelCategoryPresentations}

We recall some basics of model categories
(e.g. \cite[2]{GoerssJardine99})
of simplicial presheaves (\cite{Jardine87}\cite{Jardine96}\cite{Jardine15})
as presentations of $\infty$-toposes (\cite[A.2, A.3]{Lurie09}).

\medskip

\noindent {\bf Model categories of simplicial presheaves.}
\begin{defn}[Model category of simplicial presheaves]
\label{SimplicialPresheaves}
Let $\mathcal{C}$ be a site. We write

\noindent {\bf (i)}
\vspace{-2mm}
\begin{equation}
  \label{CategoryOfSimplicialPresheavesAsHomotopical}
  \mathrm{sPSh}(\mathcal{C})_{\mathrm{loc}}
  \;\;\in\;
  \mathrm{HomotopicalCategories}
\end{equation}
for the category of simplicial presheaves on $\mathcal{C}$,
regarded as a homotopical category with weak equivalences
the local weak homotopy equivalences of simplicial sets.

\noindent {\bf (ii)}
\vspace{-2mm}
\begin{equation}
  \label{ModelCategoryOfSimplicialPresheaves}
  \mathrm{sPSh}(\mathcal{C})_{
    {
      {\mathrm{inj}/}
      \atop
      {\mathrm{proj}}
    },
    \mathrm{loc}
  }
  \;\;\in\;
  \mathrm{ModelCategories}
\end{equation}
for the same category regarded as either the
corresponding injective or projective model category.

\noindent {\bf (iii)}
\vspace{-2mm}
\begin{equation}
  \label{SimplicialLocalizationOfSimplicialPresheaves}
  \xymatrix{
    \mathrm{sPSh}
    \ar[r]^-{ \ell }
    &
    L_{{}_{\mathrm{lwhe}}} \mathrm{sPSh}_{\mathrm{loc}}
    \;=:\;
    \mathbf{H}
  }
\end{equation}
for the corresponding simplicial localization.
\end{defn}

\begin{lemma}[Cofibrancy in projective model structure {\cite[Cor. 9.4]{Dugger01}}]
  \label{CofibrancyInTheProjectiveModelStructure}
  Let $\mathcal{C}$ be a site.
  For a simplicial presheaf
  $X_\bullet \in \mathrm{sPSh}(\mathcal{C})_{\mathrm{proj},\mathrm{loc}}$
  in the projective model structure \eqref{ModelCategoryOfSimplicialPresheaves}
  to be cofibrant it is sufficient that $X_\bullet$
  is degreewise

\noindent {\bf (i)} a coproduct of representables, such that

\noindent  {\bf (ii)} the degenerate cells split off as a direct summand.
\end{lemma}

\begin{lemma}[Simplicial presheaf represents its own hocolim {\cite[2.1]{DuggerHollanderIsaksen04}\cite[2.3.21]{dcct}}]
  \label{SimplicialPresheafRepresentsItsOwnHocolim}
  Let $\mathcal{C}$ be a site
  and $X_\bullet \in \mathrm{sPSh}(\mathcal{C})$
  a simplicial presheaf (Def. \ref{SimplicialPresheaves}).
  Then its image under simplicial localization \eqref{SimplicialLocalizationOfSimplicialPresheaves}
  is equivalently the simplicial homotopy colimit over the
  images of its component presheaves:
  \vspace{-2mm}
  $$
    \ell(X_\bullet)
    \;\simeq\;
    \underset{\longrightarrow}{\mathrm{lim}}
    (\ell X)_\bullet
    \;\;
    \in
    \;
    \mathbf{H}
    \,.
  $$
\end{lemma}

\medskip

\noindent {\bf Topological mapping stacks}

\begin{example}[Model category presentation of smooth $\infty$-groupoids]
  \label{ModelCategoryPresentationOfSmoothInfinityGroupoids}
  Let $\mathcal{C} = \mathrm{CartesianSpaces}$ (Def. \ref{CartesianSpaces}).
  Then the simplicial localization \eqref{SimplicialLocalizationOfSimplicialPresheaves}
  of $\mathrm{sPSh}(\mathcal{C})_{\mathrm{loc}}$ \eqref{ModelCategoryOfSimplicialPresheaves}
  is $\mathrm{SmoothGroupoids}_\infty$ (Example \ref{SmoothInfinityGroupoids}):
\vspace{-2mm}
  $$
    L_{{}_{\mathrm{lwhe}}}
    \mathrm{sPSh}(\mathrm{CartesianSpaces})_{\mathrm{loc}}
    \;\simeq\;
    \mathrm{SmoothGroupoids}_\infty\;.
  $$
\end{example}

\begin{lemma}[Mapping stack from delooping of discrete group to topological stack]
\label{MappingStackFromDiscreteGroupToTopologicalStack}
  In $\mathrm{SmoothGroupoids}_\infty$ (Example \ref{SmoothInfinityGroupoids})
  consider

  \noindent
  {\bf (i)} a finite group embedded via \eqref{InclusionOfFiniteGroups}
  \vspace{-2mm}
  \begin{equation}
    \label{ADiscreteGroupForHommingIntoATopologicalGroupoid}
    G \;\in\;
    \xymatrix{
      \mathrm{Groups}
      \ar[r]^-{ \mathrm{Disc} }
      &
      \mathrm{Groups}
      \big(
        \mathrm{SmoothGroupoids}_\infty
      \big),
    }
  \end{equation}

\vspace{-2mm}
  \noindent {\bf (ii)} a topological groupoid,
  embedded via \eqref{1TrunConcreteSmoothInfinityGroupoidsAreDiffeologicalSpaces}
  \vspace{-2mm}
  \begin{equation}
    \label{ATopologicalGroupoidToBeHommedIntoFromADiscreteGroupDelooping}
    \xymatrix@R=1pt{
      \mathrm{TopologicalGroupoids}
      \ar[r]^-{ \scalebox{.6}{$\mathrm{Cdfflg}$} }
      &
      \mathrm{SmoothGroupoids}_\infty
      \\
      \mathcal{X}_{\mathrm{top}}
      \ar@{}[r]|-{\mapsto}
      &
      \mathcal{X}_{\tiny\smooth}
    }
  \end{equation}

\vspace{-2mm}
  \noindent
  Then the mapping stack \eqref{InternalHomAdjunction} formed in
  $\mathrm{SmoothGroupoids}_\infty$ is the
  degreewise image under
  $\mathrm{Cdfflg}$ \eqref{AdjunctionBetweenTopologicalAndDiffeologicalSpaces}
  of the topological groupoid representing the
  mapping stack of topological groupoids
  (which exists by \cite{Noohi08} since $G$ is finite, hence compact):
  \begin{equation}
    \label{MappingStackFromDeloopingOfDiscreteToTopologicalStack}
    \mathbf{Maps}
    \big(
      \mathbf{B}G
      \,,\,
      \mathcal{X}_{\tiny\smooth}
    \big)
    \;\simeq\;
    \mathrm{Cdfflg}
    \,
    \mathbf{Maps}
    \big(
      \mathbf{B}G
      \,,\,
      \mathcal{X}_{\mathrm{top}}
    \big).
  \end{equation}
\end{lemma}

\begin{proof}
Since (by Example \ref{SmoothInfinityGroupoids})

\vspace{-.8cm}
$$
  \mathrm{SmoothGroupoids}_\infty
  \;\simeq\;
  \xymatrix{
    \mathrm{Sheaves}_\infty(\mathrm{CartesianSpaces})
    \;
    \ar@{<-}@<+5pt>[rr]^-{ L }
    \ar@<-5pt>@{^{(}->}[rr]^-{
      \raisebox{-1pt}{
        \scalebox{.8}{$
          \bot
        $}
      }
    }
    &&
   \; \mathrm{PreSheaves}_\infty(\mathrm{CartesianSpaces})
  }
$$

\vspace{-1mm}
\noindent
it is sufficient to show that we have an equivalence
of $\infty$-presheaves of the form
\begin{equation}
  \label{AnEquivalenceOfInfinityPresheaves}
  \mathbb{R}^n
    \xymatrix{
        \;   \ar@{|->}[r]
      &
} \;\;\;
  \begin{aligned}
    & \mathrm{PreSheaves}_\infty
    \big(\mathrm{CartesianSpaces}\big)
    ( \mathbb{R}^n \times \mathbf{B}G\,,\, \mathcal{X}_{\tiny\smooth} )
    \\
    &
    \simeq\;
    \mathrm{PreSheaves}_\infty
    \big(\mathrm{CartesianSpaces}\big)
    (
      \mathbb{R}^n \,,\,
      \mathrm{Cdfflg}
      \,
      \mathbf{Maps}(\mathbf{B}G, \mathcal{X}_{\mathrm{top}})
    )
  \end{aligned}
\end{equation}
By Example \ref{ModelCategoryPresentationOfSmoothInfinityGroupoids},
we may model this in the global projective model structure on
simplicial presheaves over $\mathrm{CartesianSpaces}$:

\vspace{-.7cm}

\begin{equation}
  \label{GlobalSimplicialPresheafLocalization}
  \xymatrix{
    \mathrm{sPSh}(\mathrm{CartesianSpaces})
    \ar[r]^-{ \ell }
    &
    \;\;
    L_{{}_{\mathrm{lwhe}}}
    \mathrm{sPSh}(\mathrm{CartesianSpaces})_{\mathrm{proj}}
    \;\simeq\;
    \mathrm{PreSheaves}_\infty(\mathrm{CartesianSpaces})
  }
\end{equation}

 \vspace{-2mm}
\noindent
with, by Lemma \ref{SimplicialPresheafRepresentsItsOwnHocolim},
the following models:

\noindent
{\bf (a)}
A model under $\ell$ \eqref{GlobalSimplicialPresheafLocalization} of
the Cartesian product
$\mathbb{R}^n \times \mathbf{B}G$
with the delooping $\mathbf{B}G \simeq
\underset{\longrightarrow}{\mathrm{lim}}\, G^{\times_\bullet}$
\eqref{DeloopingAsHomotopyColimit},
is given by the simplicial presheaf
\begin{equation}
  \label{RnBGAsSimplicialPresheaf}
  \mathbb{R}^n \times G^{\times_\bullet}
  \;\in\;
  \mathrm{sPSh}(\mathrm{CartesianSpaces})_{\mathrm{proj}}\;.
\end{equation}

\noindent
{\bf (b)} A model under $\ell$ \eqref{GlobalSimplicialPresheafLocalization}
for the image \eqref{ATopologicalGroupoidToBeHommedIntoFromADiscreteGroupDelooping}
of a topological groupoid $\mathcal{X}_{\mathrm{top}}$ is given by
its nerve
regarded
as a simplicial presheaf, componentwise via \eqref{ConcreteSmoothInfinityGroupoidsAreDiffeologicalSpaces}
\vspace{-2mm}
\begin{equation}
  \label{TopologicalGroupoidAsSimplicialPrsheaf}
  N_\bullet(\mathcal{X}_{\mathrm{top}})
  \;\in\;
  \mathrm{sPSh}(\mathrm{CartesianSpaces})_{\mathrm{proj}}\;.
\end{equation}

\vspace{-2mm}
\noindent Moreover:
\begin{itemize}
\vspace{-2mm}
\item The object \eqref{RnBGAsSimplicialPresheaf} is projectively cofibrant,
by Lemma \ref{CofibrancyInTheProjectiveModelStructure}, as
is its Cartesian product with an $k$ simplex $\Delta[k]$.

\vspace{-2mm}
\item  The object \eqref{TopologicalGroupoidAsSimplicialPrsheaf} is
projectively fibrant (objectwise a Kan complex) by the groupoid property
of $\mathcal{X}_{\mathrm{top}}$.
\end{itemize}

\vspace{-2mm}
\noindent Therefore to get \eqref{AnEquivalenceOfInfinityPresheaves}
it is, in turn, sufficient to exhibit
for $\mathbb{R}^n \in \mathrm{CartesianSpaces}$ a natural isomorphism
of simplicial sets of the form
\begin{equation}
  \label{AMorphismOfSimplicialPresheaves}
    \int_{[k] \in \Delta}
    \mathrm{PSh}
    \big(
      \mathbb{R}^n
        \times
      \big(
        G^{\times_k}
        \times
        \Delta(k,\bullet)
      \big)
      \,,\,
      \mathrm{Cdfflg}(
        N_k(\mathcal{X}_{\mathrm{top}})
      )
    \big)
    \;\simeq\;
    \mathrm{PSh}
    \left(
      \mathbb{R}^n
      \,,\,
      \mathrm{Cdfflg}
      \left(
       \int_{[k] \in \Delta}
        N_k(\mathcal{X}_{\mathrm{top}})^{
          (
            G^{\times_k}
            \times
            \Delta(k,\bullet)
          )
        }
      \right)
    \right),
\end{equation}
where the \emph{end} $\int_{[k] \in \Delta} (-)$
expresses the limit that computes the morphism of simplicial
sets as a subset of the product of the function spaces of components.
We obtain this as the following composite of
natural isomorphisms:
$$
  \begin{aligned}
        \int_{[k] \in \Delta}
    \mathrm{PSh}
    \Big(
      \mathbb{R}^n
        \times
      \big(
        G^{\times_k}
        \times
        \Delta(k,\bullet)
      \big)
      \,,\,
      \mathrm{Cdfflg}\big(
        N_k(\mathcal{X}_{\mathrm{top}})
      \big)
    \Big)
       & \simeq
    \int_{[k] \in \Delta}
    \mathrm{PSh}
    \left(
      \mathbb{R}^n
      \,,\,
      \big(
        \mathrm{Cdfflg}(
          N_k(\mathcal{X}_{\mathrm{top}})
        )
      \big)^{
        (
          G^{\times_k} \times \Delta(k,\bullet)
        )
      }
    \right)
    \\
    & \simeq
    \int_{[k] \in \Delta}
    \mathrm{PSh}
    \left(
      \mathbb{R}^n
      \,,\,
      \mathrm{Cdfflg}
      \left(
        {
          \big(
            N_k(\mathcal{X}_{\mathrm{top}})
          \big)
        }^{
          (
            G^{\times_k}
            \times
            \Delta(k,\bullet)
          )
        }
      \right)
    \right)
    \\
       & \simeq
    \mathrm{PSh}
    \left(
      \mathbb{R}^n
      \,,\,
      \int_{[k] \in \Delta}
      \mathrm{Cdfflg}
      \left(
        \big(
          N_k(\mathcal{X}_{\mathrm{top}})
        \big)^{
          (
            G^{\times_k}
            \times
            \Delta(k,\bullet)
          )
        }
      \right)
    \right)
    \\
    & \simeq
    \mathrm{PSh}
    \left(
      \mathbb{R}^n
      \,,\,
      \mathrm{Cdfflg}
      \left(
       \int_{[k] \in \Delta}
        \big(
          N_k(\mathcal{X}_{\mathrm{top}})
        \big)^{
          (
            G^{\times_k}
            \times
            \Delta(k,\bullet)
          )
        }
      \right)
    \right).
  \end{aligned}
$$

\vspace{0mm}
\noindent
Here the first step is the definition of function spaces $(-)^{(-)}$,
the second step uses that $\mathrm{Cdfflg}$, being a right adjoint,
preserves products (Prop. \ref{AdjointsPreserveCoLimits}).
The third step uses that the Hom-functor
preserves limits (hence ends) in its second argument,
while the fourth step uses that $\mathrm{Cdfflg}$, being a right adjoint,
preserves limits (hence ends), again by Prop. \ref{AdjointsPreserveCoLimits}.
\hfill \end{proof}

\section{Equivariant homotopy theory}
\label{GEquivariantHomotopyTheory}

For reference, we recall some basics of unstable equivariant homotopy theory
(see \cite{May96}\cite{Blu17}). We focus here on finite groups,
for simplicity and since this is what we need in the main text
(Remark \ref{NeedForDiscreteGroupsInSingularities}),
but all statements in the following, notably
Elmendorf's theorem (Prop. \ref{ElmendorfTheorem} below)
generalize to compact Lie groups.

\medskip

\begin{defn}[Topological $G$-spaces]
  \label{GSpaces}
 Let $G \in \mathrm{Groups}^{\mathrm{fin}}$ be a finite group.
\vspace{-1mm}
 \item {\bf (i)} We write
  \begin{equation}
    \label{GTopologicalSpaces}
    G \mathrm{DTopologicalSpaces}
     \xymatrix{\ar@{^{(}->}[r]&}
    G \mathrm{TopologicalSpaces}
    \;\in\; \mathrm{Categories}
  \end{equation}
  for the categories whose objects
  \begin{equation}
    \label{GAction}
    G\acts X
    \;:=\;
    \big(
      X,
      \,
      G \times X
      \overset{
        \rho
      }{\longrightarrow}
      X
    \big)
  \end{equation}
  are
  topological spaces $X$ (as in Def. \ref{TopologicalSpaces})
  or specifically
  D-topological spaces (as in Def. \ref{DTopologicalSpaces}),
  respectively,
  equipped with continuous left $G$-actions $\rho$,
  and whose morphisms are the $G$-equivariant continuous functions:
  \begin{equation}
    \label{HomSetOfGEquivariantMaps}
    \hspace{-1cm}
    G\mathrm{TopologicalSpaces}
    \big(
      G\acts X_1
      ,
      \;
      G \acts X_2
    \big)
    \;:=\;
    \left\{
    \xymatrix@C=3.5em{
      X_1
      \ar[r]_-{\rm   \scalebox{.6}{
          continuous}
        }^f
        & X_2
      }
          \;
      \left\vert
      \;
      \raisebox{16pt}{
      \xymatrix@R=1.2em{
        G \times X_1
          \ar[d]_-{f}
          \ar[rr]^-{ \rho_1 }
          &&
        X_1
          \ar[d]_-{f}
        \\
        G \times X_2
          \ar[rr]_-{ \rho_2 }
          && X_2
      }
      }
      \;
      \right.
    \right\}.
  \end{equation}

\vspace{-3mm}
 \item {\bf (ii)}  For $G\acts X_1$ a (D-)topological $G$-space
 and $H \overset{\iota}{\hookrightarrow} G$ a subgroup, we write
  \begin{equation}
    \label{HFixedPointSpaces}
    X^H
    \;:=\;
    \Big\{
      x \in X
      \;\vert\;
      \underset{h \in H \subset G}{\forall} \;\;
      \rho(h,x)
      =
      x
    \Big\}
  \end{equation}

  \vspace{-1mm}
  \noindent
  for the topological subspace of $H$-fixed points
  (which, if $X$ is D-topological, is itself
  again D-topological, by Prop. \ref{NiceDTopologicalSpaces}).

  \item {\bf (iii)} For
  $G \acts X_1$ and $G \acts X_2$
  two (D-)topological $G$-spaces, the mapping space \eqref{InternalHomInDTopologicalSpaces}
  between their underlying (D-)topological spaces canonically becomes a
  $G$-space
  via the conjugation action
  and the corresponding fixed point space \eqref{HFixedPointSpaces}
  \vspace{-2mm}
  \begin{equation}
    \label{MappingSpaceOfEquivariantFunctions}
    \xymatrix{
    \mathrm{Maps}(X_1,X_2)^G \; \ar@{^{(}->}[r]& \mathrm{Maps}(X_1,X_2)
    }
  \end{equation}

 \vspace{-1mm}
\noindent
  is the subspace on the $G$-equivariant functions
  \eqref{HomSetOfGEquivariantMaps}.
\end{defn}

\vspace{1mm}
\begin{example}[$G$-cells]
For $G \in \mathrm{Grp}_{\mathrm{fin}}$,
  $H \subset G$ a subgroup
  and $n \in \mathbb{N}$ we have the $G$-spaces
  (Def. \ref{GSpaces})
  $$
    \big(G/H \big) \times D^n
    \;\,,\;
    \big(G/H \big) \times S^{n-1}
    \;\in\;
    G \mathrm{DTopologicalSpaces}
  $$
  being the product spaces of the discrete orbit spaces
  with the standard topological unit disk and unit circle,
  respectively,
  the latter equipped with the trivial $G$-action.
  The boundary inclusions
  $
      \partial D^n = S^{n-1} \overset{\iota_n}{\hookrightarrow}  D^n
    $
  induce $G$-equivariant maps
  \vspace{-2mm}
  \begin{equation}
    \label{GCellBoundaryInclusion}
    \iota_{n,H}
    \;\colon\;
    \xymatrix{
      \big(G/H\big) \times S^{n-1}
      \; \ar@{^{(}->}[rr]^{ (\mathrm{id}, \iota_n) }
      &&
      \big(G/H\big) \times D^n
    }
  \end{equation}
  for all $n \in \mathbb{N}$, $H \subset G$.
\end{example}

\medskip

\begin{defn}[$G$-CW-complexes]
  \label{GCWComplexes}
 $ \phantom{A}$
  \vspace{-1mm}
 \item {\bf (i)}   A \emph{$G$-CW-complex} $X$ is a
 D-topological $G$-space
 (Def. \ref{GSpaces})
  which is equipped with the realization as a colimit
\vspace{-1mm}
 $$
    X
      \;\simeq\;
    \underset{\longrightarrow_{\mathrlap{n}}}{\lim}
    X_n
    \;\;
    \in
    G \mathrm{DTopologicalSpaces}
  $$

 \vspace{-3mm}
\noindent  over a sequence
  $$
    X_{-1} \longrightarrow X_0 \longrightarrow X_1 \longrightarrow X_2
    \longrightarrow \cdots
    \;\;\;
    \in
    G \mathrm{DTopologicalSpaces}
  $$
  where
  $
    X_{-1} = \varnothing
  $
  and where each $X_n \to X_{n-1}$ is given by a
  set of attachments of $G$-cells along \eqref{GCellBoundaryInclusion},
  hence by a pushout of the form:

  \newpage

 \vspace{-2mm}
  $$
    \xymatrix@R=.6em@C=3em{
      \underset{
        {H \subset G}
        \atop
        {i \in I_n}
      }{\coprod}
      G/H \times S^{n-1}
      \ar@{^{(}->}[dd]_-{ (\iota_{n,H})_{n,H} }
      \ar[rr]
      \ar@{}[ddrr]|{ \mbox{\tiny (po)} }
      &&
      X_{n-1}
      \ar@{^{(}->}[dd]
      \\
      \\
      \underset{
        {H \subset G}
        \atop
        {i \in I_n}
      }{\coprod}
      G/H \times D^n
      \ar[rr]
      &&
      X_{n}
    }
  $$

\vspace{-3mm}
\item {\bf (ii)} We write
  \begin{equation}
    \label{GSets}
    \xymatrix{
    G \mathrm{Sets}
      \; \ar@{^{(}->}[r] & \;
    G \mathrm{CWComplexes}
      \; \ar@{^{(}->}[r] & \;
    G \mathrm{DTopologicalSpaces}
    }
  \end{equation}
  for the full subcategories on those
  D-topological $G$-spaces which
  admit the structure of $G$-CW-complexes.
\end{defn}

  \begin{defn}[Homotopy theory of D-topological $G$-spaces]
    \label{HomotopyTheoryOfGSpaces}
  The \emph{homotopy theory of topological $G$-spaces}
  is the $\infty$-category
  \begin{equation}
    \label{InfinityCategoryOfGInfinityGroupoids}
    G \mathrm{Groupoids}_{\infty}
    \;\in\;
    \mathrm{Categories}_\infty
  \end{equation}
  which has the same objects as $G\mathrm{CWComplexes}$
  (Def. \ref{GCWComplexes}), and with $\infty$-groupoids the
  topological shapes (Def. \ref{ShapeOfTopologicalSpaces})
  of the mapping spaces \eqref{MappingSpaceOfEquivariantFunctions} of $G$-equivariant maps:
  \vspace{-2mm}
  \begin{equation}
    \label{GEquivariantMappingSpaces}
    G \mathrm{Groupoids}_\infty\big( G \acts X_1, G \acts X_2 \big)
    \;:=\;
    \mathrm{Shp}_{\mathrm{Top}}
    \left(
      \mathrm{Maps}\big( X_1, X_2 \big)^G
    \right)
    \,.
  \end{equation}
  \end{defn}

\begin{defn}[Shape of $G$-topological spaces]

 \item {\bf (i)}  We write

  \vspace{-.8cm}

  \begin{equation}
    \label{ShapeOfGCWComplexes}
    \mathrm{Shp}_{G\mathrm{Top}}
    \;:\;
    \xymatrix{
      G \mathrm{CWComplexes}
      \ar[rr]^-{  }
      &&
      G \mathrm{Groupoids}_\infty
    }
  \end{equation}
  for the canonical $\infty$-functor (topologically enriched functor)
  from the 1-category of $G$-CW-complexes (Def. \ref{GCWComplexes})
  to the $\infty$-category of $G$-$\infty$-groupoids
  (Def. \ref{HomotopyTheoryOfGSpaces}),
  which is the identity on objects and which on Hom-spaces is the
  continuous map given by the identity fuction
  from the discrete set of $G$-equivariant maps
  \eqref{HomSetOfGEquivariantMaps} to
  the topological space of $G$-equivariant maps \eqref{GEquivariantMappingSpaces}.

 \item {\bf (ii)}    For any choice of $G$-CW-approximation functor
 \vspace{-3mm}
  $$
    \xymatrix{
      G \mathrm{TopologicalSpaces}
      \ar[rr]^-{\;
        (-)_{\mathrm{cof}}
      \;}
      &&
      G \mathrm{CWComplex}
    }
  $$
  we get the corresponding shape functor on
  all of $G \mathrm{TopologicalSpaces}$ (Def. \ref{GSpaces})
  and hence on
  $G \mathrm{DTopologicalspaces}$, which we denote by the same
  symbol:
  \vspace{-3mm}
  \begin{equation}
    \label{ShapeOfGTopologicalSpaces}
    \mathrm{Shp}_{G\mathrm{Top}}
    \;:\;
    \xymatrix{
      G\mathrm{TopologicalSpaces}
      \ar[rr]^-{ (-)_{\mathrm{cof}} }
      &&
      G \mathrm{CWComplexes}
      \ar[rr]^-{ \mathrm{Shp}_{G \mathrm{top}} }
      &&
      G\mathrm{Groupoids}_\infty\;.
    }
  \end{equation}
\end{defn}

\begin{defn}[Proper $G$-equivariant generalized cohomology of topological $G$-spaces]
  \label{GenuineGEquivariantCohomology}
  For $G \in \mathrm{Groups}^{\mathrm{fin}}$, we say that the
  \emph{proper $G$-equivariant cohomology}
  of a topological $G$-space (Def. \ref{GSpaces})
  $
    X \;\in\; G\mathrm{TopologicalSpaces}
  $
  with coefficients in a (pointed) $G$-$\infty$-groupoid (Def. \ref{HomotopyTheoryOfGSpaces}),
  $
    A \in  G\mathrm{Groupoids}_\infty
  $,
  is
   \vspace{-2mm}
  $$
    H^{-n}_G(X, A)
    \;:=\;
    \pi_n
    \Big(
      G\mathrm{Groupoids}
      \big(
        \mathrm{Shp}_{G\mathrm{Top}}(X)
        \,,\,
        A
      \big)
    \Big)
    \,,
  $$

   \vspace{-2mm}
\noindent
  where on the right we have the $n$th homotopy group
  (at the given basepoint) of the hom-$\infty$-groupoid
  \eqref{GEquivariantMappingSpaces}
  from the $G$-topological shape of $X$ \eqref{ShapeOfGTopologicalSpaces}
  to $A$.
\end{defn}

\medskip

\noindent {\bf Elmendorf's theorem.}
\begin{defn}[Orbit of action of a finite group]
   Let $G$ be a finite group.
  If $G \acts S$ is a set
  equipped with an action by $G$, then an \emph{orbit} of $G$ in $S$ is a subset
  of points $\{g(s) | g \in G\} \subset S$
  obtained from any single point $s \in S$ by acting on it with all elements of $G$.
  \end{defn}
 \begin{defn}[Orbit category of a finite group]
  \label{OrdinaryOrbitCategory}
The \emph{category of $G$-orbits} or \emph{orbit category of $G$}
 \vspace{-1mm}
  $$
    G \mathrm{Orbits}
    \longhookrightarrow
    G\mathrm{Sets}
    \;\in\;
    \mathrm{Categories}
  $$

   \vspace{-2mm}
  \noindent
  is the category whose objects correspond to subgroup inclusions
  $H \overset{\iota}{\longhookrightarrow} G$
  and whose morphisms are $G$-equivariant functions, hence morphisms of
  $G$-sets \eqref{GSets}, between the corresponding coset spaces
  $G/H_1 \longrightarrow G/H_2$.
\end{defn}

\newpage

\begin{example}[Systems of fixed point spaces]
  \label{SystemsOfFixedPointSpaces}
  Consider a topological space equipped with a $G$-action
  $G \acts X \;\in\; G \mathrm{DTopologicalSpaces}$
  (Def. \ref{GSpaces}) and $H \subset G$ a subgroup.
  Then a $G$-equivariant function
  $
    G/H \overset{f}{\longrightarrow} X
  $
  from the corresponding $G$-orbit (Def. \ref{OrdinaryOrbitCategory}) is
  determined by its image $f\big([e]\big) \in X$ of the class of the neutral element,
  and that image has to be fixed by the action of $H \subset G$ of $X$. Therefore,
  the corresponding $G$-equivariant mapping spaces
  \eqref{MappingSpaceOfEquivariantFunctions}
  \vspace{-2mm}
  $$
    \mathrm{Maps}\big( G/H,\, X\big)^G
    \;\simeq\;
    X^H \;:=\;
    \left\{
      x \in X \,\vert\, \underset{h \in H \subset G}{\forall} (h(x) = x)
    \right\}
    \;\subset X\;
  $$
  are the topological subspaces of $H$-fixed points inside $X$ \eqref{HFixedPointSpaces}.
  By functoriality of the mapping space construction, these fixed point spaces are exhibited
  as arranging into a topological presheaf on the $G$-orbit category (Def. \ref{OrdinaryOrbitCategory}):
   \vspace{-2mm}
  $$
    \xymatrix@R=8pt{
      X^{(-)}
      \;:\;
      &
      G \mathrm{Orbits}^{\mathrm{op}}
      \ar[rr]^-{\mathrm{Maps}(-,\, X)^G}
      &&
      \mathrm{TopologicalSpaces}
    }
  $$
\end{example}

\begin{prop}[Elmendorf's theorem {\cite{Elmendorf83}, see \cite[Thm. 1.3.6 and 1.3.8]{Blu17}}]
  \label{ElmendorfTheorem}
  Let $G$ be a finite group. The functor which sends a $G$-space $G \acts X$ (Def. \ref{GSpaces})
  to its system of $H$-fixed point spaces (Example \ref{SystemsOfFixedPointSpaces})
  constitutes an equivalence of $\infty$-categories
  \vspace{-5mm}
  \begin{equation}
    \label{FixedPointsForElmendorf}
    \hspace{1cm}
    \xymatrix@R=-2pt{
      G \mathrm{Groupoids}_\infty
      \ar[rr]^-{\simeq}
      &&
      \mathrm{Sheaves}_\infty\big( G\mathrm{Orbits}  \big)
      \\
      G \acts X
      \ar@{|->}[rr]
      &&
      X^{(-)} = \mathrm{Maps}\big( -, X\big)^G
    }
  \end{equation}
\end{prop}

\end{appendices}

\vfill

\noindent {\large \bf Acknowledgements.}
The authors would like to thank
Vincent Braunack-Mayer, David Carchedi,
Dmitri Pavlov, Charles Rezk, and David Roberts
for useful discussions.

\medskip

\vspace{.5cm}
\noindent Hisham Sati, {\it Mathematics, Division of Science, New York University Abu Dhabi, UAE.}

 \medskip
\noindent Urs Schreiber,  {\it Mathematics, Division of Science, New York University Abu Dhabi, UAE;
 on leave from Czech Academy of Science, Prague.}

\newpage


\end{document}